\documentclass[envcountsame]{svproc}

\usepackage{amsmath}
\usepackage{amssymb}
\usepackage{graphicx}
\usepackage{filecontents}  
\usepackage{mathrsfs}
\usepackage{cite}

\newif\ifmytoc
\mytoctrue

\DeclareFontFamily{OT1}{fraktura}{}
\DeclareFontShape{OT1}{fraktura}{m}{n} {<5> <6> <7> <8> <9> <10>
<11> <12> <13> <14.4> [1.1] eufm10}{}
\DeclareMathAlphabet{\fraktura}{OT1}{fraktura}{m}{n}

\DeclareFontFamily{OT1}{cmss}{} \DeclareFontShape{OT1}{cmss}{m}{n}
{<5> <6> <7> <8> <9> <10> <11> <12> <13> <14.4> cmss10}{}
\DeclareMathAlphabet{\cmss}{OT1}{cmss}{m}{n}

\newcommand\MBtableofcontentsfilename{proceedings-version.mbtoc}

\newcommand\mytocentry[1]{%
\newwrite\file
\immediate\openout\file=\MBtableofcontentsfilename
\immediate\write\file{#1}
\closeout\file
}

\newcommand\mytocfilestring{}

\newcommand\mytableofcontents{%
\centerline{\bf\large Contents}
\vskip0.5cm
\IfFileExists{\MBtableofcontentsfilename}{%
\newcommand\odskok{6mm}
\newcommand\odskocek{3mm}
\newcommand\mychaptertocentry[2]{\smallskip\hbox to \textwidth{\hbox to \odskok{\bf #1}{\bf #2}\dotfill \pageref{MBchapter-label#1}}\mytoclineskip}
\newcommand\mysectiontocentry[2]{\hbox to \textwidth{\hskip \odskok#1\hskip\odskocek #2\dotfill \pageref{MBsection-label#1}}\mytoclineskip}
\newcommand\myspecialtocentry[2]{\hbox to \textwidth{{\bf #2}\dotfill \pageref{MBspecialtocentry-label#1}}\mytoclineskip}

\input{\MBtableofcontentsfilename}
\let\mychaptertocentry\relax
\let\mysectiontocentry\relax
\let\myspecialtocentry\relax

}{\typeout{File "\MBtableofcontentsfilename" does not exist!}}
}

\renewcommand{\chapter}{\secdef\sct\sect}
\newcommand{\sct}[2][default]{\refstepcounter{chapter}
\vskip0.6cm
\vbox{\hrule
\vskip0.3cm
\centerline{\bfseries\large Lecture~\arabic{chapter}}
\label{MBchapter-label\thechapter}
\vskip0.2cm
\centerline{\bfseries\large #1}
\vskip0.3cm\hrule
}\nopagebreak\vskip0.5cm\nopagebreak
\addcontentsline{toc}{chapter}{#1}
\edef\mytocfilestring{\expandafter\mytocfilestring\csname mychaptertocentry\endcsname{\thechapter}{#1}}}
\newcommand{\sect}[1]{%
\vskip1cm \centerline{\bfseries\large Lecture~\arabic{chapter}}
\nopagebreak\vskip0.5cm
\centerline{\bfseries\large #1}
\vskip0.5cm
\addcontentsline{toc}{chapter}{#1}

}

\newcommand{\specialchapter}[1]{%
\vskip0.4cm
\centerline{\bfseries\large #1}
\nopagebreak\vskip0.35cm\nopagebreak
\specialtocentry{#1}}

\newcounter{specialtocentrycounter}
\newcommand{\specialtocentry}[1]{\refstepcounter{specialtocentrycounter}
\label{MBspecialtocentry-label\thespecialtocentrycounter}
\edef\mytocfilestring{\expandafter\mytocfilestring\csname myspecialtocentry\endcsname{\thespecialtocentrycounter}{#1}}
\addcontentsline{toc}{section}{#1}}

\renewcommand{\section}{\secdef\ssct\ssect}
\newcommand{\ssct}[2][default]{\refstepcounter{section}
\vskip0.5cm
\leftline{\bfseries\large \arabic{chapter}.\arabic{section}\ \ #1}
\label{MBsection-label\thesection}\nopagebreak
\vskip0.3cm\nopagebreak
\addcontentsline{toc}{section}{#1}
\edef\mytocfilestring{\expandafter\mytocfilestring\csname mysectiontocentry\endcsname{\thesection}{#1}}
}
\newcommand{\ssect}[1]{
\nopagebreak\vskip0.5cm
\centerline{\bfseries\large #1}
\vskip0.5cm
\addcontentsline{toc}{section}{#1}
}

\newenvironment{proofsect}[1]
{\vskip0.1cm\noindent{\rmfamily\itshape#1.}}{\qed\vspace{0.15cm}}

\numberwithin{section}{chapter}
\numberwithin{equation}{chapter}
\numberwithin{theorem}{chapter}

\spnewtheorem{mytheorem}[theorem]{Theorem}{\bf}{\it}
\spnewtheorem{myproblem}[mytheorem]{Problem}{\bf}{\it}
\spnewtheorem{myexercise}[mytheorem]{Exercise}{\bf}{\it}
\spnewtheorem{mydefinition}[mytheorem]{Definition}{\bf}{\it}
\spnewtheorem{mylemma}[mytheorem]{Lemma}{\bf}{\it}
\spnewtheorem{myproposition}[mytheorem]{Proposition}{\bf}{\it}
\spnewtheorem{myconjecture}[mytheorem]{Conjecture}{\bf}{\it}
\spnewtheorem{mycorollary}[mytheorem]{Corollary}{\bf}{\it}
\spnewtheorem{myremark}[mytheorem]{Remark}{\bf}{\rm}

\renewcommand{\AA}{\mathcal A}
\newcommand{\BB}{\mathcal B}
\newcommand{\CC}{\mathcal C}

\newcommand{\EE}{\mathcal E}
\newcommand{\FF}{\mathcal F}
\newcommand{\GG}{\mathcal G}
\newcommand{\HH}{\mathcal H}
\newcommand{\II}{\mathcal I}

\newcommand{\LL}{\mathcal L}
\newcommand{\MM}{\mathcal M}
\newcommand{\NN}{\mathcal N}

\newcommand{\PP}{\mathcal P}

\newcommand{\RR}{\mathcal R}

\newcommand{\VV}{\mathcal V}

\newcommand{\XX}{\mathcal X}

\newcommand{\C}{\mathbb C}
\newcommand{\D}{\mathbb D}
\newcommand{\E}{\mathbb E}

\newcommand{\N}{\mathbb N}

\newcommand{\BbbP}{\mathbb P}
\newcommand{\Q}{\mathbb Q}
\newcommand{\R}{\mathbb R}

\newcommand{\T}{\mathbb T}

\newcommand{\Z}{\mathbb Z}

\newcommand{\scrI}{\mathscr{I}}

\newcommand{\scrX}{\mathscr{X}}

\newcommand{\twoeqref}[2]{(\ref{#1}--\ref{#2})}

\def\myffrac#1#2 in #3{\raise 2.6pt\hbox{$#3 #1$}\mkern-1.5mu\raise 0.8pt\hbox{$#3/$}\mkern-1.1mu\lower 1.5pt\hbox{$#3 #2$}}
\newcommand{\ffrac}[2]{\mathchoice%
	{\myffrac{#1}{#2} in \scriptstyle}
	{\myffrac{#1}{#2} in \scriptstyle}
	{\myffrac{#1}{#2} in \scriptscriptstyle}
	{\myffrac{#1}{#2} in \scriptscriptstyle}
}

\newcommand{\cmssP}{\cmss P}
\newcommand{\fraka}{\fraktura a}
\newcommand{\frakf}{\fraktura f}
\newcommand{\frakb}{\fraktura b}
\newcommand{\frakg}{\fraktura g}
\newcommand{\eff}{{\text{\rm eff}}}

\newcommand\independent{\protect\mathpalette{\protect\independenT}{\perp}}
\def\independenT#1#2{\mathrel{\rlap{$#1#2$}\mkern3mu{#1#2}}}

\newcommand{\ssst}{\scriptscriptstyle}
\newcommand{\PPP}{\text{\rm PPP}}

\newcommand{\dist}{\operatorname{dist}}
\newcommand{\supp}{\operatorname{supp}}
\newcommand{\diam}{\operatorname{diam}}
\newcommand{\textd}{\text{\rm d}\mkern0.5mu}

\newcommand{\texte}{\text{\rm e}}
\newcommand{\Var}{{\text{\rm Var}}}
\newcommand{\Cov}{{\text{\rm Cov}}}
\newcommand{\sgn}{\text{\rm sgn}}

\newcommand{\1}{\text{\sf 1}}
\newcommand{\cc}{{\text{\rm c}}}
\newcommand{\Leb}{\text{\rm Leb}}
\newcommand{\wt}{\widetilde}
\newcommand{\wh}{\widehat}
\newcommand{\lawarrow}{{\overset{\text{\rm law}}\longrightarrow}}
\newcommand{\laweq}{{\overset{\text{\rm law}}=}}

\newcommand{\rhomax}{\rho_{\text{\rm max}}}
\newcommand{\dmax}{\mkern1mu\mathfrak{d}}

\newcommand{\ssup}[1]{{\scriptscriptstyle(#1)}}
\newcommand{\Xiin}{\Xi^{\text{\rm in}}}
\newcommand{\Xiout}{\Xi^{\text{\rm out}}}
\newcommand{\wtXiout}{\wt\Xi^{\text{\rm out}}}
\newcommand{\whXiout}{\wh\Xi^{\text{\rm out}}}

\begin{document}
\mainmatter              
\title{Extrema of the two-dimensional Discrete Gaussian Free Field}
\toctitle{}
\titlerunning{Extrema of 2D DGFF}  
%
\author{Marek Biskup}
\authorrunning{Marek Biskup} 
\institute{Department of Mathematics, UCLA, Los Angeles, CA 90015-1555, U.S.A.\\
\email{biskup@math.ucla.edu.edu}\\ WWW home page:
\texttt{http://www.math.ucla.edu/\homedir biskup/index.html}
}

\maketitle              

\begin{abstract}
These lecture notes offer a gentle introduction to the two-dimensional Discrete Gaussian Free Field with particular attention paid to the scaling limits of the level sets at heights proportional to the absolute maximum. The bulk of the text is based on recent joint papers with O. Louidor and with J. Ding and S. Goswami. Still, new proofs of the tightness and distributional convergence of the centered DGFF maximum are presented that by-pass the use of the modified Branching Random Walk. The text contains a wealth of instructive exercises and a list of open questions and conjectures for future research.

\end{abstract}

\ifmytoc
\newcommand\mytoclineskip{\vskip0.9mm}
\mytableofcontents
\fi

\bigskip
\setcounter{chapter}{0}

\specialchapter{What this course is about (and what it is not)}{}
\addcontentsline{toc}{chapter}{What this course is about (and what it is not)}

\noindent
This is a set of lecture notes for a course on \emph{Discrete Gaussian Free Field} (DGFF) delivered at the 2017 PIMS-CRM Summer School in Probability at the University of British Columbia in June 2017. The course has been quite intense with the total of sixteen 90-minute lectures spanning over 4 weeks. Still, the subject of the DGFF has become so so developed that we  could mainly focus only on one specific aspect: \emph{extremal values}. 
The text below stays close to the actual lectures although later additions have been made to make the exposition self-contained. Each lecture contains a number of exercises that address reasonably accessible aspects of the presentation.

In Lecture~\ref{lec-1} we give an introduction to the DGFF in general spatial dimension and discuss possible limit objects. The aim here is to show that the scaling properties single out the \emph{two-dimensional} DGFF as a case of special interest to which the rest of the text is then exclusively devoted. 

Lecture~\ref{lec-2} opens up by some earlier developments that captured the leading-order behavior of the maximum of the DGFF in growing lattice domains. This sets the scale for the study (in Lectures~\ref{lec-2}--\ref{lec-4}) of what we call \emph{intermediate level sets} --- namely, those where the~DGFF takes values above a constant multiple of the absolute maximum.  The technique of proof is of interest here as it will be reused later: we first establish tightness, then extract a subsequential limit and, finally, identify the limit object uniquely, thus proving the existence of the limit. The limit is expressed using the concept of \emph{Liouville Quantum Gravity} (LQG) which we introduce via Kahane's theory of Gaussian Multiplicative Chaos.

Our next item of interest is the behavior of the DGFF maximum. For this we recall (in Lectures~\ref{lec-5} and~\ref{lec-6}) two basic, albeit technically advanced, techniques from the theory of Gaussian processes: correlation inequalities (Kahane, Borell-TIS, Sudakov-Fernique, FKG) and Fernique's majorization bound for the expected maximum. 

Once the generalities have been dispensed with, we return (in Lectures~\ref{lec-7} and~\ref{lec-8}) to the DGFF and relate the tightness of its centered maximum to that of a \emph{Branching Random Walk}~(BRW). A novel proof of tightness of the DGFF maximum is presented that avoids the so-called modified BRW; instead we rely on the Sudakov-Fernique inequality and the Dekking-Host argument applied to the DGFF coupled with a BRW. This handles the upper tail tightness; for the lower tail we develop the concept of a \emph{concentric decomposition} of the~DGFF that will be useful later as well.

In Lectures~\ref{lec-9}--\ref{lec-11} we move to the \emph{extremal level sets} --- namely, those where the DGFF is within order-unity of the absolute maximum. Similarly to the intermediate levels, we encode these via a three-coordinate point measure that records the scaled spatial position, the centered value and the ``shape'' of the configuration at the local maxima. Using the above proof strategy we then show that, in the scaling limit, this measure tends in law to a Poisson Point Process with a \emph{random} intensity measure whose spatial part can be identified with the \emph{critical} LQG. 

A key technical tool in this part is Liggett's theory of \emph{non-interacting} particle systems which we discuss in full detail (in Lecture~\ref{lec-9}); the uniqueness of the limit measure is linked to the convergence in law of the centered DGFF maximum (Lecture~\ref{lec-10}). Another key tool is the concentric decomposition which permits us to control (in Lecture~\ref{lec-11}) the local structure of the extremal points and  to give (in Lecture~\ref{lec-12}) independent proofs  of spatial tightness of the extremal level sets and convergence in law of the centered DGFF maximum. Interesting corollaries (stated at the end of Lecture~\ref{lec-11}) give existence of \emph{supercritical} LQG measures and the limit law for the Gibbs measure naturally associated with the DGFF. 

The final segment of the course (Lectures~\ref{lec-13}--\ref{lec-15}) is devoted to a \emph{random walk driven by the~DGFF}. After the statement of main theorems we proceed to develop (in Lecture~\ref{lec-13}) the main technique of the proofs: an electric network associated with the DGFF. In Lecture~\ref{lec-14} we give variational characterizations of the effective resistance/conductance using only geometric objects such as paths and cuts. Various duality relations (conductance/resistance reciprocity, path/cut planar duality) along with concentration of measure are invoked to control the effective resistivity across rectangles. These are the key ingredients for controlling (in Lecture~\ref{lec-15}) the relevant aspects of the random walk.

A standalone final lecture (Lecture~\ref{lec-16}) discusses some open research-level problems that stem directly  from the topics discussed in these notes.

\smallskip

To keep the course at a reasonable level of mathematical depth, considerable sacrifices in terms of scope had to be made. The course thus omits a large body of literature on Gaussian Multiplicative Chaos and Liouville Quantum Gravity. We do not address the level sets at heights of order unity and their connections to the Schramm-Loewner evolution and the Conformal Loop Ensemble. We ignore recent developments in Liouville First Passage percolation nor do we discuss the close cousin of our random walk, the Liouville Brownian Motion. We pass over the connections to random-walk local time. 
In many of these cases, better resources exist to obtain the relevant information elsewhere. 

Thanks to our focused approach we are able to present many difficult proofs in nearly complete detail. A patient reader will thus have a chance to learn a number of useful arguments specific to the DGFF as well as several general techniques relevant for probability at large. The hope is that interlacing specifics with generalities will make these notes interesting for newcomers as well as advanced players in this field. 
The exposition should be sufficiently self-contained to serve as the basis for a one-semester graduate course.

\newpage
\bigskip
\centerline{\bf\large Acknowledgments}
\specialtocentry{Acknowledgments}
\medskip\noindent
These notes are largely based on research (supported, in part, by NSF grant DMS-1407558 and GA\v CR project P201/16-15238S) that was done jointly with\newline

\centerline{\bf Oren Louidor} 

\noindent
and, in the last part, with
\newline

\centerline{{\bf Subhajit Goswami} and\ \bf Jian Ding} 

\medskip\smallskip\noindent
I would like to thank these gentlemen for all what I have learned from them, much of which is now reproduced here. 

Next I wish to thank the organizers of the PIMS-CRM summer school, Omer Angel, Mathav Murugan, Ed Perkins and Gordon Slade, for the opportunity to be one of the main lecturers. These notes would not have been written without the motivation (and pressure) arising from that task. Special thanks go to Ed Perkins for finding a cozy place to stay for my whole family over the duration of the school; the bulk of this text was typed over the long summer nights there. 

I have had the opportunity to lecture on various parts of this material at other schools and events; namely, in Bonn, Hejnice, Lyon, Haifa, Atlanta and Salt-Lake City. I am grateful to the organizers of these events, as well as participants thereof, for helping me find a way to parse some of the intricate details into (hopefully) a reasonably digestible form.

Finally, I wish to express my gratitude to Jian Ding, Jean-Christophe Mourrat and Xinyi Li for comments and corrections and to Sara\'{\i} Hern\'andez Torres and, particularly, Yoshi Abe for pointing out numerous errors in early versions of this manuscript. I take the blame for all the issues that went unnoticed.

\bigskip
\leftline{Marek Biskup}
\leftline{May 2019}




\newpage

\chapter{Discrete Gaussian Free Field and scaling limits}
\label{lec-1}

\noindent
In this lecture we define the main object of interest in this course: the Discrete Gaussian Free Field (henceforth abbreviated as~DGFF). By studying its limit properties we are naturally guided towards the two-dimensional case where we describe, in great level of detail, its scaling limit. The limit object, the continuum Gaussian Free Field (CGFF), will underlie, albeit often in disguise, most of the results discussed in the course.

\section{Definitions}
\noindent
For~$d\ge1$ integer, let~$\Z^d$ denote the $d$-dimensional hypercubic lattice. This is an unoriented graph with vertices at the points of~$\R^d$ with integer coordinates and an edge between any pair of vertices at unit Euclidean distance. Denoting the set of all edges (with both orientations identified) by~$\cmss E(\Z^d)$, we  put forward:

\begin{mydefinition}[DGFF, explicit formula]
\label{def-1.1}
Let~$V\subset\Z^d$ be finite. The~DGFF in~$V$ is a process $\{h_x^V\colon x\in \Z^d\}$ indexed by the vertices of~$\Z^d$ with the law given (for any measurable~$A\subseteq\R^{\Z^d}$) by
\begin{equation}
\label{E:1.1}
P(h^V\in A):=\frac1{\text{\rm norm.}}\int_A\texte^{-\frac1{4d}\sum_{(x,y)\in\cmss E(\Z^d)}(h_x-h_y)^2}\prod_{x\in V}\textd h_x\prod_{x\not\in V}\delta_0(\textd h_x)\,.
\end{equation}
Here~$\delta_0$ is the Dirac point-mass at~$0$ and ``{\rm norm.}'' is a normalization constant.
\end{mydefinition}

Notice that the definition forces the values of~$h$ outside~$V$ to zero --- we thus talk about \emph{zero boundary condition}. To see that this definition is good, we pose:

\begin{myexercise}
Prove that the integral in \eqref{E:1.1} is finite for $A:=\R^{\Z^d}$ and so the measure can be normalized to be a probability.
\end{myexercise}

\noindent
The appearance of the~$4d$ factor in the exponent is a convention used in probability; physicists would write $\frac12$ instead of~$\frac1{4d}$. Without this factor, the definition extends readily from~$\Z^d$ to any locally-finite graph but since~$\Z^d$ (in fact~$\Z^2$) will be our main focus, we keep the normalization as above.

Definition~\ref{def-1.1} gives the law of the~DGFF in the form of a \emph{Gibbs measure}, i.e., a measure of the form 
\begin{equation}
\frac1{\text{\rm norm.}}\,\texte^{-\beta H(h)}\,\nu(\textd h)
\end{equation}
where~``norm.'' is again a normalization constant, $H$ is the Hamiltonian, $\beta$ is the inverse temperature and~$\nu$ is an \emph{a priori} (typically a product) measure. Many models of statistical mechanics are cast this way; a key feature of the~DGFF is that the Hamiltonian is a positive-definite quadratic form and~$\nu$ a product Lebesgue measure which makes the law of~$h$ a multivariate Gaussian. This offers the possibility to define the law directly by prescribing its mean and covariance. 

Let~$\{X_n\colon n\ge0\}$ denote the path of a simple symmetric random walk on~$\Z^d$. For~$V\subset\Z^d$, we recall the definition of the \textit{Green function} in~$V$:
\begin{equation}
\label{E:1.3au}
G^V(x,y):=E^x\biggl(\,\sum_{n=0}^{\tau_{V^\cc}-1}1_{\{X_n=y\}}\biggr)\,,
\end{equation}
where~$E^x$ is the expectation with respect to the law of~$X$ with~$X_0=x$ a.s.\ and $\tau_{V^\cc}:=\inf\{n\ge0\colon X_n\not\in V\}$ is the first exit time of the walk from~$V$. We note:

\begin{myexercise}
\label{ex:1.3}
Prove that, for any~$x,y\in\Z^d$,
\begin{equation}
V\mapsto G^V(x,y) \text{ \rm is non-decreasing with respect to set inclusion}.
\end{equation}
In particular, for any~$V\subsetneq\Z^d$ (in any~$d\ge1$), we have $G^V(x,y)<\infty$ for all~$x,y\in\Z^d$.
\end{myexercise}

As is also easy to check, $G^V(x,y)=0$ unless~$x,y\in V$ (in fact, unless~$x$ and~$y$ can be connected by a path on~$\Z^d$ that lies entirely in~$V$). The following additional properties of~$G^V$ will be  important in the sequel:

\begin{myexercise}
Let~$\Delta$ denote the discrete Laplacian on~$\Z^d$ acting on finitely-supported functions $f\colon\Z^d\to\R$ as
\begin{equation}
\Delta f(x):=\sum_{y\colon (x,y)\in\cmss E(\Z^d)}\bigl(f(y)-f(x)\bigr).
\end{equation}
Show that for any~$V\subsetneq\Z^d$ and any~$x\in\Z^d$, $y\mapsto G^V(x,y)$ is the solution to
\begin{equation}
\label{E:1.5}
\begin{cases}
\Delta G^V(\cdot,x)=-2d\delta_x(\cdot),\qquad&\text{on }V,
\\
G^V(\cdot,x)=0,\qquad&\text{on }V^\cc,
\end{cases}
\end{equation}
where~$\delta_x$ is the Kronecker delta at~$x$. 
\end{myexercise}

\noindent
Another way to phrase this exercise is by saying that the Green function is a $(2d)$-multiple of the inverse of the negative Laplacian on~$\ell^2(V)$ --- with Dirichlet boundary conditions on~$V^\cc$. This functional-analytic representation of the Green function allows us to solve:

\begin{myexercise}
Prove that for any~$V\subsetneq\Z^d$ we have:
\begin{enumerate}
\item[(1)] for any $x,y\in\Z^d$,
\begin{equation}
G^V(x,y)=G^V(y,x)\,,
\end{equation}
\item[(2)] for any~$f\colon\Z^d\to\R$ with finite support,
\begin{equation}
\sum_{x,y\in\Z^d}G^V(x,y)f(x)f(y) \ge0.
\end{equation}
\end{enumerate}
\end{myexercise}

We remark that purely probabilistic ways to solve this (i.e., using solely considerations of random walks) exist as well. What matters for us is that properties (1-2) make~$G^V$ a covariance of a Gaussian process. This gives rise to:

\begin{mydefinition}[DGFF, via the Green function]
\label{def-1.2}
Let~$V\subsetneq\Z^d$ be given. The~DGFF in~$V$ is a multivariate Gaussian process $\{h_x^V\colon x\in \Z^d\}$ with law determined by
\begin{equation}
E(h_x^V)=0\quad\text{and}\quad
E(h_x^Vh_y^V)=G^V(x,y),\qquad x,y\in\Z^d,
\end{equation}
or, written concisely,
\begin{equation}
h^V:=\NN(0,G^V)\,.
\end{equation}
\end{mydefinition}

\noindent
Here and henceforth, $\NN(\mu,C)$ denotes the (multivariate) Gaussian with mean $\mu$ and covariance~$C$. The dimensions of these objects will usually be clear from the context.

In order to avoid accusations of logical inconsistency, we pose:

\begin{myexercise}
Prove that for~$V\subset\Z^d$ finite, Definitions~\ref{def-1.1} and~\ref{def-1.2} coincide.
\end{myexercise}

\noindent
The advantage of Definition~\ref{def-1.2} over Definition~\ref{def-1.1} is that it works for infinite~$V$ as well. The functional-analytic connection can be pushed further as follows. Given a finite set $V\subset\Z^d$, consider the Hilbert space $\HH^V:=\{f\colon\Z^d\to\R,\,\supp(f)\subset V\}$ endowed with the \emph{Dirichlet inner product}
\begin{equation}
\langle f,g\rangle_\nabla:=\frac1{2d}\sum_{x\in\Z^d}\nabla f(x)\cdot\nabla g(x)\,,
\end{equation}
where $\nabla f(x)$, the discrete gradient of~$f$ at~$x$, is the vector in~$\R^d$ whose $i$-th component is $f(x+e_i)-f(x)$, for~$e_i$ the $i$-th unit coordinate vector in~$\R^d$. (Since the supports of~$f$ and~$g$ are finite, the sum is effectively finite. The normalization is for consistence with Definition~\ref{def-1.1}.) We then have:

\begin{mylemma}[DGFF as Hilbert-space Gaussian]
\label{lemma-1.8}
For the setting as above with~$V$ finite, let~$\{\varphi_n\colon n=1,\dots,|V|\}$ be an orthonormal basis in~$\HH^V$ and let $Z_1,\dots,Z_{|V|}$ be i.i.d.\ standard normals. Then $\{\wt h_x\colon x\in\Z^d\}$, where
\begin{equation}
\wt h_x:=\sum_{n=1}^{|V|}\varphi_n(x)Z_n,\qquad x\in\Z^d,
\end{equation}
has the law of the~DGFF in~$V$.
\end{mylemma}

\begin{myexercise}
Prove the previous lemma.
\end{myexercise}

Lemma~\ref{lemma-1.8} thus gives yet another way to define~DGFF. (The restriction to finite~$V$ was imposed only for simplicity.) As we will see at the end of this lecture, this is the definition that generalizes seamlessly to  a continuum underlying space. 
Writing~$\{\psi_n\colon n=1,\dots,|V|\}$ for the orthonormal set in~$\ell^2(V)$ of eigenfunctions of the negative Laplacian with the corresponding eigenvalue written as~$2d\lambda_n$, i.e., $-\Delta\psi_n=2d\lambda_n\psi_n$ on~$V$, the choice
\begin{equation}
\varphi_n(x):=\frac1{\sqrt{\lambda_n}}\,\psi_n(x)\,,
\end{equation}
produces an orthonormal basis in~$\HH^V$. This is useful when one wishes to generate samples of the~DGFF efficiently on a computer.

The fact that~$\Z^d$ is an additive group, and the random walk is translation invariant, implies that the Green function and thus the law of $h^V$ are \emph{translation invariant} in the sense that, for all~$z\in\Z^2$ and for~$z+V:=\{z+x\colon x\in V\}$,
\begin{equation}
\label{E:G-trans}
G^{z+V}(x+z,y+z)=G^V(x,y),\quad x,y\in\Z^d,
\end{equation}
and
\begin{equation}
\label{E:h-trans}
\{h^{z+V}_{x+z}\colon x\in\Z^2\}\,\laweq\{h^V_x\colon x\in\Z^2\}
\end{equation}
A similar rule applies to rotations by multiples of~$\frac\pi2$ around any vertex of~$\Z^2$.

We finish this section by a short remark on notation: Throughout these lectures we will write $h^V$ to denote the whole configuration of the~DGFF in~$V$ and write $h^V_x$ for the value of~$h^V$ at~$x$. We may at times write $h^V(x)$ instead of~$h^V_x$ when the resulting expression is easier to parse. 

\section{Why $d=2$ only?}
\noindent
As noted above, this course will focus on the DGFF in~$d=2$. Let us therefore explain what makes the two-dimensional~DGFF special. We begin by noting:

\begin{mylemma}[Green function growth rate]
\label{lemma-1.10}
Let~$V_N:=(0,N)^d\cap\Z^d$ and, for any $\epsilon\in(0,1/2)$, denote $V_N^\epsilon:=(\epsilon N,(1-\epsilon)N)^d\cap\Z^d$. Then for any~$x\in V_N^\epsilon$,
\begin{equation}
G^{V_N}(x,x)\,\underset{N\to\infty}\sim\,\,\begin{cases}
N,\qquad&\text{}d=1,
\\
\log N,\qquad&\text{}d=2,
\\
1,\qquad&\text{}d\ge3,
\end{cases}
\end{equation}
where~``$\sim$'' means that the ratio of the left and right-hand side tends to a positive and finite number as~$N\to\infty$ (which may depend on where~$x$ is asymptotically located in~$V_N$).
\end{mylemma}

\begin{proofsect}{Proof (sketch)}
We will only prove this in~$d=1$ and~$d\ge3$ as the~$d=2$ case will be treated later in far more explicit terms. First note that a routine application of the Strong Markov Property for the simple symmetric random walk yields
\begin{equation}
\label{E:1.15wk}
G^V(x,x)=\frac1{P^x(\hat\tau_x>\tau_{V^\cc})}\,,
\end{equation}
where~$\hat\tau_x:=\inf\{n\ge1\colon X_n=x\}$ is the first \emph{return} time to~$x$. Assuming~$x\in V$, in~$d=1$ we then have
\begin{equation}
P^x(\hat\tau_x>\tau_{V^\cc}) = \frac12\Bigl[P^{x+1}(\tau_x>\tau_{V^\cc})+P^{x-1}(\tau_x>\tau_{V^\cc})\Bigr]\,,
\end{equation}
where~$\tau_x:=\inf\{n\ge0\colon X_n=x\}$ is the first \emph{hitting} time of~$x$.
For $V:=V_N$, the Markov property of the simple random walk shows that $y\mapsto P^y(\tau_x>\tau_{V_N^\cc})$ is discrete harmonic (and thus piecewise linear) on~$\{1,\dots,x-1\}\cup\{x+1,\dots,N-1\}$ with value zero at~$y=x$ and one at~$y=0$ and~$y=N$. Hence
\begin{equation}
P^{x+1}(\tau_x>\tau_{V_N^\cc}) = \frac1{N-x}\quad\text{and}\quad P^{x-1}(\tau_x>\tau_{V_N^\cc})=\frac1x.
\end{equation}
As~$x\in V_N^\epsilon$, both of these probabilities are order~$1/N$ and scale nicely when~$x$ grows proportionally to~$N$. This proves the claim in~$d=1$.

In~$d\ge3$ we note that the transience and translation invariance of the simple symmetric random walk imply
\begin{equation}
G^{V_N}(x,x)\,\underset{N\to\infty}\longrightarrow\,\frac1{P^0(\hat\tau_0=\infty)}
\end{equation}
uniformly in~$x\in V_N^\epsilon$. (Transience is equivalent to $P^0(\hat\tau_0=\infty)>0$.)
\end{proofsect}

Let us now proceed to examine the law of the whole DGFF configuration in the limit as~$N\to\infty$.
Focusing for the moment on~$d=1$, the fact that the variance blows up suggests that we normalize the~DGFF by the square-root of the variance, i.e.,~$\sqrt N$, and attempt to extract a limit. This does work and yields:

\begin{mytheorem}[Scaling to Brownian Bridge in $d=1$]
\label{thm-1.11}
Suppose~$d=1$ and let~$h^{V_N}$ be the~DGFF in~$V_N:=(0,N)\cap\Z$. Then
\begin{equation}
\Bigl\{\frac1{\sqrt N}h_{\lfloor tN\rfloor}^{V_N}\colon t\in[0,1]\Bigr\}
\,\underset{N\to\infty}\lawarrow\,\bigl\{\sqrt2\, W_t\colon t\in[0,1]\bigr\}\,,
\end{equation}
where~$W$ is the standard Brownian Bridge on~$[0,1]$.
\end{mytheorem}

We leave it to the reader to solve:

\begin{myexercise}
Prove Theorem~\ref{thm-1.11} with the convergence in the sense of finite-dimensional distributions or, if you like the challenge, in Skorokhod topology. Hint: Note that $\{h_{x+1}^{V_N}-h_x^{V_N}\colon x=0,\dots,N-1\}$ are i.i.d.\ $\NN(0,2)$ conditioned on their total sum being zero.
\end{myexercise}

We remark that the limit taken in Theorem~\ref{thm-1.11} is an example of a \emph{scaling (or continuum) limit} --- the lattice spacing is taken to zero while keeping the overall (continuum) domain fixed. In renormalization group theory, taking the scaling limit corresponds to the removal of an ultraviolet cutoff.

Moving to~$d\ge3$, in the proof of Lemma~\ref{lemma-1.10} we observed enough to conclude:

\begin{mytheorem}[Full-space limit in $d\ge3$]
Suppose~$d\ge3$ and abbreviate $\wt V_N:=(-N/2,N/2)^d\cap\Z^d$. Then for any~$x,y\in\Z^d$,
\begin{equation}
G^{\wt V_N}(x,y)\underset{N\to\infty}\longrightarrow\, G^{\Z^d}(x,y):=\int\frac{\textd k}{(2\pi)^d}\frac{\cos(k\cdot(x-y))}{\frac2d\sum_{j=1}^d\sin(k_j/2)^2}\,.
\end{equation}
In particular, $h^{\wt V_N}\to \NN(0, G^{\Z^d})$ = full space~DGFF.
\end{mytheorem}

This means that the~DGFF in large enough (finite) domains is well approximated by the full-space~DGFF as long as we are far away from the boundary. This is an example of a \emph{thermodynamic limit} --- the lattice stays fixed and the domain boundaries slide off to infinity. In renormalization group theory, taking the thermodynamic limit corresponds to the removal of an infrared cutoff.

From Lemma~\ref{lemma-1.10} it is clear that the thermodynamic limit is meaningless for the two-dimensional~DGFF (indeed, variances blow up and, since we are talking about Gaussian random variables, there is no tightness). Let us attempt to take the scaling limit just as in~$d=1$: normalize the field by the square-root of the variance (i.e.,~$\sqrt{\log N}$) and extract a distributional limit (for which it suffices to prove the limit of the covariances). Here we note that, for all~$\epsilon>0$,
\begin{equation}
\sup_{N\ge1}\,\sup_{\begin{subarray}{c}
x,y\in V_N\\|x-y|\ge\epsilon N
\end{subarray}}
G^{V_N}(x,y)<\infty,
\end{equation}
a fact that we will prove later. For any $s,t\in(0,1)^2$ we thus get
\begin{equation}
\label{E:1.22uai}
\Cov\biggl(\frac{h^{V_N}_{\lfloor tN\rfloor}}{\sqrt{\log N}},\frac{h^{V_N}_{\lfloor sN\rfloor}}{\sqrt{\log N}}\biggr)
\,\,\underset{N\to\infty}\longrightarrow\,\,\begin{cases}
c(t)>0,\qquad&\text{if }s=t,
\\
0,\qquad&\text{else},
\end{cases}
\end{equation}
where, here and henceforth,
\begin{equation}
\lfloor tN\rfloor := \text{the unique~$z\in\Z^2$ such that }\,tN\in z+[0,1)^2\,.
\end{equation}
The only way to read \eqref{E:1.22uai} is that the limit distribution is a collection of independent normals indexed by~$t\in(0,1)^2$ --- an object too irregular and generic to retain useful information from before the limit was taken. 

As we will see, the right thing to do is to take the limit of the~DGFF \emph{without} any normalization. That will lead to a singular limit as well but one that captures better the behavior of the~DGFF. Moreover, the limit object exhibits additional symmetries (e.g., conformal invariance) not present at the discrete level.

\section{Green function asymptotic}
\noindent
Let us now analyze the situation in $d=2$ in more detail. Our aim is to consider the~DGFF in \emph{sequences} of lattice domains~$\{D_N\colon N\ge1\}$ that somehow correspond to the scaled-up versions of a given continuum domain~$D\subset\C$. The assumptions on the continuum domains are the content of:

\begin{mydefinition}[Admissible domains]
\label{def-admis-domain}
Let~$\mathfrak D$ be the class of sets~$D\subset\C$ that are bounded, open and such that their topological boundary~$\partial D$ is the union of a finite number of connected components each of which has a positive (Euclidean) diameter.
\end{mydefinition}

For the sake of future use we note:

\begin{myexercise}
\label{ex:1.15}
All bounded, open and simply connected~$D\subset\C$ belong to~$\mathfrak D$.
\end{myexercise} 

\begin{figure}[t]
\vglue2mm
\centerline{\includegraphics[width=0.6\textwidth]{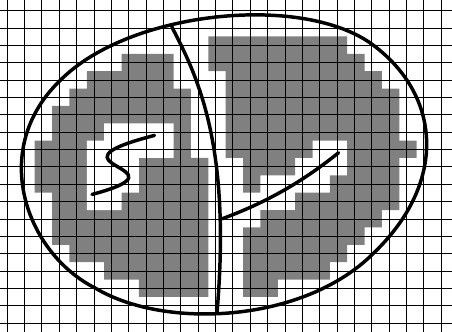}
}
\vglue2mm
\begin{quote}
\small 
\caption{\small
\label{fig-domains}
An example of an admissible discretization of a continuum domain. The continuum domain is the region in the plane bounded by the solid curves. The discrete domain is the set of lattice vertices in the shaded areas.}
\normalsize
\end{quote}
\end{figure}

As to what  sequences of discrete approximations of~$D$ we will permit, a natural choice would be to work with plain discretizations $\{x\in\Z^2\colon x/N\in D\}$. However, this is too crude because parts of~$\partial D$ could be missed out completely; see Fig.~\ref{fig-domains}. We thus have to qualify admissible discretizations more precisely: 

\begin{mydefinition}[Admissible lattice approximations]
A family of subsets $\{D_N\colon N\ge1\}$ of~$\Z^2$ is a \emph{sequence of admissible lattice approximations} of a domain~$D\in\mathfrak D$ if
\begin{equation}
\label{E:1.22}
D_N \subseteq \bigl\{x \in \Z^2 \colon \dist_\infty(x/N,D^\cc)>\ffrac1N\bigr\}\,,
\end{equation}
where $\dist_\infty$ denotes the $\ell^\infty$-distance on~$\Z^2$, and if, for each $\delta>0$,
\begin{equation}
\label{E:1.23}
D_N\supseteq\bigl\{x \in \Z^2 \colon \dist_\infty(x/N,D^\cc)>\delta\bigr\}
\end{equation}
holds once~$N$ is sufficiently large (depending possibly on~$\delta$). 
\end{mydefinition}

Note that this still ensures that $x\in D_N$ implies $x/N\in D$. For the statement of the next theorem, consider the (standard) Brownian motion $\{B_t\colon t\ge0\}$ on~$\R^2$ and let $\tau_{D^\cc}:=\inf\{t\ge0\colon B_t\not\in D\}$ be the first exit time from~$D$. Denote by $\Pi^D(x,\cdot)$ the law of the exit point from~$D$ of the Brownian motion started from~$x$; i.e.,
\begin{equation}
\label{E:1.26ua}
\Pi^D(x,A):=P^x\bigl(B_{\tau_{D^\cc}}\in A\bigr),
\end{equation}
for all Borel~$A\subseteq\R^2$.
This measure is supported on~$\partial D$ and is also known as the \emph{harmonic measure from~$x$ in~$D$}. The following then holds:

\begin{mytheorem}[Green function asymptotic in~$d=2$]
\label{thm-1.17}
Suppose~$d=2$ and let $g:=2/\pi$. There is $c_0\in\R$ such that for all domains~$D\in\mathfrak D$, all sequences~$\{D_N\colon N\ge1\}$ of admissible lattice approximations of~$D$, and all~$x\in D$,
\begin{equation}
\label{E:1.25}
G^{D_N}\bigl(\lfloor xN\rfloor,\lfloor xN\rfloor\bigr)= g\log N+c_0+g\int_{\partial D}\Pi^D(x,\textd z)\log|x-z|+o(1),
\end{equation}
where~$o(1)\to0$ as~$N\to\infty$ locally uniformly in~$x\in D$.
Moreover, for all $x,y\in D$ with~$x\ne y$, we also have
\begin{equation}
\label{E:1.26}
G^{D_N}\bigl(\lfloor xN\rfloor,\lfloor yN\rfloor\bigr)= -g\log|x-y|+g\int_{\partial D}\Pi^D(x,\textd z)\log|y-z|+o(1),
\end{equation}
where~$o(1)\to0$ as~$N\to\infty$ locally uniformly in $(x,y)\in D\times D$ with~$x\ne y$ and $|\cdot|$ denotes the Euclidean norm on~$\R^2$.
\end{mytheorem}

\begin{proofsect}{Proof (modulo two lemmas)}
The proof of the theorem starts by a convenient representation of the Green function using the \emph{potential kernel} $\fraka\colon\Z^2\to[0,\infty)$ defined, e.g., by the explicit formula
\begin{equation}
\label{E:1.27}
\fraka(x):=\int_{(-\pi,\pi)^2}\frac{\textd k}{(2\pi)^2}\,\frac{1-\cos(k\cdot x)}{\sin(k_1/2)^2+\sin(k_2/2)^2}\,.
\end{equation}
A characteristic (although not completely characterizing) property of the potential kernel is the content of:

\begin{myexercise}
Show that~$\fraka$ solves the Poisson problem
\begin{equation}
\label{E:1.28}
\begin{cases}
\Delta\fraka(\cdot)=4\delta_0(\cdot),\qquad&\text{on }\Z^2,
\\
\fraka(0)=0,\qquad&\text{}
\end{cases}
\end{equation}
where, as before, $\delta_0$ is the Kronecker delta at~$0$.
\end{myexercise}

Writing~$\partial V$ for external vertex boundary of~$V$, i.e., the set of vertices in~$V^\cc$ that have an edge to a vertex in~$V$, we now get:

\begin{mylemma}[Green function from potential kernel]
\label{lemma-1.19}
For any finite set $V\subset\Z^2$ and any vertices $x,y\in V$,
\begin{equation}
\label{E:1.29}
G^V(x,y)=-\fraka(x-y)+\sum_{z\in\partial V}H^V(x,z)\fraka(z-y)\,,
\end{equation}
where~$H^V(x,z):=P^x(X_{\tau_{V^c}}=z)$ is the probability that the simple symmetric random walk~$X$ started at~$x$ exists~$V$ at~$z\in\partial V$.
\end{mylemma}

\begin{proofsect}{Proof}
Fix~$x\in V$ and let~$X$ denote a path of the  simple symmetric random walk. In light of \eqref{E:1.5}, \eqref{E:1.28} and the translation invariance of the lattice Laplacian,
\begin{equation}
\phi(y):=G^V(x,y)+\fraka(x-y)
\end{equation}
is discrete harmonic in~$V$. Hence~$M_n:=\phi(X_{\tau_{V^\cc}\wedge n})$ is a martingale for the usual filtration $\sigma(X_0,\dots,X_n)$. The finiteness of~$V$ ensures that~$M$ is bounded and~$\tau_{V^\cc}<\infty$ a.s. The Optional Stopping Theorem then gives
\begin{equation}
\phi(y)=E^y\bigl(\phi(X_{\tau_{V^\cc}})\bigr)=\sum_{z\in\partial V}H^V(y,z)\phi(z).
\end{equation}
Since~$G^V(x,\cdot)$ vanishes on~$\partial V$, this along with the symmetry of~$x,y\mapsto G^V(x,y)$ and $x,y\mapsto\fraka(x-y)$ readily imply the claim.
\end{proofsect}

We note that the restriction to finite~$V$ was not a mere convenience. Indeed:

\begin{myexercise}
Show that $G^{\Z^2\smallsetminus\{0\}}(x,y)=\fraka(x)+\fraka(y)-\fraka(x-y)$. Use this to conclude that, in particular, \eqref{E:1.29} is generally false for infinite~$V$.
\end{myexercise}

As the next step of the proof, we invoke:

\begin{mylemma}[Potential kernel asymptotic]
\label{lemma-1.20}
There is~$c_0\in\R$ such that
\begin{equation}
\label{E:1.32}
\fraka(x)=g\log|x|+c_0+O\bigl(|x|^{-2}\bigr),
\end{equation}
where, we recall, $g=2/\pi$ and~$|x|$ is the Euclidean norm of~$x$. 
\end{mylemma}

\noindent
This asymptotic form was apparently first proved by St\"ohr~\cite{Stoehr} in 1950. In a 2004 paper, Kozma and Schreiber~\cite{Kozma-Schreiber} analyzed the behavior of the potential kernel on other lattices and identified the constants~$g$ and~$c_0$ in terms of specific geometric attributes of the underlying lattice. In our case, we can certainly attempt to compute the asymptotic explicitly:

\begin{myexercise}
Prove Lemma~\ref{lemma-1.20} by asymptotic analysis of \eqref{E:1.27}.
\end{myexercise}

\noindent
Using \eqref{E:1.29} and \eqref{E:1.32} together, for~$x,y\in D$ with~$x\ne y$ we now get
\begin{multline}
\label{E:1.33}
\quad
G^{D_N}\bigl(\lfloor xN\rfloor,\lfloor yN\rfloor\bigr)
= g\log|x-y|
\\+g\sum_{z\in\partial D_N}H^{D_N}\bigl(\lfloor xN\rfloor,z\bigr)\log|y-z/N|+O(1/N)\,,
\quad
\end{multline}
where the~$O(1/N)$ term arises from the approximation of (the various occurrences of) $\lfloor xN\rfloor$ by~$xN$ and also from the error in \eqref{E:1.32}. To get \eqref{E:1.26}, we thus need to convert the sum into the integral. Here we will need:

\begin{mylemma}[Weak convergence of discrete harmonic measures]
\label{lemma-1.22}
For any domain~$D\in\mathfrak D$ and any sequence $\{D_N\colon N\ge1\}$ of admissible lattice approximations of~$D$,
\begin{equation}
\label{E:1.34}
\sum_{z\in\partial D_N}H^{D_N}\bigl(\lfloor xN\rfloor,z\bigr)\,\delta_{z/N}(\cdot)\,\underset{N\to\infty}{\overset{\text{\rm weakly}}\longrightarrow}\,\Pi^D(x,\cdot).
\end{equation}
\end{mylemma}

We will omit the proof as it would take us too far on a tangent; the reader is instead referred to (the Appendix of)  Biskup and Louidor~\cite{BL2}. The idea is to use Donsker's Invariance Principle to extract a coupling of the simple random walk and the Brownian motion that ensures that once the random walk exits~$D_N$, the Brownian motion will exit~$D$ at a ``nearby'' point, and \textit{vice versa}. This is where we find it useful that the boundary components have a positive diameter.

Since~$u\mapsto\log|y-u|$ is bounded and continuous in a neighborhood of~$\partial D$ (whenever~$y\in D$), the weak convergence in Lemma~\ref{lemma-1.22} implies
\begin{equation}
\sum_{z\in\partial D_N}H^{D_N}\bigl(\lfloor xN\rfloor,z\bigr)\log|y-z/N|
=
\int_{\partial D}\Pi^D(x,\textd z)\log|y-z| +o(1)\,,
\end{equation}
with~$o(1)\to0$ as~$N\to\infty$. This proves \eqref{E:1.26}. The proof of \eqref{E:1.25} only amounts to the substitution of $-g\log|x-y|$ in \eqref{E:1.33} by $g\log N+c_0$; the convergence of the sum to the corresponding integral is then handled as before.
\end{proofsect}

We remark that the convergence of the discrete harmonic measure to the continuous one in Lemma~\ref{lemma-1.22} is where the above assumptions on the continuum domain and its lattice approximations crucially enter. (In fact, we could perhaps even \emph{define} admissible lattice approximations by requiring \eqref{E:1.34} to hold.) The reader familiar with conformal mapping theory may wish to think of \eqref{E:1.34} as a version of Carath\'eodory convergence for discrete planar domains.

The objects appearing on the right-hand side of \twoeqref{E:1.25}{E:1.26} are actually well known. Indeed, we set:

\begin{mydefinition}[Continuum Green function and conformal radius]
\label{def-1.24}
For bounded, open $D\subset\C$, we define the continuum Green function in~$D$ from~$x$ to~$y$ by
\begin{equation}
\wh G^D(x,y):=-g\log|x-y|+g\int_{\partial D}\Pi^D(x,\textd z)\log|y-z|\,.
\end{equation}
Similarly, for~$x\in D$ we define
\begin{equation}
\label{E:rad}
r_D(x):=\exp\Bigl\{\int_{\partial D}\Pi^D(x,\textd z)\log|x-z|\Bigr\}
\end{equation}
to be the conformal radius of~$D$ from~$x$.
\end{mydefinition}

The continuum Green function is usually defined as the fundamental solution to the Poisson equation; i.e., a continuum version of \eqref{E:1.5}. We will not need this characterization in what follows so we will content ourselves with the explicit form above. Similarly, for open simply connected~$D\subset\C$, the conformal radius of~$D$ from~$x$ is  defined as the value $|f'(x)|^{-1}$ for~$f$ any conformal bijection of~$D$ onto $\{z\in\C\colon|z|<1\}$ such that $f(x)=0$. (The result does not depend on the choice of the bijection.) The reader will find it instructive to solve:

\begin{myexercise}
Check that this coincides with~$r_D(x)$ above.
\end{myexercise}

\noindent
For this as well as later derivations it may be useful to know:

\begin{mylemma}[Conformal invariance of harmonic measure]
\label{lemma-1.26}
For any conformal bijection~$f\colon D\mapsto f(D)$ between $D,f(D)\in\mathfrak D$, 
\begin{equation}
\Pi^D(x,A)=\Pi^{f(D)}\bigl(\,f(x),f(A)\bigr)
\end{equation}
for any measurable~$A\subseteq\partial D$. 
\end{mylemma}

\noindent
This can be proved using conformal invariance of the Brownian motion, although more direct approaches to prove this exist as well. The proof is fairly straightforward when both~$f$ and~$f^{-1}$ extend continuously to the boundaries; the general case is handled by stopping the Brownian motion before it hits the boundary and invoking approximation arguments.

\section{Continuum Gaussian Free Field}
\noindent
Theorem~\ref{thm-1.17} reveals two important facts: First, a pointwise limit of the unscaled DGFF is meaningless as, in light of \eqref{E:1.25}, there is no tightness. Notwithstanding, by \eqref{E:1.26}, the off-diagonal covariances of the~DGFF do have a limit which is given by the continuum Green function. This function is singular on the diagonal, but the singularity is only logarithmic and thus relatively mild. This permits us to derive:

\begin{mytheorem}[Scaling limit of DGFF]
\label{thm-1.26}
Let~$D\in\mathfrak D$ and consider a sequence~$\{D_N\colon N\ge1\}$ of admissible lattice approximations of~$D$. For any bounded measurable function $f\colon D\to\R$, let
\begin{equation}
\label{E:1.39}
h^{D_N}(f):=\int_D \textd x \,f(x) h^{D_N}_{\lfloor xN\rfloor}.
\end{equation}
 Then
\begin{equation}
h^{D_N}(f)\,\underset{N\to\infty}\lawarrow\,\NN(0,\sigma_f^2)\,,
\end{equation}
where
\begin{equation}
\label{E:1.41}
\sigma_f^2:=\int_{D\times D}\textd x\textd y\,f(x)f(y)\wh G^D(x,y)\,.
\end{equation}
\end{mytheorem}

\begin{proofsect}{Proof}
The random variable $h^{D_N}(f)$ is Gaussian with mean zero and variance
\begin{equation}
\label{E:1.42}
E\bigl[h^{D_N}(f)^2\bigr]=\int_{D\times D}\textd x\textd y\,f(x)f(y) G^{D_N}\bigl(\lfloor xN\rfloor,\lfloor yN\rfloor\bigr). 
\end{equation}
The monotonicity from Exercise~\ref{ex:1.3} and the reasoning underlying the proof of Theorem~\ref{thm-1.17} show the existence of~$c\in\R$ such that for all~$N$ large enough and all~$x,y\in D_N$,
\begin{equation}
\label{E:1.44a}
G^{D_N}(x,y)\le g\log\biggl(\,\frac N{|x-y|\vee1}\biggr)+ c.
\end{equation}
Since~$D$ is bounded, this gives
\begin{equation}
\label{E:1.43}
G^{D_N}\bigl(\lfloor xN\rfloor,\lfloor yN\rfloor\bigr)\le -g\log|x-y| + \tilde c
\end{equation}
for some~$\tilde c\in\R$, uniformly in~$x,y\in D$. Using this bound, we can estimate (and later neglect) the contributions to the integral in \eqref{E:1.42} from pairs~$(x,y)$ with $|x-y|<\epsilon$ for any given~$\epsilon>0$. The convergence of the remaining part of the integral is then treated using the pointwise convergence \eqref{E:1.26} (which is locally uniform on the set where ~$x\ne y$) and the Bounded Convergence Theorem.
\end{proofsect}

\begin{myexercise}
\label{ex:1.28}
Give a detailed proof of \eqref{E:1.44a}.
\end{myexercise}

Formula \eqref{E:1.39} can be viewed as a projection of the field configuration onto a test function. Theorem~\ref{thm-1.26} then implies that these projections admit a joint distributional limit. This suggests we could regard the limit object as a \emph{linear functional} on a suitable space of test functions, which leads to:

\begin{mydefinition}[CGFF as a function space-indexed Gaussian]
\label{def-1.28}
A continuum Gaussian Free Field (CGFF) on a bounded, open~$D\subset\C$ is an assignment~$f\mapsto\Phi(f)$ of a random variable to each bounded measurable $f\colon D\to\R$ such that
\begin{enumerate}
\item[(1)] $\Phi$ is a.s.\ linear, i.e., 
\begin{equation}
\Phi(a f+bg)=a\Phi(f)+b\Phi(g)\quad\text{\rm a.s.}
\end{equation}
 for all bounded measurable~$f$ and~$g$ and each~$a,b\in\R$, and
\item[(2)] for all bounded measurable~$f$,
\begin{equation}
\Phi(f)\,\,\laweq\,\,\NN(0,\sigma_f^2)\,,
\end{equation}
where $\sigma_f^2$ is as in \eqref{E:1.41}.
\end{enumerate}
\end{mydefinition}

Theorem~\ref{thm-1.26} shows that the~DGFF scales to the CGFF in this sense --- which, modulo a density argument, also entails that a CGFF as defined above exists! (Linearity is immediate from \eqref{E:1.39}. An independent construction of a CGFF will be performed in Exercise~\ref{ex:2.15}.) 

We remark that definitions given in the literature usually require that the CGFF is even a.s.\ \emph{continuous} in a suitable topology over a suitable (and suitably completed) class of functions. We will not need such continuity considerations here so we do not pursue them in any detail. However, they are important when one tries to assign meaning to~$\Phi(f)$ for singular~$f$ (e.g., those vanishing Lebesgue a.e.) or for~$f$ varying continuously with respect to some parameter. 
This is for example useful in the study of the \emph{disc-average process} $t\mapsto\Phi(f_t)$, where
\begin{equation}
f_t(y):=\frac1{\pi\texte^{-2t}}1_{B(x,\texte^{-t})}(y)
\quad \text{for}\quad
B(x,r):=\bigl\{y\in\C\colon |y-x|<r\bigr\}\,.
\end{equation}
Here it is quite instructive to note:

\begin{myexercise}[Disc average process]
\label{ex:1.29}
For CGFF as defined above and for any~$x\in D$ at Euclidean distance~$r>0$ from~$D^\cc$, show that for~$t>\log(1/r)$, the process $t\mapsto\Phi(f_t)$ has independent increments with
\begin{equation}
\Var\bigl(\Phi(f_t)\bigr)=g t+c_1+g\log r_D(x)
\end{equation}
for some~$c_1\in\R$ independent of~$x$. In particular, $t\mapsto\Phi(f_t)$ admits a continuous version whose law is (for $t>\log(1/r)$) that of a Brownian motion.
\end{myexercise}

The same conclusion is obtained for the \emph{circle average process}, where~$\Phi$ is projected, via a suitable regularization procedure, onto the indicator of a circle $\{y\in\C\colon |x-y|=r\}$. This is because the circle and disc average of the continuum Green function (in one variable with the other one fixed) coincide for small-enough radii, due to the fact that the Green function is harmonic away from the diagonal.

We refer to the textbook by Janson~\cite{Janson} for a thorough discussion of various aspects of generalized Gaussian processes. The Gaussian Free Field existed in physics for a long time where it played the role of a ``trivial,'' which in the jargon of physics means ``non-interacting,'' field theory. Through various scaling limits as well as in its own right, it has recently come to the focus of a larger mathematical community as well. A pedestrian introduction to the subject of the CGFF can be found in Sheffield~\cite{Sheffield-review} and also in Chapter~5 of the recent posting by Armstrong, Kuusi and Mourrat~\cite{AKM}.

\normalsize


\chapter{Maximum and intermediate values}
\label{lec-2}
\noindent
In this lecture we begin to discuss the main topic of interest in this course: extremal properties of the DGFF sample paths. After some introduction and pictures, we focus attention on the behavior of the absolute maximum and the level sets at heights proportional to the absolute maximum. We then state the main theorem on the scaling limit of such level sets and link the limit object to the concept of Liouville Quantum Gravity. The proof of the main theorem is relegated to the forthcoming lectures.

\section{Level set geometry}
\begin{figure}[t]
\vglue-3mm
\centerline{\includegraphics[height=0.6\textwidth]{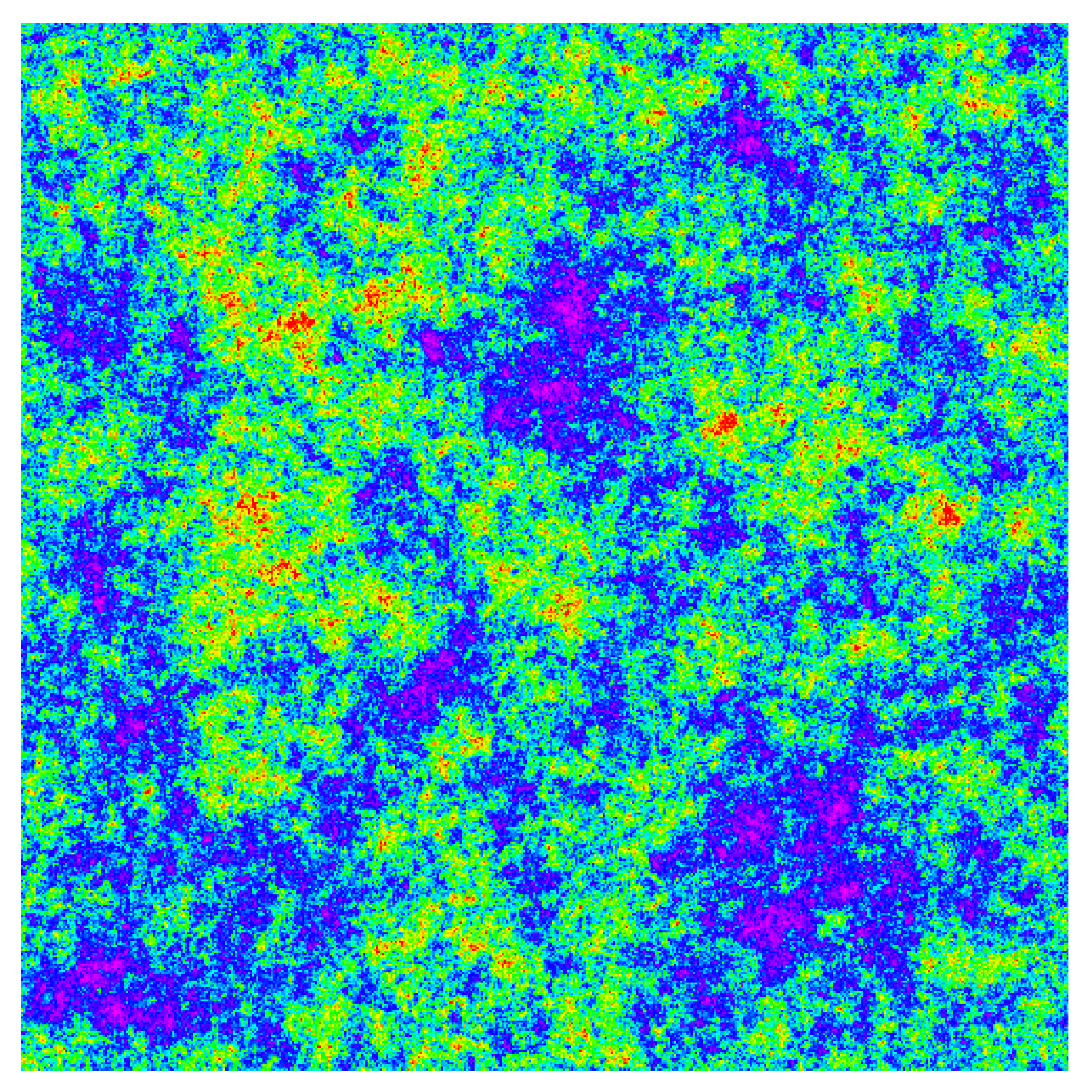}
}
\vglue0mm
\begin{quote}
\small 
\caption{\small
\label{fig-DGFF}
A sample of the DGFF on 300$\times$300 square in~$\Z^2$. The cold colors (purple and blue) indicate low values, the warm colors (yellow and red) indicate large values. The fractal nature of the sample is quite apparent.}
\normalsize
\end{quote}
\end{figure}

\noindent
The existence of the scaling limit established in Theorem~\ref{thm-1.26} indicates that the law of the DGFF is asymptotically scale invariant. Scale invariance of a random object usually entails one of the following two possibilities: 
\begin{itemize}
\item
either the object is trivial and boring (e.g., degenerate, flat, non-random), 
\item
or it is very interesting (e.g., chaotic, rough, fractal).
\end{itemize} 
As attested by Fig.~\ref{fig-DGFF}, the two-dimensional DGFF seems, quite definitely, to fall into the latter category.

Looking at Fig.~\ref{fig-DGFF} more closely, a natural first question is to understand the behavior of the (implicit) boundaries between warm and cold colors. As the field averages to zero, and should thus take both positive and negative values pretty much equally likely, this amounts to looking at the contour lines \emph{between} the regions where the field is positive and where it is negative. 
This has been done and constitutes the beautiful work started by Schramm and Sheffield~\cite{Schramm-Sheffield} and continued in Sheffield and Werner~\cite{Sheffield-Werner}, Sheffield~\cite{Sheffield-CLE} and Sheffield and Miller~\cite{MS1,MS2,MS3,MS4}. We thus know that the contour lines admit a scaling limit to a process of nested collections of loops called the \emph{Conformal Loop Ensemble} with the individual curves closely related to the \emph{Schramm-Loewner process}~SLE$_4$.

\begin{figure}[t]
\vglue1mm
\centerline{\includegraphics[width=0.32\textwidth]{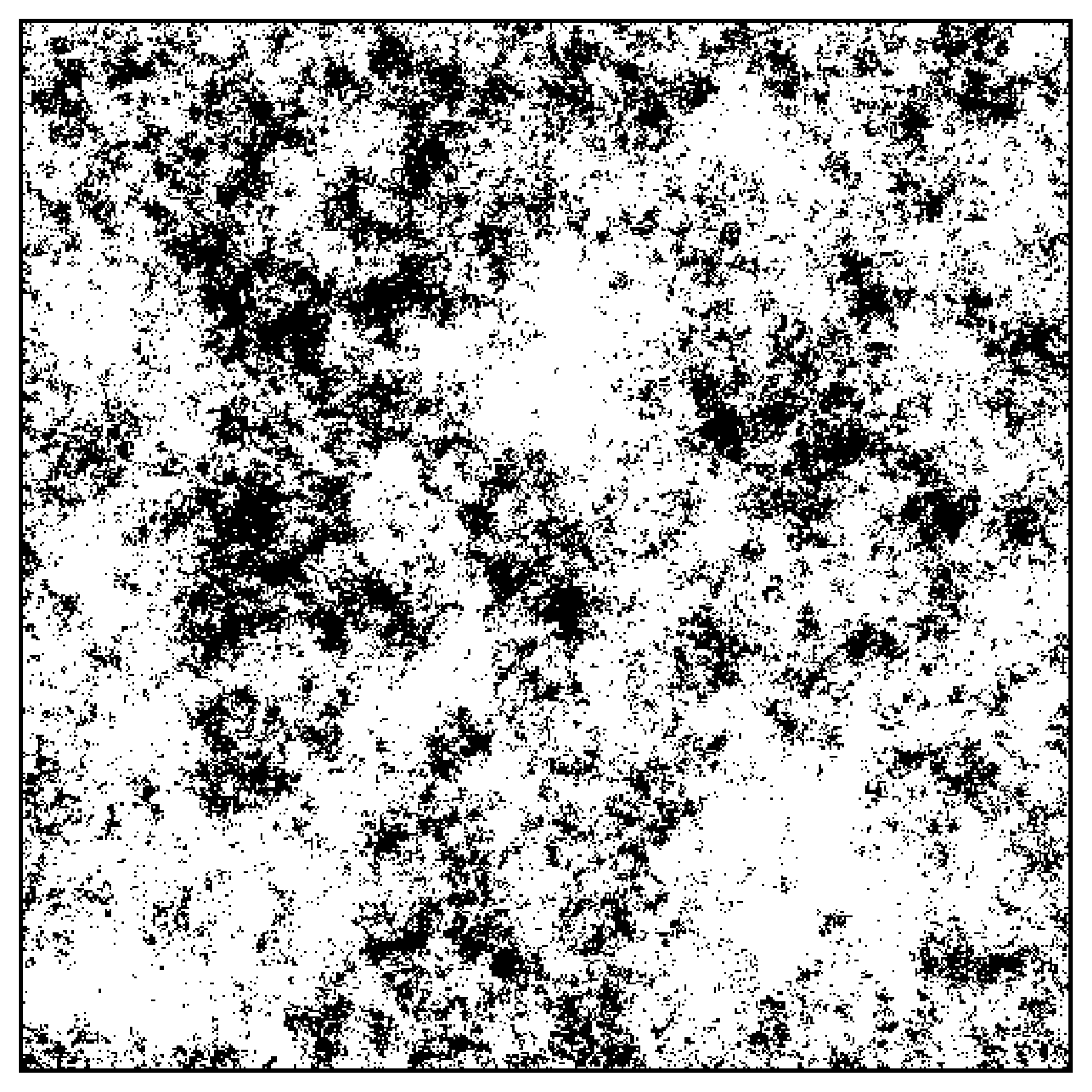}
\includegraphics[width=0.32\textwidth]{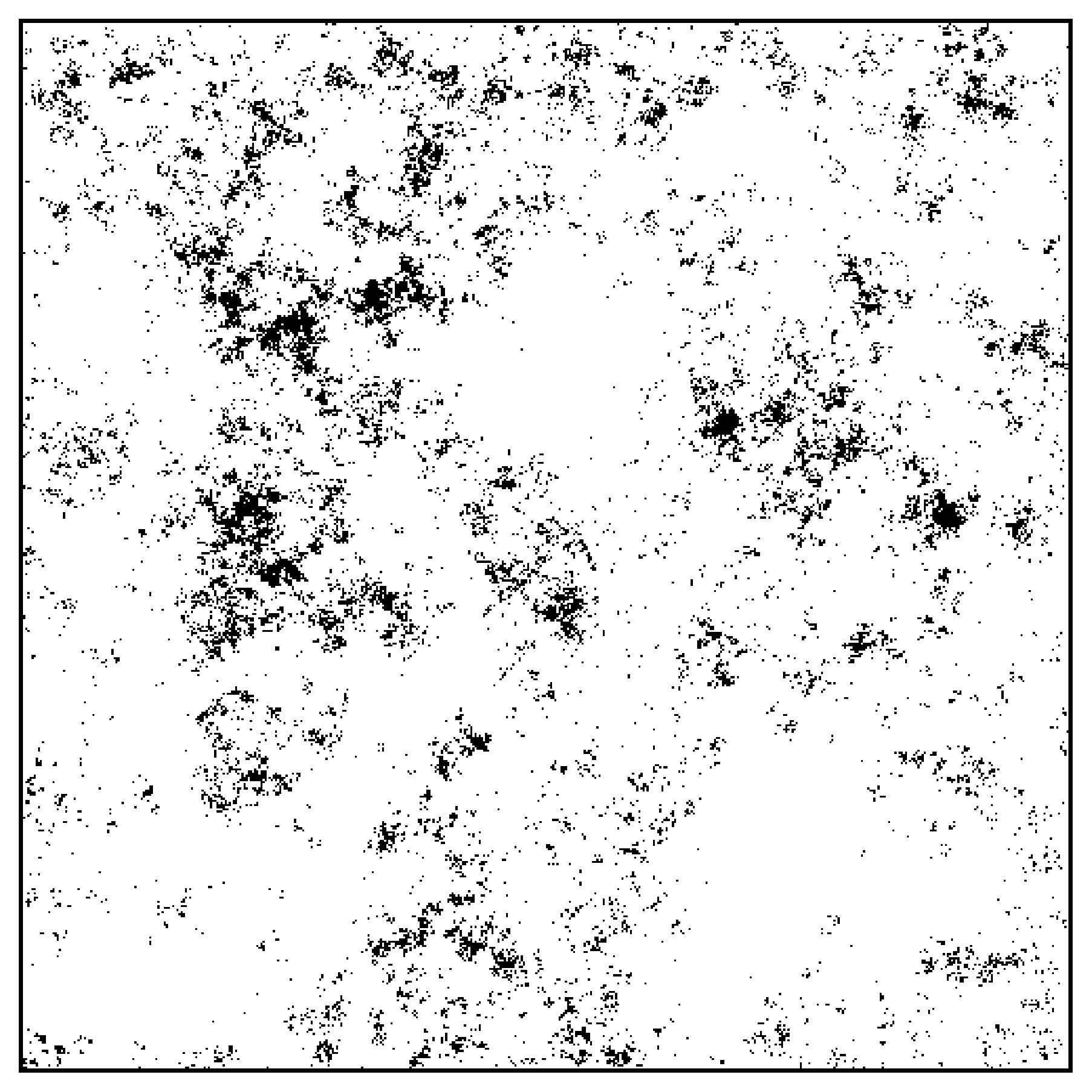}
\includegraphics[width=0.32\textwidth]{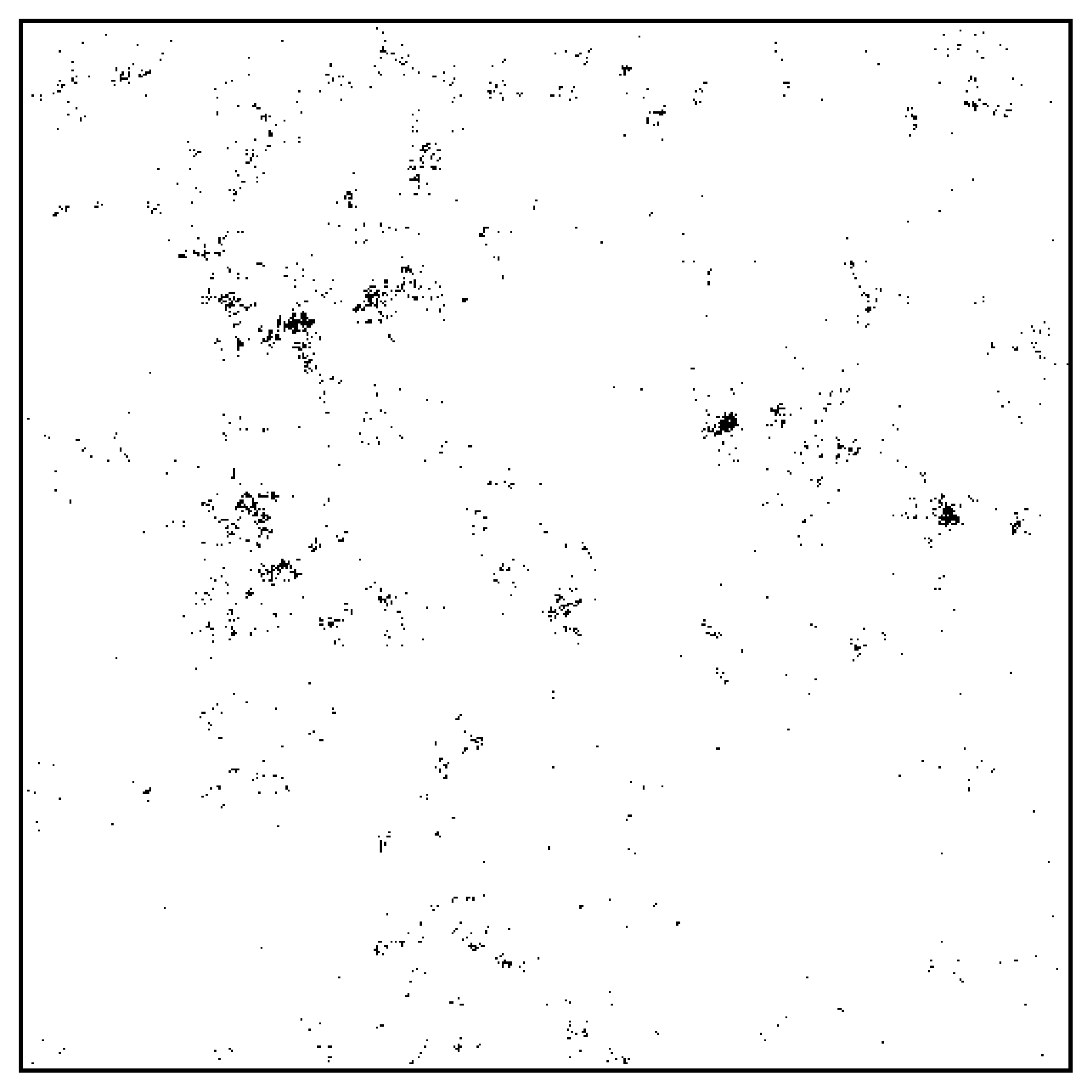}
}
\begin{quote}
\small 
\vglue1mm
\caption{\small
\label{fig-interlevel}
Plots of the points in the sample of the DGFF in Fig.~\ref{fig-DGFF} at heights (labeled left to right) above $0.1$, $0.3$ and $0.5$-multiples of the absolute maximum, respectively. Higher level sets are too sparse to produce a visible~effect.}
\normalsize
\end{quote}
\end{figure}

Our interest in these lectures is somewhat different as we wish to look at the level sets at heights that scale proportionally to the absolute maximum. We call these the \emph{intermediate level sets} although the term \emph{thick points} is quite common as well. Samples of such level sets are shown in Fig.~\ref{fig-interlevel}.

\nopagebreak


The self-similar nature of the plots in Fig.~\ref{fig-interlevel} is quite apparent. This motivates the following questions:
\begin{itemize}
\item
Is there a way to take a scaling limit of the samples in Fig.~\ref{fig-interlevel}?
\item
And if so, is there a way to characterize the limit object directly?
\end{itemize}
Our motivation for these questions stems from Donsker's Invariance Principle for random walks. There one first answers the second question by constructing a limit process; namely, the Brownian Motion. Then one proves that, under the diffusive scaling of space and time, all random walks whose increments have zero mean and finite second moment scale to that Brownian motion.

The goal of this and the next couple of lectures is to answer the above questions for the intermediate level sets of the DGFF. We focus only on one underlying discrete process so this can hardly be sold as a full-fledged analogue of Donsker's Invariance Principle. Still, we will perform the analysis over a whole class of  continuum domains thus achieving some level of \emph{universality}. The spirit of the two results is thus quite similar.

\section{Growth of absolute maximum}
\noindent
In order to set the scales for our future discussion, we first have to identify the growth rate of the absolute maximum. Here an early result of Bolthausen, Deuschel and Giacomin~\cite{BDeuG} established the leading-order asymptotic in square boxes. Their result reads: 

\begin{mytheorem}[Growth of absolute maximum]
\label{thm-2.1}
For $V_N:=(0,N)^2\cap\Z^2$,
\begin{equation}
\label{E:2.1}
\max_{x\in V_N}h^{V_N}_x = \bigl(2\sqrt g+o(1)\bigr)\log N
\end{equation}
where~$o(1)\to0$ in probability as~$N\to\infty$.
\end{mytheorem}

\begin{proofsect}{Proof of upper bound in \eqref{E:2.1}}
We start by noting the well-known tail estimate for centered normal random variables:

\begin{myexercise}[Standard Gaussian tail bound]
\label{ex:2.2} Prove that, for any $\sigma>0$,
\begin{equation}
\label{E:2.2}
Z\,\,\laweq\,\,\NN(0,\sigma^2)
\quad\Rightarrow\quad P(Z>a)\le\frac{\sigma}{\sigma+a}\,\texte^{-\frac{a^2}{2\sigma^2}},\qquad a>0.
\end{equation}
\end{myexercise}

\noindent
We want to use this for~$Z$ replaced by~$h^{V_N}_x$ but for that we need  to bound the variance of $h^{V_N}_x$ uniformly in~$x\in V_N$. Here we observe that, by the monotonicity of~$V\mapsto G^V(x,x)$ and translation invariance \eqref{E:G-trans}, denoting $\wt V_N:=(-N/2,N/2)^2\cap\Z^2$, there is a~$c\ge0$ such that
\begin{equation}
\label{E:2.3uai}
\max_{x\in V_N}\Var(h^{V_N}_x)\le \Var\bigl(h^{\wt V_{2N}}_0\bigr)\le g\log N+c,
\end{equation}
where the last bound follows from the asymptotic in Theorem~\ref{thm-1.17}. (We have just given away the solution to Exercise~\ref{ex:1.28}.) Plugging this in \eqref{E:2.2}, for any~$\theta>0$ we thus get
\begin{equation}
P\bigl(h_x^{V_N}>\theta\log N\bigr)\le\exp\Bigl\{-\frac12\frac{\theta^2(\log N)^2}{g\log N+c}\Bigr\}.
\end{equation}
Using that $(1+\lambda)^{-1}\ge1-\lambda$ for~$\lambda\in(0,1)$ we obtain
\begin{equation}
\frac1{g\log N+c}\ge\frac1{g\log N}-\frac c{(g\log N)^2}
\end{equation}
as soon as~$N$ is sufficiently large. Then
\begin{equation}
\label{E:2.6}
\max_{x\in V_N}\,P\bigl(h_x^{V_N}>\theta\log N\bigr)\le c' N^{-\frac{\theta^2}{2g}}
\end{equation}
for~$c':=\texte^{\theta^2c/(2g^2)}$ as soon as~$N$ is large enough. The union bound and the fact that~$|V_N|\le N^2$ then give
\begin{equation}
\label{E:2.7}
\begin{aligned}
P\Bigl(\,\max_{x\in V_N}h^{V_N}_x>\theta\log N\Bigr)
&\le \sum_{x\in V_N}P\bigl(h_x^{V_N}>\theta\log N\bigr)
\\
&\le c'|V_N|N^{-\frac{\theta^2}{2g}}=c'N^{2-\frac{\theta^2}{2g}}.
\end{aligned}
\end{equation}
This tends to zero as~$N\to\infty$ for any~$\theta>2\sqrt g$ thus proving~``$\le$'' in \eqref{E:2.1}.
\end{proofsect}

The proof of the complementary lower bound is considerably harder. The idea is to use the second-moment method but that requires working with a scale decomposition of the DGFF and computing the second moment under a suitable truncation. We will not perform this calculation here as the result will follow as a corollary from Theorem~\ref{thm-7.3}. (Section~\ref{sec:truncate} gives some hints how the truncation is performed.)

Building on~\cite{BDeuG}, Daviaud~\cite{Daviaud} was able to extend the control to the level sets of the form
\begin{equation}
\label{E:2.8}
\bigl\{x\in V_N\colon h^{V_N}_x\ge 2\sqrt g\,\lambda\log N\bigr\}\,,
\end{equation}
where~$\lambda\in(0,1)$. The relevant portion of his result reads:

\begin{mytheorem}[Size of intermediate level sets]
\label{thm-2.2}
For any~$\lambda\in(0,1)$,
\begin{equation}
\label{E:2.10}
\#\bigl\{x\in V_N\colon h^{V_N}_x\ge 2\sqrt g\,\lambda\log N\bigr\} = N^{2(1-\lambda^2)+o(1)}\,,
\end{equation}
where~$o(1)\to0$ in probability as~$N\to\infty$.
\end{mytheorem}

\begin{proofsect}{Proof of $``\le$'' in \eqref{E:2.10}}
Let~$L_N$ denote the cardinality of the set in \eqref{E:2.8}. Using the Markov inequality and the reasoning \twoeqref{E:2.6}{E:2.7},
\begin{equation}
\begin{aligned}
P\bigl(L_N\ge N^{2(1-\lambda^2)+\epsilon}\bigr)
&\le N^{-2(1-\lambda^2)-\epsilon}E(L_N)
\\
&\le c' N^{-2(1-\lambda^2)-\epsilon}N^{2-2\lambda^2}=c' N^{-\epsilon}.
\end{aligned}
\end{equation}
This tends to zero as~$N\to\infty$ for any~$\epsilon>0$ thus proving ``$\le$'' in \eqref{E:2.10}.
\end{proofsect}

We will not give a full proof of the lower bound for all~$\lambda\in(0,1)$ as that requires similar truncations as the corresponding bound for the maximum. However, these truncations are avoidable for~$\lambda$ sufficiently small, so we will content ourselves with:

\begin{proofsect}{Proof of ``$\ge$'' in \eqref{E:2.10} with positive probability for~$\lambda<1/\sqrt2$}
Define
\begin{equation}
Y_N:=\sum_{x\in V_N^\epsilon}\texte^{\beta h^{V_N}_x}\,,
\end{equation}
where~$\beta>0$ is a parameter to be adjusted later and~$V_N^\epsilon:=(\epsilon N,(1-\epsilon)N)^2\cap\Z^2$ for some~$\epsilon\in(0,1/2)$ to be fixed for the rest of the calculation. The quantity~$Y_N$ may be thought of as the normalizing factor (a.k.a.\ the partition function) for the Gibbs measure on~$V_N$ where the ``state''~$x$ gets ``energy''~$h^{V_N}_x$. (See Section~\ref{sec-11.4} for more on this.) Our first observation is:

\begin{mylemma}
For~$\beta>0$ such that $\beta^2g<2$ there is~$c=c(\beta)>0$ such that
\begin{equation}
\label{E:2.12new}
P\bigl(Y_N\ge cN^{2+\frac12\beta^2g}\bigr)\ge c
\end{equation}
once~$N$ is sufficiently large.
\end{mylemma}

\begin{proofsect}{Proof}
We will prove this by invoking the \emph{second moment method} whose driving force is the following inequality:

\begin{myexercise}[Second moment estimate]
\label{ex:2.5}
Let~$Y\in L^2$ be a non-negative random variable with~$EY>0$. Prove that
\begin{equation}
\label{E:2.13}
P\bigl(Y\ge q EY\bigr)\ge(1-q)^2\frac{[E(Y)]^2}{E(Y^2)},\quad q\in(0,1).
\end{equation}
\end{myexercise}

\noindent
In order to make use of \eqref{E:2.13} we have to prove that the second moments of~$Y_N$ is of the same order as the first moment squared.
We begin by the first moment of~$Y_N$. The fact that $E\texte^X=\texte^{EX+\frac12\Var(X)}$ for any~$X$ normal yields
\begin{equation}
\label{E:2.14}
EY_N=\sum_{x\in V_N^\epsilon}\texte^{\frac12\beta^2\Var(h^{V_{N}}_x)}.
\end{equation}
Writing $\wt V_N:=(-N/2,N/2)^2\cap\Z^2$, the monotonicity of~$V\mapsto G^V(x,x)$ gives
$\Var(h^{\wt V_{\epsilon N}}_0)\le\Var(h^{V_N}_x)\le \Var(h^{\wt V_{2N}}_0)$.
Theorem~\ref{thm-1.17} then implies
\begin{equation}
\label{E:2.16}
\sup_{N\ge1}\,\,\max_{x\in V_N^\epsilon}\,\,\bigl|\Var(h^{V_N}_x)-g\log N\bigr|<\infty
\end{equation}
As~$|V_N^\epsilon|$ is of order~$N^2$, using this in \eqref{E:2.14} we conclude that
\begin{equation}
\label{E:2.17u}
cN^{2+\frac12\beta^2g}\le EY_N\le c^{-1}N^{2+\frac12\beta^2g}
\end{equation}
holds for some constant~$c\in(0,1)$ and all~$N\ge1$.
 
Next we will estimate the second moment of~$Y_N$. Recalling the notation~$G^{V_N}$ for the Green function in~$V_N$, we have
\begin{equation}
E(Y_N^2)=\sum_{x,y\in V_N^\epsilon}\texte^{\frac12\beta^2[G^{V_N}(x,x)+G^{V_N}(y,y)+2G^{V_N}(x,y)]}\,.
\end{equation}
Invoking \eqref{E:2.16} and \eqref{E:1.44a} we thus get
\begin{equation}
E(Y_N^2)\le c'N^{\beta^2 g}\sum_{x,y\in V_N^\epsilon}\biggl(\,\frac N{|x-y|\vee1}\biggr)^{\beta^2g}\,.
\end{equation}
For~$\beta^2g<2$ the sum is dominated by pairs~$x$ and~$y$ with~$|x-y|$ of order~$N$. The sum is thus of order~$N^4$ and so we conclude
\begin{equation}
E(Y_N^2)\le c''N^{\beta^2 g+4}
\end{equation}
for some constant~$c''>0$. By \eqref{E:2.17u}, this bound is proportional to $[EY_N]^2$ so using this in \eqref{E:2.13} (with, e.g.,~$q:=1/2)$ readily yields \eqref{E:2.12new}.
\end{proofsect}

Next we will need to observe that the main contribution to~$Y_N$ comes from the set of points where the field roughly equals~$\beta g\log N$:

\begin{mylemma}
For any~$\delta>0$,
\begin{equation}
P\biggl(\,\sum_{x\in V_N^\epsilon}1_{\{|h^{V_N}_x-\beta g\log N|>(\log N)^{1/2+\delta}\}}\,\texte^{\beta h^{V_N}_x}\ge\delta N^{2+\frac12\beta^2g}\biggr)\,\underset{N\to\infty}\longrightarrow\,0.
\end{equation}
\end{mylemma}

\begin{proofsect}{Proof}
By \eqref{E:2.16}, we may instead prove this for $\beta g\log N$ replaced by $\beta\Var(h^{V_N}_x)$. Using the Markov inequality, the probability is then bounded by
\begin{equation}
\frac1{\delta N^{2+\frac12\beta^2g}}\sum_{x\in V_N^\epsilon}E\Bigl(1_{\{|h^{V_N}_x-\beta\Var(h^{V_N}_x)|>(\log N)^{1/2+\delta}\}}\,\texte^{\beta h^{V_N}_x}\Bigr)\,.
\end{equation}
Changing variables inside the (single-variable) Gaussian integral gives
\begin{multline}
\quad
E\Bigl(1_{\{|h^{V_N}_x-\beta\Var(h^{V_N}_x)|>(\log N)^{1/2+\delta}\}}\,\texte^{\beta h^{V_N}_x}\Bigr)
\\
=\texte^{\frac12\beta^2\Var(h^{V_N}_x)}P\bigl(|h^{V_N}_x|>(\log N)^{1/2+\delta}\bigr)
\\
\le c N^{\frac12\beta^2g}\texte^{-c''(\log N)^{2\delta}}
\quad
\end{multline}
for some $c,c''>0$, where we used again \eqref{E:2.16}. The probability in the statement is thus at most a constant times $\delta^{-1}\texte^{-c''(\log N)^{2\delta}}$, which vanishes as~$N\to\infty$.
\end{proofsect}

\noindent
Continuing the proof of ``$\ge$'' in \eqref{E:2.10},  the above lemmas yield
\begin{equation}
P\biggl(\,\sum_{x\in V_N^\epsilon}1_{\{h^{V_N}_x\ge \beta g\log N-(\log N)^{1/2+\delta}\}}
\ge\frac c2 N^{2+\frac12\beta^2g-\beta^2 g}\texte^{-\beta(\log N)^{1/2+\delta}}\biggr)
\ge\frac c2
\end{equation}
as soon as~$N$ is sufficiently large. Choosing~$\beta$ so that
\begin{equation}
\beta g\log N-(\log N)^{1/2+\delta}=2\sqrt g\lambda\log N
\end{equation}
gives $2-\frac12\beta^2 g = 2(1-\lambda^2)+O((\log N)^{-1/2+\delta})$ and so, since~$V_N^\epsilon\subset V_N$, the cardinality~$L_N$ of the level set \eqref{E:2.8} obeys
\begin{equation}
P\bigl(L_N\ge N^{2(1-\lambda^2)-c'(\log N)^{1/2+\delta}}\bigr)\ge\frac c2
\end{equation}
for some constant~$c'\in\R$ once~$N$ is large enough. This implies ``$\ge$'' in \eqref{E:2.10} with $o(1)\to0$ with a uniformly positive probability. The proof used that~$\beta^2g<2$, which means that it applies only to~$\lambda<1/\sqrt2$.
\end{proofsect}

Having the lower bound with a uniformly positive probability is actually sufficient to complete the proof of \eqref{E:2.10} as stated (when~$\lambda<1/\sqrt 2$). The key additional tool needed for this is the Gibbs-Markov decomposition of the DGFF which will be discussed in the next lecture. (See Exercise~\ref{ex:GM-bootstrap}.)

\begin{figure}[t]
\centerline{\includegraphics[width=0.58\textwidth]{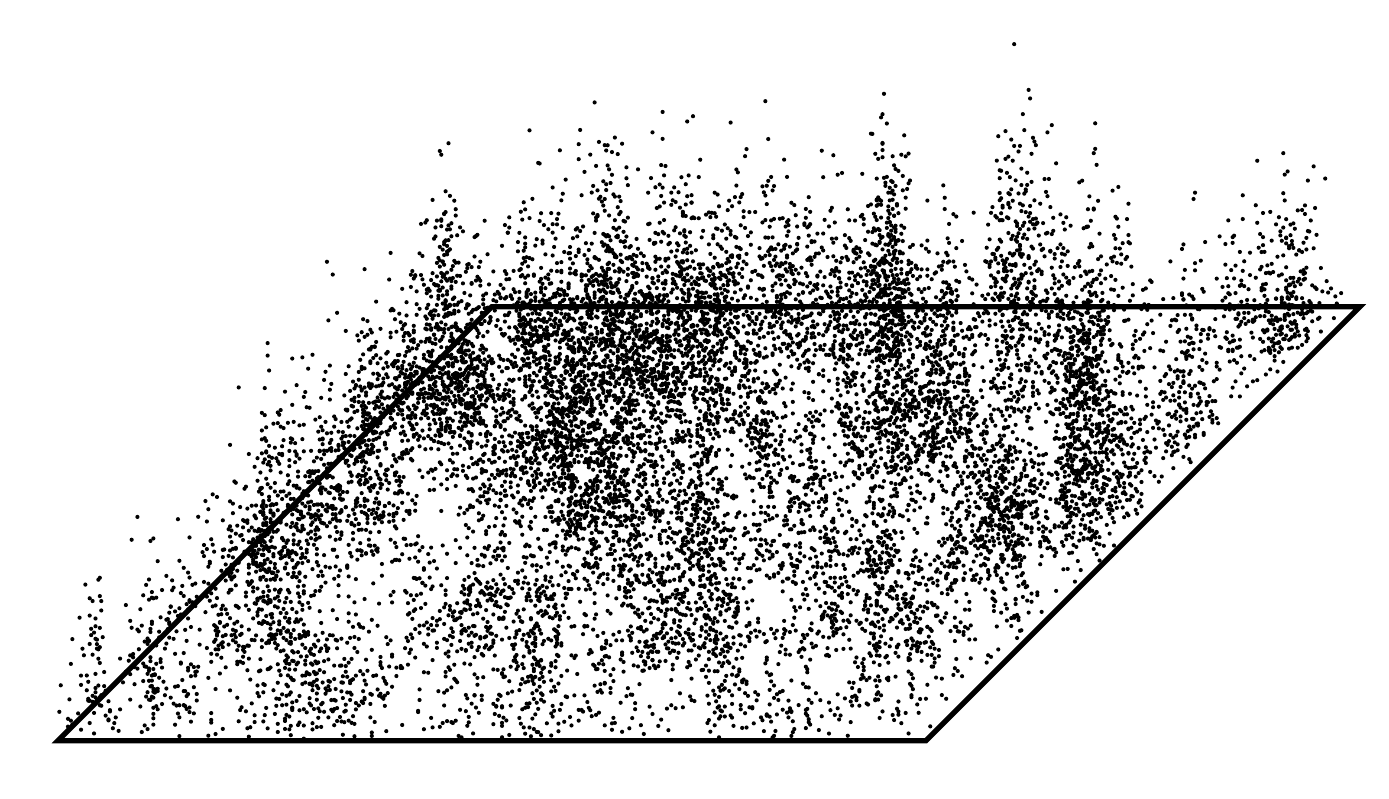}
\hglue-2.2cm
\includegraphics[width=0.58\textwidth]{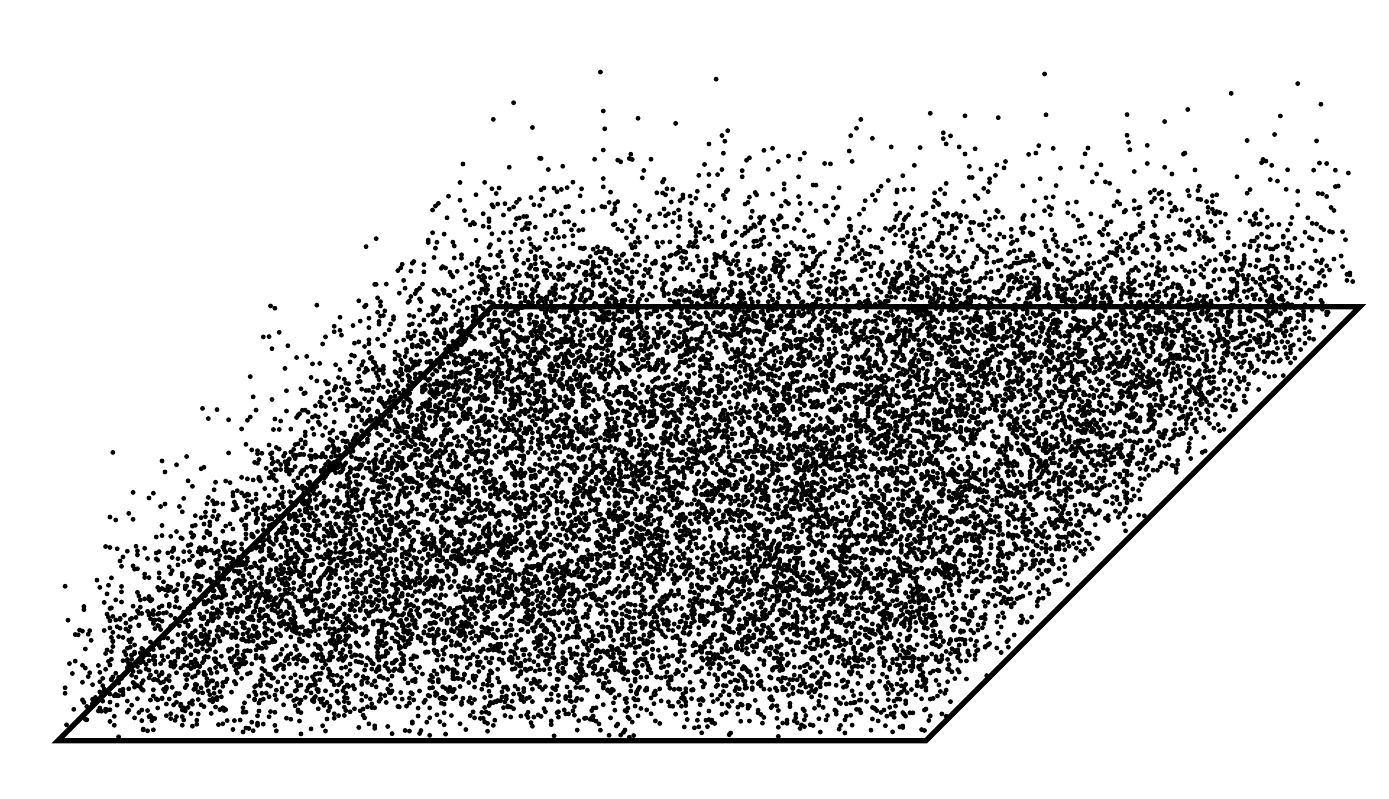}
}
\begin{quote}
\small 
\vglue-0.2cm
\caption{\small
\label{fig-PPP} 
Left: A sample of the level set \eqref{E:2.8}, or rather the point measure \eqref{E:2.27}, on a square of side~$N:=300$ with~$\lambda:=0.2$. Right: A corresponding sample for i.i.d.\ normals with mean zero and variance~$g\log N$. Although both samples live on the same ``vertical scale'', their local structure is very different.}
\normalsize
\end{quote}
\end{figure}

It is actually quite remarkable that the first-moment calculation alone is able to nail the correct leading order of the maximum as well as the asymptotic size of the level set~\eqref{E:2.8}. As that calculation did not involve correlations between the DGFF at distinct vertices, the same estimate would apply to i.i.d.\ Gaussians with the same growth rate of the variances. This (and many subsequent derivations) may lead one to think that the extreme values behave somehow like those of i.i.d.\ Gaussians. However, although some connection does exist, this is very far from the truth, as seen in Fig.~\ref{fig-PPP}.

The factor~$1-\lambda^2$ in the exponent is ubiquitous in this subject area. Indeed, it appears in the celebrated analysis of the Brownian \textit{fast points} (Orey and Taylor~\cite{Orey-Taylor}) and, as was just noted, for i.i.d.\ Gaussians with variance~$g\log N$. A paper by Chatterjee, Dembo and Ding~\cite{Chatterjee-Dembo-Ding} gives (generous) conditions under which such a factor should be expected. For suitably formulated analogue of Daviaud's level sets, called the \emph{thick points}, of the two-dimensional CGFF, this has been shown by Hu, Miller and Peres~\cite{HMP10}.

\section{Intermediate level sets}
\noindent
The main objective in this part of the course is to show that the intermediate level sets~\eqref{E:2.8} admit a non-trivial \emph{scaling limit} whose law can be explicitly characterized. We start by pondering about the formulation that makes taking a scaling limit meaningful. Not all natural ideas may work; for instance,  scaling the box down to a unit size, the set \eqref{E:2.8} becomes increasingly dense everywhere so taking its limit using, e.g., the topology of Hausdorff convergence does not seem useful. A better idea here is to encode the set into the point measure on~$[0,1]^2\times\R$ of the form
\begin{equation}
\label{E:2.27}
\sum_{x\in  V_N}\delta_{x/N}\otimes\delta_{h_x^{V_N}-a_N}\,,
\end{equation}
where, to allow for some generalizations, $a_N$ is a scale sequence such that
\begin{equation}
\label{E:2.28}
\frac{a_N}{\log N}\,\,\underset{N\to\infty}\longrightarrow\,\,2\sqrt g\,\lambda
\end{equation}
for some~$\lambda\in(0,1)$. A sample of the measure \eqref{E:2.27} is shown on the left of Fig.~\ref{fig-PPP}. 

By Theorem~\ref{thm-2.2}, the measures in \eqref{E:2.27} assign unbounded mass to bounded intervals in the second variable, and so a normalization is required. We will show that this can be done (somewhat surprisingly) by a \emph{deterministic} sequence of the form
\begin{equation}
K_N:=\frac{N^2}{\sqrt{\log N}}\,\,\texte^{-\frac{a_N^2}{2g\log N}}\,.
\end{equation}
Note that \eqref{E:2.28} implies $K_N=N^{2(1-\lambda^2)+o(1)}$ so the normalization is consistent with Theorem~\ref{thm-2.2}.
Our main result, proved in Biskup and Louidor~\cite{BL4}, is then:

\begin{mytheorem}[Scaling limit of intermediate level sets]
\label{thm-intermediate}
For each~$\lambda\in(0,1)$ and each $D\in\mathfrak D$ there is an a.s.-finite random Borel measure~$Z^D_\lambda$ on~$D$ such that for any~$a_N$ satisfying \eqref{E:2.28} and any admissible sequence $\{D_N\colon N\ge1\}$ of lattice approximations of~$D$, the normalized point measure
\begin{equation}
\label{E:2.30ua}
\eta_N^D:=\frac1{K_N}\sum_{x\in D_N}\delta_{x/N}\otimes\delta_{h^{D_N}_x-a_N}
\end{equation}
obeys
\begin{equation}
\label{E:2.31}
\eta^D_{N}\,\,\,\underset{N\to\infty}\lawarrow\,\,\,Z_\lambda^D(\textd x)\,\otimes\,\texte^{-\alpha\lambda h}\textd h,
\end{equation}
where $\alpha:=2/\sqrt g$. Moreover, $Z^D_\lambda(A)>0$ a.s.\ for every non-empty open~$A\subseteq D$.
\end{mytheorem}

A remark is perhaps in order on what it means for random measures to converge in law. The space of Radon measures on~$D\times\R$ (of which~$\eta^D_N$ is an example) is naturally endowed with the topology of vague convergence. This topology makes the space of Radon measures a Polish space which permits discussion of distributional convergence. We pose:

\begin{myexercise}
\label{ex:2.8nwt}
Let~$\scrX$ be a Polish space and~$\MM(\scrX)$ the space of Radon measures on~$\scrX$ endowed with the vague topology. Prove that a sequence of random measures~$\eta_N\in\MM(\scrX)$ converges in law to~$\eta$ if and only if for each~$f\in C_\cc(\XX)$, 
\begin{equation}
\langle \eta_N,f\rangle \,\,\,\underset{N\to\infty}\lawarrow\,\,\, \langle\eta,f\rangle,
\end{equation}
where~$\langle \eta,f\rangle$ denotes the integral of~$f$ with respect to~$\eta$ and the convergence in law is in the sense of ordinary $\R$-valued random variables.
\end{myexercise}

A subtlety of the above theorem is that the convergence actually happens over a larger ``base'' space, namely, $\overline D\times(\R\cup\{+\infty\})$ implying, in particular, $Z^D_\lambda(\partial D)=0$ a.s. This means that we get weak convergence of the said integral even for  (continuous) functions that take non-zero values on~$\partial D$ in the~$x$-variable and/or at~$+\infty$ in the~$h$-variable. This implies:

\begin{mycorollary}
\label{cor-2.6}
For the setting of Theorem~\ref{thm-intermediate},
\begin{equation}
\label{E:2.32}
\frac1{K_N}\#\bigl\{x\in D_N\colon h^{D_N}_x\ge a_N\bigr\}
\,\,\,\underset{N\to\infty}\lawarrow\,\,\, (\alpha\lambda)^{-1}\, Z_\lambda^D(D).
\end{equation}
\end{mycorollary}

\begin{proofsect}{Proof (idea)}
Use Theorem~\ref{thm-intermediate} for~$\eta_N^D$ integrated against $f(x,y):=1_{[0,\infty)}(h)$.
\end{proofsect}

\begin{myexercise}
\label{ex2.10}
Apply suitable monotone limits to check that the convergence in \eqref{E:2.31} --- which involves \emph{a priori} only integrals of these measures with respect to compactly-supported continuous functions --- can be applied to functions of the form
\begin{equation}
f(x,h):=\tilde f(x)1_{[a,b]}(h)
\end{equation}
 for~$\tilde f\colon \overline D\to\R$ continuous and~$a<b$ (including~$b:=\infty$).
\end{myexercise}

We remark that Corollary~\ref{cor-2.6} extends, quite considerably, Theorem~\ref{thm-2.2} originally proved by Daviaud~\cite{Daviaud}.
In order to get some feeling for what Theorem~\ref{thm-intermediate} says about the positions of the points in the level set, we also state:

\begin{mycorollary}
For~$a_N$ as in \eqref{E:2.28}, given a sample of~$h^{D_N}$, let~$X_N$ be a point chosen uniformly from $\{x\in D_N\colon\,h^{D_N}_x\ge a_N\}$. Then
\begin{equation}
\label{E:2.34ua}
\frac1NX_N\,\,\,\underset{N\to\infty}\lawarrow\,\,\,\wh X\quad\text{with}\quad\text{\rm law}(\wh X)=E\biggl[\frac{Z^D_\lambda(\cdot)}{Z^D_\lambda(D)}\biggr].
\end{equation}
In fact, the joint law of~$X_N$ and~$\eta^D_N$ obeys
\begin{equation}
\label{E:2.35uai}
(N^{-1}X_N,\eta_N^D)\,\,\,\underset{N\to\infty}\lawarrow\,\,\,(\wh X,\eta^D),
\end{equation}
where the marginal law of~$\eta^D$ is that of $Z^D_\lambda(\textd x)\otimes\texte^{-\alpha\lambda h}\textd h$ and the law of~$\wh X$ conditional on~$\eta^D$ is given by $Z^D_\lambda(\cdot)/Z^D_\lambda(D)$.
\end{mycorollary}

\begin{proofsect}{Proof}
We easily check that, for any~$f\colon \overline D\to\R$ continuous (and thus bounded),
\begin{equation}
\label{E:2.36ua}
E\bigl[\,f(X_N/N)\bigr]=E\Biggl[\frac{\langle\eta^D_N,f\otimes 1_{[0,\infty)}\rangle}{\langle\eta^D_N,1_{[0,\infty)}\rangle}\Biggr]\,,
\end{equation}
where $(f\otimes 1_{[0,\infty)})(x,h):=f(x)1_{[0,\infty)}(h)$ and the brackets denote the integral of the function with respect to the measure. Applying Theorem~\ref{thm-intermediate}, we get
\begin{equation}
\label{E:2.35}
\frac{\langle\eta^D_N,f\otimes 1_{[0,\infty)}\rangle}{\langle\eta^D_N,1\otimes 1_{[0,\infty)}\rangle}
\,\,\,\underset{N\to\infty}\lawarrow\,\,\,
\frac{\displaystyle\int_D Z^D_\lambda(\textd x)\,f(x)}{Z^D_\lambda(D)}.
\end{equation}
This yields \eqref{E:2.34ua}. The more general clause \eqref{E:2.35uai} is then proved by considering test functions in both variables in \eqref{E:2.36ua} and proceeding as above.
\end{proofsect}

\begin{myexercise}
The statement \eqref{E:2.35} harbors a technical caveat: we are taking the distributional limit of a ratio of two random variables, each of which converges (separately) in law. Fill the details needed to justify the conclusion.
\end{myexercise}

We conclude that the spatial part of the right-hand side of \eqref{E:2.31} thus tells us about the local ``intensity'' of the sets in the samples in Fig.~\ref{fig-PPP}. Concerning the values of the field, we may be tempted to say that these are Gumbel ``distributed'' with decay exponent~$\alpha\lambda$. This is not justified by the statement \emph{per se} as the measure on the right of \eqref{E:2.31} is not a probability (it is not even finite). Still, one can perhaps relate this to the corresponding problem for i.i.d.\ Gaussians with variance~$g\log N$. Indeed, as an exercise in extreme-value statistics, we pose:

\begin{myexercise}
Consider the measure~$\eta^D_N$ for~$h^{D_N}$ replaced by i.i.d.\ Gaussians with mean zero and variance~$g\log N$.
Prove the same limit as in \eqref{E:2.31}, with the same~$K_N$ and but with~$Z^D_\lambda$ replaced by a multiple of the Lebesgue measure on~$D$. 
\end{myexercise}

\noindent
We rush to add that (as we will explain later) $Z^D_\lambda$ is a.s.\ singular with respect to the Lebesgue measure. Compare, one more time, the two samples in Fig.~\ref{fig-PPP}.

\section{Link to Liouville Quantum Gravity}
\noindent
Not too surprisingly, the random measures $\{Z^D_\lambda\colon D\in\mathfrak D\}$ (or, rather, their laws) are very closely related. We will later give a list of properties that characterize these laws uniquely. From these properties one can derive the following transformation rule for conformal maps between admissible domains: 

\begin{mytheorem}[Conformal covariance]
\label{thm-2.4}
Let~$\lambda\in(0,1)$. Under any conformal bijection $f\colon D\to f(D)$ between the admissible domains $D,f(D)\in\mathfrak D$, the laws of the above measures transform as
\begin{equation}
\label{E:2.16a}
Z_\lambda^{f(D)}\circ f(\textd x)\,\,\laweq\,\, |f'(x)|^{2+2\lambda^2}\, Z_\lambda^D(\textd x).
\end{equation}
\end{mytheorem}

Recall that~$r_D(x)$, defined in \eqref{E:rad}, denotes the conformal radius of~$D$ from~$x$. The following is now a simple consequence of the above theorem:

\begin{myexercise}
Show that in the class of admissible~$D$, the law of
\begin{equation}
\label{E:2.37uua}
\frac1{r_D(x)^{2+2\lambda^2}}\,Z^D_\lambda(\textd x)
\end{equation}
is invariant under conformal maps.
\end{myexercise}

In light of Exercise~\ref{ex:1.15}, the law of~$Z^D_\lambda$ for any bounded, open, simply connected~$D$ can thus be reconstructed from the law on, say, the open unit disc. As we will indicate later, such a link exists for all admissible~$D$. However, that still leaves us with the task of determining the law of~$Z^D$ for at least one admissible domain. We will instead give an independent construction of the law of~$Z^D_\lambda$ that works for general~$D$. This will also give us the opportunity to review some of the ideas of Kahane's theory of multiplicative chaos. 

Let $\cmss H_0^1(D)$ denote the closure of the set of smooth functions with compact support in~$D$ in the topology of the Dirichlet inner product
\begin{equation}
\langle f,g\rangle_\nabla:= \frac14\int_D\nabla f(x)\cdot\nabla g(x)\,\textd x\,,
\end{equation}
where $\nabla f$ is now the ordinary (continuum) gradient and $\cdot$ denotes the Euclidean scalar product in~$\R^2$. For $\{X_n\colon n\ge1\}$ i.i.d.\ standard normals and $\{f_n\colon n\ge1\}$ an orthonormal basis  in $\cmss H_0^1(D)$, let
\begin{equation}
\label{E:2.17}
\varphi_n(x):=\sum_{k=1}^n X_k f_k(x).
\end{equation}
These are to be thought of as regularizations of the CGFF:

\begin{myexercise}[Hilbert space definition of CGFF]
\label{ex:2.15}
For any smooth function~$f\in\cmss H_0^1(D)$, let $\Phi_n(f):=\int_D f(x)\varphi_n(x)\textd x$. Show that~$\Phi_n(f)$ converges, as~$n\to\infty$, in~$L^2$ to a CGFF in the sense of Definition~\ref{def-1.28}. 
\end{myexercise}

Using the above random fields, for each~$\beta\in[0,\infty)$, we then define the random measure
\begin{equation}
\label{E:2.17a}
\mu_n^{D,\beta}(\textd x):=1_D(x)\texte^{\beta\varphi_n(x)-\frac{\beta^2}2E[\varphi_n(x)^2]}\,\textd x.
\end{equation}
The following observation goes back to Kahane~\cite{Kahane} in 1985:

\begin{mylemma}[Gaussian Multiplicative Chaos]
\label{lemma-GMC}
There exists a random, a.s.-finite (albeit possibly trivial) Borel measure $\mu_\infty^{D,\beta}$ on~$D$ such that for  each Borel~$A\subseteq D$,
\begin{equation}
\label{E:2.38}
\mu_n^{D,\beta}(A)\,\,\underset{n\to\infty}\longrightarrow\,\,\mu_\infty^{D,\beta}(A),\quad \text{\rm a.s.}
\end{equation}
\end{mylemma}

\begin{proofsect}{Proof}
Pick~$A\subseteq D$ Borel measurable. First we will show that~$\mu_n^{D,\beta}(A)$ converge a.s. To this end, for each~$n\in\N$, define 
\begin{equation}
M_n:=\mu_n^{D,\beta}(A) \quad\text{and}\quad \FF_n:=\sigma(X_1,\dots,X_n).
\end{equation}
We claim that $\{M_n\colon n\ge1\}$ is a martingale with respect to~$\{\FF_n\colon n\ge1\}$. Indeed, using the regularity of the underlying measure space (to apply Fubini-Tonelli)
\begin{multline}
\label{E:2.39}
\quad
E(M_{n+1}|\FF_n)=E\bigl(\mu_{n+1}^{D,\beta}(A)\,\big|\,\FF_n\bigr)
\\
=\int_A\textd x \,E\bigl(\texte^{\beta\varphi_{n+1}(x)-\frac{\beta^2}2E[\varphi_{n+1}(x)^2]}\,\big|\,\FF_n\bigr)\,.
\quad
\end{multline}
The additive structure of~$\varphi_n$ now gives
\begin{multline}
\quad
E\bigl(\texte^{\beta\varphi_{n+1}(x)-\frac{\beta^2}2E[\varphi_{n+1}(x)^2]}\,\big|\,\FF_n\bigr)
\\
=\texte^{\beta\varphi_{n}(x)-\frac{\beta^2}2E[\varphi_{n}(x)^2]}
E\bigl(\texte^{\beta f_{n+1}(x)X_{n+1}-\frac12\beta^2 f_{n+1}(x)^2 E[X_{n+1}^2]}\bigr)
\\
=\texte^{\beta\varphi_{n}(x)-\frac{\beta^2}2E[\varphi_{n}(x)^2]}.
\quad
\end{multline}
Using this in \eqref{E:2.39}, the right-hand side then wraps back into~$\mu_{n}^{D,\beta}(A)=M_n$ and so $\{M_n\colon n\ge1\}$ is a martingale as claimed.

The martingale~$\{M_n\colon n\ge1\}$ is non-negative and so the Martingale Convergence Theorem yields~$M_n\to M_\infty$ a.s. (with the implicit null event depending on~$A$). In order to identify the limit in terms of a random measure, we have to  rerun the above argument as follows: For any bounded measurable~$f\colon D\to\R$ define
\begin{equation}
\phi_n(f):=\int f\,\textd \mu_{n}^{D,\beta}\,.
\end{equation}
 Then the same argument as above shows that~$\phi_n(f)$ is a bounded martingale and so~$\phi_n(f)\to \phi_\infty(f)$ a.s. Specializing to continuous~$f$, the bound
 \begin{equation}
\bigl|\phi_n(f)\bigr|\le\,\mu_{n}^{D,\beta}(D)\,\Vert f\Vert_{C(\overline D)}
\end{equation}
along with the above a.s.\ convergence~$M_n\to M_\infty$ (for~$A:=D$) yields
\begin{equation}
\bigl|\phi_\infty(f)\bigr|\le M_\infty\Vert f\Vert_{C(\overline D)}.
\end{equation}
Restricting to a countable dense subclass of~$f\in C(\overline D)$ to manage the proliferation of null sets, $f\mapsto \phi_\infty(f)$ extends to a continuous linear functional $\phi_\infty'$ on~$C(\overline D)$ a.s.\ such that $\phi_\infty(f)=\phi_\infty'(f)$ a.s.\ for each~$f\in C(\overline D)$ (the null set possibly depending on~$f$). On the event that~$\phi_\infty'$ is well-defined and continuous, the Riesz Representation Theorem yields existence of a (random) Borel measure~$\mu_\infty^{D,\beta}$ such that
\begin{equation}
\label{E:2.50uai}
\phi_\infty(f)\,\overset{\text{a.s.}}=\,\phi_\infty'(f)=\int f\,\textd \mu_\infty^{D,\beta}
\end{equation}
for each~$f\in C(\overline D)$. 

To identify the limit in \eqref{E:2.38} (which we already showed to exist a.s.) with $\mu_\infty^{D,\beta}(A)$ we proceed as follows. First, given a~$\GG_\delta$-set~$A\subseteq D$, we can find functions $f_k\in C(\overline D)$ such that~$f_k\downarrow 1_A$ as~$k\to\infty$. Since $\mu_n^{D,\beta}(A)\le\phi_n(f_k)\to\phi_\infty'(f_k)$ a.s.\ as~$n\to\infty$ and $\phi_\infty'(f_k)\downarrow\mu_\infty^{D,\beta}(A)$ as~$k\to\infty$ by the Bounded Convergence Theorem, we get
\begin{equation}
\label{E:conv}
\lim_{n\to\infty}\mu_n^{D,\beta}(A)\le\mu_\infty^{D,\beta}(A)\quad\text{\rm a.s.}
\end{equation}
Next, writing $L^1(D)$ for the space of Lebesgue absolutely integrable~$f\colon D\to\R$, Fatou's lemma shows
\begin{equation}
E\bigl|\phi_\infty'(f)\bigr|\le\Vert f\Vert_{L^1(D)},\quad f\in C(\overline D).
\end{equation}
The above approximation argument, along with the fact the Lebesgue measure is outer regular now shows that $\Leb(A)=0$ implies $\mu_\infty^{D,\beta}(A)=0$ a.s. This gives \eqref{E:conv} for all Borel~$A\subseteq D$. The string of equalities
\begin{equation}
\mu_n^{D,\beta}(A)+\mu_n^{D,\beta}(A^\cc)=\phi_n(1)\,\,\overset{\text{a.s.}}{\underset{n\to\infty}\longrightarrow}\,\,\phi_\infty'(1)=\mu_\infty^{D,\beta}(A)+\mu_\infty^{D,\beta}(A^\cc)
\end{equation}
now shows that equality must hold in \eqref{E:conv}.
\end{proofsect}

An interesting question is how the limit object~$\mu_\infty^{D,\beta}$ depends on~$\beta$ and, in particular, for what~$\beta$ it is non-trivial. We pose:

\begin{myexercise}
\label{ex:2.17ua}
Use Kolmogorov's zero-one law to show that, for each~$A\subseteq D$ measurable,
\begin{equation}
\bigl\{\mu_\infty^{D,\beta}(A)=0\bigr\}\text{ is a zero-one event}.
\end{equation}
\end{myexercise}

\begin{myexercise}
Prove that there is~$\beta_\cc\in[0,\infty]$ such that
\begin{equation}
\mu_\infty^{D,\beta}(D)\begin{cases}
>0,\qquad&\text{if }0\le \beta<\beta_\cc,
\\
=0,\qquad&\text{if }\beta>\beta_\cc.
\end{cases}
\end{equation}
Hint: Use Gaussian interpolation and conditional Fatou's lemma to show that the Laplace transform of $\mu_n^{D,\beta}(A)$ is non-decreasing in~$\beta$.
\end{myexercise}

We will show later that, in our setting (and with~$\alpha$ from Theorem~\ref{thm-intermediate}),
\begin{equation}
\beta_\cc:=\alpha=2/\sqrt g.
\end{equation}
As it turns out, the law of the measure~$\mu_\infty^{D,\beta}$ is independent of the choice of the underlying basis in~$\cmss H_0^1(D)$. This has been proved gradually starting with somewhat restrictive Kahane's theory~\cite{Kahane} (which we will review later) and culminating in a recent paper by Shamov~\cite{Shamov}. (See Theorems~\ref{thm-GMC-unique} and~\ref{thm-5.5}.)
 
\begin{figure}[t]
\vglue0.2cm
\centerline{\includegraphics[width=0.7\textwidth]{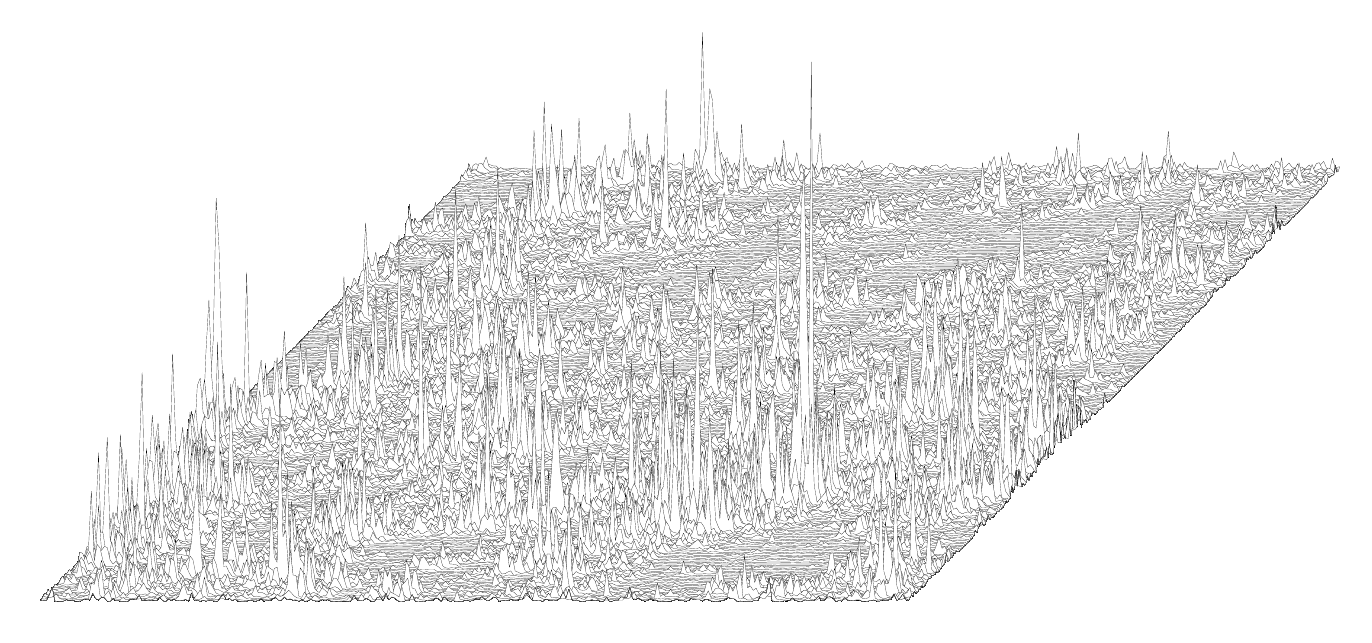}
}
\begin{quote}
\small 
\caption{\small
\label{fig3}
A sample of the LQG measure $r_D(x)^{2\lambda^2}\mu^{D,\lambda\alpha}_\infty(\textd x)$ for~$D$ a unit square and~$\lambda:=0.3$. The high points indicate places of high local intensity.}
\normalsize
\end{quote}
\end{figure}
 
The measure~$\mu_\infty^{D,\beta}$ is called the \emph{Gaussian multiplicative chaos} associated with the continuum Gaussian Free Field.  We now claim:

\begin{mytheorem}[$Z^D_\lambda$-measure as LQG measure]
\label{thm-LQG}
Assume the setting of Theorem~\ref{thm-intermediate} with~$\lambda\in(0,1)$.
Then there is $\hat c\in(0,\infty)$ such that for all~$D\in\mathfrak D$,
\begin{equation}
\label{E:2.43}
Z_\lambda^D(\textd x) \,\laweq\,\hat c\,r_D(x)^{2\lambda^2}\,\mu_\infty^{D,\,\lambda\alpha}(\textd x).
\end{equation} 
where, we recall,~$r_D(x)$ denotes the conformal radius of~$D$ from~$x$.
\end{mytheorem}

The measure on the right of \eqref{E:2.43} (without the constant~$\hat c$) is called the \emph{Liouville Quantum Gravity} (LQG) measure in~$D$ for parameter $\beta:=\lambda\alpha$. This object is currently heavily studied in connection with random conformally-invariant geometry (see, e.g., Miller and Sheffield~\cite{MS5,MS6}).

An alternative construction of the LQG measure was given by Duplantier and Sheffield~\cite{Duplantier-Sheffield} using disc/circle averages (cf Exercise~\ref{ex:1.29}). This construction is technically more demanding (as it is not based directly on martingale convergence theory) but, as a benefit, one gets some regularity of the limit. 

We note (and will prove this for~$\lambda<1/\sqrt2$ in Corollary~\ref{cor-3.17}) that,  
\begin{equation}
E\mu_\infty^{D,\lambda\alpha}(A)=\Leb(A),\quad \lambda\in[0,1).
\end{equation}
For each~$\lambda<1$, the total mass $\mu_\infty^{D,\lambda\alpha}(D)$ has moments up to $1+\epsilon(\lambda)$, for~$\epsilon(\lambda)>0$ that tends to zero as~$\lambda\uparrow1$. 

We will not discuss Gaussian Multiplicative Chaos and/or the LQG measures as a stand-alone topic much in these lectures (although these objects will keep popping up in our various theorems) and instead refer the reader to Berestycki~\cite{Berestycki-GMC} and the review by Rhodes and Vargas~\cite{RV-review}. The proofs of the Theorems~\ref{thm-intermediate}, \ref{thm-2.4} and~\ref{thm-LQG} will be given in the forthcoming lectures.


\chapter{Intermediate level sets: factorization}
\label{lec-3}
\noindent
The aim of this and the following lecture is to give a fairly detailed account of the proofs of the above theorems  on the scaling limit of the intermediate level sets. We will actually do this only in the regime where the second-moment calculations work without the need for truncations; this requires restricting to~$\lambda<1/\sqrt2$. We  comment on the changes that need to be made for the complementary set of~$\lambda$'s at the end of the next lecture.

\section{Gibbs-Markov property of DGFF}
\noindent
A number of forthcoming proofs will use a special property of the DGFF that addresses the behavior of the field restricted, via conditioning, to a subdomain. This property is the spatial analogue of the Markov property in one-parameter stochastic processes and is a consequence of the Gaussian decomposition into orthogonal subspaces along with the fact that the law of the DGFF is a Gibbs measure for a nearest-neighbor Hamiltonian (cf Definition~\ref{def-1.1}). For this reason, we will attach the adjective \emph{Gibbs-Markov} to this property, although the literature uses the term \emph{domain-Markov} as well. Here is the precise statement:

\begin{mylemma}[Gibbs-Markov property]
\label{lemma-GM}
For $U\subsetneq V\subsetneq\Z^2$, denote
\begin{equation}
\varphi^{V,U}_x:=E\bigl(h^V_x\,\big|\,\sigma(h^V_z\colon z\in V\smallsetminus U)\bigr).
\end{equation}
where~$h^V$ is the DGFF in~$V$. Then we have:
\begin{enumerate}
\item[(1)] A.e.\ sample of $x\mapsto\varphi^{V,U}_x$ is discrete harmonic on~$U$ with ``boundary values'' determined by $\varphi^{V,U}_x=h^V_x$ for each $x\in V\smallsetminus U$.
\item[(2)] The field $h^V-\varphi^{V,U}$ is independent of~$\varphi^{V,U}$ and
\begin{equation}
h^V-\varphi^{V,U} \,\,\laweq\,\, h^{U}\,.
\end{equation}
\end{enumerate}
\end{mylemma}

\begin{proofsect}{Proof}
Assume that~$V$ is finite for simplicity.
Conditioning a multivariate Gaussian on part of the values preserves the multivariate Gaussian nature of the law. Hence~$\varphi^{V,U}$ and~$h^V-\varphi^{V,U}$ are multivariate Gaussians that are, by the properties of the conditional expectation, uncorrelated. It follows that $\varphi^{V,U}\independent h^V-\varphi^{V,U}$.

Next let us prove that~$\varphi^{V,U}$ has discrete-harmonic sample paths in~$U$. To this end pick any~$x\in U$ and note that the ``smaller-always-wins'' principle for nested conditional expectations yields
\begin{equation}
\label{E:3.3}
\varphi^{V,U}_x=E\Bigl(E\bigl(h^V_x\,\big|\,\sigma(h^V_{z'}\colon z'\ne x)\bigr)\,\Big|\,\sigma(h^V_z\colon z\in V\smallsetminus U)\Bigr)\,.
\end{equation}
In light of Definition~\ref{def-1.1}, the inner conditional expectation admits the explicit form
\begin{equation}
\label{E:3.4}
E\bigl(h^V_x\,\big|\,\sigma(h^V_{z'}\colon z'\ne x)\bigr)
=\frac{\displaystyle\int_\R h_x\,\texte^{-\frac18\sum_{y\colon y\sim x}(h_y^V-h_x)^2}\textd h_x}
{\displaystyle\int_\R \texte^{-\frac18\sum_{y\colon y\sim x}(h_y^V-h_x)^2}\textd h_x},
\end{equation}
where~$y\sim x$ abbreviates $(x,y)\in E(\Z^2)$. Now
\begin{equation}
\begin{aligned}
\frac14\sum_{y\colon y\sim x}(h_y-h_x)^2
&=h_x^2-2\frac14\sum_{y\colon y\sim x}h_y+\frac14\sum_{y\colon y\sim x}h_y^2
\\
&=\Bigl(h_x-\frac14\sum_{y\colon y\sim x}h_y\Bigr)^2+\frac14\sum_{y\colon y\sim x}h_y^2-\Bigl(\frac14\sum_{y\colon y\sim x}h_y\Bigr)^2\,.
\end{aligned}
\end{equation}
The last two terms factor from both the numerator and denominator on the right of \eqref{E:3.4}. Shifting~$h_x$ by the average of the neighbors then gives
\begin{equation}
E\bigl(h^V_x\,\big|\,\sigma(h^V_{z'}\colon z'\ne x)\bigr) = \frac14\sum_{y\colon y\sim x}h_y^V\,.
\end{equation}
Using this in \eqref{E:3.3} shows that~$\varphi^{V,U}$ has the mean-value property, and is thus discrete harmonic, on~$U$.

Finally, we need to show that $\wt h^U:=h^V-\varphi^{V,U}$ has the law of~$h^U$. The mean of~$\wt h^U$ is zero so we just need to verify that the covariances match. Here we note that, using $H^U$ to denote the discrete harmonic measure on~$U$, the mean-value property of~$\varphi^{V,U}$ yields
\begin{equation}
\wt h^U_x = h^V_x-\sum_{z\in\partial U}H^U(x,z)h^V_z,\quad x\in U.
\end{equation}
For any~$x,y\in U$, this implies
\begin{multline}
\label{E:3.8}
\quad
\Cov(\wt h^U_x,\wt h^U_y)=G^V(x,y)-\sum_{z\in\partial U}H^U(x,z)G^V(z,y)
\\
-\sum_{z\in\partial U}H^U(y,z)G^V(z,x)+\sum_{z,\tilde z\in\partial U}H^U(x,z)H^U(y,\tilde z)G^V(z,\tilde z)\,.
\quad
\end{multline}
Now recall the representation \eqref{E:1.29} which casts $G^V(x,y)$ as $-\fraka(x-y)+\phi(y)$ with~$\phi$ harmonic on~$V$. Plugging this in \eqref{E:3.8}, the fact that
\begin{equation}
\sum_{z\in\partial U}H^U(x,z)\phi(z)=\phi(x),\quad x\in U,
\end{equation}
 shows that all occurrences of~$\phi$ in \eqref{E:3.8} cancel out. As $x\mapsto\fraka(z-x)$ is discrete harmonic on~$U$ for any~$z\in\partial U$, replacing~$G^V(\cdot,\cdot)$ by~$-\fraka(\cdot-\cdot)$ in the last two sums on the right of \eqref{E:3.8} makes these sums cancel each other as well. 

We are thus left with the first two terms on the right of \eqref{E:3.8} in which~$G^V(\cdot,\cdot)$ is now replaced by~$-\fraka(\cdot-\cdot)$. The representation \eqref{E:1.29} then tells us that
\begin{equation}
\Cov(\wt h^U_x,\wt h^U_y) = G^U(x,y),\qquad x,y\in U.
\end{equation}
Since both~$\wt h^U$ and~$h^U$ vanish outside~$U$, we have $\wt h^U\,\,\laweq\,\,h^U$ as desired.
\end{proofsect}

\begin{myexercise}
Supply the missing (e.g., limiting) arguments to prove the Gibbs-Markov decomposition applies even to the situation when~$U$ and~$V$ are allowed to be infinite.
\end{myexercise}

\nopagebreak
\begin{figure}[t]
\vglue-2mm
\centerline{\includegraphics[width=0.5\textwidth]{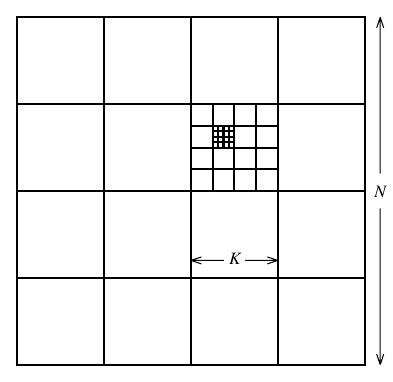}
}
\vglue-1mm
\begin{quote}
\small 
\caption{\small
\label{fig-GM}
A typical setting for the application of the Gibbs-Markov property. The box $V_N=(0,N)^2\cap\Z^2$ is partitioned into~$(N/K)^2$ translates of boxes $V_K:=(0,K)^2\cap\Z^2$ of side~$K$ (assuming~$K$ divides~$N$). This leaves a ``line of sites'' between any two adjacent translates of~$V_K$. The DGFF on~$V_N$ is then partitioned as $h^{V_N}=h^{V_N^\circ}+\varphi^{V_N,V_N^\circ}$ with $h^{V_N^\circ}\independent\varphi^{V_N,V_N^\circ}$, where~$V_N^\circ$ is the union of the shown translates of~$V_K$ and~$\varphi^{V_N,V_N^\circ}$ has the law of the harmonic extension to~$V_N^\circ$ of the values of~$h^{V_N}$ on~$V_N\smallsetminus V_N^\circ$. The translates of~$V_K$ can be further  subdivided to produce a hierarchical description of the DGFF.}
\normalsize
\end{quote}
\end{figure}
\vglue-1mm

A short way to write the Gibbs-Markov decomposition is as
\begin{equation}
\label{E:3.14}
h^V\,\,\laweq\,\, h^U+\varphi^{V,U}\quad\text{where}\quad h^U\independent\varphi^{V,U}.
\end{equation}
with the law of~$h^U$ and~$\varphi^{V,U}$ (often implicitly) as above.

We have seen that the monotonicity~$V\mapsto G^V(x,y)$ allows for control of the variance of the DGFF in general domains by that in more regular ones. One of the important consequences of the Gibbs-Markov property is to give similar comparisons for various probabilities involving a finite number of vertices. The following examples will turn out to be useful: 

\begin{myexercise} 
\label{ex:3.3}
Suppose~$\emptyset\ne U\subseteq V\subsetneq\Z^2$. Prove that for every~$a\in\R$,
\begin{equation}
\label{E:3.4a}
P\bigl(h^U_x\ge a\bigr)\le 2 P\bigl(h^V_x\ge a\bigr),\qquad x\in U.
\end{equation}
Similarly, for any binary relation~$\RR\subseteq \Z^d\times\Z^d$, show that also
\begin{multline}
\label{E:3.12}
\qquad
P\bigl(\exists x,y\in U\colon (x,y)\in \RR,\,h^U_x,h^U_y\ge a\bigr)
\\
\le 4 P\bigl(\exists x,y\in V\colon (x,y)\in \RR,\,h^V_x,h^V_y\ge a\bigr)\,.
\qquad
\end{multline}
\end{myexercise}

Similar ideas lead to:

\begin{myexercise}
\label{ex:3.4}
Prove that for any~$\emptyset\ne U\subseteq V\subsetneq\Z^2$ and any~$a\in\R$,
\begin{equation}
P\bigl(\max_{x\in U}h^U_x>a\bigr)\le 2 P\bigl(\max_{x\in V}h^V_x>a\bigr).
\end{equation}
For finite~$V$ we get
\begin{equation}
E\bigl(\,\max_{x\in U}h^U_x\bigr)\le E\bigl(\,\max_{x\in V}h^V_x\bigr),
\end{equation}
and so~$U\mapsto E(\max_{x\in U}h^U_x)$ is non-decreasing with respect to the set inclusion.
\end{myexercise}

These estimates squeeze the maximum in one domain between that in smaller and larger domains. If crude bounds are sufficient, this permits reduction to simple domains, such as squares.

\begin{figure}[t]
\vglue-2mm
\centerline{\includegraphics[width=0.7\textwidth]{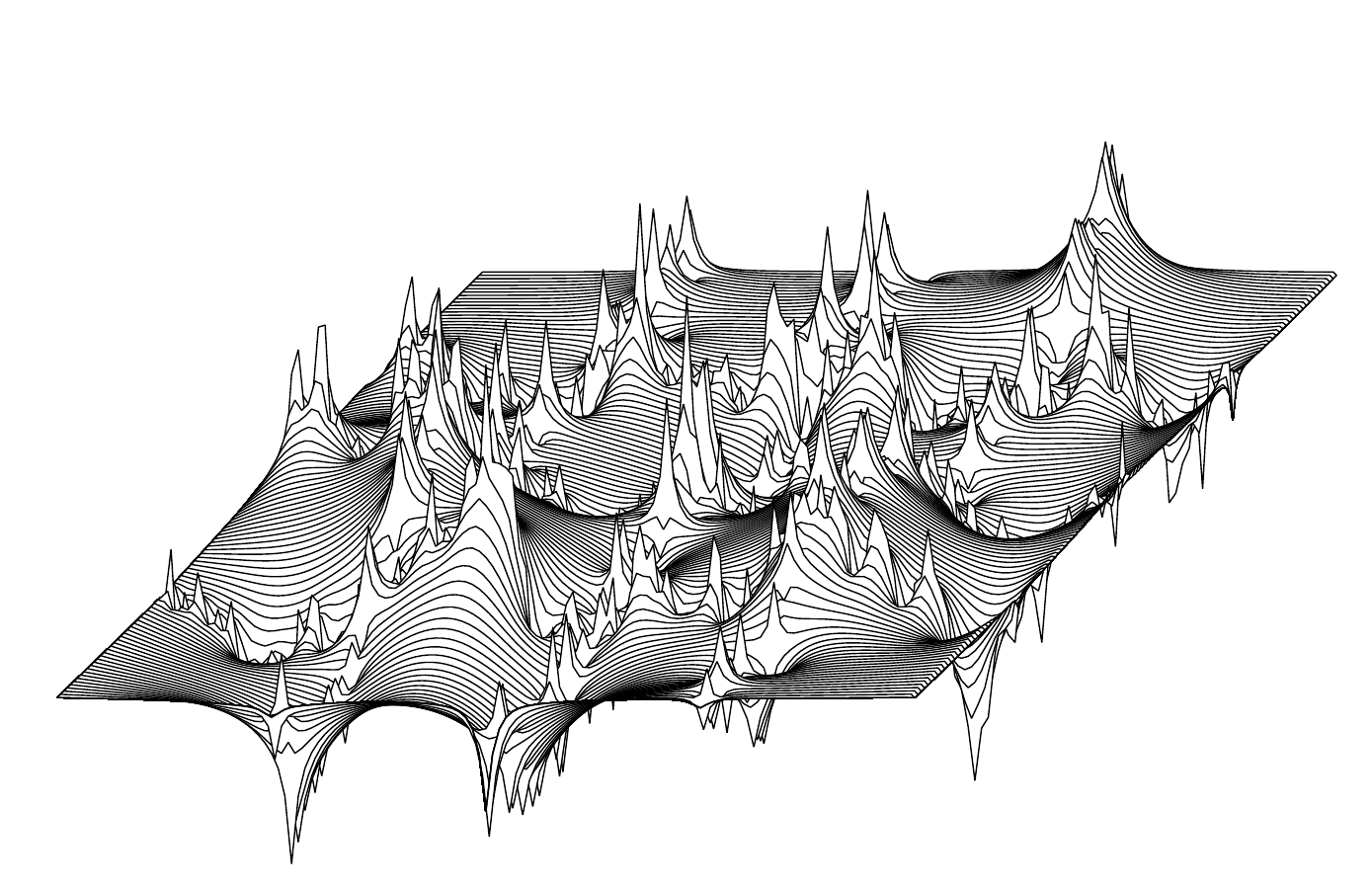}
}
\vglue-1mm
\begin{quote}
\small 
\caption{\small
\label{fig-binding-4x4}
A sample of the binding field~$\varphi^{V_{N},V_{N}^\circ}$ for the (first-level) partition depicted in Fig.~\ref{fig-GM} with $N:=4K$. Here~$V_{N}^\circ$ is the union of~$16$ disjoint translates of~$V_K$. Note that while samples of the field are discrete harmonic inside the individual squares, they become quite rough on the dividing lines of sites.}
\normalsize
\end{quote}
\end{figure}

A typical setting for the application of the Gibbs-Markov property is depicted in~Fig.~\ref{fig-GM}. There each of the small boxes (the translates of~$V_K$) has its ``private'' independent copy of the DGFF. By \eqref{E:3.14}, to get~$h^{V_N}$ these copies are ``bound together'' by an independent Gaussian field $\varphi^{V_N,V_N^\circ}$ that, as far as its law is concerned, is just the harmonic extension of the values of~$h^{V_N}$ on the dividing lines that separate the small boxes from each other. For this reason we sometimes refer to~$\varphi^{V_N,V_N^\circ}$ as the \emph{binding field}. Note that~$\varphi^{V_N,V_N^\circ}$ has discrete-harmonic sample paths on~$V_N^\circ$ yet it becomes quite singular on $V_N\smallsetminus V_N^\circ$; cf Fig.~\ref{fig-binding-4x4}.

Iterations of the partitioning sketched in~Fig.~\ref{fig-GM} lead to a hierarchical description of the DGFF on a square of side~$N:=2^n$ as the sum (along root-to-leaf paths of length~$n$) of a family of tree-indexed binding fields. If these binding fields could be regarded as constant on each of the ``small'' square, this would cast the DGFF as a Branching Random Walk. Unfortunately, the binding fields are not constant on relevant squares so this representation is only approximate. Still, it is extremely useful; see Lecture~\ref{lec-7}. 

The Gibbs-Markov property can be used to bootstrap control from positive probability (in small boxes) to probability close to one (in a larger box, perhaps for a slightly modified event). An example of such a statement is:

\begin{myexercise}
\label{ex:GM-bootstrap}
Prove that ``$\ge$'' holds in \eqref{E:2.10} with probability tending to one as~$N\to\infty$ (still assuming that $\lambda<1/\sqrt 2$).
\end{myexercise}

\noindent
The term \emph{sprinkling technique} is sometimes used for such bootstrap arguments in the literature although we prefer to leave it reserved for parameter manipulations of independent Bernoulli or Poisson random variables.

\nopagebreak
\vglue-1mm

\section{First moment of level-set size}
\noindent
Equipped with the Gibbs-Markov property, we are now ready to begin the proof of the scaling limit of the measures in Theorem~\ref{thm-intermediate}. 
The key point is to estimate, as well as compute the asymptotic of, the first two moments of the size of the level set
\begin{equation}
\Gamma_N^D(b):=\bigl\{x\in D_N\colon h^{D_N}_x\ge a_N+b\bigr\}.
\end{equation}
We begin with the first-moment calculation. Assume that $\lambda\in(0,1)$, an admissible domain~$D\in\mathfrak D$ and an admissible sequence $\{D_N\colon N\ge1\}$ of domains approximating~$D$ are fixed. Our first lemma is then:

\begin{mylemma}[First moment upper bound]
\label{lemma-3.2}
For each~$\delta\in(0,1)$ there is $c\in(0,\infty)$ such that for all $N\ge1$, all $b\in\R$ with $|b|\le\log N$ and all~$a_N$ with $\delta\log N\le a_N\le \delta^{-1}\log N$, and all~$A\subseteq D_N$,
\begin{equation}
\label{E:3.13}
E\bigl|\Gamma_N^D(b)\cap A\bigr|\le cK_N\,\frac{|A|}{N^2}\,\texte^{-\frac{a_N}{g\log N} b}\,.
\end{equation}
\end{mylemma}

\begin{proofsect}{Proof}
Similarly as in the proof of the upper bound in Theorem~\ref{thm-2.1}, the claim will follow by summing over~$x\in A$ the inequality
\begin{equation}
\label{E:3.18uai}
P\bigl(h^{D_N}_x\ge a_N+b\bigr)\le c\frac1{\sqrt{\log N}}\,\,\texte^{-\frac{a_N^2}{2g\log N}}\,\,\texte^{-\frac{a_N}{g\log N} b},
\end{equation}
which (we claim) holds with some~$c>0$
uniformly in~$x\in D_N$ and $b$ with~$|b|\le \log N$. By \eqref{E:3.4a} and translation invariance of the DGFF, it suffices to prove \eqref{E:3.18uai} for~$x:=0$ and~$D_N$ replaced by the box~$\wt D_N$ of side length $4\diam_\infty(D_N)$ centered at the origin. (We will still write~$x$ for the vertex in question though.)

For this setting, Theorem~\ref{thm-1.17} ensures that the variance of~$h_x^{\wt D_N}$ is within a constant~$\tilde c>0$ of~$g\log N$. Hence we get
\begin{equation}
P\bigl(h^{\wt D_N}_x\ge a_N+b\bigr)
\le\frac1{\sqrt{2\pi}}\frac1{\sqrt{g\log N-\tilde c}}\int_b^\infty\texte^{-\frac12\frac{(a_N+s)^2}{g\log N+\tilde c}}\,\textd s.
\end{equation}
Bounding $(a_N+ s)^2\ge a_N^2+2a_N s$ and noting that the assumptions on~$a_N$ (and the inequality $(1+r)^{-1}\ge1-r$ for~$0<r<1$) imply
\begin{equation}
\frac{a_N^2}{g\log N+c}\ge \frac{a_N^2}{g\log N}-\frac c{g^2\delta}\,,
\end{equation}
we get
\begin{equation}
\int_b^\infty\texte^{-\frac12\frac{(a_N+s)^2}{g\log N+\tilde c}}\,\textd s
\le c'\texte^{-\frac{a_N^2}{2g\log N}}\,\,\texte^{-\frac{a_N}{g\log N + \tilde c} b}
\end{equation}
for some constant~$c'>0$. As~$a_N\le\delta^{-1}\log N$ and~$|b|\le\log N$, the constant~$\tilde c$ in the exponent can be dropped at the cost of another multiplicative (constant) term popping up in the front. The claim follows.
\end{proofsect}

The fact that the estimate in Lemma~\ref{lemma-3.2} holds uniformly for all subsets of~$D_N$ will be quite useful for the following reason:

\begin{myexercise}
\label{ex:3.7}
Show that for each~$D\in\mathfrak D$, each~$b\in\R$ and each~$\epsilon>0$ there is~$\delta>0$ such that for all~$N$ sufficiently large,
\begin{equation}
E\bigl|\{x\in \Gamma_N^D(b)\colon \dist_\infty(x,D_N^\cc)<\delta N\}\bigr|\le\epsilon K_N\,.
\end{equation}
\end{myexercise}

\noindent
We remark that no properties of~$\partial D$ other than those stated in Definition~\ref{def-admis-domain} should  be assumed in the solution. Next let us note the following fact:

\begin{myexercise}
\label{ex:3.8}
A sequence $\{\mu_n\colon n\ge1\}$ of random Borel measures on a topological space~$\scrX$ is tight with respect to the vague topology if and only if the sequence of random variables $\{\mu_N(K)\colon n\ge1\}$ is tight for every compact~$K\subseteq\scrX$.
\end{myexercise}

\noindent
Lemma~\ref{lemma-3.2} then gives:

\begin{mycorollary}[Tightness]
The family $\{\eta^D_N\colon N\ge1\}$, regarded as measures on~$\overline D\times(\R\cup\{+\infty\})$, is a tight sequence in the topology of vague convergence.
\end{mycorollary}

\begin{proofsect}{Proof}
Every compact set in $\overline D\times(\R\cup\{+\infty\})$ is contained in $K_b:=\overline D\times[b,\infty]$ for some~$b\in\R$. The definition of~$\eta^D_N$ shows
\begin{equation}
\eta^D_N(K_b)=\frac1{K_N}\bigl|\Gamma^D_N(b)\bigr|.
\end{equation}
Lemma~\ref{lemma-3.2} shows these have uniformly bounded expectations and so are tight as ordinary random variables. The claim follows by Exercise~\ref{ex:3.8}.
\end{proofsect}

In probability, tightness is usually associated with ``mass not escaping to infinity'' or ``the total mass being conserved.'' However, for convergence of random unnormalized measures in the vague topology, tightness does not prevent convergence to zero measure. In order to rule that out, we will need: 

\begin{mylemma}[First moment asymptotic]
\label{lemma-3.3}
Assume that~$a_N$ obeys \eqref{E:2.28} and let~$c_0$ be as in \eqref{E:1.32}. Then for all~$b\in\R$ and all open $A\subseteq D$,
\begin{equation}
\label{E:3.6}
E\bigl|\{x\in\Gamma_N^D(b)\colon x/N\in A\}\bigr| = 
\frac{\texte^{2c_0\lambda^2/g}}{\lambda\sqrt{8\pi}}\texte^{-\alpha\lambda b}\Bigl[\,o(1)+\int_A\textd x\, r_D(x)^{2\lambda^2}\Bigr]\,K_N\,,
\end{equation}
where $o(1)\to0$ as~$N\to\infty$ uniformly on compact sets of~$b$.
\end{mylemma}

\begin{proofsect}{Proof}
Thanks to Exercise~\ref{ex:3.7} we may remove a small neighborhood of~$\partial D$ from~$A$ and thus assume that~$\dist_\infty(A,D^\cc)>0$. We will proceed by extracting an asymptotic expression for $P(h^{D_N}_x\ge a_N+b$) with~$x$ such that~$x/N\in A$.
For such~$x$, Theorem~\ref{thm-1.17}  gives
\begin{equation}
G^{D_N}(x,x) = g\log N+\theta_N(x),
\end{equation}
where
\begin{equation}
\label{G-asymp}
\theta_N(x) = c_0+g\log r_D(x/N)+o(1)\,,
\end{equation}
with $o(1)\to0$ as $N\to\infty$ uniformly in~$x\in D_N$ with~$x/N\in A$. Using this in the formula for the probability density of~$h^{D_N}_x$ yields
\begin{equation}
P\bigl(h^{D_N}_x\ge a_N+b\bigr)=\frac1{\sqrt{2\pi}}\frac1{\sqrt{g\log N+\theta_N(x)}}\int_b^\infty\texte^{-\frac12\frac{(a_N+s)^2}{g\log N+\theta_N(x)}}\,\textd s\,.
\end{equation}
The first occurrence of~$\theta_N(x)$ does not affect the overall asymptotic as this quantity is bounded uniformly for all~$x$ under consideration. Expanding the square $(a_N+s)^2=a_N^2+2a_Ns+s^2$
 and noting that (by decomposing the integration domain into~$s\ll{\log N}$ and its complement) the~$s^2$ term has negligible effect on the overall asymptotic of the integral, we find out
\begin{equation}
\int_b^\infty\texte^{-\frac12\frac{(a_N+s)^2}{g\log N+\theta_N(x)}}\,\textd s
=\bigl(1+o(1)\bigr)(\alpha\lambda)^{-1}\texte^{-\frac12\frac{a_N^2}{g\log N+\theta_N(x)}}\,\texte^{-\alpha\lambda b+o(1)}\,.
\end{equation}
We now use Taylor's Theorem (and the asymptotic of~$a_N$) to get
\begin{equation}
\frac{a_N^2}{g\log N+\theta_N(x)}
=\frac{a_N^2}{g\log N}-\frac{4\lambda^2}g\theta_N(x)+o(1)
\end{equation}
with~$o(1)\to0$ again uniformly in all~$x$ under consideration. This yields
\begin{equation}
P\bigl(h^{D_N}_x\ge a_N+b\bigr) = \bigl(1+o(1)\bigr)\frac{\texte^{2c_0\lambda^2/g}}{\lambda\sqrt{8\pi}}\,\texte^{-\alpha\lambda b}\,r_D(x/N)^{2\lambda^2}\,\frac{K_N}{N^2}.
\end{equation}
The result follows by summing this probability over~$x$ with~$x/N\in A$ and using the continuity of~$r_D$ to convert the resulting Riemann sum into an integral.
\end{proofsect}

\section{Second moment estimate}
\noindent
Our next task is to perform a rather tedious estimate on the second moment of the size of~$\Gamma_N^D(b)$. It is here where we need to limit the range of possible~$\lambda$. 

\begin{mylemma}[Second moment bound]
\label{lemma-3.3b}
Suppose $\lambda\in(0,1/\sqrt2)$. For each $b_0>0$ and each~$D\in\mathfrak D$ there is $c_1\in(0,\infty)$ such that for each $b\in[-b_0,b_0]$ and each~$N\ge1$,
\begin{equation}
\label{e:2.3}
E \bigl(|\Gamma_N^D (b)|^2 \bigr)
\le c_1 K_N^2\,.
\end{equation}
\end{mylemma}

\begin{proofsect}{Proof}
Assume $b:=0$ for simplicity (or absorb~$b$ into~$a_N$). Writing
\begin{equation}
\label{E:3.9a}
E \bigl(|\Gamma_N^D (0)|^2 \bigr) = \sum_{x,y\in D_N}P\bigl(h^{ D_N}_x\ge a_N,\,h^{ D_N}_y\ge a_N\bigr).
\end{equation}
we need to derive a good estimate for the probability on the right-hand side.  In order to ensure uniformity, let~$\wt D_N$ be a neighborhood of~$D_N$ of diameter twice the diameter of~$D_N$. Exercise~\ref{ex:3.3} then shows
\begin{equation}
\label{E:3.9ua}
P\bigl(h^{ D_N}_x\ge a_N,\,h^{ D_N}_y\ge a_N\bigr)
\le 4P\bigl(h^{\wt D_N}_x\ge a_N,\,h^{\wt D_N}_y\ge a_N\bigr)\,.
\end{equation}
We will now estimate the probability on the right by conditioning on $h^{\wt D_N}_x$.

First note that the Gibbs-Markov property yields a pointwise decomposition
\begin{equation}
\label{E:3.10}
h^{\wt D_N}_y = \frakg_x(y)h^{\wt D_N}_x+\hat h^{\wt D_N\smallsetminus\{x\}}_y,
\end{equation} 
where 
\settowidth{\leftmargini}{(1111)}
\begin{enumerate}
\item[(1)] $h^{\wt D_N}_x$ and $\hat h^{\wt D_N\smallsetminus\{x\}}$ are independent,
\item[(2)] $\hat h^{\wt D_N\smallsetminus\{x\}}$ has the law of the DGFF in~$\wt D_N\smallsetminus\{x\}$, and
\item[(3)] $\frakg_x$ is harmonic in~$\wt D_N\smallsetminus\{x\}$, vanishes outside~$\wt D_N$ and obeys $\frakg_x(x)=1$.
\end{enumerate}
Using \eqref{E:3.10}, the above probability is recast as
\begin{multline}
\label{E:3.11}
P\bigl(h^{\wt D_N}_x\ge a_N,\,h^{\wt D_N}_y\ge a_N\bigr)
\\
=\int_0^\infty P\Bigl(\,\hat h^{\wt D_N\smallsetminus\{x\}}_y\ge a_N(1-\frakg_x(y))-s\frakg_x(y)\Bigr)P\bigl(h^{\wt D_N}_x-a_N\in\textd s\bigr).
\end{multline}
Given~$\delta>0$ we can always bound the right-hand side by $P(h^{\wt D_N}_x\ge a_N)$ when $|x-y|\le\delta\sqrt{K_N}$. This permits us to assume that $|x-y|>\delta\sqrt{K_N}$ from now on. The $s\ge a_N$ portion of the integral is similarly bounded by $P(h^{\wt D_N}_x\ge 2a_N)$, so we will henceforth focus on $s\in[0,a_N]$.

Since~$x,y$ lie ``deep'' inside~$\wt D_N$ and $|x-y|>\delta\sqrt{K_N}=N^{1-\lambda^2+o(1)}$, we have
\begin{equation}
\begin{aligned}
\label{E:3.13ua}
\frakg_x(y)=\frac{G^{\wt D_N}(x,y)}{G^{\wt D_N}(x,x)}
&\le \frac{\log\frac{N}{|x-y|}+c}{\log N-c}
\\
&\le1-(1-\lambda^2)+o(1)=\lambda^2+o(1),
\end{aligned}
\end{equation}
where $o(1)\to0$ uniformly in~$x,y\in D_N$. For~$s\in[0,a_N]$, $\lambda<1/\sqrt 2$ implies the existence of an~$\epsilon>0$ such that, for~$N$ is large enough,
\begin{equation}
\epsilon a_N\le a_N\bigl(1-\frakg_x(y)\bigr)-s\frakg_x(y)\le a_N,\quad x,y\in D_N.
\end{equation}
The Gaussian bound $P(X\ge t)\le\sigma t^{-1}\texte^{-\frac{t^2}{2\sigma^2}}$ for~$X=\NN(0,\sigma^2)$ and any~$t>0$ along with $G^{\wt D_N\smallsetminus\{x\}}(y,y)\le g\log N+c$ uniformly in~$y\in D_N$ show
\begin{multline}
P\Bigl(\,\hat h^{\wt D_N\smallsetminus\{x\}}_y\ge a_N\bigl(1-\frakg_x(y)\bigr)-s\frakg_x(y)\Bigr)
\\
\le \frac{\sqrt{G(y,y)}}{\epsilon a_N}\,\texte^{-\frac{[a_N(1-\frakg_x(y))-s\frakg_x(y)]^2}{2G(y,y)}}
\le c\frac{K_N}{N^2}\,\texte^{\frakg_x(y)\frac{a_N^2}{g\log N}+\frac{a_N}{g\log N}\frakg_x(y)s}\,,
\end{multline}
where $G(y,y)$ abbreviates $G^{\wt D_N\smallsetminus\{x\}}(y,y)$. The first inequality in \eqref{E:3.13ua} gives
\begin{equation}
\label{E:3.18e}
\texte^{\frakg_x(y)\frac{a_N^2}{g\log N}}
\le c\biggl(\frac N{|x-y|}\biggr)^{4\lambda^2+o(1)}
\end{equation}
with~$o(1)\to0$ uniformly in~$x,y\in D_N$ with~$|x-y|>\delta\sqrt{K_N}$.

The uniform upper bound on $G^{\wt D_N}(x,x)$ allows us to dominate the law of~$h^{\wt D_N}_x-a_N$ on~$[0,\infty)$ by
\begin{equation}
\label{E:3.17e}
P\bigl(h^{\wt D_N}_x-a_N\in\textd s\bigr)
\le c\frac{K_N}{N^2}\,\texte^{-\frac{a_N}{g\log N} s}\textd s.
\end{equation}
As $1-\frakg_x(y)=1-\lambda^2+o(1)$, the $s\in[0,a_N]$ part of the integral in \eqref{E:3.11} is then readily performed to yield
\begin{equation}
\label{E:3.16e}
P\bigl(h^{\wt D_N}_x\ge a_N,\,h^{\wt D_N}_y\ge a_N\bigr)
\le P\bigl(h^{\wt D_N}_x\ge 2a_N\bigr)
+c\Bigl(\frac {K_N}{N^2}\Bigr)^2\biggl(\frac N{|x-y|}\biggr)^{4\lambda^2+o(1)}
\end{equation}
uniformly in~$x,y\in D_N$ withÊ$|x-y|>\delta\sqrt{K_N}$.

In order to finish the proof, we now use \eqref{E:3.9ua} to write
\begin{equation}
\label{E:3.17}
\begin{aligned}
E \bigl(|\Gamma_N^D (0)|^2 \bigr)
&\le\sum_{\begin{subarray}{c}
x,y\in D_N\\|x-y|\le\delta\sqrt{K_N}
\end{subarray}}
P\bigl(h^{ D_N}_x\ge a_N\bigr)
\\
&\qquad\qquad\qquad+\sum_{\begin{subarray}{c}
x,y\in D_N\\|x-y|>\delta\sqrt{K_N}
\end{subarray}}
4P\bigl(h^{\wt D_N}_x\ge a_N,\,h^{\wt D_N}_y\ge a_N\bigr).
\end{aligned}
\end{equation}
Summing over~$y$ and invoking Lemma~\ref{lemma-3.2} bounds the first term by a factor of order $(\delta K_N)^2$. The contribution of the first term on the right of \eqref{E:3.16e} to the second sum is bounded via Lemma~\ref{lemma-3.2} as well:
\begin{equation}
\begin{aligned}
P\bigl(h^{D_N}_x\ge 2a_N\bigr)
&\le \frac{c}{\sqrt{\log N}}\texte^{-2\frac{a_N^2}{g\log N}}
\\
&= c\Bigl(\frac{K_N}{N^2}\Bigr)^2 \texte^{-\frac{a_N^2}{g\log N}}\sqrt{\log N}\le c \delta\Bigl(\frac{K_N}{N^2}\Bigr)^2.
\end{aligned}
\end{equation}
Plugging in also the second term on the right of \eqref{E:3.16e}, we thus get
\begin{equation}
\label{E:3.22e}
E \bigl(|\Gamma_N^D (0)|^2 \bigr)
\le 2c\delta(K_N)^2+c\Bigl(\frac {K_N}{N^2}\Bigr)^2\sum_{\begin{subarray}{c}
x,y\in D_N\\|x-y|>\delta\sqrt{K_N}
\end{subarray}}
\biggl(\frac N{|x-y|}\biggr)^{4\lambda^2+o(1)}\,.
\end{equation}
Dominating the sum by $cN^4\int_{D\times D}|x-y|^{-4\lambda^2+o(1)}\textd x\textd y$, with the integral convergent due to $4\lambda^2<2$, we find that also the second term on the right is of order~$(K_N)^2$. 
\end{proofsect}

As a corollary we now get:

\begin{mycorollary}[Subsequential limits are non-trivial]
Let~$\lambda\in(0,1/\sqrt2)$. Then every subsequential limit~$\eta^D$ of $\{\eta^D_N\colon N\ge1\}$ obeys
\begin{equation}
P\bigl(\eta^D(A\times[b,b'])>0\bigr)>0
\end{equation}
for any open and non-empty~$A\subseteq D$ and every~$b<b'$.
\end{mycorollary}

\begin{proofsect}{Proof}
Abbreviate~$X_N:=\eta^D_N(A\times[b,b'])$. Then Lemma~\ref{lemma-3.3} implies
\begin{equation}
E(X_N)\,\underset{N\to\infty}\longrightarrow\,\hat c\Bigl[\int_A\textd x\, r_D(x)^{2\lambda^2}\Bigr]\bigl(\texte^{-\lambda\alpha b}-\texte^{-\lambda\alpha b'}\bigr)\,,
\end{equation}
where $\hat c:=\texte^{2c_0\lambda^2/g}/(\lambda\sqrt{8\pi})$. This is positive and finite for any~$A$ and~$b,b'$ as above. On the other hand, Lemma~\ref{lemma-3.3b} shows that~$\sup_{N\ge1}E(X_N^2)<\infty$. The second-moment estimate (Exercise~\ref{ex:2.5}) then yields the claim.
\end{proofsect}

\section{Second-moment asymptotic and factorization}
\noindent
At this point we know that the subsequential limits exist and are non-trivial (with positive probability). The final goal of this lecture is to prove:

\begin{myproposition}[Factorization]
\label{prop-factorization}
Suppose~$\lambda\in(0,1/\sqrt2)$. Then every subsequential limit~$\eta^D$ of $\{\eta^D_N\colon N\ge1\}$ takes the form
\begin{equation}
\label{E:3.46}
\eta^D(\textd x\,\textd h) = Z^D_\lambda(\textd x)\otimes\texte^{-\alpha\lambda h}\textd h,
\end{equation}
where~$Z^D_\lambda$ is a random, a.s.-finite Borel measure on~$D$ with $P(Z^D_\lambda(D)>0)>0$.
\end{myproposition}

The proof relies on yet another (and this time quite lengthy) second-moment calculation. The result of this calculation is the content of:

\begin{mylemma}
\label{lemma-aux}
For any $\lambda\in(0,1/\sqrt2)$, any open $A\subseteq D$, any $b\in\R$, and
\begin{equation}
A_N:=\{x\in \Z^2\colon x/N\in A\}
\end{equation}
we have
\begin{equation}
\lim_{N\to\infty}\,
\frac1{K_N}E\Bigl|\bigl|\Gamma_N^D(0)\cap A_N\bigr|-\texte^{\alpha\lambda b}\bigl|\Gamma_N^D(b)\cap A_N\bigr|\Bigr|=0.
\end{equation}
\end{mylemma}

\begin{proofsect}{Proof (modulo a computation)}
By Lemma~\ref{lemma-3.2} we may assume~$\dist_\infty(A,D^\cc)>\epsilon$ for some~$\epsilon>0$. We will  prove
\begin{equation}
\label{E:3.51uai}
\lim_{N\to\infty}\,
\frac1{K_N^2}E\biggl(\Bigl|\bigl|\Gamma_N^D(0)\cap A_N\bigr|-\texte^{\alpha\lambda b}\bigl|\Gamma_N^D(b)\cap A_N\bigr|\Bigr|^2\biggr)=0
\end{equation}
which implies the claim via the Cauchy-Schwarz inequality.
Plugging
\begin{equation}
\bigl|\Gamma_N^D(\cdot)\cap A_N\bigr|=\sum_{x\in A_N}1_{\{h^{D_N}_x\ge a_N+\cdot\}}
\end{equation}
into \eqref{E:3.51uai} 
we get a sum of pairs of (signed) products of the various combinations of these indicators. The argument in the proof of Lemma~\ref{lemma-3.3b} allows us to estimate the pairs where $|x-y|\le\delta N$ by a quantity that vanishes as~$N\to\infty$ and~$\delta\downarrow0$. It will thus suffice to show
\begin{equation}
\label{E:3.27}
\begin{aligned}
\max_{\begin{subarray}{c}
x,y\in A_N\\|x-y|>\delta N
\end{subarray}}
\biggl(&P\bigl(h^{D_N}_x\ge a_N,\,h^{D_N}_y\ge a_N\bigr)
\\*[-4mm]
&-\texte^{\alpha\lambda b}P\bigl(h^{D_N}_x\ge a_N+b,\,h^{D_N}_y\ge a_N\bigr)
\\*[1mm]
&\quad-\texte^{\alpha\lambda b}P\bigl(h^{D_N}_x\ge a_N,\,h^{D_N}_y\ge a_N+b\bigr)
\\
&\qquad
+
\texte^{2\alpha\lambda b}P\bigl(h^{D_N}_x\ge a_N+b,\,h^{D_N}_y\ge a_N+b\bigr)\biggr)
= o\Bigl(\frac{K_N^2}{N^4}\Bigr)
\end{aligned}
\end{equation}
as~$N\to\infty$. A computation refining the argument in the proof of Lemma~\ref{lemma-3.3b} to take into account the precise asymptotic of the Green function (this is where we get aided by the fact that~$|x-y|>\delta N$ and~$\dist_\infty(A,D^\cc)>\epsilon$) now shows that, for any $b_1,b_2\in\{0,b\}$,
\begin{multline}
\label{E:3.32ua}
\qquad
P\bigl(h^{D_N}_x\ge a_N+b_1,\,h^{D_N}_y\ge a_N+b_2\bigr)
\\=\bigl(\texte^{-\alpha\lambda (b_1+b_2)}+o(1)\bigr)P\bigl(h^{D_N}_x\ge a_N,\,h^{D_N}_y\ge a_N\bigr)
\qquad
\end{multline}
with~$o(1)\to0$ as~$N\to\infty$ uniformly in~$x,y\in A_N$.
This then implies \eqref{E:3.27} and thus the whole claim.
\end{proofsect}

\begin{myexercise}
Supply a detailed proof of \eqref{E:3.32ua}. (Consult \cite{BL2} if lost.)
\end{myexercise}

From Lemma~\ref{lemma-aux} we get:

\begin{mycorollary}
\label{cor-3.12}
Suppose~$\lambda\in(0,1/\sqrt2)$. Then any subsequential limit~$\eta^D$ of the processes~$\{\eta^D_N\colon N\ge1\}$ obeys the following: For any open~$A\subseteq D$ and any~$b\in\R$,
\begin{equation}
\eta^D\bigl(A\times[b,\infty)\bigr)=\texte^{-\alpha\lambda b}\eta^D\bigl(A\times[0,\infty)\bigr),\quad\text{\rm a.s.}
\end{equation}
\end{mycorollary}

\begin{proofsect}{Proof}
In the notation of Lemma~\ref{lemma-aux}, 
\begin{equation}
\eta^D_N\bigl(A\times[b,\infty)\bigr) = \frac1{K_N}\bigl|\Gamma_N^D(b)\cap A_N\bigr|.
\end{equation}
Taking a \emph{joint} distributional limit of $\eta^D_N(A\times[b,\infty))$ and $\eta^D_N(A\times[0,\infty))$ 
along the given subsequence, Lemma~\ref{lemma-aux} along with Fatou's lemma show 
\begin{equation}
E\Bigl|\eta^D\bigl(A\times[0,\infty)\bigr)-\texte^{\alpha\lambda b}\eta^D\bigl(A\times[b,\infty)\bigr)\Bigr|=0.
\end{equation}
(This requires a routine approximation of indicators of these events by continuous functions as in Exercise~\ref{ex2.10}.) The claim follows.
\end{proofsect}

We now give:

\begin{proofsect}{Proof of Proposition~\ref{prop-factorization}}
For each Borel $A\subseteq D$ define
\begin{equation}
Z^D_\lambda(A):=(\alpha\lambda)\eta^D\bigl(A\times[0,\infty)\bigr)
\end{equation}
Then~$Z^D_\lambda$ is an a.s.-finite (random) Borel measure on~$D$.
Letting~$\AA$ be the set of all half-open dyadic boxes entirely contained in~$D$, Corollary~\ref{cor-3.12} and a simple limiting argument show that, for any~$A\in\AA$ and any $b\in\Q$,
\begin{equation}
\label{E:3.55}
\eta^D\bigl(A\times[b,\infty)\bigr) = (\alpha\lambda)^{-1}Z^D_\lambda(A)\texte^{-\alpha\lambda b},\quad \text{a.s.}
\end{equation}
Since~$\AA\times\Q$ is countable, the null set in \eqref{E:3.55} can be chosen so that the equality in \eqref{E:3.55} holds for all~$A\in\AA$ and all~$b\in\Q$ simultaneously, a.s. Note that the sets $\{A\times[b,\infty)\colon A\in\AA, b\in\Q\}$ constitute a $\pi$-system that generates all Borel sets in~$D\times\R$. In light of
\begin{equation}
(\alpha\lambda)^{-1}Z^D_\lambda(A)\texte^{-\alpha\lambda b} = \int_{A\times[b,\infty)}Z^D_\lambda(\textd x)\otimes\texte^{-\alpha\lambda h}\textd h
\end{equation}
the claim follows from Dynkin's $\pi$-$\lambda$-theorem.
\end{proofsect}

We also record an important observation:

\begin{mycorollary}
\label{cor-3.17}
Assume~$\lambda\in(0,1/\sqrt 2)$ and denote
\begin{equation}
\label{E:3.58}
\hat c:=\frac{\texte^{2c_0\lambda^2/g}}{\lambda\sqrt{8\pi}}
\end{equation}
 for~$c_0$ as in~\eqref{E:1.32}. Then~$Z^D_\lambda$ from \eqref{E:3.46} obeys
\begin{equation}
\label{E:3.59}
E\bigl[Z^D_\lambda(A)\bigr] = \hat c\int_A\textd x\, r_D(x)^{2\lambda^2}
\end{equation}
for each Borel~$A\subseteq D$. Moreover, there is~$c\in(0,\infty)$ such that for any open square~$S\subset\C$,
\begin{equation}
\label{E:3.60}
E\bigl[Z^S_\lambda(S)^2\bigr]\le c \diam(S)^{4+4\lambda^2}.
\end{equation}
\end{mycorollary}

\begin{proofsect}{Proof (sketch)}
Thanks to the uniform square integrability proved in Lemma~\ref{lemma-3.3b}, the convergence in probability is accompanied by convergence of the first moments. Then \eqref{E:3.59} follows from Lemma~\ref{lemma-3.3}. To get also \eqref{E:3.60} we need a uniform version of the bound in Lemma~\ref{lemma-3.3b}. We will not perform the requisite calculation, just note that for a~$c'\in(0,\infty)$ the following holds for all~$D\in\mathfrak D$,
\begin{equation}
\label{e:2.3b}
\limsup_{N\to\infty}\frac1{K_N^2}\,E \bigl(|\Gamma_N^D (0)|^2\bigr)\le
c'\int_{D\times D}\textd x\otimes\textd y\,\biggl(\frac{[\diam D]^2}{|x-y|}\biggr)^{4\lambda^2}\,,
\end{equation}
where $\diam D$ is the diameter of~$D$ in the Euclidean norm. We leave further details of the proof to the reader.
\end{proofsect}

This closes the first part of the proof of Theorem~\ref{thm-intermediate} which showed that every subsequential limit of the measures of interest factors into the desired product form. The proof continues in the next lecture.


\chapter{Intermediate level sets: nailing the limit}
\label{lec-4}
\noindent
The goal of this lecture is to finish the proof of Theorem~\ref{thm-intermediate} and the results that follow thereafter. This amounts to proving a list of properties that the~$Z^D_\lambda$-measures (still tied to a specific subsequence) satisfy and showing that these properties characterize the law of the $Z^D_\lambda$-measures uniquely. As part of the proof, we obtain a conformal transformation rule for~$Z^D_\lambda$ and a representation thereof as a Liouville Quantum Gravity measure. All proofs remain restricted to $\lambda<1/\sqrt2$; the final section comments on necessary changes in the complementary regime of~$\lambda$'s.

\section{Gibbs-Markov property in the scaling limit}
\noindent
We have shown so far that every subsequential limit of the family of point measures $\{\eta_N^D\colon N\ge1\}$ takes the form
\begin{equation}
Z^D_\lambda(\textd x)\otimes\texte^{-\alpha\lambda h}\textd h
\end{equation}
 for \emph{some} random Borel measure~$Z^D_\lambda$ on~$D$. Our next goal is to identify properties of these measures that will ultimately nail their law uniquely. The most important of these is the behavior under restriction to a subdomain which arises from the Gibbs-Markov decomposition of the DGFF (which defines the concept of the binding field used freely below). However, as the~$Z^D_\lambda$-measure appears only in the scaling limit, we have to first describe the scaling limit of the Gibbs-Markov decomposition itself.

The main observation to be made below is that, although the DGFF has no pointwise scaling limit, the binding field does. This is facilitated (and basically implied) by the fact that the binding field has discrete-harmonic sample paths. To define the relevant objects, let~$\wt D,D\in\mathfrak D$ be two domains satisfying~$\wt D\subseteq D$. For each~$x,y\in\wt D$, set
\begin{multline}
\label{E:4.1}
\quad
C^{D,\wt D}(x,y) 
:= g\int_{\partial D}\Pi^D(x,\textd z)\log|y-z|
\\-g\int_{\partial \wt D}\Pi^{\wt D}(x,\textd z)\log|y-z|\,,
\quad
\end{multline}
where~$\Pi^D$ is the harmonic measure from \eqref{E:1.26ua}.
Given any admissible approximations $\{D_N\colon N\ge1\}$ and~$\{\wt D_N\colon N\ge1\}$ of domains~$D$ and~$\wt D$, respectively, of which we also assume that~$\wt D_N\subseteq D_N$ for each~$N\ge1$, we now observe:

\begin{mylemma}[Convergence of covariances]
\label{lemma-4.1}
Locally uniformly in~$x,y\in\wt D$,
\begin{equation}
\label{E:4.2}
G^{D_N}\bigl(\lfloor xN\rfloor,\lfloor yN\rfloor\bigr)
-G^{\wt D_N}\bigl(\lfloor xN\rfloor,\lfloor yN\rfloor\bigr)\,\underset{N\to\infty}\longrightarrow\, C^{D,\wt D}(x,y).
\end{equation}
\end{mylemma}

We leave it to the reader to solve:

\begin{myexercise}
Prove Lemma~\ref{lemma-4.1} while noting that this includes uniform convergence on the diagonal $x=y$. Hint: Use the representation in Lemma~\ref{lemma-1.19}.
\end{myexercise}

From here we get:

\begin{mylemma}[Limit binding field]
For any~$D$ and~$\wt D$ as above, $x,y\mapsto C^{D,\wt D}(x,y)$ is a symmetric, positive semi-definite kernel on~$\wt D\times\wt D$. In particular, there is a Gaussian process $x\mapsto\Phi^{D,\wt D}(x)$ on~$\wt D$ with zero mean and covariance
\begin{equation}
\Cov\bigl(\Phi^{D,\wt D}(x),\Phi^{D,\wt D}(y)\bigr)  =C^{D,\wt D}(x,y),\qquad x,y\in\wt D.
\end{equation}
\end{mylemma}

\begin{proofsect}{Proof}
Let~$U\subset V$ be non-empty and finite. The Gibbs-Markov decomposition implies
\begin{equation}
\label{E:4.4ua}
\Cov(\varphi^{V,U}_x,\varphi^{V,U}_y) = G^V(x,y)-G^U(x,y).
\end{equation}
Hence, $x,y\mapsto G^V(x,y)-G^U(x,y)$ is symmetric and positive semi-definite on $U\times U$. In light of~\eqref{E:4.2}, this extends to~$C^{D,\wt D}$ on~$\wt D\times\wt D$ by a limiting argument. Standard arguments then imply the existence of the Gaussian process~$\Phi^{D,\wt D}$. 
\end{proofsect}

\nopagebreak
\begin{figure}[t]
\vglue-2mm
\centerline{\includegraphics[width=0.7\textwidth]{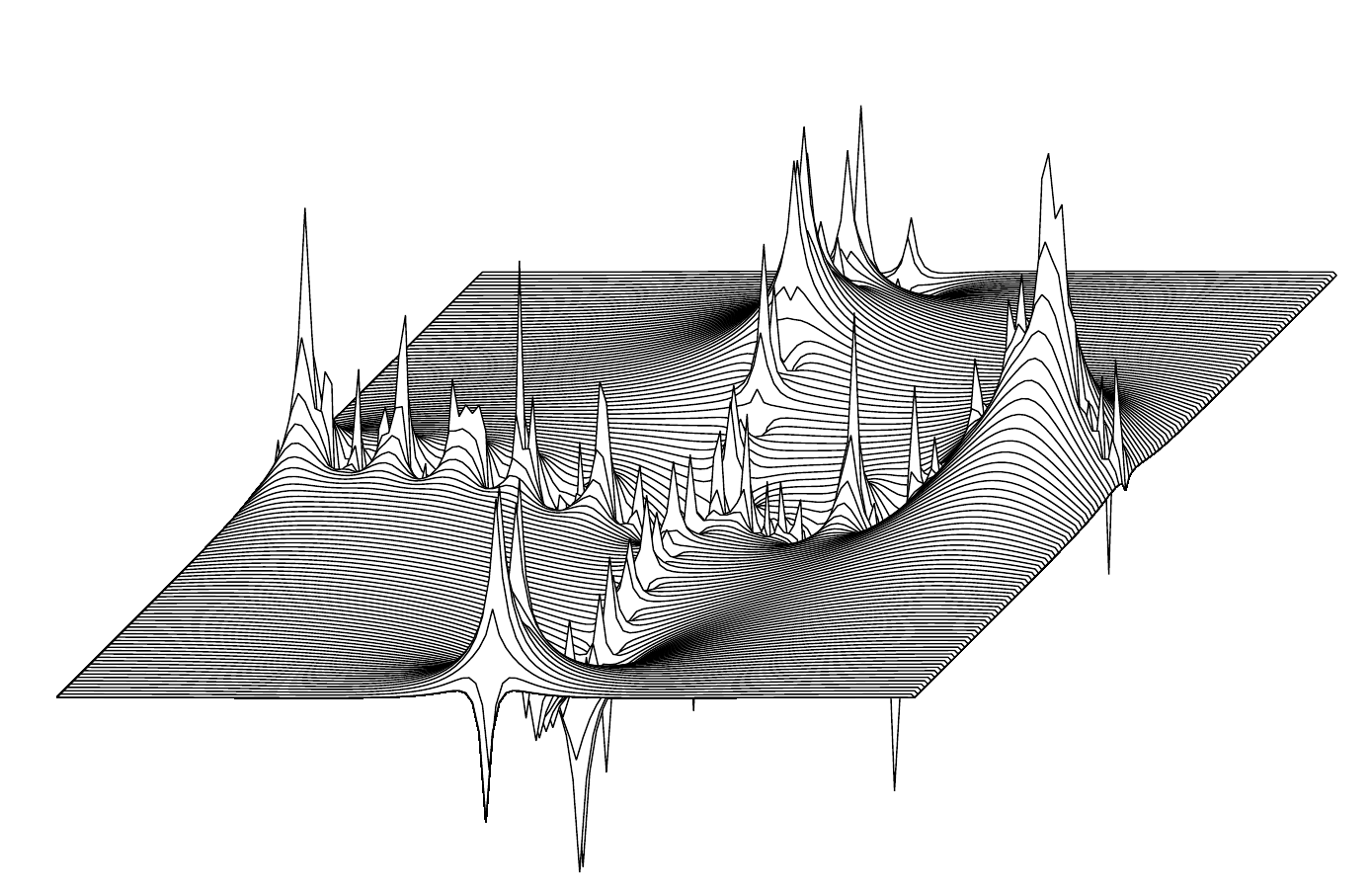}
}
\vglue-1mm
\begin{quote}
\small 
\caption{
\label{fig-binding}
\small
A sample of~$\varphi_N^{D,\wt D}$ where~$D:=(-1,1)^2$ and~$\wt D$ obtained from~$D$ by removing points on the coordinate axes.}
\normalsize
\end{quote}
\end{figure}
\vglue-1mm

In light of Exercise~\ref{ex:1.3} and \eqref{E:4.4ua} we see that $C^{D,\wt D}(x,y)\ge0$ for all~$x,y\in\wt D$. When $\wt D\subsetneq D$, we will even have $C^{D,\wt D}(x,x)>0$ for some~$x\in\wt D$. We will call $\Phi^{D,\wt D}$ the \emph{continuum binding field}. To justify this name, observe:

\begin{mylemma}[Coupling of binding fields]
\label{lemma-4.4}
$\Phi^{D,\wt D}$ has a version with continuous sample paths on~$\wt D$. Moreover, for each $\delta>0$ and each $N\ge1$ there is a coupling of $\varphi^{D_N,\wt D_N}$ and~$\Phi^{D,\wt D}$ such that
\begin{equation}
\sup_{\begin{subarray}{c}
x\in \wt D\\\dist(x, \partial\wt D)>\delta
\end{subarray}}
\bigl|\,\Phi^{D,\wt D}(x)-\varphi^{D_N,\wt D_N}_{\lfloor xN\rfloor}\bigr|\,\underset{N\to\infty}\longrightarrow\,0,\qquad\text{\rm in probability}.
\end{equation}
\end{mylemma}

\begin{proofsect}{Proof (assuming regularity of the fields)}
Abbreviate~$\varphi^{D,\wt D}_N(x):=\varphi^{D_N,\wt D_N}_{\lfloor xN\rfloor}$.
The convergence of the covariances from Lemma~\ref{lemma-4.1} implies $\varphi_N^{D,\wt D}(x)\to \Phi^{D,\wt D}(x)$ in law for each~$x\in\wt D$, so the point is to extend this to the convergence of these fields as random functions. Fix $\delta>0$ and denote
\begin{equation}
\wt D^\delta:=\bigl\{x\in\C\colon\dist(x,\partial\wt D)>\delta\bigr\}\,.
\end{equation}
Fix~$r>0$ small and let $x_1,\dots,x_k$ be an $r$-net in $\wt D^\delta$. As convergence of the covariances implies convergence in law, and convergence in law on~$\R^k$ can be turned into convergence in probability under a suitable coupling measure, for each~$N\ge1$ there is a coupling of $\varphi^{D_N,\wt D_N}$ and $\Phi^{D,\wt D}$ such that
\begin{equation}
P\Bigl(\,\,\max_{i=1,\dots,k}\,\,\bigl|\,\Phi^{D,\wt D}(x_i)-\varphi_N^{D,\wt D}(x_i)\bigr|>\epsilon\Bigr)
\,\underset{N\to\infty}\longrightarrow\,0.
\end{equation}
The claim will then follow if we can show that
\begin{equation}
\lim_{r\downarrow0}\,P\biggl(\,\,\sup_{\begin{subarray}{c}
x,y\in \wt D^\delta\\|x-y|<r
\end{subarray}}
\bigl|\Phi^{D,\wt D}(x)-\Phi^{D,\wt D}(y)\bigr|>\epsilon\biggr)=0
\end{equation}
and similarly (with an additional \textit{limes superior} as~$N\to\infty$ preceding the limit~$r\downarrow0$) for~$\Phi^{D,\wt D}(\cdot)$ replaced by~$\varphi_N^{D,\wt D}$. This, along with the continuity of~$\Phi^{D,\wt D}$, will follow from the regularity estimates on Gaussian processes to be proved in Lecture~\ref{lec-6} (specifically, in Lemma~\ref{lemma-6.17}). 
\end{proofsect}

For future reference, we suggest that the reader solve the following exercises:

\begin{myexercise}[Projection definition]
\label{ex:4.5}
Prove (by comparing covariances) that $\Phi^{D,\wt D}$ is the projection of the CGFF on~$D$ onto the subspace of functions in~$\cmss H^1_0(D)$ that are harmonic on (the connected components of)~$\wt D$.
\end{myexercise}

We will also need to note that the binding field has a nesting property:

\begin{myexercise}[Nesting property]
\label{ex:4.6}
Show that if~$U\subset V\subset W$, then
\begin{equation}
\varphi^{W,U}\,\,\laweq\,\,\varphi^{W,V}+\varphi^{V,U}\quad\text{with}\quad \varphi^{W,V}\independent\varphi^{V,U}.
\end{equation}
Similarly, if~$D''\subset D'\subset D$ are admissible, then (on~$D''$)
\begin{equation}
\Phi^{D,D''} \,\,\laweq\,\, \Phi^{D,D'}+\Phi^{D',D''}\quad\text{with}\quad \Phi^{D,D'}\independent\Phi^{D',D''}\,.
\end{equation}
\end{myexercise}

To demonstrate the benefits of taking the scaling limit we note that continuous binding field behaves nicely under conformal maps of the underlying domain. Indeed, we have:

\begin{myexercise}
Show that under a conformal bijection~$f\colon D\to f(D)$, we have
\begin{equation}
C^{f(D),f(\wt D)}\bigl(\,f(x),f(y)\bigr) = C^{D,\wt D}(x,y),\qquad x,y\in\wt D,
\end{equation}
for any $\wt D\subset D$. Prove that this implies
\begin{equation}
\Phi^{f(D),f(\wt D)}\circ f\,\,\laweq\,\, \Phi^{D,\wt D}.
\end{equation}
Hint: Use Lemma~\ref{lemma-1.26}.
\end{myexercise}

\section{Properties of $Z^D_\lambda$-measures}
\noindent
We are now ready to move to the discussion of the properties of (the laws of) the~$Z^D_\lambda$-measures. These will often relate the measures in different domains and so we need to ensure that the subsequential limit can be taken for all of these domains simultaneously. Applying Cantor's diagonal argument, we can definitely use the same subsequence for any \emph{countable} collection~$\mathfrak D_0\subset\mathfrak D$ of admissible domains and thus assume that a subsequential limit~$\eta^D$, and hence also the measure~$Z^D_\lambda$, has been extracted for each~$D\in\mathfrak D_0$. All of our subsequent statements will then be restricted to the domains in~$\mathfrak D_0$. We will always assume that~$\mathfrak D_0$ contains all finite unions of dyadic open and/or half-open squares plus any finite number of additional domains of present (or momentary) interest.

Some of the properties of the measures~$\{Z^D_\lambda\colon D\in\mathfrak D_0\}$ are quite elementary consequences of the above derivations and so we relegate them to:

\begin{myexercise}[Easy properties]
\label{ex:prop-Z}
Prove that:
\settowidth{\leftmargini}{(1111)}
\begin{enumerate}
\item[(1)] for each~$D\in\mathfrak D_0$, $Z^D_\lambda$ is supported on~$D$; i.e., $Z^D_\lambda(\partial D)=0$ a.s.,
\item[(2)] if $A\subset D\in\mathfrak D_0$ is measurable with $\Leb(A)=0$, then $Z_\lambda^D(A)=0$ a.s.,
\item[(3)] if $D,\wt D\in\mathfrak D_0$ obey $D\cap \wt D=\emptyset$, then 
\begin{equation}
\label{E:2.13uq}
Z_\lambda^{D\cup\wt D}(\textd x)\,\,\laweq\,\, Z_\lambda^{D}(\textd x)+Z_\lambda^{\wt D}(\textd x),
\end{equation}
with the measures $Z_\lambda^{D}$ and $Z_\lambda^{\wt D}$ on the right regarded as independent, and
\item[(4)] the law of $Z^D_\lambda$ is translation invariant; i.e.,
\begin{equation}
Z^{a+D}_\lambda(a+\textd x)\,\,\laweq\,\, Z^D_\lambda(\textd x)
\end{equation}
 for each $a\in\C$ such that~$D,a+D\in\mathfrak D_0$.
\end{enumerate}
\end{myexercise}

\noindent
As already mentioned, a key point for us is to prove the behavior under the restriction to a subdomain. We formulate this as follows:

\begin{myproposition}[Gibbs-Markov for~$Z^D_\lambda$-measures]
\label{prop-GM-Z}
For any $D,\wt D\in\mathfrak D_0$ satisfying $\wt D\subseteq D$ and~$\Leb(D\smallsetminus\wt D)=0$,
\begin{equation}
\label{E:4.15nw}
Z^D_\lambda(\textd x)\,\,\laweq\,\,\,\texte^{\alpha\lambda\Phi^{D,\wt D}(x)}\,Z^{\wt D}_\lambda(\textd x),
\end{equation}
where $\Phi^{D,\wt D}\independent Z^{\wt D}_\lambda$ with the laws as above.
\end{myproposition}

Note that, by Exercise~\ref{ex:prop-Z}(2), both sides of the expression assign zero mass to~$D\smallsetminus\wt D$ a.s. It therefore does not matter that~$\Phi^{D,\wt D}$ is not really defined there. (The decomposition in Fig.~\ref{fig-GM} shows that, even with $D$ and~$\wt D$ restricted to~$\mathfrak D_0$, non-trivial examples exist with~$\wt D\subset D$ and~$\Leb(D\smallsetminus\wt D)=0$.)

\begin{proofsect}{Proof}
Suppose $D,\wt D\in\mathfrak D_0$ obey $\wt D\subseteq D$ and $\Leb(D\smallsetminus\wt D)=0$.
Thanks to Exercise~\ref{ex:prop-Z}(2) and the Monotone Convergence Theorem, it suffices to prove equality in law for the integrals of the two measures with respect to any function bounded and continuous~$f\colon D\times\R\to\R$ with support in~$K\times[-b,b]$ for a compact~$K\subset\wt D$ and some~$b>0$.  

We start by invoking the Gibbs-Markov decomposition
\begin{equation}
h^{D_N} \,\,\laweq \,\,h^{\wt D_N}+\varphi^{D,\wt D}_N\quad\text{where}\quad h^{\wt D_N}\independent\varphi^{D,\wt D}_N.
\end{equation}
 A calculation then shows
\begin{equation}
\label{E:3.35}
\langle\eta^D_N,f\rangle \,\,\laweq\,\,\langle\eta^{\wt D}_N,f_\varphi\rangle\,,
\end{equation}
where
\begin{equation}
f_\varphi(x,h) := f\bigl(x,h+\varphi^{D_N,\wt D_N}_{\lfloor xN\rfloor}\bigr)\,.
\end{equation}
Next consider the coupling of~$\varphi^{D_N,\wt D_N}$ and~$\Phi^{D,\wt D}$ from Lemma~\ref{lemma-4.4} where we may and will assume $\Phi^{D,\wt D}\independent \eta^{\wt D}_N$. Our aim is to replace~$f_\varphi$ by
\begin{equation}
\label{E:4.20ue}
f_\Phi(x,h):=f\bigl(x,h+\Phi^{D,\wt D}(x)\bigr)
\end{equation}
in \eqref{E:3.35}. Given~$\epsilon>0$ let~$\delta>0$ be such that for all~$x\in \wt D$ and all~$h,h'\in[-b,b]$ with~$|h-h'|<\delta$ we have~$|f(x,h)-f(x,h')|<\epsilon$. Then, on the event
\begin{equation}
\label{E:4.20uai}
\bigl\{\sup_{x\in K}|\Phi^{D,\wt D}(x)|\le M\bigr\}\cap\bigl\{\sup_{x\in K}|\varphi^{D_N,\wt D_N}_{\lfloor xN\rfloor}-\Phi^{D,\wt D}(x)|<\delta\bigr\},
\end{equation}
we get 
\begin{equation}
\label{E:4.21uai}
\bigl|\langle\eta^{\wt D}_N,f_\varphi\rangle -\langle\eta^{\wt D}_N,f_\Phi\rangle\bigr|
\le\eta^{\wt D}_N\bigl(\wt D\times[-b-M-\delta,\infty)\bigr)\epsilon.
\end{equation}
By Lemmas~\ref{lemma-3.2} and~\ref{lemma-4.4}, the probability of the event in \eqref{E:4.20uai} tends to one while the right-hand side of \eqref{E:4.21uai} tends to zero in probability as $N\to\infty$ followed by $\epsilon\downarrow0$ and~$M\to\infty$. Any simultaneous subsequential limits $\eta^D$, resp.,~$\eta^{\wt D}$ of $\{\eta^D_N\colon N\ge1\}$, resp., $\{\eta^{\wt D}_N\colon N\ge1\}$ therefore obey
\begin{equation}
\label{E:4.21ue}
\langle \eta^D,f\rangle \,\,\laweq\,\, \langle\eta^{\wt D},f_\Phi\rangle\,,
\end{equation}
where $\Phi^{D,\wt D}$ (implicitly contained in $f_\Phi$) is independent of $\eta^{\wt D}$. 

The representation from Proposition~\ref{prop-factorization} now permits us to write
\begin{equation}
\begin{aligned}
\int_{D\times\R}Z^D_\lambda(\textd x)&\otimes\textd h \,\,\texte^{-\alpha\lambda h}\,f(x,h)
=\langle\eta^D,f\rangle
\,\,\laweq\,\,
\langle\eta^{\wt D},f_\Phi\rangle 
\\
&= \int_{D\times\R}Z^{\wt D}_\lambda(\textd x)\otimes\textd h \,\,\texte^{-\alpha\lambda h}\,f\bigl(x,h+\Phi^{D,\wt D}(x)\bigr)
\\
&=\int_{D\times\R}Z^{\wt D}_\lambda(\textd x)\otimes\textd h \,\,\texte^{\alpha\lambda\Phi^{D,\wt D}(x)}\,\texte^{-\alpha\lambda h}\,f(x,h).
\end{aligned}
\end{equation}
As noted at the beginning of the proof, this implies the claim.
\end{proofsect}

\section{Representation via Gaussian multiplicative chaos}
\noindent
We will now show that the above properties determine the laws of $\{Z^D_\lambda\colon D\in\mathfrak D_0\}$ uniquely. To this end we restrict our attention to dyadic open squares, i.e., those of the form
\begin{equation}
2^{-n}z+(0,2^{-n})^2\quad\text{for}\quad z\in\Z^2\,\,\text{and}\,\, n\ge0.
\end{equation}
 For a fixed~$m\in\Z$, set $S:=(0,2^{-m})^2$ and let $\{S_{n,i}\colon i=1,\dots,4^n\}$ be an enumeration of the dyadic squares of side~$2^{-(n+m)}$ that have a non-empty intersection with~$S$. (See Fig.~\ref{fig-GM} for an illustration of this setting.) Recall that we assumed that~$\mathfrak D_0$ contains all these squares which makes  $Z^D_\lambda$ defined on all of them. 

Abbreviating
\begin{equation}
\tilde S^n:=\bigcup_{i=1}^{4^n}S_{n,i}
\end{equation}
the Gibbs-Markov decomposition \eqref{E:4.15nw} gives
\begin{equation}
\label{E:4.21}
Z^S_\lambda(\textd x)\,\,\laweq\,\,\sum_{i=1}^{4^n}\texte^{\alpha\lambda\Phi^{S,\tilde S^n}(x)}\,Z^{S_{n,i}}_\lambda(\textd x),
\end{equation}
where the measures $\{Z^{S_{n,i}}_\lambda\colon i=1,\dots,4^n\}$ and the binding field~$\Phi^{S,\tilde S^n}$ on the right-hand side are regarded as independent. The expression \eqref{E:4.21} links~$Z_\lambda$-measure in one set in terms of~$Z_\lambda$-measures in a scaled version thereof. This suggests we think of \eqref{E:4.21} as a fixed point of a \emph{smoothing transformation}; cf e.g., Durrett and Liggett~\cite{Durrett-Liggett}. Such fixed points are found by studying the variant of the fixed point equation where the object of interest on the right is replaced by its expectation. This leads to the consideration of the measure
\begin{equation}
\label{E:4.22}
Y_n^S(\textd x):=\hat c\sum_{i=1}^{4^n}1_{S_{n,i}}(x)\,\texte^{\alpha\lambda\Phi^{S,\tilde S^n}(x)}\,r_{\tilde S^n}(x)^{2\lambda^2}\,\textd x. 
\end{equation}
where~$\hat c$ is as in \eqref{E:3.58}. Indeed, in light of the fact that
\begin{equation}
r_{\tilde S^n}(x) = r_{S_{n,i}}(x)\quad\text{for}\quad
x\in S_{n,i}
\end{equation}
we have 
\begin{equation}
\label{E:4.23}
Y_n^S(A) = E\bigl[Z^S_\lambda(A)\,\big|\,\sigma(\Phi^{S,\tilde S^n})\bigr]
\end{equation}
for any Borel $A\subset\C$. The next point to observe is that these measures can be interpreted in terms of \emph{Gaussian multiplicative chaos} (see Lemma~\ref{lemma-GMC}):

\begin{mylemma}
\label{lemma-4.9}
There is a random measure~$Y_\infty^S$ such that for all Borel~$A\subset\C$,
\begin{equation}
\label{E:4.30nw}
Y^S_n(A)\,\,\underset{n\to\infty}\lawarrow\,\, Y^S_\infty(A).
\end{equation}
\end{mylemma}

\begin{proofsect}{Proof}
Denoting~$\tilde S^0:=S$, the nesting property of the binding field (Exercise~\ref{ex:4.6}) allows us to represent   $\{\Phi^{S,\tilde S^n}\colon n\ge1\}$ on the same probability space via
\begin{equation}
\label{E:4.25}
\Phi^{S,\tilde S^n} :=\sum_{k=0}^{n-1}\Phi^{\tilde S^k,\tilde S^{k+1}}\,,
\end{equation}
where the fields $\{\Phi^{\tilde S^k,\tilde S^{k+1}}\colon k\ge0\}$ are independent with their corresponding laws. In this representation, the measures~$Y_n$ are defined all on the same probability space and so we can actually prove the stated convergence in almost-sure sense. 
Indeed, \eqref{E:rad}, \eqref{E:4.1} and $\alpha^2g=4$ imply
\begin{equation}
\label{E:4.33nwwwt}
\wt D\subseteq D\quad\Rightarrow\quad r_{\wt D}(x)^{2\lambda^2}=r_D(x)^{2\lambda^2}
\texte^{\frac12\alpha^2\lambda^2\Var[\Phi^{D,\wt D}(x)]},\quad x\in \wt D.
\end{equation}
This permits us to rewrite \eqref{E:4.22} as
\begin{equation}
\label{E:4.32a}
Y_n^S(\textd x):=\hat c\,r_S(x)^{2\lambda^2}
\sum_{i=1}^{4^n}1_{S_{n,i}}(x)\,\texte^{\alpha\lambda\Phi^{S,\tilde S^n}(x)-\frac12\alpha^2\lambda^2\Var[\Phi^{D,\wt D}(x)]}\,\textd x
\end{equation}
and thus cast~$Y_n^S$ in the form we encountered in the definition of the Gaussian Multiplicative Chaos. Adapting the proof of Lemma~\ref{lemma-GMC} (or using it directly with the help of Exercise~\ref{ex:4.5}), we get \eqref{E:4.30nw} any Borel~$A\subset\C$.
\end{proofsect}

We now claim:

\begin{myproposition}[Characterization of $Z^D_\lambda$ measure]
\label{prop-4.12}
For any dyadic squa\-re $S\subset\C$ and any bounded and continuous function~$f\colon S\to[0,\infty)$, we have
\begin{equation}
\label{E:4.29}
E\bigl(\texte^{-\langle Z^S_\lambda,f\rangle}\bigr)=E\bigl(\texte^{-\langle Y^S_\infty,f\rangle}\bigr).
\end{equation}
In particular,
\begin{equation}
\label{E:4.30}
Z^S_\lambda(\textd x)\,\,\laweq\,\,Y^S_\infty(\textd x).
\end{equation}
\end{myproposition}

\begin{proofsect}{Proof of ``$\ge$'' in \eqref{E:4.29}}
Writing~$Z^S_\lambda$ via \eqref{E:4.21} and invoking conditional expectation given $\Phi^{S,\tilde S^n}$ with the help of \eqref{E:4.23}, the conditional Jensen inequality shows
\begin{equation}
\begin{aligned}
E\bigl(\texte^{-\langle Z^S_\lambda,f\rangle}\bigr)
&=E\Bigl(E\bigl(\texte^{-\langle Z^S_\lambda,f\rangle}\,\big|\,\sigma(\Phi^{S,\tilde S^n})\bigr)\Bigr)
\\
&\ge E\bigl(\texte^{-E[\langle Z^S_\lambda,f\rangle\,|\,\sigma(\Phi^{S,\tilde S^n}]}\bigr)
=E\bigl(\texte^{-\langle Y^S_n,f\rangle}\bigr)
\end{aligned}
\end{equation}
for any continuous $f\colon S\to[0,\infty)$.
The convergence in Lemma~\ref{lemma-4.9} implies
\begin{equation}
E(\texte^{-\langle Y^S_n,f\rangle})\,\underset{n\to\infty}\longrightarrow\, E(\texte^{-\langle Y^S_\infty,f\rangle})
\end{equation}
and so we get~``$\ge$'' in \eqref{E:4.29}.
\end{proofsect}

For the proof of the opposite inequality in \eqref{E:4.29} we first note:

\begin{mylemma}[Reversed Jensen's inequality]
\label{lemma-3.11a}
If $X_1,\dots,X_n$ are independent non-negative  random variables, then for each~$\epsilon>0$,
\begin{equation}
\label{E:4.32}
E\biggl(\exp\Bigl\{-\sum_{i=1}^n X_i\Bigr\}\biggr)\le \exp\Bigl\{-\texte^{-\epsilon}\sum_{i=1}^n E(X_i\,;\, X_i\le \epsilon)\Bigr\}\,.
\end{equation}
\end{mylemma}

\begin{proofsect}{Proof}
In light of assumed independence, it suffices to prove this for~$n=1$. This is checked by bounding $E(\texte^{-X})\le E(\texte^{-\wt X})$, where $\wt X :=X\1_{\{X\le\epsilon\}}$, writing 
\begin{equation}
-\log E(\texte^{-\wt X})=\int_0^1\textd s\,\, \frac{E(\wt X\texte^{-s\wt X})}{E(\texte^{-s \wt X})}
\end{equation}
and invoking $E(\wt X\texte^{-s\wt X})\ge \texte^{-\epsilon} E(\wt X)$ and $E(\texte^{-s \wt X})\le1$.
\end{proofsect}

We are now ready to give:

\begin{proofsect}{Proof of ``$\le$'' in \eqref{E:4.29}}
Pick~$n$ large and assume~$Z^S_\lambda$ is again represented via \eqref{E:4.21}. We first invoke an additional truncation: Given~$\delta>0$, let~$S_{n,i}^\delta$ be the translate of~$(\delta2^{-(n-m)},(1-\delta)2^{-(n-m)})$ centered at the same point as~$S_{n,i}$. Denote
\begin{equation}
\tilde S^n_\delta:=\bigcup_{i=1}^{4^n}S_{n,i}^\delta
\quad\text{and}\quad f_{n,\delta}(x):=f(x)1_{\tilde S^n_\delta}(x).
\end{equation}
Denoting also
\begin{equation}
X_i:=\int_{S_{n,i}} f_{n,\delta}(x)\,\texte^{\alpha\lambda\Phi^{S,\tilde S^n}}\,Z^{S_{n,i}}_\lambda(\textd x)
\end{equation}
from~$f\ge f_{n,\delta}$ we then have
\begin{equation}
E\bigl(\texte^{-\langle Z^S_\lambda,f\rangle}\bigr)
\le E\bigl(\texte^{-\langle Z^S_\lambda,f_{n,\delta}\rangle}\bigr)
=E\biggl(\exp\Bigl\{-\sum_{i=1}^n X_i\Bigr\}\biggr).
\end{equation}
Conditioning on~$\Phi^{S,\tilde S^n}$, the bound \eqref{E:4.32} yields
\begin{equation}
\label{E:4.37}
E\bigl(\texte^{-\langle Z^S_\lambda,f\rangle}\bigr)
\le E\biggl(\exp\Bigl\{-\texte^{-\epsilon}\sum_{i=1}^{4^n} E\bigl(X_i1_{\{X_i\le\epsilon\}}\,\big|\,\sigma(\Phi^{S,\tilde S^n})\bigr)\Bigr\}\biggr)\,.
\end{equation}
Since \eqref{E:4.23} shows
\begin{equation}
\sum_{i=1}^{4^n}E\bigl(X_i\,\big|\,\sigma(\Phi^{S,\tilde S^n})\bigl) = \langle Y^S_n,\,f_{n,\delta}\rangle\,,
\end{equation}
we will also need:

\begin{mylemma}
\label{lemma-4.13}
Assume $\lambda\in(0,1/\sqrt2)$. Then for each~$\epsilon>0$,
\begin{equation}
\label{E:4.39}
\lim_{n\to\infty}\,\sum_{i=1}^{4^n} E\bigl(X_i;\,X_i>\epsilon\bigr)=0\,.
\end{equation}
\end{mylemma}

\noindent
Postponing the proof until after that of Proposition~\ref{prop-4.12}, from \eqref{E:4.37} and \eqref{E:4.39} we now get
\begin{equation}
E\bigl(\texte^{-\langle Z^S_\lambda,f\rangle}\bigr)
\le \limsup_{n\to\infty}E(\texte^{-\texte^{-\epsilon}\langle Y^S_n,f_{n,\delta}\rangle})\,.
\end{equation}
But 
\begin{equation}
\langle Y^S_n,f_{n,\delta}\rangle \ge \langle Y^S_n,f\rangle - \Vert f\Vert Y^S_n(S\smallsetminus\tilde S^n_\delta)
\end{equation}
and so
\begin{equation}
E(\texte^{-\texte^{-\epsilon}\langle Y^S_n,f_{n,\delta}\rangle})
\le \texte^{\epsilon\Vert f\Vert}\,E(\texte^{-\texte^{-\epsilon}\langle Y^S_n,f\rangle})+P\bigl(Y^S_n(S\smallsetminus\tilde S^n_\delta)>\epsilon\bigr)\,.
\end{equation}
A calculation based on \eqref{E:4.32a} shows
\begin{equation}
P\bigl(Y^S_n(S\smallsetminus\tilde S^n_\delta)>\epsilon\bigr)\le c\epsilon^{-1} \Leb(S\smallsetminus\tilde S^n_\delta)\le c'\epsilon^{-1}\delta.
\end{equation}
Invoking also \eqref{E:4.30nw} we get
\begin{equation}
E\bigl(\texte^{-\langle Z^S_\lambda,f\rangle}\bigr)
\le
\lim_{\epsilon\downarrow0}\,\limsup_{\delta\downarrow0}\,\limsup_{n\to\infty}E(\texte^{-\texte^{-\epsilon}\langle Y^S_n,f_{n,\delta}\rangle})
\le E(\texte^{-\langle Y^S_\infty,f\rangle}).
\end{equation}
This completes the proof of~\eqref{E:4.29};  \eqref{E:4.30} then directly follows.
\end{proofsect}

It remains to give:

\begin{proofsect}{Proof of Lemma~\ref{lemma-4.13}}
First we note
\begin{equation}
\label{E:4.40}
\begin{aligned}
\sum_{i=1}^{4^n} &E\bigl(X_i;\,X_i>\epsilon\bigr)
\le\frac1\epsilon 
\sum_{i=1}^{4^n} E(X_i^2)
\\
&\le\frac{\Vert f\Vert^2}\epsilon\sum_{i=1}^{4^n}E\biggl(\,\int_{S_{n,i}^\delta\times S_{n,i}^\delta}E\bigl(\texte^{\alpha\lambda[\Phi^{S,\tilde S^n}(x)+\Phi^{S,\tilde S^n}(y)]}\bigr)Z^{S_{n,i}}_\lambda(\textd x)Z^{S_{n,i}}_\lambda(\textd y)\biggr).
\end{aligned}
\end{equation}
Denote~$L:=2^n$. In light of the fact that, for some constant~$c$ independent of~$n$,
\begin{equation}
\Var(\Phi^{S,\tilde S^n}(x))=g\log\frac{r_S(x)}{r_{S_{n,i}}(x)}\le g\log(L) + c
\end{equation}
holds uniformly in~$x\in\tilde S^n_\delta$, \eqref{E:4.40} is bounded by $\Vert f\Vert^2/\epsilon^2$ times
\begin{multline}
\quad
c' \texte^{4\frac12\alpha^2\lambda^2g\log(L)}\sum_{i=1}^{4^n}\,E\bigl[Z^{S_{n,i}}_\lambda(S_{n,i})^2\bigr]
\\
\le c'' L^{8\lambda^2+2-(4+4\lambda^2)}=c'' L^{-2(1-2\lambda^2)},
\quad
\end{multline}
where we also used $\alpha^2g=4$, invoked \eqref{E:3.60} and noted that there are~$L^2$ terms in the sum. In light of $\lambda<1/\sqrt2$, this tends to zero as~$L\to\infty$.
\end{proofsect}

\section{Finishing touches}
\noindent
We are now ready to combine the above observations to get the main conclusion about the existence of the limit of processes~$\{\eta^D_N\colon N\ge1\}$:

\begin{proofsect}{Proof of Theorem~\ref{thm-intermediate} for~$\lambda<1/\sqrt2$}
Pick~$D\in\mathfrak D$ and assume, as discussed before, that~$\mathfrak D_0$ used above contains~$D$. For any subsequence of~$N$'s for which the limit of the measures in question exists for all domains in~$\mathfrak D_0$ we then have the representation \eqref{E:3.46} as well as the properties stated in Exercise~\ref{ex:prop-Z} and Proposition~\ref{prop-GM-Z}. We just have to show that the limit measure $Z^D_\lambda$ is determined by these properties and the observation from Proposition~\ref{prop-4.12}. This will also prove the existence of the limit of $\{\eta^D_N\colon N\ge1\}$. 

By Exercise~\ref{ex:prop-Z}(1,2) it suffices to show that the law of $\langle Z^D_\lambda,f\rangle$ is determined for any continuous~$f\colon D\times\R\to\R$ with compact support. Let~$D^n$ be the union of all the open dyadic  squares of side~$2^{-n}$ entirely contained in~$D$. Letting~$n$ be so large that $\supp(f)\subseteq D^n\times[-n,n]$ and denoting $\wt D^n:=D\smallsetminus\partial D^n$, Proposition~\ref{prop-GM-Z} and Exercise~\ref{ex:prop-Z}(3) then imply 
\begin{equation}
\label{E:4.51}
\langle Z^D_\lambda,f\rangle\,\,\laweq\,\,\bigl\langle Z^{D^n}_\lambda,\,\texte^{\alpha\lambda\Phi^{D,\wt D^n}}f\bigr\rangle\quad\text{with}\quad  Z^{D^n}_\lambda\independent\Phi^{D,\wt D^n}.
\end{equation}
It follows that the law of the left-hand side is determined once the law of~$Z^{D^n}_\lambda$ is determined. But part~(3) of Exercise~\ref{ex:prop-Z} also shows (as in \eqref{E:4.21}) that the measure $Z^{D^n}_\lambda$ is the exponential of the continuum binding field times the sum of independent copies of~$Z^S_\lambda$ for~$S$ ranging over the dyadic squares constituting~$D^n$. The laws of these~$Z^S_\lambda$ are determined uniquely by Proposition~\ref{prop-4.12} and hence so are those of~$Z^{D^n}_\lambda$ and~$Z^D_\lambda$ as well.
\end{proofsect}

Concerning the proof of conformal invariance and full characterization by the LQG measure, we will need the following result:

\begin{mytheorem}[Uniqueness of the GMC/LQG measure]
\label{thm-GMC-unique}
The law of the Gaussian Multiplicative Chaos measure $\mu^{D,\beta}_\infty$ does not depend on the choice of the orthonormal basis in~$\cmss H_0^1(D)$ that was used to define it (see Lemma~\ref{lemma-GMC}).
\end{mytheorem}

We will not prove this theorem in these notes as that would take us on a tangent that we do not wish to follow. We remark that the result has a rather neat proof due to Shamov~\cite{Shamov} which is made possible by his ingenious characterization of the GMC measures using Cameron-Martin shifts. An earlier work of Kahane~\cite{Kahane} required uniform convergence of the covariances of the approximating fields; we state and prove this version in Theorem~\ref{thm-5.5}. This version also suffices to establish a conformal transformation rule for the limit (cf Exercise~\ref{ex:Kahane-CFI}).

Equipped with the uniqueness claim in Theorem~\ref{thm-GMC-unique}, let us now annotate the steps that identify~$Z^D_\lambda$ with the LQG-measure $\hat cr_D(x)^{2\lambda^2}\mu^{D,\lambda\alpha}_\infty(\textd x)$. Let us start with a unit square~$S$. For each~$k\ge0$, let~$\{S_{k,i}\colon i=1,\dots,n(k)\}$ be the collection of open dyadic squares of side~$2^{-k}$ that are entirely contained in~$S$.  Denote
\begin{equation}
D^k:=\bigcup_{j=1}^{n(k)}S_{k,j},\quad k\ge0,
\end{equation}
with $D^{-1}:=S$ and observe that $D^k\subset D^{k-1}$ for each~$k\ge0$. Observe that~$\partial D^k$ is a collection of horizontal and vertical lines. Note also that~$\partial D^{k-1}\subset\partial D^k$.

\begin{myexercise}
\label{ex:4.15}
For each~$k\ge0$, let~$\cmss H_k$ denote the subspace of functions in~$\cmss H_0^1(S)$ that are harmonic in~$D^k$ and vanish on~$\partial D^{k-1}$. Prove that
\begin{equation}
\cmss H_0^1(S) = \bigoplus_{k=0}^\infty\cmss H_k\,.
\end{equation}
\end{myexercise}

\noindent
Next, for each~$k\ge1$, let $\{\wt f_{k,j}\colon j\ge1\}$ be an orthonormal basis in~$\cmss H_k$ with respect to the Dirichlet inner product. Then show:

\begin{myexercise}
Prove that, for $\{X_{k,j}\colon k,j\ge1\}$ i.i.d.\ standard normals,
\begin{equation}
\label{E:4.53}
\Phi^{D^{k-1},D^k}\,\,\laweq\,\,\sum_{j\ge1}X_{k,j}\wt f_{k,j}\quad\text{on}\,\,D^k
\end{equation}
holds for all $k\ge1$ with the sums converging locally uniformly in~$D^k$, a.s. Conclude that for all $m\ge1$,
\begin{equation}
\Phi^{D,D^m}\,\,\laweq\,\,\sum_{k=1}^m\sum_{j\ge1}X_{k,j}\wt f_{k,j}\quad\text{on}\,\,D^m\,.
\end{equation}
Hint: Define the~$X_{k,j}$'s by suitable inner products and check the covariances.
\end{myexercise}

\noindent
From here and Theorem~\ref{thm-GMC-unique} we now infer the following:

\begin{myexercise}
\label{ex:Z-mu}
Using suitable test functions, and a convenient enumeration of the above orthogonal basis in~$\cmss H_0^1(S)$, show that
\begin{equation}
\label{E:3.59uai}
Z^S_\lambda(\textd x) \,\,\laweq\,\,
Y_\infty^S(\textd x)\,\,\laweq\,\,\hat c\,r_S(x)^{2\lambda^2}\mu^{S,\lambda\alpha}_\infty(\textd x).
\end{equation}
Hint: Use either the argument based on Jensen's inequality from the previous section or invoke Kahane's convexity inequality from Proposition~\ref{prop-5.6}.
\end{myexercise}

\noindent
The point here is that although the right-hand side \eqref{E:4.53} casts the binding field from~$D$ to~$D^m$ in the form akin to \eqref{E:2.17}, the sum over~$j$ is infinite. One thus has to see that a suitable truncation to a finite sum will do as well.

Exercise~\ref{ex:Z-mu} identifies $Z^D_\lambda$ with the LQG measure for~$D$ a dyadic square. To extend this to general domains, which can be reduced to dyadic squares using the Gibbs-Markov property, we also need to solve:

\begin{myexercise}[Gibbs-Markov for LQG measure]
Suppose that domains $\wt D\subset D$ obey $\Leb(D\smallsetminus\wt D)=0$. Use \eqref{E:4.33nwwwt} to prove that, for all $\beta\in(0,\alpha)$,
\begin{equation}
\label{E:4.60uai}
r_D(x)^{2(\beta/\alpha)^2}\mu^{D,\beta}_\infty(\textd x)\,\,\laweq\,\,\texte^{\beta\Phi^{D,\wt D}(x)}\,\,r_{\wt D}(x)^{2(\beta/\alpha)^2}\mu^{\wt D,\beta}_\infty(\textd x)\,,
\end{equation}
where $\Phi^{D,\wt D}$ and $\mu^{\wt D,\beta}_\infty$ are regarded as independent.
\end{myexercise}

We remark that, in the regime~$\lambda<1/\sqrt2$, formula \eqref{E:4.60uai} would be enough to complete \eqref{E:3.59uai}; indeed, our Proposition~\ref{prop-4.12} only required the Gibbs-Markov property, the expectation formula \eqref{E:3.59} and the bounds on the second moment in \twoeqref{E:3.60}{e:2.3b}, which are known to hold for the LQG measure as well. However, this does not apply in the complementary regime of~$\lambda$'s where we use a truncation on the underlying field that is hard to interpret in the language of the above random measures.

Once we identify the limit measure with the LQG-measure --- and accept Theorem~\ref{thm-GMC-unique} without proof  --- the proof of the conformal transformation rule in Theorem~\ref{thm-2.4} is quite easy. A key  point is to solve:

\begin{myexercise}[Conformal transform of GMC measure]
\label{ex:4.18}
Let~$f\colon D\mapsto f(D)$ be a conformal map between bounded and open domains in~$\C$. Show that if~$\{g_n\colon n\ge1\}$ is an orthonormal basis in~$\cmss H_0^1(D)$ with respect to the Dirichlet inner product, then $\{g_n\circ f^{-1}\colon n\ge1\}$ is similarly an orthonormal basis in~$\cmss H_0^1(f(D))$. Prove that this implies
\begin{equation}
\mu^{f(D),\beta}_\infty\circ f(\textd x)\,\,\laweq\,\,\bigl|f'(x)\bigr|^2\mu^{D,\beta}_\infty(\textd x).
\end{equation}
\end{myexercise}

\noindent
To get the proof of Theorem~\ref{thm-2.4}, one then needs to observe:

\begin{myexercise}
For any conformal bijection~$f\colon D\to f(D)$,
\begin{equation}
r_{f(D)}\bigl(f(x)\bigr)=\bigl|f'(x)\bigr|\,r_D(x),\qquad x\in D.
\end{equation}
\end{myexercise}

\noindent
We reiterate that a proof of conformal invariance avoiding the full strength of Theorem~\ref{thm-GMC-unique} will be posed as Exercise~\ref{ex:Kahane-CFI}.

\section{Dealing with truncations}
\label{sec:truncate}\noindent
The above completes the proof of our results in the regime where second-moment calculations can be applied without truncations. To get some feeling for what happens in the the complementary regime, $1/\sqrt2\le\lambda<1$, let us at least introduce the basic definitions and annotate the relevant steps. 

\begin{figure}[t]
\vglue-2mm
\centerline{\includegraphics[width=0.5\textwidth]{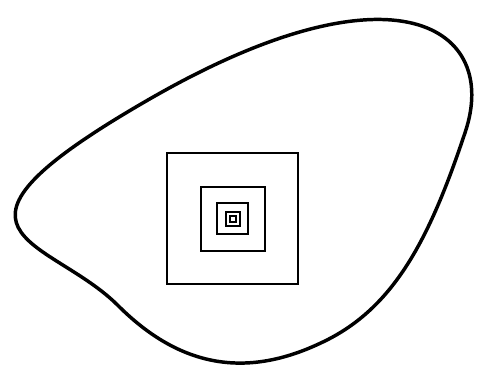}
}
\vglue-1mm
\begin{quote}
\small 
\caption{
\label{fig-squares}
\small
An illustration of the collection of sets~$\Delta^k(x)$ above. The domain~$D_N$ corresponds to the region marked by the outer curve.
}
\normalsize
\end{quote}
\end{figure}

Denote~$\Lambda_r(x):=\{y\in\Z^d\colon\dist_\infty(x,y)\le r\}$. Given a discretized version~$D_N$ of a continuum domain~$D$, for each~$x\in D_N$ define
\begin{equation}
\label{e:4.2new}
\Delta^k(x) := \begin{cases}
	\emptyset			  & \text{for } k = 0\,,\\
	\Lambda_{\texte^k}(x) & \text{for } k= 1, \dots, n(x)-1 \,, \\
	D_N 				  & \text{for } k = n(x) \,,
\end{cases} 
\end{equation}
where $n(x):=\max\{n\ge0\colon\Lambda_{\texte^{n+1}}(x)\subseteq D_N\}$. See Fig.~\ref{fig-squares} for an illustration.

\nopagebreak
\vglue-1mm

Using the definition of~$\varphi^{D_N,\Delta^k}(x)$ as a conditional expectation of~$h^{D_N}$, we now define the truncation event
\begin{equation}
T_{N,M}(x):=\bigcap_{k= k_N}^{n(x)}\biggl\{\Bigl|\varphi^{D_N,\Delta^k}_x- a_N\frac{n(x)-k}{n(x)}\Bigr|\le M\bigl[n(x)-k\bigr]^{3/4}\biggr\}\,,
\end{equation}
 where $M$ is a parameter and $k_N:=\frac18\log(K_N)\approx\frac14(1-\lambda^2)\log N$.
Next we introduce the \emph{truncated point measure}
\begin{equation}
\label{E:3.19anew}
\wh\eta_N^{D,M}:=\frac1{K_N}\sum_{x\in D_N}\1_{T_{N,M}(x)}\,\delta_{x/N}\otimes\delta_{h^{D_N}_x-a_N}.
\end{equation}
The following elementary inequality will be quite useful:
\begin{equation}
\label{E:3.20anew}
\eta^D_N\ge\wh\eta^{D,M}_N,\qquad M\in(0,\infty).
\end{equation}
Indeed, the tightness of~$\{\wh\eta_N^D\colon N\ge1\}$ is thus inherited from $\{\eta_N^D\colon N\ge1\}$ and, as can be shown, the limit points of the former increase to those of the latter as~$M\to\infty$. The requisite (now really ugly) second moment calculations are then performed which yield the following conclusions for all~$M<\infty$ and all~$\lambda\in(0,1)$:
\begin{enumerate}
\item[(1)] Defining $\wh\Gamma^{D,M}_N(b):=\{x\in D_N\colon h^{D_N}_x\ge a_N+b,\, T_{N,M}(x)\text{ occurs}\}$, we have
\begin{equation}
\sup_{N\ge1}\frac1{K_N^2}E\bigl(|\wh\Gamma^{D,M}_N(b)|^2\bigr)<\infty.
\end{equation}
By a second-moment argument, the limits of $\{\wh\eta_N^{D,M}\colon N\ge1\}$ are non-trivial. 
\item[(2)] The factorization property proved in Proposition~\ref{prop-factorization} applies to limit points of $\{\wh\eta_N^{D,M}\colon N\ge1\}$ with~$Z^D_\lambda$ replaced by some~$\wh Z^{D,M}_\lambda$ instead. 
\end{enumerate}
The property \eqref{E:3.20anew} now implies that $M\mapsto \wh Z^{D,M}_\lambda$ is pointwise increasing and so we may define
\begin{equation}
\label{E:4.67}
Z^D_\lambda(\cdot):=\lim_{M\to\infty}\wh Z^{D,M}_\lambda(\cdot).
\end{equation}
We then check that this measure has the properties in Exercise~\ref{ex:prop-Z} as well as the Gibbs-Markov property from Proposition~\ref{prop-GM-Z}. However, although the limit in \eqref{E:4.67} exists in~$L^1$ for all~$\lambda\in(0,1)$, it does not exist in~$L^2$ for~$\lambda\ge1/\sqrt2$ --- because, e.g., by Fatou's Lemma, the limit LQG measure would then be in $L^2$ as well which is known to be false --- and so we have to keep using $\wh Z^{D,M}_\lambda$ whenever estimates involving second moments are needed. 

This comes up only in one proof: ``$\le$'' in \eqref{E:4.29}. There we use the inequality $Z^D_\lambda(\cdot)\ge \wh Z^{D,M}_\lambda(\cdot)$ and the fact that~$Z^D_\lambda$ satisfies the Gibbs-Markov property  to dominate $Z^D_\lambda(\cdot)$ from below by the measure
\begin{equation}
\label{E:4.21b}
\wt Z^S_\lambda(\textd x):=\sum_{i=1}^{4^n}\texte^{\alpha\lambda\Phi^{S,\tilde S^n}}\,\wh Z^{S_{n,i},M}_\lambda(\textd x).
\end{equation}
Then we perform the calculation after \eqref{E:4.21} with this measure instead of $Z^S_\lambda$ modulo one change: In the proof of Lemma~\ref{lemma-4.13} we truncate to the event
\begin{equation}
\sup_{x\in \wt S^n_\delta}\Phi^{S,\wt S^n}(x) <  2\sqrt g\log(2^n)+c\sqrt{\log(2^n)}\,,
\end{equation}
which has probability very close to one. On this event, writing again~$L:=2^{n}$, the sum on the right-hand side of \eqref{E:4.40} is thus bounded by 
\begin{equation}
c'\texte^{2\sqrt g\alpha\lambda\log(L)+c\sqrt{\log(L)}}\,\texte^{\frac12\alpha^2\lambda^2 g\log(L)}L^2L^{-2(2+2\lambda^2)}.
\end{equation}
Using the definition of~$\alpha$, this becomes $L^{-2(1-\lambda)^2+o(1)}$ which vanishes as~$n\to\infty$ for all~$\lambda\in(0,1)$. The rest of the proof is then more or less the same.

The full proof of Theorem~\ref{thm-intermediate} is carried out in~\cite{BL4} to which we refer the reader for further details.


\chapter{Gaussian comparison inequalities}
\label{lec-5}
\noindent
In our discussion we have so far managed to get by using only elementary facts about Gaussian processes. The forthcoming derivations will require more sophisticated techniques and so it is time we addressed them properly. In this lecture we focus on Gaussian comparison inequalities, starting with Kahane's inequality and its corollaries called the Slepian Lemma and the Sudakov-Fernique inequality. We give an application of Kahane's inequality to uniqueness of the Gaussian Multiplicative Chaos, a subject touched upon before. In the last section, we will review the concepts of stochastic domination and the FKG inequality. The presentation draws on Adler~\cite{Adler}, Adler and Taylor~\cite{Adler-Taylor} and Liggett~\cite{Liggett-IPS}.

\section{Kahane's inequality}
\noindent
Let us first make a short note on terminology: We say that a (multivariate) Gaussian $X=(X_1,\dots,X_n)$ is \emph{centered} if~$E(X_i)=0$ for all~$i=1,\dots,n$. A function $f\colon\R^n\to\R$ is said to have a \emph{subgaussian growth} if for each~$\epsilon>0$ there is~$C>0$ such that $|f(x)|\le C\texte^{\epsilon|x|^2}$ holds for all~$x\in\R^n$.

In his development of the theory of Gaussian Multiplicative Chaos, Kahane made convenient use of inequalities that, generally, give comparison estimates of expectation of functions (usually convex in appropriate sense) of Gaussian random variables whose covariances can be compared in a pointwise sense. One version of this inequality is as follows:

\begin{myproposition}[Kahane inequality]
\label{thm-5.1}
Let~$X,Y$ be centered Gaussian vectors on~$\R^n$ and~$f\in C^2(\R^n)$ a function whose second derivatives have a subgaussian growth. Assume
\begin{equation}
\forall i,j=1,\dots,n:\quad\left\{\begin{alignedat}{3}
E(Y_iY_j)&> E(X_iX_j)\quad&\Rightarrow\quad &\frac{\partial f}{\partial x_i\partial x_j}(x)\ge 0,\quad x\in\R^n
\\*[2mm]
E(Y_iY_j)&< E(X_iX_j)\quad&\Rightarrow\quad &\frac{\partial f}{\partial x_i\partial x_j}(x)\le 0,\quad x\in\R^n
\end{alignedat}
\right.
\end{equation}
Then
\begin{equation}
E f(Y)\ge E f(X).
\end{equation}
\end{myproposition}

\noindent
Note also that for pairs~$i,j$ such that $E(Y_iY_j)=E(X_iX_j)$ the sign of $\frac{\partial f}{\partial x_ix_j}$ is not constrained. For the proof we will need the following standard fact:

\begin{mylemma}[Gaussian integration by parts]
\label{lemma-5.2}
Let~$X$ be a centered Gaussian vector on~$\R^n$ and suppose $f\in C^1(\R^n)$ is such that~$f,\nabla f$ have subgaussian growth. Then for each~$i=1,\dots,n$,
\begin{equation}
\label{E:5.3nw}
\Cov\bigl(f(X),\,X_i\bigr)=\sum_{j=1}^n \Cov(X_i,X_j)E\Bigl(\frac{\partial f}{\partial x_j}(X)\Bigr).
\end{equation}
\end{mylemma}

\begin{myexercise}
Prove Lemma~\ref{lemma-5.2}. Hint: The proof is an actual integration by parts for~$X$ one-dimensional. For the general case use the positive definiteness of the covariance to find an $n\times n$ matrix~$A$ such that $X=AZ$ for~$Z$ i.i.d.\ $\NN(0,1)$. Then apply the one-dimensional results to each coordinate of~$Z$. 
\end{myexercise}

Gaussian integration can also be proved by approximating~$f$ by polynomials and invoking the following identity:

\begin{myexercise}[Wick pairing formula]
Let~$(X_1,\dots,X_{2n})$ be a centered multivariate Gaussian (with some variables possibly repeating). Show that
\begin{equation}
E(X_1\dots X_{2n})=\sum_{\pi\colon\text{\rm\ pairing}}\,\,\prod_{i=1}^n \Cov\bigl(X_{\pi_1(i)}X_{\pi_2(i)}\bigr),
\end{equation}
where a ``pairing'' is a pair~$\pi=(\pi_1,\pi_2)$ of functions $\pi_1,\pi_2\colon\{1,\dots,n\}\to\{1,\dots,2n\}$ such that
\begin{equation}
\pi_1(i)<\pi_2(i),\quad i=1,\dots,n,
\end{equation}
\begin{equation}
\pi_1(1)<\pi_1(2)\dots<\pi_1(n)
\end{equation}
and
\begin{equation}
\{1,\dots,2n\}= \bigcup_{i=1}^n\bigl\{\pi_1(i),\pi_2(i)\bigr\}\,.
\end{equation}
(Note that these force~$\pi_1(1)=1$.)
\end{myexercise}

\noindent
The pairing formula plays an important role in computations involving Gaussian fields; in fact, it is the basis of perturbative calculations of functionals of Gaussian processes and their organization in terms of Feynman diagrams.

\begin{proofsect}{Proof of Proposition~\ref{thm-5.1}}
Suppose that~$X$ and~$Y$ are realized on the same probability space so that $X\independent Y$. Consider the Gaussian interpolation
\begin{equation}
Z_t:=\sqrt{1-t^2}\,X+t Y,\qquad t\in[0,1].
\end{equation}
Then $Z_0=X$ and~$Z_1=Y$ and so
\begin{equation}
\label{E:5.5}
E f(Y)-Ef(X) = \int_0^1\frac{\textd}{\textd t}E f(Z_t)\,\textd t.
\end{equation}
Using elementary calculus along with Lemma~\ref{lemma-5.2},
\begin{equation}
\begin{aligned}
\frac{\textd}{\textd t}E f(Z_t)
&=\sum_{i=1}^nE\biggl(\Bigl(\frac{-t}{\sqrt{1-t^2}}X_i+Y_i\Bigr) \frac{\partial f}{\partial x_i}(Z_t)\biggr)
\\
&=t\sum_{i,j=1}^nE\biggl(\bigl[E(Y_iY_j)-E(X_iX_j)\bigr] \frac{\partial^2 f}{\partial x_i\partial x_j}(Z_t)\biggr).
\end{aligned}
\end{equation}
Based on our assumptions, the expression under the expectation is non-negative for every realization of~$Z_t$. Using this in \eqref{E:5.5} yields the claim.
\end{proofsect}

\section{Kahane's theory of Gaussian Multiplicative Chaos}
\noindent
We will find Theorem~\ref{thm-5.1} useful later but Kahane's specific interest in Gaussian Multiplicative Chaos actually required a version that is not directly obtained from the one above. Let us recall the setting more closely.

Let~$D\subset\R^d$ be a bounded open set and let~$\nu$ be a finite Borel measure on~$D$. Assume that~$C\colon D\times D\to\R\cup\{\infty\}$ is a symmetric, positive semi-definite kernel in~$L^2(\nu)$ which means that
\begin{equation}
\int_{D\times D}\nu(\textd x)\otimes\nu(\textd y)\,C(x,y)f(y)f(x)\ge0
\end{equation}
holds for every bounded measurable~$f\colon D\to\R$. If~$C$ is finite everywhere, then one can define a Gaussian process $\varphi=\NN(0,C)$. Our interest is, however, in the situation when~$C$ is allowed to diverge on the diagonal $\{(x,x)\colon x\in D\}\subset D\times D$ which means that the Gaussian process exists only in a generalized sense --- e.g., as a random distribution on a suitable space of test functions. 

We will not try to specify the conditions on~$C$ that would make this setting fully meaningful; instead, we will just assume that~$C$ can be written as
\begin{equation}
C(x,y)=\sum_{k=1}^\infty C_k(x,y),\qquad x,y\in D,
\end{equation}
where~$C_k$ is a continuous (and thus finite) covariance kernel for each~$k$ and the sum converges pointwise everywhere (including, possibly, to infinity when~$x=y$). We then consider the Gaussian processes 
\begin{equation}
\varphi_k=\NN(0,C_k)\quad\text{with}\quad\{\varphi_k\colon k\ge1\}\,\,\text{independent}\,.
\end{equation}
Letting
\begin{equation}
\label{E:5.12}
\Phi_n(x):=\sum_{k=1}^n\varphi_k(x)
\end{equation}
we define
\begin{equation}
\mu_n(\textd x):=\texte^{\Phi_n(x)-\frac12\Var[\Phi_n(x)]}\nu(\textd x).
\end{equation}
Lemma~\ref{lemma-GMC} (or rather its proof) gives the existence of a random Borel measure~$\mu_\infty$ such that
for each~$A\subset D$ Borel,
\begin{equation}
\label{E:5.14}
\mu_n(A)\,\underset{n\to\infty}\longrightarrow\,\,\mu_\infty(A)\quad\text{a.s.}
\end{equation}
As the covariances $\Cov(\Phi_n(x),\Phi_n(y))$ converge to~$C(x,y)$, we take $\mu_\infty$ as our interpretation of the measure
\begin{equation}
\text{``}\,\,\texte^{\Phi_\infty(x)-\frac12\Var[\Phi_\infty(x)]}\nu(\textd x)\,\,\text{''}
\end{equation}
for~$\Phi_\infty$ being the centered generalized Gaussian field with covariance~$C$. A key problem that Kahane had to deal with was the dependence of the limit measure on the above construction, and the uniqueness of the law of~$\mu_\infty$ in general. This is, at least partially, resolved in:

\begin{mytheorem}[Kahane's Uniqueness Theorem]
\label{thm-5.5}
For~$D\subset\R^d$ bounded and open, suppose there are covariance kernels $C_k,\wt C_k\colon D\times D\to\R$ such that
\settowidth{\leftmargini}{(11)}
\begin{enumerate}
\item[(1)] both $C_k$ and~$\wt C_k$ is continuous and non-negative everywhere on~$D\times D$,
\item[(2)] for each~$x,y\in D$,
\begin{equation}
\label{E:5.16}
\sum_{k=1}^\infty C_k(x,y) = \sum_{k=1}^\infty \wt C_k(x,y)
\end{equation}
with both sums possibly simultaneously infinite, and
\item[(3)] 
the fields $\varphi_k=\NN(0,C_k)$ and~$\wt\varphi_k=\NN(0,\wt C_k)$, with $\{\varphi_k,\wt\varphi_k\colon k\ge1\}$ all independent of one-another, have versions with continuous paths for each~$k\ge1$.
\end{enumerate}
Define, via \twoeqref{E:5.12}{E:5.14}, the random measures~$\mu_\infty$ and~$\wt \mu_\infty$ associated with these fields. Then
\begin{equation}
\label{E:5.18uai}
\mu_\infty(\textd x)\,\,\laweq\,\,\wt\mu_\infty(\textd x).
\end{equation}
\end{mytheorem}

In order to prove Theorem~\ref{thm-5.5} we will need the following variation on Proposition~\ref{thm-5.1} which lies at the heart of Kahane's theory:

\begin{myproposition}[Kahane's convexity inequality]
\label{prop-5.6}
Let~$D\subset\R^n$ be bounded and open and let~$\nu$ be a finite Borel measure on~$D$. Let~$C,\wt C\colon D\times D\to\R$ be covariance kernels such that~$\varphi=\NN(0,C)$ and~$\wt\varphi=\NN(0,\wt C)$ have continuous paths a.s. If
\begin{equation}
\label{E:5.18}
\wt C(x,y)\ge C(x,y),\quad x,y\in D,
\end{equation}
then for each convex~$f\colon[0,\infty)\to\R$ with at most polynomial growth at infinity,
\begin{equation}
E\,f\Bigl(\,\int_D\texte^{\wt\varphi(x)-\frac12\Var[\wt\varphi(x)]}\nu(\textd x)\Bigr)
\ge E\,f\Bigl(\,\int_D\texte^{\varphi(x)-\frac12\Var[\varphi(x)]}\nu(\textd x)\Bigr).
\end{equation}
\end{myproposition}

\begin{proofsect}{Proof}
By approximation we may assume that~$f\in C^2(\R)$ (still convex). By the assumption of the continuity of the fields, it suffices to prove this for~$\nu$ being the sum of a finite number of point masses, $\nu=\sum_{i=1}^n p_i\delta_{x_i}$ where~$p_i>0$. (The general case then follows by the weak limit of such measures to~$\nu$.)

Assume that the fields~$\varphi$ and~$\wt\varphi$ are realized on the same probability space so that $\varphi\independent\wt\varphi$. Consider the interpolated field
\begin{equation}
\varphi_t(x):=\sqrt{1-t^2}\,\varphi(x)+t\wt\varphi(x),\qquad t\in[0,1].
\end{equation}
Since~$\varphi_0(x)=\varphi(x)$ and~$\varphi_1(x)=\wt\varphi(x)$, it suffices to show
\begin{equation}
\label{E:5.21}
\frac{\textd}{\textd t}
E\,f\Bigl(\,\sum_{i=1}^np_i\texte^{\varphi_t(x_i)-\frac12\Var[\varphi_t(x_i)]}\Bigr)\ge0\,.
\end{equation}
For this we abbreviate $W_t(x):=\texte^{\varphi_t(x)-\frac12\Var[\varphi_t(x)]}$ and use elementary calculus to get
\begin{multline}
\quad
\frac{\textd}{\textd t} E\,f\Bigl(\,\sum_{i=1}^np_iW_t(x_i)\Bigr)
=\sum_{i=1}^n p_iE\biggl(\Bigl[-\frac{t}{\sqrt{1-t^2}}\varphi(x_i)+\wt\varphi(x_i)
\\
+t\Var\bigl(\varphi(x_i)\bigr)-t\Var\bigl(\wt\varphi(x_i)\bigr)\Bigr] W_t(x_i)f'(\cdots)\biggr)
\end{multline}
Next we integrate by parts (cf Lemma~\ref{lemma-5.2}) the terms involving~$\varphi(x_i)$, which results in the~$\varphi(x_j)$-derivative of~$W_t(x_i)$ or~$f(\cdots)$. A similar process is applied to the term~$\wt\varphi(x_i)$. A key point is that the contribution from differentiating~$W_t(x_i)$ exactly cancels that coming from the variances. Hence we get
\begin{multline}
\quad
\frac{\textd}{\textd t} E\,f\Bigl(\,\sum_{i=1}^np_iW_t(x_i)\Bigr)
\\
=\sum_{i,j=1}^n p_ip_j\bigl[\wt C(x_i,x_j)-C(x_i,x_j)\bigr]\,E\Bigl(W_t(x_i)W_t(x_j)f''(\cdots)\Bigr).
\quad
\end{multline}
As~$f''\ge0$ by assumption and $W_t(x)\ge0$ and~$p_i,p_j\ge0$ by inspection, \eqref{E:5.18} indeed implies \eqref{E:5.21}. The claim follows by integration over~$t$.
\end{proofsect}

This permits us to give:

\begin{proofsect}{Proof of Theorem~\ref{thm-5.5}}
The claim \eqref{E:5.18uai} will follow once we show
\begin{equation}
\label{E:5.24}
\int g(x)\mu_\infty(\textd x)\,\,\laweq\,\,\int g(x)\wt\mu_\infty(\textd x)
\end{equation}
for any continuous function $g\colon D\to[0,\infty)$ supported in a compact set $A\subset D$. Let~$\{C_k\colon k\ge1\}$ and~$\{\wt C_k\colon k\ge1\}$ be the covariances in the statement. We now invoke a reasoning underlying the proof of Dini's Theorem: For each~$\epsilon>0$ and each~$n\in\N$, there is~$m\in\N$ such that
\begin{equation}
\label{E:5.25}
\sum_{k=1}^n C(x,y)<\epsilon+\sum_{k=1}^m\wt C(x,y),\qquad x,y\in A.
\end{equation}
Indeed, fix $n\in\N$  and let~$F_m$ be the set of pairs~$(x,y)\in A\times A$ where \eqref{E:5.25} fails. Then~$F_m$ is closed (and thus compact) by the continuity of the covariances; their non-negativity in turn shows that $m\mapsto F_m$ is decreasing with respect to set inclusion. The equality \eqref{E:5.16} translates into $\bigcap_{m\ge1}F_m=\emptyset$ and so, by Heine-Borel, we must have~$F_m=\emptyset$ for~$m$ large enough thus giving us \eqref{E:5.25}.

Interpreting the~$\epsilon$ term on the right-hand side of \eqref{E:5.25} as the variance of the random variable~$Z_\epsilon=\NN(0,\epsilon)$ that is independent of~$\wt\varphi$, Proposition~\ref{prop-5.6} with the choice~$f(x):=\texte^{-\lambda x}$ for some~$\lambda\ge0$ gives us
\begin{equation}
E\bigl(\texte^{-\lambda\,\texte^{Z_\epsilon-\epsilon/2}\,\int g\,\textd\wt\mu_m}\bigr)
\ge E\bigl(\texte^{-\lambda\int g\,\textd\mu_n}\bigr)\,. 
\end{equation}
Invoking the limit \eqref{E:5.14} and taking~$\epsilon\downarrow0$ afterwards yields
\begin{equation}
\label{E:5.29nw}
E\bigl(\texte^{-\lambda\int g\,\textd\wt\mu_\infty}\bigr)
\ge E\bigl(\texte^{-\lambda\int g\,\textd\mu_\infty}\bigr)\,.
\end{equation}
Swapping~$C$ and~$\wt C$ shows that equality holds in \eqref{E:5.29nw} and since this is true for every~$\lambda\ge0$, we get \eqref{E:5.24} as desired.
\end{proofsect}

As far as the GMC associated with the two-dimensional Gaussian Free Field is concerned, Theorem~\ref{thm-5.5} shows that any decomposition of the continuum Green function $\wh G^D(x,y)$ into the sum of positive covariance kernels will yield, through the construction in Lemma~\ref{lemma-GMC}, the same limiting measure~$\mu_\infty^{D,\beta}$. One example of such a decomposition is that induced by the Gibbs-Markov property upon reductions to a subdomain; this is the content of Exercises~\ref{ex:4.15}--\ref{ex:Z-mu}. Another example is the \emph{white noise decomposition} that will be explained in Section~\ref{sec-10.5}.

We remark that (as noted above) uniqueness of the GMC measure has now been proved in a completely general setting  by Shamov~\cite{Shamov}. Still, Kahane's approach is sufficient to prove conformal transformation rule for the $Z^D_\lambda$-measures discussed earlier in these notes:

\begin{myexercise}
\label{ex:Kahane-CFI}
Prove that $\{Z^D_\lambda\colon D\in\mathfrak D\}$ obey the conformal transformation rule stated in Theorem~\ref{thm-2.4}. Unlike Exercise~\ref{ex:4.18}, do not assume the full uniqueness stated in Theorem~\ref{thm-GMC-unique}.
\end{myexercise}

\section{Comparisons for the maximum}
\noindent
Our next task is to use Kahane's inequality from Theorem~\ref{thm-5.1} to provide comparisons between the maxima of two Gaussian vectors with point-wise ordered covariances. We begin with a corollary to Theorem~\ref{thm-5.1}:

\begin{mycorollary}
\label{cor-5.8}
Suppose that~$X$ and~$Y$ are centered Gaussians on~$\R^n$ such that
\begin{equation}
\label{E:5.28}
E(X_i^2)=E(Y_i^2),\qquad i=1,\dots,n
\end{equation}
and
\begin{equation}
\label{E:5.29}
E(X_iX_j)\le E(Y_iY_j),\qquad i,j=1,\dots,n.
\end{equation}
Then for any~$t_1,\dots,t_n\in\R$,
\begin{equation}
\label{E:5.31uai}
P\bigl(X_i\le t_i\colon i=1,\dots,n\bigr)\le P\bigl(Y_i\le t_i\colon i=1,\dots,n\bigr).
\end{equation}
\end{mycorollary}

\begin{proofsect}{Proof}
Consider any collection~$g_1,\dots,g_n\colon\R\to\R$ of non-negative bounded functions that are smooth and non-increasing. Define
\begin{equation}
f(x_1,\dots,x_n):=\prod_{i=1}^ng_i(x_i).
\end{equation}
Then $\frac{\partial^2f}{\partial x_i\partial x_j}\ge0$ for each~$i\ne j$. Hence, by Theorem~\ref{thm-5.1}, conditions \twoeqref{E:5.28}{E:5.29} imply $E f(Y)\ge E f(X)$. The claim follows by letting~$g_i$ decrease to~$1_{(-\infty,t_i]}$.
\end{proofsect}

From here we now immediately get:

\begin{mycorollary}[Slepian's lemma]
Suppose~$X$ and~$Y$ are centered Gaussians on~$\R^n$ with
\begin{equation}
\label{E:5.32}
E(X_i^2)=E(Y_i^2),\qquad i=1,\dots,n
\end{equation}
and
\begin{equation}
\label{E:5.33}
E\bigl((X_i-X_j)^2\bigr)\le E\bigl((Y_i-Y_j)^2\bigr),\qquad i,j=1,\dots,n.
\end{equation}
Then for each~$t\in\R$,
\begin{equation}
\label{E:5.35ua}
P\Bigl(\,\max_{i=1,\dots,n}X_i>t\Bigr)\le P\Bigl(\,\max_{i=1,\dots,n}Y_i>t\Bigr).
\end{equation}
\end{mycorollary}

\begin{proofsect}{Proof}
Set~$t_1=\dots=t_n:=t$ in the previous corollary.
\end{proofsect}

Slepian's lemma (proved originally in~\cite{Slepian}) has a nice verbal formulation using the following concept: Given a Gaussian process~$\{X_t\colon t\in T\}$ on a set~$T$,
\begin{equation}
\label{E:5.36ua}
\rho_X(t,s):=\sqrt{E\bigl((X_t-X_s)^2\bigr)}
\end{equation}
defines a pseudometric on~$T$. Indeed, we pose:

\begin{myexercise}
Verify that~$\rho_X$ is indeed a pseudo-metric on~$T$.
\end{myexercise}

\noindent
Disregarding the prefix ``pseudo'', we will call $\rho_X$ the \emph{canonical}, or \emph{intrinsic}, metric associated with Gaussian processes. Slepian's lemma may then be verbalized as follows: \textit{For two Gaussian processes with equal variances, the one with larger intrinsic distances has a stochastically larger maximum.}

The requirement of equal variances is often too much to ask for. One way to compensate for an inequality there is by adding suitable independent Gaussians to~$X$ and~$Y$. However, it turns out that this inconvenience disappears altogether if we contend ourselves with the comparison of expectations only (which is, by way of integration, implied by \eqref{E:5.35ua}):

\begin{myproposition}[Sudakov-Fernique inequality]
\label{prop-SF}
Suppose that~$X$ and~$Y$ are centered Gaussians in~$\R^n$ such that
\begin{equation}
\label{E:5.36}
E\bigl((X_i-X_j)^2\bigr)\le E\bigl((Y_i-Y_j)^2\bigr),\qquad i,j=1,\dots,n.
\end{equation}
Then
\begin{equation}
\label{E:5.37}
E\bigl(\,\max_{i=1,\dots,n}X_i\bigr)\le E\bigl(\,\max_{i=1,\dots,n}Y_i\bigr)
\end{equation}
\end{myproposition}

\begin{proofsect}{Proof}
Consider the function
\begin{equation}
\label{E:f-beta}
f_\beta(x_1,\dots,x_n):=\frac1\beta\log\Bigl(\,\sum_{i=1}^n\texte^{\beta x_i}\Bigr).
\end{equation}
For readers familiar with statistical mechanics, $f_\beta$ can be thought of as a free energy. H\"older's inequality implies that~$x\mapsto  f_\beta(x)$ is convex. In addition, we also get
\begin{equation}
\lim_{\beta\to\infty}f_\beta(x)=\max_{i=1,\dots,n}x_i.
\end{equation}
Using Dominated Convergence, it therefore suffices to show that
\begin{equation}
\label{E:5.40}
E f_\beta(X)\le E f_\beta(Y),\qquad \beta\in(0,\infty).
\end{equation}
The  proof of this inequality will be based on a re-run of the proof of Kahane's inequality. Assuming again~$X\independent Y$ and letting $Z_t:=\sqrt{1-t^2}\,X+tY$, differentiation yields
\begin{equation}
\label{E:5.41}
\frac{\textd}{\textd t}E f_\beta(Z_t)=t\sum_{i,j=1}^n E\Bigl(\, \bigl[E(Y_iY_j)-E(X_iX_j)\bigr]\frac{\partial^2 f_\beta}{\partial x_i\partial x_j}(Z_t)\Bigr).
\end{equation}
Now
\begin{equation}
\frac{\partial f_\beta}{\partial x_i}=\frac{\texte^{\beta x_i}}{\sum_{j=1}^n\texte^{\beta x_j}}=:p_i(x)
\end{equation}
where~$p_i\ge0$ with~$\sum_{i=1}^np_i(x)=1$. For the second derivatives we get
\begin{equation}
\frac{\partial^2 f_\beta}{\partial x_i\partial x_j}=\beta\bigl[p_i(x)\delta_{ij}-p_i(x)p_j(x)\bigr]\,.
\end{equation}
Plugging this on the right of \eqref{E:5.41} (and omitting the argument~$Z_t$ of the second derivative as well as the~$p_i$'s) we then observe
\begin{equation}
\begin{aligned}
\sum_{i,j=1}^n\bigl[E(Y_i&Y_j)-E(X_iX_j)\bigr]\frac{\partial^2 f_\beta}{\partial x_i\partial x_j}
\\
&=\beta\sum_{i,j=1}^n\bigl[E(Y_iY_j)-E(X_iX_j)\bigr]\bigl[p_i\delta_{ij}-p_ip_j\bigr]
\\
&=\beta\sum_{i,j=1}^n \bigl[E(Y_i^2)+E(X_i^2)\bigr]p_ip_j+\beta\sum_{i,j=1}^n\bigl[E(Y_iY_j)-E(X_iX_j)\bigr]p_ip_j
\\
&=\frac12\beta\sum_{i,j=1}^n\Bigl[E\bigl((Y_i-Y_j)^2\bigr)-E\bigl((X_i-X_j)^2\bigr)\Bigr]p_ip_j,
\end{aligned}
\end{equation}
where we used that $\{p_i\colon i=1,\dots,n\}$ are probabilities in the second line and then symmetrized the first sum under the exchange of~$i$ for~$j$ to wrap the various parts into the expression in the third line.
Invoking \eqref{E:5.36}, this is non-negative (pointwise) and so we get \eqref{E:5.40} by integration. The claim follows.
\end{proofsect}

The Sudakov-Fernique inequality may be verbalized as follows: \textit{For two Gaussian processes, the one with larger intrinsic distances has a larger expected maximum.}
Here are some other, rather elementary, facts related to the same setting:

\begin{myexercise}
\label{ex:max-E-nonzero}
Show that, for any centered Gaussians~$X_1,\dots,X_n$,
\begin{equation}
E\bigl(\,\max_{i=1,\dots,n}X_i\bigr)\ge0.
\end{equation}
Prove that equality occurs if and only if $X_i=X_1$ for all~$i=1,\dots,n$ a.s.
\end{myexercise}

\begin{myexercise}
Suppose that~$X$, resp.,~$Y$ are centered Gaussian vectors on~$\R^n$ with covariances~$C$, resp.,~$\wt C$. Show that if~$\wt C-C$ is positive semi-definite, then \eqref{E:5.37} holds.
\end{myexercise}

We will use the Sudakov-Fernique inequality a number of times in these notes. Unfortunately, the comparison between the expected maxima is often insufficient for the required level of precision; indeed, \eqref{E:5.37} only tells us that~$\max_iX_i$ is not larger than a \emph{large multiple} of $E\max_i Y_i$ with significant probability. In Lemma~\ref{lemma-SF-extension}, which is a kind of cross-breed between the Slepian Lemma and the Sudakov-Fernique inequality, we will show how this can be boosted to control the upper tail of $\max_i X_i$ by that of~$\max_i Y_i$, assuming that the latter maximum is concentrated. 

\section{Stochastic domination and FKG inequality}
\noindent
Yet another set of convenient inequalities that we will use on occasion is related to the (now) classical notions of stochastic domination. This topic has been covered in the concurrent class by Hugo Duminil-Copin~\cite{Hugo-notes}, but we will still review its salient features for future (and independent) reference. 

We have already invoked the following version of stochastic domination: Given real-valued random variables~$X$ and~$Y$, we say that~\emph{$X$ is stochastically larger than~$Y$} or, equivalently, that \emph{$X$ stochastically dominates~$Y$} if the cumulative distribution function of~$Y$ exceeds that of~$X$ at all points, i.e.,
\begin{equation}
P(X\le t)\le P(Y\le t),\quad t\in\R.
\end{equation}
 Although this may seem just a statement about the laws of these random variables, there is a way to realize the ``domination'' pointwise:

\begin{myexercise}
\label{ex:5.14}
Suppose that~$X$ stochastically dominates~$Y$. Prove that there is a coupling of these random variables (i.e., a realization of both of them on the same probability space) such that $P(Y\le X)=1$.
\end{myexercise}

The complete order of the real line plays a crucial role here. A natural question is what to do about random variables that take values in a set which admits only a (non-strict) \emph{partial order}; i.e., a reflexive, antisymmetric and transitive binary relation~$\preccurlyeq$. An instance of this is the product space $\R^A$ where \begin{equation}
\label{E:5.48}
x\preccurlyeq y\text{ means that }x_i\le y_i\text{ for every }i\in A.
\end{equation}
In the more general context, we will rely on the following notion: A real-valued function~$f$ is called \emph{increasing} if~$x\preccurlyeq y$ implies~$f(x)\le f(y)$. The same concept applies for Cartesian products of arbitrary ordered spaces (not just~$\R$).
 
\begin{mydefinition}[Stochastic domination]
We say that~$X$ stochastically dominates~$Y$, writing $Y\preccurlyeq X$, if $E f(Y)\le E f(X)$ holds for every bounded measurable increasing~$f$.
\end{mydefinition}

It is an easy exercise to check that, for~$X$ and~$Y$ real valued and~$\preccurlyeq$ being the usual order of the reals, this coincides with the definition given above. Exercise~\ref{ex:5.14} then turns into an instance of:

\begin{mytheorem}[Strassen's lemma]
Suppose~$X$ and~$Y$ take values in a par\-ti\-ally-ordered compact metric space $\scrX$ such that $\{(x,y)\colon x\preccurlyeq y\}$ is closed in $\scrX\times\scrX$. Then~$X$ stochastically dominates~$Y$ if and only if there is a coupling of these random variables such that~$P(Y\preccurlyeq X)=1$.
\end{mytheorem}

For a proof we refer to Liggett~\cite[Theorem~2.4]{Liggett-IPS}. Stochastic domination is often interpreted as a property of probability measures; i.e., we write $\mu\preccurlyeq\nu$ if the random variables~$X$ and~$Y$ with respective laws~$\nu$ and $\mu$ obey $Y\preccurlyeq X$. 

Stochastic domination is closely related to the concept of \emph{positive correlations}, also known as \emph{positive association} or \emph{(weak) FKG inequality}. Let us call an event~$A$ increasing if its indicator~$1_A$ is an increasing function, and decreasing if $1_A$ is a decreasing function.
The following concept inherits its acronym from the authors of Fortuin, Kasteleyn and Ginibre~\cite{FKG}:

\begin{mydefinition}[Positive correlations a.k.a.\ FKG inequality]
A probability measure~$\mu$ on a partially-order space is said to have positive correlations, or satisfy the (weak) FKG inequality, if
\begin{equation}
\label{E:5.49}
\mu(A\cap B)\ge \mu(A)\mu(B)
\end{equation}
holds for any pair of increasing events~$A$ and~$B$.
\end{mydefinition}

Positive correlations can be interpreted using the concept of stochastic domination as follows:
\begin{equation}
\forall B \text{ increasing with }\mu(B)>0\colon\quad
\mu\preccurlyeq\mu(\cdot|B).
\end{equation}
In other words, a probability measure has positive correlations if (and only if) conditioning on an increasing event ``increases'' the measure. 

As is readily checked, if \eqref{E:5.49} holds for pairs of increasing events then it also holds for pairs of decreasing events, and that the opposite inequality applies when one event is increasing and the other decreasing. A seemingly stronger, albeit equivalent, formulation of positive correlations is via increasing functions:

\begin{myexercise}
\label{ex:5.18}
Show that~$\mu$ has positive correlations if and only if
\begin{equation}
\label{E:5.51}
E_\mu(\,fg)\ge E_\mu(f)E_\mu(g)
\end{equation}
holds true for any pair of increasing (or decreasing) functions~$f,g\in L^2(\mu)$. Hint: Write every such~$f$ as a limit of linear combinations of increasing indicators.
\end{myexercise}

That positive correlations and stochastic domination are related is seen from:

\begin{myexercise}
Let~$\mu$ and~$\nu$ be probability measures on a partially-ordered space. If~$\mu$ and~$\nu$ have positive correlations and $\mu\preccurlyeq\nu$ (or~$\nu\preccurlyeq\mu$), then their convex combination $t\mu+(1-t)\nu$ has positive correlations for all~$t\in[0,1]$.
\end{myexercise}

To see how the above concepts are related, let us recall the situation of independent random variables where these connections were observed first:

\begin{mylemma}[Harris' inequality]
\label{lemma-Harris}
For any set~$A$, any product law on~$\R^A$ (endowed with the product $\sigma$-algebra and the partial order \eqref{E:5.48}) has positive correlations.
\end{mylemma}

\begin{proofsect}{Proof}
Let~$\mu$ be a product law which (by the assumed product structure) we may think of as the distribution of independent random variables~$\{X_i\colon i\in A\}$. We first prove \eqref{E:5.51} for~$f,g\in L^2(\mu)$ that depend only on one of these random variables, say~$X_i$. Let~$\wt X_i$ be an independent copy of~$X_i$. If $f,g$ are increasing, then
\begin{equation}
\bigl[f(X_i)-f(\wt X_i)\bigr]\bigl[g(X_i)-g(\wt X_i)\bigr]\ge0.
\end{equation}
Taking expectation then yields \eqref{E:5.51}.

Assuming \eqref{E:5.51} holds for~$f,g\in L^2(\mu)$ that depend on random variables $X_1,\dots,X_k$, we will now show that it holds for any~$f,g\in L^2(\mu)$ that depend on~$X_1,\dots,X_{k+1}$. Denote
\begin{equation}
\FF_k:=\sigma(X_1,\dots,X_k),\quad f_k:=E(f|\FF_k)\quad\text{and}\quad g_k:=E(g|\FF_k)
\end{equation}
and write $\mu_k$ for the regular conditional probability $\mu(\cdot|\FF_k)$. Then \eqref{E:5.51} for the one-dimensional case yields $E_{\mu_k}(fg)\ge f_kg_k$. Moreover, thanks to the product structure and Fubini-Tonelli, $f_k,g_k$ are both increasing. They are also functions of~$X_1,\dots,X_k$ only and so the induction assumption shows
\begin{equation}
\begin{aligned}
E_\mu(fg)&=E_\mu\bigl(E_{\mu_k}(fg)\bigr)
\\
&\ge E_\mu\bigl(f_kg_k\bigr)\ge E_\mu(f_k)E_\mu(g_k)=E_\mu(f)E_\mu(g).
\end{aligned}
\end{equation}
We conclude that \eqref{E:5.51} applies to all~$f,g\in L^2(\mu)$ depending only on a finite number of coordinates. 

Now take a general~$f\in L^2(\mu)$. By elementary measurability considerations,~$f$ is a function of at most a countable (sub)collection $\{X_1,X_2,\dots\}$ of the above random variables; Levy's Forward Theorem ensures $E_\mu(f|\FF_k)\to f$ in~$L^2(\mu)$ as~$k\to\infty$. Since $E_\mu(f|\FF_k)$ is also increasing, \eqref{E:5.51} for any~$f,g\in L^2(\mu)$ follows from the finite-dimensional case by usual approximation arguments.
\end{proofsect}

The above proof actually showed more than \eqref{E:5.51}; namely, that any product law has the following property:

\begin{mydefinition}[strong-FKG property]
We say that a probability measure on~$\R^A$ is strong-FKG if the conditional law given the values for any finite number of coordinates has positive correlations.
\end{mydefinition}

We remark that the expression ``positive correlations'' is sometimes used in the context when~$\mu$ is a law of random variables~$\{X_i\colon i\in I\}$ and ``positivity of correlations'' refers to $\Cov(X_i,X_j)\ge0$ --- namely, a special case of \eqref{E:5.51} with $f(x):=x_i$ and~$g(x):=x_j$. This is the reason why the term ``positive association'' is generally preferred to capture the full strength of \eqref{E:5.51}. Notwithstanding, this is all the same for Gaussian random variables:

\begin{myproposition}[strong-FKG for Gaussians]
\label{prop-5.24}
Suppose that~$\mu$ is the law of a Gaussian vector~$X$ on~$\R^n$. Then
\begin{equation}
\text{$\mu$ is strong-FKG}\quad\Leftrightarrow\quad \Cov(X_i,X_j)\ge0,\quad i,j=1,\dots,n.
\end{equation}
\end{myproposition}

\begin{proofsect}{Proof}
To get $\Rightarrow$ we just use the definition of (weak) FKG along with the fact that~$f(X):=X_i$ is increasing. Moving to the~$\Leftarrow$ part, assume $\Cov(X_i,X_j)\ge0$ for all~$i,j=1,\dots,n$. Conditioning a multivariate Gaussian on part of the variables preserves the multivariate Gaussian structure as well as the covariances. It thus suffices to prove that~$\mu$ satisfies the weak FKG inequality for which, by Exercise~\ref{ex:5.18} and routine approximation arguments, it suffices to show $\Cov(f(X),g(X))\ge0$ for any non-decreasing smooth functions~$f,g\colon\R^n\to\R$ with bounded gradients. This follows from an enhanced version of Gaussian integration by parts in Lemma~\ref{lemma-6.2} (to be proved soon) and the fact that the first partial derivatives of~$f$ and~$g$ are non-negative.
\end{proofsect}

We note that, since the DGFF covariance is given by the Green function which is non-negative everywhere, Proposition~\ref{prop-5.24} shows that the DGFF is a strong-FKG process.

We close this lecture by noting that the above discussion of stochastic domination focused only on the topics that are needed for a full understanding of the arguments carried out in these notes. The reader is referred to, e.g., Liggett~\cite{Liggett-IPS} or Grimmett~\cite{Grimmett-RCM} for a comprehensive treatment of this subject including its (interesting) history.


\chapter{Concentration techniques}
\label{lec-6}
\noindent
In this lecture we will establish bounds on the maximum of Gaussian random variables that are not be based on comparisons but rather on the metric properties of the covariance kernel. The first result to be proved here is the Borell-Tsirelson-Ibragimov-Sudakov inequality on the concentration of the maximum. Any use of this inequality will inevitably entail estimates on the expected maximum which we do via the Fernique majorization technique. Once these are stated and proved, we will infer some standard but useful consequences concerning boundedness and continuity of centered Gaussian processes. The presentation draws on that in Ledoux~\cite{Ledoux}, Adler~\cite{Adler} and Adler and Taylor~\cite{Adler-Taylor}.

\section{Inheritance of Gaussian tails}
\noindent
Much of the present day probability hinges on the phenomenon of \emph{concentration of measure}. For Gaussian random variables this is actually a very classical subject. The following inequality will come up frequently in the sequel:

\begin{mytheorem}[Borell-TIS inequality]
\label{thm-BTIS}
Let~$X$ be a centered Gaussian on~$\R^n$ and set
\begin{equation}
\label{E:6.1}
\sigma_X^2:=\max_{i=1,\dots,n}E(X_i^2).
\end{equation}
Then for each~$t>0$,
\begin{equation}
P\Bigl(\bigl|\max_{i=1,\dots,n}X_i-E(\,\max_{i=1,\dots,n} X_i)\bigr|>t\Bigr)\le2\texte^{-\frac{t^2}{2\sigma_X^2}}.
\end{equation}
\end{mytheorem}

\noindent
This result may be verbalized as: \textit{The maximum of Gaussian random variables has a tail no heavier than the heaviest tail seen among these random variables}. Of course, the maximum is no longer centered (cf Exercise~\ref{ex:max-E-nonzero}) and so any use of this bound requires information on its expectation as well.

The original proof of Theorem~\ref{thm-BTIS} was given by Borell~\cite{Borell} using a Gaussian isoperimetric inequality; the inequality was discovered independently in the Eastern block by Tsirelson, Ibragimov and Sudakov~\cite{TIS}. We will  instead proceed using analytic techniques based on hypercontractivity. The following lemma offers a generalization of the formula on the Gaussian integration by parts:

\begin{mylemma}
\label{lemma-6.2}
Let~$X$ be a Gaussian vector on~$\R^n$ and let~$f,g\in C^1(\R^n)$ have subgaussian growth. Then
\begin{multline}
\label{E:6.3nw}
\quad
\Cov\bigl(\,f(X),g(X)\bigr)
\\
=\int_0^1\textd t \,\sum_{i,j=1}^n\Cov(X_i,X_j) E\biggl(\frac{\partial f}{\partial x_i}(X) \frac{\partial g}{\partial x_j}\bigl(\,tX+\sqrt{1-t^2}\, Y\bigr)\biggr),
\quad
\end{multline}
where~$Y\,\laweq\,X$ with~$Y\independent X$ on the right-hand side.
\end{mylemma}

\begin{proofsect}{Proof}
Since \eqref{E:6.3nw} is an equality between bilinear expressions for two functions of finitely-many variables, we may try to prove it by first checking it for a sufficiently large class of functions (e.g., the exponentials $x\mapsto\texte^{k\cdot x}$) and then using extension arguments. We will instead rely on Gaussian integration by parts.

For~$X$ and~$Y$ as above and $t\in[0,1]$, abbreviate
\begin{equation}
Z_t:=tX+\sqrt{1-t^2}\, Y.
\end{equation}
 Approximation arguments permit us to assume~$g\in C^2$ along with all the second partial derivatives having subgaussian growth. Then
\begin{equation}
\begin{aligned}
\Cov\bigl(\,f(X),g(X)\bigr)
&=E\Bigl( f(X)\bigl[g(Z_1)-g(Z_0)\bigr]\Bigr)
\\
&=\int_0^1\textd t\, E\Bigl(\,f(X)\frac{\textd}{\textd t}g(Z_t)\Bigr)
\\
&=\int_0^1\textd t\,\sum_{i=1}^n\,E\biggl(\Bigl[X_i-\frac{t}{\sqrt{1-t^2}}Y_i\Bigr]f(X)\frac{\partial g}{\partial x_i}(Z_t)\biggr).
\end{aligned}
\end{equation}
The integration by parts (cf Lemma~\ref{lemma-5.2}) will eliminate the square bracket and yield two contributions: one from the derivative hitting~$f$ and the other from the derivative hitting the partial derivative of~$g$. The latter term equals the sum over~$j$ of $t\Cov(X_i,X_j)-t\Cov(Y_i,Y_j)$ times $E[f(X)\frac{\partial^2g}{\partial x_i\partial x_j}(Z_t)]$. As~$Y\,\laweq\,X$, this term vanishes identically. The term where the derivative hits~$f$ produces the integrand in \eqref{E:6.3nw}.
\end{proofsect}

As a side note, we notice that this implies:

\begin{mycorollary}[Gaussian Poincar\'e inequality]
For~$X_1,\dots,X_n$ i.i.d.\ copies of~$\NN(0,1)$ and any~$f\in C^1(\R^n)$ with~$f,\nabla f\in L^2(\texte^{-|x|^2/2}\textd x)$,
\begin{equation}
\Var\bigl(f(X)\bigr)\le E\bigl(|\nabla f(X)|^2\bigr).
\end{equation}
\end{mycorollary}

\begin{proofsect}{Proof}
Apply Cauchy-Schwarz on the right-hand side of \eqref{E:6.3nw} while noting also that $tX+\sqrt{1-t^2}\,Y\,\,\laweq\,\, X$. An alternative is to use Gaussian integration by parts formula instead of \eqref{E:6.3nw}.
\end{proofsect}

Note that this bound is of dimension-less nature --- meaning: with no~$n$ dependence of the (implicit) constant on the right-hand side. This is quite in contrast to the Poincar\'e inequality on~$\R^d$. (A generalization to non-i.i.d.\ Gaussian vectors is straightfoward.)

Moving along with the proof of the Borell-TIS inequality, next we will prove:

\begin{mylemma}[Concentration for Lipschitz functions]
\label{lemma-6.4}
Let~$X_1,\dots,X_n$ be i.i.d.\ copies of~$\NN(0,1)$ and let~$f\colon\R^n\to\R$ be Lipschitz in the sense that, for some~$M\in(0,\infty)$,
\begin{equation}
\label{E:6.6nw}
\bigl|f(x)-f(y)\bigr|\le M|x-y|,\quad x,y\in\R^n,
\end{equation}
where~$|\cdot|$ on the right-hand side is the Euclidean norm. Then for each~$t>0$,
\begin{equation}
\label{E:6.7}
P\bigl(\,f(X)-Ef(X)>t\bigr)\le\texte^{-\frac{t^2}{2M^2}}.
\end{equation}
\end{mylemma}

\begin{proofsect}{Proof}
By approximation we may assume that~$f\in C^1$ with~$\nabla f$ having the Euclidean norm at most~$M$. By adding a suitable constant to~$f$ we may assume~$E f(X)=0$. The exponential Chebyshev inequality then shows
\begin{equation}
\label{E:6.8}
P\bigl(\,f(X)-Ef(X)>t\bigr)
\le\,\texte^{-\lambda t}E\bigl(\texte^{\lambda f(X)}\bigr)
\end{equation}
for any~$\lambda\ge0$ and so just we need to bound the expectation on the right.

Here we note that Lemma~\ref{lemma-6.2} with~$g(x):=\texte^{\lambda f(x)}$ and \eqref{E:6.6nw} imply
\begin{multline}
\quad
E\bigl(\,f(X)\texte^{\lambda f(X)}\bigr)
\\
=\int_0^1\textd t\,\lambda E\Bigl(\nabla f(X)\cdot\nabla f(Z_t)\texte^{\lambda f(Z_t)}\Bigr)
\overset{\lambda\ge0}\le \lambda M^2 E\bigl(\texte^{\lambda f(X)}\bigr).
\end{multline}
The left-hand side is the derivative of the expectation on the right-hand side. It follows that the function
\begin{equation}
h(\lambda):=E\bigl(\texte^{\lambda f(X)}\bigr)
\end{equation}
obeys the differential inequality
\begin{equation}
h'(\lambda)\le \lambda M^2 h(\lambda),\qquad \lambda\ge0.
\end{equation}
As~$h(0)=1$, this is readily solved to give
\begin{equation}
E\bigl(\texte^{\lambda f(X)}\bigr)\le\texte^{\frac12\lambda^2M^{2}}.
\end{equation}
Inserting this into \eqref{E:6.8} and optimizing over~$\lambda\ge0$ then yields the claim.
\end{proofsect}

In order to prove the Borell-TIS inequality, we will also need:

\begin{myexercise}
\label{ex:6.5}
Denote $f(x):=\max_{i=1,\dots,n}x_i$. Prove that for any $n\times n$-matrix~$A$,
\begin{equation}
\label{E:6.13}
\bigl|f(Ax)-f(Ay)\bigr|\le \sqrt{\max_{i=1,\dots,n}(A^{\text{T}}A)_{ii}}\,\,
|x-y|,\quad x,y\in\R^n,
\end{equation}
with $|x-y|$ denoting the Euclidean norm of~$x-y$ on the right-hand side. 
\end{myexercise}

\begin{proofsect}{Proof of Theorem~\ref{thm-BTIS}}
Let~$X$ be the centered Gaussian on~$\R^n$ from the statement and let~$C$ denote its covariance matrix. In light of the symmetry and positive semi-definiteness of~$C$, there is an $n\times n$-matrix~$A$ such that $C=A^{\text{T}}A$. If~$Z=(Z_1,\dots,Z_n)$ are i.i.d.\ copies of~$\NN(0,1)$, then
\begin{equation}
\label{E:6.14}
X\,\,\laweq\,\,AZ.
\end{equation}
Denoting~$f(x):=\max_{i=1,\dots,n}x_i$, Exercise~\ref{ex:6.5} shows that $x\mapsto f(Ax)$ is Lipschitz with Lipschitz constant~$\sigma_X$. The claim follows from \eqref{E:6.7} and a union bound.
\end{proofsect}

For a future reference, note that using \eqref{E:6.14}, Theorem~\ref{thm-BTIS} generalizes to all functions that are Lipschitz with respect to the~$\ell^\infty$-norm:

\begin{mycorollary}[Gaussian concentration, a general case]
\label{cor-6.6}
Let~$f\colon\R^n\to\R$ be such that for some~$M>0$ and all~$x,y\in\R^n$,
\begin{equation}
\label{E:6.15a}
\bigl|f(y)-f(x)\bigr|\le M\max_{i=1,\dots,n}|x_i-y_i|.
\end{equation}
Then for any centered Gaussian~$X$ on~$\R^n$ with~$\sigma_X$ as in \eqref{E:6.1} and any~$t\ge0$,
\begin{equation}
P\bigl(\,f(X)-E f(X)>t\bigr)\le \texte^{-\frac{t^2}{2M^2\sigma_X^2}}.
\end{equation}
\end{mycorollary}

\begin{proofsect}{Proof}
Let~$A$ be the $n\times n$ matrix such that \eqref{E:6.14} holds. From \eqref{E:6.15a} and \eqref{E:6.13} we get
\begin{equation}
\bigl|f(Ay)-f(Ax)\bigr|\le M\max_{i=1,\dots,n}\bigl|(A(x-y))_i\bigr|\le M\sigma_X|x-y|\,.
\end{equation}
Now apply Lemma~\ref{lemma-6.4}.
\end{proofsect}

\section{Fernique majorization}
\noindent
As noted before, the Borell-TIS inequality is of little use unless we have a way to control the expected maximum of a large collection of Gaussian random variables. Our next task is to introduce a method for this purpose. We will actually do this for the supremum over a countable family of such variables as that requires no additional effort. A principal notion here is that of the canonical (pseudo) metric $\rho_X$ associated via \eqref{E:5.36ua} with the Gaussian process~$\{X_t\colon t\in T\}$ on any set~$T$. Our principal result here is:

\begin{mytheorem}[Fernique majorization]
\label{thm-F}
There is~$K\in(0,\infty)$ such that the following holds for any Gaussian process $\{X_t\colon t\in T\}$ over a countable set~$T$ for which~$(T,\rho_X)$ is totally bounded: For any probability measure~$\mu$ on~$T$, we have
\begin{equation}
\label{E:6.15}
E\bigl(\,\sup_{t\in T}X_t\bigr)\le K\sup_{t\in T}\int_0^\infty\textd r\sqrt{\log\frac1{\mu(B(t,r))}}\,,
\end{equation}
where $B(t,r):=\{s\in T\colon\rho_X(t,s)<r\}$.
\end{mytheorem}

A measure~$\mu$ for which the integral in \eqref{E:6.15} converges is called a \emph{majorizing measure}. Note that the integral exists because the integrand is non-increasing and left-continuous. Also note that the domain of integration is effectively bounded because $\mu(B(t,r))=1$ whenever~$r$ exceeds the $\rho_X$-diameter of~$T$, which is finite by the assumed total boundedness.

The above theorem takes its origin in Dudley's work~\cite{Dudley} whose main result is the following theorem:

\begin{mytheorem}[Dudley's inequality]
\label{thm-Dudley}
For the same setting as in the previous theorem, there is a universal constant~$K\in(0,\infty)$ such that
\begin{equation}
\label{E:6.16}
E\bigl(\,\sup_{t\in T}X_t\bigr)\le K\int_0^\infty\textd r\sqrt{\log N_X(r)}\,,
\end{equation}
where~$N_X(r)$ is the minimal number of $\rho_X$-balls of radius~$r$ that are needed to cover~$T$.
\end{mytheorem} 

We will prove Dudley's inequality by modifying a couple of last steps in the proof of Fernique's estimate. Dudley's inequality is advantageous as it sometimes easier to work with. To demonstrate its use we note that the setting of the above theorems is so general that they fairly seamlessly connect boundedness of Gaussian processes to sample-path continuity. Here is an exercise in this vain:

\begin{myexercise}
\label{ex:6.8}
Apply Dudley's inequality to the process $X_{t,s}:=X_t-X_s$ with $t,s\in T$ restricted by $\rho_X(t,s)\le R$ to prove, for $K'$ a universal constant,
\begin{equation}
\label{E:6.17}
E\Bigl(\,\sup_{\begin{subarray}{c}
t,s\in T\\\rho_X(t,s)\le R
\end{subarray}}
|X_t-X_s|\Bigr)\le K'\int_0^R\textd r\sqrt{\log N_X(r)}\,.
\end{equation}
Conclude that if~$r\mapsto \sqrt{\log N_X(r)}$ is integrable, then $t\mapsto X_t$ has a version with (uniformly) $\rho_X$-continuous sample paths a.s.
\end{myexercise}

If this exercise seems too hard at first, we suggest that the reader first reads the proof of Theorem~\ref{thm-6.13} and solves Exercise~\ref{ex:6.14}. To see \eqref{E:6.17} in action, it is instructive to solve:

\begin{myexercise}
Use Dudley's inequality \eqref{E:6.17} to prove the existence of a $\rho_X$-continuous version for the following Gaussian processes:
\begin{enumerate}
\item[(1)] the standard Brownian motion, i.e., a centered Gaussian process $\{B_t\colon t\in[0,1]\}$ with $E(B_tB_s)=t\wedge s$,
\item[(2)] the Brownian sheet, i.e., a centered Gaussian process $\{W_t\colon t\in[0,1]^d\}$ with
\begin{equation}
E(W_tW_s)=\prod_{i=1}^d(t_i\wedge s_i),
\end{equation}
\item[(3)] any centered Gaussian process $\{X_t\colon t\in[0,1]\}$ such that
\begin{equation}
E\bigl([X_t-X_s]^2\bigr)\le c[\log(1/|t-s|)]^{-1-\delta}
\end{equation}
for some~$\delta>0$ and~$c>0$ and~$|t-s|$ sufficiently small.
\end{enumerate}
\end{myexercise}

\noindent
The continuity of these processes can as well be proved via the Komogorov-\v Censtov condition. Both techniques play small probability events against the entropy arising from their total count. (Roughly speaking, this is why the logarithm of $N_X(r)$ appears; the square root arises from Gaussian tails.) Both techniques offer an extension to the proof of uniform H\"older continuity.

Notwithstanding the above discussion, an advantage of Fernique's bound over Dudley's inequality is that it allows optimizing over the probability measure~$\mu$. This is in fact all one needs to get a \emph{sharp} estimate. Indeed, a celebrated result of Talagrand~\cite{Talagrand} shows that a choice of~$\mu$ exists such that the corresponding integral bounds the expectation from below modulo a universal multiplicative constant. This is known to fail for Dudley's inequality. The optimizing~$\mu$ may be thought of as the distribution of the point~$t$ where the maximum of~$X_t$ is achieved, although this has not been made rigorous.

\section{Proof of Fernique's estimate}
\noindent
We will now give a proof of Fernique's bound but before we get embroiled in detail, let us outline the main idea. The basic strategy is simple: We identify an auxiliary centered Gaussian process $\{Y_t\colon t\in T\}$ whose intrinsic distance function~$\rho_Y$ dominates~$\rho_X$. The Sudakov-Fernique inequality then bounds the expected supremum of~$X$ by that of~$Y$. 

For the reduction to be useful, the~$Y$-process must be constructed with a lot of independence built in from the outset. This is achieved by a method called \emph{chaining}. First we organize the elements of~$T$ in a kind of tree structure by defining, for each~$n\in\N$, a map~$\pi_n\colon T\to T$ whose image is a finite set such that the $\rho_X$-distance between~$\pi_{n-1}(t)$ and~$\pi_n(t)$ tends to zero \emph{exponentially} fast with~$n\to\infty$, uniformly in~$t$. Assuming that~$\pi_0(T)$ is a singleton,~$\pi_0(T)=\{t_0\}$, the Borel-Cantelli estimate then allows us to write
\begin{equation}
X_t-X_{t_0}=\sum_{n=1}^\infty \bigl[X_{\pi_n(t)}-X_{\pi_{n-1}(t)}\bigr]\,,
\end{equation}
with the sum converging a.s.\ for each~$t\in T$. We then define~$Y$ by replacing the increments $X_{\pi_n(t)}-X_{\pi_{n-1}(t)}$ by \emph{independent} random variables with a similar variance. The intrinsic distances for~$Y$ can be computed quite explicitly and shown, thanks to careful choices in the definition of~$\pi_n$, to dominate those for~$X$. A computation then bounds the expected supremum of~$Y_t$ by the integral in \eqref{E:6.15}.
We will now begin with the actual proof:

\begin{proofsect}{Proof of Theorem~\ref{thm-F}}
Assume the setting of the theorem and fix a probability measure~$\mu$ on~$T$. The proof (whose presentation draws on Adler~\cite{Adler}) comes in five steps. 

\medskip\noindent\textsl{STEP~1: Reduction to unit diameter.}
If~$D:=\diam(T)$ vanishes, $T$ is effectively a singleton and the statement holds trivially. So we may assume $D>0$. The process~$\wt X_t:=D^{-1/2}X_t$ has a unit diameter. In light of $\rho_{\wt X}(s,t)=D^{-1/2}\rho_X(s,t)$, the $\rho_{\wt X}$-ball of radius~$r$ centered at~$t$ coincides with~$B(t,D^{-1/2} r)$. Passing from~$X$ to~$\wt X$ in \eqref{E:6.15}, both sides scale by factor~$\sqrt D$.

\medskip\noindent\textsl{STEP~2: Construction of the tree structure.}
Next we will define the aforementioned maps~$\pi_n$ subject to properties that will be needed later:

\begin{mylemma}
\label{lemma-6.10}
For each~$n\in\N$ there is~$\pi_n\colon T\to T$ such that
\begin{enumerate}
\item[(1)] $\pi_n(T)$ is finite,
\item[(2)] for each~$t\in T$, we have $\rho(t,\pi_n(t))<2^{-n}$,
\item[(3)] for each~$t\in T$,
\begin{equation}
\mu\bigl(B(\pi_n(t),2^{-n-2})\bigr)\ge\mu\bigl(B(t,2^{-n-3})\bigr),
\end{equation}
\item[(4)] the sets $\{B(t,2^{-n-2})\colon t\in\pi_n(T)\}$ are (pairwise) disjoint.
\end{enumerate}
\end{mylemma}

\begin{proofsect}{Proof}
Fix~$n\in\N$ and, using the assumption of total boundedness, let $t_1,\dots,t_{r_n}$ be points such that 
\begin{equation}
\label{E:6.25nw}
\bigcup_{i=1}^{r_n}B(t_i,2^{-n-3})=T.
\end{equation}
Assume further that these points have been labeled so that
\begin{equation}
\label{E:6.22}
i\mapsto\mu\bigl(B(t_i,2^{-n-2})\bigr)\quad\text{ is non-increasing}.
\end{equation}
We will now identify $\{C_1,C_2,\dots\}\subseteq\{\emptyset\}\cup\{B(t_i,2^{-n-2})\colon i=1,\dots, r_n\}$ by progressively dropping balls that have a non-empty intersection with one of the lesser index. Formally, we set
\begin{equation}
C_1:=B(t_1,2^{-n-2})
\end{equation}
and, assuming that $C_1,\dots,C_i$ have already been defined, let
\begin{equation}
C_{i+1}:=\begin{cases}
B(t_{i+1},2^{-n-2}),\qquad&\text{if }\,\,\displaystyle B(t_{i+1},2^{-n-2})\cap\bigcup_{j=1}^i C_j=\emptyset,
\\
\emptyset,\qquad&\text{else}.
\end{cases}
\end{equation}
Now we define $\pi_n$ as the composition of two maps described informally as follows: Using the ordering induced by \eqref{E:6.22}, first assign~$t$ to the point~$t_i$ of smallest index~$i$ such that~$t\in B(t_i,2^{-n-3})$. Then assign this~$t_i$ to the~$t_j$ with the largest $j\in\{1,\dots,i\}$ such that~$B(t_{i},2^{-n-2})\cap C_j\ne\emptyset$. In summary,
\begin{equation}
\begin{aligned}
i=i(t)&:=\min\bigl\{i=1,\dots,r_n\colon t\in B(t_i,2^{-n-3})\bigr\}
\\
j=j(t)&:=\max\bigl\{j=1,\dots,i(t)\colon B(t_{i(t)},2^{-n-2})\cap C_j\ne\emptyset\bigr\},
\end{aligned}
\end{equation}
where we notice that, by the construction of~$\{C_k\}$, the set in the second line is always non-empty. We then define
\begin{equation}
\pi_n(t):=t_j\quad\text{for}\quad j:=j(t).
\end{equation}
This implies~$\pi_n(T)\subseteq\{t_1,\dots,t_{r_n}\}$ and so~$\pi_n(T)$ is indeed finite, proving~(1). For~(2), using~$i$ and~$j$ for the given~$t$ as above, the construction gives
\begin{equation}
\begin{aligned}
\rho_X(t,\pi_n(t))&=\rho_X(t,t_j)
\\&
\le\rho_X(t,t_i)+\rho_X(t_i,t_j)
\\&\le 2^{-n-3}+2\,2^{-n-2}<2^{-n}.
\end{aligned}
\end{equation}
For~(3) we note that
\begin{equation}
B(t,2^{-n-3})\subseteq B(t_i,2^{-n-2})
\end{equation}
and, by \eqref{E:6.22} and $j\le i$,
\begin{equation}
\mu\bigl(B(t_i,2^{-n-2})\bigr)\le \mu\bigl(B(t_j,2^{-n-2})\bigr).
\end{equation}
Finally, $t_j\in\pi_n(T)$ only if~$C_j\ne\emptyset$ and, when that happens, $C_j=B(t_j,2^{-n-2})$. The construction ensures that the~$C_j$'s are disjoint, thus proving~(4).
\end{proofsect}

\medskip\noindent\textsl{STEP~3: Auxiliary process.}
We are now ready to define the aforementioned process $\{Y_t\colon t\in T\}$. For this, consider a collection $\{Z_n(t)\colon n\in\N,\,t\in\pi_n(T)\}$ of i.i.d.\ standard normals and set
\begin{equation}
Y_t:=\sum_{n\ge1}2^{-n}Z_n\bigl(\pi_n(t)\bigr).
\end{equation}
The sum converges absolutely a.s.\ for each~$t$ due to the fact that the maximum of the first~$n$ terms in a sequence of i.i.d.\ standard normals grows at most like a constant times $\sqrt{\log n}$. We now state:

\begin{mylemma}
For any~$t,s\in T$,
\begin{equation}
\label{E:6.31}
E\bigl([X_t-X_s]^2\bigr)\le 6 E\bigl([Y_t-Y_s]^2\bigr)\,.
\end{equation}
In particular,
\begin{equation}
E\bigl(\sup_{t\in T}X_t\bigr)\le\sqrt6\,E\bigl(\sup_{t\in T}Y_t\bigr)\,.
\end{equation}
\end{mylemma}

\begin{proofsect}{Proof}
We may assume~$\rho_X(t,s)>0$ as otherwise there is nothing to prove. Since $\diam(T)=1$, there is an integer~$N\ge1$ such that $2^{-N}<\rho_X(t,s)\le 2^{-N+1}$. Lemma~\ref{lemma-6.10}(2) and the triangle inequality then show $\pi_n(t)\ne\pi_n(s)$ for all $n\ge N+1$. 
This is quite relevant because the independence built into~$Y_t$ yields
\begin{equation}
E\bigl([Y_t-Y_s]^2\bigr)
=\sum_{n\ge1}2^{-2n}E\Bigl(\bigl[Z_n(\pi_n(t))-Z_n(\pi_n(s))\bigr]^2\Bigr)
\end{equation}
and the expectation on the right vanishes unless~$\pi_n(t)\ne\pi_n(s)$. As that expectation is either zero or 2, we get
\begin{equation}
E\bigl([Y_t-Y_s]^2\bigr)
\ge2\sum_{n\ge N+1}2^{-2n} =2\frac{4^{-(N+1)}}{3/4}=\frac16 4^{-N+1}
\ge\frac16E\bigl([X_t-X_s]^2\bigr)\,,
\end{equation}
where the last inequality follows from the definition of~$N$. This is \eqref{E:6.31}; the second conclusion then follows from the Sudakov-Fernique inequality.
\end{proofsect}

\medskip\noindent\textsl{STEP~4: Majorizing $E(\sup_{t\in T}Y_t)$.}
For the following argument it will be convenient to have a random variable~$\tau$, taking values in~$T$, that identifies the maximizer of~$t\mapsto Y_t$. Such a random variable can certainly be defined when~$T$ is finite. For~$T$ infinite, one has to work with approximate maximizers only. To this end we pose:

\begin{myexercise}
\label{ex:6.12}
Suppose there is~$M\in(0,\infty)$ such that $E(Y_\tau)\le M$ holds for any~$T$-valued random variable~$\tau$ that is measurable with respect to~$\sigma(Y_t\colon t\in T)$. Prove that then also $E(\sup_{t\in T}Y_t)\le M$.
\end{myexercise}

It thus suffices to estimate~$E(Y_\tau)$ for any~$T$-valued random variable~$\tau$. For this we first partition the expectation according to the values of~$\pi_n(\tau)$ as 
\begin{equation}
E(Y_\tau)=\sum_{n\ge1}2^{-n}\sum_{t\in\pi_n(T)}E\bigl(Z_n(t)1_{\{\pi_n(\tau)=t\}}\bigr)\,.
\end{equation}
We now estimate the expectation on the right as follows: Set $g(a):=\sqrt{2\log(1/a)}$ and note that, for~$Z=\NN(0,1)$ and any~$a>0$,
\begin{equation}
E\bigl(Z\1_{\{Z>g(a)\}}\bigr)=\frac1{\sqrt{2\pi}}\int_{g(a)}^\infty x\,\texte^{-\frac12x^2}\textd x
=\frac1{\sqrt{2\pi}}\texte^{-\frac12g(a)^2}=\frac{a}{\sqrt{2\pi}}.
\end{equation}
Therefore, 
\begin{equation}
\label{E:6.38}
\begin{aligned}
E\bigl(Z_n(t)1_{\{\pi_n(\tau)=t\}}\bigr)
&\le E\bigl(Z_n(t)1_{\{Z_n(t)>g(a)\}}\bigr)+g(a)P\bigl(\pi_n(\tau)=t\bigr)
\\
&=\frac{a}{\sqrt{2\pi}}+g(a)P\bigl(\pi_n(\tau)=t\bigr)\,.
\end{aligned}
\end{equation}
Now set $a:=\mu(B(t,2^{-n-2}))$ and perform the sum over~$t$ and~$n$. In the first term we use the disjointness claim from Lemma~\ref{lemma-6.10}(4) to get
\begin{equation}
\label{E:6.42nw}
\frac1{\sqrt{2\pi}}\sum_{t\in\pi_n(T)}\mu(B(t,2^{-n-2}))\le\frac1{\sqrt{2\pi}}
\end{equation}
while in the second term we invoke
\begin{equation}
g\bigl(\mu(B(\pi_n(t),2^{-n-2}))\bigr)\le g\bigl(\mu(B(\tau,2^{-n-3}))\bigr),
\end{equation}
as implied by Lemma~\ref{lemma-6.10}(3), and the fact that~$g$ is non-increasing to get
\begin{multline}
\label{E:6.41}
\quad
\sum_{n\ge1}2^{-n}\sum_{t\in\pi_n(T)}g\bigl(\mu(B(t,2^{-n-2}))\bigr)P\bigl(\pi_n(\tau)=t\bigr)
\\
=E\biggl[\,\sum_{n\ge1}2^{-n}g\bigl(\mu(B(\pi_n(\tau),2^{-n-2}))\bigr)\biggr]
\\
\quad\le E\biggl[\,\sum_{n\ge1}2^{-n}g\bigl(\mu(B(\tau,2^{-n-3}))\bigr)\biggr]
\\
\le\sup_{t\in T}\,\,\sum_{n\ge1}2^{-n}g\bigl(\mu(B(t,2^{-n-3}))\bigr)\,.
\quad
\end{multline}
Using the monotonicity of~$g$,
\begin{equation}
\label{E:6.45nw}
2^{-n}g\bigl(\mu(B(t,2^{-n-3}))\bigr)\le16\,\int_{2^{-n-4}}^{2^{-n-3}}g\bigl(\mu(B(t,r))\bigr)\,\textd r, 
\end{equation}
and so the last sum in \eqref{E:6.41} can now be dominated by $16$-times the integral in the statement of the theorem. Putting the contribution of both terms on the right of \eqref{E:6.38} together, we thus conclude
\begin{equation}
\label{E:6.42}
E(Y_\tau)\le \frac1{\sqrt{2\pi}}+16\sup_{t\in T}\,\int_0^1 g\bigl(\mu(B(t,r))\bigr)\,\textd r.
\end{equation}
Exercise~\ref{ex:6.12} now extends this to a bound on the expected supremum.

\medskip\noindent\textsl{STEP~5: A final touch.}
In order to finish the proof, we need to show that the term~$1/\sqrt{2\pi}$ is dominated by, and can thus be absorbed into, the integral. Here we use the fact that, since~$\diam(T)=1$, there is~$t\in T$ such that $\mu(B(t,1/2))\le1/2$. The supremum on the right of \eqref{E:6.42} is then at least~$\frac12\sqrt{2\log 2}$. The claim follows with $K:=[\frac1{\sqrt{2\pi}}\frac2{\sqrt{\log 2}}+16\sqrt2]\sqrt6$.
\end{proofsect}

The proof of Dudley's inequality requires only minor adaptations:

\begin{proofsect}{Proof of Theorem~\ref{thm-Dudley}}
We follow the previous proof \emph{verbatim} (while disregarding all statements concerning~$\mu$) until \eqref{E:6.38} at which point we instead choose
\begin{equation}
a:=\frac1{N_X(2^{-n-3})}.
\end{equation}
As the number of balls in \eqref{E:6.25nw} could be assumed minimal for the given radius, we have $|\pi_n(T)|\le N_X(2^{-n-3})$ and so the analogue of \eqref{E:6.42nw} applies. Hence,
\begin{equation}
E(Y_\tau)\le\frac1{\sqrt{2\pi}}+\sum_{n\ge1}2^{-n}g\bigl(1/N_X(2^{-n-3})\bigr).
\end{equation}
The sum is converted to the desired integral exactly as in \eqref{E:6.45nw}. The additive prefactor is absorbed by noting that $N_X(1/2)\ge2$ because $\diam(X)=1$.
\end{proofsect}

\section{Consequences for continuity}
\noindent
As already alluded to after the statement of Dudley's inequality, the generality of the setting in which Fernique's inequality was proved permits a rather easy extension to a criterion for continuity. The relevant statement is as follows:

\begin{mytheorem}
\label{thm-6.13}
There is a universal constant~$K'\in(0,\infty)$ such that the following holds for every Gaussian process $\{X_t\colon t\in T\}$ on a countable set~$T$ such that $(T,\rho_X)$ is totally bounded: For any probability measure~$\mu$ on~$T$ and any~$R>0$,
\begin{equation}
\label{E:6.49nw}
E\Bigl(\,\sup_{\begin{subarray}{c}
t,s\in T\\\rho_X(t,s)\le R
\end{subarray}}
|X_t-X_s|
\Bigr)\le
K'\sup_{t\in T}\int_0^R\textd r\sqrt{\log\frac1{\mu(B(t,r))}}\,.
\end{equation}
\end{mytheorem}

\begin{proofsect}{Proof}
We will reduce this to Theorem~\ref{thm-F} but that requires preparations. Let
\begin{equation}
U\subseteq\bigl\{(t,s)\in T\times T\colon\rho_X(t,s)\le R\bigr\}
\end{equation}
be a finite and symmetric set. Denote $Y_{s,t}:=X_t-X_s$ and notice that
\begin{equation}
\rho_Y\bigl((s,t),(s',t')\bigr):=\sqrt{E\bigl([Y_{s,t}-Y_{s',t'}]^2\bigl)}
\end{equation}
obeys
\begin{equation}
\label{E:6.52nw}
\rho_Y\bigl((s,t),(s',t')\bigr)
\le\begin{cases}
\rho(s,s')+\rho(t,t'),
\\
\rho(s,t)+\rho(s',t').
\end{cases}
\end{equation}
Writing~$B_Y$ for the balls (in~$T\times T$) in the $\rho_Y$-metric and~$B_X$ for the balls (in~$T$) in~$\rho_X$-metric, the first line in \eqref{E:6.52nw} then implies
\begin{equation}
\label{E:6.49}
B_Y\bigl((s,t),r\bigr)\supseteq B_X(s,r/2)\times B_X(t,r/2)
\end{equation}
while the second line shows $\diam_{\rho_Y}(U)\le 2R$.
Now define $f\colon T\times T\to U$ by
\begin{equation}
f(y):=\begin{cases}
y,\qquad&\text{if }y\in U,
\\
\operatornamewithlimits{argmin}_U \rho_Y(y,\cdot),\qquad&\text{else},
\end{cases}
\end{equation}
where in the second line the minimizer exists because~$U$ is finite and, in case of ties, is chosen minimal in some \emph{a priori} complete ordering of~$U$.  Then~$f$ is clearly measurable and so, given a probability measure~$\mu$ on~$T$\begin{equation}
\nu(A):=\mu\otimes\mu\bigl(f^{-1}(A)\bigr)\,.
\end{equation}
defines a probability measure on~$U$. Theorem~\ref{thm-F} then yields
\begin{equation}
\label{E:6.56nw}
E\Bigl(\,\,\sup_{(s,t)\in U}
Y_{s,t}
\Bigr)\le K\sup_{(s,t)\in U}\int_0^{2R}\sqrt{\log\frac1{\nu(B_Y((t,s),r))}}\,\textd r\,.
\end{equation}
Our next task is to bring the integral on the right to the form in the statement.

First observe that if~$x\in U$ and $y\in B_Y(x,r)$, then
\begin{equation}
\rho_Y\bigl(x,f(y)\bigr)\le \rho_Y(x,y)+\rho_Y\bigl(y,f(y)\bigr)
\overset{x\in U}\le 2\rho_Y(x,y)\,.
\end{equation}
Hence we get
\begin{equation}
B_Y(x,r)\subseteq f^{-1}\bigl(B_Y(x,2r)\bigr),\qquad x\in U,
\end{equation}
and so, in light of \eqref{E:6.49},
\begin{equation}
\begin{aligned}
\nu\bigl(B_Y((s,t),2r)\bigr)&=\mu\otimes\mu\bigl(f^{-1}\bigl(B_Y((s,t),2r)\bigr)\bigr)
\\
&\ge \mu\otimes\mu\bigl(B_Y((s,t),r)\bigr)
\\
&\ge \mu\bigl(B_X(s,r/2)\bigr)\mu\bigl(B_X(t,r/2)\bigr)\,.
\end{aligned}
\end{equation}
Plugging this in \eqref{E:6.56nw} and invoking $\sqrt{a+b}\le\sqrt a+\sqrt b$, elementary calculus gives\begin{equation}
E\Bigl(\,\,\sup_{(s,t)\in U}
Y_{t,s}
\Bigr)\le 4K\sup_{t\in T}\int_0^R\sqrt{\log\frac1{\mu(B_X(t,r))}}\,\textd r\,.
\end{equation}
Increasing~$U$ to~$U_R:=\{(s,t)\in T\times T\colon\rho_X(s,t)\le R\}$ and applying the Monotone Convergence Theorem, the bound holds for~$U:=U_R$ as well. To connect this to the expectation on the left of \eqref{E:6.49nw},  the symmetry of $\rho_X(t,s)$ and antisymmetry of $|X_t-X_s|$ under the exchange of~$t$ and~$s$ shows
\begin{equation}
E\Bigl(\,\sup_{\begin{subarray}{c}
t,s\in T\\\rho_X(t,s)\le R
\end{subarray}}
|X_t-X_s|
\Bigr)=E\Bigl(\,\,\sup_{(s,t)\in U_R}
Y_{s,t}
\Bigr).
\end{equation}
The claim follows with~$K':=4K$, where~$K$ is as in Theorem~\ref{thm-F}.
\end{proofsect}

Theorem~\ref{thm-6.13} gives us means to prove continuity with respect to the intrinsic metric. However, more often than not, $T$ has its own private metric structure and continuity is desired in the topology thereof. Here the following exercise --- a version of which we already asked for Dudley's inequality \eqref{E:6.16} --- helps:

\begin{myexercise}
\label{ex:6.14}
Suppose $(T,\rho)$ is a metric space, $\{X_t\colon t\in T\}$ a Gaussian process and~$\rho_X$ the intrinsic metric on~$T$ induced thereby. Assume
\settowidth{\leftmargini}{(1111)}
\begin{enumerate}
\item[(1)] $(T,\rho)$ is totally bounded, and
\item[(2)] $s,t\mapsto\rho_X(s,t)$ is uniformly $\rho$-continuous on~$T\times T$.
\end{enumerate}
Prove that, if there is a probability measure~$\mu$ on~$T$ such that
\begin{equation}
\lim_{R\downarrow0}\,\sup_{t\in T}\int_0^R\sqrt{\log\frac1{\mu(B(t,r))}}\,\textd r=0,
\end{equation}
then $X$ admits (uniformly) $\rho$-continuous sample paths on~$T$, a.s.
\end{myexercise}

We note that condition~(2) is necessary for sample path continuity, but definitely not sufficient. To see this, solve:

\begin{myexercise}
Given a measure space $(\scrX,\FF,\nu)$ with~$\nu$ finite, consider the (centered) Gaussian white-noise process $\{W(A)\colon A\in\FF\}$ defined by
\begin{equation}
E\bigl(W(A)W(B)\bigr)=\nu(A\cap B).
\end{equation}
This corresponds to the intrinsic metric $\rho_W(A,B)=\sqrt{\nu(A\triangle B)}$.
Give a (simple) example of $(\scrX,\FF,\nu)$ for which $A\mapsto W(A)$ does not admit $\rho_W$-continuous sample paths. 
\end{myexercise}

\section{Binding field regularity}
\noindent
As our last item of concern in this lecture, we return to the problem of uniform continuity of the binding field for the DGFF and its continuum counterpart. (We used these in the proof of the coupling of the two processes in Lemma~\ref{lemma-4.4}.) The relevant bounds are stated in:

\begin{mylemma}[Binding field regularity]
\label{lemma-6.17}
Let~$\wt D,D\in\mathfrak D$ be such that~$\wt D\subset D$ and $\Leb(D\smallsetminus\wt D)=0$. For~$\delta>0$, denote $\wt D^\delta:=\{x\in D\colon \dist(x,D^\cc)>\delta\}$. Then for each~$\epsilon,\delta>0$,
\begin{equation}
\label{E:6.61}
\lim_{r\downarrow0}\,P\biggl(\,\sup_{\begin{subarray}{c}
x,y\in \wt D^\delta\\|x-y|<r
\end{subarray}}
\bigl|\Phi^{D,\wt D}(x)-\Phi^{D,\wt D}(y)\bigr|>\epsilon\biggr)=0.
\end{equation}
Similarly, given an admissible sequence $\{D_N\colon N\ge1\}$ of approximating domains and denoting~$\wt D_N^\delta:=\{x\in \wt D_N\colon\dist(x,\wt D_N^\cc)>\delta N\}$, for each~$\epsilon,\delta>0$,
\begin{equation}
\label{E:6.62}
\lim_{r\downarrow0}\,\limsup_{N\to\infty}\,P\biggl(\,\sup_{\begin{subarray}{c}
x,y\in \wt D^\delta_N\\|x-y|<rN
\end{subarray}}
\bigl|\varphi^{D_N,\wt D_N}_x-\varphi^{D_N,\wt D_N}_y\bigr|>\epsilon\biggr)=0.
\end{equation}
\end{mylemma}

\begin{proofsect}{Proof of \eqref{E:6.61}}
Consider the set $\wt D^{\delta}_1:=\{x\in\C\colon\dist(x,\wt D^\delta)<\delta/2\}$. 
The intrinsic metric associated with~$\{\Phi^{D,\wt D}(x)\colon x\in \wt D^\delta_1\}$ is given by
\begin{equation}
\rho_\Phi(x,y)=\sqrt{C^{D,\wt D}(x,x)+C^{D,\wt D}(y,y)-2C^{D,\wt D}(x,y)}.
\end{equation}
Since $x\mapsto C^{D,\wt D}(x,y)$ is harmonic on~$\wt D$, it is continuously differentiable and thus uniformly Lipschitz on~$\wt D^\delta_1$. It follows that, for some constant $L=L(\delta)<\infty$,
\begin{equation}
\rho_\Phi(x,y)\le L\sqrt{|x-y|},\qquad x,y\in\wt D^\delta_1.
\end{equation}
Let~$B(x,r):=\{y\in\C\colon|x-y|<r\}$, denote $B_\Phi(x,r):=\{y\in\wt D^\delta_1\colon\rho_\Phi(x,y)<r\}$ and let~$\mu$ be the normalized Lebesgue measure on~$\wt D^\delta_1$. Then
\begin{equation}
B_\Phi(x,L\sqrt r)\supseteq B(x,r),\qquad x\in\wt D^\delta_1,
\end{equation}
while (by the choice of~$\wt D^\delta_1$),
\begin{equation}
\mu\bigl(B(x,r)\bigr)\ge cr^2,\quad x\in\wt D^\delta_1,
\end{equation}
for some~$c=c(\delta)>0$. Combining these observations, we get
\begin{equation}
\mu(B_\Phi(x,r))\ge cL^{-2}r^4.
\end{equation}
 As~$r\mapsto\log(1/r^4)$ is integrable at zero, \eqref{E:6.61} follows from Theorem~\ref{thm-6.13}, Exercise~\ref{ex:6.14} and Markov's inequality.
\end{proofsect}

The above argument could be improved a bit by noting that~$\rho_\Phi$ is itself Lipschitz although that does not change the main conclusion.

\begin{myexercise}
Using an analogous argument with the normalized counting measure replacing the Lebesgue measure, prove \eqref{E:6.62}.
\end{myexercise}


\chapter{Connection to Branching Random Walk}
\label{lec-7}\noindent
In this lecture we return to the two-dimensional DGFF and study the behavior of its absolute maximum beyond the leading order discussed in Lecture~\ref{lec-2}. We begin by recounting the so-called Dekking-Host argument which yields, rather effortlessly, tightness of the maximum (away from its expectation) along a subsequence. Going beyond this will require development of a connection to Branching Random Walk and proving sharp concentration for the maximum thereof. These will serve as the main ingredients for our proof of the tightness of the DGFF maximum in Lecture~\ref{lec-8}.

\section{Dekking-Host argument for DGFF}
\noindent
In Lecture~\ref{lec-2} we already noted that the maximum of the DGFF in a box of side-length~$N$,
\begin{equation}
M_N:=\max_{x\in V_N}h^{V_N}_x\,,
\end{equation}
grows as~$M_N\sim 2\sqrt g\log N$ in probability, with the same leading-order growth rate for~$EM_N$. The natural follow-up questions are then:
\begin{enumerate}
\item[(1)] What is the growth rate of
\begin{equation}
EM_N-2\sqrt g\log N
\end{equation}
with~$N$; i.e., what are the lower-order corrections?
\item[(2)] What is the size of the fluctuations, i.e., the growth rate of~$M_N-EM_N$?
\end{enumerate}
As observed by Bolthausen, Deuschel and Zeitouni~\cite{BDZ} in 2011, an argument that goes back to Dekking and Host~\cite{Dekking-Host} from 1991 shows that, for the DGFF, these seemingly unrelated questions are tied closely together:

\begin{mylemma}[Dekking-Host argument]
\label{lemma-DH}
For~$M_N$ as above and any~$N\ge2$,
\begin{equation}
\label{E:7.2}
E\bigl|M_N-EM_N\bigr|\le 2\bigl(EM_{2N}-EM_N\bigr).
\end{equation}
\end{mylemma}

\nopagebreak
\begin{figure}[t]
\vglue-1mm
\centerline{\includegraphics[width=0.4\textwidth]{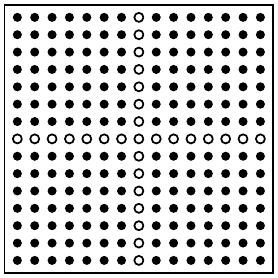}
}
\vglue0mm
\begin{quote}
\small 
\caption{
\label{fig-box}
\small
The partition of box~$V_{2N}$ (both kinds of bullets) into four translates of~$V_N$ for~$N:=8$ and two lines of sites (empty bullets) in the middle. The set $V_{2N}^\circ$ is the collection of all fat bullets.}
\normalsize
\end{quote}
\end{figure}

\begin{proofsect}{Proof}
We will use an idea underlying the solution of the second part of Exercise~\ref{ex:3.4}.
Note that the box~$V_{2N}$ embeds four translates $V_N^{\ssst(1)},\dots,V_N^{\ssst(4)}$ of~$V_N$ that are separated by two lines of sites in-between; see Fig.~\ref{fig-box}. Denoting by $V_{2N}^\circ$ the union of the four translates  of~$V_N$, the Gibbs-Markov property tells us that
\begin{equation}
h^{V_{2N}}:= h^{V_{2N}^\circ}+\varphi^{V_{2N},V_{2N}^\circ},\quad\text{with}\quad  h^{V_{2N}^\circ}\independent \varphi^{V_{2N},V_{2N}^\circ},
\end{equation}
has the law of the DGFF in~$V_{2N}$.
Writing $M_{N}$ for the maximum of $h^{V_{N}}$, using~$X$ to denote the (a.s.-unique) vertex where $h^{V_{2N}^\circ}$ achieves its maximum on~$V_{2N}^\circ$ and abbreviating
\begin{equation}
M_N^{\ssst(i)}:=\max_{x\in V_N^{\ssst(i)}}h_x^{V_{2N}^\circ},\quad i=1,\dots,4,
\end{equation}
it follows that
\begin{equation}
\label{E:7.6uai}
\begin{aligned}
EM_{2N}&=E\Bigl(\,\max_{x\in V_{2N}}(h^{V_{2N}^\circ}_x+\varphi^{V_{2N},V_{2N}^\circ}_x)\Bigr)
\\
&\ge E\Bigl(\,\max_{x\in V_{2N}^\circ}h^{V_{2N}^\circ}_x+\varphi^{V_{2N},V_{2N}^\circ}_X\Bigr) 
\\
&=E\Bigl(\,\max_{x\in V_{2N}^\circ}h^{V_{2N}^\circ}_x\Bigr) 
=E\bigl(\,\max_{i=1,\dots,4}M_N^{\ssst(i)}\bigr)\,,
\end{aligned}
\end{equation}
where we used that $\varphi^{V_{2N},V_{2N}^\circ}$ is independent of $h^{V_{2N}^\circ}$ and thus also of~$X$ to infer that $E\varphi^{V_{2N},V_{2N}^\circ}_X=0$. The portions of $h^{V_{2N}^\circ}$ restricted to $V_N^{\ssst(1)},\dots, V_N^{\ssst(4)}$ are i.i.d.\ copies of~$h^{V_N}$ and so $\{M_N^{\ssst(i)}\colon i=1,\dots,4\}$ are i.i.d.\ copies of~$M_N$. Dropping two out of the four terms from the last maximum in \eqref{E:7.6uai} then yields
\begin{equation}
M_N'\,\,\laweq\,\,M_N,\,M_N'\independent M_N\quad\Rightarrow\quad EM_{2N}\ge E\max\{M_N,M_N'\}\,.
\end{equation}
Now use $2\max\{a,b\}=a+b+|a-b|$ to turn this into
\begin{equation}
\begin{aligned}
E\bigl|M_N-M_N'\bigr|
&=2E\max\{M_N,M_N'\}-E(M_N+M_N')
\\&\le 2EM_{2N}-2EM_N\,.
\end{aligned}
\end{equation}
To get \eqref{E:7.2}, we apply Jensen's inequality to pass the expectation over~$M_N'$ inside the absolute value.
\end{proofsect}

From the growth rate of $EM_N$ we then readily conclude:

\begin{mycorollary}[Tightness along a subsequence]
There is a (deterministic) sequence $\{N_k\colon k\ge 1\}$ of integers with $N_k\to\infty$ such that $\{M_{N_k}-EM_{N_k}\colon k\ge1\}$ is tight.
\end{mycorollary}

\begin{proofsect}{Proof}
Denote $a_n:=EM_{2^n}$. From \eqref{E:7.2} we know that $\{a_n\}$ is non-decreasing. The fact that $EM_N\le c\log N$ (proved earlier by simple first-moment calculations) reads as $a_n\le c'n$ for~$c':=c\log2$. The increments of an increasing sequence with at most a linear growth cannot tend to infinity, so there must be $\{n_k\colon k\ge1\}$ such that $n_k\to\infty$ and $a_{n_{k}+1}-a_{n_k}\le 2c'$. Setting $N_k:=2^{n_k}$, from Lemma~\ref{lemma-DH} we get $E|M_{N_k}-EM_{N_k}|\le 4c'$. This gives tightness via Markov's inequality. 
\end{proofsect}

Unfortunately, tightness along an (existential) subsequence seems to be all one is able to infer from the leading-order asymptotic of~$EM_N$. If we hope to get any better along this particular line of reasoning, we need to control the asymptotic of $EM_N$ up to terms of order unity. This was achieved by Bramson and Zeitouni~\cite{BZ} in 2012. Their main result reads:

\begin{mytheorem}[Tightness of DGFF maximum]
\label{thm-7.3}
Denote
\begin{equation}
\label{E:7.8}
m_N:=2\sqrt g\,\log N-\frac34\sqrt g\,\log\log (N\vee\texte)\,.
\end{equation}
Then
\begin{equation}
\sup_{N\ge1}\,E\bigl|M_N-m_N\bigr|<\infty.
\end{equation}
As a consequence, $\{M_N-m_N\colon N\ge1\}$ is tight.
\end{mytheorem}

\noindent
This and the next lecture will be spent on proving Theorem~\ref{thm-7.3} using, however, a different (and, in the lecturer's view, easier) approach than that of~\cite{BZ}.

\section{Upper bound by Branching Random Walk}
\begin{figure}[t]
\centerline{\includegraphics[width=0.92\textwidth]{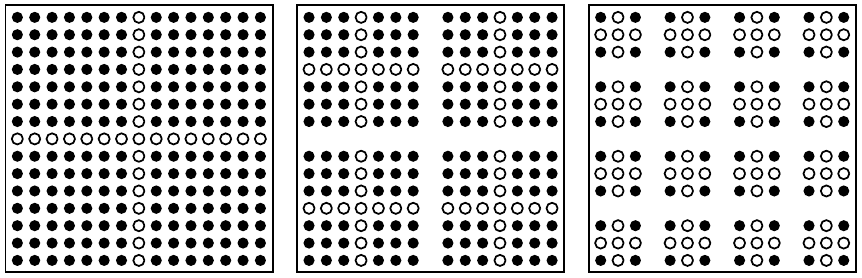}
}
\vglue0mm
\begin{quote}
\small 
\vglue2mm
\caption{
\label{fig-box-hierarchy}
\small
The sets $V_N^{\ssst (0)}=V_N$, $V_N^{\ssst(1)}$ and $V_N^{\ssst(2)}$ underlying the hierarchical representation of the DGFF on~$V_N$ with~$N:=16$. The empty bullets in~$V_N^{\ssst(i)}$ mark the vertices that are being removed to define $V_N^{\ssst (i+1)}$. Boundary vertices (where the fields are set to zero by default) are not depicted otherwise. The binding field $\varphi^{V_N^{\ssst(1)},V_N^{\ssst(2)}}$ is independent on each of the four squares constituting~$V_N^{\ssst(1)}$, but is neither constant nor independent on the squares constituting~$V_N^{\ssst(2)}$.}
\normalsize
\end{quote}
\end{figure}

\nopagebreak\noindent
The Gibbs-Markov decomposition underlying the proof of Lemma~\ref{lemma-DH} can be iterated as follows: Consider a box~$V_N:=(0,N)^2\cap\Z^2$ of side~$N:=2^n$ for some large~$n\in\N$. As illustrated in Fig.~\ref{fig-box-hierarchy}, the square~$V_{2^n}$ then contains four translates of~$V_{2^{n-1}}$ separated by a ``cross'' of ``lines of sites'' in-between, and each of these squares contains four translates of~$V_{2^{n-2}}$, etc. Letting $V_N^{\ssst(i)}$, for $i=1,\dots,n-1$, denote the union of the resulting~$4^i$ translates of~$V_{2^{n-i}}$, and setting $V_N^{\ssst(0)}:=V_N$ and $V_N^{\ssst (n)}:=\emptyset$, we can then write
\begin{equation}
\label{E:7.10}
\begin{aligned}
h^{V_N}&\,\,\laweq\,\, h^{V_N^{\ssst(1)}}+\varphi^{V_N^{\ssst(0)},V_N^{\ssst(1)}}
\\
&\,\,\laweq\,\, h^{V_N^{\ssst(2)}}+\varphi^{V_N^{\ssst(0)},V_N^{\ssst(1)}}+\varphi^{V_N^{\ssst(1)},V_N^{\ssst(2)}}
\\
&\,\,\,\,\,\vdots\qquad\ddots\qquad\quad\ddots
\\
&\,\,\laweq\,\,\varphi^{V_N^{\ssst(0)},V_N^{\ssst(1)}}+\dots+\varphi^{V_N^{\ssst(n-1)},V_N^{\ssst(n)}}\,,
\end{aligned}
\end{equation}
where $V_N^{\ssst (n)}:=\emptyset$ gives $h^{V_N^{\ssst(n-1)}}=\varphi^{V_N^{\ssst(n-1)},V_N^{\ssst(n)}}$. 

The fields in each line on the right-hand side of \eqref{E:7.10} are independent. Moreover, the binding field $\varphi^{\ssst V_N^{\ssst(i)},V_N^{\ssst(i+1)}}$ is a concatenation of $4^i$ independent copies of the binding field $\varphi^{U,V}$ for $U:=V_{2^{n-i}}$ and~$V:=V_{2^{n-i-1}}$, with one for each translate of~$V_{2^{n-i}}$ constituting~$V_N^{\ssst(i)}$. A significant nuisance is that $\varphi^{\ssst V_N^{\ssst(i)},V_N^{\ssst(i+1)}}$ is not constant on these translates. If it were, we would get a representation of the DGFF by means of a Branching Random Walk that we will introduce next. 

For an integer~$b\ge2$, consider a $b$-ary tree~$T^b$ which is a connected graph without cycles where each vertex except one (to be denoted by~$\varnothing$) has exactly $b+1$ neighbors. The distinguished vertex~$\varnothing$ is called the root; we require that the degree of the root is~$b$. We will write~$L_n$ for the set of vertices at graph-theoretical distance~$n$ from the root --- these are the \emph{leaves} at depth~$n$.

Every vertex~$x\in L_n$ can be identified with a sequence
\begin{equation}
(x_1,\dots,x_n)\in\{1,\dots,b\}^n,
\end{equation}
 where~$x_i$ can be thought of as an instruction which ``turn'' to take at the~$i$-th step on the (unique) path from the root to~$x$. The specific case of~$b=4$ can be linked with a binary decomposition of~$V_N:=(0,N)^2\cap\Z^2$ with~$N:=2^n$ as follows: Every~$x\in V_{2^n}$ has non-negative coordinates so it can be written in $\R^2$-vector notation as
\begin{equation}
x = \Bigl(\,\,\sum_{i=0}^{n-1}\sigma_i2^i,\sum_{i=0}^{n-1}\tilde\sigma_i2^i\Bigr)\,,
\end{equation}
for some uniquely-determined~$\sigma_i,\tilde\sigma_i\in\{0,1\}$. Now set the $i$-th instruction for the sequence $(x_1,\dots,x_n)$ as
\begin{equation}
x_i:=2\sigma_{n-i+1}+\tilde\sigma_{n-i+1}+1
\end{equation}
to identify~$x\in V_{2^n}$ with a point in~$L_n$. Since~$V_N$ contains only~$(N-1)^2$ vertices while, for $b:=4$, the cardinality of~$L_n$ is~$4^n$, we only get a subset of~$L_n$.
 
\begin{figure}[t]
\centerline{\includegraphics[width=0.55\textwidth]{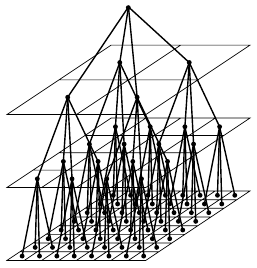}
}
\vglue0mm
\begin{quote}
\small 
\vglue2mm
\caption{
\label{fig-pyramid}
\small
A schematic picture of how the values of the BRW on~$T^4$ are naturally interpreted as a field on~$\Z^2$. There is a close link to the Gibbs-Markov decomposition of the DGFF from \eqref{E:7.10}.
}
\normalsize
\end{quote}
\end{figure}

We now come to:

\begin{mydefinition}[Branching random walk]
Given integers $b\ge2$,~$n\ge1$ and a random variable~$Z$, let $\{Z_x\colon x\in T^b\}$ be i.i.d.\ copies of~$Z$ indexed by the vertices of~$T^b$. The \emph{Branching Random Walk (BRW) on~$T^b$ of depth~$n$ with step distribution~$Z$} is the family of random variables $\{\phi_x^{T^b}\colon x\in L_n\}$ where for~$x=(x_1,\dots,x_n)\in L_n$ we set
\begin{equation}
\label{E:7.14ueu}
\phi_x^{T^b}:=\sum_{k=0}^{n-1} Z_{(x_1,\dots,x_k)},
\end{equation}
with the~$k=0$ term corresponding to the root value~$Z_\varnothing$.
\end{mydefinition}

We write~$n-1$ in \eqref{E:7.14ueu} to ensure that the restriction of $z\mapsto\phi^{T^b}_z-\phi^{T^b}_x$ to the subtree of $T^b$ rooted at~$x$ is independent of~$\phi^{T^b}_x$ with the same law (after proper relabeling) as $\phi^{T^b}$. This is, in fact, a statement of the Gibbs-Markov property for the BRW.

The specific case of interest for us is the \emph{Gaussian Branching Random Walk} where we take~$Z$ normal. The value of the BRW at a given point~$x\in L_n$ is then very much like the last line in~\eqref{E:7.10} --- the sum of~$n$ independent Gaussians along the unique path from the root to~$x$. As already noted, the correspondence is not perfect because of the more subtle covariance structure of the DGFF compared to the BRW and also because~$L_n$ has more active vertices than~$V_{2^n}$. Still, we can use this fruitfully to get:

\begin{mylemma}[Domination of DGFF by BRW]
\label{lemma-UB-BRW}
Consider a BRW $\phi^{T^4}$ on~$T^4$ with step distribution $\NN(0,1)$ and identify~$V_N$ for~$N:=2^n$ with a subset of~$L_n$ as above. There is~$ c>0$ such that for each~$n\ge1$ and each~$x,y\in L_n$,
\begin{equation}
\label{E:7.13}
E\bigl([h^{V_N}_x-h^{V_N}_y]^2\bigr)\le c+(g\log2)\,E\Bigl(\bigl[\phi^{T^4}_x-\phi^{T^4}_y\bigl]^2\Bigr)\,.
\end{equation}
In particular, there is~$k\in\N$ such that for each~$n\ge1$ (and $N:=2^n$),
\begin{equation}
\label{E:7.14}
E\Bigl(\,\max_{x\in V_N}h^{V_N}_x\Bigr)\le \sqrt{g\log2\,}\,E\Bigl(\,\max_{x\in L_{n+k}}\phi^{T^4}_x\Bigr).
\end{equation}
\end{mylemma}

\begin{proofsect}{Proof}
Since~$V\mapsto E([h^V_x-h^V_y]^2)$ is non-decreasing under the set inclusion, the representation of the Green function from Lemma~\ref{lemma-1.19} along with the asymptotic for the potential kernel from Lemma~\ref{lemma-1.20} shows that, for some constant~$\tilde c>0$ and all distinct $x,y\in L_n$,
\begin{equation}
\label{E:7.15}
E\bigl([h^{V_N}_x-h^{V_N}_y]^2\bigr)\le \tilde c+2g\log|x-y|.
\end{equation}
Denoting by $d_n(x,y)$ the \emph{ultrametric distance} between~$x,y\in L_n$, which is defined as the graph-theoretical distance on~$T^b$ from~$x$ to the nearest common ancestor with~$y$, from \eqref{E:7.14ueu} we readily infer 
\begin{equation}
\label{E:7.16}
E\Bigl(\bigl[\phi^{T^4}_x-\phi^{T^4}_y\bigl]^2\Bigr)=2[d_n(x,y)-1]
\end{equation}
for any two distinct $x,y\in L_n$. We now pose:

\begin{myexercise}
There is~$\tilde c'\in(0,\infty)$ such that for each~$n\ge1$ and each~$x,y\in L_n$,
\begin{equation}
\label{E:7.18uai}
|x-y|\le \tilde c'\, 2^{d_n(x,y)}\,.
\end{equation}
\end{myexercise}

\noindent
Combining this with \twoeqref{E:7.15}{E:7.16}, we then get \eqref{E:7.13}. 

To get \eqref{E:7.14}, let $k\in\N$ be so large that~$c$ in \eqref{E:7.13} obeys $c\le 2k(g\log2)$. Now, for each~$x=(x_1,\dots,x_n)\in L_n$ let $\theta(x):=(x_1,\dots,x_n,1,\dots,1)\in L_{n+k}$. Then \eqref{E:7.13} implies
\begin{equation}
\label{E:7.20uai}
E\bigl([h^{V_N}_x-h^{V_N}_y]^2\bigr)\le (g\log2)\,E\Bigl(\bigl[\phi^{T^4}_{\theta(x)}-\phi^{T^4}_{\theta(y)}\bigl]^2\Bigr),\qquad x,y\in L_n.
\end{equation}
The Sudakov-Fernique inequality then gives
\begin{equation}
E\Bigl(\,\max_{x\in V_N}h^{V_N}_x\Bigr)\le \sqrt{g\log2\,}\,E\Bigl(\,\max_{x\in L_{n}}\phi^{T^4}_{\theta(x)}\Bigr).
\end{equation}
The claim now follows by extending the maximum on the right from~$\theta(L_n)$ to all vertices in~$L_{n+k}$.
\end{proofsect}

\section{Maximum of Gaussian Branching Random Walk}
\noindent
In order to use Lemma~\ref{lemma-UB-BRW} to bound the expected maximum of the DGFF, we need a good control of the expected maximum of the BRW. This is a classical subject with strong connections to large deviation theory. (Indeed, as there are~$b^n$ branches of the tree, the maximum will be determined by events whose probability decays exponentially with~$n$. See, e.g., Zeitouni's notes~\cite{Zeitouni-BRW-notes}.) For Gaussian BRWs, we can rely on explicit calculations and so the asymptotic is completely explicit as well:

\begin{mytheorem}[Maximum of Gaussian BRW]
\label{thm-7.7}
For~$b\ge2$, let $\{\phi^{T^b}_x\colon x\in T^b\}$ be the Branching Random Walk on~$b$-ary tree with step distribution~$\NN(0,1)$. Then
\begin{equation}
\label{E:7.20}
E\Bigl(\,\max_{x\in L_n}\phi^{T^b}_x\Bigr)=\sqrt{2\log b\,}\,n-\frac3{2\sqrt{2\log b\,}}\log n+O(1),
\end{equation}
where $O(1)$ is a quantity that is bounded uniformly in~$n\ge1$. Moreover,
\begin{equation}
\Bigl\{\,\max_{x\in L_n}\phi^{T^b}_x-E(\max_{x\in L_n}\phi^{T^b}_x)\colon n\ge1\Bigr\}
\end{equation}
is tight.
\end{mytheorem}

Let us henceforth abbreviate the quantity on the right of \eqref{E:7.20} as
\begin{equation}
\label{E:7.22}
\wt m_n:=\sqrt{2\log b\,}\,n-\frac3{2\sqrt{2\log b}}\log n\,.
\end{equation}
The proof starts by showing that the maximum exceeds~$\wt m_n-O(1)$ with a uniformly positive probability. This is achieved by a second moment estimate of the kind we employed for the intermediate level sets of the DGFF. However, as we are dealing with the absolute maximum, a truncation is necessary. Thus, for~$x=(x_1,\dots,x_n)\in L_n$, let
\begin{equation}
G_n(x):=\bigcap_{k=0}^{n-1}\Bigl\{\phi^{T^b}_{(x_1,\dots,x_{k-1})}\le \frac kn \wt m_n+2\Bigr\}
\end{equation}
be the ``good'' event that curbs the growth the BRW on the unique path from the root to~$x$. Now define
\begin{equation}
\Gamma_n:=\bigl\{x\in L_n\colon\phi^{T^b}_x\ge\wt m_n,\, G_n(x)\text{ occurs}\bigr\}
\end{equation}
as the analogue of the truncated level set $\wh\Gamma^{D,M}_N(b)$ from our discussion of intermediate levels of the DGFF. We now claim:

\begin{mylemma}
\label{lemma-7.8}
For the setting as above,
\begin{equation}
\label{E:7.24}
\inf_{n\ge1}E|\Gamma_n|>0
\end{equation}
while
\begin{equation}
\label{E:7.25}
\sup_{n\ge1}\,E\bigl(|\Gamma_n|^2\bigr)<\infty.
\end{equation}
\end{mylemma}

Let us start with the first moment calculations:

\begin{proofsect}{Proof of \eqref{E:7.24}}
Fix~$x\in L_n$ and, for $k=1,\dots,n$, abbreviate $Z_k:=Z_{(x_1,\dots,x_{k-1})}$ (with~$Z_1:=Z_\varnothing$). Set
\begin{equation}
\label{E:7.40ueu}
S_k:=Z_1+\dots+Z_k,\quad k=1,\dots, n.
\end{equation}
Then
\begin{multline}
\label{E:7.28ueu}
\qquad
P\bigl(\phi^{T^b}_x\ge\wt m_n,\,G_n(x)\text{ occurs}\bigr)
\\
=P\biggl(\{S_n\ge\wt m_n\}\cap \bigcap_{k=1}^n\Bigl\{S_k\le \frac kn\wt m_n+2\Bigr\}\biggr).
\qquad
\end{multline}
In what follows we will make frequent use of:

\begin{myexercise}
\label{ex:Gauss-shift}
Prove that, for $Z_1,\dots,Z_n$ i.i.d.\ normal and any $a\in\R$,
\begin{multline}
\label{E:7.32nw}
\qquad
P\biggl((Z_1,\dots,Z_n)\in\cdot\,\Big|\,\sum_{i=1}^nZ_i=a\biggr)
\\
=P\biggl(\Bigl(Z_1+\frac an,\dots,Z_n+\frac an\Bigr)\in\cdot\,\Big|\,\sum_{i=1}^nZ_i=0\biggr).
\qquad
\end{multline}
\end{myexercise}

\noindent
To see this in action, denote
\begin{equation}
\label{E:7.27}
\mu_n(\textd s):=P(S_n-\wt m_n\in\textd s)
\end{equation}
and use \eqref{E:7.32nw} to express the probability in \eqref{E:7.28ueu} as
\begin{equation}
\label{E:7.29}
\int_0^\infty
\mu_n(\textd s)
P\biggl(\,\bigcap_{k=1}^{n}\Bigl\{S_k\le -\frac kn s+2\Bigr\}\,\bigg|\,S_n=0\biggr).
\end{equation}
As a lower bound, we may restrict the integral to $s\in[0,1]$ which yields $\frac kn s\le1$.
Realizing~$Z_k$ as the increment of the standard Brownian motion on~$[k-1,k)$, the giant probability on the right is bounded from below by the probability that the standard Brownian motion on~$[0,n]$, conditioned on~$B_n=0$, stays below~$1$ for all times in~$[0,n]$. For this we observe:

\begin{myexercise}
\label{ex:7.9}
Let~$\{B_t\colon t\ge0\}$ be the standard Brownian motion started from~$0$. Prove that for all~$a,b>0$ and all~$r>0$,
\begin{equation}
\label{E:7.31ueu}
P^a\Bigl(B_t\ge 0\colon t\in[0,r]\,\Big|\,B_r=b\Bigr)=1-\exp\Bigl\{-2\frac{ab}r\Bigr\}.
\end{equation}
\end{myexercise}

\noindent
Invoking \eqref{E:7.31ueu} with~$a,b:=1$ and~$r:=n$ and applying the shift invariance of the Brownian motion, the giant probability in \eqref{E:7.29} is at most $1-\texte^{-2/n}$. A calculation shows that, for some constant~$c>0$,
\begin{equation}
\label{E:7.36nw}
\mu_n\bigl([0,1]\bigr)
\ge \frac{c}{\sqrt{n}}\,\texte^{-\frac{\wt m_n^2}{2n}}=c\,\texte^{O(n^{-1}\log n)} b^{-n}n
\end{equation}
thanks to our choice of~$\wt m_n$. The product~$n(1-\texte^{-2/n})$ is  uniformly positive and so we conclude that the probability in \eqref{E:7.28ueu} is at least a constant times~$b^{-n}$. Since $|L_n|=b^n$, summing over~$x\in L_n$ we get \eqref{E:7.24}.
\end{proofsect}

We remark that \eqref{E:7.36nw} (and later also \eqref{E:7.47nw}) is exactly what determines the precise constant in the subleading term in \eqref{E:7.22}.
Next we tackle the second moment estimate which is somewhat harder:

\begin{proofsect}{Proof of \eqref{E:7.25}}
Pick distinct~$x,y\in L_n$ and let~$k\in\{1,\dots,n-1\}$ be such that the paths from the root to these vertices have exactly~$k$ vertices (including the root) in common. Let~$S_1,\dots,S_n$ be the values of the BRW on the path to~$x$ and let~$S_1',\dots,S_n'$ be the values on the path to~$y$. Then $S_i=S_i'$ for~$i=1,\dots,k$ while
\begin{equation}
\{S_{k+j}-S_k\colon j=1,\dots,n-k\}\quad\text{and}\quad \{S_{k+j}'-S_k'\colon j=1,\dots,n-k\}
\end{equation}
are independent and each equidistributed to $S_1,\dots,S_{n-k}$. Denoting
\begin{equation}
\mu_{n,k}(\textd s)=P\Bigl(S_k-\frac kn\wt m_n\in\textd s\Bigr),
\end{equation}
conditioning on $S_k-\frac kn\wt m_n$ then implies
\begin{multline}
\qquad
\label{E:7.34uei}
P\Bigl(\phi^{T^b}_x\vee\phi^{T^b}_y\ge\wt m_n,\,G_n(x)\cap G_n(y)\text{ occurs}\Bigr)
\\
=\int_{-\infty}^{2}\mu_{n,k}(\textd s)\,f_k(s)\,g_{k,n-k}(s)^2\,,
\qquad
\end{multline}
where
\begin{equation}
f_k(s):=P\biggl(\,\bigcap_{j=1}^{k}\Bigl\{S_j\le\frac{j}n\wt m_n+2\Bigr\}\,\bigg|\, S_k-\frac kn\wt m_n=s\biggr)
\end{equation}
and
\begin{equation}
g_{k,r}(s):=P\biggl(\,\bigcap_{j=1}^{r}\Bigl\{S_j\le\frac j n\wt m_n+2-s\Bigr\}\cap\Bigl\{S_{r}\ge\frac{r}n\wt m_n-s\Bigr\}\biggr)\,.
\end{equation}
We now need good estimates on both~$f_k$ and~$g_{k,r}$. For this we will need a version of Exercise~\ref{ex:7.9} for random walks:

\begin{myexercise}[Ballot problem]
\label{ex:7.15a}
Let~$S_k:=Z_1+\dots+Z_k$ be the above Gaussian random walk. Prove that there is~$c\in(0,\infty)$ such that for any~$a\ge1$ and any~$n\ge1$,
\begin{equation}
P\biggl(\,\,\bigcap_{k=1}^{n-1}\bigl\{S_k\le a\bigr\}\,\bigg|\,S_n=0\biggr)
\le c\frac{a^2}n.
\end{equation}
\end{myexercise}

\noindent
In what follows we will use~$c$ to denote a positive constant whose meaning may change line to line. 

Concerning an upper bound on~$f_k$,
Exercise~\ref{ex:Gauss-shift} gives
\begin{equation}
f_k(s)=P\biggl(\,\bigcap_{j=1}^{k}\Bigl\{S_j\le2-\frac jk\,s\Bigr\}\,\bigg|\, S_k=0\biggr)\,.
\end{equation}
Exercise~\ref{ex:7.15a} (and~$s\le2$) then yields
\begin{equation}
\label{E:7.38uei}
f_k(s)\le c\frac{1+s^2}k\,.
\end{equation}
As for~$g_{k,r}$, we again invoke Exercise~\ref{ex:Gauss-shift} to write
\begin{equation}
g_{k,r}(s)=\int_{-s}^{2-s}\mu_{n,r}(\textd u) \,P\biggl(\,\bigcap_{j=1}^{r}\Bigl\{S_j\le2-s-\frac j{r}u\Bigr\}\,\bigg|\,S_r=0\biggr).
\end{equation}
Exercise~\ref{ex:7.15a} (and~$s\le2$) again gives
\begin{equation}
\label{E:7.40uei}
g_{k,r}(s)\le c\,\frac{1+s^2}r\,P\Bigl(S_r\ge\frac rn\wt m_n -s\Bigr).
\end{equation}
In order to plug this into the integral in \eqref{E:7.34uei}, we invoke the standard Gaussian estimate \eqref{E:2.2} and the fact that $\wt m_n/n$ is uniformly positive and bounded to get, for all $k=1,\dots,n-1$ and all $s\le2$,
\begin{equation}
\label{E:7.42uei}
\begin{aligned}
\mu_{n,k}(\textd s)P\Bigl(S_{n-k}&\ge\frac {n-k}n\wt m_n -s\Bigr)^2
\\
&\le\frac{c}{\sqrt k}\,\texte^{-\frac{(\frac kn\wt m_n+s)^2}{2k}}\Bigl[\frac1{\sqrt{n-k}}\,\texte^{-\frac{(\frac {n-k}n\wt m_n-s)^2}{2(n-k)}}\Bigr]^2\,\textd s
\\
&\le\frac{ c}{\sqrt k\,(n-k)}\,\texte^{-\frac12(\wt m_n/n)^2(2n-k)+(\wt m_n/n)s}\,\textd s\,.
\end{aligned}
\end{equation}
(When $\frac{n-k}n\wt m_n-s$ is not positive, which for $n$ large happens only when $\log b\le2$, $k=n-1$ and $s$ is close to~$2$, the probability is bounded simply by one.)
The explicit form of~$\wt m_n$ then gives
\begin{equation}
\label{E:7.47nw}
\texte^{-\frac12(\wt m_n/n)^2(2n-k)}\le c b^{-(2n-k)}n^{3-\frac32 k/n}.
\end{equation}
The exponential factor in~$s$ in \eqref{E:7.42uei} ensures that the integral in \eqref{E:7.34uei}, including all~$s$-dependent terms arising from \eqref{E:7.38uei} and \eqref{E:7.40uei}, is bounded. Collecting the denominators from \eqref{E:7.38uei} and \eqref{E:7.40uei}, we get that
\begin{equation}
\label{E:7.49nwwt}
P\Bigl(\phi^{T^b}_x\vee\phi^{T^b}_y\ge\wt m_n,\,G_n(x)\cap G_n(y)\text{ occurs}\Bigr)
\le c b^{-(2n-k)}\frac{n^{3-\frac32 k/n}}{k^{3/2}(n-k)^3}
\end{equation}
holds whenever $x,y\in L_n$ are distinct and~$k$ is as above.

The number of distinct pairs~$x,y\in L_n$ with the same~$k$ is~$b^{2n-k}$. Splitting off the term corresponding to~$x=y$, from \eqref{E:7.49nwwt} we obtain
\begin{equation}
\label{E:7.50nw}
E\bigl(|\Gamma_n|^2\bigr)\le E|\Gamma_n|+c\sum_{k=1}^{n-1}\frac{n^{3-\frac32 k/n}}{k^{3/2}(n-k)^3}.
\end{equation}
For $k< n/2$ the expression under the sum is bounded by $8k^{-3/2}$ which is summable on all~$k\ge1$. For the complementary~$k$ we use~$k^{-3/2}\le \sqrt8 n^{-3/2}$ and then change variables to~$j:=n-k$ to get 
\begin{equation}
\sum_{n/2\le k<n}\frac{n^{3-\frac32 k/n}}{k^{3/2}(n-k)^3}\le\sqrt8\sum_{1\le j\le n/2}\frac{\texte ^{\frac32 (j/n)\log n}}{j^3}.
\end{equation}
Using~$\ell$ to denote the unique integer such that $\frac n{\log n}\ell\le j<\frac n{\log n}(\ell+1)$, the sum on the right is further bounded by
\begin{equation}
\sum_{j\ge1}\frac{\texte^{3/2}}{j^3}+2\frac{(\log n)^2}{n^2}
\sum_{1\le\ell<\log n}\frac{\texte^{\frac32(\ell+1)}}{\ell^3}.
\end{equation}
The first sum is convergent and, by $\texte^{\frac32(\ell+1)}\le c n^{3/2}$ for $\ell<\log n$, the second part of the expression is at most order $(\log n)^2/\sqrt n$. Using this in \eqref{E:7.50nw} we get $E(|\Gamma_n|^2)\le E|\Gamma_n|+c$ for all~$n\ge1$. The Cauchy-Schwarz inequality and some algebra then show $E(|\Gamma_n|^2)\le1+\sqrt c+c$ thus proving \eqref{E:7.25}.
\end{proofsect}

As a consequence of Lemma~\ref{lemma-7.8} we get:

\begin{mycorollary}
Using the notation~$\wt m_n$ from \eqref{E:7.22},
\begin{equation}
\label{E:7.30}
\inf_{n\ge1}\,P\Bigl(\,\max_{x\in L_n}\phi^{T^b}_x\ge\wt m_n\Bigr)>0.
\end{equation}
\end{mycorollary}

\begin{proofsect}{Proof}
The probability in \eqref{E:7.30} is bounded below by $P(|\Gamma_n|>0)$ which by~\eqref{E:2.13} is bounded from below by $(E|\Gamma_n|)^2/E(|\Gamma_n|^2)$. Thanks to \twoeqref{E:7.24}{E:7.25}, this ratio is positive uniformly in~$n\ge1$.
\end{proofsect}

\section{Bootstrap to exponential tails}
\noindent
Our next task in this lecture is to boost the uniform lower bound \eqref{E:7.30} to an exponential tail estimate. Our method of proof will for convenience be restricted to~$b>2$ (remember that we are interested in $b=4$) and so this is what we will assume in the statement:

\begin{mylemma}[Lower tail]
\label{lemma-7.11}
For each integer~$b>2$ there is $a=a(b)>0$ such~that
\begin{equation}
\label{E:7.31}
\sup_{n\ge1}\,P\Bigl(\,\max_{x\in L_n}\phi^{T^b}_x<\wt m_n-t\Bigr)\le\frac1a\texte^{-at},\quad t>0.
\end{equation}
In particular, ``$\ge$'' holds in \eqref{E:7.20}.
\end{mylemma}

\begin{proofsect}{Proof}
The proof will be based on a percolation argument. Recall that the threshold for site percolation on~$T^b$ is~$p_\cc(b)=1/b$. (This is also the survival threshold of a branching process with offspring distribution Bin($b,p$).) Since $P(Z_x\ge0)=1/2$, for any~$b>2$ there is~$\epsilon>0$ such that the set $\{x\in T^b\colon Z_x\ge\epsilon\}$ contains an infinite connected component a.s. We denote by~$\CC$ the one closest to the origin (breaking ties using an arbitrary \emph{a priori} ordering of the vertices). 

Noting that~$\CC$, if viewed from the point closest to the origin, contains a supercritical branching process that survives forever, the reader will surely be able to solve:
\begin{myexercise}
Show that there are~$\theta>1$ and $c>0$ such that for all~$r\ge1$,
\begin{equation}
\label{E:7.33}
P\Bigl(\exists n\ge r\colon |\CC\cap L_n|<\theta^n\Bigr)\le\texte^{-c r}
\end{equation}
\end{myexercise}

\noindent
Writing again~$c$ for a generic positive constant, we claim that this implies
\begin{equation}
\label{E:7.34}
P\Bigl(|\{x\in L_k\colon\phi^{T^b}_x\ge0\}|<\theta^k\Bigr)\le\texte^{-ck},\quad k\ge1.
\end{equation}
Indeed, a crude first moment estimate shows
\begin{equation}
\label{E:7.37ueu}
P\Bigl(\,\min_{x\in L_r}\phi^{T^b}_x\le- 2\sqrt{\log b\,}r\Bigr)\le cb^{-r},\quad r\ge0.
\end{equation}
Taking~$r:=\delta n$, with~$\delta\in(0,1)$, on the event that $\min_{x\in L_r}\phi^{T^b}_x>- 2\sqrt{\log b\,}r$ and~$\CC\cap L_r\ne\emptyset$, we then have
\begin{equation}
\phi^{T^b}_x\ge -2\sqrt{\log b\,}\delta n+\epsilon(n-\delta n),\qquad x\in \CC\cap L_n.
\end{equation}
Assuming~$\delta>0$ is small enough so that $\epsilon(1-\delta)>2\sqrt{\log b\,}\,\delta$, this gives~$\phi^{T^b}_x\ge0$ for all~$x\in\CC\cap L_n$. Hence \eqref{E:7.34} follows from~\eqref{E:7.33} and \eqref{E:7.37ueu}.

Moving to the proof of \eqref{E:7.31}, fix~$t>0$ and let~$k$ be the largest integer less than~$n$ such that $\wt m_n-t\le \wt m_{n-k}$. Denote $A_k:=\{x\in L_k\colon \phi_x^{T^b}\ge0\}$. On the event in \eqref{E:7.31}, the maximum of the BRW of depth~$n-k$ started at any vertex in~$A_k$ must be less than~$\wt m_{n-k}$. Conditional on~$A_k$, this has probability at most $(1-q)^{|A_k|}$, where~$q$ is the infimum in \eqref{E:7.30}. On the event that $|A_k|\ge\theta^k$, this decays double exponentially with~$k$ and so the probability in \eqref{E:7.31} is dominated by that in \eqref{E:7.34}. The claim follows by noting that $t\approx\sqrt{2\log b\,}\,k$.
\end{proofsect}

With the lower-tail settled, we can address the upper bound as well:

\begin{mylemma}[Upper tail]
\label{lemma-7.12}
For each~$b>0$ there is~$\tilde a=\tilde a(b)>0$ such that
\begin{equation}
\label{E:7.39ueu}
\sup_{n\ge1}\,P\Bigl(\,\max_{x\in L_n}\phi^{T^b}_x>\wt m_n+t\Bigr)\le\frac1{\tilde a}\texte^{-\tilde at},\quad t>0.
\end{equation}
\end{mylemma}

\nopagebreak
\begin{figure}[t]
\vglue4mm
\centerline{\includegraphics[width=0.55\textwidth]{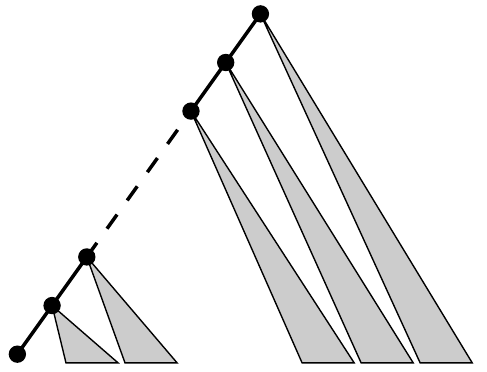}
}
\vglue0mm
\begin{quote}
\small 
\vglue2mm
\caption{
\label{fig-tree}
\small
The picture demonstrating the geometric setup for the representation in~\eqref{E:7.40a}. The bullets mark the vertices on the path from the root (top vertex) to~$x$ (the vertex on the bottom left). The union of the relevant subtrees of these vertices are marked by shaded triangles. The maximum of the field in the subtrees of~$\ell$-th vertex on the path is the quantity in \eqref{E:7.39a}. }
\normalsize
\end{quote}
\end{figure}

\begin{proofsect}{Proof}
The continuity of the involved Gaussian random variables ensures that the maximum occurs at a unique $x=(x_1,\dots,x_n)\in L_n$ a.s. Write $Z_k:=Z_{(x_1,\dots,x_{k-1})}$ (with $Z_1:=Z_\varnothing$) and recall the notation~$S_k$ from \eqref{E:7.40ueu}.
Note that each vertex $(x_1,\dots,x_k)$ on the path from the root to~$x$ has $b-1$ ``children'' $y_1,\dots,y_{b-1}$ not lying on this path. Let~$\wt M_{\ell}^{\ssst(i)}$ denote the maximum of the BRW of depth~$\ell$ rooted at~$y_i$ and abbreviate
\begin{equation}
\label{E:7.39a}
\wh M_{\ell}:=\max_{i=1,\dots,b-1}\wt M_{\ell-1}^{\ssst(i)},
\end{equation}
see Fig.~\ref{fig-tree}. Since the maximum occurs at~$x$ and equals~$\wt m_n+u$ for some~$u\in\R$, we must have $S_k+\wh M_{n-k}\le\wt m_n+u$ for all~$k=1,\dots,n-1$. The symmetries of the BRW then allow us to write the probability in \eqref{E:7.39ueu} as
\begin{equation}
\label{E:7.40a}
b^n\int_t^\infty\mu_n(\textd u)P\biggl(\,\bigcap_{k=1}^{n-1}\bigl\{S_k+\wh M_{n-k}\le\wt m_n+u\bigr\}\,\bigg|\,S_n=\wt m_n+u\biggr),
\end{equation}
where~$\mu_n$ is the measure in \eqref{E:7.27} and where $\wh M_1,\dots,\wh M_n$ are independent of each other and of the random variables $Z_1,\dots,Z_n$ that define the~$S_k$'s. 

We will now estimate \eqref{E:7.40a} similarly as in the proof of Lemma~\ref{lemma-7.8}. First, shifting the normals by their arithmetic mean, the conditional probability in the integral is recast as
\begin{equation}
\label{E:7.37}
P\biggl(\,\bigcap_{k=1}^{n-1}\Bigl\{S_k+\wh M_{n-k}\le \frac{n-k}n(\wt m_n+u)\Bigr\}\,\bigg|\,S_n=0\biggr).
\end{equation}
Letting $\theta_n(k):=(k\wedge(n-k))^{1/5}$, we readily check that
\begin{equation}
\frac{n-k}n\wt m_n\le\wt m_{n-k}+\theta_n(k),\quad k=1,\dots,n,
\end{equation}
as soon as~$n$ is sufficiently large. Introducing 
\begin{equation}
\Theta_n:=\max_{k=1,\dots,n}\bigl[\wt m_{n-k}-\wh M_{n-k}-\theta_n(k)\bigr]_+,
\end{equation}
the probability in \eqref{E:7.37} is thus bounded from above by
\begin{equation}
\label{E:7.43}
P\biggl(\,\,\bigcap_{k=1}^{n-1}\Bigl\{S_k\le \Theta_n+2\theta_n(k)+u\Bigr\}\,\bigg|\,S_n=0\biggr).
\end{equation}
We now observe:

\begin{myexercise}[Inhomogenous ballot problem]
\label{ex:7.15}
Let~$S_k:=Z_1+\dots+Z_k$ be the above Gaussian random walk. Prove that there is~$c\in(0,\infty)$ such that for any~$a\ge1$ and any~$n\ge1$,
\begin{equation}
P\biggl(\,\,\bigcap_{k=1}^{n-1}\Bigl\{S_k\le a+2\theta_n(k)\Bigr\}\,\bigg|\,S_n=0\biggr)
\le c\frac{a^2}n.
\end{equation}
(This is in fact quite hard. Check the Appendix of \cite{BL3} for ideas and references.)
\end{myexercise}

\noindent
Noting that~$\Theta_n$ is independent of~$Z_1,\dots,Z_n$, Exercise~\ref{ex:7.15} bounds the probability in \eqref{E:7.43} by a constant times $n^{-1}E([\Theta_n+u]^2)$. The second moment of~$\Theta_n$ is bounded uniformly in~$n\ge1$ because, by Lemma~\ref{lemma-7.11} and a union bound,
\begin{equation}
P(\Theta_n>u)\le\sum_{k=1}^n \texte^{-a(\theta_n(k)+u)},\quad u>0.
\end{equation}
The probability in \eqref{E:7.43} is thus at most a constant times $(1+u^2)/n$. Since
\begin{equation}
\mu_n\bigl([u,u+1]\bigr)\le c\,\texte^{-(\wt m_n/n) u}\,nb^{-n},\quad u\ge0,
\end{equation}
the claim follows by a routine calculation.
\end{proofsect}

\begin{myremark}
Exercise~\ref{ex:7.15} is our first excursion to the area of ``random walks above slowly-varying curves'' or ``Inhomogenous Ballot Theorems'' which we will encounter several times in these notes. We will not supply detailed proofs of these estimates as these are quite technical and somewhat detached from the main theme of these notes. The reader is encouraged to consult Bramson's seminal work~\cite{Bramson} as well as the Appendix of~\cite{BL3} or the recent posting by Cortines, Hartung and Louidor~\cite{CHL} for a full treatment.
\end{myremark}

We now quickly conclude:

\begin{proofsect}{Proof of Theorem~\ref{thm-7.7}}
Combining Lemmas~\ref{lemma-7.11}--\ref{lemma-7.12}, the maximum has exponential tails away from~$\wt m_n$, uniformly in~$n\ge1$. This yields the claim.
\end{proofsect}

Notice that the underlying idea of the previous proof is to first establish a bound on the lower tail of the maximum and then use it (via a bound on the second moment of~$\Theta_n$) to control the upper tail. A similar strategy, albeit with upper and lower tails interchanged, will be used in the next lecture to prove tightness of the DGFF maximum.

\chapter{Tightness of DGFF maximum}
\label{lec-8}\noindent
We are now ready to tackle the tightness of the DGFF maximum stated in Theorem~\ref{thm-7.3}. The original proof due to Bramson and Zeitouni~\cite{BZ} was based on comparisons with the so called modified Branching Random Walk. We bypass this by proving tightness of the upper tail directly using a variation of the Dekking-Host argument and controlling the lower tail via a concentric decomposition of the DGFF. This brings us closer to what we have done for the Branching Random Walk. The concentric decomposition will be indispensable later as well; specifically, in the analysis of the local structure of nearly-maximal local maxima and the proof of distributional convergence of the DGFF maximum.

\section{Upper tail of DGFF maximum}
\noindent
Recall the notation~$m_N$ from \eqref{E:7.8} and $\wt m_n$ from \eqref{E:7.22}. As is easy to check, for~$b:=4$ and~$N:=2^n$ we have
\begin{equation}
\label{E:7.42}
\sqrt{g\log 2\,}\,\wt m_n = m_N+O(1)
\end{equation}
and so \eqref{E:7.14} and \eqref{E:7.20} yield $EM_N\le m_N+O(1)$. Unfortunately, this does not tell us much by itself (indeed, the best type of bound we can extract from this is that $P(M_N>2m_N)$ is at most about a half.) Notwithstanding, the argument can be enhanced to yield tightness of the upper tail of~$M_N$:

\begin{mylemma}[Upper tail tightness]
\label{lemma-7.13}
We have
\begin{equation}
\sup_{N\ge1}E\bigl((M_N-m_N)_+\bigr)<\infty.
\end{equation}
\end{mylemma}

For the proof we will need the following general inequality:

\begin{mylemma}[Between Slepian and Sudakov-Fernique]
\label{lemma-SF-extension}
Suppose~$X$ and~$Y$ are centered Gaussians on~$\R^n$ such that
\begin{equation}
\label{E:8.3a}
E\bigl((X_i-X_j)^2\bigr)\le E(\bigl(Y_i-Y_j)^2\bigr),\quad i,j=1,\dots,n
\end{equation}
and
\begin{equation}
\label{E:8.3b}
E(X_i^2)\le E(Y_i^2),\quad i=1,\dots,n.
\end{equation}
Abbreviate
\begin{equation}
M_X:=\max_{i=1,\dots,n}X_i\quad\text{and}\quad
M_Y:=\max_{i=1,\dots,n}Y_i
\end{equation}
and let~$M_Y'\,\laweq\,M_Y$ be such that~$M_Y'\independent M_Y$. Then
\begin{equation}
E\bigl((M_X-EM_Y)_+\bigr)\le E\bigl(\max\{M_Y,M_Y'\}\bigr)-EM_Y.
\end{equation}
\end{mylemma}

\begin{proofsect}{Proof}
Let~$Y'$ be a copy of~$Y$ and assume~$X,Y,Y'$ are realized as independent on the same probability space. Define random vectors $Z,\wt Z\in\R^{2n}$  as
\begin{equation}
Z_i:=X_i\quad\text{and}\quad Z_{n+i}:=Y_i',\quad i=1,\dots,n
\end{equation}
and
\begin{equation}
\wt Z_i:=Y_i\quad\text{and}\quad \wt Z_{n+i}:=Y_i',\quad i=1,\dots,n\,.
\end{equation}
Since $E((X_i-Y_j')^2)=E(X_i^2)+E((Y_j')^2)$, from \twoeqref{E:8.3a}{E:8.3b} we readily get
\begin{equation}
E\bigl((Z_i-Z_j)^2\bigr)\le E\bigl((\wt Z_i-\wt Z_j)^2\bigr),\quad i,j=1,\dots,2n.
\end{equation}
Writing $M_Y':=\max_{i=1,\dots,n}Y'_i$, from the Sudakov-Fernique inequality we infer
\begin{equation}
E\bigl(\max\{M_X,M_Y'\}\bigr)\le E\bigl(\max\{M_Y,M_Y'\}\bigr).
\end{equation}
Invoking $(a-b)_+=\max\{a,b\}-b$ and using Jensen's inequality to pass the expectation over~$Y'$ inside the positive-part function, the claim follows.
\end{proofsect}

\begin{proofsect}{Proof of Lemma~\ref{lemma-7.13}}
Let~$n$ be the least integer such that~$N\le 4^n$, assume that~$V_N$ is naturally embedded into~$L_n$ and consider the map~$\theta\colon L_n\to L_{n+k}$ for~$k$ as in Lemma~\ref{lemma-UB-BRW}. Assume~$h^{V_N}$ and~$\phi^{T^4}$ are realized independently on the same probability space, denote
\begin{equation}
\wt M_n:=\sqrt{g\log 2}\,\max_{x\in L_n}\phi^{T^4}_x\,,
\end{equation}
and observe that
\begin{equation}
\label{E:wtMineq1}
\sqrt{g\log 2}\,\,\max_{x\in L_n}\phi^{T^4}_{\theta(x)}\le M_{n+k}
\end{equation}
and 
\begin{equation}
\label{E:wtMineq2}
\sqrt{g\log 2}\,\,E\Bigl(\,\max_{x\in L_n}\phi^{T^4}_{\theta(x)}\Bigr)\ge E\wt M_n\,,
\end{equation}
where the second inequality follows by the same reasoning as \eqref{E:7.6uai} in the Dekking-Host argument.  
In light of \eqref{E:7.20uai} and
\begin{equation}
E\bigl([h^{V_N}_x]^2\bigr)\le (g\log2)\,E\bigl([\phi^{T^4}_{\theta(x)}]^2\bigr),\qquad x\in L_n,
\end{equation}
(proved by the same computation as \eqref{E:7.20uai}) the conditions \twoeqref{E:8.3a}{E:8.3b} are satisfied for~$X:=h^{V_N}$ and~$Y:=\sqrt{g\log2}\,\phi^{T^4}_{\theta(\cdot)}$ (both indexed by~$L_n$).
Lemma~\ref{lemma-SF-extension} along with \twoeqref{E:wtMineq1}{E:wtMineq2} and downward monotonicity of~$b\mapsto(a-b)_+$
show
\begin{equation}
E\bigl((M_N-E\wt M_{n+k})_+\bigr)\le E\bigl(\max\{ \wt M_{n+k},\wt M_{n+k}'\}\bigr)-E(\wt M_n),
\end{equation}
where
\begin{equation}
\wt M_{n+k}'\,\laweq\, \wt M_{n+k}\quad\text{with}\quad
\wt M_{n+k}'\independent \wt M_{n+k}.
\end{equation}
The maximum on the right-hand side is now dealt with as in the Dekking-Host argument (see the proof of Lem\-ma~\ref{lemma-DH}). Indeed, the definition of the BRW gives
\begin{equation}
M_{n+k+1}\,\,\overset{\text{law}}\ge\,\, Z+\max\{ \wt M_{n+k},\wt M_{n+k}'\}\,,
\end{equation}
 where~$Z\,\laweq\,\NN(0,g\log 2)$ is independent of~$\wt M_{n+k}$ and~$\wt M_{n+k}'$. It follows that
\begin{equation}
E\bigl((M_N-E\wt M_{n+k})_+\bigr)\le E(\wt M_{n+k+1})-E(\wt M_n).
\end{equation}
By Theorem~\ref{thm-7.7} and \eqref{E:7.42}, both the right-hand side and $E\wt M_{n+k}-m_N$ are bounded uniformly in~$N\ge1$, thus proving the claim.
\end{proofsect}

Once we know that~$(M_N-m_N)_+$ cannot get too large, we can bootstrap this to an exponential upper tail, just as in Lemma~\ref{lemma-7.12} for the BRW:

\begin{mylemma}
\label{lemma-8.2}
There are~$\tilde a>0$ and~$t_0>0$ such that
\begin{equation}
\label{E:8.8ueu}
\sup_{N\ge1}\,P\bigl(M_N\ge m_N+t\bigr)\le\texte^{-\tilde a t},\qquad t\ge t_0.
\end{equation}
\end{mylemma}

\nopagebreak
\begin{figure}[t]
\vglue-1mm
\centerline{\includegraphics[width=0.45\textwidth]{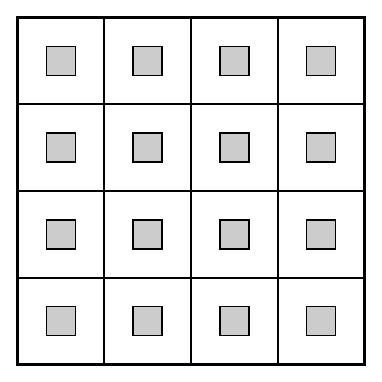}
}
\vglue0mm
\begin{quote}
\small 
\caption{
\label{fig-N-3N}
\small
The geometric setup underlying the proof of Lemma~\ref{lemma-8.2}. The large square marks the domain~$V_{3KN}$ for $K:=4$. The shaded squares are the translates $V_{N}^{\ssst(i)}$, with~$i=1,\dots,K^2$, of $V_N$. }
\normalsize
\end{quote}
\end{figure}

\begin{proofsect}{Proof}
Lemma~\ref{lemma-7.13} and the Markov inequality ensure that, for some $r>0$,
\begin{equation}
\label{E:8.9ueu}
\inf_{N\ge1}P\bigl(\,M_{N}\le m_{N}+r\bigr)\ge\frac12.
\end{equation}
Fix an even integer~$K\ge1$ and consider the DGFF in~$V_{3KN}$. 
Identify~$K^2$ disjoint translates of~$V_{3N}$ inside~$V_{3KN}$ such that any pair of adjacent translates is separated by a line of sites. Denote these translates~$V_{3N}^{\ssst(i)}$, with~$i=1,\dots,K^2$, and abusing our earlier notation slightly, write $V_{3KN}^\circ:=\bigcup_{i=1}^{ K^2}V_{3N}^{\ssst(i)}$. Moreover, let~$V_N^{\ssst(i)}$  be a translate of~$V_N$ centered at the same point as $V_{3N}^{\ssst(i)}$; see Fig.~\ref{fig-N-3N}. Using the Gibbs-Markov decomposition, we then have
\begin{equation}
\label{E:8.10uai}
M_{3KN}\,\,\overset{\text{law}}\ge\,\,\max_{i=1,\dots,K^2}\max_{x\in V_N^{\ssst(i)}}\Bigl(h^{V_{3N}^{\ssst(i)}}_x+\varphi^{V_{3KN},V_{3KN}^\circ}_x\Bigr)\,.
\end{equation}
Consider the event
\begin{equation}
A_K:=\biggl\{\#\Bigl\{i\in\{1,\dots, K^2\}\colon \min_{x\in V_N^{\ssst(i)}}\varphi^{V_{3KN},V_{3KN}^\circ}_x\ge -\sqrt{t}\,\sqrt{\log K}\Bigr\}\ge K^2/2\biggr\}.
\end{equation}
Since $\Var(\varphi^{V_{3KN},V_{3KN}^\circ}_x)\le c\log K$, a combination of Borell-TIS inequality with Fernique's majorization permits us to solve:

\begin{myexercise}
\label{ex:8.4a}
Prove that there are~$a>0$ and~$t_0\ge0$ such that for all~$t\ge t_0$,
\begin{equation}
\sup_{N\ge1}\,\max_{i=1,\dots,K}P\Bigl(\,\min_{x\in V_N^{\ssst(i)}}\varphi^{V_{3KN},V_{3KN}^\circ}_x<-\sqrt{t}\,\sqrt{\log K}\Bigr)\le\texte^{-at}.
\end{equation}
\end{myexercise} 

\noindent
Assuming~$t\ge t_0$, the Markov inequality shows $P(A_K^\cc)\le 2\texte^{-a t}$. The Gibbs-Markov decomposition \eqref{E:8.10uai} (and translation invariance of the DGFF) then yields
\begin{multline}
\label{E:8.13}
\quad
P\bigl(\,M_{3KN}\le m_{3KN}+r\bigr)
\\
\le2\texte^{-a t}+P\Bigl(\,\max_{x\in V_N'}h^{V_{3N}}_x\le m_{KN}+r+\sqrt{t}\,\sqrt{\log K}\Bigr)^{ K^2/2},
\quad
\end{multline}
where~$V_N'$ is the translate of~$V_N$ centered at the same point as~$V_{3N}$. We now invoke the bound $P(X\le s)=1-P(X>s)\le\texte^{-P(X>s)}$ along with Exercise~\ref{ex:3.4} to replace~$h^{V_{3N}}$ by~$h^{V_N'}$ in the maximum over~$x\in V_N'$ at the cost of another factor of $1/2$ popping in front of the probability. This turns \eqref{E:8.9ueu} and \eqref{E:8.13} into
\begin{equation}
\label{E:8.24uea}
\frac12-2\texte^{-a t}\le\exp\biggl\{-\frac14 K^2P\bigl(\,M_N>m_{3NK}+r+\sqrt{t}\,\sqrt{\log K}\Bigr)\biggr\}.
\end{equation}
Assuming~$t$ is so large that~$2\texte^{-a t}\le 1/4$ and applying the inequality
\begin{equation}
m_{3KN}\le m_N+2\sqrt g\log K,
\end{equation}
this proves the existence of~$c'>0$ such that, for all even~$K\ge1$,
\begin{equation}
\sup_{N\ge1}P\Bigl(M_N>m_N+r+2\sqrt{g}\log K+\sqrt{t\log K}\Bigr)\le c' K^{-2}\,.
\end{equation}
Substituting~$t:=c\log K$ for~$c$ large enough that $r+2\sqrt{g}\log K+\sqrt{t\log K}\le t$ then gives \eqref{E:8.8ueu}.
\end{proofsect}

\section{Concentric decomposition}
\noindent
\label{sec:8.2}\noindent
Although the above conclusions seem to be quite sharp, they are not inconsistent with~$M_N$ being concentrated at values much smaller than~$m_N$. (Indeed, comparisons with the BRW will hardly get us any further because of the considerable defect in \eqref{E:7.18uai}. This is what the modified BRW was employed in~\cite{BZ} for but this process is then much harder to study than the BRW.) To address this deficiency we now develop the technique of \emph{concentric decomposition} that will be useful in several parts of these notes.

In order to motivate the definitions to come, note that, to rule out~$M_N\ll m_N$, by Lemma~\ref{lemma-8.2} it suffices to show~$EM_N\ge m_N+O(1)$. In the context of BRW, this was reduced (among other things) to calculating the asymptotic of a probability that for the DGFF takes the form
\begin{equation}
P\bigl(h^{D_N}\le m_N+t\,\big|\, h^{D_N}_0=m_N+t\bigr),
\end{equation}
where we assumed that~$0\in D_N$. For the BRW it was useful (see Exercise~\ref{ex:Gauss-shift}) that the conditional event can be transformed into (what for the DGFF is) $h^{D_N}_0=0$ at the cost of subtracting a suitable linear expression in~$h^{D_N}_0$ from all fields. Such a reduction is possible here as well and yields:

\begin{mylemma}[Reduction to pinned DGFF]
\label{lemma-8.3}
Suppose~$D_N\subset\Z^2$ is finite with~$0\in D_N$. Then for all~$t\in\R$ and~$s\ge0$,
\begin{multline}
\label{E:8.18}
\qquad
P\bigl(h^{D_N}\le m_N+t+s\,\big|\, h^{D_N}_0=m_N+t\bigr)
\\=P\bigl(h^{D_N}\le (m_N+t)(1-\frakg^{D_N})+s\,\big|\, h^{D_N}_0=0\bigr)\,,
\qquad
\end{multline}
where $\frakg^{D_N}\colon\Z^2\to[0,1]$ is discrete-harmonic on~$D_N\smallsetminus\{0\}$ with $\frakg^{D_N}(0)=1$ and~$\frakg^{D_N}=0$ on~$D_N^\cc$. In particular, the probability is non-decreasing in~$s$ and~$t$.
\end{mylemma}

\begin{proofsect}{Proof}
The Gibbs-Markov decomposition of~$h^{D^N}$ reads
\begin{equation}
h^{D_N}\,\,\laweq\,\, h^{D_N\smallsetminus\{0\}}+\varphi^{D_N,D_N\smallsetminus\{0\}}.
\end{equation}
Now $\varphi^{D_N,D_N\smallsetminus\{0\}}$ has the law of the discrete-harmonic extension of~$h^{D_N}$ on~$\{0\}$ to~$D_N\smallsetminus\{0\}$. This means $\varphi^{D_N,D_N\smallsetminus\{0\}} = \frakg^{D_N}h^{D_N}(0)$. Using this, the desired probability can be written as
\begin{equation}
P\bigl(h^{D_N\smallsetminus\{0\}}\le (m_N+t)(1-\frakg^{D_N})+s\bigr).
\end{equation}
The claim then follows from the next exercise.
\end{proofsect}

\begin{myexercise}[Pinning is conditioning]
\label{ex:8.4}
For any finite~$D\subset\Z^2$ with~$0\in D$,
\begin{equation}
(h^D\,|\,h^D_0=0)\,\,\laweq\,\,h^{D\smallsetminus\{0\}}.
\end{equation}
\end{myexercise}

The conditioning the field to be zero is useful for the following reason:

\begin{myexercise}[Pinned field limit]
\label{ex:8.5}
Prove that, for $0\in D_n\uparrow\Z^2$,
\begin{equation}
h^{D_N\smallsetminus\{0\}}\,\,\underset{N\to\infty}\lawarrow\,\, h^{\Z^2\smallsetminus\{0\}}
\end{equation}
 in the sense of finite dimensional distributions.
\end{myexercise}

Let us now inspect the event $h^{D_N}\le m_N(1-\frakg^{D_N})$ in \eqref{E:8.18} --- with~$t$ and~$s$ dropped for simplicity. The following representation using the Green function~$G^{D_N}$ will be useful
\begin{equation}
\label{E:8.22}
m_N\bigl(1-\frakg^{D_N}(x)\bigr)=m_N\frac{G^{D_N}(0,0)-G^{D_N}(0,x)}{G^{D_N}(0,0)}.
\end{equation}
Now $m_N=2\sqrt g\log N+o(\log N)$ while (for~$0$ deep inside~$D_N$) for the Green function we get $G^{D_N}(0,0)=g\log N+O(1)$. With the help of the relation between the Green function and the potential kernel~$\fraka$ as well as the large-scale asymptotic form of~$\fraka$ (see Lemmas~\ref{lemma-1.19} and \ref{lemma-1.20}) we then obtain
\begin{equation}
\label{E:8.24}
m_N\bigl(1-\frakg^{D_N}(x)\bigr)=\frac2{\sqrt g}\fraka(x)+o(1)
=2\sqrt g\,\log|x|+O(1).
\end{equation}
The event $\{h^{D_N}\le m_N(1-\frakg^{D_N})\}\cap\{h^{D_N}_0=0\}$ thus forces the field to stay below the logarithmic cone $x\mapsto 2\sqrt g\log|x|+O(1)$.
 Notice that this is the analogue of the event that the BRW on \emph{all subtrees} along the path to the maximum stay below a linear curve; see~\eqref{E:7.40a}.

In order to mimic our earlier derivations for the BRW, we need to extract as much independence from the DGFF as possible. The Gibbs-Markov property is the right tool to use here. We will work with a decomposition over a sequence of domains defined, with a slight abuse of our earlier notation, by
\begin{equation}
\label{E:8.24ai}
\Delta^k:=\begin{cases}
\{x\in\Z^2\colon|x|_\infty< 2^k\},\qquad&\text{if }k=0,\dots,n-1,
\\
D_N,\qquad&\text{if }k=n,
\end{cases}
\end{equation}
where~$n$ is the largest integer such that $\{x\in\Z^2\colon|x|_\infty\le 2^{n+1}\}\subseteq D_N$; see Fig.~\ref{fig-squares} in Section~\ref{sec:truncate}. The Gibbs-Markov property now gives
\begin{equation}
\label{E:8.25}
h^{D_N}=h^{\Delta^n}\,\,\laweq\,\, h^{\Delta^{n-1}}+h^{\Delta^n\smallsetminus\overline{\Delta^{n-1}}}+\varphi^{\Delta^n,\Delta^n\smallsetminus\partial\Delta^{n-1}}\,,
\end{equation}
where~$\partial D$ stands for the set of vertices on the external boundary of~$D\subset\Z^2$ and~$\overline D:=D\cup\partial D$.
This relation can be iterated to yield:

\begin{mylemma}
For the setting as above,
\begin{equation}
\label{E:8.26}
h^{D_N}\,\,\laweq\,\,\sum_{k=0}^n\bigl(\varphi_k+h_k'\bigr)
\end{equation}
where all the fields on the right are independent with
\begin{equation}
\label{E:8.27}
\varphi_k\,\,\laweq\,\,\varphi^{\Delta^k,\Delta^k\smallsetminus\partial\Delta^{k-1}}
\quad\text{and}\quad h'_k\,\,\laweq\,\, h^{\Delta^k\smallsetminus\overline{\Delta^{k-1}}}
\end{equation}
for~$k=1,\dots,n$ and
\begin{equation}
\label{E:8.28}
\varphi_0\,\,\laweq\,\,h^{\{0\}}\quad\text{and}\quad h'_0=0.
\end{equation}
\end{mylemma}

\begin{proofsect}{Proof}
Apply induction to \eqref{E:8.25} while watching for the provisos at~$k=0$.
\end{proofsect}

The representation \eqref{E:8.26} is encouraging in that it breaks~$h^{D_N}$ into the sum of independent contributions of which one (the~$\varphi_k$'s) is ``smooth'' and timid and the other (the~$h'_k$'s) is, while ``rough'' and  large, localized to the annulus $\Delta^k\smallsetminus\overline{\Delta^{k-1}}$. In order to make the correspondence with the BRW closer, we need to identify an analogue of the Gaussian random walk from \eqref{E:7.40ueu} in this expression. Here we use that, since~$\varphi_k$ is harmonic on~$\Delta^k\smallsetminus\partial\Delta^{k-1}$, its typical value is well represented by its value at the origin. This gives rise to:

\nopagebreak
\begin{figure}[t]
\vglue-1mm
\centerline{\includegraphics[width=0.75\textwidth]{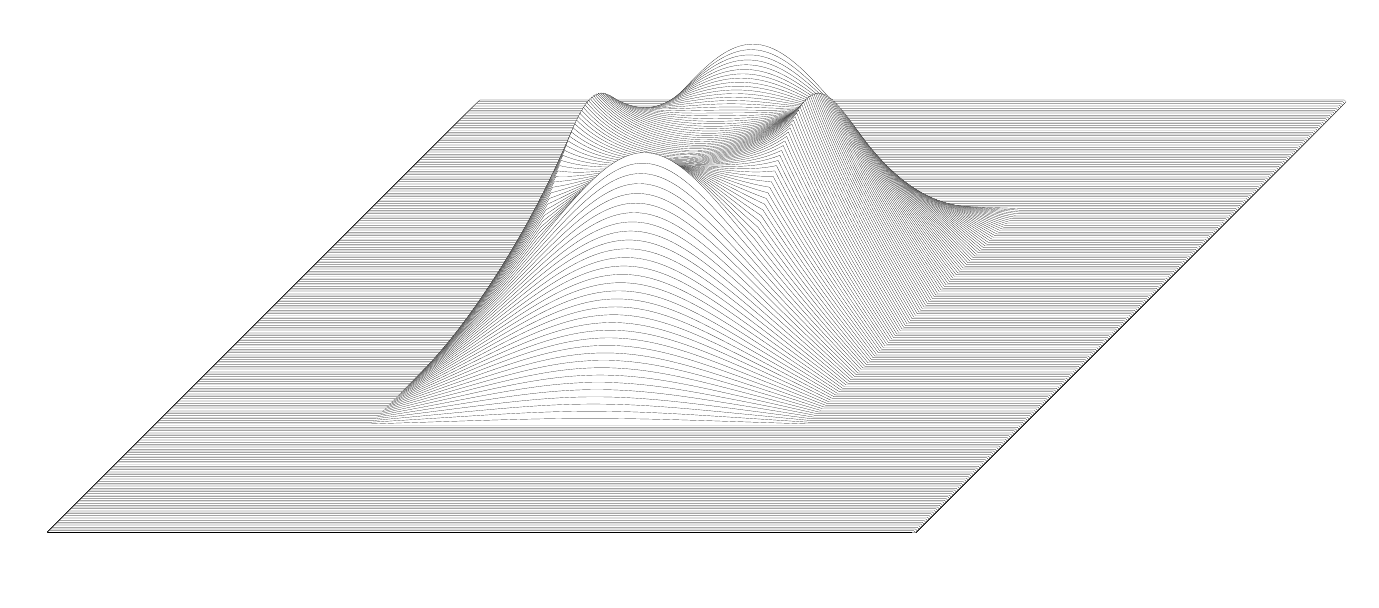}
}
\vglue0mm
\begin{quote}
\small 
\caption{
\label{fig-b-function}
\small
A plot of~$\frakb_k$ on the set~$\Delta^{k+1}$ for~$k$ large. The function equals~$-1$ outside~$\Delta^k$ and vanishes at the origin. It is discrete harmonic on~$\Delta^k\smallsetminus\partial\Delta^{k-1}$.}
\normalsize
\end{quote}
\end{figure}

\begin{myproposition}[Concentric decomposition of DGFF]
For the setting as above,
\begin{equation}
\label{E:8.29}
h^{D_N}\,\,\laweq\,\,\sum_{k=0}^n\Bigl(\bigl(1+\frakb_k\bigr)\varphi_k(0)+\chi_k+h'_k\Bigr),
\end{equation}
where all the random objects in $\bigcup_{k=0}^n\{\varphi_k(0),\chi_k,h_k'\}$ are independent of one another with the law of~$\varphi_k(0)$ and $h_k'$ as in \twoeqref{E:8.27}{E:8.28} and with
\begin{equation}
\label{E:8.30}
\chi_k(\cdot)\,\,\laweq\,\,\varphi_k(\cdot)-E\bigl(\varphi_k(\cdot)\,\big|\,\sigma(\varphi_k(0))\bigr)\,.
\end{equation}
The function $\frakb_k\colon\Z^2\to\R$ is defined by
\begin{equation}
\frakb_k(x):=\frac{E\bigl([\varphi_k(x)-\varphi_k(0)]\varphi_k(0)\bigr)}{E\bigl(\varphi_k(0)^2\bigr)}.
\end{equation}
\end{myproposition}

\begin{proofsect}{Proof}
Define~$\chi_k$ from~$\varphi_k$ by the right-hand side of \eqref{E:8.30}. Then~$\chi_k$ and~$\varphi_k(0)$ are uncorrelated and, being Gaussian, independent. Moreover, the fact that conditional expectation is a projection in~$L^2$ ensures that $E\bigl(\varphi_k(\cdot)\,\big|\,\sigma(\varphi_k(0))\bigr)$ is a linear function of~$\varphi_k(0)$. The fact that these fields have zero mean then implies
\begin{equation}
E\bigl(\varphi_k(x)\,\big|\,\sigma(\varphi_k(0))\bigr) = \frakf_k(x)\varphi_k(0)
\end{equation}
for some deterministic~$\frakf_k\colon\Z^2\to\R$. A covariance computation shows $\frakf_k=1+\frakb_k$. Substituting
\begin{equation}
\varphi_k=(1+\frakb_k)\varphi_k(0)+\chi_k,
\end{equation}
which, we note, includes the case~$k=0$, into \eqref{E:8.26} then gives the claim.
\end{proofsect}

\section{Bounding the bits and pieces}
\noindent
One obvious advantage of \eqref{E:8.29} is that it gives us a representation of DGFF as the sum of \emph{independent}, and reasonably localized, objects. However, in order to make use of this representation, we need estimates on the sizes of these objects as well. The various constants in the estimates that follow will depend on the underlying set~$D_N$ but only via the smallest~$k_1\in\N$ such that 
\begin{equation}
\label{E:8.34}
D_N\subseteq \{x\in\Z^2\colon|x|_\infty\le 2^{n+1+k_1}\}
\end{equation}
with~$n$ as above. We thus assume this~$k_1$ to be fixed; all estimates are then uniform in the domains satisfying \eqref{E:8.34}. We begin with the~$\varphi_k(0)$'s:

\begin{mylemma}
\label{lemma-8.8}
For each~$\epsilon>0$ there is~$k_0\ge0$ such that for all~$n\ge k_0+1$,
\begin{equation}
\max_{k=k_0,\dots,n-1}\,\Bigl|\Var\bigl(\varphi_k(0)\bigr)-g\log2\Bigr|<\epsilon
\end{equation} 
Moreover, $\Var(\varphi_k(0))$ is bounded away from zero and infinity by positive constants that depend only on~$k_1$ from \eqref{E:8.34}.
\end{mylemma}

\begin{proofsect}{Proof (sketch)}
For~$k<n$ large, $\varphi_k(0)$ is close in law to the value at zero of the continuum binding field $\Phi^{B_2,B_2\smallsetminus\partial B_1}$, where $B_r:=[-r,r]^2$. A calculation shows $\Var(\Phi^{B_2,B_2\smallsetminus\partial B_1}(0))=g\log 2$.
\end{proofsect}

\begin{mylemma}
\label{lemma-8.10a}
The function $\frakb_k$ is bounded uniformly in~$k$, is discrete-harmonic on~$\Delta^k\smallsetminus\partial\Delta^{k-1}$ and obeys
\begin{equation}
\label{E:8.48nwt}
\frakb_k(0)=0\quad\text{\rm and}\quad\frakb_k(\cdot)=-1\text{\rm\ \ on }\Z^2\smallsetminus\Delta^k\,.
\end{equation}
Moreover, there is~$c>0$ such that for all~$k=0,\dots,n$,
\begin{equation}
\label{E:8.36}
\bigl|\frakb_k(x)\bigr|\le c\frac{\dist(0,x)}{\dist(0,\partial\Delta^k)},\qquad x\in\Delta^{k-2}.
\end{equation}
\end{mylemma}

\begin{proofsect}{Proof (sketch)}
The harmonicity of~$\frakb_k$ follows from harmonicity of~$\varphi_k$. The overall boundedness is checked by representing~$\frakb_k$ using the covariance structure of~$\varphi_k$. The bound \eqref{E:8.36} then follows from uniform Lipschitz continuity of the (discrete) Poisson kernels in square domains. 
See Fig.~\ref{fig-b-function}.
\end{proofsect}

\begin{mylemma}
\label{lemma-8.11a}
For~$k=0,\dots,n$ and $\ell=0,\dots,k-2$,
\begin{equation}
E\biggl(\max_{x\in\Delta^\ell}\bigl|\chi_k(x)\bigr|\biggr)\le c2^{\ell-k}
\end{equation}
and
\begin{equation}
P\biggl(\Bigl|\max_{x\in\Delta^\ell}\chi_k(x)-E\max_{x\in\Delta^\ell}\chi_k(x)\Bigr|>\lambda\biggr)
\le\texte^{-c4^{k-\ell}\lambda^2}.
\end{equation}
\end{mylemma}

\begin{proofsect}{Proof (idea)}
These follow from Fernique majorization and Borell-TIS inequality and Lipschitz property of the covariances of~$\varphi_k$ (which extend to~$\chi_k$).
\end{proofsect}

\nopagebreak
\begin{figure}[t]
\vglue-1mm
\centerline{\includegraphics[width=0.75\textwidth]{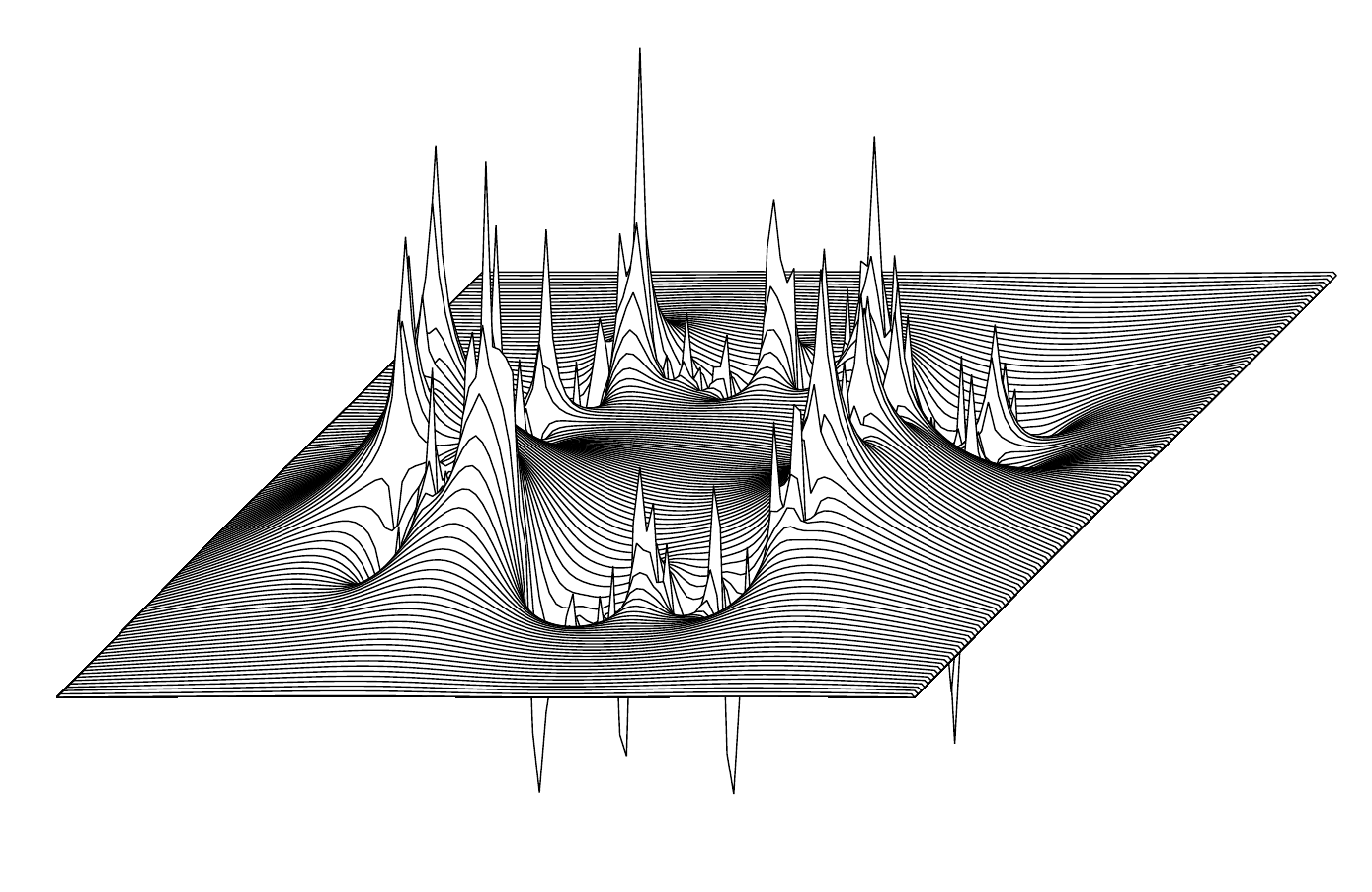}
}
\vglue0mm
\begin{quote}
\small 
\caption{
\label{fig-chi}
\small
A plot of a sample of~$\chi_k$ on~$\Delta^k$ for~$k:=7$. The sample function is discrete harmonic on~$\Delta^k\smallsetminus\partial\Delta^{k-1}$ but quite irregular on and near~$\partial \Delta^{k-1}$.}
\normalsize
\end{quote}
\end{figure}

The case of~$\ell=k-1,k$ has intentionally been left out of Lemma~\ref{lemma-8.11a} because~$\chi_k$ ceases to be regular on and near~$\partial\Delta^{k-1}$, being essentially equal (in law) to the DGFF there; see Fig.~\ref{fig-chi}. Combining~$\chi_k$ with~$h_k'$ and~$\chi_{k+1}$ we get:

\begin{mylemma}[Consequence of upper-tail tightness of~$M_N$]
\label{lemma-8.11}
There exists $a>0$ such that each~$k=1,\dots,n$ and each~$t\ge1$,
\begin{equation}
\label{E:8.39}
P\Bigl(\,\max_{x\in\Delta^k\smallsetminus\Delta^{k-1}}\bigl[\chi_{k+1}(x)+\chi_k(x)+h'_k(x)\bigr] - m_{2^k}\ge t\Bigr)
\le\texte^{-a t}.
\end{equation}
\end{mylemma}

\begin{proofsect}{Proof (sketch)}
Recalling how the concentric decomposition was derived, 
\begin{equation}
\varphi_{k+1}+\varphi_k+h'_k\,\,\laweq\,\,h^{\Delta^k}\quad\text{on }\Delta^k\smallsetminus\Delta^{k-1}.
\end{equation}
Lemma~\ref{lemma-8.2} along with the first half of Exercise~\ref{ex:3.4} show that this field has exponential upper tail above~$m_{2^k}$. But this field differs from the one in the statement by the term $(1+\frakb_k)\varphi_k(0)+ (1+\frakb_{k+1})\varphi_{k+1}(0)$ which has even a Gaussian tail. The claim follows from a union bound.
\end{proofsect}

\begin{myremark}
Once we prove the full tightness of the DGFF maximum, we will augment \eqref{E:8.39} to an estimate on the maximal \emph{absolute value} of the quantity in~\eqref{E:8.39}, see Lemma~\ref{lemma-10.4a}. However, at this point  we claim only a bound on the upper tail in~\eqref{E:8.39}.
\end{myremark}

In light of~$\chi_k(0)=0$, $h_k'(0)=0$ and $\frakb_k(0)=0$ we have $h^{\Delta^n}_0=\sum_{k=0}^n\varphi_k(0)$. This leads to a representation of the field at the prospective maximum by the $(n+1)$st member of the sequence
\begin{equation}
\label{E:8.42a}
S_k:=\sum_{\ell=0}^{k-1}\varphi_\ell(0)\,,
\end{equation}
which we think of as an analogue of the random walk in \eqref{E:7.40ueu} albeit this time with time inhomogeneous steps.
The observation
\begin{equation}
\label{E:8.42}
h^{\Delta^n}_0 = 0\quad\Leftrightarrow\quad S_{n+1}=0
\end{equation}
along with the fact that, for any $k=0,\dots,n-1$, neither $\frakb_k$ nor the laws of $\varphi_k$ and~$h'_k$ depend on~$n$ and~$D_N$ then drive:

\begin{myexercise}[Representation of pinned DGFF]
\label{ex:8.13}
The DGFF on~$\Z^2\smallsetminus\{0\}$ can be represented as the a.s.-convergent sum
\begin{equation}
h^{\Z^2\smallsetminus\{0\}}\,\,\laweq\,\,\sum_{k=0}^\infty\Bigl(\varphi_k(0)\frakb_k+\chi_k+h'_k\Bigr),
\end{equation}
where the objects on the right are independent with the laws as above for the sequence $\{\Delta^k\colon k\ge0\}$ from \eqref{E:8.24ai} with~$n:=\infty$. [Hint: Use Exercise~\ref{ex:8.5}.]
\end{myexercise}

We remark that, besides the connection to the random walk from \eqref{E:7.40ueu} associated with the Gaussian BRW, the random walk in \eqref{E:8.42a} can also be thought of as an analogue of circle averages of the CGFF; see Exercise~\ref{ex:1.29}. As we will show next, this random walk will by and large determine the behavior of the DGFF in the vicinity of a point where the field has a large local maximum.

\section{Random walk representation}
\label{sec8.4}\noindent
We will now move to apply the concentric decomposition in the proof of the lower bound on~$EM_N$. A key technical step in this will be the proof of:

\begin{myproposition}
\label{prop-8.13}
For all~$\epsilon\in(0,1)$ there is~$c=c(\epsilon)>1$ such that for all naturals~$N>2$ and all sets $D_N\subset\Z^2$ satisfying
\begin{equation}
[-\epsilon N,\epsilon N]^2\cap\Z^2\subseteq D_N\subseteq[-\epsilon^{-1}N,\epsilon^{-1}N]^2\cap\Z^2
\end{equation}
we have
\begin{equation}
\label{E:8.46}
P\bigl(\,h^{D_N}\le m_N\,\big|\,h^{D_N}_0=m_N\bigr)\ge\frac {c^{-1}}{\log N}\bigl[1-2P(M_N> m_N-c)\bigr],
\end{equation}
where, abusing our earlier notation, $M_N:=\max_{x\in D_N}h^{D_N}_x$.
\end{myproposition}

In order to prove this, we will need to control the growth of the various terms on the right-hand side of \eqref{E:8.29}. This will be achieved using a single \emph{control variable}~$K$ that we define next:

\begin{mydefinition}[Control variable]
\label{def-control-var}
For~$k,\ell,n$ positive integers denote
\begin{equation}
\label{E:def-theta}
\theta_{n,k}(\ell):=\bigl[\log(k\vee(\ell\wedge(n-\ell)))]^2.
\end{equation}
Then define~$K$ as the smallest~$k\in\{1,\dots,\lfloor \frac n2\rfloor\}$ such that for all~$\ell=0,\dots,n$:
\settowidth{\leftmargini}{(1111)}
\begin{enumerate}
\item[(1)] $|\varphi_\ell(0)|\le\theta_{n,k}(\ell)$,
\item[(2)] for all~$r=1,\dots,\ell-2$,
\begin{equation}
\max_{x\in\Delta^r}\,\bigl|\chi_\ell(x)\bigr|\le 2^{(r-\ell)/2}\theta_{n,k}(\ell)\,,
\end{equation}
\item[(3)] 
\begin{equation}
\label{E:8.62nwwt}
\max_{x\in\Delta^\ell\smallsetminus\Delta^{\ell-1}}\bigl[\chi_\ell(x)+\chi_{\ell+1}(x)+h_\ell'(x)- m_{2^\ell}\bigr]\le\theta_{n,k}(\ell)\,.
\end{equation}
\end{enumerate}
If no such~$K$ exists, we set $K:=\lfloor \frac n2\rfloor+1$.
We call~$K$ the control variable.
\end{mydefinition}

Based on the above lemmas, one readily checks:

\begin{myexercise}
\label{ex:8.17}
For some~$c,c'>0$, all $n\ge1$ and all $k=1,\dots,\lfloor\frac n2\rfloor+1$,
\begin{equation}
P(K=k)\le c'\,\texte^{-c(\log k)^2},\qquad k\ge1.
\end{equation}
\end{myexercise}

As we will only care to control events up to probabilities of order~$1/n$, this permits us to disregard the situations when~$K$ is at least $n^\epsilon$, for any~$\epsilon>0$. Unfortunately, for smaller~$k$ we will need to control the growth of the relevant variables on the background of events whose probability is itself small (the said order~$1/n$). The key step is to link the event in Proposition~\ref{prop-8.13} to the behavior of the above random walk. This is the content of:

\begin{mylemma}[Reduction to a random walk event]
\label{lemma-8.16}
Assume~$h^{D_N}$ is realized as the sum in \eqref{E:8.29}. There is a numerical constant~$C>0$ such that, uniformly in the above setting, the following holds for each~$k=0,\dots,n$:
\begin{multline}
\qquad
\{h^{D_N}_0=0\}\cap\bigl\{h^{D_N}\le m_N(1-\frakg^{D_N})\text{\rm\ on }\Delta^k\smallsetminus\Delta^{k-1}\bigr\}
\\
\supseteq\{S_{n+1}=0\}\cap\bigl\{S_k\ge C[1+\theta_{n,K}(k)]\bigr\}.
\qquad
\end{multline}
\end{mylemma}

\begin{proofsect}{Proof}
Fix~$k$ as above and let~$x\in\Delta^k\smallsetminus\Delta^{k-1}$. In light of \twoeqref{E:8.42a}{E:8.42}, on the event $\{h^{D_N}_0=0\}$ we can drop the ``1'' in the first term on the right-hand side of \eqref{E:8.29} without changing the result. Noting that~$\frakb_\ell(x)=-1$ for~$\ell<k$, on this event we then get
\begin{multline}
\label{E:8.65nwwt2}
\quad
h^{D_N}_x-m_{2^k}=-S_k+\sum_{\ell=k}^n\frakb_\ell(x)\varphi_\ell(0)
\\
+\biggl(\,\sum_{\ell=k+2}^n\chi_\ell(x)\biggr)+\bigl[\chi_{k+1}(x)+\chi_k(x)+h_k'(x)-m_{2^k}\bigr].
\quad
\end{multline}
The bounds in the definition of the control variable permit us to estimate all terms after~$-S_k$ by $C\theta_{n,K}(k)$ from above, for~$C>0$ a numerical constant independent of~$k$ and~$x$. Adjusting~$C$ if necessary, \eqref{E:8.22} along with the approximation of the Green function using the potential kernel and invoking the downward monotonicity of $N\mapsto\frac{\log\log N}{\log N}$ for~$N$ large enough shows
\begin{equation}
\label{E:8.65nwwt}
m_N\bigl(1-\frakg^{D_N}(x)\bigr)\ge m_{2^k}-C.
\end{equation}
Hence,
\begin{multline}
\quad
m_N(1-\frakg^{D_N})-h^{D_N}(x)
\\
\ge m_{2^k}-h^{D_N}(x)-C
\ge S_k-C[1+\theta_{n,K}(k)].
\quad
\end{multline}
This, along with \eqref{E:8.42}, now readily yields the claim.
\end{proofsect}

We are ready to give:

\begin{proofsect}{Proof of Proposition~\ref{prop-8.13}}
We first use Lemma~\ref{lemma-8.3} to write
\begin{equation}
P\bigl(\,h^{D_N}\le m_N\,\big|\,h^{D_N}_0=m_N\bigr) = P\bigl(\,h^{D_N\smallsetminus\{0\}}\le m_N(1-\frakg^{D_N})\bigr)\,.
\end{equation}
Next pick a~$k\in\{1,\dots,\lfloor n/2\rfloor\}$ and note that, by the FKG inequality for the DGFF (see Proposition~\ref{prop-5.24}),
\begin{equation}
\label{E:8.56}
P\bigl(\,h^{D_N\smallsetminus\{0\}}\le m_N(1-\frakg^{D_N})\bigr)\ge P(A_{n,k}^1)P(A_{n,k}^2)P(A_{n,k}^3)\,,
\end{equation}
where we used that
\begin{equation}
\begin{aligned}
A_{n,k}^1&:=\bigl\{h^{D_N\smallsetminus\{0\}}\le m_N(1-\frakg^{D_N})\text{ on }\Delta^k\bigr\}
\\
A_{n,k}^2&:=\bigl\{h^{D_N\smallsetminus\{0\}}\le m_N(1-\frakg^{D_N})\text{ on }\Delta^{n-k}\smallsetminus\Delta^k\bigr\}
\\
A_{n,k}^3&:=\bigl\{h^{D_N\smallsetminus\{0\}}\le m_N(1-\frakg^{D_N})\text{ on }\Delta^n\smallsetminus\Delta^{n-k}\bigr\}
\end{aligned}
\end{equation}
are increasing events. We will now estimate the three probabilities on the right of \eqref{E:8.56} separately.

First we observe that, for any~$k$ fixed,
\begin{equation}
\label{E:8.58}
\inf_{n\ge1} P(A_{n,k}^1)>0
\end{equation}
in light of the fact that $h^{D_N\smallsetminus\{0\}}$ tends in law to~$h^{\Z^2\smallsetminus\{0\}}$ (see Exercise~\ref{ex:8.5}), while $m_N(1-\frakg^{D_N})$ tends to~$\frac2{\sqrt g}\fraka$, see \eqref{E:8.24}. For $A_{n,k}^3$ we note (similarly as in \eqref{E:8.65nwwt}) that, for some $c$ depending only on~$k$,
\begin{equation}
m_N(1-\frakg^{D_N})\ge m_N - c\quad\text{on }\Delta^n\smallsetminus\Delta^{n-k}\,.
\end{equation}
Denoting $M_N:=\max_{x\in D_N}h^{D_N}_x$ and $M_N^0:=\max_{x\in D_N}h^{D_N\smallsetminus\{0\}}_x$, Exercise~\ref{ex:3.3} then yields
\begin{equation}
\label{E:8.60}
P(A_{n,k}^3)\ge P\bigl(M_N^0\le m_N-c)\ge 1-2P\bigl(M_N> m_N-c\bigr)\,.
\end{equation}
We may and will assume the right-hand side to be positive in what follows, as there is nothing to prove otherwise.

It remains to estimate $P(A_{n,k}^2)$. Using Lemma~\ref{lemma-8.16} and the fact that~$k\mapsto\theta_{n,k}(\ell)$ is non-decreasing, we bound this probability as
\begin{equation}
\label{E:8.61}
\begin{aligned}
P(&A_{n,k}^2)\ge P\biggl(\{K\le k\}\cap\bigcap_{\ell=k+1}^{n-k-1}\bigl\{S_\ell\ge C[1+\theta_{n,k}(\ell)]\bigr\}\,\bigg|\, S_{n+1}=0\biggr)
\\
&\ge P\biggl(\,\bigcap_{\ell=k+1}^{n-k-1}\{S_\ell\ge C[1+\theta_{n,k}(\ell)]\bigr\}\cap \bigcap_{\ell=1}^{n}\{S_\ell\ge-1\}\,\bigg|\, S_{n+1}=0\biggr)
\\
&\qquad\qquad\qquad\qquad\quad
-P\biggl(\{K>k\}\cap \bigcap_{\ell=1}^{n}\{S_\ell\ge-1\}\,\bigg|\, S_{n+1}=0\biggr)\,.
\end{aligned}
\end{equation}
To estimate the right-hand side, we invoke the following lemmas:

\begin{mylemma}[Entropic repulsion]
\label{lemma-8.18}
There is a constant~$c_1>0$ such that for all $n\ge1$ and all~$k=1,\dots,\lfloor \ffrac n2\rfloor$
\begin{equation}
P\biggl(\,\bigcap_{\ell=k+1}^{n-k-1}\{S_\ell\ge C[1+\theta_{n,k}(\ell)]\bigr\}\,\bigg|\, \bigcap_{\ell=1}^{n}\{S_\ell\ge-1\}\cap\{S_{n+1}=0\}\biggr)\ge c_1.
\end{equation}
\end{mylemma}

\begin{mylemma}
\label{lemma-8.19}
There is~$c_2>0$ such that for all~$n\ge1$ and all~$k=1,\dots,\lfloor \ffrac n2\rfloor$,
\begin{equation}
P\biggl(\{K>k\}\cap \bigcap_{\ell=1}^{n}\{S_\ell\ge-1\}\,\bigg|\, S_{n+1}=0\biggr)\le\frac1n\texte^{-c_2(\log k)^2}.
\end{equation}
\end{mylemma}

As noted before, we will not supply a proof of these lemmas as that would take us on a detour into the area of ``Inhomogenous Ballot Theorems;'' instead, the reader is asked to consult~\cite{BL3} where these statements are given a detailed proof. We refer to Lemma~\ref{lemma-8.18} using the term ``entropic repulsion'' in reference to the following observation from statistical mechanics: An interface near a barrier appears to be pushed away as that increases the entropy of its available fluctuations. 

Returning to the proof of Proposition~\ref{prop-8.13}, we note that Lemmas~\ref{lemma-8.18} and~\ref{lemma-8.19} reduce \eqref{E:8.61} to a lower bound on the probability of $\bigcap_{\ell=1}^{n}\{S_\ell\ge-1\}$ conditional on~$S_{n+1}=0$. We proceed by embedding the random walk into a path of the standard Brownian motion $\{B_t\colon t\ge0\}$ via
\begin{equation}
\label{E:embed}
S_k:=B_{t_k}\quad\text{where}\quad t_k:=\Var(S_k)=\sum_{\ell=0}^{k-1}\Var\bigl(\varphi_\ell(0)\bigr).
\end{equation}
 This readily yields
\begin{multline}
\quad
P\biggl(\,\bigcap_{\ell=1}^{n}\{S_\ell\ge-1\}\,\bigg|\,S_{n+1}=0\biggr)
\\
\ge P^0\Bigl(B_t\ge-1\colon t\in[0,t_{n+1}]\,\Big|\,B_{t_{n+1}}=0\Bigr).
\quad
\end{multline}
Lemma~\ref{lemma-8.8} ensures that~$t_{n+1}$ grows proportionally to~$n$ and the Reflection Principle then bounds the last probability by~$c_3/n$ for some~$c_3>0$ independent of~$n$; see Exercise~\ref{ex:7.9}. Lemmas~\ref{lemma-8.18}--\ref{lemma-8.19} then show 
\begin{equation}
P(A_{n,k}^2)\ge \frac{c_1c_3}n-\frac1n\texte^{-c_2(\log k)^2}.
\end{equation}
For~$k$ sufficiently large, this is at least a constant over~$n$. Since $n\approx\log_2N$, we get \eqref{E:8.46} from \eqref{E:8.56}, \eqref{E:8.58} and~\eqref{E:8.60}.
\end{proofsect}

\section{Tightness of DGFF maximum: lower tail}
\noindent
We will now harvest the fruits of our hard labor in the previous sections and prove tightness of the maximum of DGFF. First we claim:

\begin{mylemma}
For the DGFF in~$V_N$, we have
\begin{equation}
\label{E:8.65}
\inf_{N\ge1}\,P\bigl(M_N\ge m_N)>0.
\end{equation}
\end{mylemma}

\begin{proofsect}{Proof}
We may and will assume~$N\ge10$ without loss of generality. Let~$V_{N/2}'$ denote the square of side~$\lfloor N/2\rfloor$ centered (roughly) at the same point as~$V_N$. For each~$x\in V_{N/2}'$ and denoting~$D_N:=-x+V_N$, the translation invariance of the DGFF gives
\begin{equation}
\label{E:8.68ueu}
\begin{aligned}
P\bigl(&h^{V_N}_x\ge m_N,\,h^{V_N}\le h^{V_N}_x\bigr)
=P\bigl(h^{D_N}_0\ge m_N,\,h^{D_N}\le h^{D_N}_0\bigr)
\\
&=\int_0^\infty P\bigl(h^{D_N}_0- m_N\in\textd s\bigr)\,P\bigl(h^{D_N}\le m_N+s\,\big|\,h^{D_N}_0=m_N+s\bigr)\,.
\end{aligned}
\end{equation}
Rewriting the conditional probability using the DGFF on~$D_N\smallsetminus\{0\}$, the monotonicity in Lemma~\ref{lemma-8.3} along with Proposition~\ref{prop-8.13} show, for any~$s\ge0$,
\begin{multline}
P\bigl(h^{D_N}\le m_N+s\,\big|\,h^{D_N}_0=m_N+s\bigr)
\\
\ge P\bigl(h^{D_N}\le m_N\,\big|\,h^{D_N}_0=m_N\bigr)\ge\frac{c^{-1}}{\log N}\bigl[1-2P(M_N\le m_N-c)\bigr].
\end{multline}
Plugging this in \eqref{E:8.68ueu} yields
\begin{multline}
\label{E:8.82nwwt}
\qquad
P\bigl(h^{V_N}_x\ge m_N,\,h^{V_N}\le h^{V_N}_x\bigr)
\\
\ge\frac{c^{-1}}{\log N}P\bigl(h^{V_N}_x\ge m_N\bigr)\bigl[1-2P(M_N\le m_N-c)\bigr].
\qquad
\end{multline}
If $P(M_N\le m_N-c)>1/2$, we may skip directly to \eqref{E:8.71ueu}. Otherwise, invoking
\begin{equation}
P\bigl(h^{V_N}_x\ge m_N\bigr)\ge c_1(\log N)N^{-2}
\end{equation}
with some~$c_1>0$ uniformly for all~$x\in V_{N/2}'$, summing \eqref{E:8.82nwwt} over~$x\in V_{N/2}'$ and using that $|V_{N/2}'|$ has order~$N^2$ vertices  shows
\begin{equation}
P(M_N\ge m_N)\ge c_2 \bigl[1-2P(M_N> m_N-c)\bigr]
\end{equation}
for some~$c_2>0$. As~$c>0$, this implies
\begin{equation}
\label{E:8.71ueu}
P(M_N> m_N-c)\ge c_2/(1+2c_2).
\end{equation}
This is almost the claim except for the constant~$c$ in the event. 

Consider the translate~$V_{N/2}'$ of~$V_{\lfloor N/2\rfloor}$ centered at the same point as~$V_N$ and use the Gibbs-Markov property to write $h^{V_N}$ as~$h^{V'_{N/2}}+\varphi^{V_N,V'_{N/2}}$. 
Denoting $M_{N/2}':=\max_{x\in V'_{N/2}}h^{V'_{N/2}}_x$, a small variation on Exercise~\ref{ex:3.4} shows
\begin{multline}
\label{E:8.83uai}
\qquad
P(M_N\ge m_N)
\ge P(M_{N/2}'> m_{N/2}-c)
\\
\times
\min_{x\in V'_{N/2}} P\bigl(\varphi^{V_N,V'_{N/2}}_x\ge c+m_N-m_{N/2}\bigr).
\qquad
\end{multline}
Since~$m_N-m_{N/2}$ remains bounded as~$N\to\infty$ and 
\begin{equation}
\label{E:8.72ueu}
\inf_{N\ge10}\,\min_{x\in V'_{N/2}}\Var\bigl(\varphi^{V_N,V'_{N/2}}_x\bigr)>0\,,
\end{equation}
the minimum on the right of \eqref{E:8.83uai} is positive uniformly in~$N\ge10$. The claim follows from \eqref{E:8.71ueu}.
\end{proofsect}

As our final step, we boost \eqref{E:8.65} to a bound on the lower tail:

\begin{mylemma}[Tightness of lower tail]
\label{lemma-8.21}
There is $a>0$ and~$t_0>0$ such that
\begin{equation}
\label{E:8.70}
\sup_{N\ge1}\,P\bigl(M_N<m_N-t\bigr)\le \texte^{-a t},\quad t>t_0.
\end{equation}
\end{mylemma}

Before we delve into the proof, let us remark that the bound \eqref{E:8.70} is not sharp even as far its overall structure is concerned. Indeed, Ding and Zeitouni~\cite{DZ} showed that the lower tails of~$M_N$ are in fact doubly exponential. However, the proof of the above is easier and fully suffices for our later needs.

\begin{proofsect}{Proof of Lemma~\ref{lemma-8.21}}
Consider the setting as in the proof of Lemma~\ref{lemma-8.2}. The reasoning leading up to \eqref{E:8.13} plus the arguments underlying \eqref{E:8.24uea} show
\begin{equation}
\label{E:bd-u}
P\bigl(\,M_{3KN}< m_{3KN}-t\bigr)
\le2\texte^{-a t}+\texte^{-\frac14 K^2P(\,M_N\ge m_{3NK}-t+\sqrt{t}\,\sqrt{\log K})}.
\end{equation}
Now link~$K$ to~$t$ by setting $t:=c\log K$ with~$c>0$ so large that $m_{3NK}-t+\sqrt{t}\,\sqrt{\log K}\le m_N$.
With the help of \eqref{E:8.65}, the second term on the right of \eqref{E:bd-u} is doubly exponentially small in~$t$, thus proving the claim (after an adjustment of~$a$) for~$N\in (3K)\N$. The claim follows from the next exercise.
\end{proofsect}

\begin{myexercise}
Use a variation of Exercise~\ref{ex:3.4} to show that if~$c$ relating~$K$ to~$t$ is chosen large enough, then \eqref{E:8.70} for~$N\in (3K)\N$ extends to a similar bound for all~$N\in\N$. Hint: Note that $\max_{0\le r<3K}[m_{3Kn+r}-m_{3Kn}]\le 2\sqrt g\,\log(3K)$.
\end{myexercise}

From here we finally get:

\begin{proofsect}{Proof of Theorem~\ref{thm-7.3}}
The tightness of $\{M_N-m_N\colon N\ge1\}$ follows from Lemmas~\ref{lemma-8.2} and~\ref{lemma-8.21}. Alternatively, use these lemmas to infer~$EM_N=m_N+O(1)$ and then apply Lemma~\ref{lemma-DH}.
\end{proofsect}


\chapter{Extremal local extrema}
\label{lec-9}\noindent
Having established tightness of the DGFF maximum, we turn our attention to the structure of the extremal level sets --- namely, the sets of vertices where the DGFF takes values within order unity of the absolute maximum. As for the intermediate levels, we cast these in the form of a point process. In this lecture we state the main result proved in a sequence of joint papers with O.~Louidor~\cite{BL1,BL2,BL3} and then discuss one of the main ingredients of the proof: Liggett's theory of invariant measures for point processes evolving (or Dysonized) by independent Markov chains. This characterizes the subsequential limits of the resulting point process; uniqueness (and thus existence) of the limit will be addressed later.

\section{Extremal level sets}
\noindent
In the previous lectures we showed that the maximum of the DGFF in~$V_N$ is tight around the sequence~$m_N$ defined in \eqref{E:7.8}. By Exercise~\ref{ex:3.4}, this applies to any sequence $\{D_N\colon N\ge1\}$ of admissible discretizations of a continuum domain~$D\in\mathfrak D$. Once the tightness of the maximum is in place, additional conclusions of interest can be derived concerning the structure of the \emph{extremal level set}
\begin{equation}
\label{E:9.1}
\Gamma^D_N(t):=\bigl\{x\in D_N\colon h^{D_N}_x\ge m_N-t\bigr\}.
\end{equation}
We will need two theorems, both proved originally in Ding and Zeitouni~\cite{DZ}:

\begin{mytheorem}[Size of extremal level set]
\label{thm-DZ1}
There are constants~$c,C\in(0,\infty)$ such that
\begin{equation}
\label{E:9.2nwt}
\lim_{t\to\infty}\liminf_{N\to\infty}\,P\bigl(\texte^{ct}\le|\Gamma^D_N(t)|\le\texte^{Ct}\bigr)=1.
\end{equation}
\end{mytheorem}

\begin{mytheorem}[Geometry of extremal level set]
\label{thm-DZ2}
For each~$t\in\R$,
\begin{equation}
\lim_{r\to\infty}\,\limsup_{N\to\infty}\,
P\Bigl(\exists x,y\in \Gamma^D_N(t)\colon r<|x-y|<\frac Nr\Bigr)=0\,.
\end{equation}
\end{mytheorem}

The upshot of these theorems is the following picture: \textit{With high probability, the extremal level set $\Gamma_N^D(t)$ is the union of a finite number of ``islands'' of bounded size at distances of order~$N$}. For this reason, we like to think of Theorems~\ref{thm-DZ1}--\ref{thm-DZ2} as tightness results; the former for the size and the latter for the spatial distribution of~$\Gamma_N^D(t)$. Our proofs of these theorems will follow a different route than~\cite{DZ} and are therefore postponed to Lecture~\ref{lec-12}. 

Our focus on scaling limits naturally suggests considering subsequential distributional limits of the scaled version of~$\Gamma_N^D(t)$, namely, the set $\{x/N\colon x\in\Gamma_N^D(t)\}$ regarded as a random subset of the continuum domain~$D$.
As is common in the area of extreme-order statistics, these are best encoded in terms of the empirical point measure on $D\times\R$ defined by
\begin{equation}
\label{E:eta-ND}
\eta_N^D:=\sum_{x\in D_N}\delta_{x/N}\otimes\delta_{h^{D_N}_x-m_N}\,.
\end{equation}
Indeed, we have
\begin{equation}
\eta^D_N\bigl(D\times[-t,\infty)\bigr)=\bigl|\Gamma_N^D(t)\bigr|
\end{equation}
 and so~$\eta^D_N$ definitely captures at least the size of the level set. In addition,~$\eta^D_N$ also keeps track of how the points are distributed in the level set and what the value of the field is there, very much like the measures discussed in Lectures~\ref{lec-2}--\ref{lec-5}. However, unlike for the intermediate level sets, no normalization needs to be imposed here thanks to Theorem~\ref{thm-DZ1}.

\nopagebreak
\begin{figure}[t]
\vglue-1mm
\centerline{\includegraphics[width=0.65\textwidth]{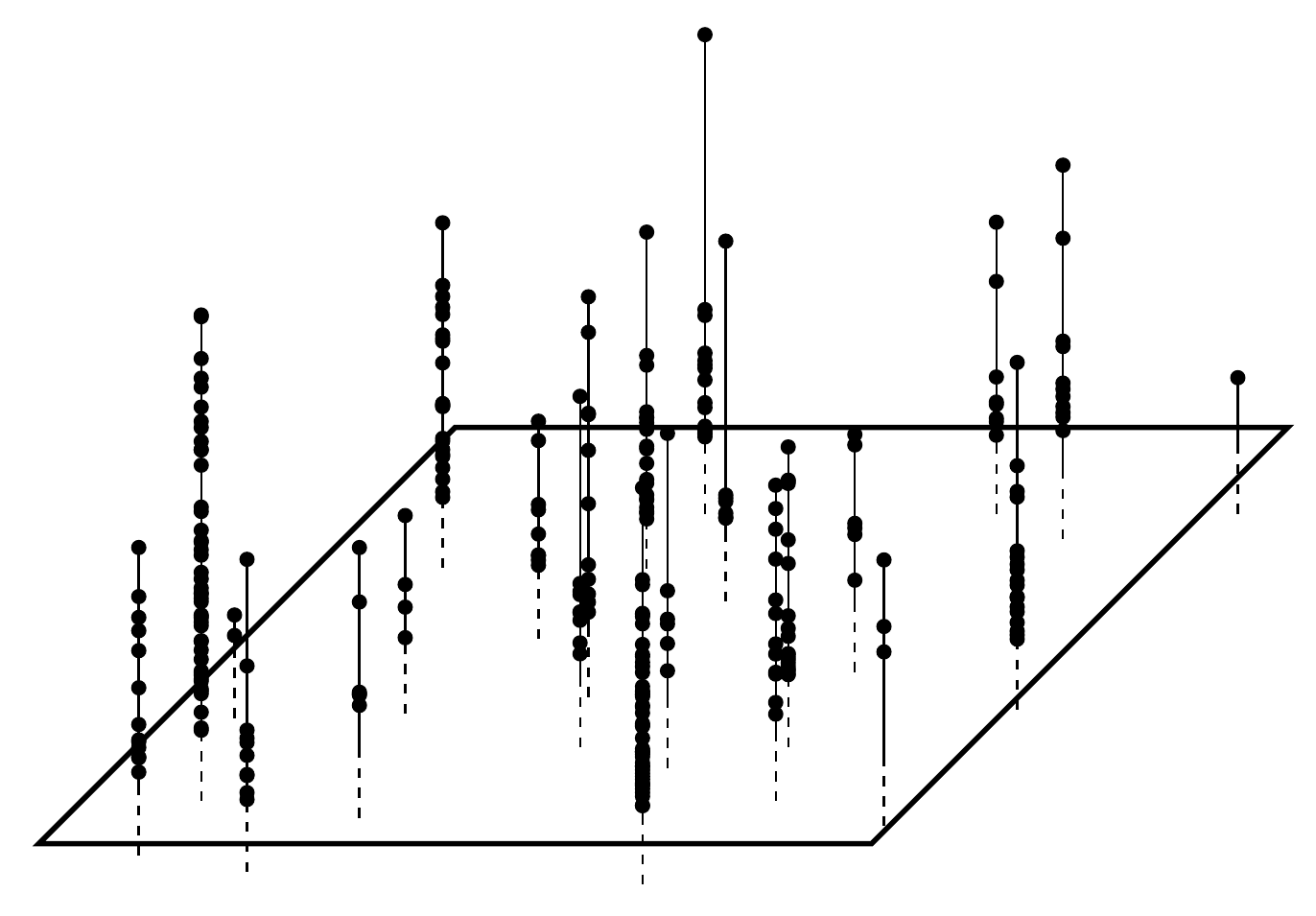}}
\vglue0mm
\begin{quote}
\small 
\caption{
\label{fig-cluster-process}
\small
A sample of the limit process obtained from~$\eta_N^D$ in \eqref{E:eta-ND} in the limit~$N\to\infty$. Each ``leading'' local maximum is accompanied by a cluster of points that, in the limit, all line up over the same spatial location.}
\normalsize
\end{quote}
\end{figure}

Notwithstanding all that was just said, $\eta^D_N$ are actually not the most natural \emph{primary} objects to work with. This is because each high value of the DGFF comes along with a whole \emph{cluster} of high values; namely, the values of the field at the points of~$\Z^2$ in close vicinity thereof. In the limit~$N\to\infty$, this whole cluster collapses to a single spatial location; see Fig.~\ref{fig-cluster-process}. As the field-values within a single cluster are heavily correlated (and remain so through $N\to\infty$ limit) depending on their relative lattice positions, it is advantageous to track them at the (original) lattice scale; i.e., without scaling space by~$N$.

Thus, denoting (with a slight abuse of our earlier notation) by
\begin{equation}
\Lambda_r(x):=\bigl\{y\in\Z^2\colon |x-y|<r\bigr\}
\end{equation}
the $r$-neighborhood of~$x$ in (say) Euclidean norm,
instead of~$\eta_N^D$ we will consider its \emph{structured} version
\begin{equation}
\label{E:9.6}
\eta_{N,r}^D:=\sum_{x\in D_N}\1_{\{h^{D_N}_x=\max_{y\in\Lambda_r(x)}h^{D_N}_y\}}\,
\delta_{x/N}\otimes\delta_{h^{D_N}_x-m_N}\otimes\delta_{\{h^{D_N}_x-h^{D_N}_{x+z}\colon z\in\Z^2\}}\,.
\end{equation}
This measure  picks one reference point in each ``cluster'' --- namely, the ``$r$-local maximum'' which is a point~$x$ where~$h^{D_N}$ dominates all values in~$\Lambda_r(x)$  --- and records the scaled position and reduced field-value at this point along with the ``shape'' of the field thereabout.

Let us write~$\PPP(\mu)$ to denote the Poisson point process with ($\sigma$-finite) intensity measure~$\mu$. Abbreviate $\overline\R:=\R\cup\{+\infty,-\infty\}$. Our main result on the structured extremal point processes \eqref{E:9.6} is the following theorem, derived jointly with O.~Louidor in~\cite{BL1,BL2,BL3}:

\begin{mytheorem}[DGFF extremal process]
\label{thm-extremal-vals}
There is a probability measure~$\nu$ on $[0,\infty)^{\Z^2}$ and, for each~$D\in\mathfrak D$, a random Borel measure~$Z^D$ on~$D$ with $Z^D(D)\in(0,\infty)$ a.s.\ such that the following holds for any sequence $\{D_N\colon N\ge1\}$ of admissible approximations of~$D$ and any sequence $\{r_N\colon N\ge1\}$ of numbers with~$r_N\to0$ and~$N/r_N\to\infty$:
\begin{equation}
\label{E:9.8uai}
\eta^D_{N,r_N}\,\,\underset{N\to\infty}\lawarrow\,\, \PPP\bigl(Z^D(\textd x)\otimes\texte^{-\alpha h}\textd h\otimes\nu(\textd\phi)\bigr),
\end{equation}
where the convergence in law is with respect to the vague convergence of Radon measures on~$\overline D\times(\R\cup\{\infty\})\times\overline\R^{\Z^d}$ and $\alpha:=2/\sqrt g$.
\end{mytheorem}

A couple of remarks are in order:
\begin{enumerate}
\item[(1)] By Exercise~\ref{ex:2.8nwt} with $\scrX:=\overline D\times(\R\cup\{\infty\})\times\overline\R^{\Z^d}$, the statement \eqref{E:9.8uai} is equivalent to the distributional convergence of the sequence of real-valued random variables $\{\langle\eta^D_{N,r_N},f\rangle\colon N\ge1\}$
 for all functions $f\in C_\cc(\scrX)$, or even just such~$f$'s that depend only on a finite number of the coordinates in the third (i.e.,``cluster'') variable.
\item[(2)] 
The $\PPP$ process with a random intensity is to be understood as a random sample from a class of Poisson processes. To generate $\PPP(\mu)$ for~$\mu$ random, we thus first sample~$\mu$ and then generate the Poisson process conditional on~$\mu$. The term \emph{Cox process} is sometimes used for this object as well.
\item[(3)] 
Somewhat surprisingly, the expression we saw in the limit of the intermediate level sets appears in \eqref{E:9.8uai} as well, albeit now as the intensity measure of a Poisson process. See Conjectures~\ref{cor-12.13ua} and~\ref{conj-15.5} for some explanation.
\item[(4)] Using local maxima as the reference points in \eqref{E:9.6} is the most natural choice although other choices might presumably work as well. (For instance, one could take~$x$ to be a sample from a suitably normalized measure $z\mapsto\texte^{\beta h^{D_N}_z}$ on~$\Lambda_r(y)$ for~$y$ being an $r$-local maximum of~$h^{D_N}$.)
\end{enumerate}

Let us elaborate on remark~(2) above by offering a more explicit way of sampling the limit process in \eqref{E:9.8uai}. Writing $Z^D(\cdot)$ as~$Z^D(D)\wh Z^D(\cdot)$, where
\begin{equation}
\label{E:10.18ueu}
\wh Z^D(A):=\frac{Z^D(A)}{Z^D(D)},
\end{equation}
the limit process in \eqref{E:9.8uai} can  be sampled as follows: Let $\{h_i'\colon i\in\N\}$ be points in a sample of the Poisson point process~$\R$ with (Gumbel) intensity $\texte^{-\alpha h}\textd h$. Given an independent sample of the~$Z^D$ measure, set
\begin{equation}
\label{E:9.11nwwt}
h_i:=h_i'+\alpha^{-1}\log Z^D(D),\quad i\in\N,
\end{equation}
and, conditionally on~$Z^D$, sample independently:
\settowidth{\leftmargini}{(1111)}
\begin{enumerate}
\item[(1)] $\{x_i\colon i\in \N\}\,\laweq$ i.i.d.\ with law~$\wh Z^D$,
\item[(2)] $\{\phi_i\colon i\in \N\}\,\laweq$  i.i.d.\ with law~$\nu$.
\end{enumerate}
The point measure $\sum_{i\in\N}\delta_{x_i}\otimes\delta_{h_i}\otimes\delta_{\phi_i}$ then has the law on the right-hand side of \eqref{E:9.8uai}. See again Fig.~\ref{fig-cluster-process}.

We remark that a conclusion analogous to \eqref{E:9.8uai} is known to hold also for the Branching Brownian Motion. This is thanks to the work of McKean~\cite{McKean}, Bramson~\cite{Bramson,Bramson2}, Lalley and Sellke~\cite{LS} culminating in Arguin, Bovier and Kistler~\cite{ABK1,ABK2,ABK3}, A\"idekon, Berestycki, Brunet and Shi~\cite{ABBS} and Bovier and Hartung~\cite{Bovier-Hartung}; see also the review by Bovier~\cite{Bovier-review}. The corresponding problem was solved for the DGFF in~$d\ge3$ as well (Chiarini, Cipriani and Hazra~\cite{CCH1,CCH2,CCH3}) although (since these DGFFs are no longer log-correlated) there the limit process has a non-random intensity. In fact, the same applies to the field of i.i.d.\ standard normals.

\section{Distributional invariance}
\noindent
The proof of Theorem~\ref{thm-extremal-vals} will take the total of four lectures. (The actual proof appears in Sections~\ref{sec:10.2} and~\ref{sec-full-convergence} with all needed technicalities settled only in Lecture~\ref{lec-12}.) The basic strategy will be very much like that used for the intermediate level sets although, chronologically, these ideas started with extremal level sets. Namely, we extract a subsequential limit of the processes of interest and then derive enough properties to identify the limit law uniquely.

\smallskip
We start by noting a tightness statement:

\begin{myexercise}
\label{ex:9.4}
Let $f\in C_\cc(\overline D\times(\R\cup\{\infty\})\times\overline\R^{\Z^2})$. Use Theorems~\ref{thm-DZ1}--\ref{thm-DZ2} to show that the sequence of random variables $\{\langle \eta^D_{N,r_N},f\rangle\colon N\ge1\}$ is tight. 
\end{myexercise}

\noindent
The use of subsequential convergence is then enabled by:

\begin{myexercise}
\label{ex:9.5}
Suppose~$\scrX$ is a locally compact, separable Hausdorff space and let~$\eta_N$ be a sequence of random (non-negative) integer-valued Radon measures on~$\scrX$. Assume $\{\langle\eta_N,f\rangle\colon N\ge1\}$ is tight for each~$f\in C_\cc(\scrX)$. Prove that there is a sequence $N_k\to\infty$ such that $\langle \eta_{N_k},f\rangle$ converges in law to a random variable that takes the form~$\langle\eta,f\rangle$ for some random Radon measure~$\eta$ on~$\scrX$. Prove that~$\eta$ is integer valued.
\end{myexercise}

\noindent

In this lecture we go only part of the way to the full proof of Theorem~\ref{thm-extremal-vals} by showing that any (subsequential) limit of the measures~$\{\eta^D_{N,r_N}\colon N\ge1\}$ must take the form on the right of \eqref{E:9.8uai}. We will work only with the first two coordinates of the process for a while so we will henceforth regard~$\eta^D_{N,r_N}$ as a random measure on~$D\times\R$. The main focus of this lecture is thus:

\begin{mytheorem}[Poisson structure of subsequential limits]
\label{prop-subseq}
Any weak subsequential limit~$\eta^D$ of the processes~$\{\eta^D_{N,r_N}\colon N\ge1\}$ restricted to just the first two coordinates takes the form
\begin{equation}
\eta^D\,\,\laweq\,\,\PPP\bigl(Z^D(\textd x)\otimes\texte^{-\alpha h}\textd h\bigr)
\end{equation}
for some random Borel measure~$Z^D$ on~$\overline D$ with~$Z^D(\overline D)\in(0,\infty)$ a.s.
\end{mytheorem}

As we will demonstrate, this arises from the fact that \emph{every} subsequential limit measure~$\eta^D$ has the following distributional symmetry:

\begin{myproposition}[Invariance under Dysonization]
\label{prop-Dyson}
For any~$\eta^D$ as above (projected on the first two coordinates) and any function $f\in C_\cc(\overline D\times\R)$, we have
\begin{equation}
E\bigl(\texte^{-\langle\eta,f\rangle}\bigr)=E\bigl(\texte^{-\langle\eta,f_t\rangle}\bigr),\quad t>0,
\end{equation}
where
\begin{equation}
f_t(x,h):=-\log E^0\bigl(\texte^{-f(x,h+B_t-\frac\alpha2 t)}\bigr)
\end{equation}
with $\{B_t\colon t\ge0\}$ denoting the standard Brownian motion.
\end{myproposition}

Let us pause to explain why we refer to this as ``invariance under Dysonization.'' Exercises~\ref{ex:9.4}-\ref{ex:9.5} ensure that~$\eta^D$ is a point measure, i.e.,
\begin{equation}
\eta^D=\sum_{i\in\N}\delta_{x_i}\otimes\delta_{h_i}.
\end{equation}
Given a collection $\{B_t^{(i)}\colon t\ge0\}_{i=1}^\infty$ of independent standard Brownian motions, we then set
\begin{equation}
\eta^D_t:=\sum_{i\in\N}\delta_{x_i}\otimes\delta_{h_i+B^{(i)}_t-\frac\alpha2t}\,.
\end{equation}
Of course, for~$t>0$ this may no longer be a ``good'' point measure as we cannot \emph{a priori} guarantee that~$\eta^D_t(C)<\infty$ for any compact set. Nonetheless, we can use Tonelli's theorem to perform the following computations:
\begin{equation}
\label{E:9.13}
\begin{aligned}
E\bigl(\texte^{-\langle\eta_t,f\rangle}\bigr) 
&= E\Bigl(\,\prod_{i\in\N}\texte^{-f(x_i,h_i+B^{(i)}_t-\frac\alpha2t)}\Bigr)
\\
&=E\Bigl(\,\prod_{i\in\N}\texte^{-f_t(x_i,h_i)}\Bigr)
=E\bigl(\texte^{-\langle\eta,f_t\rangle}\bigr)\,,
\end{aligned}
\end{equation}
where in the middle equality we used the Bounded Convergence Theorem to pass the expectation with respect to each Brownian motion inside the infinite product. Proposition~\ref{prop-Dyson} then tells us
\begin{equation}
E\bigl(\texte^{-\langle\eta_t,f\rangle}\bigr) =E\bigl(\texte^{-\langle\eta,f\rangle}\bigr),\quad t\ge0,
\end{equation}
and, since this holds for all~$f$ as above,
\begin{equation}
\eta_t\,\,\laweq\,\,\eta,\quad t\ge0.
\end{equation}
Thus, attaching to each point of~$\eta^D$ an independent diffusion~$t\mapsto B_t-\frac\alpha2t$ in the second coordinate preserves the law of~$\eta^D$. We call this operation ``Dysonization'' in analogy to Dyson's proposal~\cite{Dyson} to consider dynamically-evolving random matrix ensembles; i.e., those where the static matrix entries are replaced by stochastic processes.

\begin{proofsect}{Proof of Proposition~\ref{prop-Dyson} (main idea)}
The proof is based on the following elementary computation.
Writing $h$ for the DGFF on~$D_N$, let $h',h''$ be independent copies of~$h$. For any~$s\in[0,1]$, we can then realize~$h$ as
\begin{equation}
h = \sqrt{1-s}\,\,h'+\sqrt s\,h''.
\end{equation}
Choosing~$s:=t/(g\log N)$, we get
\begin{equation}
\label{E:9.17}
h = \sqrt{1-\frac{t}{g\log N}\,}\,\,h'+\frac{\sqrt t}{\sqrt{g\log N}\,}\,h''\,.
\end{equation}
Now pick~$x$ at or near a local maximum of~$h'$ of value~$m_N+O(1)$. Expanding the first square-root into the first-order Taylor polynomial and using that $h'_x=2\sqrt g\,\log N+O(\log\log N)$ yields
\begin{equation}
\begin{aligned}
h_x&=h'_x-\frac12\frac t{g\log N}\,h'_x+\frac{\sqrt t}{\sqrt{g\log N}\,}\,h''_x+O\Bigl(\frac1{\log N}\Bigr)
\\
&=h'_x-\frac t{\sqrt g}+\frac{\sqrt t}{\sqrt{g\log N}\,}\,h''_x+O\Bigl(\frac{\log\log N}{\log N}\Bigr)\,.
\end{aligned}
\end{equation}
As noted in \eqref{E:1.22uai}, the covariance of the second field on the right of~\eqref{E:9.17} satisfies
\begin{equation}
\text{Cov}\Bigl(\frac{\sqrt t}{\sqrt{g\log N}\,}\,h''_x,\,\frac{\sqrt t}{\sqrt{g\log N}\,}\,h''_y\Bigr)
=\begin{cases}
t+o(1),\quad&\text{if }|x-y|<r_N,
\\
o(1),\quad&\text{if }|x-y|>N/r_N.
\end{cases}
\end{equation}
Thus, the second field behaves as a \emph{random constant} (with the law of~$\NN(0,t)$) on the whole island of radius~$r_N$ around~$x$, and these random constants on distinct islands can be regarded as more or less independent.

The technical part of the proof (which we skip) requires showing that, if~$x$ is a local maximum of~$h'_x$ at height~$m_N+O(1)$, the \emph{gap} between~$h'_x$ and the second largest value in an~$r_N$-neighborhood of~$x$ stays positive with probability tending to~$1$ as~$N\to\infty$. Once this is known, the errors in all approximations can be fit into this gap and~$x$ will also be the local maximum of the field~$h$ --- and so the locations of the relevant local maxima of~$h$ and~$h'$ coincide with high probability. 

The values of the field $\frac{\sqrt t}{\sqrt{g\log N}}\,h''_x$ for~$x$ among the local maxima of~$h$ (or~$h'$) then act as independent copies of~$\NN(0,t)$, and so they can be realized as values of independent Brownian motions at time~$t$. The drift term comes from the linear shift by~$t/\sqrt g$ by noting that~$1/\sqrt g=\alpha/2$.  
\end{proofsect}

The full proof of Proposition~\ref{prop-Dyson} goes through a series of judicious approximations for which we refer the reader to~\cite{BL1}. One of them addresses the fact that~$f_t$ no longer has compact support. We pose this as:

\begin{myexercise}
Show that for each $f\in C_\cc(\overline D\times\R)$, each~$t>0$ and each~$\epsilon>0$ there is $g\in C_\cc(\overline D\times\R)$ such that
\begin{equation}
\limsup_{N\to\infty}\,P\Bigl(\bigl|\langle \eta_{N,r_N}^D,f_t\rangle-\langle \eta_{N,r_N}^D,g\rangle\bigr|>\epsilon\Bigr)<\epsilon.
\end{equation}
Hint: Use Theorem~\ref{thm-DZ1}.
\end{myexercise}

\noindent
We will also revisit the ``gap estimate'' in full detail later; see Exercise~\ref{ex:11.9}.

\section{Dysonization-invariant point processes}
\noindent
Our next task is to use the distributional invariance articulated in Proposition~\ref{prop-Dyson} to show that the law of~$\eta^D$ must take the form in Theorem~\ref{prop-subseq}. A first natural idea is to look at the moments of~$\eta^D$, i.e., measures~$\mu_n$ on~$\R^n$ of the form
\begin{equation}
\mu_n(A_1\times\dots\times A_n):=E\bigl[\eta^D(D\times A_1)\dots\eta^D(D\times A_n)\bigr].
\end{equation}
As is readily checked, by the smoothness of the kernel associated with $t\mapsto B_t-\frac\alpha2t$, these would have a density $f(x_1,\dots,x_n)$ with respect to the Lebesgue measure on~$\R^n$ and this density would obey the PDE
\begin{equation}
\alpha (1,\dots,1)\cdot\nabla f +\Delta f = 0.
\end{equation}
This seems promising since all full-space solutions of this PDE can be classified. Unfortunately, the reasoning fails right at the beginning because, as it turns out, all positive-integer moments of~$\eta^D$ are infinite. We thus have to proceed using a different argument. Fortunately, a 1978 paper of T.~Liggett~\cite{Liggett-ZWVG} does all what we need to do, so this is what we will discuss next.

Liggett's interest in \cite{Liggett-ZWVG} was in \emph{non-interacting} particle systems which he interpreted as families of independent Markov chains. His setting is as follows: Let~$\scrX$ be a locally compact, separable Hausdorff space and, writing~$\BB(\scrX)$  for the class of Borel sets in~$\scrX$, let $\cmss P\colon\scrX\times\BB(\scrX)\to[0,1]$ be a transition kernel of a Markov chain on~$\scrX$. Denoting $\N^\star:=\N\cup\{0\}\cup\{\infty\}$, we are interested in the evolution of point processes on~$\scrX$, i.e., random elements of
\begin{equation}
\MM:=\bigl\{\text{$\N^\star$-valued Radon measures on~$(\scrX,\BB(\scrX))$}\bigr\},
\end{equation}
under the (independent) Markovian ``dynamics.''

Recall that being ``Radon'' means that the measure is inner and outer regular with a finite value on every compact subset of~$\scrX$. Calling the measures in~$\MM$ ``point processes'' is meaningful in light of:

\begin{myexercise}
Every~$\eta\in\MM$ takes the form
\begin{equation}
\label{E:9.22}
\eta = \sum_{i=1}^N\delta_{x_i},
\end{equation}
where $N\in\N^\star$ and where (for $N>0$) $\{x_i\colon i=1,\dots,N\}$ is a multiset of points from~$\scrX$ such that $\eta(C)<\infty$ a.s.\ for every~$C\subset\scrX$ compact.
\end{myexercise}

To describe the above ``dynamics'' more precisely, given~$\eta\in\MM$ of the form \eqref{E:9.22}, we consider a collection of independent Markov chains $\{X_n^{(i)}\colon n\ge0\}_{i=1}^N$ such that
\begin{equation}
P\bigl(X_0^{(i)}=x_i\bigr)=1,\quad i=1,\dots,N.
\end{equation}
Then we define
\begin{equation}
\eta_n:=\sum_{i=1}^N\delta_{X_n^{(i)}},\quad n\ge 0,
\end{equation}
while noting that, due to the independence of the chains $X^{(i)}$, the law of this object does not depend on the initial labeling of the atoms of~$\eta$.

Note that we cannot generally claim that~$\eta_n\in\MM$ for~$n\ge1$; in fact, easy counterexamples can be produced to the contrary whenever~$N=\infty$. Nonetheless, for processes whose law is \emph{invariant} under the above time evolution $\eta\mapsto\{\eta_n\colon n\ge0\}$, i.e., those in
\begin{equation}
\label{E:I}
\scrI:=\bigl\{\eta\colon \text{random element of~$\MM$ such that~$\eta_1\,\,\laweq\,\,\eta$}\bigr\},
\end{equation}
this is guaranteed \emph{ex definitio} in light of 
\begin{equation}
\eta\in\MM\quad\text{and}\quad\eta_1\,\,\laweq\,\,\eta\quad\Rightarrow\quad \eta_n\in\MM \text{ a.s. } \forall n\ge0.
\end{equation}
Let~$\cmss P^n$ denote the $n$-th (convolution) power of the kernel~$\cmss P$ defined, inductively, by $\cmss P^{n+1}(x,\cdot)=\int_{\scrX}\cmss P(x,\textd y)\cmss P^n(y,\cdot)$. 
The starting point of Liggett's paper is the following observation which, although attributed to folklore, is still quite ingenious:

\begin{mytheorem}[Folklore theorem]
\label{thm-Liggett-folk}
Suppose~$\cmss P$ obeys the following ``strong dispersivity'' assumption
\begin{equation}
\label{E:9.27}
\forall C\subset\scrX\text{\rm\ compact}\colon\quad
\sup_{x\in\scrX}\cmss P^n(x,C)\,\,\underset{n\to\infty}\longrightarrow\,\,0\,.
\end{equation}
Then the set of invariant measures $\scrI$ defined in \eqref{E:I} is given as
\begin{equation}
\scrI=\bigl\{\PPP(M)\colon M = \text{\rm (random) Radon measure on~$\scrX$ s.t. $M\,\,\laweq\,\,M\cmss P$}\bigr\},
\end{equation}
where $M\cmss P(\cdot):=\int M(\textd x)\cmss P(x,\cdot)$.
\end{mytheorem}

We refer to the condition \eqref{E:9.27} using the physics term ``dispersivity'' as that is often used to describe the rate of natural spread (a.k.a.\ dispersion) of waves or wave packets in time. Note that this condition rules out existence of stationary \emph{distributions} for the Markov chain.

Before we prove Theorem~\ref{thm-Liggett-folk}, let us verify its easier part --- namely, that all the Poisson processes~$\PPP(M)$ with $M\,\,\laweq\,\,M\cmss P$ lie in~$\scrI$. This follows from:

\begin{mylemma}
\label{lemma-9.9}
Let~$M$ be a Radon measure on~$\scrX$. Then
\begin{equation}
\eta\,\,\laweq\,\,\PPP(M)\quad\Rightarrow\quad \eta_n\,\, \laweq\,\,\PPP(M\cmssP^n),\quad\forall n\ge0.
\end{equation}
\end{mylemma}

\begin{proofsect}{Proof}
We start by noting that $\eta$ being a sample from~$\PPP(M)$ for a given Radon measure~$M$ is equivalent to saying that, for every~$f\in C_\cc(\scrX)$ with~$f\ge0$,
\begin{equation}
\label{E:9.30}
E\bigl(\texte^{-\langle\eta,f\rangle}\bigr) = \exp\Bigl\{-\int M(\textd x)(1-\texte^{-f(x)})\Bigr\}.
\end{equation}
The argument in \eqref{E:9.13} shows
\begin{equation}
\label{E:9.31}
E\bigl(\texte^{-\langle\eta_n,f\rangle}\bigr)=E\bigl(\texte^{-\langle\eta,f_n\rangle}\bigr)\,,
\end{equation}
where 
\begin{equation}
\label{E:9.36uea}
f_n(x):=-\log \bigl[(\cmss P^n\texte^{-f})(x)\bigr].
\end{equation}
 This equivalently reads as~$\texte^{-f_n}=\cmss P^n\texte^{-f}$ and so from \eqref{E:9.30} for~$f$ replaced by~$f_n$ and the fact that~$\cmss P^n1=1$ we get
\begin{equation}
E\bigl(\texte^{-\langle\eta,f_n\rangle}\bigr) = \exp\Bigl\{-\int M(\textd x)\,\cmss P^n(1-\texte^{-f})(x)\Bigr\}.
\end{equation}
Tonelli's theorem identifies the right-hand side with that of \eqref{E:9.30} for~$M$ replaced by~$M\cmss P^n$. Since \eqref{E:9.30} characterizes $\PPP(M)$, the claim follows from \eqref{E:9.31}.
\end{proofsect}

We are now ready to give:

\begin{proofsect}{Proof of Theorem~\ref{thm-Liggett-folk}}
In light of Lemma~\ref{lemma-9.9} we only need to verify that every element of~$\scrI$ takes the form~$\PPP(M)$ for some $M$ satisfying $M\,\,\laweq\,\,M\cmss P$. Let~$\eta\in\scrI$ and pick $f\in C_\cc(\scrX)$ with~$f\ge0$. Since~$\texte^{-f}$ equals~$1$ outside a compact set, the strong dispersivity condition \eqref{E:9.27} implies
\begin{equation}
\sup_{x\in\scrX}\bigl|(\cmss P^n\texte^{-f})(x)-1\bigr|\,\underset{n\to\infty}\longrightarrow\,0.
\end{equation}
Recalling the definition of~$f_n(x)$ from \eqref{E:9.36uea}, this and~$\cmss P^n1=1$ yields the existence of an~$\epsilon_n\downarrow0$ such that
\begin{equation}
(1-\epsilon_n)\cmss P^n(1-\texte^{-f})(x)\le f_n(x)\le (1+\epsilon_n)\cmss P^n(1-\texte^{-f})(x),\quad x\in\scrX.
\end{equation}
By the Intermediate Value Theorem there is a (random) $\wt\epsilon_n\in[-\epsilon_n,\epsilon_n]$ such that
\begin{equation}
\langle \eta,f_n\rangle=(1+\wt\epsilon_n)\bigl\langle\eta,\cmss P^n(1-\texte^{-f})\bigr\rangle.
\end{equation}
Denoting $\eta\cmss P^n(\cdot):=\int\eta(\textd x)\cmss P^n(x,\cdot)$, from $\eta_n\,\,\laweq\,\,\eta$ and \eqref{E:9.31} we get
\begin{equation}
\label{E:9.36}
E\bigl(\texte^{-\langle\eta,f\rangle}\bigr) = E\bigl(\texte^{-(1+\wt\epsilon_n)\langle\eta\cmss P^n,(1-\texte^{-f})\rangle}\bigr)\,.
\end{equation}
Noting that every~$g\in C_\cc(\scrX)$ can be written as $g=\lambda(1-\texte^{-f})$ for some $f\in C_\cc(\scrX)$ and $\lambda>0$, we now ask the reader to solve:

\begin{myexercise}
Prove that $\{\langle \eta\cmss P^n,g\rangle\colon n\ge0\}$ is tight for all~$g\in C_\cc(\scrX)$.
\end{myexercise} 

Using this along with Exercise~\ref{ex:9.5} permits us to extract a subsequence $\{n_k\colon k\ge1\}$ with~$n_k\to\infty$ and a random Borel measure~$M$ on~$\scrX$ such that
\begin{equation}
\label{E:9.42ue}
\langle \eta\cmss P^{n_k},g\rangle\,\,\underset{k\to\infty}\lawarrow\,\,\langle M,g\rangle,\quad g\in C_\cc(\scrX).
\end{equation}
As~$\wt\epsilon_n\to0$ in~$L^\infty$, from \twoeqref{E:9.36}{E:9.42ue} we conclude
\begin{equation}
\label{E:9.37}
E\bigl(\texte^{-\langle\eta,f\rangle}\bigr) = E\bigl(\texte^{-\langle M,(1-\texte^{-f})\rangle}\bigr),\quad f\in C_\cc(\scrX).
\end{equation}
A comparison with \eqref{E:9.30} proves that $\eta\,\,\laweq\,\,\PPP(M)$. Replacing~$n$ by~$n+1$ in \eqref{E:9.36} shows that \eqref{E:9.37} holds with~$M$ replaced by~$M\cmss P$. From \eqref{E:9.37} we infer that~$M\cmss P$ is equidistributed to~$M$.
\end{proofsect}

\section{Characterization of subsequential limits}
\noindent
By Proposition~\ref{prop-Dyson}, for the case at hand --- namely, the problem of the extremal process of the DGFF restricted to first two coordinates --- the relevant setting is $\scrX:=\overline D\times\R$ and, fixing any~$t>0$, the transition kernel
\begin{equation}
\cmss P_t\bigl((x,h),A\bigr):=P^0\bigl((x,h+B_t-\tfrac\alpha2 t)\in A\bigr).
\end{equation}
We leave it to the reader to verify:

\begin{myexercise}
\label{ex:disperse}
Prove that (for all~$t>0$) the kernel $\cmss P_t$ has the strong dispersivity property~\eqref{E:9.27}.
\end{myexercise}

Hence we get:

\begin{mycorollary}
\label{cor-9.12}
Every subsequential limit~$\eta^D$ of processes $\{\eta_{N,r_N}^D\colon N\ge1\}$ (projected on the first two coordinates and taken with respect to the vague topology on Random measures on~$\overline D\times\R$) takes the form
\begin{equation}
\eta^D\,\,\laweq\,\,\PPP(M),
\end{equation}
where $M=M(\textd x\textd h)$ is a Radon measure on~$\overline D\times\R$ such that
\begin{equation}
\label{E:9.40}
M\cmss P_t\,\,\laweq\,\, M,\quad t>0.
\end{equation}
\end{mycorollary}

\begin{proofsect}{Proof}
Just combine Exercise~\ref{ex:disperse} with Theorem~\ref{thm-Liggett-folk}.
\end{proofsect}

Moving back to the general setting of the previous section, Theorem~\ref{thm-Liggett-folk} reduces the problem of classifying the invariant measures on point processes evolving by independent Markov chains to a question involving a \emph{single} Markov chain only:
\begin{equation*}
\text{Characterize (random) Radon measures~$M$ on~$\scrX$ satisfying $M\cmss P\,\,\laweq\,\, M$.}
\end{equation*}
Note that if~$M$ is a random sample from the invariant measures for the Markov chain~$\cmss P$, then $M\cmss P=M$ a.s.\ and so, in particular, $M\cmss P\,\,\laweq\,\, M$. This suggests that we recast the above question as:
\begin{equation*}
\text{When does }M\cmss P\,\,\laweq\,\, M\text{ imply } M\cmss P =M\text{ a.s.?}
\end{equation*}
In his 1978 paper, Liggett~\cite{Liggett-ZWVG} identified a number of examples when this is answered in the affirmative. The salient part of his conclusions is condensed~into:

\begin{mytheorem}[Cases when $M\cmss P\,\,\laweq\,\, M$ implies $M\cmss P =M$ a.s.]
\label{thm-9.13}
Let~$M$ be a random Radon measure on~$\scrX$. Suppose that
\settowidth{\leftmargini}{(11)}
\begin{enumerate}
\item[(1)] either $\cmss P$ is an irreducible, aperiodic, Harris recurrent Markov chain, or
\item[(2)] $\cmss P$ is a random walk on an Abelian group such that~$\cmss P(0,\cdot)$, where~$0$ is the identity, is not supported on a translate of a proper closed subgroup.
\end{enumerate}
Then $M\cmss P\,\,\laweq\,\, M$ implies $M\cmss P=M$ a.s.
\end{mytheorem}

Liggett's proof of Theorem~\ref{thm-9.13} uses sophisticated facts from the theory of Markov chains and/or random walks on Abelian groups and so we will not reproduce it in full generality here. Still, the reader may find it instructive to~solve:

\begin{myexercise}
Find examples of Markov chains that miss just one of the three required attributes in~(1) and for which the conclusion of the theorem fails.
\end{myexercise}

Returning back to the specific case of extremal point process associated with the DGFF, here the first alternative in Theorem~\ref{thm-9.13} does not apply as our Markov chain --- namely, the Brownian motion with a constant negative drift evaluated at integer multiplies of some~$t>0$ --- is definitely not Harris recurrent. Fortunately, the second alternative does apply and hence we get:

\begin{mycorollary}
For each~$t>0$, any~$M$ satisfying \eqref{E:9.40} obeys $M\cmss P_t=M$ a.s.
\end{mycorollary}

\noindent
We will provide an independent proof of this corollary, and thus also highlight the main ideas behind the second part of Theorem~\ref{thm-9.13}, by showing:

\begin{mylemma}
\label{lemma-9.17}
Any~$M$ satisfying \eqref{E:9.40} takes the form
\begin{equation}
\label{E:9.42}
M(\textd x\textd h)\,\,\laweq\,\,Z^D(\textd x)\otimes\texte^{-\alpha h}\textd h+\wt Z^D(\textd x)\otimes\textd h,
\end{equation}
where~$(Z^D,\wt Z^D)$ is a pair of random Borel measures on~$\overline D$.
\end{mylemma}

\begin{proofsect}{Proof}
Let~$M$ be a random Radon measure on~$\overline D\times\R$ with $M\cmss P_t \laweq M$ and let $A\subseteq\overline D$ be a Borel set. Since the Markov kernel~$\cmss P_t$ does not move the spatial coordinate of the process, we project that coordinate out by considering the measure
\begin{equation}
M^A(C):=M(A\times C),\quad C\subseteq\R\text{ Borel},
\end{equation}
and the kernel
\begin{equation}
\cmss Q_t(h,C):=P^0\bigl(h+B_t-\tfrac\alpha2 t\in C\bigr).
\end{equation}
Note that the only invariant sets for~$\cmss Q_t$ are~$\emptyset$ and~$\R$ and so~$\cmss Q_t$ falls under alternative~(2) in Theorem~\ref{thm-9.13}. We will nevertheless give a full proof in this~case.  

By assumption, the sequence $\{M^A\cmss Q_t^n\colon n\ge0\}$ is stationary with respect to the left shift. The Kolmogorov Extension Theorem (or weak-limit arguments) then permits us to embed it into a \emph{two-sided} stationary sequence $\{M^A_n\colon n\in\Z\}$ such that
\begin{equation}
\label{E:9.46uai}
M^A_n\,\,\laweq\,\,M^A\quad\text{and}\quad M^A_{n+1}=M_n^A\cmss Q_t \text{ a.s.}
\end{equation}
hold for each~$n\in\Z$. This makes~$\{M^A_n\colon n\in\Z\}$ a random instance of:

\begin{mydefinition}[Entrance law]
Given a Markov chain on~$S$ with transition kernel~$\cmss P$, a family of measures $\{\pi_n\colon n\in\Z\}$ on~$S$ satisfying 
\begin{equation}
\pi_{n+1}=\pi_n\cmss P,\quad n\in\Z.
\end{equation}
is called an \emph{entrance law}.
\end{mydefinition}

\noindent
From~$\cmss Q_t(h,\cdot)\ll\Leb$ we infer $M_n^A\ll\Leb$ for each~$n\in\Z$. Hence, there are (random) densities $h\mapsto f(n,h)$ such that
$M^A_n(\textd h)=f(n,h)\textd h$. Denoting by
\begin{equation}
k_t(h):=\frac1{\sqrt{2\pi t}}\,\texte^{-\frac{h^2}{2t}}
\end{equation}
the probability density of~$\NN(0,t)$, the second relation in \eqref{E:9.46uai} then reads
\begin{equation}
\label{E:9.46}
f(n+1,h)=f(n,\cdot)\star k_t\bigl(h+\tfrac\alpha2t\bigr)\,,
\end{equation}
where~$\star$ denotes the convolution of functions over~$\R$. Strictly speaking, this relation holds only for Lebesgue a.e.~$h$ but this can be mended by solving:

\begin{myexercise}
Show that, for each~$n\in\N$, $h\mapsto f(n,h)$ admits a continuous version such that $f(n,h)\in(0,\infty)$ for each~$h\in\R$. The identity \eqref{E:9.46} then holds for all~$h\in\R$.
\end{myexercise}

A key observation underlying Liggett's proof of~(2) in Theorem~\ref{thm-9.13} is then:

\begin{myexercise}
Let $X$ be a random variable on an Abelian group~$S$. Prove that every entrance law $\{\pi_n\colon n\in\Z\}$ for the random walk on~$S$ with step distribution~$X$ is a stationary measure for the random walk on $\Z\times S$ with step distribution $(1,X)$. 
\end{myexercise}

\noindent
This is relevant because it puts us in the setting to which the \emph{Choquet-Deny theory} applies (see Choquet and Deny~\cite{CD} or Deny~\cite{Deny}) with the following conclusion: Every $f\colon\Z\times\R\to(0,\infty)$ that obeys \eqref{E:9.46} takes the form
\begin{equation}
\label{E:CD1}
f(n,h)=\int \lambda(\kappa)^n\,\texte^{\kappa h}\,\nu(\textd \kappa)
\end{equation}
for some Borel measure~$\nu$ on~$\R$ and for
\begin{equation}
\label{E:CD2}
\lambda(\kappa):=\texte^{\frac12\kappa(\kappa+\alpha)t}\,.
\end{equation}
One important ingredient that goes in the proof is:

\begin{myexercise}
Prove that the random walk on (the Abelian group)~$\Z\times\R$ with step distribution $(1,\NN(t,-\frac\alpha2t))$ has no non-trivial closed invariant subgroup.
\end{myexercise}

Underlying the Choquet-Deny theory is Choquet's Theorem, which states that every compact, convex subset of a Banach space is the closed convex hull of its extreme points --- i.e., those that cannot be written as a convex combination of distinct points from the set. It thus suffices to classify the extreme points: 

\begin{myexercise}
Observe that if~$f$ lies in
\begin{equation}
\CC:=\Bigl\{f\colon \Z\times\R\to(0,\infty)\colon \eqref{E:9.46}\text{ holds},\, f(0,0)=1\Bigr\},
\end{equation}
then
\begin{equation}
f_s(n,h):=\frac{f(n-1,h+s)}{f(-1,s)}
\end{equation}
obeys $f_s\in\CC$ for each~$s\in\R$. Use this to prove that every extreme $f\in\CC$ takes the form $f(n,h)=\lambda(\kappa)^n\texte^{\kappa h}$ for some~$\kappa\in\R$ and~$\lambda(\kappa)$ as in \eqref{E:CD2}. 
\end{myexercise}

\noindent
With \twoeqref{E:CD1}{E:CD2} in hand, elementary integration (and Tonelli's theorem) gives
\begin{equation}
\label{E:9.56ueu}
M^A_n\bigl([-1,1]\bigr)=\int_{[-1,1]} f(n,h)\,\textd h = \int \lambda(\kappa)^n\,\frac{\sinh(\kappa)}\kappa\,\nu(\textd \kappa)\,,
\end{equation}
where, in light of \eqref{E:CD1},~$\nu$ is determined by (and thus a function of) the realization of~$M^A$.
Since $\{M^A_n\colon n\in\Z\}$ is stationary, the last integral in \eqref{E:9.56ueu} cannot diverge to infinity as $n\to\infty$ or $n\to-\infty$. This means that $\nu$ must be supported on~$\{\kappa\in\R\colon\lambda(\kappa)=1\}$ which, via \eqref{E:CD2}, forces~$\nu$ to be of the form
\begin{equation}
\nu = X^A\delta_{-\alpha}+Y^A\delta_0
\end{equation}
for some non-negative~$X^A$ and~$Y^A$. Hence we get
\begin{equation}
\label{E:9.52}
M^A(\textd h) = X^A\texte^{-\alpha h}\textd h+Y^A\textd h.
\end{equation}
But $A\mapsto M^A$ is a Borel measure and so $Z^D(A):=X^A$ and $\wt Z^D(A):=Y^A$
defines two random Borel measures for which \eqref{E:9.42} holds.
\end{proofsect}

\begin{myexercise}
A technical caveat in the last argument is that \eqref{E:9.52} holds only a.s.\ with the null set depending possibly on~$A$. Prove that, thanks to the fact that the Borel sets in~$\R$ are countably generated, the conclusion does hold as stated. (Remember to resolve all required measurability issues.).
\end{myexercise}

We are now ready to establish the desired Poisson structure for any weak subsequential limit of the processes $\{\eta_{N,r_N}^D\colon N\ge1\}$:

\begin{proofsect}{Proof of Theorem~\ref{prop-subseq}}
First we note that, if $f\in C_\cc(\overline D\times(\R\cup\{\infty\}))$ is such that $f\ge0$ and $\supp(f)\subseteq\overline D\times[t,\infty]$, then
\begin{equation}
\label{E:9.65nwwt}
\bigl\{\max_{x\in D_N}h^{D_N}_x<m_N+t\bigr\}=\Bigl\{\eta_{N,r_N}^D\bigl(D\times[t,\infty)\bigr)=0\Bigr\}
\end{equation}
implies
\begin{equation}
P\bigl(\langle\eta^D_{N,r_N},f\rangle>0\bigr)\le P\Bigl(\,\max_{x\in D_N}h^{D_N}_x\ge m_N+t\Bigr).
\end{equation}
The upper-tail tightness of the centered maximum (cf Lemma~\ref{lemma-8.2}) shows that the right-hand side tends to zero as~$N\to\infty$ and~$t\to\infty$. As a consequence we get that every subsequential limit of $\{\eta^D_{N,r_N}\colon N\ge1\}$ (in the first two coordinates) is concentrated on~$\overline\D\times\R$ and so may focus on convergence in this space.

Consider a sequence $N_k\to\infty$ such that $\eta^D_{N_k,r_{N_k}}$ converges in law with respect to the vague topology on the space of Radon measures on~$\overline D\times\R$. By Corollary~\ref{cor-9.12} and Lemma~\ref{lemma-9.17}, the limit is then of the form $\PPP(M)$ for some Radon measure~$M$ on~$\overline D\times\R$ of the form \eqref{E:9.42}. The probability that $\PPP(M)$ has no points in the set~$A$ is $\texte^{-M(A)}$. Standard approximation arguments and \eqref{E:9.65nwwt} then show
\begin{equation}
\label{E:9.53}
P\Bigl(\,\max_{x\in D_{N_k}}h^{D_{N_k}}_x<m_{N_k}+t\Bigr)\,\,\underset{k\to\infty}\longrightarrow\,\,E\Bigl(\texte^{-M(\overline D\times[t,\infty))}\Bigr)\,.
\end{equation}
(The convergence \emph{a priori} holds only for a dense set of~$t$'s but, since we already know that~$M$ has a density in the~$h$ variable, it extends to all~$t$.)
The upper-tail tightness of the maximum bounds the left-hand side by~$1-\texte^{-\tilde at}$ from below once~$t>t_0$. Taking~$t\to\infty$ thus yields
\begin{equation}
M\bigl(\overline D\times[t,\infty)\bigr)\,\,\underset{t\to\infty}\longrightarrow\,\,0\text{ a.s.}
\end{equation}
which forces~$\wt Z^D(\overline D)=0$ a.s. 

To get that~$Z^D(\overline D)>0$ a.s.\ we instead invoke the lower-tail tightness of the maximum (cf Lemma~\ref{lemma-8.21}), which bounds the left-hand side of \eqref{E:9.53} by~$\texte^{at}$ from above once~$t$ is large negative. The right-hand side in turn equals
\begin{equation}
E\bigl(\texte^{-\alpha^{-1}\texte^{-\alpha t}Z^D(\overline D)}\bigr),
\end{equation}
which tends to zero as~$t\to-\infty$ only if $Z^D(\overline D)>0$ a.s.
\end{proofsect}

We will continue the proof of Theorem~\ref{thm-extremal-vals} in the upcoming lectures. To conclude the present lecture we note that, although the subject of entrance laws has been studied intensely (e.g., by Cox~\cite{Cox}, the aforementioned paper of Liggett~\cite{Liggett-ZWVG} and also the diploma thesis of Secci~\cite{Secci}), a number of interesting open questions remains; particularly, for transient Markov chains. Ruzmaikina and Aizenman~\cite{AR} applied Liggett's theory to classify quasi-stationary states for competing particle systems on the real line where the limiting distribution ends up to be also of Gumbel extreme-order type.


\chapter{Nailing the intensity measure}
\label{lec-10}\noindent
In this lecture we continue the proof of extremal point process convergence in Theorem~\ref{thm-extremal-vals}. We first observe that the intensity measure associated with a subsequential limit of the extremal point process is in close correspondence with the limit distribution of the DGFF maximum. The existence of the latter limit is supplied by a theorem of Bramson, Ding and Zeitouni~\cite{BDingZ} which we will prove in Lecture~\ref{lec-12}, albeit using a different technique. Next we state properties that link the intensity measure in various domains; e.g., under the restriction to a subdomain (Gibbs-Markov) and conformal maps between domains. These give rise to a set of conditions that ultimately identify the intensity measure uniquely and, in fact, link them to a version of the critical Liouville Quantum Gravity.

\section{Connection to the DGFF maximum}
\noindent
On the way to the proof of Theorem~\ref{thm-extremal-vals} we have so far shown that any \emph{subsequential} limit of the measures~$\{\eta_{N,r_N}^D\colon N\ge1\}$, restricted to the first two coordinates, is a Poisson Point Process on $D\times \R$ with intensity
\begin{equation}
Z^D(\textd x)\otimes\texte^{-\alpha h}\textd h
\end{equation}
 for some random Borel measure~$Z^D$ on~$D$. Our next task is to prove the existence of the limit. Drawing on the proofs for the intermediate level sets, a natural strategy would be to identify~$Z^D$ through its properties uniquely. Although this strategy now seems increasingly viable, it was not available at the time when these results were first derived (see the remarks at the end of this lecture).  We thus base our proof of uniqueness on the connection with the DGFF maximum, along the lines of the original proof in~\cite{BL1}.

We start by formally recording an observation used in the last proof of the previous lecture:

\begin{mylemma}
\label{lemma-10.1ua}
Suppose $N_k\to\infty$ is such that
\begin{equation}
\eta^D_{N_k,r_{N_k}}\,\,\underset{k\to\infty}\lawarrow\,\,\,\PPP\bigl(Z^D(\textd x)\otimes\texte^{-\alpha h}\textd h\bigr).
\end{equation}
Then for each~$t\in\R$,
\begin{equation}
P\Bigl(\,\max_{x\in D_{N_k}}h^{D_{N_k}}_x<m_{N_k}+t\Bigr)\,\,\underset{k\to\infty}\longrightarrow\,\,E\Bigl(\texte^{-Z^D(D)\alpha^{-1}\texte^{-\alpha t}}\Bigr)
\end{equation}
\end{mylemma}

\begin{proofsect}{Proof}
Apply \eqref{E:9.53} along with the known form of the intensity measure.
\end{proofsect}

Hence we get:

\begin{mycorollary}
\label{cor-10.2}
If $\max_{x\in D_N}h^{D_N}_x-m_N$ converges in distribution, then the law of~$Z^D(D)$ is the same for every subsequential limit of $\{\eta^D_{N,r_N}\colon N\ge1\}$.
\end{mycorollary}

\begin{proofsect}{Proof}
The function~$t\mapsto\alpha^{-1}\texte^{-\alpha t}$ sweeps through $(0,\infty)$ as~$t$ varies through~$\R$. The limit distribution of the maximum thus determines the Laplace transform of the random variable~$Z^D(D)$ which in turn determines the law of~$Z^D(D)$.
\end{proofsect}

\noindent
In the original argument in~\cite{BL1}, the premise of Corollary~\ref{cor-10.2} was supplied by the main result of Bramson, Ding and Zeitouni~\cite{BDingZ} (which conveniently appeared while the first version of~\cite{BL1} was being written):

\begin{mytheorem}[Convergence of DGFF maximum]
\label{thm-BDZ}
As $N\to\infty$, the centered maximum $\max_{x\in V_N}h^{V_N}_x-m_N$ of the DGFF in~$V_N:=(0,N)^2\cap\Z^2$ converges  in law to a non-degenerate random variable.
\end{mytheorem}

The proof of Theorem~\ref{thm-BDZ} in \cite{BDingZ} is (similarly as the proof of tightness in~\cite{BZ}) based on an intermediate process between the BRW and the DGFF called the modified Branching Random Walk. In Lecture~\ref{lec-12} we will give a different proof that instead relies on the concentric decomposition and entropic-repulsion techniques encountered already in our proof of tightness. Our proof of Theorem~\ref{thm-BDZ} runs somewhat logically opposite to Corollary~\ref{cor-10.2} as it yields directly the uniqueness of the~$Z^D$ measure; the limit of the DGFF maximum then follows from Lemma~\ref{lemma-10.1ua}. The following text highlights some of the main ideas shared by both proofs.

\begin{figure}[t]
\vglue2mm
\centerline{\includegraphics[width=0.5\textwidth]{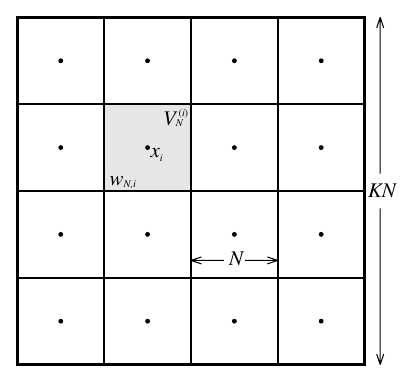}
}
\vglue2mm
\begin{quote}
\small 
\caption{
\label{fig-N-K}
\small
The decomposition of~$V_{NK}$ into~$K^2$ translates of~$V_N$. The center of (the shaded box)~$V_N^{\ssup i}$ is denoted by~$x_i$; the lower-left corner thereof by~$w_{N,i}$.}
\normalsize
\end{quote}
\end{figure}

Pick~$N,K\in\N$ large and consider the DGFF on~$V_{KN}$. Building towards the use of the Gibbs-Markov decomposition, identify within~$V_{KN}$ the (unique) collection of~translates $\{V_N^{\ssup{i}}\colon i=1,\dots,K^2\}$ of~$V_N$ separated by ``lines of sites'' in-between; cf Fig.~\ref{fig-N-K}. Realizing the DGFF $h^{V_{KN}}$ on~$V_{KN}$ by way of a binding field and~$K^2$ independent copies of the DGFF on~$V_N$'s, i.e.,  
\begin{equation}
h^{V_{KN}}:=h^{V_{KN}^\circ}+\varphi^{V_{KN},V_{KN}^\circ}\quad\text{with}\quad h^{V_{KN}^\circ}\independent\varphi^{V_{KN},V_{KN}^\circ}\,,
\end{equation}
where $V_{KN}^\circ:=\bigcup_{i=1}^{K^2}V_N^{\ssup{i}}$, we will now study what the occurrence of the absolute maximum of~$h^{V_{KN}}$ at some~$x\in V_N^{\ssup{i}}$ means for the DGFF on~$V_N^{\ssup{i}}$.

We first observe that we may safely ignore the situations when~$x$ falls near the said ``lines of sites.'' Indeed, in Lecture~\ref{lec-12} we will prove:

\begin{mylemma}
\label{lemma-10.4}
There is a constant~$c>0$ such that for all~$N\ge1$, all $t$ with $0\le |t|\le (\log N)^{1/5}$ and all non-empty sets $A\subseteq D_N$,
\begin{equation}
P\bigl(\,\max_{x\in A}h^{D_N}_x\ge m_N+t\bigr)\le c(1+t^2)\texte^{-\alpha t}\frac{|A|}{N^2}
\biggl[\log\Bigl(1+\frac{N^2}{|A|}\Bigr)\biggr]^5.
\end{equation}
\end{mylemma}

\begin{myremark}
We note that, although we expect the fifth power to be suboptimal, the example of~$A$ being a singleton shows that a logarithmic correction cannot be eliminated completely. The prefactor~$t^2$ is suboptimal as well as it can be omitted for~$t<0$ and reduced to~$t$ for~$t>0$; see Theorem~\ref{prop-10.7} below.
\end{myremark}

Next we note that, once~$x$ is at least~$\epsilon N$ away from the boundaries of boxes~$V_N^{(i)}$, the binding field~$\varphi^{V_{KN},V_{KN}^\circ}$ is well behaved. This might suggest that we could replace the whole binding field by its value at the center~$x_i$ of~$V_N^{\ssup{i}}$. However, this is correct only to the leading order as the field
\begin{equation}
\label{E:10.6}
x\mapsto \varphi^{V_{KN},V_{KN}^\circ}_{x}-\varphi^{V_{KN},V_{KN}^\circ}_{x_i}\quad\text{for~$x$ with }
\dist(x,(V_N^{\ssup{i}})^\cc)>\epsilon N,
\end{equation}
 retains variance of order unity (after all, it scales to a non-trivial Gaussian field with smooth, but non-constant, sample paths). Notwithstanding, the covariances at the center points satisfy
\begin{multline}
\qquad
\Cov\Bigl(\varphi^{V_{KN},V_{KN}^\circ}_{x_i},\varphi^{V_{KN},V_{KN}^\circ}_{x_j}\Bigr)
\\=
\Cov\bigl(h^{V_{KN}}_{x_i},h^{V_{KN}}_{x_j}\bigr)-\Cov\bigl(h^{V_{N}^{\ssup{i}}}_{x_i},h^{V_N^{(j)}}_{x_j}\bigr)
\qquad
\end{multline}
from which, checking the cases $i=j$ and~$i\ne j$ separately, we get
\begin{equation}
\Cov\Bigl(\varphi^{V_{KN},V_{KN}^\circ}_{x_i},\varphi^{V_{KN},V_{KN}^\circ}_{x_j}\Bigr)
=g\log\Bigl(\frac{KN}{|x_i-x_j|\vee N}\Bigr)+O(1).
\end{equation}
In particular, $\{\varphi^{V_{KN},V_{KN}^\circ}_{x_i}\colon i=1,\dots,K^2\}$ behaves very much like the DGFF in~$V_K$. This is further corroborated by solving:
 
\begin{myexercise}
Prove that for $c>0$ small enough,
\begin{equation}
\label{ex:10.6}
\sup_{N\ge1}\,P\Bigl(\,\max_{i=1,\dots,K^2}\varphi^{V_{KN},V_{KN}^\circ}_{x_i}>2\sqrt g\log K-c\log\log K\Bigr)
\,\,\underset{K\to\infty}\longrightarrow\,\,0\,.
\end{equation}
\end{myexercise}

The factor $c\log\log K$ is quite important because, ignoring the variations of the field in \eqref{E:10.6}, the assumption $h^{V_{KN}}_x\ge m_{KN}+t$ for  $x\in V_N^{\ssup{i}}$ at least~$\epsilon N$ away from the boundary of the box implies (on the complementary event to \eqref{ex:10.6})
\begin{equation}
\label{E:10.9ueu}
\begin{aligned}
h^{V_N^{\ssup{i}}}_x
&\ge m_{KN}-2\sqrt g\log K+c\log\log K+t
\\
&\ge m_N+c\log\log K+t+O(1),
\end{aligned}
\end{equation}
i.e., (for~$K$ large) $h^{V_N^{\ssup{i}}}$ takes an \emph{unusually high} value at~$x$. 
Since the binding field has a well-defined (and smooth) scaling limit, it thus appears that, to prove convergence in law for the maximum the DGFF in a large box~$V_{KN}$ as~$N\to\infty$, it suffices to prove the convergence in law of the maximum in~$V_N$ conditioned to exceed~$m_N+t$, in the limits as~$N\to\infty$ followed by~$t\to\infty$. 

What we just described will indeed be our strategy except that, since the spatial variations of the field \eqref{E:10.6} remain non-trivial in the said limits, we have to control \emph{both} the value of the maximum \emph{and} the position of the maximizer of the DGFF in~$V_N$. This is the content of the following theorem proved originally (albeit in a  different form) as \cite[Proposition~4.1]{BDingZ}:

\begin{mytheorem}[Limit law for large maximum]
\label{prop-10.7}
There is~$c_\star>0$ such that
\begin{equation}
\label{E:10.11}
P\bigl(\,\max_{x\in V_N}h^{V_N}_x\ge m_N+t\bigr)=\bigl[c_\star+o(1)\bigr] t\texte^{-\alpha t}\,,
\end{equation}
where $o(1)\to0$ in the limit~$N\to\infty$ followed by~$t\to\infty$. Furthermore, there is a continuous function $\psi\colon(0,1)^2\to[0,\infty)$ such that for all $A\subseteq(0,1)^2$ open,
\begin{multline}
\label{E:10.11ueu}
\quad
P\Bigl(\,N^{-1}\operatornamewithlimits{argmax}_{V_N}h^{V_N}\in A\,\Big|\,\max_{x\in V_N}h^{V_N}_x\ge m_N+t\Bigr)
\\
=o(1)+\int_A\psi(x)\textd x\,,
\quad
\end{multline}
where $o(1)\to0$ in the limit~$N\to\infty$ followed by~$t\to\infty$.
\end{mytheorem}

Our proof of Theorem~\ref{prop-10.7} will come in Lecture~\ref{lec-12}, along with the proof of Theorem~\ref{thm-BDZ} and the finishing touches for the proof of Theorem~\ref{thm-extremal-vals}. We will now comment on how Theorem~\ref{prop-10.7} is used in the proof of Theorem~\ref{thm-BDZ} in~\cite{BDingZ}. A key point to observe is that (except for the $o(1)$ term) the right-hand side of \eqref{E:10.11ueu} is independent of~$t$; the position of the conditional maximum thus becomes \emph{asymptotically independent} of its value. This suggest that we define an auxiliary process of triplets of {independent} random variables,
\begin{equation}
\bigl\{(\wp_i,h_i,X_i)\colon i=1,\dots,K^2\bigr\},
\end{equation}
which encode the limiting centered values and scaled positions of the excessive maxima in the boxes~$V_N^{\ssup{i}}$, as follows: Fix a ``cutoff scale''
 \begin{equation}
t_K:=(c/2)\log\log K
\end{equation}
 with~$c$ as in \eqref{ex:10.6}. Noting that $\psi$ is a probability density on $(0,1)^2$, set
\settowidth{\leftmargini}{(11111)}
\begin{enumerate}
\item[(1)] $\wp_i\in\{0,1\}$ with $P(\wp_i=1):=c_\star\, t_K\,\texte^{-\alpha t_K}$,
\item[(2)] $h_i\ge0$ with $P(h_i\ge t):=\frac{t_K+t}{t_K}\,\texte^{-\alpha t}$ for all $t\ge0$, and
\item[(3)] $X_i\in[0,1]^2$ with $P\bigl(X_i\in A):=\int_A\psi(x)\textd x$.
\end{enumerate}
The meaning of~$\wp_i$ is that $\wp_i=1$ designates that the maximum of the DGFF in~$V_N^{\ssup i}$ is at least~$m_N+t_K$. The excess of this maximum above~$m_N+t_K$ is then distributed  as~$h_i$ while the position of the maximizer (scaled and shifted to a unit box) is distributed as~$X_i$. See \twoeqref{E:10.11}{E:10.11ueu}.

Writing~$w_{K,i}$ for the lower-left corner of~$V_N^{\ssup{i}}$, from Theorem~\ref{thm-DZ2} (which ensures that only one value per~$V_N^{\ssup i}$ needs to be considered) it follows that the absolute maximum of~$h^{V_{KN}}-m_{NK}$ is for~$N$ large well approximated by
\begin{multline}
\label{E:10.12}
\quad
2\sqrt g\,\log K+t_K
\\
+\max\Bigl\{h_i+\varphi^{V_{KN},V_{KN}^\circ}_{w_{N,i}+\lfloor NX_i\rfloor}\colon i=1,\dots,K^2,\,\wp_i=1\Bigr\} \,, \quad
\end{multline}
where the very first term arises as the $N\to\infty$ limit of $m_{NK}-m_N$ and where the binding field is independent of the auxiliary process. The choice of~$t_K$ guarantees that, for~$N$ large, $\{i=1,\dots,K^2\colon \wp_i=1\}\ne\emptyset$ with high probability. 

Next we note that all $N$-dependence in \eqref{E:10.12} now comes through the binding field which, for~$K$ fixed and~$N\to\infty$, converges to
\begin{equation}
\Phi^{S,S^\circ_K}:=\NN(0,C^{S,S^\circ_K})\,,
\end{equation}
 where
\begin{equation}
S:=(0,1)^2\quad\text{and}\quad S^\circ:=\bigcup_{i=1}^{K^2}[\hat w_i+(0,1/K)^2]
\end{equation}
with~$\hat w_i$ denoting the position of the lower-left corner of the box that, for finite~$N$, has the lower-left corner at~$w_{N,i}$. It follows that \eqref{E:10.12} approximates \emph{every} subsequential limit of the absolute maximum by 
\begin{multline}
\label{E:10.18nwt}
\quad
2\sqrt g\,\log K+t_K
\\+\max\Bigl\{h_i+\Phi^{S,S^\circ_K}(\hat w_i+X_i/K)\colon i=1,\dots,K^2,\,\wp_i=1\Bigr\}\,,
\quad
\end{multline}
with all errors washed out when~$K\to\infty$. The maximum thus converges in law to the $K\to\infty$ limit of the random variable in \eqref{E:10.18nwt} (whose existence is a by-product of the above reasoning).

\section{Gumbel convergence}
\label{sec:10.2}\noindent\nopagebreak
Moving back to the convergence of the two-coordinate extremal process, let us explain on how the above is used to prove the uniqueness of the subsequential limit of the processes~$\{\eta^D_{N,r_N}\colon N\ge1\}$. First off, all that was stated above for square boxes~$V_N$ extends readily to any admissible sequence~$D_N$ of lattice approximations of~$D\in\mathfrak D$. Defining, for $A\subseteq D$ with a non-empty interior, 
\begin{equation}
h^{D_N}_{A,\star}:=\max_{\begin{subarray}{c}
x\in D_N\\x/N\in A
\end{subarray}}
h^{D_N}_x.
\end{equation}
the methods underlying the proof of Theorem~\ref{thm-BDZ} also give:

\begin{mylemma}
\label{lemma-mult-max}
Let~$A_1,\dots,A_k\subseteq D$ be disjoint open sets.  Then the joint law of
\begin{equation}
\bigl\{h^{D_N}_{A_i,\star}-m_N\colon i=1,\dots,k\bigr\}
\end{equation}
admits a non-degenerate weak limit as~$N\to\infty$.
\end{mylemma}

\noindent
This lemma was originally proved as~\cite[Theorem~5.1]{BL1}; we will prove it in Lecture~\ref{lec-12} via the uniqueness of the limiting~$Z^D$ measure. This will be aided by:

\begin{mylemma}
\label{lemma-10.8}
For any subsequential limit~$\eta^D=\PPP(Z^D(\textd x)\otimes\texte^{-\alpha h}\textd h)$ of $\{\eta^D_{N,r_N}\colon N\ge1\}$, any disjoint open sets~$A_1,\dots,A_k\subseteq D$ and any~$t_1,\dots,t_k\in\R$,
\begin{multline}
\quad
P\bigl(h^{D_N}_{A_i,\star}<m_N+t_i\colon i=1,\dots,k\bigr)
\\\underset{N\to\infty}\longrightarrow\,\,E\bigl(\texte^{-\sum_{i=1}^k Z^D(A_i)\alpha^{-1}\texte^{-\alpha t_i}}\bigr).
\quad
\end{multline}
\end{mylemma}

\begin{proofsect}{Proof (sketch)}
The event on the left-hand side can be written as $\langle \eta^D_{N,r_N},f\rangle=0$, where $f:=\sum_{i=1}^k1_{A_i}\otimes 1_{[t_i,\infty)}$. Since the right-hand side equals the probability that $\langle \eta^D,f\rangle=0$ the claim follows by approximating~$f$ by bounded continuous functions with compact support in $\overline D\times(\R\cup\{\infty\})$. The next exercise helps overcome (unavoidable) boundary issues.
\end{proofsect}

\begin{myexercise}
\label{ex:10.5}
Use Lemma~\ref{lemma-10.4} to prove that, for every subsequential limit~$\eta^D$ of the processes of interest, the associated $Z^D$ measure obeys
\begin{equation}
\forall A\subseteq \overline{D} \text{\ \rm\ Borel}\colon
\quad
\Leb(A)=0\quad\Rightarrow\quad Z^D(A)=0 \quad\text{\rm a.s.}
\end{equation}
In addition, we get $Z^D(\partial D)=0$ a.s.\ so~$Z^D$ is concentrated on~$D$. (You should not need to worry about whether~$\Leb(\partial D)=0$ for this.)
\end{myexercise}

These observations permit us to give:

\begin{proofsect}{Proof of Theorem~\ref{thm-extremal-vals} (first two coordinates)}
Lemmas~\ref{lemma-mult-max} and~\ref{lemma-10.8} imply that the joint law of $(Z^D(A_1),\dots,Z^D(A_k))$ is the same for every subsequential limit~$\eta^D$ of our processes of interest. This means that we know the law of $\langle Z^D,f\rangle$ for any~$f$ of the form $f=\sum_{i=1}^ka_i1_{A_i}$ with~$A_i$ open disjoint and $a_1,\dots,a_k>0$. Every bounded and continuous~$f$ can be approximated by a function of the form $\sum_{i=1}^k a_i1_{\{a_{i-1}\le f<a_i\}}$ with $\Leb(f=a_i)=0$ for every~$i=1,\dots,k$. With the help of Exercise~\ref{ex:10.5} we then get uniqueness of the law of $\langle Z^D,f\rangle$ for any~$f\in C_\cc(D\times\R)$. Thanks to Exercise~\ref{ex:10.5}, this identifies the law of~$Z^D$ on~$\overline D$ uniquely. That $Z^D(D)\in(0,\infty)$ a.s.\ was already shown in the proof of Theorem~\ref{prop-subseq}. 
\end{proofsect}

The structure of the limit process gives rise to interesting formulas that are worth highlighting at this point. First we ask the reader to solve:

\begin{myexercise}[Connection to Gumbel law]
Let $\{(x_i,h_i)\colon i\in\Z\}$ label the points in a sample of~$\PPP(Z^D(\textd x)\otimes\texte^{-\alpha h}\textd h)$. Show that
\begin{equation}
\label{E:10.23nwt}
\sum_{i\in\N}\delta_{h_i-\alpha^{-1}\log Z^D(D)}\,\,\laweq\,\,\PPP(\texte^{-\alpha h}\textd h).
\end{equation}
In particular, the absolute maximum $\max_{i\in\N}h_i$ has the law of a randomly shifted Gumbel random variable. See also \eqref{E:9.11nwwt}.
\end{myexercise}

By \eqref{E:10.23nwt}, the \emph{gap} between the $i$-th and~$i+1$-st local maximum are distributed as those in~$\PPP(\texte^{-\alpha h}\textd h)$. (This will largely be responsible for the simple form of the limit in Corollary~\ref{cor-11.16}.) We remark that proving convergence of the maximum to a randomly-shifted Gumbel random variable has been one of the holy grails of the subject area (stimulated by the corresponding result for the Branching Brownian Motion; cf McKean~\cite{McKean} and Lalley and Sellke~\cite{LS}).
For the so called star-scale invariant Gaussian fields (in any spatial dimension), this was proved by Madaule~\cite{Madaule} very soon after~\cite{BDingZ} and \cite{BL1} were posted. A corresponding result for general Branching Random Walks is due to A\"idekon~\cite{Aidekon}.

Another, perhaps more interesting, aspect to showcase is:

\begin{mylemma}
\label{lemma-10.11}
Let~$X_N$ be the (a.s.-unique) vertex where $h_{X_N}^{D_N}=\max_{x\in D_N}h^{D_N}_x$. Then for any $A\subseteq D$ open with~$\Leb(\partial A)=0$ and any~$t\in\R$,
\begin{multline}
\label{E:10.17}
\quad
P\Bigl(\tfrac1N X_N\in A,\,\max_{x\in D_N}h^{D_N}_x<m_N+t\Bigr)
\\\underset{N\to\infty}\longrightarrow\, E\bigl(\wh Z^D(A)\texte^{-\alpha^{-1}\texte^{-\alpha t}Z^D(D)}\bigr),
\quad
\end{multline}
where $\wh Z$ is the measure from \eqref{E:10.18ueu}.
\end{mylemma}

\begin{proofsect}{Proof}
Let~$A\subseteq D$ be open with~$\Leb(\partial A)=0$ (which implies $Z^D(\partial A)=0$ a.s.). Lemma~\ref{lemma-mult-max} along with a continuity argument based on Exercise~\ref{ex:10.5} show
\begin{equation}
\bigl(h^{D_N}_{A,\star}-m_N,\,h^{D_N}_{A^\cc,\star}-m_N\bigr)\,\,\underset{N\to\infty}\lawarrow\,\,(h^\star_A,h^\star_{A^\cc})\,.
\end{equation}
The continuity of the law of the DGFF yields
\begin{equation}
\begin{aligned}
P\Bigl(\tfrac1N X_N\in A,\,\max_{x\in D_N}&h^{D_N}_x<m_N+t\Bigr) 
\\
&= P\bigl(h^{D_N}_{A^\cc,\star}<h^{D_N}_{A,\star},\,h^{D_N}_{A,\star}-m_N<t\bigr)
\\
&= P\bigl(h^{D_N}_{A^\cc,\star}\le h^{D_N}_{A,\star},\,h^{D_N}_{A,\star}-m_N\le t\bigr).
\end{aligned}
\end{equation}
Since the event in the first line is open while that in the second line is closed, and the limit distribution are continuous thanks to Theorem~\ref{prop-subseq} and Exercise~\ref{ex:10.12nwt} below, the standard facts about the convergence in law imply
\begin{multline}
\label{E:10.21}
\quad
P\Bigl(\tfrac1N X_N\in A,\,\max_{x\in D_N}h^{D_N}_x<m_N+t\Bigr) 
\\\underset{N\to\infty}\longrightarrow\, P\bigl(h^\star_{A^\cc}<h^\star_A,\,h_A^\star<t\bigr).
\quad
\end{multline}
Now we invoke:

\begin{myexercise}
\label{ex:10.12nwt}
Let~$A\subseteq D$ be open with~$\Leb(\partial A)=0$. Prove that
\begin{equation}
h_A^\star\,\,\laweq\,\,\inf\bigl\{t\in\R\colon\eta^D(A\times[t,\infty))=0\}
\end{equation}
and that this in fact applies jointly to $(h^\star_A,h^\star_{A^\cc})$.
\end{myexercise}

\noindent
This means that we can now rewrite the probability on the right of \eqref{E:10.21} in terms of the maximal points in the sample of~$\eta^D$. For $\eta = \PPP(M\otimes\texte^{-\alpha h}\textd h)$ with a fixed~$M$, the point with the maximal $h$-value in~$A$ has probability density $\alpha^{-1}\texte^{-\alpha h}M(A)\texte^{-\alpha^{-1}\texte^{-\alpha h}M(A)}$ with respect to the Lebesgue measure, while the probability that no point in~$A^\cc$ has $h$-value above~$h$ is~$\texte^{-\alpha^{-1}\texte^{-\alpha h}M(A^\cc)}$. The underlying Poisson structure (conditional on~$Z^D$) therefore gives
\begin{equation}
P\bigl(h^\star_{A^\cc}<h^\star_A,\,h_A^\star<t\bigr) = \int_{-\infty}^t \alpha^{-1}\texte^{-\alpha h}\,E\Bigl(Z^D(A)\texte^{-\alpha^{-1}\texte^{-\alpha h}Z^D(D)}\Bigr)\textd h.
\end{equation}
The result now follows by integration (and Tonelli's Theorem).
\end{proofsect}

Hereby we get:

\begin{mycorollary}
\label{cor-10.13}
The measure $A\mapsto E(\wh Z^D(A))$ is the $N\to\infty$ weak limit of the marginal law of the (a.s.-unique) maximizer of~$h^{D_N}$ scaled by~$N$.
\end{mycorollary}

\section{Properties of $Z^D$-measures}
\noindent
Once the convergence issue has been settled (modulo proofs deferred till later), the next natural question is: What is~$Z^D$? Or, more precisely: Can the law of~$Z^D$ be independently characterized?
As was the case of intermediate level sets, although the laws of $\{Z^D\colon D\in\mathfrak D\}$ are defined individually, they are very much interrelated. The following theorem articulates these relations explicitly: 

\begin{mytheorem}[Properties of~$Z^D$-measures]
\label{thm-10.14}
The family~$\{Z^D\colon D\in\mathfrak D\}$ satisfies the following:
\settowidth{\leftmargini}{(1111)}
\begin{enumerate}
\item[(1)] $Z^D(A)=0$ a.s.\ when $A=\partial D$ and for any $A\subseteq D$ with $\Leb(A)=0$,
\item[(2)] for any $a\in\C$ and any $b>0$, 
\begin{equation}
Z^{a+b D}(a+b\textd x)\,\,\laweq\,\, b^4 Z^D(\textd x),
\end{equation}
\item[(3)] if $D\cap\widetilde D=\emptyset$, then
\begin{equation}
Z^{D\cup\wt D}(\textd x)\,\,\laweq\,\, Z^D(\textd x)+Z^{\wt D}(\textd x),
\quad Z^D\independent Z^{\wt D},
\end{equation}
\item[(4)]
if $\widetilde D\subseteq D$ and $\Leb(D\smallsetminus\widetilde D)=0$, then for $\Phi^{ D, \wt{ D}}:=\NN(0,C^{D,\wt D})$ and $\alpha:=2/\sqrt g$,
\begin{equation}
Z^D(\textd x) \,\,\laweq\,\, Z^{\wt{ D}}(\textd x)\, \texte^{\alpha  \Phi^{ D, \wt{ D}}(x)},\quad
Z^{\wt{ D}}\independent \Phi^{ D, \wt{ D}},
\end{equation}
\item[(5)] there is~$\hat c\in(0,\infty)$ such that for all open $A\subseteq D$,
\begin{equation}
\label{E:1.25ue}
\lim_{\lambda\downarrow0}\frac{E(Z^D(A)\texte^{-\lambda Z^D(D)})}{\log(\ffrac1\lambda)}= \hat c\int_A r_D(x)^2\,\textd x,
\end{equation}
where $r_D(x)$ is the conformal radius of~$D$ at~$x$.
\end{enumerate}
\end{mytheorem}

Before we get to the proof of this theorem, let us make a few remarks. First off, with the exception of~(5), these properties are shared by the whole family of measures~$\{Z^D_\lambda\colon D\in\mathfrak D\}$ introduced in Lectures~\ref{lec-2}--\ref{lec-5}. The condition~(5) is different than for the measures $Z^D_\lambda$ because (by~(5)) $E Z^D(A)=\infty$ for any non-empty open~$A\subseteq D$. Incidentally, this is also what stands in the way of proving uniqueness of the law of~$Z^D$ by the argument underlying Proposition~\ref{prop-4.12}.

\begin{proofsect}{Proof of Theorem~\ref{thm-10.14}, (1) and (3)}
Part (1) is proved by way of reference to Exercise~\ref{ex:10.5}.
For (3) we just observe that the DGFF on domains separated by at least two lattice steps are independent. (This is one place where the first condition \eqref{E:1.22} of admissible approximations is needed.)
\end{proofsect}

Next we will establish the \emph{Gibbs-Markov property}:

\begin{proofsect}{Proof of Theorem~\ref{thm-10.14}(4)}
Most of the argument can be borrowed from the proof of Proposition~\ref{prop-GM-Z} on the Gibbs-Markov property for the measures arising from the intermediate level sets. Let~$\wt D$ and~$D$ as in the statement, pick~$f\in C_\cc(\wt D\times\R)$ with~$f\ge0$ and recall the notation
\begin{equation}
f_\Phi(x,h):=f\bigl(x,\,h+\Phi^{D,\wt D}(x)\bigr)\,.
\end{equation}
Writing~$\eta^D$ and~$\eta^{\wt D}$ for the limit processes in the respective domains, the argument leading up to \eqref{E:4.21} then ensures
\begin{equation}
\langle \eta^D,f\rangle \,\,\laweq\,\, \langle\eta^{\wt D},f_\Phi\rangle, \quad\Phi^{D,\wt D}\independent\eta^{\wt D}\,.
\end{equation}
The Poisson nature of the limit process then gives, via a routine change of variables,
\begin{equation}
\begin{aligned}
E\bigl(\texte^{- \langle\eta^{\wt D},f_\Phi\rangle}\bigr)
&=E\Bigl(\exp\Bigl\{-\int Z^{\wt D}(\textd x)\texte^{-\alpha h}\textd h\,\bigl(1-\texte^{-f(x,\,h+\Phi^{D,\wt D}(x))}\bigr)\Bigr\}\Bigr)
\\
&=E\Bigl(\exp\Bigl\{-\int Z^{\wt D}(\textd x)\texte^{\alpha\Phi^{D,\wt D}(x)}\texte^{-\alpha h}\textd h\,\bigl(1-\texte^{-f(x,h)}\bigr)\Bigr\}\Bigr)\,.
\end{aligned}
\end{equation}
The claim follows by comparing this with the expression one would get for the Laplace transform of $\langle \eta^D,f\rangle$.
\end{proofsect}

With the help of the Gibbs-Markov property we then readily solve:

\begin{myexercise}[Support of~$Z^D$ and its non-atomicity]
Prove that~$Z^D$ charges every non-empty open subset of~$D$ with probability one. In particular, we have $\supp Z^D=\overline D$ a.s. Prove also that $Z^D$ is a.s.\ non-atomic. [Hint: Recall that Theorem~\ref{thm-DZ2} rules out nearby occurrence of near-maximal local maxima.]
\end{myexercise} 

Next let us address the behavior under shifts and dilations:

\begin{proofsect}{Proof of Theorem~\ref{thm-10.14}(2) (sketch)}
The proof will require working with approximations of the given domain. We thus first ask the reader to solve:

\begin{myexercise}
\label{ex:10.16}
Let~$D\in\mathfrak D$ and assume that $\{D^n\colon n\ge1\}\in\mathfrak D$ are such that~$D^n\uparrow D$ 
with $C^{D,D^n}(x,y)\to0$ locally uniformly on~$D$. Then
\begin{equation}
Z^{D^n}(\textd x)\,\,\underset{n\to\infty}\lawarrow\,\,Z^D(\textd x).
\end{equation}
Hint: Use the Gibbs-Markov property.
\end{myexercise}

Thanks to this exercise we may assume that~$a\in\Q^2$ and~$b\in\Q$. As $aN$ and~$bN$ will then be integer for an infinite number of~$N$'s, the existence of the limit permits us to even assume that~$a\in\Z^2$ and~$b\in\N$. The invariance of the law of~$Z^D$ under integer-valued shifts is a trivial consequence of the corresponding invariance of the DGFF. Concerning the behavior under scaling, here we note that if~$\{D_N\colon N\ge1\}$ approximates~$D$, then $\{D_{bN}\colon N\ge1\}$ approximates~$bD$. The only item to worry about is the centering sequence for which we get
\begin{equation}
m_{bN}-m_N = 2\sqrt g\,\log b+o(1).
\end{equation}
Following this change through the limit procedure yields
\begin{equation}
Z^{b D}(b\,\textd x)\,\,\laweq\,\, \texte^{\alpha 2\sqrt g\,\log b} Z^D(\textd x).
\end{equation}
The claim follows by noting that $\alpha2\sqrt g = 4$.
\end{proofsect}

\begin{proofsect}{Proof of Theorem~\ref{thm-10.14}(5)}
The starting point is an extension of Theorem~\ref{prop-10.7} to all $D\in\mathfrak D$. Considering two domains, $D$ and~$\wt D$, with~$\wt D\subseteq D$ and~$\wt D$ a square; i.e., $\wt D=a+(0,b)^2$ for some $a\in\C$ and $b>0$. We will for a while assume that~$\wt D$ is in fact a unit square (i.e., $b=1$); by Exercise~\ref{ex:10.16} this can be achieved by redefining~$N$. Pick approximating domains $\{D_N\}$ and~$\{\wt D_N\}$ respectively, and couple the DGFFs therein via the Gibbs-Markov property to get 
\begin{equation}
h^{D_N}=h^{\wt D_N}+\varphi^{D_N,\wt D_N}\quad\text{where}\quad h^{\wt D_N}\independent\varphi^{D_N,\wt D_N}\,.
\end{equation}
We then claim the following intuitive fact:
 
\begin{mylemma}
\label{lemma-10.17}
Conditional on the maximum of $h^{D_N}$ to exceed $m_N+t$, the position of the (a.s.-unique) maximizer will, with  probability tending to one as~$N\to\infty$ and~$t\to\infty$, lie within $o(N)$ distance of the position of the maximizer of $h^{\wt D_N}$.
\end{mylemma}

\noindent
We will will not prove this lemma here; instead,  we refer the reader to~\cite[Proposition~5.2]{BL2}. (The proof can also be gleaned from that for a closely-related, albeit far stronger, Lemma~\ref{lemma-same-max} that we give in Lecture~\ref{lec-12}.) We now extract the limit law of the scaled maximizer/centered maximum for both domains and write $(X_\star,h_\star)$ for the limiting pair for~$D$ and $(\wt X_\star,\wt h_\star)$ for the limiting pair for~$\wt D$ (these exist by Lemma~\ref{lemma-10.11}). Let~$A\subset \wt D$ be a non-empty open set with~$\overline A\subset\wt D$ and $\Leb(\partial A)=0$. Lemma~\ref{lemma-10.17} and Exercise~\ref{ex:10.5} then give
\begin{equation}
P\bigl(X_\star\in A,\,h_\star>t\bigr)=\bigl(1+o(1)\bigr)P\Bigl(\wt X_\star\in A,\,h_\star+\Phi^{D,\wt D}(\wt X_\star)>t\Bigr),
\end{equation}
where~$o(1)\to0$ in the limit~$t\to\infty$. Since~$\wt D$ is a unit square and~$\Phi^{D,\wt D}$ has uniformly continuous sample paths on~$A$, a routine approximation argument combined with conditioning on~$\Phi^{D,\wt D}$ and Theorem~\ref{prop-10.7} show
\begin{multline}
\quad
t^{-1}\texte^{\alpha t}P\Bigl(\wt X_\star\in A,\,\wt h_\star+\Phi^{D,\wt D}(\wt X_\star)>t\Bigr)
\\\underset{t\to\infty}\longrightarrow\,c_\star\,E\Bigl(\,\int_A\texte^{\alpha\Phi^{D,\wt D}(x)}\psi(x)\textd x\Bigr)\,.
\quad
\end{multline}
Hereby we get
\begin{equation}
\label{E:10.36}
t^{-1}\texte^{\alpha t}
P\bigl(X_\star\in A,\,h_\star>t\bigr)
\,\,\underset{t\to\infty}\longrightarrow\,\int_A\psi^D(x)\textd x\,,
\end{equation}
where
\begin{equation}
\psi^D(x):=c_\star\psi(x)\texte^{\frac12\alpha^2\,C^{D,\wt D}(x,x)},\quad x\in\wt D.
\end{equation}
Note that, since~$\wt D$ is a unit square, we get $\psi^{\wt D}(x):=c_\star\psi(x)$ directly from Theorem~\ref{prop-10.7}.

From~$\frac12\alpha^2 g = 2$ and the relation between~$C^{D,\wt D}$ and the conformal radius (see Definition~\ref{def-1.24} and \eqref{E:4.1}) we then get
\begin{equation}
\label{E:10.41}
\frac{\psi^D(x)}{r_D(x)^2}=\frac{\psi^{\wt D}(x)}{r_{\wt D}(x)^2},\quad x\in\wt D.
\end{equation}
This was derived for~$\wt D$ a unit square but the following exercise relaxes that:

\begin{myexercise}
Suppose that $\wt D$ is a unit square centered at~$0$. Show that, for any~$b\ge 1$, \eqref{E:10.36} holds for a domain~$D\in\mathfrak D$ with~$\wt D\subseteq D$ if and only if it holds for domain~$bD$. Moreover, 
\begin{equation}
\psi^{bD}(bx)=b^2\psi^D(x).
\end{equation}
Conclude that \eqref{E:10.41} holds if~$\wt D$ is a square of any size (of the form $a+(0,b)^2$).
\end{myexercise}

\noindent
We now claim that, for all~$\wt D,D\in\mathfrak D$,
\begin{equation}
\label{E:10.41x}
\frac{\psi^D(x)}{r_D(x)^2}=\frac{\psi^{\wt D}(x)}{r_{\wt D}(x)^2},\quad x\in D\cap\wt D.
\end{equation} 
Indeed, if~$x\in D\cap\wt D$ then (as $D\cap\wt D$ is open) there is a square~$\wt D'\subseteq D\cap\wt D$ containing~$x$ and so we get \eqref{E:10.41x} by iterating \eqref{E:10.41}. Using \eqref{E:10.41x} for~$\wt D$ a translate of~$D$ shows that $x\mapsto\psi^D(x)/r_D(x)^2$ is constant on~$D$; one more use of \eqref{E:10.41x} shows that this constant is the same for all admissible domains. So
\begin{equation}
\psi^D(x)=c r_D(x)^2,\quad x\in D,
\end{equation}
for some~$c>0$ independent of~$D$.

Having extended Theorem~\ref{prop-10.7} to all~$D\in\mathfrak D$, writing $\lambda:=\alpha^{-1}\texte^{-\alpha t}$, \eqref{E:10.36}, Corollary~\ref{cor-10.13} and \eqref{E:10.17} then show
\begin{equation}
\label{E:10.43ueu}
\frac{E(\wh Z^D(A)[1-\texte^{-\lambda Z^D(D)}])}{\lambda\log(\ffrac1\lambda)}\,\,\underset{\lambda\downarrow0}\longrightarrow\,\,\alpha c_\star\,c\int_A r_D(x)^2\textd x
\end{equation}
for any~$A\subseteq D$ open with $\Leb(\partial A)=0$. The claim follows with~$\hat c:=\alpha c_\star\, c$ from the  next exercise.
\end{proofsect}

\begin{myexercise}
\label{ex:10.19}
Prove that \eqref{E:10.43ueu} implies \eqref{E:1.25ue}.
\end{myexercise}

We note that, in Lecture~\ref{lec-12}, we will establish the explicit relation between the asymptotic density~$\psi^D$ and the square of the conformal radius directly.

\section{Uniqueness up to overall constant}
\noindent
As our next item of business, we wish to explain that the properties listed in Theorem~\ref{thm-10.14} actually determine the laws of the~$Z^D$'s uniquely. 

\begin{figure}[t]
\vglue2mm
\centerline{\includegraphics[width=0.6\textwidth]{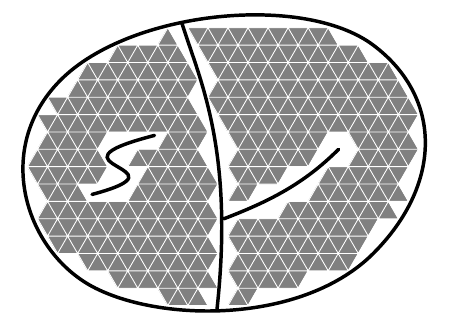}
}
\vglue2mm
\begin{quote}
\small 
\caption{
\label{fig-triangles}
\small
A tiling of domain~$D$ by equilateral triangles of side-length~$K^{-1}$.}
\normalsize
\end{quote}
\end{figure}

Given an integer~$K\ge1$, consider a tiling of the plane by equilateral triangles of side-length~$K^{-1}$. For a domain~$D\in\mathfrak D$, let $T^1,\dots,T^{m_K}$ be the triangles entirely contained in~$D$, cf Fig.~\ref{fig-triangles}. Abbreviate
\begin{equation}
\label{E:10.45aa}
\wt D:=\bigcup_{i=1}^{m_K}T^i.
\end{equation}
Given $\delta\in(0,1)$, label the triangles so that $i=1,\dots,n_K$, for some $n_K\le m_K$, enumerate the triangles that are at least distance~$\delta$ away from~$D^\cc$. Define $T^1_\delta,\dots,T^{n_K}_\delta$ as the equilateral triangles of side length $(1-\delta)K^{-1}$ that have the same orientation and centers as $T^1,\dots,T^{n_K}$, respectively. 
Recall that the oscillation of a function~$f$ on a set~$A$ is given by
\begin{equation}
\text{osc}_A f:=\sup_{x\in A}f(x)-\inf_{x\in A}f(x).
\end{equation}
We then claim:

\begin{mytheorem}
\label{thm-10.16}
Consider a family $\{M^D\colon D\in\mathfrak D\}$ of random Borel measures satisfying (1-5) in Theorem~\ref{thm-10.14} with some~$\hat c\in(0,\infty)$. Define events $A^i_{K,R}$, $i=1,\dots,n_K$ by
\begin{equation}
\label{E:10.52nwwt}
A^i_{K,R}:=\bigl\{\text{\rm osc}_{T^i_\delta}\Phi^{D,\wt D}\le R\bigr\}\cap\bigl\{\max_{{T^i_\delta}}\Phi^{D,\wt D}\le2\sqrt g\log K-R\bigr\}.
\end{equation}
Then for any~$D\in\mathfrak D$, for~$\wt D$ related to~$K$ as in \eqref{E:10.45aa} and $\Phi^{D,\wt D}:=\NN(0,C^{D,\wt D})$, the random measure
\begin{multline}
\label{E:10.45}
\quad\alpha \hat c\, r_D(x)^2\sum_{i=1}^{n_K}\1_{A_{K,R}^i}\bigl(\alpha\Var(\Phi^{D,\wt D}(x))-\Phi^{D,\wt D}(x)\bigr)
\\
\times\texte^{\alpha \Phi^{D,\wt D}(x)-\frac12\alpha^2\Var(\Phi^{D,\wt D}(x))}\,\1_{T_\delta^i}(x)\,\textd x
\quad
\end{multline}
tends in law to $M^D$ in the limit as $K\to\infty$, $R\to\infty$ and~$\delta\downarrow0$ (in this order). This holds irrespective of the orientation of the triangular grid.
\end{mytheorem}

Before we delve into the proof we note that, by computations involving the covariance kernel $C^{D,\wt D}$, we get
\begin{equation}
\min_{i=1,\dots, n_K}\,\inf_{x\in T_\delta^i}\,\Var(\Phi^{D,\wt D}(x))\ge g\log K-c
\end{equation}
for some $c=c(\delta)>0$ and all $K\ge1$. Hence, for~$R$ sufficiently large (depending only on~$\delta>0$), $\alpha g = 2\sqrt g$ implies
\begin{equation}
\label{E:10.53uia}
\alpha\Var(\Phi^{D,\wt D}(x))-\Phi^{D,\wt D}(x)>0\quad\text{on }A_{K,R}^i.
\end{equation}
In particular, \eqref{E:10.45} is a positive measure. We also note that, by combining Theorems~\ref{thm-10.14} and~\ref{thm-10.16}, we get:

\begin{mycorollary}[Characterization of~$Z^D$-measures]
\label{cor-10.21}
The laws of the random measures $\{Z^D\colon D\in\mathfrak D\}$ are determined uniquely by conditions (1-5) of Theorem~\ref{thm-10.14}.
\end{mycorollary}

In order to avoid confusion about how the various results above depend on each other, we put forward:

\begin{myremark}
Although Corollary~\ref{cor-10.21} may seem to produce the desired independent characterization of the~$Z^D$ measures by their underlying properties, we remind the reader that our derivation of these properties in Theorem~\ref{thm-10.14} assumed (or came after proving, if we take Theorem~\ref{thm-BDZ} for granted) that the limit of $\{\eta^D_{N,r_N}\colon N\ge1\}$ exists, which already entails some version of their uniqueness. We will amend this in Lecture~\ref{lec-12} in our proofs of Theorems~\ref{prop-10.7} and~\ref{thm-BDZ} which draw heavily on the proof of Theorem~\ref{thm-10.16}.
\end{myremark}

The main thrust of Theorem~\ref{thm-10.16} is that it gives a representation of~$Z^D$ as the limit of the measures in \eqref{E:10.45} that are determined solely by the binding fields~$\Phi^{D,\wt D}$ (where~$\wt D$ depends on~$K$). By Exercise~\ref{ex:4.15}, we may think of $\Phi^{D,\wt D}$ as the orthogonal projection of the CGFF onto the subspace of functions in~$\cmss H_0^1(D)$  that are harmonic in each of the triangles constituting~$\wt D$. The representation by measures in \eqref{E:10.45} is thus akin to that of the~$Z^D_\lambda$-measures by way of the measures in~\eqref{E:4.22} except that instead of using
\begin{equation}
\label{E:10.49}
\texte^{\beta \Phi^{D,\wt D}(x)-\frac12\beta^2\Var(\Phi^{D,\wt D}(x))}
\end{equation}
as the density with respect to the Lebesgue measure, we use its \emph{derivative}
\begin{multline}
\label{E:10.50}
\quad
-\frac{\textd}{\textd \beta}\,\texte^{\beta \Phi^{D,\wt D}(x)-\frac12\beta^2\Var(\Phi^{D,\wt D}(x))}\Bigl|_{\beta=\alpha}
\\
=\bigl(\alpha\Var(\Phi^{D,\wt D}(x))-\Phi^{D,\wt D}(x)\bigr)\texte^{\alpha \Phi^{D,\wt D}(x)-\frac12\alpha^2\Var(\Phi^{D,\wt D}(x))}.
\quad
\end{multline}
Note that the key fact underlying the derivations in Lecture~\ref{lec-4} is still true: the expression in \eqref{E:10.50} is a martingale with respect to the filtration induced by the fields~$\Phi^{D,\wt D}$ along nested subdivisions of~$D$ into finite unions of equilateral triangles (i.e., with~$K$ varying along powers of~$2$). However, the lack of positivity makes manipulations with this object considerably more difficult.

\begin{proofsect}{Proof of Theorem~\ref{thm-10.16} (sketch)} 
The proof is based on a number of relatively straightforward observations (and some technical calculations that we will mostly skip). Denote $\wt D_\delta:=\bigcup_{i=1}^{n_K}T^i_\delta$ and, for $f\in C_\cc(D)$, let~$f_\delta:=f1_{\wt D_\delta}$. Consider a family $\{M^D\colon D\in\mathfrak D\}$ of measures in the statement. Property~(1) then ensures
\begin{equation}
\label{E:10.51uai}
\langle M^D,f_\delta\rangle\,\,\underset{\delta\downarrow0}\lawarrow\,\,\langle M^D,f\rangle\,,
\end{equation}
and so we may henceforth focus on~$f_\delta$. Properties~(3-4) then give
\begin{equation}
\label{E:10.47}
1_{\wt D_\delta}(x)M^D(\textd x)\,\,\laweq\,\,\sum_{i=1}^{n_K}\texte^{\alpha\Phi^{D,\wt D}(x)}1_{T^i_\delta}(x)\,M^{T^i}(\textd x)\,,
\end{equation}
with $M^{T^1},\dots, M^{T^{n_K}}$ and~$\Phi^{D,\wt D}$ all independent.
Let~$x_1,\dots,x_{n_K}$ be the center points of the triangles $T^1,\dots,T^{n_K}$, respectively. A variation on Exercise~\ref{ex:10.6} then shows
\begin{equation}
\label{E:10.53ua}
\limsup_{K\to\infty}\,P\Bigl(\,\max_{i=1,\dots,n_K}\Phi^{D,\wt D}(x_i)>2\sqrt g\log K-c\log\log K\Bigr)=0
\end{equation}
for some~$c>0$. The first technical claim is the content of:

\begin{myproposition}
\label{prop-10.18}
For any~$\delta\in(0,1)$ and any $\epsilon>0$,
\begin{equation}
\label{EE:1.2}
\lim_{R\to\infty}\,\limsup_{K\to\infty}\,
P\Bigl(\,\,\sum_{i=1}^{n_K}\int_{T^i_\delta}M^{T^i}(\textd x)\texte^{\alpha\Phi^{D,\wt D}(x)}\1_{\{\text{\rm osc}_{T^i_\delta}\Phi^{D,\wt D}>R\}}>\epsilon\Bigr)=0.
\end{equation}
\end{myproposition}

\begin{proofsect}{Proof (main ideas)}
We will not give a detailed proof for which we refer the reader to \cite[Proposition~6.5]{BL2}. Notwithstanding, we will give the main ideas and thus also explain why we have resorted to triangle partitions. 

For triangle partitions the subspace of functions in $\cmss H_0^1(D)$ that are piece-wise harmonic on~$\wt D$ naturally decomposes into a direct sum of the space $\cmss H^{\triangle}$ of functions that are \emph{affine} on each~$T^i$ and its orthogonal complement~$\cmss H^\perp$ in~$\cmss H_0^1(D)$. The projection onto~$\cmss H^\perp$ can be controlled uniformly in~$K$ thanks to:

\begin{myexercise}
For~$D$ and~$\wt D$ (which depends on~$K$) as above, let~$\Phi^\perp_K$ denote the orthogonal projection of~$\Phi^{D,\wt D}$ onto~$\cmss H^\perp$. Show that, for each~$\delta>0$ small,
\begin{equation}
\label{E:10.60nw}
\sup_{K\ge1}\sup_{x\in\wt D_\delta}\Var\bigl(\Phi^\perp_K(x)\bigr)<\infty.
\end{equation}
Hint: Use the observation (made already in Sheffield's review~\cite{Sheffield-review}) that
\begin{equation}
\Phi_K^\triangle:=\Phi^{D,\wt D}-\Phi_K^\perp
\end{equation} 
is piece-wise affine and so it is determined by its values at the vertices of the triangles where it has the law of (a scaled) DGFF on the triangular lattice. See the proof of \cite[Lemma~6.10]{BL2}.
\end{myexercise}

The uniform bound \eqref{E:10.60nw} implies (via Borell-TIS inequality and Fernique majorization) a uniform Gaussian tail for the oscillation of~$\Phi^\perp_K$. This gives
\begin{equation}
\label{E:6.41nwt}
\epsilon_R:=\sup_{K\ge1}\,\,\max_{i=1,\dots,n_K}\,\,\sup_{x\in T_\delta^i}\,
E\bigl(\texte^{\alpha\Phi_K^\perp(x)}\1_{\{\text{\rm osc}_{T^i_\delta}\Phi_K^{\perp}>R\}}\bigr)\,\underset{R\to\infty}\longrightarrow\,0
\end{equation}
(see \cite[Corollary~6.11]{BL2} for details). To get the claim from this, let $M_{K,R}^\perp$ denote the giant sum  in (the event in) \eqref{EE:1.2}. The independence of~$\Phi^\perp_K$ and~$\Phi_K^\triangle$
then gives, for each~$\lambda>0$ and for~$Y$ a centered normal independent of~$\Phi^\perp_K$, $\Phi^\triangle_K$ and of~$M^{T^i}$'s, and thus of $M_{K,R}^\perp$,
\begin{equation}
E\texte^{-\lambda \texte^{\alpha Y}M_{K,R}^{\perp}}
\ge E\biggl(\exp\Bigl\{-\lambda\sum_{i=1}^{n}\int_{T^i_\delta}M^{T^i}(\textd x)\texte^{\alpha\Phi_K^\triangle(x)+Y}\epsilon_R\bigr)\Bigr\}\biggr).
\end{equation}
If $\Var(Y)$ exceeds the quantity in \eqref{E:10.60nw}, Kahane's convexity inequality (cf Proposition~\ref{prop-5.6}) permits us to replace $\alpha\Phi_K^\triangle(x)+Y$ by $\alpha\Phi^{D,\wt D}+\frac12\Var(Y)$ and wrap the result (as a further lower bound) into $E\texte^{-\lambda c\epsilon_R M^D(D)}$, for some~$c>0$. As~$\epsilon_R\to0$ when~$R\to\infty$ we conclude that $M_{K,R}^{\perp}\to0$ in probability as~$K\to\infty$ and~$R\to\infty$. This gives \eqref{EE:1.2}.
\end{proofsect}

Returning to the proof of Theorem~\ref{thm-10.16}, \twoeqref{E:10.53ua}{EE:1.2} permit us to restrict attention only to those triangles where $A_{K,R}^i$ occurs. This is good because, in light of the piece-wise harmonicity of the sample paths of~$\Phi^{D,\wt D}$, the containment in $A_{K,R}^i$ forces
\begin{equation}
x\mapsto\Phi^{D,\wt D}(x)-\Phi^{D,\wt D}(x_i)
\end{equation}
 to be bounded and uniformly Lipschitz on~$T^i_\delta$, for each~$i=1,\dots,_K$. Let $\FF_{R,\beta,\delta}^T$, for $R,\beta,\delta>0$, denote the class of continuous functions $\phi\colon T\to\R$ on an equilateral triangle~$T$ such that
\begin{equation}
\label{E:5.8c}
\phi(x)\ge\beta\quad\text{and}\quad |\phi(x)-\phi(y)|\le R|x-y|,\qquad x,y\in T_\delta.
\end{equation}
For such functions, property~(5) in Theorem~\ref{thm-10.14} we assume for~$M^D$ yields:

\begin{myproposition}
\label{prop-3}
Fix~$\beta>0$ and~$R>0$. For each $\epsilon>0$ there are $\delta_0>0$ and $\lambda_0>0$ such that, for all $\lambda\in(0,\lambda_0)$, all $\delta\in(0,\delta_0)$ and all~$f\in\FF_{R,\beta,\delta}^T$,
\begin{multline}
\qquad
(1-\epsilon)\hat c\int_{T_\delta}f(x)r_T(x)^2\,\textd x\le
\frac{\log E(\texte^{-\lambda M^T(f\1_{T_\delta})})}{\lambda\log\lambda}
\\
\le(1+\epsilon)\hat c\int_{T_\delta}f(x)r_T(x)^2\,\textd x,
\qquad
\end{multline}
where $M^T(f\1_{T_\delta}):=\int_{T_\delta}M^T(\textd x)\,f(x)$.
\end{myproposition}

\begin{proofsect}{Proof (some ideas)}
This builds on Exercise~\ref{ex:10.19} and a uniform approximation of~$f$ by functions that are piecewise constant on small subtriangles of~$T$. Thanks to the known scaling properties of the $Z^D$ measure (see Theorem~\ref{thm-10.14}(2)) and the conformal radius, it suffices to prove this just for a unit triangle. We refer the reader to \cite[Lemma~6.8]{BL2} for further details and calculations. 
\end{proofsect}

As noted above, whenever $f\in\FF_{R,\beta,\delta}$ and $A_{K,R}^i$ occurs, Proposition~\ref{prop-3} may be applied (with perhaps slightly worse values of~$R$ and~$\beta$) to the test function
\begin{equation}
x\mapsto f(x)\texte^{\alpha(\Phi^{D,\wt D}(x)-\Phi^{D,\wt D}(x_i))}
\end{equation}
because the harmonicity of $x\mapsto \Phi^{D,\wt D}(x)$ turns the uniform bound on oscillation into a Lipschitz property. Thus, for $\lambda:=K^{-4}\texte^{\alpha\Phi^{D,\wt D}(x_i)}$, on $A_{K,R}^i$ we then get
\begin{equation}
\label{E:10.58ui}
\begin{aligned}
E\biggl(&\,\exp\Bigl\{-\texte^{\alpha\Phi^{D,\wt D}(x_i)}M^{T^i}(f\1_{T_\delta^i}\texte^{\alpha(\Phi^{D,\wt D}-\Phi^{D,\wt D}(x_i))})\Bigr\}\,\bigg|\,\Phi^{D,\wt D}\biggr)
\\
&=\exp\biggl\{\hat c(1+\tilde\epsilon)\,\log\bigl(K^{-4}\texte^{\alpha\Phi^{D,\wt D}(x_i)}\bigr)\int_{T_\delta^i}\textd x\,f(x)\,\texte^{\alpha\Phi^{D,\wt D}(x)}\,r_{T^i}(x)^2\biggr\}
\end{aligned}
\end{equation}
for some random~$\tilde\epsilon\in[-\epsilon,\epsilon]$ depending only on~$\Phi^{D,\wt D}$, where we also used that, by \eqref{E:10.53ua} and $\alpha 2\sqrt g=4$, we have $\lambda\le(\log K)^{-\alpha c}$ which tends to zero as~$K\to\infty$.

Denote by $Z^D_{K,R,\delta}$ the measure in \eqref{E:10.45} and write $M^D_{K,R,\delta}$ for the expression on the right of \eqref{E:10.47} with the sum restricted to~$i$ where~$A_{K,R}^i$ occurs. Note that, by \eqref{E:10.53uia} and standard (by now) estimates for~$C^{D,\wt D}$,
\begin{equation}
\label{E:10.59ui}
\log(K^{-4}\texte^{\alpha\Phi^{D,\wt D}(x_i)})
=\bigl(1+\tilde\epsilon'(x)\bigr)\alpha\bigl(\Phi^{D,\wt D}(x)-\alpha\Var(\Phi^{D,\wt D}(x))\bigr)
\end{equation}
with~$\tilde\epsilon'(x)\in[-\epsilon,\epsilon]$ for all~$x\in T^i_\delta$, provided~$R$ is large enough. Recalling also that
\begin{equation}
\label{E:10.60nwwt}
r_{T^i}(x)^2=r_D(x)^2\texte^{-\frac12\alpha^2\Var(\Phi^{D,\wt D}(x))},\qquad x\in T^i,
\end{equation}
from \twoeqref{E:10.58ui}{E:10.59ui} we obtain
\begin{equation}
E\bigl(\texte^{-(1+2\epsilon) Z^D_{K,R,\delta}(f)}\bigr)\le
E\bigl(\texte^{-M^D_{K,R,\delta}(f)}\bigr)\le E\bigl(\texte^{-(1-2\epsilon) Z^D_{K,R,\delta}(f)}\bigr).
\end{equation}
Since $M^D_{K,R,\delta}(f)$ tends in distribution to $M^D(f)$ in the stated limits, the law of~$M^D(f)$ is given by the corresponding limit law of~$ Z^D_{K,R,\delta}(f)$, which must therefore exist as well.
\end{proofsect}

As an immediate consequence of Theorem~\ref{thm-10.16}, we then get:

\begin{mycorollary}[Behavior under rigid rotations]
For each~$a,b\in\C$,
\begin{equation}
Z^{a+bD}(a+b\textd x)\,\,\laweq\,\,|b|^4\,Z^D(\textd x).
\end{equation}
\end{mycorollary}

\begin{proofsect}{Proof}
By Theorem~\ref{thm-10.14}(2), we just need to prove this for~$a:=0$ and~$|b|=1$. This follows from Theorem~\ref{thm-10.16} and the fact that both the law of~$\Phi^{D,\wt D}$ and the conformal radius~$r_D$ are invariant under the rigid rotations of~$D$ (and~$\wt D$).
\end{proofsect}

Somewhat more work is required to prove:

\begin{mytheorem}[Behavior under conformal maps]
Let $f\colon D\to f(D)$ be a conformal bijection between admissible domains $D,f(D)\in\mathfrak D$. Then
\begin{equation}
Z^{f(D)}\circ f(\textd x)\,\,\laweq\,\,\bigl|f'(x)\bigr|^4\,Z^D(\textd x).
\end{equation}
In particular, the law of $r_D(x)^{-4}\,Z^D(\textd x)$ is invariant under conformal maps.
\end{mytheorem}

This follows, roughly speaking, by the fact that the law of~$\Phi^{D,\wt D}$ is conformally invariant and also by the observation that a conformal map is locally the composition of a dilation with a rotation. In particular, the triangles~$T^i$ map to near-triangles~$f(T^i)$ with the deformation tending to zero with the size of the triangle. See the proof of \cite[Theorem~7.2]{BL2}.

\section{Connection to Liouville Quantum Gravity}
\label{sec-10.5}\noindent
We will now move to the question of direct characterization of the law of~$Z^D$-measures. As we will show, $Z^D$ has the law of a critical Liouville Quantum Gravity associated with the continuum GFF. These are versions of the measures from Lemma~\ref{lemma-GMC} for a critical value of parameter~$\beta$ which in our case corresponds to~$\beta_\cc=\alpha$. At this value, the martingale argument from Lemma~\ref{lemma-GMC} produces a vanishing measure and so additional manipulations are needed to obtain a non-trivial limit object. 

One way to extract a non-trivial limit measure at~$\beta_\cc$ is to replace the exponential factor in \eqref{E:2.17a} by its negative derivative; cf~\twoeqref{E:10.49}{E:10.50}. The existence of the corresponding limit was proved in~\cite{DRSV1}. Another possibility is to re-scale the approximating measures so that a non-trivial limit is achieved. This goes by the name \emph{Seneta-Heyde norming} as discussed in a 2014 paper by Duplantier, Rhodes, Sheffield and Vargas~\cite{DRSV2}. A technical advantage of the scaling over taking the derivative is the fact that Kahane's convexity inequality (Proposition~\ref{prop-5.6}) remains applicable in that case.

We will use approximating measures based on the \emph{white-noise decomposition} of the CGFF. The specifics are as follows: For~$\{B_t\colon t\ge0\}$ the standard Brownian motion, let $p^D_t(x,y)$ be the transition density from~$x$ to~$y$ before exiting~$D$. More precisely, letting $\tau_{D^\cc}:=\inf\{t\ge0\colon B_t\not\in D\}$ we have
\begin{equation}
p^D_t(x,y)\textd y:=P^x\bigl(B_t\in\textd y,\,\tau_{D^\cc}>t\bigr).
\end{equation}
Writing $W$ for the white noise on $D\times(0,\infty)$ with respect to the Lebesgue measure,  consider the Gaussian process~$t,x\mapsto\varphi_t(x)$ defined by
\begin{equation}
\varphi_t(x):=\int_{D\times[\texte^{-2t},\infty)}p^D_{s/2}(x,z)W(\textd z\,\textd s).
\end{equation}
The Markov property of~$p^D$ along with reversibility give
\begin{equation}
\begin{aligned}
\Cov\bigl(\varphi_t(x),\varphi_t(y)\bigr)&=\int_{D\times[\texte^{-2t},\infty)}\textd z\otimes\textd s\,\,p^D_{s/2}(x,z)p^D_{s/2}(y,z)
\\
&=\int_{[\texte^{-2t},\infty)}\textd s\,\,p^D_s(x,y)
\,\,\underset{t\to\infty}\longrightarrow\,\wh G^D(x,y)\end{aligned}
\end{equation}
and so~$\varphi_t$ tends in law to the CGFF.

Define the random measure (and compare with~\eqref{E:2.17a})
\begin{equation}
\mu^{D,\alpha}_t(\textd x):=\sqrt{t}\,1_D(x)\,\texte^{\alpha\varphi_t(x)-\frac12\alpha^2\Var[\varphi_t(x)]}\,\textd x\,.
\end{equation}
The key point to show that the scale factor~$\sqrt t$ ensures that the $t\to\infty$ limit will produce a non-vanishing and yet finite limit object. This was proved in Duplantier, Rhodes, Sheffield and Vargas~\cite{DRSV2}:

\begin{mytheorem}[Critical GMC measure]
\label{thm-DRSV}
There is a non-vanishing a.s.-finite random Borel measure~$\mu^{D,\alpha}_\infty$ such that for every Borel set~$A\subseteq D$,
\begin{equation}
\mu^{D,\alpha}_t(A)\,\,\underset{t\to\infty}\longrightarrow\,\, \mu^{D,\alpha}_\infty(A),\quad\text{\rm in probability.}
\end{equation}
\end{mytheorem}

\begin{figure}[t]
\centerline{\includegraphics[width=0.65\textwidth]{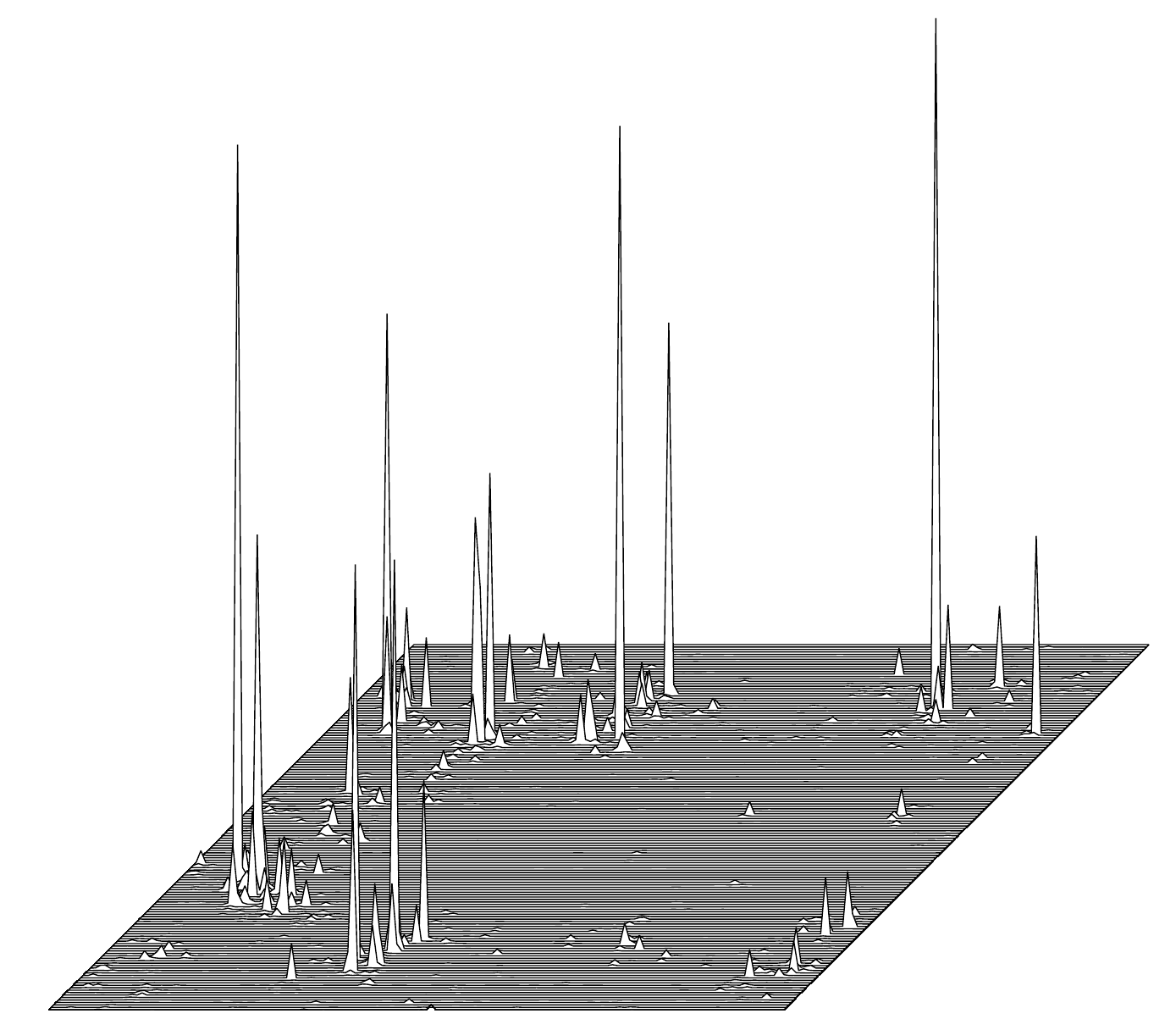}
}
\vglue2mm
\begin{quote}
\small 
\caption{
\label{fig-critical-LQG}
\small
A sample of critical Liouville Quantum Gravity measure over a square domain. Despite the occurence of prominent spikes, the measure is non-atomic although just barely so as it is (conjecturally) supported on a set of zero Hausdorff dimension.}
\normalsize
\end{quote}
\end{figure}

\noindent
We call the measure $\mu^{D,\alpha}_\infty$ the \emph{critical GMC} in~$D$. It is a fact that $E\mu^{D,\alpha}_\infty(A)=\infty$ for any non-empty open subset of~$D$ so, unlike the subcritical cases $\beta=\lambda\alpha$ for~$\lambda\in(0,1)$, the overall normalization of this measure cannot be fixed by its expectation. 

Alternative definitions of $\mu^{D,\alpha}_\infty$ have appeared in the meantime which, through the contributions of Rhodes and Vargas~\cite{RV-review},  Powell~\cite{Powell} and Junnila and Saksman~\cite{JS}, are now known to be all equal up to a multiplicative (deterministic) constant. The measure in \eqref{E:10.45} is yet another example of this kind, although the abovementioned uniqueness theorems do not apply to this case. Notwithstanding, one is able to call on Theorem~\ref{thm-10.14} instead and get:

\begin{mytheorem}[Identification with critical LQG measure]
\label{thm-BLnew}
The family of measures
\begin{equation}
\label{E:10.77}
r_D(x)^2\mu^{D,\alpha}_\infty(\textd x),\quad D\in\mathfrak D\,,
\end{equation}
 constructed in Theorem~\ref{thm-DRSV} obeys conditions~(1-5) of Theorem~\ref{thm-10.14} for some $\hat c>0$. In particular, there is a constant~$c\in(0,\infty)$ such that
\begin{equation}
Z^D(\textd x)\,\,\laweq\,\, c\, r_D(x)^2\,\mu^{D,\alpha}_\infty(\textd x),\quad D\in\mathfrak D.
\end{equation}
\end{mytheorem}

\noindent
This appears as \cite[Theorem~2.9]{BL2}. The key technical challenge is to prove that the measure in \eqref{E:10.77} obeys the Laplace-transform asymptotic \eqref{E:1.25ue}. This is based on a version of the concentric decomposition and approximations via Kahane's convexity inequality. We refer the reader to \cite{BL2} for further details.

\smallskip
We conclude by reiterating that many derivations in the present lecture were originally motivated by the desire to prove directly the uniqueness of the subsequential limit of $\eta^D_{N,r_N}$ (as we did for the intermediate level sets) based only on the easily obtainable properties of the~$Z^D$-measure. This program will partially be completed in our proof of Theorem~\ref{thm-Z-unique} although there we still rely on the Laplace transform tail \eqref{E:1.25ue} which also underpins the limit of the DGFF maximum. We believe that even this step can be bypassed and the reference to convergence of the DGFF maximum avoided; see Conjecture~\ref{conj-Laplace-tail}.


\chapter{Local structure of extremal points}
\label{lec-11}\noindent
In this lecture we augment the conclusions obtained for the point processes associated with extremal local maxima to include information about the local behavior. The proofs are based on the concentric decomposition of the DGFF and entropic-repulsion arguments developed for the proof of tightness of the absolute maximum. Once this is done, we give a formal proof of our full convergence result from Theorem~\ref{thm-extremal-vals}. A number of  corollaries are presented that concern the cluster-process structure of the extremal level sets, a Poisson-Dirichlet limit of the Gibbs measure associated with the DGFF, the Liouville Quantum Gravity in the  glassy phase and the freezing phenomenon.

\section{Cluster at absolute maximum}
\noindent
Our focus in this lecture is on the local behavior of the field near its near-maximal local extrema. For reasons mentioned earlier, we will refer of these values as a \emph{cluster}. As it turns out, all that will technically be required is the understanding of the cluster associated with the absolute maximum:

\begin{mytheorem}[Cluster law]
\label{thm-11.1}
Let~$D\in\mathfrak D$ with~$0\in D$ and let~$\{D_N\colon N\ge1\}$ be an admissible sequence of approximating domains. Then for each~$t\in\R$ and each function $f\in C_\cc(\overline\R^{\Z^2})$ depending only on a finite number of coordinates,
\begin{equation}
\label{E:11.1}
E\Bigl(\,f\bigl(h^{D_N}_0-h^{D_N}\bigr)\,\Big|\, h^{D_N}_0=m_N+t,\,h^{D_N}\le h^{D_N}_0\Bigr)
\,\,\underset{N\to\infty}\longrightarrow\, E_\nu(f)\,,
\end{equation}
where~$\nu$ is a probability measure on~$[0,\infty)^{\Z^2}$ defined from~$\phi:=$ DGFF on~$\Z^2\smallsetminus\{0\}$ via the weak limit
\begin{equation}
\label{E:11.2}
\nu(\cdot):=\lim_{r\to\infty} P\biggl(\phi+\frac2{\sqrt g}\fraka\in\cdot\,\bigg|\,\phi_x+\frac2{\sqrt g}\fraka(x)\ge0\colon |x|\le r\biggr)
\end{equation}
in which (we recall)~$\fraka$ denotes the potential kernel on~$\Z^2$.
\end{mytheorem}

The existence of the weak limit in \eqref{E:11.2} is part of the statement of the theorem. Remarkably, the convergence alone may be inferred by much softer~means:

\begin{myexercise}
Let~$\nu_r$ be the conditional measure on the right of \eqref{E:11.2}. Prove that $r\mapsto\nu_r$ is stochastically increasing. [Hint: Under this type of conditioning, the same holds true for any strong-FKG  measure.]
\end{myexercise}

\noindent
The exercise shows that~$r\mapsto\nu_r(A)$ is increasing on increasing events and so the limit in~\eqref{E:11.2} exists for any event depending on a finite number of coordinates. The problem is that~$\nu_r$ is \emph{a priori} a measure on $[0,\infty]^{\Z^2}$ and the interpretation of the limit as a distribution on~$[0,\infty)^{\Z^2}$ requires a proof of tightness. This  additional ingredient will be supplied by our proof of Theorem~\ref{thm-11.1}.

\nopagebreak
\begin{figure}[t]
\vglue-1mm
\centerline{\includegraphics[width=0.75\textwidth]{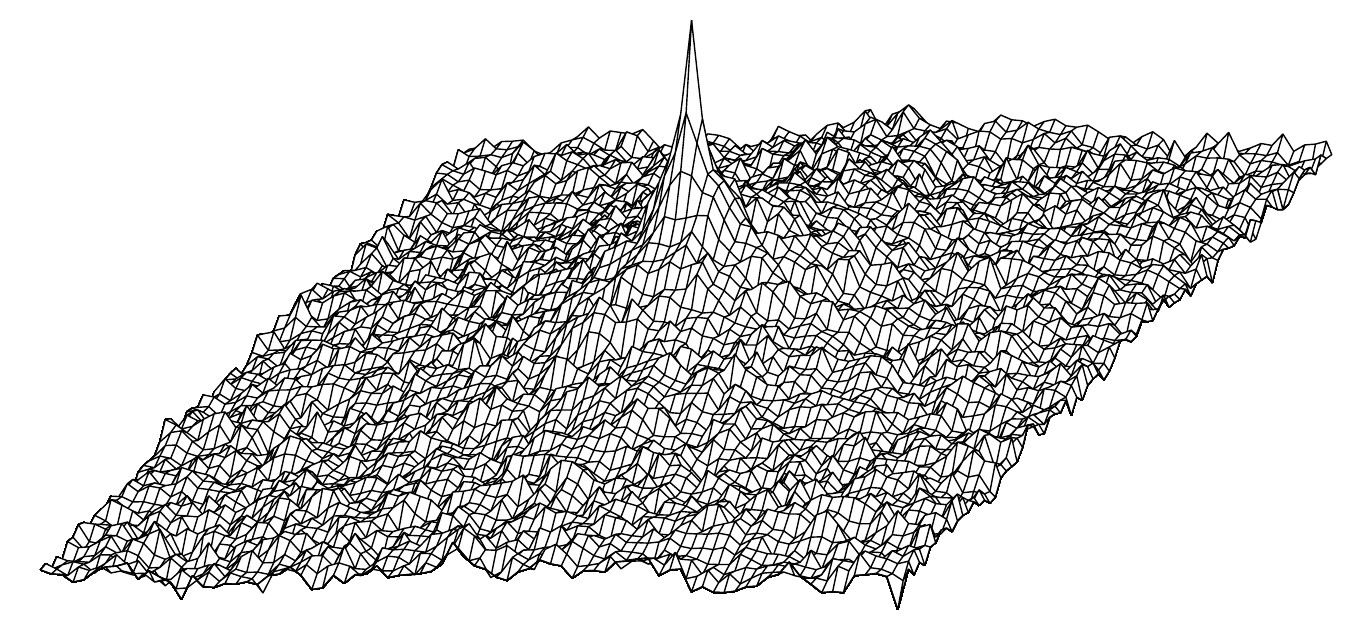}
}
\vglue0mm
\begin{quote}
\small 
\vglue2mm
\caption{
\label{fig-cluster}
\small
A sample of the configuration of the DGFF in the vicinity of its (large) local maximum.}
\normalsize
\end{quote}
\end{figure}

That the limit takes the form in \eqref{E:11.2} may be explained via a simple heuristic calculation. Indeed, by Lemma~\ref{lemma-8.3}, conditioning the field $h^{D_N}$ on $h^{D_N}_0=m_N+t$ effectively shifts the mean of $h^{D_N}_0-h^{D_N}_x$ by a quantity with the asymptotic
\begin{equation}
\label{E:11.3}
(m_N+t)\bigl(1-\frakg^{D_N}(x)\bigr)\,\underset{N\to\infty}\longrightarrow\frac2{\sqrt g}\fraka(x).
\end{equation}
A variance computation then shows that the law of $x\mapsto h^{D_N}_0-h^{D_N}_x$ tends, in the sense of finite-dimensional distributions, to $x\mapsto\phi_x+\frac2{\sqrt g}\fraka(x)$, where~$\phi$ is the pinned DGFF; see Fig.~\ref{fig-cluster}. The conditioning on the origin being the maximum then translates into the conditioning in \eqref{E:11.2}. This would more or less suffice to prove the desired result were it not for the following fact:

\begin{mytheorem}
\label{thm-11.3}
There exists~$c_\star\in(0,\infty)$ such that
\begin{equation}
\label{E:11.4}
P\biggl(\phi_x+\frac2{\sqrt g}\fraka(x)\ge0\colon |x|\le r\biggr)=\frac{c_\star}{\sqrt{\log r}}\bigl(1+o(1)\bigr),\quad r\to\infty.
\end{equation}
\end{mytheorem}

\noindent
Since the right-hand side of \eqref{E:11.4} vanishes in the limit~$r\to\infty$, the conditioning in \eqref{E:11.2} is increasingly singular and so it is hard to imagine that one could control the limit in \eqref{E:11.1} solely via manipulations based on weak convergence.

\section{Random walk based estimates} 
\noindent
Our proof of Theorem~\ref{thm-1.11} relies on the concentric decomposition of the DGFF developed in Sections~\ref{sec:8.2}--\ref{sec8.4}. As these sections were devoted to the proof of tightness of the lower tail of the maximum, we were not allowed to assume any bounds on the lower tails in estimates there. However, with the tightness settled, Lemma~\ref{lemma-8.21} augments Lemma~\ref{lemma-8.11} to a two-sided estimate:

\begin{mylemma}
\label{lemma-10.4a}
There is~$a>0$ such that each~$k=1,\dots,n$ and each~$t\ge0$,
\begin{equation}
\label{E:8.39a}
P\biggl(\,\Bigl|\max_{x\in\Delta^k\smallsetminus\Delta^{k-1}}\bigl[\chi_{k+1}(x)+\chi_k(x)+h'_k(x)\bigr]- m_{2^k}\Bigr|\ge t\biggr)
\le\texte^{-a t}.
\end{equation}
\end{mylemma}

\noindent
This allows for control of the deviations of the field $h^{D_N}$ from the random walk~$-S_k$ in both directions which upgrades Lemma~\ref{lemma-8.16} to the form: 

\begin{mylemma}[Reduction to random walk events]
\label{lemma-10.5}
Assume~$h^{D_N}$ is realized as the sum on the right of \eqref{E:8.29}. There is a numerical constant~$C>0$ such that uniformly in the above setting, the following holds for each~$k=1,\dots,n$ and each~$t\in\R$:
\begin{multline}
\label{E:11.6}
\quad
\{S_{n+1}=0\}\cap\bigl\{S_k\ge R_K(k)+|t|\bigr\}
\\
\subseteq
\{h^{D_N}_0=0\}\cap\bigl\{h^{D_N}\le (m_N+t)(1-\frakg^{D_N})\text{\rm\ on }\Delta^k\smallsetminus\Delta^{k-1}\bigr\}
\\
\subseteq\{S_{n+1}=0\}\cap\bigl\{S_k\ge -R_K(k)-|t|\bigr\}\,,
\quad
\end{multline}
where $K$ is the control variable from Definition~\ref{def-control-var} with the absolute value signs added around the quantity on the left of \eqref{E:8.62nwwt} and
\begin{equation}
R_k(\ell):=C[1+\theta_{n,k}(\ell)],
\end{equation}
with~$\theta_{n,k}$ as in \eqref{E:def-theta} and~$C$ as in Lemma~\ref{lemma-8.16}. (We recall that, here and henceforth, $n$ is the largest integer such that $\{x\in\Z^2\colon|x|_\infty\le 2^{n+1}\}\subseteq D_N$.)
\end{mylemma}

\begin{proofsect}{Proof}
From \eqref{E:8.65nwwt2} we get, for all $k=0,\dots,n$ (and with~$\Delta^{-1}:=\emptyset$),
\begin{multline}
\quad
\Bigl|\,\max_{x\in\Delta^k\smallsetminus\Delta^{k-1}}\bigl[h^{D_N}_x-(m_N+t)(1-\frakg^{D_N}(x))\bigr]+S_k\Bigr|
\\
\le \sum_{\ell=k}^n \max_{x\in\Delta^k\smallsetminus\Delta^{k-1}}\bigl|\frakb_\ell(x)\bigr|\bigl|\varphi_\ell(0)\bigr|
+\biggl(\,\sum_{\ell=k+2}^n\max_{x\in\Delta^k\smallsetminus\Delta^{k-1}}\bigl|\chi_\ell(x)\bigr|\biggr)
\\+
\Bigl|\,\max_{x\in\Delta^k\smallsetminus\Delta^{k-1}}\bigl[\chi_{k+1}(x)+\chi_k(x)+h_k'(x)-m_{2^k}\bigr]\Bigr|
\\+\max_{x\in\Delta^k\smallsetminus\Delta^{k-1}}\Bigl|\,m_N(1-\frakg^{D_N}(x))-m_{2^k}\Bigr|+|t|.
\quad
\end{multline}
The definition of~$K$ then bounds the first three terms on the right by a quantity of order~$1+\theta_{n,K}(k)$. For the second to last term, here instead of \eqref{E:8.65nwwt} we need:

\begin{myexercise}
\label{ex:11.6}
There is $c>0$ such that for all $n\ge1$ and all~$k=0,\dots,n$,
\begin{equation}
\max_{x\in\Delta^k\smallsetminus\Delta^{k-1}}\Bigl|\,m_N(1-\frakg^{D_N}(x))-m_{2^k}\Bigr|\le c\bigl[1+\log(1+k\wedge(n-k))\bigr].
\end{equation}
\end{myexercise}
\noindent
The inclusions \eqref{E:11.6} follow readily from this.
\end{proofsect}

\noindent
We will now use the random walk $\{S_1,\dots,S_n\}$ to control all important aspects of the conditional expectation in the statement of Theorem~\ref{thm-11.1}.

First note that the event $\bigcap_{k=1}^n\{S_k\ge -R_K(k)-|t|\}$ encases all of the events of interest and so we can use it as the basis for estimates of various undesirable scenarios. (This is necessary because the relevant events will have probability tending to zero proportionally to~$1/n$.) In particular, we need to upgrade Lemma~\ref{lemma-8.19} to the form:

\begin{mylemma}
\label{lemma-11.7}
There are constants~$c_1,c_2>0$ such that for all~$n\ge1$, all~$t$ with~$0\le t\le n^{1/5}$ and all~$k=1,\dots,\lfloor \ffrac n2\rfloor$,
\begin{equation}
P\biggl(\{K>k\}\cap \bigcap_{\ell=1}^{n}\{S_\ell\ge -R_k(\ell)-t\}\,\bigg|\, S_{n+1}=0\biggr)\le c_1\frac{1+t^2}n\texte^{-c_2(\log k)^2}.
\end{equation}
\end{mylemma}

\noindent
Since the target decay is order-$1/n$, this permits us to assume~$\{K\le k\}$ for~$k$ sufficiently large but independent of~$n$ whenever need arises in the forthcoming derivations. Lemma~\ref{lemma-8.18} then takes the form:

\begin{mylemma}[Entropic repulsion]
\label{lemma-11.6}
For each~$t\ge0$ there is~$c>0$ such that for all $n\ge1$ and all~$k=1,\dots,\lfloor \ffrac n2\rfloor$
\begin{multline}
P\biggl(\,\{S_k,S_{n-k}\ge k^{1/6}\}\cap\bigcap_{\ell=k+1}^{n-k-1}\{S_\ell\ge R_k(\ell)+t\bigr\}
\\\bigg|\, \bigcap_{\ell=1}^{n}\{S_\ell\ge -R_k(\ell)-t\}\cap\{S_{n+1}=0\}\biggr)\ge 1-ck^{-\frac1{16}}.
\end{multline}
\end{mylemma}

\noindent
We will not give formal proofs of Lemmas~\ref{lemma-11.7}--\ref{lemma-11.6} here; instead we refer the reader to~\cite[Section~4]{BL3}.

Consider the expectation in the statement of Theorem~\ref{thm-11.1}. Lemma~\ref{lemma-8.3} permits us to shift the conditioning event to $h^{D_N}_0=0$ at the cost of adding $(m_N+t)\frakg^{D^N}$ to all occurrences of the field. Abbreviating
\begin{equation}
\label{E:mNtx}
m_N(t,x):=(m_N+t)(1-\frakg^{D_N}(x)),
\end{equation}
the expectation in \eqref{E:11.1} can thus be written as
\begin{equation}
\label{E:11.11}
\frac{E\Bigl(\,f\bigl(m_N(t,\cdot)-h^{D_N}\bigr)\1_{\{h^{D_N}\le m_N(t,\cdot)\}}\,\Big|\, h^{D_N}_0=0\,\Bigr)}
{E\Bigl(\,\1_{\{h^{D_N}\le m_N(t,\cdot)\}}\,\Big|\, h^{D_N}_0=0\,\Bigr)}.
\end{equation}
We will control this ratio by deriving, separately, asymptotic expressions for the numerator and the denominator which (in light of Lemmas~\ref{lemma-10.4a} and~\ref{lemma-11.7}) will both decay as a constant over~$n$. As both the numerator and the denominator have the same structure, it suffices to focus on the numerator. We claim:

\begin{myproposition}
\label{prop-11.8}
For each~$\epsilon>0$ and each~$t_0>0$ there is~$k_0\ge1$ such that for all~$k$ with~$k_0\le k\le n^{1/6}$ and all $t\in[-t_0,t_0]$,
\begin{multline}
\label{E:11.12ueu}
\Biggl|E\Bigl(\,f\bigl(m_N(t,\cdot)-h^{D_N}\bigr)\1_{\{h^{D_N}\le m_N(t,\cdot)\}}\,\Big|\, h^{D_N}_0=0\,\Bigr)
\\
-E\biggl(\,f\bigl(\tfrac2{\sqrt g}\fraka+\phi_k\bigr)\1_{\{\phi_k+\frac2{\sqrt g}\fraka\ge0\text{\rm\ in }\Delta^k\}}\1_{\{S_k,S_{n-k}\in[k^{1/6},\,k^2]\}}\Bigl(\,\prod_{\ell=k}^{n-k}1_{\{S_\ell\ge0\}}\Bigr)
\\
\times\1_{\{h^{D_N}\le m_N(t,\cdot)\text{\rm\ in }D_N\smallsetminus\Delta^{n-k}\}}\bigg|\, h^{D_N}_0=0\,\biggr)
\Biggr|\le\frac{\epsilon}n\,,
\end{multline}
where we used the shorthand
\begin{equation}
\label{E:11.13}
\phi_k(x):=h^{\Delta^k}_0-h^{\Delta^k}_x.
\end{equation}
\end{myproposition}

\begin{proofsect}{Proof (sketch)}
The proof is based on a sequence of replacements that gradually convert one expectation into the other.
First we separate scales by noting that,  given~$k\in\{1,\dots,\lfloor n/2\rfloor\}$, we can write the ``hard'' event in the expectation as the intersection of ``inner'', ``middle'' and ``outer'' events,
\begin{multline}
\quad
\1_{\{h^{D_N}\le m_N(t,\cdot)\}}
=\1_{\{h^{D_N}\le m_N(t,\cdot)\text{ in }\Delta^k\}}
\\
\times
\1_{\{h^{D_N}\le m_N(t,\cdot)\text{ in }\Delta^{n-k}\smallsetminus\Delta^k\}}
\1_{\{h^{D_N}\le m_N(t,\cdot)\text{ in }D_N\smallsetminus\Delta^{n-k}\}}.
\quad
\end{multline}
Plugging this in we can also use Lemma~\ref{lemma-11.7} to insert the indicator of $\{K\le k\}$ into the expectation. The inclusions in \eqref{E:11.6} permit us to replace the ``middle'' event
\begin{equation}
\bigl\{h^{D_N}\le m_N(t,\cdot)\text{ in }\Delta^{n-k}\smallsetminus\Delta^k\bigr\}
\end{equation}
by
\begin{equation}
\bigcap_{\ell=k}^{n-k}\{S_\ell\ge\pm(R_k(\ell)+|t|)\}
\end{equation}
with the sign depending on whether we aim to get upper or lower bounds. Lem\-ma~\ref{lemma-11.6} then tells us that the difference between these upper and lower bounds is negligible, and so we may further replace $\{S_\ell\ge\pm(R_k(\ell)+|t|)\}$ by $\{S_\ell\ge0\}$. In addition, by Lemma~\ref{lemma-11.6} we may also assume $S_k,S_{n-k}\ge k^{1/6}$. The bounds $S_k,S_{n-k}\le k^2$ then arise (for~$k$ large enough) from the restriction to~$\{K\le k\}$ and the inequalities in Definition~\ref{def-control-var}.

At this point we have replaced the first expectation in \eqref{E:11.12ueu} by $1+O(k^{-\frac1{16}})$-times the conditionalexpectation
\begin{multline}
E\biggl(f\bigl(m_N(t,\cdot)-h^{D_N}\bigr)
\1_{\{K\le k\}}\1_{\{h^{D_N}\le m_N(t,\cdot)\text{ in }\Delta^k\}}\1_{\{S_k,S_{n-k}\in[k^{1/6},\,k^2]\}}
\\
\times\Bigl(\prod_{\ell=k}^{n-k}1_{\{S_\ell\ge0\}}\Bigr)
\1_{\{h^{D_N}\le m_N(t,\cdot)\text{ in }D_N\smallsetminus\Delta^{n-k}\}}\,\bigg|\, h^{D_N}_0=0\,\biggr)
\end{multline}
plus a quantity of order $\frac{1+t^2}n\texte^{-c_2(\log k)^2}$.
Next we will use the continuity of~$f$ to replace $m_N(t,\cdot)$ in the argument of~$f$ by its limit value $\frac2{\sqrt g}\fraka$. For this we assume that~$k$ is much larger than the diameter of the set of vertices that~$f$ depends on. Conditional on $h^{D_N}_0=0$ we then have
\begin{equation}
\label{E:11.17}
h^{D_N}_x = -\phi_k(x)+\sum_{\ell>k}\bigl[\frakb_\ell(x)\varphi_\ell(0)+\chi_\ell\bigr],\quad x\in \Delta^k.
\end{equation}
The bounds arising from the restriction $K\le k$ then let us replace~$h^{D_N}$ in the argument of~$f$ by~$-\phi_k$. 

It remains to deal with the indicator of the ``inner'' event
\begin{equation}
\bigl\{h^{D_N}\le m_N(t,\cdot)\text{ in }\Delta^k\bigr\}
\end{equation}
which we want to replace by
\begin{equation}
\bigl\{\phi_k+ \tfrac2{\sqrt g}\fraka\ge0\text{ in }\Delta^k\bigr\}\,.
\end{equation}
Here continuity arguments cannot be used; instead, we have to show that entropic repulsion creates a sufficiently large gap between $h^{D_N}$ and $m_N(t,\cdot)$ in the first expectation in \eqref{E:11.12ueu}, and between $\phi_k+\frac2{\sqrt g}\fraka$  and zero in the second expectation, to absorb the sum on the right of \eqref{E:11.17} and the difference between $m_N(t,\cdot)$ and~$\frac2{\sqrt g}\fraka$ everywhere on~$\Delta^k$. We refer to \cite[Lemma~4.22]{BL3} (and Exercise~\ref{ex:11.9} below) for details of this step.

Even after all the replacements above have been done, the quantity under expectation remains concentrated on $\bigcap_{\ell=1}^n\{S_\ell\ge-R_k(\ell)-|t|\}$. Lemma~\ref{lemma-11.6} then permits us to drop the restriction to $\{K\le k\}$ and get the desired result.
\end{proofsect}

The entropic repulsion step from \cite[Lemma~4.22]{BL3} relies on the following general fact which constitutes a nice exercise for the reader:

\begin{myexercise}[Controlling the gap]
\label{ex:11.9}
Prove that, if $\phi$ is a field on a set~$\Lambda$ with the strong-FKG property, then for all~$\delta>0$ and all $f\colon \Lambda\to\R$ (and writing $\{\phi\not\ge f\}$ for $\{\phi(x)\ge f(x)\colon x\in\Lambda\}^\cc$)
\begin{equation}
\label{E:11.19nwt}
P\bigl(\phi\not\ge f+\delta\,\big|\,\phi\ge f\bigr)
\le\sum_{x\in\Lambda}P\bigl(\phi(x)< f(x)+\delta\,\big|\,\phi(x)\ge f(x)\bigr).
\end{equation}
For $\phi(x)=\NN(0,\sigma^2)$ with $\sigma^2>0$, and $f(x)\le0$, show also that 
\begin{equation}
\label{E:11.20nwt}
P\bigl(\phi(x)\not\ge f(x)+\delta\,\big|\,\phi(x)\ge f(x)\bigr)\le\frac2{\sqrt{2\pi}}\frac\delta\sigma\le\frac\delta\sigma.
\end{equation}
\end{myexercise}

\noindent
Thanks to \twoeqref{E:11.19nwt}{E:11.20nwt}, in order to show the existence of a gap it suffices to prove a uniform lower bound on the variance of the relevant fields. (Note that the gap can be taken to zero slowly with~$n\to\infty$.)

Moving back to the main line of argument underlying the proof of Theorem~\ref{thm-11.1}, next we observe:

\begin{mylemma} 
\label{lemma-11.10ueu}
Conditionally on~$S_k$ and~$S_{n-k}$ and the event~$S_{n+1}=0$, 
\begin{enumerate}
\item[(1)] the ``inner'' field $\phi_k$,
\item[(2)] the random variables $\{S_\ell\colon \ell=k,\dots,n-k\}$, and
\item[(3)] the ``outer'' field $\{h^{D_N}_x\colon x\in D_N\smallsetminus\Delta^{n-k}\}$
\end{enumerate}
are independent of one another whenever~$n>2k$.
\end{mylemma}

\begin{proofsect}{Proof}
Inspecting \eqref{E:11.13}, $\phi_k$ is measurable with respect to
\begin{equation}
\sigma\bigl(\varphi_1(0),\dots,\varphi_k(0),\chi_1,\dots,\chi_k,h_0,\dots,h_k'\bigr),
\end{equation}
while the random variables in~(2) are measurable with respect to
\begin{equation}
\sigma\bigl(\{S_k\}\cup\{\varphi_\ell(0)\colon \ell=k+1,\dots,n-k\}\bigr).
\end{equation}
The definition of the concentric decomposition in turn ensures that the random variables in~(3), as well as the event~$\{S_{n+1}=0\}$, are measurable with respect to
\begin{equation}
\sigma\Bigl(\{S_{n-k}\}\cup\bigcup_{\ell=n-k}^n\{\varphi_\ell(0),\chi_\ell,h'_\ell\}\Bigr).
\end{equation}
The claim follows from the independence of the random objects constituting the concentric decomposition.
\end{proofsect}

We note that the Lemma~\ref{lemma-11.10ueu} is the prime reason why we had to replace $h^{\Delta^n}-h^{\Delta^n}_0$ by~$\phi_k$ in the ``inner'' event. (Indeed, $h^{\Delta^n}-h^{\Delta^n}_0$ still depends on the fields~$\chi_\ell$, $\ell=n-k,\dots,n$.) Our next step is to take expectation conditional on~$S_k$ and~$S_{n-k}$. Here we will use:

\begin{mylemma}
\label{lemma-11.9}
For each~$t_0>0$ there is~$c>0$ such that for all~$1\le k\le n^{1/6}$,
\begin{equation}
\label{E:11.24nw}
\biggl|\,P\Bigl(\,\bigcap_{\ell=k}^{n-k}\{S_\ell\ge0\}\,\Big|\,\sigma(S_k,S_{n-k})\Bigr)-\frac2{g\log 2}\frac{S_kS_{n-k}}n\biggr|\le c\frac{k^4}n\frac{S_kS_{n-k}}n
\end{equation}
holds everywhere on $\{S_k,S_{n-k}\in[k^{1/6},k^2]\}$.
\end{mylemma}

\begin{proofsect}{Proof (idea)}
We will only explain the form of the leading term leaving the error to a reference to \cite[Lemma~5.6]{BL3}. Abbreviating $x:=S_k$ and~$y:=S_{n-k}$, the conditional probability in \eqref{E:11.24nw} is lower bounded by
\begin{equation}
P\bigl(B_t\ge0\colon t\in[t_k,t_{n-k}]\,\big|\,B_{t_k}=x,B_{t_{n-k}}=y\bigr)\,,
\end{equation}
where we used the embedding \eqref{E:embed} of the walk into a path of Brownian motion. In light of Lemma~\ref{lemma-8.8}, we know that 
\begin{equation}
t_\ell-t_k = \bigl(g\log 2+o(1)\bigr)(\ell-k),\quad \ell\ge k,
\end{equation}
with $o(1)\to0$ as~$k\to\infty$ uniformly in~$\ell\ge k$.
Exercise~\ref{ex:7.9} then gives
\begin{equation}
\begin{aligned}
P\Bigl(\,\bigcap_{\ell=k}^{n-k}\{S_\ell\ge0\}\,\Big|\,\sigma(S_k,S_{n-k})\Bigr)
&\gtrsim
\frac{2S_kS_{n-k}}{t_{n-k}-t_k}
\\
&=\frac2{g\log 2}\frac{S_kS_{n-k}}n\bigl(1+o(1)\bigr)
\end{aligned}
\end{equation}
on $\{S_k,S_{n-k}\in[k^{1/6},k^2]\}$ whenever~$k^4\ll n$.

To get a similar upper bound, one writes the Brownian motion on interval $[t_\ell,t_{\ell+1}]$ as a linear curve connecting~$S_\ell$ to~$S_{\ell+1}$ plus a Brownian bridge. Then we observe that the entropic repulsion pushes the walk far away from the positivity constraint so that the Brownian bridges do not affect the resulting probability much. See \cite[Lemma~4.15]{BL3} for details.
\end{proofsect}

The main consequence of the above observations is an asymptotic formula for the numerator (and thus also the denominator) in \eqref{E:11.11}. Indeed, denote
\begin{equation}
\label{E:Xi-in}
\Xiin_k(f):=E\Bigl(\,f\bigl(\phi_k+\tfrac2{\sqrt g}\fraka\bigr)\1_{\{\phi_k+\frac2{\sqrt g}\fraka\,\ge\,0\text{\rm\ in }\Delta^k\}}\1_{\{S_k\in[k^{1/6},\,k^2]\}}S_k\Bigr)
\end{equation}
and
\begin{equation}
\label{E:Xi-out}
\Xiout_{N,k}(t):=E\biggl(\1_{\{h^{D_N}\le m_N(t,\cdot)\text{\rm\ in }D_N\smallsetminus\Delta^{n-k}\}}
\1_{\{S_{n-k}\in[k^{1/6},\,k^2]\}}S_{n-k}\,\bigg|\, h^{D_N}_0=0\,\biggr)\,.
\end{equation}
As a result of the above manipulations, we then get:

\begin{mylemma}[Main asymptotic formula]
\label{lemma-11.12ueu}
For $o(1)\to0$ in the limits as~$N\to\infty$ and~$k\to\infty$, uniformly on compact sets of~$t$ and compact families of~$f\in C_\cc(\overline\R^{\Z^2})$ depending on a given finite number of variables,
\begin{multline}
\label{E:11.28ueu}
\quad
E\Bigl(\,f\bigl(m_N(t,\cdot)-h^{D_N}\bigr)\1_{\{h^{D_N}\le m_N(t,\cdot)\}}\,\Big|\, h^{D_N}_0=0\,\Bigr)
\\
=\frac2n\frac{\Xiin_k(f)\Xiout_{N,k}(t)}{g\log 2}+\frac{o(1)}n\,.
\quad
\end{multline}
\end{mylemma}

\begin{proofsect}{Proof}
This follows by plugging Lemma~\ref{lemma-11.9} in the second expectation in \eqref{E:11.12ueu} and using the independence stated in Lemma~\ref{lemma-11.10ueu}.
\end{proofsect}

In order to control the $N\to\infty$ and $k\to\infty$ limits of the right-hand side of \eqref{E:11.28ueu}, we will also need to observe:

\begin{mylemma}
\label{lemma-11.13ueu}
For each~$t_0>0$ there are~$c_1,c_2\in(0,\infty)$ such that for all~$t\in[-t_0,t_0]$, all~$N\ge1$ and all~$k\le n^{1/6}$ (with~$N$ and~$n$ related as above),
\begin{equation}
\label{E:11.33ueu}
c_1<\Xiout_{N,k}(t)<c_2
\end{equation}
and
\begin{equation}
\label{E:11.34ueu}
c_1<\Xiin_k(1)<c_2.
\end{equation}
\end{mylemma}

\noindent
We postpone the proof of these bounds until the next section. Instead, we note the following consequence thereof:

\begin{mycorollary}
\label{cor-11.11}
Uniformly in~$f$ and~$t$ as stated in Lemma~\ref{lemma-11.12ueu}, 
\begin{equation}
\lim_{N\to\infty}\,
E\Bigl(\,f\bigl(h^{D_N}_0-h^{D_N}\bigr)\,\Big|\, h^{D_N}_0=m_N+t,\,h^{D_N}\le h^{D_N}_0\Bigr)
=\lim_{k\to\infty}\,\frac{\Xiin_k(f)}{\Xiin_k(1)}\,,
\end{equation}
where, in particular, both limits exist.
\end{mycorollary}

\begin{proofsect}{Proof}
The bounds \eqref{E:11.33ueu} allow us to write the right-hand side of \eqref{E:11.28ueu} as
\begin{equation}
\frac{2\Xiout_{N,k}(t)}{g\log 2}\,\frac{\Xiin_k(f)+o(1)}n\,.
\end{equation}
The ratio in \eqref{E:11.11} thus simplifies into the ratio of $\Xiin_k(f)+o(1)$ and $\Xiin_k(1)+o(1)$. This depends on~$N$ only through the $o(1)$ terms which tend to zero as $N\to\infty$ and $k\to\infty$. Using the lower bound in \eqref{E:11.34ueu}, the claim follows.
\end{proofsect}

\section{Full process convergence}
\label{sec-full-convergence}\noindent
Moving towards the proof of of Theorem~\ref{thm-11.1}, thanks to the representation of the pinned DGFF in Exercise~\ref{ex:8.13}, the above derivation implies, even in a somewhat simpler form, also the limit of the probabilities in~\eqref{E:11.2}. The difference is that here the random walk is not constrained to~$S_{n+1}=0$ and also there is no~$t$ to worry about. This affects the asymptotics of the relevant probability as follows: 

\begin{mylemma}
\label{lemma-11.12}
For~$f\in C_\cc(\overline\R^{\Z^2})$ depending on a finite number of coordinates,
\begin{equation}
\label{E:11.35nwt}
E\Bigl(f\bigl(\phi+\tfrac2{\sqrt g}\fraka\bigr)\1_{\{\phi+\frac2{\sqrt g}\fraka\ge0\text{\rm\ in }\Delta^r\}}\Bigr)
=\frac1{\sqrt{\log 2}}\frac{\Xiin_k(f)}{\sqrt r}+\frac{o(1)}{\sqrt r}\,,
\end{equation}
where~$o(1)\to0$ as~$r\to\infty$ followed by~$k\to\infty$.
\end{mylemma}

Similarly to~$1/n$ in \eqref{E:11.28ueu} arising from the asymptotic of the probability that a random-walk bridge of time length~$n$ to stay positive, the asymptotic $1/\sqrt{r}$ stems from the probability that an (unconditioned) random walk stays positive for the first~$r$ steps. Indeed, the reader will readily check:

\begin{myexercise}
\label{ex:11.16}
Let~$\{B_t\colon t\ge0\}$ be the standard Brownian motion with~$P^x$ denoting the law started from~$B_0=x$. Prove that for all $x>0$ and all~$t>0$,
\begin{equation}
\sqrt{\frac2\pi}\,\frac x{\sqrt t}\Bigl(1-\frac{x^2}{2t}\Bigr)\le P^x\bigl(B_s\ge0\colon 0\le s\le t\bigr)\le\sqrt{\frac2\pi}\,\frac x{\sqrt t}\,.
\end{equation}
\end{myexercise}

\noindent
We will not give further details concerning the proof of Lemma~\ref{lemma-11.12} as that would amount to repetitions that the reader may not find illuminating. (The reader may consult \cite[Proposition~5.1]{BL3} for that.) Rather we move on to:

\begin{proofsect}{Proof of Lemma~\ref{lemma-11.13ueu}}
We begin with \eqref{E:11.34ueu}. We start by introducing a variant~$\wt K$ of the control variable~$K$. Let $\wt\theta_k(\ell):=1+[\log(k\vee\ell)]^2$ and define~$\wt K$ to be the least positive natural~$k$ such that Definition~\ref{def-control-var}(1,2) --- with~$\theta_{n,k}(\ell)$ replaced by~$\wt\theta_k(\ell)$ --- as well as
\begin{equation}
\Bigl|\,\min_{x\in\Delta^\ell\smallsetminus\Delta^{\ell-1}}\,\bigl[\chi_\ell(x)+\chi_{\ell+1}(x)+h'_\ell(x)\bigr]+\frac2{\sqrt g}\fraka(x)\Bigr|\le\theta_k(\ell)
\end{equation}
and
\begin{equation}
\label{E:11.42nwut}
\Bigl|\,\max_{x\in\Delta^\ell\smallsetminus\Delta^{\ell-1}}\,\bigl[\chi_\ell(x)+\chi_{\ell+1}(x)+h'_\ell(x)\bigr]-\frac2{\sqrt g}\fraka(x)\Bigr|\le\theta_k(\ell)
\end{equation}
hold true for all~$\ell\ge k$. (That~$\wt K<\infty$ a.s.\ follows by our earlier estimates and the Borel-Cantelli lemma. The condition \eqref{E:11.42nwut} is introduced for later convenience.) 

Recall that Exercise~\ref{ex:8.13} expresses the pinned DGFF using the objects in the concentric decomposition. Hence we get
\begin{equation}
\bigcap_{k=1}^r\bigl\{S_k\ge C\wt\theta_{\wt K}(k)\bigr\}\subseteq 
\Bigl\{\phi+\frac2{\sqrt g}\fraka\ge0\text{\rm\ in }\Delta^r\Bigr\}
\subseteq \bigcap_{k=1}^r\bigl\{S_k\ge-C\wt\theta_{\wt K}(k)\bigr\},
\end{equation}
for some sufficiently large absolute constant~$C>0$. Ballot problem/entropy repulsion arguments invoked earlier then show that the probability of both sides decays proportionally to that of the event $\bigcap_{k=1}^r\{S_\ell\ge-1\}$. Lemma~\ref{lemma-11.12} and Exercise~\ref{ex:11.16} then give \eqref{E:11.34ueu}. 

Plugging \eqref{E:11.34ueu} to \eqref{E:11.28ueu}, the inclusions in Lemma~\ref{lemma-10.5} along with the bounds in Lemmas~\ref{lemma-11.7}--\ref{lemma-11.6} then imply \eqref{E:11.33ueu} as well.
\end{proofsect}

We are finally ready to give:

\begin{proofsect}{Proof of Theorem~\ref{thm-11.1}}
From Lemmas~\ref{lemma-11.12} and~\ref{lemma-11.13ueu} we have
\begin{equation}
\label{E:11.31}
\lim_{r\to\infty} E\Bigl(f\bigl(\phi+\tfrac2{\sqrt g}\fraka\bigr)\,\Big|\,\phi_x+\frac2{\sqrt g}\fraka(x)\ge0\colon |x|\le r\Bigr)=\lim_{k\to\infty}\,\frac{\Xiin_k(f)}{\Xiin_k(1)}.
\end{equation}
Jointly with Corollary~\ref{cor-11.11}, this proves equality of the limits in the statement. 
To see that~$\nu$ concentrates on~$\R^{\Z^2}$ (and, in fact, on $[0,\infty)^{\Z^2}$) we observe that, since $\ell\mapsto\wt\theta_k(\ell)$ from Lemma~\ref{lemma-11.13ueu} grows only polylogarithmically while~$\fraka$ is order~$k$ on~$\Delta^k\smallsetminus\Delta^{k-1}$, once~$k$ is sufficiently large we get
\begin{equation}
\label{E:11.39nw}
\Bigl\{\phi+\frac2{\sqrt g}\fraka\not\le k^2\text{\rm\ on }\Delta^k\Bigr\}
\subseteq \{\wt K> k\}.
\end{equation}
(It is here where we make use of \eqref{E:11.42nwut}.)
By an analogue of Lemma~\ref{lemma-11.7} for random variable~$\wt K$, we have $\nu(\wt K>k)\le c_1\texte^{-c_2(\log k)^2}$. The Borel-Cantelli lemma then implies $\nu(\R^{\Z^2})=1$, as desired.
\end{proofsect}

The above techniques also permit us to settle the asymptotic \eqref{E:11.4}:

\begin{proofsect}{Proof of Theorem~\ref{thm-11.3}}
Observe that after multiplying both sides of \eqref{E:11.35nwt} by~$\sqrt r$, the left-hand side is independent of~$k$ while the right-hand side depends on~$r$ only through the~$o(1)$-term. 
In light of \eqref{E:11.34ueu}, ~$r\to\infty$ and~$k\to\infty$ can be taken independently (albeit in this order) to get that
\begin{equation}
\label{E:11.40nwt}
\Xiin_\infty(1):=\lim_{k\to\infty} \Xiin_k(1)
\end{equation}
exists, is positive and finite. Since $r\log 2$ is, to the leading order, the logarithm of the diameter of~$\Delta^r$, the claim follows from \eqref{E:11.35nwt} with~$c_\star:=\Xiin_\infty(1)$.
\end{proofsect}

We are now finally ready to move to the proof of Theorem~\ref{thm-extremal-vals} dealing with the weak convergence of the full, three coordinate process~$\eta_{N,r}^D$ defined in \eqref{E:9.6}. In light of Theorem~\ref{thm-11.1}, all we need to do is to come up with a suitable localization method that turns a large local maximum (in a large domain) to the actual maximum (in a smaller domain). We will follow a different line of reasoning than~\cite{BL3} as the arguments there seem to contain flaws whose removal would take us through lengthy calculations that we prefer to avoid. 

\begin{proofsect}{Proof of Theorem~\ref{thm-extremal-vals}}
Every~$f\in C_\cc(D\times\R\times\overline\R^{\Z^2})$ is uniformly close to a compactly supported function that depends only on a finite number of ``cluster'' coordinates. So let us assume that~$x,h,\phi\mapsto f(x,h,\phi)$ is continuous and depends only on~$\{\phi_y\colon y\in\Lambda_r(0)\}$ for some~$r\ge1$ and vanishes when~$|h|,\max_{x\in\Lambda_r}|\phi_y|\ge \lambda$ for some~$\lambda>0$ or if~$\dist(x,D^\cc)<\delta$ for some~$\delta>0$. 

The proof hinges on approximation of~$\eta^D_N$ by three auxiliary processes. Given an integer~$K\ge1$, let $\{S^i\colon i=1,\dots,m\}$ be the set of all the squares of the form $z/K+(0,1/K)^2$ with $z\in\Z^2$ that fit entirely into~$D^\delta$.  For each~$i=1,\dots,m$, let~$S^i_N$ be the lattice box of side-length (roughly) $N/K$ contained in $NS^i$ so that every pair of neighboring squares $S^i_N$ and $S^j_N$ keep a ``line of sites'' in-between. Given $\delta>0$ small, let $S^i_{N,\delta}:=\{x\in S^i_N\colon \dist_\infty(x,S^i_N)>\delta N\}$. For~$x\in D_N$ such that $x\in S^i_{N,\delta}$ for some $i=1,\dots,m$, set
\begin{equation}
\Theta_{N,K,\delta}(x):=\Bigl\{h^{D_N}_x=\max_{y\in S^i_N}h^{D_N}_y\Bigr\}
\end{equation}
and let $\Theta_{N,K,\delta}(x):=\emptyset$ for $x\in D_N$ where no such~$i$ exists. Setting
\begin{equation}
\eta^{D,K,\delta}_N:=\sum_{x\in D_N}1_{\Theta_{N,K,\delta}(x)}\,\delta_{x/N}\otimes\delta_{h^{D_N}_x-m_N}\otimes\delta_{\{h^{D_N}_x-h^{D_N}_{x+z}\colon z\in\Z^2\}},
\end{equation}
Lemma~\ref{lemma-10.4} and Theorem~\ref{thm-DZ2} show that, for any~$f$ as above,
\begin{equation}
\lim_{\delta\downarrow0}\,\limsup_{K\to\infty}\,\limsup_{N\to\infty}\,
\Bigl|E(\texte^{-\langle\eta_{N,r_N}^D,f\rangle})-E(\texte^{-\langle\eta^{D,K,\delta}_N,f\rangle})\Bigr|=0.
\end{equation}
Next we consider the Gibbs-Markov decomposition~$h^{D_N}=h^{\wt D_N}+\varphi^{D_N,\wt D_N}$, where $\wt D_N:=\bigcup_{i=1}^m S^i_N$.
Denoting, for~$x\in \wt D_N$ and~$i=1,\dots,m$ such that~$x\in S^i_{N,\delta}$,
\begin{equation}
\wt\Theta_{N,K,\delta}(x):=\Bigl\{h^{\wt D_N}_x=\max_{y\in S^i_N}h^{\wt D_N}_y\Bigr\},
\end{equation}
we then define~$\wt\eta^{D,K,\delta}_N$ by the same formula as~$\eta^{D,K,\delta}_N$ but with~$\Theta_{N,K,\delta}(x)$ replaced by~$\wt\Theta_{N,K,\delta}(x)$ and the sum running only over~$x\in\wt D_N$. The next point to observe is that, for~$K$ sufficiently large, the maximum of~$h^{D_N}$ in~$S^i_N$ coincides (with high probability) with that of~$h^{\wt D_N}$. (We used a weaker version of this already in Lemma~\ref{lemma-10.17}.) This underpins the proof of:

\begin{mylemma}
\label{lemma-same-max}
For any~$f$ as above,
\begin{equation}
\lim_{\delta\downarrow0}\,\limsup_{K\to\infty}\,\limsup_{N\to\infty}\,
\Bigl|E(\texte^{-\langle\eta^{D,K,\delta}_N,f\rangle})-E(\texte^{-\langle\wt\eta^{D,K,\delta}_N,f\rangle})\Bigr|=0.
\end{equation}
\end{mylemma}

\noindent
Postponing the proof to Lecture~\ref{lec-12}, we now define a third auxiliary process $\wh\eta^{D,K,\delta}_N$ by replacing $h^{D_N}$ in the cluster variables of $\wt\eta^{D,K,\delta}_N$ by~$h^{\wt D_N}$,
\begin{equation}
\wh\eta^{D,K,\delta}_N:=\sum_{x\in \wt D_N}1_{\wt\Theta_{N,K,\delta}(x)}\,\delta_{x/N}\otimes\delta_{h^{D_N}_x-m_N}\otimes\delta_{\{h^{\wt D_N}_x-h^{\wt D_N}_{x+z}\colon z\in\Z^2\}}.
\end{equation}
By the uniform continuity of~$f$ and the fact that~$\varphi^{D_N,\wt D_N}$ converges locally-uni\-form\-ly to the (smooth) continuum binding field~$\Phi^{D,\wt D}$ (see Lemma~\ref{lemma-4.4}) the reader will readily verify:

\begin{myexercise}
Show that for each~$K\ge1$, each~$\delta>0$ and each~$f$ as above,
\begin{equation}
\lim_{N\to\infty}\,
\Bigl|E(\texte^{-\langle\wh\eta^{D,K,\delta}_N,f\rangle})-E(\texte^{-\langle\wt\eta^{D,K,\delta}_N,f\rangle})\Bigr|=0.
\end{equation}
\end{myexercise}

Now comes the key calculation of the proof. Let~$X_i$, for $i=1,\dots,m$, denote the (a.s.-unique) maximizer of~$h^{\wt D_N}$ in~$S^i_N$. For the given~$f$ as above,~$x\in S^i$ for some~$i=1,\dots,m$ and any~$t\in\R$, abbreviate $x_N:=\lfloor xN\rfloor$ and let
\begin{multline}
f_{N,K}(x,t)
\\
:=-\log E\Bigl(\texte^{-f(x,t,h^{\wt D_N}_{x_N}-h^{\wt D_N}_{x_N+z}\colon z\in\Z^2)}\,\Big|\, X_i=x_N,\,h_{x_N} = m_N+t\Bigr).
\end{multline}
Thanks to the independence of the field $h^{\wt D_N}$ over the boxes~$S^i_N$, $i=1,\dots,m$, by conditioning on the binding field~$\varphi^{D_N,\wt D_N}$ and using its independence of~$h^{\wt D_N}$ (and thus also of the~$X_i$'s), and then also conditioning on the~$X_i$'s and the values~$h^{\wt D_N}_{X_i}$, we get
\begin{equation}
\label{E:11.49nwt}
\begin{aligned}
&E(\texte^{-\langle\wh\eta^{D,K,\delta}_N,f\rangle})
=E\biggl(\,\prod_{i=1}^m \texte^{-f(X_i/N,\,h^{\wt D_N}_{X_i}-m_N+\varphi^{D_N,\wt D_N}_{X_i},\{h^{\wt D_N}_{X_i}-h^{\wt D_N}_{X_i+z}\colon z\in\Z^2\})}\biggr)
\\
&=
E\biggl(\,\prod_{i=1}^m \texte^{-f_{N,K}(X_i/N,\,h^{\wt D_N}_{X_i}-m_N+\varphi^{D_N,\wt D_N}_{X_i})}\biggr)
=E(\texte^{-\langle\wh\eta^{D,K,\delta}_N,f_{N,K}\rangle}).
\end{aligned}
\end{equation}
Since~$X_i$ marks the location of the absolute maximum of~$h^{\wt D_N}$ in~$S^i_{N}$, recalling the notation~$\nu$ for the cluster law from \eqref{E:11.2}, Theorem~\ref{thm-11.1} yields
\begin{equation}
f_{N,K}(x,t)\,\underset{N\to\infty}\longrightarrow\, f_\nu(x,t):=-\log E_\nu(\texte^{-f(x,t,\phi)})
\end{equation}
uniformly in~$t$ and~$x\in\bigcup_{i=1}^m 
\wt S^i_\delta$, where~$\wt S^i_\delta$ is the shift of~$(0,(1-\delta)/K)^2$ centered the same point as~$S^i$. Using this in \eqref{E:11.49nwt}, the series of approximations invoked earlier  shows
\begin{equation}
E\bigl(\texte^{-\langle\eta^D_{N,r_N},\,f\rangle}\bigr)
= E\bigl(\texte^{-\langle\eta^D_{N,r_N},\,f_\nu\rangle}\bigr)+o(1)
\end{equation}
with~$o(1)\to0$ as~$N\to\infty$. 
As~$f_\nu\in C_\cc(D\times\R)$, the convergence of the two-coordinate process proved in Section~\ref{sec:10.2} yields
\begin{multline}
\label{E:11.58ueu}
\quad
E\bigl(\texte^{-\langle\eta^D_{N,r_N},\,f_\nu\rangle}\bigr)
\\\underset{N\to\infty}\longrightarrow\,E\biggl(\exp\Bigl\{-\int_{D\times\R}Z^D(\textd x)\otimes\texte^{-\alpha h}\textd h\,\bigl(1-\texte^{-f_\nu(x,h)}\bigr) \Bigr\}\biggr)\,.
\quad
\end{multline}
To conclude, it remains to observe that
\begin{multline}
\quad\int_{D\times\R}Z^D(\textd x)\otimes\texte^{-\alpha h}\textd h\,\bigl(1-\texte^{-f_\nu(x,h)}\bigr) 
\\
= \int_{D\times\R\times\R^{\Z^2}}Z^D(\textd x)\otimes\texte^{-\alpha h}\textd h\otimes\nu(\textd\phi)\bigl(1-\texte^{-f(x,h,\phi)}\bigr)
\quad
\end{multline}
turns \eqref{E:11.58ueu} into the Laplace transform of $\PPP(Z^D(\textd x)\otimes\texte^{-\alpha h}\textd h\otimes\nu(\textd\phi))$.
\end{proofsect}

\section{Some corollaries}
\label{sec-11.4}\noindent
Having established the limit of the structured point measure, we proceed to state a number of corollaries of interest. We begin with the limit of the ``ordinary'' extreme value process \eqref{E:eta-ND}:

\begin{mycorollary}[Cluster process]
\label{cor-cluster-process}
For~$Z^D$ and~$\nu$ as in Theorem~\ref{thm-extremal-vals},
\begin{equation}
\label{E:11.47}
\sum_{x\in D_N}\delta_{x/N}\otimes\delta_{h^{D_N}_x-m_N}\,\,\,\underset{N\to\infty}\lawarrow\,\,\,\sum_{i\in\N}\sum_{z\in\Z^2}\delta_{(x_i,\,h_i-\phi_z^{(i)})},
\end{equation}
where the right-hand side is defined using the following independent objects:
\settowidth{\leftmargini}{(1111)}
\begin{enumerate}
\item[(1)] $\{(x_i,h_i)\colon i\in\N\}:=$ sample from $\PPP(Z^D(\textd x)\otimes\texte^{-\alpha h}\textd h)$, 
\item[(2)] $\{\phi^{(i)}\colon i\in\N\}:=$ i.i.d.\ samples from~$\nu$.
\end{enumerate}
The measure on the right is locally finite on $D\times\R$ a.s. 
\end{mycorollary}

We relegate the proof to:

\begin{myexercise}
\label{ex:11.20}
Derive \eqref{E:11.47} from the convergence statement in Theorem~\ref{thm-extremal-vals} and the tightness bounds in Theorems~\ref{thm-DZ1}--\ref{thm-DZ2}.
\end{myexercise}

The limit object on the right of \eqref{E:11.47} takes the form of a \emph{cluster process}. This term generally refers to a collection of random points obtained by taking a sample of a Poisson point process and then associating with each point thereof an independent cluster of (possibly heavily correlated) points. See again Fig.~\ref{fig-cluster-process}. We note that a cluster process naturally appears in the limit description of the extreme-order statistics of the Branching Brownian Motion~\cite{ABK1,ABK2,ABK3}.

Another observation that is derived from the above limit law concerns the Gibbs measure on~$D_N$ associated with the DGFF on~$D_N$ as follows:
\begin{equation}
\label{E:11.48}
\mu_{\beta,N}^D\bigl(\{x\}\bigr):=\frac1{\mathfrak Z_N(\beta)}\,\texte^{\beta h_x^{D_N}}
\quad\text{where}\quad\mathfrak Z_N(\beta):=\sum_{x\in D_N}\texte^{\beta h_x^{D_N}}.
\end{equation}
In order to study the scaling limit of this object, we associate the value~$\mu_{\beta,N}^D(\{x\})$ with a point mass at~$x/N$. From the convergence of the suitably-normalized measure $\sum_{x\in D_N}\texte^{\beta h_x^{D_N}}\delta_{x/N}$
to the Liouville Quantum Gravity for $\beta<\beta_\cc:=\alpha$ it is known (see, e.g., Rhodes and Vargas~\cite[Theorem 5.12]{RV-review}) that
\begin{equation}
\label{E:11.49}
\sum_{z\in D_N}\mu_{\beta,N}^D\bigl(\{z\}\bigr)\delta_{z/N}(\textd x)
\,\,\,\underset{N\to\infty}\lawarrow\,\,\,\frac{Z^D_\lambda(\textd x)}{Z^D_\lambda(D)}\,,
\end{equation}
where~$\lambda:=\beta/\beta_\cc$ and where $Z^D_\lambda$ is the measure we saw in the discussion of the intermediate level sets (for~$\lambda<1$). The result extends (see \cite[Therem~5.13]{RV-review}, although the proof details seem scarce) to the case $\beta=\beta_\cc$, where thanks to Theorem~\ref{thm-BLnew} we get~$\wh Z^D(\textd x)$ from \eqref{E:10.18ueu} instead.

The existence of the limit in the supercritical cases $\beta>\beta_\cc$ has been open for quite a while (and was subject to specific conjectures; e.g.,~\cite[Conjecture~11]{DRSV1}). Madaule, Rhodes and Vargas~\cite{MRV} first proved it for star-scale invariant fields as well as certain specific cutoffs of the CGFF. For the DGFF considered here, Arguin and Zindy~\cite{AZ} proved convergence of the \emph{overlaps} of~$\mu^D_{\beta,N}$ to those of Poisson-Dirichlet distribution $\text{\rm PD}(s)$ with parameter $s:=\beta_\cc/\beta$. This only identified the law of the \emph{sizes} of the atoms in the limit measure; the full convergence including the spatial distribution was settled in~\cite{BL3}:

\begin{mycorollary}[Poisson-Dirichlet limit for the Gibbs measure]
\label{cor-11.16}
For all $\beta>\beta_\cc:=\alpha$ we then have
\begin{equation}
\label{E:11.51}
\sum_{z\in D_N}\mu_{\beta,N}^D\bigl(\{z\}\bigr)\delta_{z/N}(\textd x)
\,\,\,\underset{N\to\infty}\lawarrow\,\,\,
\sum_{i\in \N}p_i\,\delta_{X_i},
\end{equation}
where 
\begin{enumerate}
\item[(1)] $\{X_i\}$ are (conditionally on~$Z^D$) i.i.d.\ with law~$\widehat Z^D$, while
\item[(2)] $\{p_i\}_{i\in\N}$ is independent of $Z^D$ and $\{X_i\}$ with $\{p_i\}_{i\in\N} \,\laweq\,\text{\rm PD}(\beta_\cc/\beta)$.
\end{enumerate}
\end{mycorollary}

We remark that $\text{\rm PD}(s)$ is a law on non-increasing sequences of non-negative numbers with unit total sum obtained as follows: Take a sample $\{q_i\}_{i\in\N}$ from the Poisson process on $(0,\infty)$ with intensity $t^{-1-s}\textd t$, order the points decreasingly and normalize them by their total sum (which is finite a.s.\ when~$s<1$). 

Corollary~\ref{cor-11.16} follows from the description of the Liouville Quantum Gravity measure for~$\beta>\beta_\cc$ that we will state next. For~$s>0$ and~$Q$ a Borel probability measure on~$\C$, let $\{q_i\}_{i\in\N}$ be a sample from the Poisson process on $(0,\infty)$ with intensity $t^{-1-s}\textd t$ and let $\{X_i\}_{i\in\N}$ be independent i.i.d.\ samples from~$Q$. Use these to define the random measure
\begin{equation}
\label{E:2.17nwt}
\Sigma_{s,Q}(\textd x):=\sum_{i\in\N}q_i\,\delta_{X_i}\,.
\end{equation}
We then have:

\begin{mytheorem}[Liouville measure in the glassy phase]
\label{thm-11.7}
Let~$Z^D$ and~$\nu$ be as in Theorem~\ref{thm-extremal-vals}. For each $\beta>\beta_\cc:=\alpha$ there is~$c(\beta)\in(0,\infty)$ such that
\begin{equation}
\label{E:11.53}
\sum_{z\in D_N}\texte^{\beta (h_z-m_N)}\delta_{z/N}(\textd x)
\,\,\,\underset{N\to\infty}\lawarrow\,\,\,
c(\beta) \,Z^D(D)^{\beta/\beta_\cc}\,\,\,\Sigma_{\beta_\cc/\beta,\,\widehat Z^D}(\textd x),
\end{equation}
where~$\Sigma_{\beta_\cc/\beta,\,\widehat Z^D}$ is defined conditionally  on~$Z^D$. Moreover,
\begin{equation}
\label{E:11.56}
c(\beta)=\beta^{-\beta/\beta_\cc}\bigl[E_\nu(Y^\beta(\phi)^{\beta_\cc/\beta})\bigr]^{\beta/\beta_\cc}\quad\text{\rm with}\quad
Y^\beta(\phi):=\sum_{x\in\Z^2}\texte^{-\beta \phi_x}.
\end{equation}
In particular, $E_\nu(Y^\beta(\phi)^{\beta_\cc/\beta})<\infty$ for each~$\beta>\beta_\cc$. 
\end{mytheorem}

Note that the limit laws in \eqref{E:11.51} and \eqref{E:11.53} are \emph{purely atomic}, in contrast to the limits of the subcritical measures \eqref{E:11.49} which, albeit a.s.\ singular with respect to the Lebesgue measure, are non-atomic a.s.

\begin{proofsect}{Proof of Theorem~\ref{thm-11.7} (main computation)}
We start by noting that the Laplace transform of the above measure $\Sigma_{s,Q}$ is explicitly computable:

\begin{myexercise}
Show that for any measurable~$f\colon \C\to[0,\infty)$,
\begin{equation}
\label{E:11.62nwt}
E\bigl(\texte^{-\langle \Sigma_{s,Q},f\rangle}\bigr)=\exp\biggl\{-\int_{\C\times(0,\infty)}Q(\textd x)\otimes t^{-1-s}\textd t\,\,(1-\texte^{-t f(x)})\biggr\}.
\end{equation}
\end{myexercise}

\noindent
Pick a continuous~$f\colon \C\to[0,\infty)$ with support in~$D$ and (abusing our earlier notations) write~$M_N$ to denote the measure on the left of \eqref{E:11.53}. Through suitable truncations, the limit proved in Theorem~\ref{thm-extremal-vals} shows that $E(\texte^{-\langle M_N,f\rangle})$ tends to
\begin{equation}
\label{E:11.66ueu}
E\biggl(\exp\Bigl\{-\int Z^D(\textd x)\otimes \texte^{-\alpha h}\textd h\otimes\nu(\textd\phi)(1-\texte^{-g(x,h,\phi)})\Bigr\}\biggr)\,,
\end{equation}
where $g(x,h,\phi):=f(x)\,\texte^{\beta h}\,Y^\beta(\phi)$. The change of variables $t:=\texte^{\beta h}$ then gives
\begin{equation}
\begin{aligned}
\int\textd h\,\texte^{-\alpha h}(1-\texte^{-g(x,h,\phi)})
&=\int_0^\infty\textd t\frac{1}{\beta t}\, t^{-\alpha/\beta}\bigl(1-\texte^{-f(x)Y^\beta(\phi)t}\bigr)
\\
&=\frac1\beta \bigl(Y^\beta(\phi)\bigr)^{\alpha/\beta}\int_0^\infty\textd t\,t^{-1-\alpha/\beta}(1-\texte^{-tf(x)}).
\end{aligned}
\end{equation}
The integral with respect to~$\nu$ affects only the term $(Y^\beta(\phi))^{\alpha/\beta}$ in the front; scaling~$t$ by $c(\beta) Z^D(D)^{\beta/\alpha}$ then absorbs all prefactors and identifies the integral on the right of \eqref{E:11.66ueu} with that in \eqref{E:11.62nwt} for $Q:=\wh Z^D$ and~$f(x)$ replaced by $c(\beta) Z^D(D)^{\beta/\alpha}f(x)$. This now readily gives the claim \eqref{E:11.53}. (The finiteness of the expectation of $Y^\beta(\phi)^{\beta_\cc/\beta}$ requires a separate tightness argument.)
\end{proofsect}

The reader might wonder at this point how it is possible that the rather complicated (and correlated) structure of the cluster law~$\nu$ does not at all appear in the limit measure on the right of \eqref{E:11.53} --- that is, not beyond the expectation of $Y^\beta(\phi)^{\beta_\cc/\beta}$ in~$c(\beta)$. This may, more or less, be traced to the following property of the Gumbel~law:

\begin{myexercise}[Derandomizing i.i.d.\ shifts of Gumbel PPP]
\label{ex:11.18}
Consider a sample $\{h_i\colon i\in\N\}$ from $\PPP(\texte^{-\alpha h}\textd h)$ and let~$\{T_i\colon i\in\N\}$ be independent, i.i.d.\ random variables. Prove that
\begin{equation}
\sum_{i\in\N}\delta_{h_i+T_i}\,\,\laweq\,\,\sum_{i\in\N}\delta_{h_i+\alpha^{-1}\log c}
\end{equation}
whenever $c:=E\texte^{\alpha T_1}<\infty$.
\end{myexercise}

Our final corollary concerns the behavior of the function
\begin{equation}
 G_{N,\beta}(t):=E\biggl(\,\exp\Bigl\{-\texte^{-\beta t}\sum_{x\in D_N}\texte^{\beta h_x^{D_N}}\Bigr\}\biggr),
\end{equation}
which, we observe, is a re-parametrization of the Laplace transform of the normalizing constant $\mathfrak Z_N(\beta)$ from \eqref{E:11.48}. In their work on the Branching Brownian Motion, Derrida and Spohn~\cite{Derrida-Spohn} and later Fyodorov and Bouchaud~\cite{Fyodorov-Bouchaud} observed that, in a suitable limit, an analogous quantity ceases to depend on~$\beta$ once~$\beta$ crosses a critical threshold. They referred to this as \emph{freezing}. Our control above is sufficient to yield the same phenomenon for the quantity arising from DGFF:

\begin{mycorollary}[Freezing]
For all $\beta>\beta_\cc:=\alpha$ there is $\tilde c(\beta)\in\R$ such that
\begin{equation}
\label{E:2.26}
G_{N,\beta}\bigl(\,t+m_N+\tilde c(\beta)\bigr)\,\underset{N\to\infty}\longrightarrow\,E\bigl(\texte^{-Z^D(D)\,\texte^{-\alpha t}}\bigr).
\end{equation}
\end{mycorollary}

\begin{proofsect}{Proof}
Noting that $\texte^{-\beta m_N}\mathfrak Z_N(\beta)$ is the total mass of the measure on the left of \eqref{E:11.53}, from Theorem~\ref{thm-11.7} we get (using the notation of \eqref{E:2.17nwt} and \eqref{E:11.56})
\begin{equation}
\label{E:11.68nwt}
\texte^{-\beta m_N}\mathfrak Z_N(\beta)\,\,\,\underset{N\to\infty}\lawarrow\,\,\, c(\beta) Z^D(D)^{\beta/\alpha}\sum_{i\in\N}q_i\,.
\end{equation}
The Poisson nature of the $\{q_i\}_{i\in\N}$ shows, for any~$\lambda>0$,
\begin{equation}
\label{E:11.69nwt}
\begin{aligned}
E\bigl(\texte^{-\lambda\sum_{i\in\N}q_i}\bigr)
&=\exp\Bigl\{-\int_0^\infty\textd t\,\, t^{-1-\beta/\alpha}(1-\texte^{-\lambda t})\Bigr\}
\\
&=\exp\Bigl\{-\lambda^{\alpha/\beta}\int_0^\infty\textd t\,\, t^{-1-\beta/\alpha}(1-\texte^{-t})\Bigr\}.
\end{aligned}
\end{equation}
Thanks to the independence of $\{q_i\}_{i\in\N}$ and~$Z^D$, the Laplace transform of the random variable on the right of \eqref{E:11.68nwt} at parameter~$\texte^{-\beta t}$ is equal to the right-hand side of \eqref{E:2.26}, modulo a shift in~$t$ by a constant that depends on~$c(\beta)$ and the integral on the second line in \eqref{E:11.69nwt}.
\end{proofsect}

We refer to~\cite{BL3} for further consequences of the above limit theorem and additional details. The (apparently quite deep) connection between freezing and Gumbel laws has recently been further explored by Subag and Zeitouni~\cite{Subag-Zeitouni}.


\chapter{Limit theory for DGFF maximum}
\label{lec-12}\noindent
In this final lecture on the extremal process associated with the DGFF we apply the concentric decomposition developed in Lectures~\ref{lec-8} and~\ref{lec-11} to give the proofs of various technical results scattered throughout Lectures~\ref{lec-9}--\ref{lec-11}. We start by Theorems~\ref{thm-DZ1}-\ref{thm-DZ2} and Lemma~\ref{lemma-10.4} that we used to control the spatial tightness of the extremal level sets. Then we establish the convergence of the centered maximum from Theorem~\ref{thm-BDZ} by way of proving Theorem~\ref{prop-10.7} and extracting from this the uniqueness (in law) of the limiting $Z^D$-measures from Theorem~\ref{prop-subseq}. At the very end we state (without proof) a local limit result for both the position and the value of the absolute maximum.

\section{Spatial tightness of extremal level sets}
\noindent
Here we give the proofs of Theorems~\ref{thm-DZ1}-\ref{thm-DZ2} and Lemma~\ref{lemma-10.4} that we relied heavily on in the previous lectures.
Recall the notation~$\Gamma_N^D(t)$ from \eqref{E:9.1} for the extremal level set ``$t$ units below~$m_N$'' and the definition $D_N^\delta:=\{x\in D_N\colon\dist_\infty(x,D_N^\cc)>\delta N\}$. We start by:

\begin{mylemma}
\label{lemma-4.9new}
For each $D\in\mathfrak D$,~$\delta>0$ and~$c'>0$ there is $c>0$ such that for all $N\ge1$, all $x\in D_N^\delta$ and all $t$ and~$r$ with $|t|,r\in[0,c'(\log N)^{1/5}]$,
\begin{equation}
\label{E:4.24z}
P\Bigl(\,h^{D_N}\le m_N+r
\,\Big|\,h^{D_N}_x=m_N+t\Bigr)\le\frac{c}{\log N}\bigl(1+r+|t|\bigr)^2.
\end{equation}
\end{mylemma}

\begin{proofsect}{Proof}
Assume without loss of generality that~$0\in D$ and let us first address the case~$x=0$. Lemma~\ref{lemma-8.3} shifts the conditioning to $h^{D_N}_0=0$ at the cost of reducing the field by $(m_N+t)\frakg^{D_N}$. Invoking the concentric decomposition and Lemma~\ref{lemma-10.5}, the resulting event can be covered as
\begin{multline}
\bigl\{h^{D_N}\le m_N(1-\frakg^{D_N})+r-t\frakg^{D_N}\bigr\}
\subseteq\{K=\lfloor\ffrac n2\rfloor+1\}
\\
\cup \bigcap_{k=1}^n\Bigl(\bigl\{S_k-S_{n+1}\ge -R_K(k)-r-|t|\bigr\}\cap\bigl\{K\le\lfloor\ffrac n2\rfloor\bigr\}\Bigr)\,,
\end{multline}
where we used $0\le\frakg^{\Delta^n}\le1$. Summing the bound in Lemma~\ref{lemma-11.7} we then get \eqref{E:4.24z}. (The restriction to~$t\le n^{1/5}$ can be replaced by~$t\le c'n^{1/5}$ at the cost of changing the constants.) Noting that constant~$c$ depends only on~$k_1$ from \eqref{E:8.34}, shifting~$D_N$ around  extends the bound to all~$x\in D_N^\delta$. 
\end{proofsect}

With Lemma~\ref{lemma-4.9new} at our disposal, we can settle the tightness of the cardinality of the extremal level set:

\begin{proofsect}{Proof of Theorem~\ref{thm-DZ1} (upper bound)}
Building on the idea underlying Exercises~\ref{ex:3.3}--\ref{ex:3.4}, we first ask the reader to solve:

\begin{myexercise}
Let~$\wt D\subseteq D$. Then, assuming~$\wt D_N\subseteq D_N$, for any~$n\in\N$,
\begin{equation}
P\Bigl(\bigl|\Gamma^{\wt D}_N(t)\bigr|\ge 2n\Bigr)\le 2P\Bigl(\bigl|\Gamma^D_N(t)\cap\wt D_N\bigr|\ge n\Bigr).
\end{equation}
\end{myexercise}

\noindent
In light of this, it suffices to bound the size of~$\Gamma^D_N(t)\cap D_N^\delta$. Here we note that, by the Markov inequality
\begin{multline}
\quad
P\Bigl(|\Gamma^D_N(t)\cap D_N^\delta|\ge \texte^{Ct}\Bigr)\le P\Bigl(\,\max_{x\in D_N}h^{D_N}_x>m_N+t\Bigr)
\\+\texte^{-Ct}\sum_{x\in D_N^\delta}P\Bigl(h^{D_N}\le m_N+t,\,h^{D_N}_x\ge m_N-t\Bigr).
\quad
\end{multline}
The first probability is bounded by~$\texte^{-\tilde a t}$ using the exponential tightness of absolute maximum; cf Lemma~\ref{lemma-8.2}.
The probability under the sum is in turn bounded by conditioning on $h^{D_N}_x$ and applying Lemma~\ref{lemma-4.9new} along with the standard Gaussian estimate (for the probability of $h^{D_N}_x\ge m_N-t$). This bounds the probability by a quantity of order $(1+t)^2t\texte^{\alpha t}/N^2$. Since $|D_N^\delta|=O(N^2)$, for~$C>\alpha$, we get the upper bound in~\eqref{E:9.2nwt}.
\end{proofsect}

We leave the lower bound to:

\begin{myexercise}
Use the geometric setting, and the ideas underlying the proof of Lemma~\ref{lemma-8.2} to prove the lower bound on~$|\Gamma_N^D(t)|$ in \eqref{E:9.2nwt} for any $c<\alpha/2$.
\end{myexercise}

\noindent
We remark that the bound proposed in this exercise is far from optimal. Indeed, once Theorem~\ref{thm-extremal-vals} has been established in full (for which we need only the upper bound in Theorem~\ref{thm-DZ1}), we get $|\Gamma_N(t)|=\texte^{(\alpha+o(1))t}$. See Conjecture~\ref{cor-12.13ua} for even a more precise claim.

\smallskip
Next we will give:

\begin{proofsect}{Proof of Lemma~\ref{lemma-10.4}}
First we note that, thanks to Exercise~\ref{ex:3.4} we may assume that~$A\subset D_N^\delta$ for some~$\delta>0$ independent of~$N$. Lemma~\ref{lemma-4.9new}, the union bound, the definition of~$m_N$ and the exponential tightness of the absolute maximum from Lemma~\ref{lemma-8.2} then shows, for all $r,|t|\in[0,(\log N)^{1/5}]$,
\begin{equation}
\begin{aligned}
P\bigl(\,\max_{x\in A}&\,h_x^{D_N}\ge m_N+t\bigr)
\\
&\le P\bigl(\,\max_{x\in D_N}h_x^{D_N}\ge m_N+r\bigr)
\\&\qquad\quad+P\bigl(\,\max_{x\in A}h_x^{D_N}\ge m_N+t,\,\max_{x\in D_N}h_x^{D_N}\le m_N+r\bigr)
\\
&\le\texte^{-\tilde ar}+c\frac{|A|}{N^2}\texte^{-\alpha t}\bigl(1+r+|t|\bigr)^2.
\end{aligned}
\end{equation}
Setting $r:=5\tilde a^{-1}\log(1+N^2/|A|)+2\alpha\tilde a^{-1}|t|$ then proves the claim for all~$A$ such that $|A|\ge N^2\texte^{-\tilde a(\log N)^{1/5}}$ (as that ensures $r\le c'(\log N)^{1/5}$). For the complementary set of~$A$'s we apply the standard Gaussian tail estimate along with the uniform bound on the Green function and a union bound to get
\begin{equation}
P\bigl(\,\max_{x\in A}h_x^{D_N}\ge m_N+t\bigr)\le c\texte^{-\alpha t}\frac{|A|}{N^2}\log N.
\end{equation}
Then we note that $\log N\le c''[\log(1+N^2/|A|)]^5$ for $c''>0$ sufficiently large.
\end{proofsect}

Our next item of business is the proof of Theorem~\ref{thm-DZ2} dealing with spatial tightness of the set~$\Gamma^D_N(t)$. Here we will also need:

\begin{mylemma}
\label{lemma-12.2ua}
For each $D\in\mathfrak D$, each~$\delta>0$ and each~$t\ge0$ there is $c>0$ such that for all $t'\in[-t,t]$, all $N\ge1$, all~$x\in D_N^\delta$ and all~$r$ with $1/\delta<r<\delta N$, 
\begin{equation}
P\Bigl(\,\max_{y\in A_{r,N/r}(x)}h^{D_N}_y\ge h^{D_N}_x,\,h^{D_N}\le m_N+t\,\Big|\, h^{D_N}_x=m_N+t'\Bigr)
\le c\frac{\,(\log r)^{-\frac1{16}}}{\log N},
\end{equation}
where we abbreviated~$A_{a,b}(x):=\{y\in\Z^2\colon a\le \dist_\infty(x,y)\le b\}$.
\end{mylemma}

\begin{proofsect}{Proof}
Let us again assume~$0\in D$ and start with~$x=0$. The set under the maximum is contained in~$\Delta^{n-k}\smallsetminus\Delta^k$ for a natural~$k$ proportional to~$\log r$.
Lemma~\ref{lemma-8.3} shifts the conditioning to $h^{D_N}_0=0$ at the cost of adding $(m_N+t')\frakg^{D_N}$ to all occurrences of~$h^{D_N}$. This embeds the event under consideration into
\begin{multline}
\qquad
\bigcup_{y\in \Delta^{n-k}\smallsetminus\Delta^k}
\bigl\{h^{D_N}_y\ge (m_N+t')(1-\frakg^{D_N}(y))\bigr\}
\\
\cap\bigl\{h^{D_N}\le (m_N+t)(1-\frakg^{D_N})+(t-t')\frakg^{D_N}\bigr\}.
\qquad
\end{multline}
For $|t'|\le t$, the event is the largest when~$t'=-t$ so let us assume that from now on.
Using that $\frakg^{D_N}\le1$ and invoking the estimates underlying Lemma~\ref{lemma-10.5}, the resulting event is then covered by the union of
\begin{equation}
\label{E:12.7ua}
\{K>k\}\cap\bigcap_{\ell=1}^{n}\{S_\ell\ge -R_K(\ell)-2t\}
\end{equation}
and
\begin{equation}
\label{E:12.8ua}
\biggl(\,\bigcup_{\ell=k+1}^{n-k}\{S_\ell\le R_k(\ell)+t\bigr\}\biggr)
\cap\bigcap_{\ell=1}^{n}\{S_\ell\ge -R_k(\ell)-2t\}\,.
\end{equation}
By Lemma~\ref{lemma-11.7}, the conditional probability (given~$S_{n+1}=0$) of the event in \eqref{E:12.7ua} is at most order $\texte^{-c_2(\log k)^2}/\log N$. Lemma~\ref{lemma-11.6} in turn bounds the corresponding conditional probability of the event in \eqref{E:12.8ua} by a constant times $k^{-\frac1{16}}/\log N$. The constants in these bounds are uniform on compact sets of~$t$ and over the shifts of~$D_N$ such that $0\in D_N^\delta$. Using the translation invariance of the DGFF, the claim thus  extends to all~$x\in D_N^\delta$.
\end{proofsect}

With Lemma~\ref{lemma-12.2ua} in hand, we are ready to give:

\begin{proofsect}{Proof of Theorem~\ref{thm-DZ2}}
Given~$\delta>0$ and~$s\ge t\ge0$, the event in the statement is the subset of the union of
\begin{equation}
\label{E:12.9ua}
\bigl\{\Gamma_N^D(t)\smallsetminus D_N^{\delta}\ne\emptyset\bigr\}\cup\bigl\{\,\max_{x\in D_N}h^{D_N}_x>m_N+s\bigr\}
\end{equation}
and the event
\begin{equation}
\label{E:12.10ua}
\bigcup_{x\in D_N^\delta}\biggl\{\,
\max_{y\in A_{r,N/r}(x)}h^{D_N}_y\ge h^{D_N}_x\ge m_N-t,\,h^{D_N}\le m_N+s\biggr\}
\end{equation}
provided~$1/\delta<r<\delta N$. Lemmas~\ref{lemma-10.4} and~\ref{lemma-8.2} show that the probability of the event in \eqref{E:12.9ua} tends to zero in the limits~$N\to\infty$, $\delta\downarrow0$ and~$s\to\infty$. The probability of the event in \eqref{E:12.10ua} is in turn estimated by the union bound and conditioning on~$h^{D_N}_x=m_N+t'$ for~$t'\in[-t,s]$. Lemma~\ref{lemma-12.2ua} then dominates the probability by an $s$-dependent constant times $(\log r)^{-\frac1{16}}$. The claim follows by taking~$N\to\infty$, $r\to\infty$ followed by $\delta\downarrow0$ and~$s\to\infty$.
\end{proofsect}

The arguments underlying Lemma~\ref{lemma-12.2ua} also allow us to prove that large local extrema are preserved (with high probability) upon application of the Gibbs-Markov property:

\begin{proofsect}{Proof of Lemma~\ref{lemma-same-max}}
Let~$D_N$ be a lattice approximation of a continuum domain $D\in\mathfrak D$ and let~$\wt D_N$ be the union of squares of side length (roughly) $N/K$ that fit into~$D_N^\delta$ with distinct squares at least two lattice steps apart. For $x\in\wt D_N$, let $S_{N,K}(x)$ denote the square (of side-length roughly $N/K$) in~$\wt D_N$ containing~$x$ and let $S_{N,K}^\delta(x)$ be the square of side-length $(1-\delta)(N/K)$ centered at the same point as~$S_{N,K}(x)$. Abusing our earlier notation, let $\wt D_N^\delta:=\bigcup_{x\in\wt D_N}S_{N,K}^\delta(x)$.  

Consider the coupling of $h^{D_N}$ and~$h^{\wt D_N}$ in the Gibbs-Markov property. The claim in Lemma~\ref{lemma-same-max} will follow by standard approximation arguments once we prove that, for every $t\ge0$,
\begin{equation}
\label{E:12.13nwxt}
P\biggl(
\exists x\in\Gamma_N^D(t)\cap \wt D_N^\delta\colon \!\!\max_{y\in S_{N,K}(x)}h^{ D_N}_y\le h^{ D_N}_x,\,
 \!\max_{y\in S_{N,K}(x)}h^{\wt D_N}_y> h^{\wt D_N}_x\biggr)
\end{equation}
tends to zero as~$N\to\infty$ and~$K\to\infty$.
For this it suffices to show that for all~$t>0$ there is~$c>0$ such that for all $N\gg K\ge1$, all $s\in[-t,t]$, all~$x\in\wt D_N^\delta$ and all~$y\in S_{N,K}(x)$,
\begin{multline}
\label{E:12.13nwwt}
P\Bigl(\,\max_{y\in S_{N,K}(x)}h^{ D_N}_y\le h^{ D_N}_x,\,
 \!\max_{y\in S_{N,K}(x)}h^{\wt D_N}_y> h^{\wt D_N}_x,
 \\h^{D_N}\le m_N+t\,\Big|\, h^{D_N}_x=m_N+s\Bigr)\le c\frac{(\log K)^{-\frac1{16}}}{\log N}.
\end{multline}
Indeed, multiplying \eqref{E:12.13nwwt} by the probability density of $h^{D_N}_x-m_N$, which is of order $N^{-2}\log N$, integrating over~$s\in[-t,t]$ and summing over~$x\in \wt D_N^\delta$, we get a bound on \eqref{E:12.13nwxt} with $\Gamma_N^D(t)\cap \wt D_N$ replaced by $\Gamma_N^D(t)\cap \wt D_N^\delta$ and the event restricted to~$h^{D_N}\le m_N+t$. These defects are handled via Lemma~\ref{lemma-10.4}, which shows that $P(\Gamma^D_N(t)\smallsetminus\wt D_N^\delta=\emptyset)\to1$ as~$N\to\infty$, $K\to\infty$ and~$\delta\downarrow0$, and Lemma~\ref{lemma-8.2}, which gives $P(h^{D_N}\not\le m_N+t)\le\texte^{-\tilde a t}$ uniformly in~$N\ge1$.

In order to prove \eqref{E:12.13nwwt}, we will shift the domains so that~$x$ is at the origin. We will again rely on the concentric decomposition with the following caveat: $S_{N,K}(0)$ is among the $\Delta^i$'s. 
More precisely, we take
\begin{equation}
\Delta^j:=\begin{cases}
\{x\in\Z^2\colon|x|_\infty\le 2^j\},\quad&\text{if }j=0,\dots,r-1,
\\
S_{N,K}^\delta(0),\quad&\text{if }j=r,
\\
\{x\in\Z^2\colon|x|_\infty\le 2^{j+\wt r-r}\},\quad&\text{if }j=r+1,\dots,n-1,
\\
D_N,\quad&\text{if }j=n,
\end{cases}
\end{equation}
where
\begin{equation}
\begin{aligned}
r&:=\max\bigl\{j\ge0\colon \{x\in\Z^2\colon|x|_\infty\le 2^{j}\}\subseteq S_{N,K}^\delta(0)\bigr\},
\\
\wt r&:=\min\bigl\{j\ge0\colon S_{N,K}(0)\subseteq \{x\in\Z^2\colon|x|_\infty\le 2^{j}\}\bigr\},
\\
n&:=\max\bigl\{j\ge0\colon \{x\in\Z^2\colon|x|_\infty\le 2^{j+\wt m-m+1}\}\subseteq D_N\bigr\}.
\end{aligned}
\end{equation}
We leave it to the reader to check that all the estimates pertaining to the concentric decomposition remain valid, albeit with constants that depend on~$\delta$.

Denoting~$k:=n-r$, we have $k=\log_2(K)+O(1)$ and $n=\log_2(N)+O(1)$. Assuming $\{h^{\Delta^j}\colon j=0,\dots,n\}$ are defined via the concentric decomposition, the probability in \eqref{E:12.13nwwt} becomes
\begin{multline}
\label{E:12.15nwwt}
\quad
P\Bigl(h^{\Delta^{n-k}}\not\le h^{\Delta^{n-k}}_0,\,h^{\Delta^n}\le h^{\Delta^n}_0\text{ in }\Delta^{n-k},
\\h^{\Delta^n}\le m_N+t\,\Big|\, h^{\Delta^n}_0=m_N+s\Bigr).
\quad
\end{multline}
Observe that, for all~$x\in\Delta^{n-k}$,
\begin{equation}
h^{\Delta^n}_x-h^{\Delta^n}_0=
h^{\Delta^{n-k}}_x-h^{\Delta^{n-k}}_0+\sum_{j=n-k+1}^n\bigl[\frakb_j(x)\varphi_j(0)+\frak\chi_j(x)\bigr].
\end{equation}
Assuming~$K\le k$, the estimates in the concentric decomposition bound the~$j$-th term in the sum by $(\log(k\vee j))^2\texte^{-c(n-k-j)}$ on~$\Delta^j$ for all $j=0,\dots,n-k-1$. This will suffice for~$j\le k$; for the complementary~$j$ we instead get
\begin{equation}
\sum_{j=n-k+1}^n\bigl[\frakb_j(x)\varphi_j(0)+\frak\chi_j(x)\bigr]+m_N(s,x)
\ge -R_k(j)-t
\end{equation}
for~$x\in\Delta^j\smallsetminus\Delta^{j-1}$, where we used that~$|s|\le t$ and invoked the shorthand \eqref{E:mNtx}. The event in \eqref{E:12.15nwwt} is thus contained in the intersection  of $\{h^{\Delta^n}\le m_N+t\}$ with the union of three events:~$\{K>k\}$,
\begin{equation}
\label{E:12.20nut}
\bigl\{h^{\Delta^n}_0\ge h^{\Delta^n}\text{ on~$\Delta^k$}\bigr\}\cap\bigcup_{x\in\Delta^k}
\bigl\{h^{\Delta^n}_x> h^{\Delta^n}_0-(\log k)^2\texte^{-c(2n-k)}\bigr\}
\end{equation}
and
\begin{equation}
\label{E:12.21nut}
\{K\le k\}\cap\bigcup_{x\in\Delta^{n-k}\smallsetminus\Delta^k}\Bigl\{h^{\Delta^n}_x-h^{\Delta^n}_0
+m_N(s,x)\ge -R_k(j)-t\Bigr\}. 
\end{equation}
Lemma~\ref{lemma-11.7} bounds the contribution of the first event by $c_1\texte^{-c_2(\log k)^2}/\log N$. For  \eqref{E:12.20nut} we drop the event~$h^{D^n}\le m_N+t$ and then invoke the ``gap estimate'' in Exercise~\ref{ex:11.9} to bound the requisite probability by a quantity of order $|\Delta^k|(\log k)^2\texte^{-c(2n-k)}$. For \eqref{E:12.21nut} we shift the conditioning to~$h^{\Delta^n}_0=0$, which turns the event under the union to $\{h^{\Delta^n}_x-h^{\Delta^n}_0
\ge -R_k(j)-t\}$. The resulting event is then contained in \eqref{E:12.8ua} which has probability less than order $k^{-\frac1{16}}/\log N$. As $k\approx\log_2(K)$ and $n\gg k$ we get \eqref{E:12.13nwwt}.
\end{proofsect}

\section{Limit of atypically large maximum}
\noindent
Our next item of business here is the proof of Theorem~\ref{prop-10.7} dealing with the limit law of the maximum (and the corresponding maximizer) conditioned on atypically large values.
This will again rely on the concentric decomposition of the DGFF; specifically, the formula
\begin{equation}
\label{E:12.11}
P\Bigl(\,h^{D_N}\le h^{D_N}_0\,\Big|\, h^{D_N}_0=m_N+t\,\Bigr)
=\frac{\Xiin_k(1)\Xiout_{N,k}(t)}{g\log N}\bigl(2+o(1)\bigr)
\end{equation}
from Lemma~\ref{lemma-11.12ueu}, where the quantities $\Xiin_k(1)$ and $\Xiout_{N,k}(t)$ are as in \twoeqref{E:Xi-in}{E:Xi-out} and where we incorporated the $o(1)$-term to the expression thanks to the bounds in Lemma~\ref{lemma-11.13ueu}. The novel additional input is:

\begin{myproposition}[Asymptotic contribution of outer layers]
\label{prop-12.3}
$\Xiout_{N,k}(t)\sim t$ in the limit $N\to\infty$, $k\to\infty$ and~$t\to\infty$. More precisely,
\begin{equation}
\label{E:12.14nwt}
\lim_{t\to\infty}\,\limsup_{k\to\infty}\,\limsup_{N\to\infty}\,\biggl|\,\frac{\Xiout_{N,k}(t)}t-1\,\biggr|=0.
\end{equation}
For each~$\delta>0$, the limit is uniform in the shifts of~$D_N$ such that $0\in D_N^\delta$.
\end{myproposition}

Deferring the proof of this proposition to later in this lecture, we will now show how this implies Theorem~\ref{prop-10.7}. Recall that $r_D(x)$ denotes the conformal radius of~$D$ from~$x$ and that~$\Xiin_\infty(1)$ is the limit of~$\Xiin_k(1)$ whose existence was established in the proof of Theorem~\ref{thm-11.3}; see \eqref{E:11.40nwt}. We start by:

\begin{mylemma}
For any $\delta>0$ and any $q>1$,
\begin{multline}
\label{E:12.13uau}
\quad
\lim_{t\to\infty}\,\limsup_{N\to\infty}\,\max_{x\in D^\delta_N}\,
\biggl|\,N^2\,\frac{\texte^{\alpha t}}t\, P\Bigl(h^{D_N}\le h^{D_N}_x,\,h^{D_N}_x-m_N\in[t,qt]\Bigr)
\\
-\sqrt{\frac\pi8}\,\Xiin_\infty(1)\,r_D(x/N)^2\biggr|=0\,.
\quad
\end{multline}
\end{mylemma}

\begin{proofsect}{Proof}
Denote by $f_{N,x}$ the density with respect to the Lebesgue measure of the distribution of~$h^{D_N}_x-m_N$ and write
\begin{multline}
\qquad
P\Bigl(h^{D_N}\le h^{D_N}_x,\,h^{D_N}_x-m_N\in[t,qt]\Bigr)
\\
=\int_t^{qt}\textd s\, f_{N,x}(s)\,
P\Bigl(h^{D_N}\le h^{D_N}_x\,\Big|\, h^{D_N}_x=m_N+s\Bigr)\,.
\qquad
\end{multline}
A calculation based on the Green function asymptotic in Theorem~\ref{thm-1.17} now shows that, with $o(1)\to0$ as $N\to\infty$ uniformly in~$x\in D^\delta$ and $s\in[t,qt]$,
\begin{equation}
\frac{N^2}{\log N}\, f_{N,\lfloor Nx\rfloor}(s)= \bigl(1+o(1)\bigr)\frac1{\sqrt{2\pi g}}\,\texte^{-\alpha s}\,r_D(x)^2.
\end{equation}
Combining this with \eqref{E:12.11},  \eqref{E:12.14nwt} and  \eqref{E:11.40nwt} gives
\begin{multline}
\qquad
N^2\,\frac{\texte^{\alpha t}}t\, P\Bigl(h^{D_N}\le h^{D_N}_{x},\,h^{D_N}_{x}-m_N\in[t,qt]\Bigr)
\\
= \bigl(1+o(1)\bigr)\,\frac2g\,\Xiin_\infty(1)\,\frac1{\sqrt{2\pi g}}\,r_D(x/N)^2\,\Bigl[\,\frac{\texte^{\alpha t}}t\int_t^{qt}\textd s\,\texte^{-\alpha s}s\Bigr],
\qquad
\end{multline}
where~$o(1)\to0$ as~$N\to\infty$ followed by $t\to\infty$ uniformly in~$x\in\Z^2$ such that $x/N\in D^\delta$. The expression inside the square brackets tends to~$\alpha^{-1}$ as~$t\to\infty$. 
Since $(2/g)(2\pi g)^{-1/2}\alpha^{-1} = (\sqrt{2\pi}\,g)^{-1}=\sqrt{\pi/8}$, we get \eqref{E:12.13uau}. 
\end{proofsect}

In order to rule out the occurrence of multiple close points with excessive values of the field, and thus control certain boundary issues in the proof of Theorem~\ref{prop-10.7}, we will also need:

\begin{mylemma}
\label{lemma-12.8nqqt}
For each $\delta>0$ small enough there is~$c>0$ such that for all $N\ge 1$, all $x\in D_N$ and all $0\le s, t\le\log N$,
\begin{equation}
P\Bigl(\,\max_{\begin{subarray}{c}
y\in D_N\\\dist_\infty(x,y)>\delta N
\end{subarray}}
h^{D_N}_y\ge m_N+s\,\Big|\,h^{D_N}\le h^{D_N}_x=m_N+t\Bigr)
\le c\,\texte^{-\tilde a s},
\end{equation}
where~$\tilde a>0$ is as in Lemma~\ref{lemma-8.2}.
\end{mylemma}

\begin{proofsect}{Proof}
The FKG inequality for $h^{D_N}_x$ conditioned on $h^{D_N}_x=m_N+t$ bounds the desired probability by
\begin{equation}
P\Bigl(\,\max_{\begin{subarray}{c}
y\in D_N\\\dist_\infty(x,y)>\delta N
\end{subarray}}
h^{D_N}_y\ge m_N+s\,\Big|\,h^{D_N}_x=m_N+t\Bigr).
\end{equation}
Lemma~\ref{lemma-8.3} along with Exercise~\ref{ex:8.4} then permit us to rewrite this as the probability that there exists $y\in D_N$ with $\dist_\infty(x,y)>\delta N$ such that
\begin{equation}
h^{D_N\smallsetminus\{0\}}_y\ge m_N+s-(m_N+t)\frakg^{-x+D_N}(y-x).
\end{equation}
The Maximum Principle shows that $A\mapsto\frakg^A$ is non-decreasing with respect to the set inclusion and so $\frakg^{-x+D_N}(\cdot)\le\frakg^{\wt D_N}(\cdot)$, where~$\wt D_N$ is a box of side $4\diam_\infty(D_N)$ centered at the origin.
Under our assumptions on~$t$, the term $(m_N+t)\frakg^{\wt D_N}(y-x)$ is at most a constant times $\log(1/\delta)$. The claim then follows from the upper-tail tightness of the maximum in Lemma~\ref{lemma-8.2}.
\end{proofsect}

We are ready to give:

\begin{proofsect}{Proof of Theorem~\ref{prop-10.7}}
We may as well give the proof for a general sequence of approximating domains~$D_N$. Set
\begin{equation}
\label{E:12.17uai}
\bar c:=\sqrt{\frac\pi8}\,\,\Xiin_\infty(1).
\end{equation}
Let~$A\subseteq D$ be open and denote by $X^\star_N$ the (a.s.-unique) maximizer of~$h^{D_N}$ over~$D_N$. Abbreviate $M_N:=\max_{x\in D_N}h^{D_N}_x$. Summing \eqref{E:12.13uau} then yields
\begin{multline}
\label{E:12.18uau}
\quad
\lim_{\delta\downarrow0}\,\limsup_{t\to\infty}\,\limsup_{N\to\infty}\,\biggl|\frac{\texte^{\alpha t}}t\,P\Bigl(\frac1NX_N^\star\in D^\delta\cap A,\,M_N-m_N\in[t,qt]\Bigr)
\\
-\bar c\int_A\textd x\,r_D(x)^2\,\biggr|=0\,.
\qquad
\end{multline}
Assuming~$q>1$ obeys $q\tilde a>2\alpha$, for~$\tilde a$ as in Lemma~\ref{lemma-8.2}, the probability of $M_N\ge m_N+qt$  decays as $\texte^{-2\alpha t}$ and so \eqref{E:12.18uau} holds with $M_N-m_N\in[t,qt]$ replaced by $M_N\ge m_N+t$. This would already give us a lower bound in \eqref{E:10.11} but, in order to prove the complementary upper bound as well as \eqref{E:10.11ueu}, we need to address the possibility of the maximizer falling outside~$D^\delta_N$. 

Let~$\wt D_N$ be the square of side $4\diam_\infty(D_N)$ centered at an arbitrary point of~$D_N$.
Exercise~\ref{ex:3.4} then gives
\begin{equation}
\label{E:12.20uau}
P\Bigl(\,\max_{x\in D_N\smallsetminus D_N^\delta}h^{D_N}_x>m_N+t\Bigr)
\le 2P\Bigl(\,\max_{x\in D_N\smallsetminus D_N^\delta}h^{\wt D_N}_x>m_N+t\Bigr).
\end{equation}
Denote $B_N^\delta:=\{x\in\wt D_N\colon \dist_\infty(x,D_N\smallsetminus D_N^\delta)<\delta N\}$
and write $(\wt X_N^\star,\wt M_N)$ for the maximizer and the maximum of $h^{\wt D_N}$ in~$\wt D_N$. The probability on the right-hand side of \eqref{E:12.20uau} is then bounded by the sum of
\begin{equation}
\label{E:12.32nqqt}
P\bigl( \wt X_N^\star\in B_N^\delta,\,\wt M_N>m_N+t\bigr)
\end{equation}
and
\begin{equation}
\label{E:12.33nqqt}
\sum_{x\in\wt D_N}\,P\biggl(h^{\wt D_N}\le h^{\wt D_N}_x,\,h^{\wt D_N}_x>m_N+t,\!\!\!\!
\max_{\begin{subarray}{c}
y\in \wt D_N\\
\dist_\infty(x,y)>\delta N
\end{subarray}}\!
h^{D_N}_y> m_N+t\biggr).
\end{equation}
Once~$\delta$ is sufficiently small, we have~$B_N^\delta\subset\wt D_N^{\delta'}$ for some~$\delta'>0$ and so we may bound \eqref{E:12.32nqqt} using \eqref{E:12.18uau} by a quantity of order $t\texte^{-\alpha t}|B_N^\delta|/N^2$. Since $|B_N^\delta|/N^2$ is bounded by a constant multiple of $\Leb(D\smallsetminus D^\delta)$, this term is $o(t\texte^{-\alpha t})$ in the limits $N\to\infty$, $t\to\infty$ and~$\delta\downarrow0$. 

The term in \eqref{E:12.33nqqt} is bounded using Lemma~\ref{lemma-12.8nqqt} by $c\texte^{-\tilde a t}$ times the probability $P(\wt M_N> m_N+t)$. By Lemma~\ref{lemma-10.4}, this probability is (for~$t\ge1$) at most a constant times~$t^2\texte^{-\alpha t}$ and so \eqref{E:12.33nqqt} is $o(t\texte^{-\alpha t})$ as well. Hence, for $o(1)\to0$ as~$N\to\infty$ and~$t\to\infty$,
\begin{equation}
P\Bigl(\frac1N X_N^\star\in A,\,M_N\ge m_N+t\Bigr)=t\texte^{-\alpha t} \Bigl(o(1)+\bar c\int_A\textd x\,r_D(x)^2\Bigr).
\end{equation}
This yields the claim with~$\psi$ equal to the normalized conformal radius squared (and $c_\star:=\bar c\int_S \textd x\,r_S(x)^2$ for~$S:=(0,1)^2$).
\end{proofsect}

As noted earlier, our proof identifies~$\psi$ with (a constant multiple of) $r_D(x)^2$ directly; unlike the proof of part~(5) of Theorem~\ref{thm-10.14} which takes the existence of some~$\psi$ as an input and infers equality with the conformal radius squared using the Gibbs-Markov property of the~$Z^D$-measures (whose uniqueness we have yet to prove).
For future use, we pose a minor extension of Theorem~\ref{prop-10.7}:

\begin{myexercise}
\label{ex:12.5uai}
For any bounded continuous~$f\colon\overline D\to[0,\infty)$,
\begin{multline}
\quad
P\Bigl(\,h^{D_N}_x\le m_N+t-\alpha^{-1}\log f(x/N)\colon x\in D_N\Bigr)
\\=\exp\biggl\{-t\texte^{-\alpha t}\Bigl(o(1)+\bar c\int_D \textd x\,\,r_D(x)^2\,f(x)\Bigr)\biggr\},
\quad
\end{multline}
where~$o(1)\to0$ in the limits as~$N\to\infty$ and~$t\to\infty$ uniformly on compact families of non-negative~$f\in C(\overline D)$.
\end{myexercise}

\section{Precise upper tail asymptotic}
\noindent
We now move to the proof of Proposition~\ref{prop-12.3}.
We begin by removing the conditioning on $h^{D_N}_0=0$ from the quantity of interest. To this end we note that
\begin{equation}
h^{D_N}-h^{D_N}_0\frakg^{D_N}\,\,\laweq\,\, h^{D_N\smallsetminus\{0\}}
\end{equation}
which (by Exercise~\ref{ex:8.4}) also has the law of~$h^{D_N}$ conditioned on~$h^{D_N}_0=S_{n+1}=0$. Since $\frakg^{D_N}(x)\le c/\log N$ on~$D_N\smallsetminus\Delta^{n-k}$ with some~$c=c(k)$ while (by a first moment bound) $P(h^{D_N}_0\ge \epsilon\log N)\to0$ as~$N\to\infty$, the term~$h^{D_N}_0\frakg^{D_N}$ can be absorbed into a negligible change of~$t$. Since $\Xiout_{N,k}(t)$ is non-decreasing in~$t$, this gives the main idea how to solve:

\begin{myexercise}
Abbreviate $\wt m_N(x,t):=m_N(1-\frakg^{D_N}(x))+t$ and let
\begin{equation}
\label{E:Xi-out2}
\wtXiout_{N,k}(t):=E\biggl(\1_{\{h^{D_N}\le \wt m_N(t,\cdot)\text{\rm\ in }D_N\smallsetminus\Delta^{n-k}\}}
\1_{\{\wt S_k\in[k^{1/6},\,k^2]\}}\wt S_k\biggr),
\end{equation}
where we abbreviated
\begin{equation}
\wt S_\ell:=S_{n-\ell}-S_{n+1},\quad \ell=0,\dots,k.
\end{equation}
Prove that if, for some~$\delta>0$,
\begin{equation}
\label{E:12.14nwt2}
\lim_{t\to\infty}\,\limsup_{k\to\infty}\,\limsup_{N\to\infty}\,\biggl|\,\frac{\wtXiout_{N,k}(t)}t-1\,\biggr|=0
\end{equation}
holds uniformly in the shifts of~$D_N$ such that~$0\in D_N^\delta$, then the same applies to $\Xiout_{N,k}(t)$ --- meaning that \eqref{E:12.13uau} holds.
\end{myexercise}

We thus have to prove the claim for $\wtXiout_{N,k}(t)$ instead of $\Xiout_{N,k}(t)$. A key point is that the removal of the conditioning on $S_{n+1}=0$ makes the increments of~$\wt S$ independent (albeit still slightly time-inhomogeneous) and so we can analyze~$\wtXiout_{N,k}(t)$ using the ballot problem/entropic-repulsion methods used earlier.

Recall the basic objects $\varphi_j$, $\chi_j$ and~$h_j'$ from the concentric decomposition and define a (new) control random variable~$L$ as the minimal index in $\{1,\dots,k-1\}$ such that for all $j=0,\dots,k$,
\begin{equation}
\bigl|\varphi_{n-j}(0)\bigr|\le[\log(L\vee j)]^2,
\end{equation}
\begin{equation}
\max_{j+1\le r\le k}\,\max_{x\in \Delta^{n-r}}\bigl|\chi_{n-j}(x)\bigr|2^{(j-r)/2}\le[\log(L\vee j)]^2
\end{equation}
and
\begin{equation}
\Bigl|\,\max_{x\in\Delta^{n-j}\smallsetminus\Delta^{n-j-1}}\bigl[\chi_{n-j}(x)+\chi_{n-j+1}(x)+h_{n-j}'(x)\bigr]-m_{2^{n-j}}\Bigr|\le[\log(L\vee j)]^2.
\end{equation}
If no such index exists, we set~$L:=k$. We then have:

\begin{mylemma}
\label{lemma-12.9nwt}
There is an absolute constant~$C\in(0,\infty)$ such that for all $t\in\R$, all $n\ge1$, all $k=1,\dots,n$ obeying
\begin{equation}
\label{E:tkn}
0\le tk/n\le 1
\end{equation}
and all $\ell=0,\dots,k$, on the event $\{L<k\}$ we have
\begin{equation}
\Bigl|\,\max_{x\in\Delta^{n-\ell}\smallsetminus\Delta^{n-\ell-1}}
\bigl[h^{D_N}_x- \wt m_N(t,x)\bigr]+(\wt S_\ell+t)\Bigr|
\le\zeta(\ell\vee L),
\end{equation}
where
\begin{equation}
\zeta(s):=C\bigl[1+\log(s)]^2.
\end{equation}
\end{mylemma} 

\begin{proofsect}{Proof}
From \eqref{E:8.29}, \eqref{E:8.48nwt} and the fact that~$h'_j$ is supported on~$\Delta^j\smallsetminus\Delta^{j-1}$ while~$\chi_j$ is supported on~$\Delta^j$ we get, for all $\ell=0,\dots,k$ and all $x\in\Delta^{n-\ell}\smallsetminus\Delta^{n-\ell-1}$,
\begin{equation}
h^{D_N}_x = S_{n+1}-S_{n-\ell}+h_{n-\ell}'(x)+\sum_{j=n-\ell-1}^n\bigl(\frakb_j(x)\varphi_j(0)+\chi_j(x)\bigr).
\end{equation}
The definition of~$L$ along with the decay of~$j\mapsto\frakb_j(x)$ in \eqref{E:8.36} bound the sum by a quantity proportional to~$1+\log(L\vee\ell)^2$. The claim follows from the definition of $\wt m_N(t,\cdot)$, the tightness of the absolute maximum and Exercise~\ref{ex:11.6}.
\end{proofsect}

As in our earlier use of the concentric decomposition, Lemma~\ref{lemma-12.9nwt} permits us to squeeze the event of interest between two random-walk events:

\begin{myexercise}
\label{ex:12.7}
Prove that, under \eqref{E:tkn},
\begin{multline}
\label{E:12.33nwt}
\qquad
\{L< k\}\cap\bigcap_{\ell=0}^{k}\bigl\{\wt S_\ell\ge\zeta(\ell\vee L)-t\bigr\}
\\\subseteq\
\{L< k\}
\cap\bigl\{h^{D_N}\le \wt m_N(t,\cdot)\text{\rm\ in }D_N\smallsetminus\Delta^{n-k}\bigr\}
\\
\subseteq\bigcap_{\ell=0}^{k}\bigl\{\wt S_\ell\ge-\zeta(\ell\vee L)-t\bigr\}.
\qquad
\end{multline}
\end{myexercise}

\noindent
Notice that, unlike in Lemma~\ref{lemma-8.16}, the $t$-dependent term has the same form and sign on both ends of this inclusion. The next item to address are suitable estimates on the relevant probabilities for the control variable~$L$ and the random walk~$\wt S$.
First, a variation on Exercise~\ref{ex:8.17} shows
\begin{equation}
\label{E:12.42uau}
P\bigl(L=\ell)\le c_1\texte^{-c_2(\log\ell)^2},\quad\ell=1,\dots,k,
\end{equation}
for some~$c_1,c_2>0$. In addition, we will need:

\begin{mylemma}
\label{lemma-12.10}
There is~$c>0$ such that for all $\ell,k\in\N$ with $\ell\le k^{\frac1{12}}$ and~$k\le n/2$ and all~$0\le t\le k^{1/6}$,
\begin{equation}
P\biggl(\,\bigcap_{j=\ell}^{k}\bigl\{\wt S_j-\wt S_\ell\ge-\zeta(j)-t-\ell^2\bigr\}\biggr)\le c\frac{\ell^2+t}{\sqrt k}.
\end{equation}
\end{mylemma}

\noindent
This is showed by invoking again estimates on the inhomogeneous ballot problem encountered earlier in these lectures. We refer to \cite{BL3} for details. We use Lemma~\ref{lemma-12.10} to prove:

\begin{mylemma}[Controlling the control variable]
\label{lemma-12.11uau}
There are $c_1,c_2>0$ such that for all $\ell,k,n\in\N$ satisfying $k\le n/2$ and $1\le\ell\le k^{\frac1{12}}$,  and for all~$1\le t\le k^{1/6}$,
\begin{equation}
\label{E:12.44uau}
E\Bigl(\,(\wt S_k\vee k^{1/6})1_{\{L=\ell\}}\prod_{j=\ell}^{k}1_{\{\wt S_j\ge-\zeta(j)-t\bigr\}}\Bigr)
\le c_1(\ell^2+t)\texte^{-c_2(\log \ell)^2}
\end{equation}
and also
\begin{equation}
\label{E:12.45uau}
E\Bigl(\,\wt S_k1_{\{0\le\wt S_k\le k^{1/6}\}}\,1_{\{L\le\ell\}}\prod_{j=\ell}^{k}1_{\{\wt S_j\ge-\zeta(j)-t\bigr\}}\Bigr)
\le c_1(\ell^2+t)k^{-1/3}.
\end{equation}
\end{mylemma}

\begin{proofsect}{Proof}
Since~$|\wt S_\ell|\le\ell^2$ when $L=\ell$, we have
\begin{equation}
\wt S_k\vee k^{1/6}\le \ell^2+(\wt S_k-\wt S_\ell)\vee k^{1/6},\quad\text{on }\{L=\ell\}
\end{equation}
and
\begin{equation}
1_{\{L=\ell\}}1_{\{\wt S_j\ge-\zeta(j)-t\}}\le 1_{\{L=\ell\}}1_{\{\wt S_j-\wt S_\ell\ge-\zeta(j)-t-\ell^2\}}
\end{equation}
Substituting these in \eqref{E:12.44uau}, the indicator of~$\{L=\ell\}$ is independent of the rest of the expression. Invoking \eqref{E:12.42uau}, to get \eqref{E:12.44uau} we thus have to show
\begin{equation}
\label{E:12.47uau}
E\Bigl(\,
\bigl((\wt S_k-\wt S_\ell)\vee k^{1/6}\bigr)\,\prod_{j=\ell}^{k}1_{\{\wt S_j-\wt S_\ell\ge-\zeta(j)-t-\ell^2\}}\Bigr)
\le c(\ell^2+t)
\end{equation}
for some constant~$c>0$. Abbreviate
\begin{equation}
\wh S_j:=\wt S_{\ell+j}-\wt S_\ell.
\end{equation}
We will need:

\begin{myexercise}
Prove that the law of $\{\wh S_j\colon j=1,\dots,k-\ell\}$ on~$\R^{k-\ell}$ is strong FKG. Prove the same for the law of the standard Brownian motion on~$\R^{[0,\infty)}$.
\end{myexercise}

\begin{myexercise}
\label{ex:12.13}
Prove that if the law of random variables~$X_1,\dots,X_n$ is strong FKG, then for any increasing function~$Y$,
\begin{equation}
a_1,\dots,a_n\mapsto E\bigl(Y\,|\,X_1\ge a_1,\dots,X_n\ge a_n\bigr)
\end{equation}
is increasing in each variable.
\end{myexercise}

\noindent
Embedding the random walk into the standard Brownian motion~$\{B_s\colon s\ge0\}$ via $\wh S_j:=B_{s_j}$, where~$s_j:=\Var(\wh S_j)$, and writing $P^0$ is the law of the Brownian motion started at~$0$ and~$E^0$ is the associated expectation, Exercise~\ref{ex:12.13} dominates the expectation in \eqref{E:12.47uau} by
\begin{multline}
\label{E:12.50uau}
\quad
E^0\Bigl((B_{s_{k-\ell}}\vee k^{1/6})\,1_{\{B_s\ge-1\colon s\in[0,s_{k-\ell}]\}}\Bigr)
\\
\times\frac
{P\bigl(\,\bigcap_{j=0}^{k-\ell+1}\{\wh S_j\ge-\zeta(j+\ell)-t-\ell^2\}\bigr)}
{P^0\bigl(B_s\ge-1\colon s\in[0,s_{k-\ell}]\bigr)}.
\quad
\end{multline}
The expectation in \eqref{E:12.50uau} is bounded by conditioning on~$B_t$ and invoking Exercise~\ref{ex:7.9} along with the bound $1-\texte^{-a}\le a$ to get 
\begin{multline}
\label{E:12.51uau}
\quad
E^0\bigl((B_{s_{k-\ell}}\vee k^{1/6})\,1_{\{B_s\ge-1\colon s\in[0,s_{k-\ell}]\}}\bigr)
\\
\le E^0\biggl((B_{s_{k-\ell}}\vee k^{1/6})\frac{2 (B_{s_{k-\ell}}+1)}{s_{k-\ell}}\biggr).
\quad
\end{multline}
As~$s_{k-\ell}$ is of the order of~$k-\ell$, the expectation on the right is bounded uniformly in~$k>\ell\ge0$.
Thanks to Lemma~\ref{lemma-12.10} and Exercise~\ref{ex:11.16} (and the fact that $\zeta(\ell)\le\ell^2$ and that $j\mapsto\zeta(j+\ell)-\zeta(\ell)$ is slowly varying) the ratio of probabilities in \eqref{E:12.50uau} is at most a constant times $\ell^2+t$. 

The proof of \eqref{E:12.45uau} easier; we bound~$\wt S_k$ by~$k^{1/6}$ and thus estimate the probability by passing to $\wh S_j$ and invoking Lemma~\ref{lemma-12.10}.
\end{proofsect}

The above proof shows that the dominating contribution to~$\wtXiout_{N,k}(t)$ comes from the event when~$\wt S_{k}$ is order~$\sqrt k$. This is just enough to cancel the factor~$1/\sqrt{k}$ arising from the probability that~$\wt S$ stays above a slowly varying curve for~$k$ steps. This is what yields the desired asymptotic as~$t\to\infty$ as well:

\begin{proofsect}{Proof of Proposition~\ref{prop-12.3}}
Note that Lemma~\ref{lemma-12.11uau} and the fact that~$|S_k|\le k^2$ on~$\{L<k\}$ permits us to further reduce the computation to the asymptotic of
\begin{equation}
\label{E:Xi-out3}
\whXiout_{N,k}(t):=E\biggl(\1_{\{h^{D_N}\le \wt m_N(t,\cdot)\text{\rm\ in }D_N\smallsetminus\Delta^{n-k}\}}
(\wt S_k\vee0)\biggr).
\end{equation}
Indeed, writing~$\epsilon_\ell(t)$ for the right-hand side of \eqref{E:12.44uau}, this bound along with \eqref{E:12.33nwt} yield
\begin{multline}
\label{E:12.45nwt}
\quad
E\Bigl(\,(\wt S_k\vee 0)\prod_{j=0}^{k}1_{\{\wt S_{j}\ge\zeta(j\vee \ell)-t\}}\Bigr)
-\epsilon_\ell(t)
\\
\le
\whXiout_{N,k}(t)
\le
E\Bigl(\,(\wt S_k\vee 0)\prod_{j=0}^{k}1_{\{\wt S_{j}\ge-\zeta(j\vee \ell)-t\}}\Bigr)
+\epsilon_\ell(t).
\quad
\end{multline}
Similarly as in \eqref{E:12.50uau}, and denoting $\tilde s_k:=\Var(\wh S_k)$, Exercise~\ref{ex:12.13} bounds the expectation on the right of \eqref{E:12.45nwt} by
\begin{equation}
\label{E:12.54uau}
E^0\bigl((B_{\tilde s_k}\vee0)1_{\{B_s\ge -t\colon s\in[0,\tilde s_k]\}}\bigr)
\end{equation}
times the ratio
\begin{equation}
\label{E:12.55uau}
\frac{P\bigl(\bigcap_{j=1}^{k+1}\{\wt S_{j}\ge-\zeta(j\vee \ell)-t\bigr)}{P^0\bigl(B_s\ge -t\colon s\in[0,\tilde s_k]\bigr)}
\end{equation}
and that on the left by \eqref{E:12.54uau} times
\begin{equation}
\label{E:12.56uau}
\frac{P^0\bigl(B_s\ge \tilde\zeta(s\vee \tilde s_\ell)-t\colon s\in[0,\tilde s_k]\bigr)}{P^0\bigl(B_s\ge -t\colon s\in[0,\tilde s_k]\bigr)},
\end{equation}
where~$\tilde\zeta$ is a piecewise-linear (non-decreasing) function such that~$\tilde\zeta(\tilde s_j)=\zeta(j)$ for each~$j\ge1$.
The precise control of inhomogeneous Ballot Theorem (cf~\cite[Proposition~4.7 and~4.9]{BL3}) then shows that the ratios \twoeqref{E:12.55uau}{E:12.56uau} converge to one in the limits~$k\to\infty$ followed by~$t\to\infty$. Exercise~\ref{ex:7.9} shows that the expectation in \eqref{E:12.54uau} is asymptotic to~$t$ and so the claim follows.
\end{proofsect}

\section{Convergence of the DGFF maximum}
\noindent
We are now ready to address the convergence of the centered DGFF maximum by an argument that will also immediately yield the convergence of the extremal value process. The key point is the proof of:

\begin{mytheorem}[Uniqueness of $Z^D$-measure]
\label{thm-Z-unique} 
Let~$D\in\mathfrak D$ and let~$Z^D$ be the measure related to a subsequential weak limit of $\{\eta^D_{N,r_N}\colon N\ge1\}$ as in Theorem~\ref{thm-extremal-vals}. Then $Z^D$ is the weak limit of slight variants of the measures in \eqref{E:10.45}. In particular, its law is independent of the subsequence used to define it and the processes $\{\eta^D_{N,r_N}\colon N\ge1\}$ thus converge in law.
\end{mytheorem}

The proof will follow closely the proof of Theorem~\ref{thm-10.16} that characterizes the $Z^D$-measures by a list of natural conditions. These were supplied in Theorem~\ref{thm-10.14} for all domains of interest through the existence of the limit of processes $\{\eta^D_{N,r_N}\colon N\ge1\}$. Here we allow ourselves only subsequential convergence and so we face the following additional technical difficulties:
\settowidth{\leftmargini}{(1111)}
\begin{enumerate}
\item[(1)] We can work only with a countable set of domains at a time, and so the Gibbs-Markov property needs to be restricted accordingly.
\item[(2)] We lack the transformation rule for the behavior of the $Z^D$ measures under dilations of~$D$ in property~(3) of Theorem~\ref{thm-10.14}.
\item[(3)] We also lack the Laplace transform asymptotic in property~(5) of Theorem~\ref{thm-10.14}.
\end{enumerate}
The most pressing are the latter two issues so we start with them. It is the following lemma where the majority of the technical work is done:

\begin{mylemma}
\label{lemma-12.7ua}
Fix~$R>1$ and~$D\in\mathfrak D$ and for any positive~$n\in\N$ abbreviate~$D^n:=n^{-1}D$. Suppose that (for a given sequence of approximating domains for each~$n\in\N$), there is~$N_k\to\infty$ such that
\begin{equation}
\eta^{D^n}_{N_k,r_{N_k}}\lawarrow\,\, \eta^{D^n},\quad n\ge1.
\end{equation}
Let~$Z^{D^n}$ be the measure related to~$\eta^{D_n}$ as in Theorem~\ref{prop-subseq}. Then for any choice of $\lambda_n>0$ satisfying
\begin{equation}
\label{E:12.23uau}
\lambda_n n^{-4}\,\underset{n\to\infty}\longrightarrow\,0
\end{equation}
and any choice of $Rn$-Lipschitz functions $f_n\colon D^n\to[0,R]$, we have
\begin{multline}
\label{E:12.54nwwt}
\quad
\log E\bigl(\texte^{-\lambda_n\langle Z^{D^n},f_n\rangle}\bigr)
\\
=-\lambda_n\log(n^{4}/\lambda_n)\biggl[o(n^{-4})+\bar c\int_{D^n} \textd x\,\,r_{D^n}(x)^2f_n(x)\biggr]\,,
\quad
\end{multline}
where~$n^4o(n^{-4})\to0$ as $n\to\infty$ uniformly in the choice of the~$f_n$'s. The constant~$\bar c$ is as in \eqref{E:12.17uai}.
\end{mylemma}

\begin{proofsect}{Proof}
We start by some observations for a general domain~$D$. Let~$\varrho\colon \overline D\to[0,\infty)$ be measurable and denote~$A_\varrho:=\{(x,h)\colon x\in D,\,h>-\alpha^{-1}\log\varrho(x)\}$. Then for any~$a>0$, any subsequential limit process~$\eta^D$ obeys
\begin{equation}
(1-\texte^{-a})P\bigl(\eta^D(A_\varrho)>0\bigr)\le 1-E\bigl(\texte^{-a\eta^D(A_\varrho)}\bigr)\le P\bigl(\eta^D(A_\varrho)>0\bigr).
\end{equation}
Assuming the probability on the right-hand side is less than one, elementary estimates show
\begin{equation}
\label{E:12.56nwwt}
\Bigl|\,\log E\bigl(\texte^{-a\eta^D(A_\varrho)}\bigr)-\log P\bigl(\eta^D(A_\varrho)=0\bigr)
\Bigr|\le\frac{\texte^{-a}}{P(\eta^D(A_\varrho)=0)}.
\end{equation}
We also note that, setting
\begin{equation}
\label{E:12.53nwt}
f(x):=(1-\texte^{-a})\alpha^{-1}\varrho(x),
\end{equation}
implies 
\begin{equation}
\label{E:12.58nwwt}
E(\texte^{-a\eta^D(A_\varrho)}) = E(\texte^{-\langle Z^D,f\rangle})
\end{equation}
thanks to the Poisson structure of~$\eta^D$ proved in Theorem~\ref{prop-subseq}.
 
Given~$D\in\mathfrak D$, for each~$n\in\N$ let~$D^n:=n^{-1}D$ and let $f_n\colon D^n\to[0,R]$ be an $Rn$-Lipschitz function. Fix an auxiliary sequence $\{a_n\}$ of numbers $a_n\ge1$ such that $a_n\to\infty$ subject to a minimal growth condition to be stated later. Define $\varrho_n\colon D^n\to[0,\infty)$ so that
\begin{equation}
\label{E:12.59nwwt}
f_n(x)=(1-\texte^{-a_n})\alpha^{-1}\varrho_n(x)
\end{equation}
holds for each~$n\in\N$ and note that $\varrho_n(x)$ is $R'n$-Lipschitz for some~$R'$ depending only on~$R$. Let~$\{\lambda_n\colon n\ge1\}$ be positive numbers such that \eqref{E:12.23uau} holds and 
note that, by \eqref{E:12.56nwwt} and \eqref{E:12.58nwwt}, in order to control the asymptotic of $E(\texte^{-\lambda_n\langle Z^{D^n},f_n\rangle})$ in the limit as~$n\to\infty$ we need to control the asymptotic of $P(\eta^{D^n}(A_{\lambda_n\varrho_n})=0)$.

Writing~$\{D^n_N\colon N\ge1\}$ for the approximating domains for~$D^n$, the assumed convergence yields
\begin{multline}
\label{E:12.61nwwt} 
P\bigl(\eta^{D^n}(A_{\lambda_n\varrho_n})=0\bigr) 
\\= \lim_{k\to\infty}P\Bigl(h^{D^n_{N}}_x\le m_{N}-\alpha^{-1}\log\bigl(\lambda_n\varrho_n(x/N)\bigr)\colon x\in D^n_{N}\Bigr)\Bigl|_{N=N_k}.
\end{multline}
A key point is that, by Theorem~\ref{prop-10.7} (or, more precisely, Exercise~\ref{ex:12.5uai}) the limit of the probabilities on the right of \eqref{E:12.61nwwt} exists for all~$N$; not just those in the subsequence $\{N_k\}$.
We may thus regard~$D^n_N$ as an approximating domain of~$D$ at scale~$\lfloor N/n\rfloor$ and note that, in light of $2\sqrt g = 4\alpha^{-1}$, 
\begin{equation}
m_{N}-\alpha^{-1}\log\lambda_n = m_{\lfloor N/n\rfloor} + \alpha^{-1}\log(n^4/\lambda_n)+o(1),
\end{equation}
where $o(1)\to0$ as~$N\to\infty$. Note also that $\tilde \varrho_n\colon D\to\R$ defined by $\tilde \varrho_n(x):=\varrho_n(x/n)$ is $R'$-Lipschitz for some~$R'$ depending only on~$R$. Using \eqref{E:12.23uau}, Exercise~\ref{ex:12.5uai} with~$t=\alpha^{-1}\log(n^4/\lambda_n)+o(1)$ shows
\begin{multline}
\log P\Bigl(h^{D^n_{N}}_x\le m_{N}-\alpha^{-1}\log\lambda_n-\alpha^{-1}\log\varrho_n(x/N)\colon x\in D^n_{N}\Bigr)
\\
=-\alpha^{-1}\lambda_n\,n^{-4}\log(n^{4}/\lambda_n)\biggl[\,o(1)+\bar c\int_D\textd x\,r_D(x)^2\tilde \varrho_n(x)\biggr],
\end{multline}
where $o(1)\to0$ as~$N\to\infty$ and~$n\to\infty$.
Changing the variables back to~$x\in D^n$ absorbs the term~$n^{-4}$ into the integral and replaces~$\tilde \varrho_n$ by~$\varrho_n$, which can be then related back to~$f_n$ via \eqref{E:12.59nwwt}. Since $a_n\to\infty$, we get
\begin{multline}
\label{E:12.63nwwt}
\log P\bigl(\eta^{D^n}(A_{\lambda_n\varrho_n})=0\bigr)
\\
=-\,\frac{\lambda_n\log(n^{4}/\lambda_n)}{1-\texte^{-a_n}}
\biggl[\,o(n^{-4})+\bar c\int_{D^n}\textd x\,\,r_{D^n}(x)^2f_n(x)\biggr],
\end{multline}
where $o(n^{-4})n^4\to0$ as~$n\to\infty$ uniformly in the choice of $\{f_n\}$ with the above properties.

In remains to derive \eqref{E:12.54nwwt} from \eqref{E:12.63nwwt}. Note that the boundedness of $f_n(x)$ ensures that the right-hand side of \eqref{E:12.63nwwt} is of order $\lambda_n n^{-4}\log(n^{4}/\lambda_n)$ which tends to zero by \eqref{E:12.23uau}. In particular, the probability on the left of \eqref{E:12.63nwwt} tends to one. If~$a_n$ grows so fast that
\begin{equation}
\texte^{a_n}\lambda_n n^{-4}\log(n^{4}/\lambda_n)\,\underset{n\to\infty}\longrightarrow\,\infty,
\end{equation}
then \eqref{E:12.56nwwt} and \eqref{E:12.58nwwt} equate the leading order asymptotic of the quantities on the left-hand side of \eqref{E:12.54nwwt} and \eqref{E:12.63nwwt}. The claim follows.
\end{proofsect}

We are now ready to give:

\begin{proofsect}{Proof of Theorem~\ref{thm-Z-unique}}
Let~$\mathfrak D_0\subseteq\mathfrak D$ be a countable set containing the domain of interest along with the collection of all open equilateral triangles of side length~$2^{-n}$, $n\in\Z$, with (at least) two vertices at points with rational coordinates, and all finite disjoint unions thereof. Let $\{N_k\colon k\ge1\}$ be a subsequence along which $\eta_{N,r_N}^D$ converges in law to some~$\eta^D$ for all $D\in\mathfrak D_0$. We will never need the third coordinate of $\eta^D_{N,r_N}$ so we will henceforth think of these measures as two-coordinate processes only.

Let $Z^{D}$, for each $D\in\mathfrak D_0$, denote the random measure associated with the limit process~$\eta^{D}$. Pick~$f\in C(\overline D)$ and assume that~$f$ is positive and Lipschitz. As we follow closely the proof of Theorem~\ref{thm-10.16}, let us review it geometric setting: Given a $K\in\N$, which we will assume to be a power of two, and~$\delta\in(0,1/100)$ and let~$\{T^i\colon i=1,\dots,n_K\}$ be those triangles in the tiling of~$\R^2$ by triangles of side length~$K^{-1}$ that fit into~$D^\delta$. Denote by~$x_i$ the center of~$T^i$ and let $T^i_\delta$ be the triangle of side length~$(1-\delta)K^{-1}$ centered at~$x_i$ and oriented the same way as~$T^i$. Write~$\wt D_\delta:=\bigcup_{i=1}^{n_K}T^i_\delta$ and let $\chi_\delta\colon D\to[0,1]$ be a function that is one on~$\wt D_{2\delta}$, zero on~$D\smallsetminus \wt D_\delta$ and is $(2K\delta^{-1})$-Lipschitz elsewhere. 

Set~$f_\delta:=f\chi_\delta$. In light of Lemma~\ref{lemma-10.4} (proved above) and Exercise~\ref{ex:10.5}, we then have \eqref{E:10.51uai} and so it suffices to work with~$\langle Z^D,f_\delta\rangle$. First we ask the reader to solve:

\begin{myexercise}
Prove that the measures $\{Z^{\wt D}\colon \wt D\in\mathfrak D_0\}$ satisfy the Gibbs-Markov property in the form \eqref{E:10.47} (for~$M^D$ replaced by~$Z^D$).
\end{myexercise}

\noindent
Moving along the proof of Theorem~\ref{thm-10.16}, next we observe that \eqref{E:10.53ua} and Proposition~\ref{prop-10.18} apply because they do not use any particulars of the measures~$M^D$ there. Letting $c>0$ is the constant from \eqref{E:10.53ua} and writing
\begin{equation}
\wt A^i_{K,R}:=\bigl\{\text{\rm osc}_{T^i_\delta}\Phi^{D,\wt D}\le R\bigr\}\cap\bigl\{\Phi^{D,\wt D}(x_i)\le2\sqrt g\log K-c\log\log K\bigr\}.
\end{equation}
for a substitute of event $A^i_{K,R}$ from \eqref{E:10.52nwwt}, we thus have
\begin{multline}
\label{E:12.66nwwt}
\langle Z^D,f_\delta\rangle =o(1)
\\
+\sum_{i=1}^{n_K}1_{\wt A^i_{K,R}}\texte^{\alpha\Phi^{D,\wt D}(x_i)}\int_{T^i}Z^{T^i}(\textd x)\,\texte^{\alpha[\Phi^{D,\wt D}(x)-\Phi^{D,\wt D}(x_i)]}\,f_\delta(x),
\end{multline}
where $o(1)\to0$ in probability as $K\to\infty$ and $R\to\infty$ and where the family of measures $\{Z^{T^i}\colon i=1,\dots,n_K\}$ are independent on one another (and equidistributed, modulo shifts) and independent of $\Phi^{D,\wt D}$. 

Proceeding along the argument from the proof of Theorem~\ref{thm-10.16}, we now wish to compute the negative exponential moment of the sum on the right of \eqref{E:12.66nwwt}, conditional on~$\Phi^{D,\wt D}$. The indicator of $\wt A^i_{K,R}$ enforces that
\begin{equation}
\lambda_K:=\texte^{\alpha\Phi^{D,\wt D}(x_i)}\quad\text{obeys}\quad K^4\lambda_K\le (\log K)^{-\alpha c}
\end{equation}
while, in light of  the regularity assumptions on~$f$ and~$\chi_\delta$ and the harmonicity and boundedness of oscillation of~$\Phi^{D,\wt D}$,
\begin{equation}
f_K(x):=f(x)\chi_\delta(x)\,\texte^{\alpha[\Phi^{D,\wt D}(x)-\Phi^{D,\wt D}(x_i)]}
\end{equation}
is $R'K$-Lipschitz for some~$R'$ depending only on~$R$. Lemma~\ref{lemma-12.7ua} (which is our way to by-pass Proposition~\ref{prop-3}) then yields the analogue of \eqref{E:10.58ui}: For any $i=1,\dots,n_K$ where $\wt A^i_{K,R}$ occurs,
\begin{equation}
\begin{aligned}
E\biggl(&\,\exp\Bigl\{-\texte^{\alpha\Phi^{D,\wt D}(x_i)}Z^{T^i}(f_K)\Bigr\}\,\bigg|\,\Phi^{D,\wt D}\biggr)
\\
&=\exp\biggl\{\log\bigl(K^{-4}\texte^{\alpha\Phi^{D,\wt D}(x_i)}\bigr)\,\bar c\int_{T^i}\textd x\,f_\delta(x)\,\texte^{\alpha\Phi^{D,\wt D}(x)}\,r_{T^i}(x)^2\biggr\}
\\
&\qquad\qquad\times\exp\Bigl\{\wt\epsilon_K K^{-4}\texte^{\alpha\Phi^{D,\wt D}(x_i)}\log\bigl(K^{-4}\texte^{\alpha\Phi^{D,\wt D}(x_i)}\bigr)\Bigr\},
\end{aligned}
\end{equation}
where $\wt\epsilon_K$ is a random ($\Phi^{D,\wt D}$-measurable) quantity taking values in some $[\epsilon_K,\epsilon_K]$ for some deterministic $\epsilon_K$ with $\epsilon_K\to0$ as $K\to\infty$.

Assuming (at the cost of a routine limit argument at the end of the proof) that~$f\ge\delta$ on~$D$, the bound on the oscillation of $\Phi^{D,\wt D}$ permits us to absorb the error term into an additive modification of~$f$ by a term of order $\epsilon_K\texte^{\alpha R}$. Invoking \twoeqref{E:10.59ui}{E:10.60nwwt}, we conclude that the measure
\begin{multline}
\quad\alpha \bar c\, r_D(x)^2\sum_{i=1}^{n_K}\1_{\wt A_{K,R}^i}\bigl(\alpha\Var(\Phi^{D,\wt D}(x))-\Phi^{D,\wt D}(x)\bigr)
\\
\times\texte^{\alpha \Phi^{D,\wt D}(x)-\frac12\alpha^2\Var(\Phi^{D,\wt D}(x))}\,\1_{T_\delta^i}(x)\,\textd x
\quad
\end{multline}
tends to $Z^D$ in the limit $K\to\infty$, $R\to\infty$ and $\delta\downarrow0$. In particular,  the law of~$Z^D$ is independent of the subsequence $\{N_k\colon k\ge1\}$ used to define it and the processes $\{\eta^D_{N,r_N}\colon N\ge1\}$ thus converge in law.
\end{proofsect}

Theorem~\ref{thm-Z-unique} implies, via Lemma~\ref{lemma-10.1ua}, the convergence of the centered maximum from Theorem~\ref{thm-BDZ} and, by Lemma~\ref{lemma-10.8}, also the joint convergence of maxima in disjoint open subsets as stated in Lemma~\ref{lemma-mult-max}. The convergence statement of the full extremal process in Theorem~\ref{thm-extremal-vals} is now finally justified.

\nopagebreak
\begin{figure}[t]
\vglue-1mm
\centerline{\includegraphics[width=0.99\textwidth]{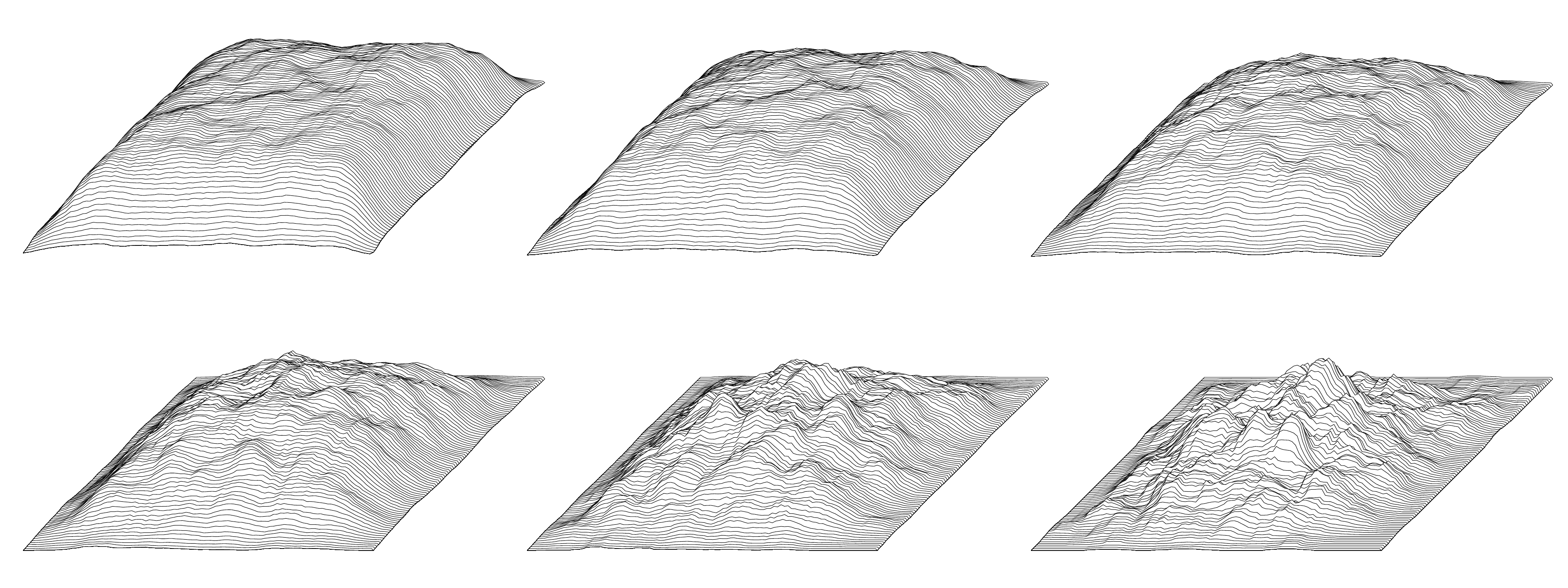}
}

\vglue0mm
\begin{quote}
\small 
\caption{
\label{fig-maxima-plot}
\small
Empirical plots of~$x\mapsto\rho^D(x,t)$ obtained from a set of about 100000 samples of the maximum of the DGFF on a 100$\times$100 square. The plots (labeled left to right starting with the top row) correspond to~$t$ increasing by uniform amounts over an interval of length 3 with~$t$ in the fourth figure set to the empirical mean.}
\normalsize
\end{quote}
\end{figure}

\section{The local limit theorem}
\noindent
Our way of control of the absolute maximum by conditioning on its position can be used to give also a local limit theorem for both the position and value of the maximum:

\begin{mytheorem}[Local limit theorem for absolute maximum]
\label{thm-LLT}
Using the setting of Theorem~\ref{thm-extremal-vals}, for any $x\in D$ and any $a<b$,
\begin{equation}
\lim_{N\to\infty}
N^2 P\Bigl(\operatornamewithlimits{argmax}_{D_N}h^{D_N}=\lfloor xN\rfloor,\,\max_{x\in D_N}h^{D_N}_x-m_N\in(a,b)\Bigr)=\int_a^b\rho^D(x,t)\textd t\,,
\end{equation}
where $x\mapsto\rho^D(x,t)$ is the Radon-Nikodym derivative of the measure
\begin{equation}
\label{E:density}
A\mapsto\texte^{-\alpha t}E\bigl(Z^D(A)\texte^{-\alpha^{-1}\texte^{-\alpha t}Z^D(D)}\bigr)
\end{equation}
with respect to the Lebesgue measure on~$D$.
\end{mytheorem}

\noindent
We only state the result here and refer the reader to~\cite{BL3} for details of the proof. Nothing explicit is known about~$\rho^D$ except its~$t\to\infty$ asymptotics; see \eqref{E:rhoD-asymptotic}.
Empirical plots of~$\rho^D$ for different values of~$t$ are shown in Fig.~\ref{fig-maxima-plot}.


\chapter{Random walk in DGFF landscape}
\label{lec-13}\noindent
In this lecture we turn our attention to a rather different problem than discussed so far: a random walk in a random environment. The connection with the main theme of these lectures is through the specific choice of the random walk dynamics which we interpret as the motion of a random particle in a DGFF landscape. We first state the results, obtained in a recent joint work with Jian Ding and Subhajit Goswami~\cite{BDG}, on the behavior of this random walk. Then we proceed to develop the key method of the proof, which is based on recasting the problem as a random walk among random conductances and applying methods of electrostatic theory.

\section{A charged particle in an electric field}
\noindent
Let us start with some physics motivation.
Suppose we are given the task to describe the motion of a charged particle in a rapidly varying electric field. A natural choice is to fit this into the framework of the theory of random walks in random environment (RWRE) as follows. The particle is confined to the hypercubic lattice~$\Z^d$ and the electric field is represented by a configuration $h=\{h_x\colon x\in\Z^d\}$ with~$h_x$ denoting the \emph{electrostatic potential} at~$x$. Given a realization of~$h$, the charged particle then performs a ``random walk'' which, technically, is a discrete-time Markov chain on~$\Z^d$ with the transition probabilities
\begin{equation}
\label{eq-def-transition}
\cmss P_h(x,y):=\frac{\texte^{\beta(h_y-h_x)}}{\displaystyle\sum_{z\colon(x,z)\in E(\Z^d)}\texte^{\beta(h_z-h_x)}} \,1_{(x,y)\in E(\Z^d)},
\end{equation}
where~$\beta$ is a parameter playing the role of the \emph{inverse temperature}. We will assume $\beta>0$ which means that the walk is  more likely to move in the direction where the electrostatic potential increases. Alternatively, we may think of~$\beta$ as the value of the charge of the moving particle.

Let~$X=\{X_n\colon n\ge0\}$ denote a sample path of the Markov chain. We will write~$P_h^x$ for the law of~$X$ with~$P_h^x(X_0=x)=1$, use~$E^x_h$ to denote expectation with respect to~$P_h^x$ and write~$\BbbP$ to denote the law of the DGFF on~$\Z^2\smallsetminus\{0\}$. As usual in RWRE theory, we will require that 
\begin{equation}
\label{E:12.2a}
\{\cmss P_h(x,\cdot)\colon x\in\Z^d\}\text{ is stationary, ergodic under the shifts of~$\Z^d$}
\end{equation}
as that is typically the minimal condition needed to extract a limit description of the path properties.
However, since~$\cmss P_h(x,\cdot)$ depends only on the differences of the field, in our case \eqref{E:12.2a} boils down to the requirement:
\begin{equation}
\label{E:12.2}
\bigl\{h_x-h_y\colon (x,y)\in E(\Z^d)\bigr\}\text{ is stationary, ergodic under the shifts of~$\Z^d$}\,.
\end{equation}
A number of natural examples may be considered, with any i.i.d.\ random field or, in fact, any stationary and ergodic random field obviously satisfying \eqref{E:12.2}.  However, our desire in the lectures is to work with the fields that exhibit  logarithmic correlations. A prime example of such a field is the two-dimensional~DGFF. 

\nopagebreak
\begin{figure}[t]
\vglue-1mm
\centerline{\includegraphics[width=0.27\textwidth]{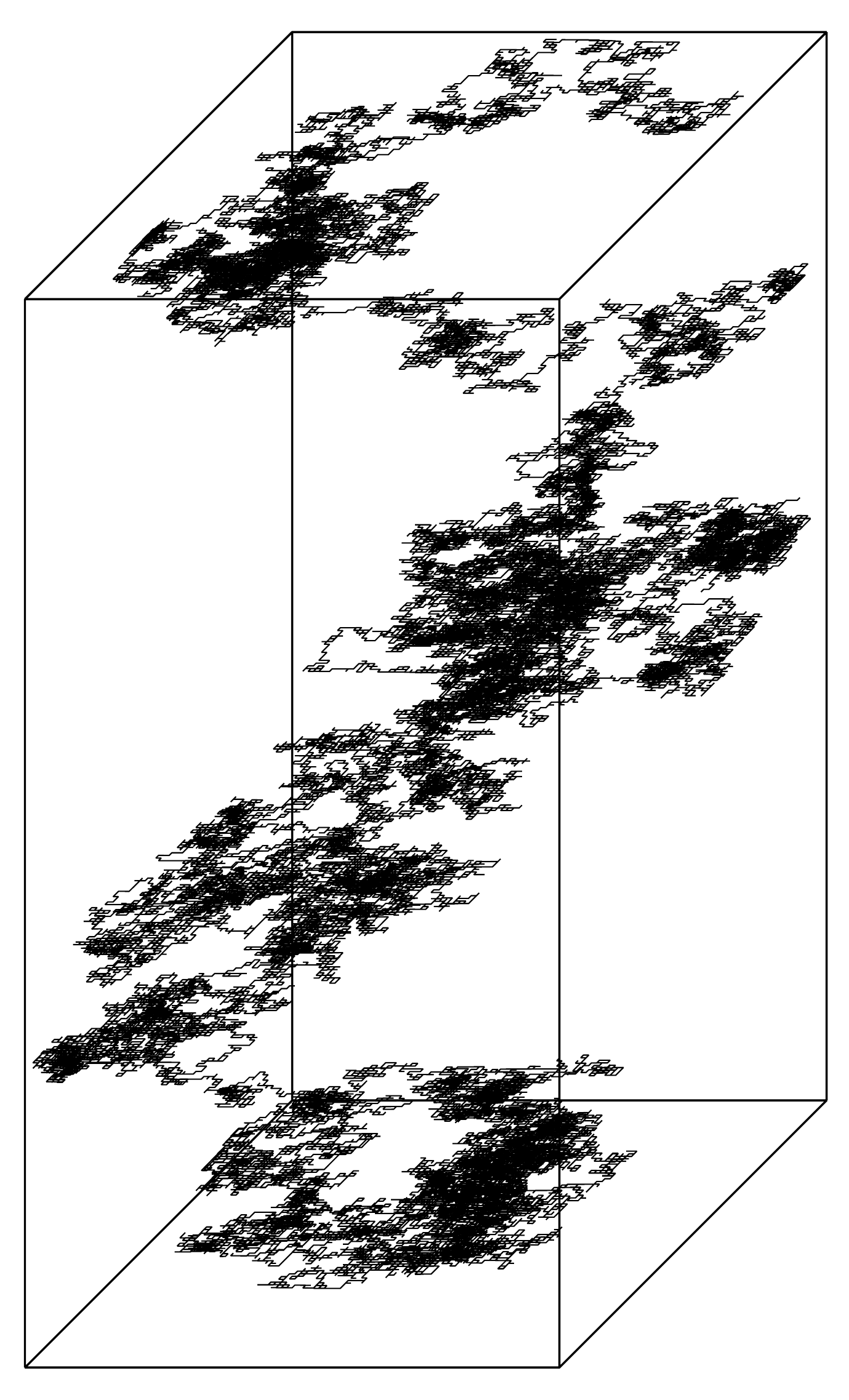}
\includegraphics[width=0.27\textwidth]{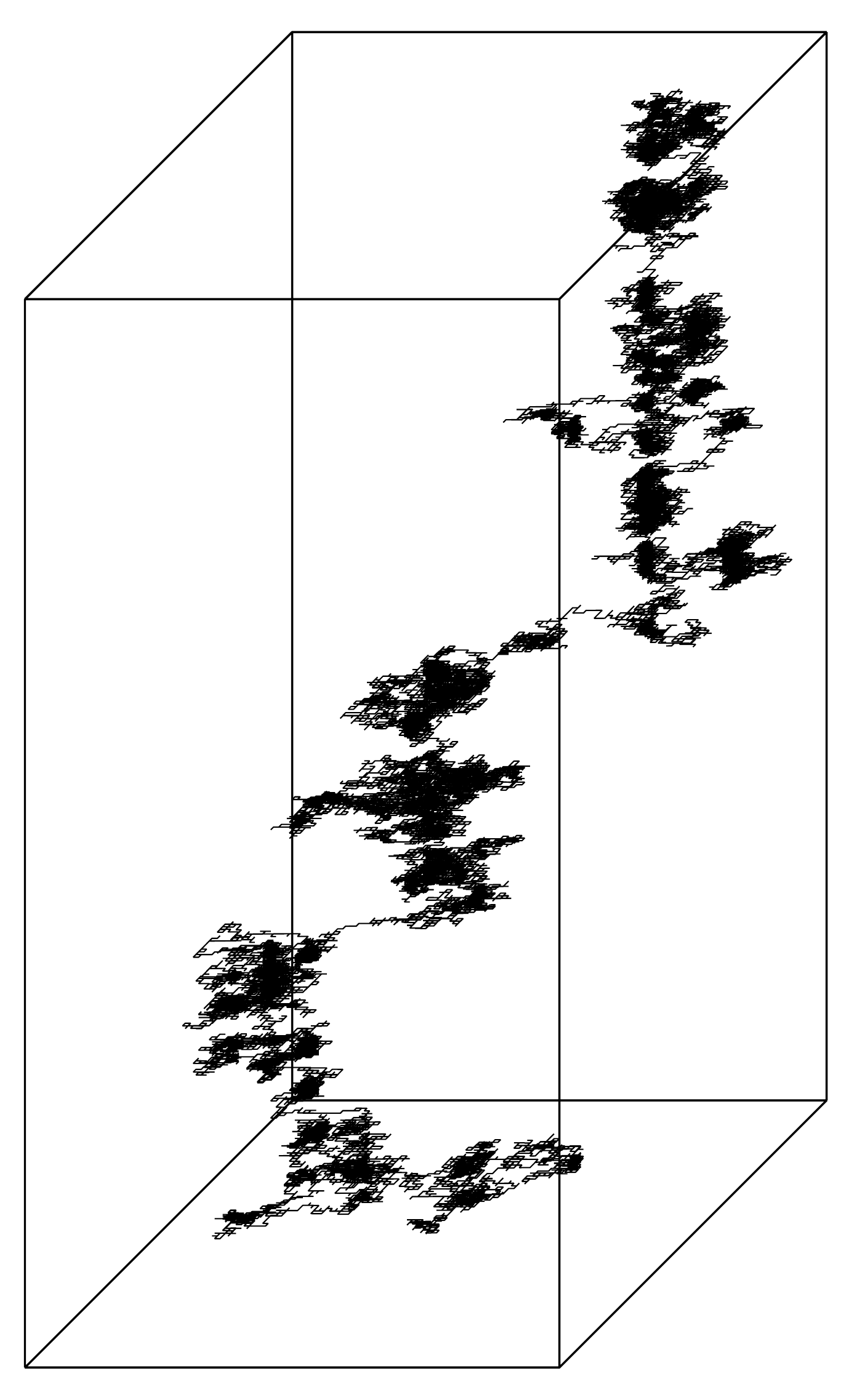}
\includegraphics[width=0.27\textwidth]{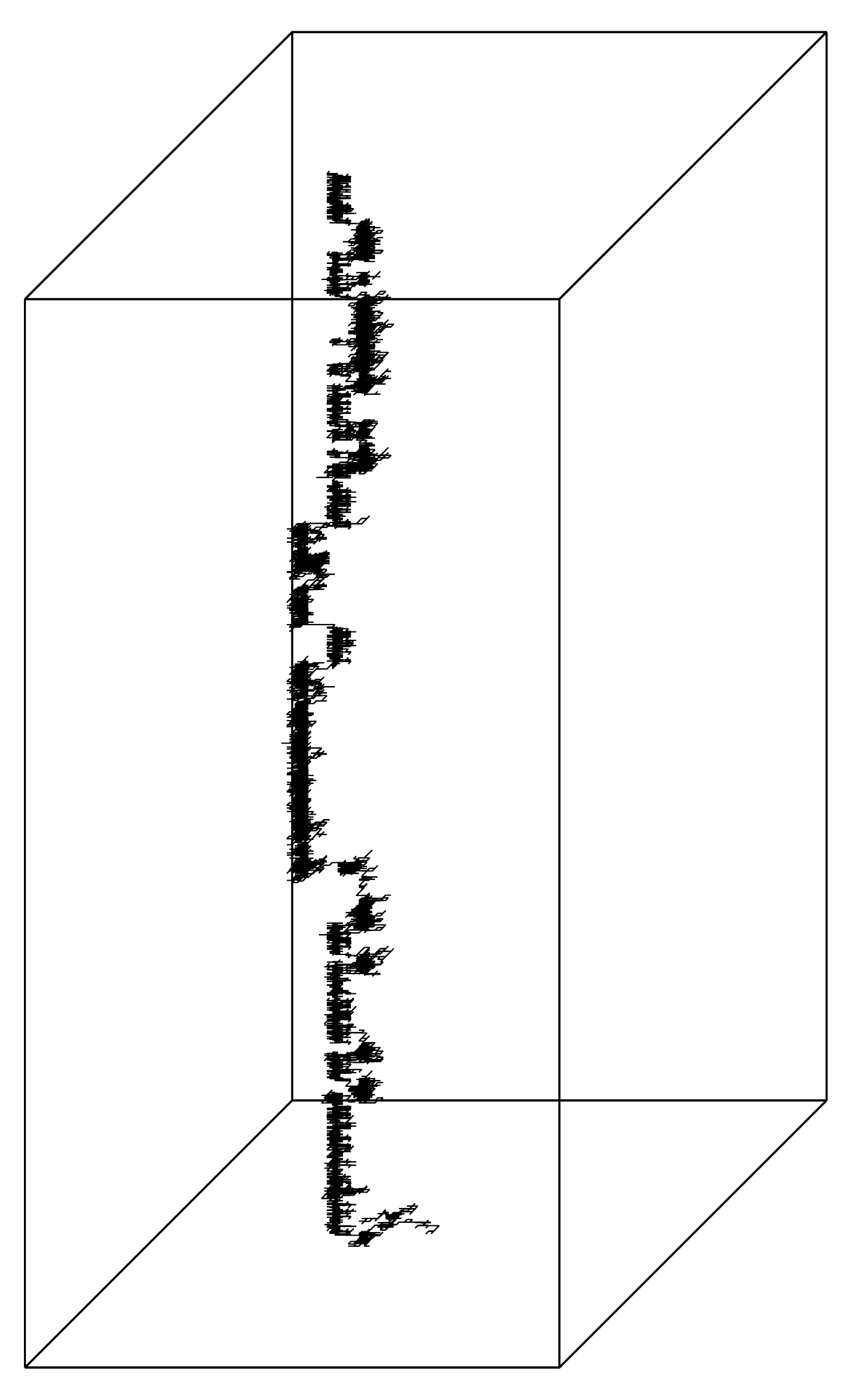}
}

\vglue0mm
\begin{quote}
\small 
\caption{
\label{fig-walks}
\small
Runs of $100000$ steps of the random walk with transition probabilities~\eqref{eq-def-transition} and~$\beta$ equal to~$0.2$, $0.6$ and~$1.2$ multiples of~$\tilde\beta_\cc$. Time runs upwards along the vertical axis. Trapping effects are quite apparent.}
\normalsize
\end{quote}
\end{figure}

The motivation for our focus on log-correlated fields comes from the 2004 paper of Carpentier and Le Doussal~\cite{Carpentier-LeDoussal}, who discovered, on the basis of physics arguments, that such environments exhibit the following phenomena:
\settowidth{\leftmargini}{(1111)}
\begin{enumerate}
\item[(1)] trapping effects make the walk behave subdiffusively with the diffusive exponent~$\nu$, defined via $|X_n|= n^{\nu+o(1)}$, depending non-trivially on~$\beta$, and
\item[(2)] $\beta\mapsto\nu(\beta)$ undergoes a phase transition (i.e., a change in analytic dependence) as~$\beta$ varies through a critical point~$\tilde\beta_\cc$.
\end{enumerate}
Log-correlated fields are in fact deemed critical for the above phase transition to occur (with weaker correlations generically corresponding to~$\beta\downarrow0$ regime and stronger correlations to~$\beta\to\infty$) although we will not try to make this precise. The  purpose of these lectures it to demonstrate that (1-2) indeed happen in at least one example; namely, the two-dimensional DGFF.

We will thus henceforth focus on $d=2$. To get around the fact that the DGFF on~$\Z^2$ does not exist, we  will  work with
\begin{equation}
\label{E:12.4}
h:=\text{ DGFF in~$\Z^2\smallsetminus\{0\}$}
\end{equation}
as the driving ``electric'' field for the rest of these lectures.
This does fall into the class of systems introduced above; indeed, we have:

\begin{myexercise}[Pinned DGFF has stationary gradients]
\label{ex:13.1}
Show that $h$ from \eqref{E:12.4} obeys \eqref{E:12.2}.
\end{myexercise}

\section{Statement of main results}
\noindent
Having elucidated the overall context of the problem at hand, we are now ready to describe the results for the particular choice \eqref{E:12.4}. These have all been proved in a joint 2016 paper with J.~Ding and S.~Goswami~\cite{BDG}. Our first result concerns the decay of return probabilities (a.k.a.\ the \emph{heat kernel}):

\begin{mytheorem}[Heat-kernel decay]
\label{thm-heatkernel}
For each~$\beta>0$ and each~$\delta>0$,
\begin{equation}
\label{E:12.5a}
\BbbP\biggl(\,\frac1T\texte^{-(\log T)^{1/2+\delta}}\le P^0_h(X_{2T}=0)\le\frac1T\texte^{(\log T)^{1/2+\delta}}\biggr)
\,\,\underset{T\to\infty}\longrightarrow\,\,1.
\end{equation}
\end{mytheorem}

A noteworthy point is that the statement features no explicit dependence on~$\beta$. (In fact, it applies even to~$\beta=0$ when~$X$ is just the simple random walk on~$\Z^2$.) Hence, as far as the leading order of the return probabilities  is concerned, the walk behaves just as the simple random walk. Notwithstanding, the $\texte^{\pm(\log T)^{1/2+\delta}}$ terms are too large to let us decide whether~$X$ is recurrent or transient, a question that we will eventually resolve as well albeit by different~means.

Although the propensity of the walk to move towards larger values of the field does not seem to affect the (leading order) heat kernel decay, the effect on the path properties is more detectable. For each set~$A\subset\Z^2$, define
\begin{equation}
\tau_A:=\inf\bigl\{n\ge0\colon X_n\in A\bigr\}.
\end{equation}
Denote also
\begin{equation}
\label{E:BN}
B(N):=[-N,N]^2\cap\Z^2.
\end{equation}
Then we have:

\begin{mytheorem}[Subdiffusive expected exit time]
\label{thm-exittime}
For each~$\beta>0$ and each~$\delta>0$,
\begin{equation}
\label{E:12.8}
\BbbP\biggl(\,N^{\theta(\beta)}\texte^{-(\log N)^{1/2+\delta}}\le E^0_h(\tau_{B(N)^\cc})\le N^{\theta(\beta)}\texte^{(\log N)^{1/2+\delta}}\biggr)
\,\,\underset{N\to\infty}\longrightarrow\,\,1\,,
\end{equation}
where, for~$\tilde\beta_\cc:=\sqrt{\pi/2}$\,,
\begin{equation}
\label{E:12.9}
\theta(\beta):=\begin{cases}
2+2(\beta/\tilde\beta_\cc)^2,\qquad&\text{if }\beta\le\tilde\beta_\cc,
\\
4\beta/\tilde\beta_\cc,\qquad&\text{if }\beta\ge\tilde\beta_\cc.
\end{cases}
\end{equation}
\end{mytheorem}

The functional dependence of~$\theta(\beta)$ on~$\beta$ actually takes a rather familiar form: Denoting $\lambda:=\beta/\tilde\beta_\cc$, for~$\lambda\in(0,1)$ we have~$\theta(\beta)=2+2\lambda^2$, which is the scaling exponent associated with the intermediate level set at height~$\lambda$-multiple of the absolute maximum. The dependence on~$\beta$ changes at~$\tilde\beta_\cc$, just as predicted by Carpentier and Le Doussal~\cite{Carpentier-LeDoussal}. Moreover, for~$\beta>0$, we have~$\theta(\beta)>0$. The walk thus takes considerably longer (in expectation) to exit a box than the simple random walk. This can be interpreted as a version of \emph{subdiffusive behavior}. 

The standard definition of subdiffusive behavior is via the typical spatial scale of the walk at large times. Here we can report only a one-way bound:

\begin{mycorollary}[Subdiffusive lower bound]
\label{cor-12.4}
For all~$\beta>0$ and all~$\delta>0$,
\begin{equation}
P_h^0\Bigl(|X_T|\ge T^{1/\theta(\beta)}\texte^{-(\log N)^{1/2+\delta}}\Bigr)
\,\,\underset{N\to\infty}\longrightarrow\,\,1\,,\quad\text{\rm in $\BbbP$-probability}.
\end{equation}
\end{mycorollary}

\noindent
Unfortunately, the more relevant upper bound is elusive at this point (although we believe that our methods can be boosted to include a matching leading-order upper bound as well).

Our method of proof of the above results relies on the following elementary observation: The transition probabilities from \eqref{eq-def-transition} can be recast as
\begin{equation}
\label{E:12.11ua}
\cmss P_h(x,y):=\frac{\texte^{\beta(h_y+h_x)}}{\pi_h(x)}\,1_{(x,y)\in E(\Z^d)}\,,
\end{equation}
where
\begin{equation}
\label{E:12.12}
\pi_h(x):=\sum_{z\colon(x,z)\in E(\Z^d)}\texte^{\beta(h_z+h_x)}.
\end{equation}
This, as we will explain in the next section, phrases the problem as a \emph{random walk among random conductances}, with the conductance of edge $(x,y)\in E(\Z^d)$ given by
\begin{equation}
\label{E:12.13}
c(x,y):=\texte^{\beta(h_y+h_x)}.
\end{equation}
As is readily checked,~$\pi_h$ is a reversible, and thus stationary, measure for~$X$.  

Over the past couple of decades, much effort went into the understanding of random walks among random conductances (see, e.g., a review by the present author~\cite{Biskup-RCM} or Kumagai~\cite{Kumagai-RCM}). A standing assumption of these is that the conductance configuration, 
\begin{equation}
\bigl\{c(x,y)\colon (x,y)\in E(\Z^d)\bigr\}\,,
\end{equation}
is stationary and ergodic with respect to the shifts of~$\Z^d$. For~$\Z^d$ in~$d\ge3$, the full-lattice DGFF  (with zero boundary conditions) is stationary and so the problem falls under this umbrella. The results of Andres, Deuschel and Slowik~\cite{ADS} then imply scaling of the random walk to a non-degenerate Brownian motion. 

When~$d=2$ and~$h=$ DGFF on~$\Z^2\smallsetminus\{0\}$ the conductances \eqref{E:12.13} are no longer stationary. Thus we have a choice to make: either work with stationary transition probabilities with no help from reversibility techniques or give up on stationarity and earn reversibility. (A similar situation occurs for Sinai's RWRE on~$\Z$~\cite{Sinai}.) It is the latter option that has (so far) led to results.

The prime benefit of reversibility is that it makes the Markov chain amenable to analysis via the methods of \emph{electrostatic theory}. This is a standard technique; cf Doyle and Snell~\cite{Doyle-Snell} or Lyons and Peres~\cite{Lyons-Peres}. We thus interpret the underlying graph as an electric network with resistance $r(x,y):=1/c(x,y)$ assigned to edge~$(x,y)$. The key notion to consider is the \emph{effective resistance} $R_\eff(0,B(N)^\cc)$ from~$0$ to~$B(N)^\cc$. We will define this quantity precisely in the next section; for the knowledgeable reader we just recall that $R_\eff(0,B(N)^\cc)$ is the voltage difference one needs to put between~$0$ and~$B(N)^\cc$ to induce unit (net) current through the network. Using sophisticated techniques, for the effective resistance we then get:

\begin{mytheorem}[Effective resistance growth]
\label{thm-12.5}
For each~$\beta>0$,
\begin{equation}
\label{E:1.11ua}
\limsup_{N\to\infty}\,\,\frac{\log R_\eff(0,B(N)^\cc)}{(\log N)^{1/2}(\log\log N)^{1/2}}<\infty,
\qquad\BbbP\text{\rm-a.s.}
\end{equation}
and
\begin{equation}
\label{E:1.12ua}
\liminf_{N\to\infty}\,\,\frac{\log R_\eff(0,B(N)^\cc)}{(\log N)^{1/2}/(\log\log\log N)^{1/2}}>0,
\qquad\BbbP\text{\rm-a.s.}
\end{equation}
\end{mytheorem}

Both conclusions of Theorem~\ref{thm-12.5} may be condensed into one (albeit weaker) statement as
\begin{equation}
R_\eff(0,B(N)^\cc) = \texte^{(\log N)^{1/2+o(1)}},\quad N\to\infty.
\end{equation}
In particular, $R_\eff(0,B(N)^\cc)\to\infty$ as~$N\to\infty$, a.s. The standard criteria of recurrence and transience of Markov chains (to be discussed later) then yield:

\begin{mycorollary}
For $\BbbP$-a.e.\ realization of~$h$, the Markov chain~$X$ is recurrent.
\end{mycorollary}

\noindent
With the help of Exercise~\ref{ex:13.1} we get:

\begin{myexercise}
Show that the limits in \twoeqref{E:1.11ua}{E:1.12ua} depend only on the gradients of the pinned DGFF and are thus constant almost surely.
\end{myexercise}

The remainder of these lectures will be spent on describing the techniques underpinning the above results. Several difficult proofs will only be sketched or outright omitted; the aim is to communicate the ideas rather then give a full fledged account that the reader may just as well get by reading the original~paper.

\section{A crash course on electrostatic theory}
\noindent
We begin by a review of the connection between Markov chains and electrostatic theory.
Consider a finite, unoriented, connected graph~$\mathfrak G=(\cmss V,\cmss E)$ where both orientations of edge~$e\in\cmss E$ are identified as one. An assignment of \emph{resistance} $r_e\in(0,\infty)$ to each edge~$e\in\cmss E$ then makes~$\mathfrak G$ an electric network. An alternative description uses \emph{conductances}~$\{c_e\colon e\in \cmss E\}$ where
\begin{equation}
\label{E:12.17}
c_e:=\frac1{r_e}.
\end{equation}
We will exchangeably write $r(x,y)$ for~$r_e$ when~$e=(x,y)$, and similarly for~$c(x,y)$. Note that these are symmetric quantities, $r(x,y)=r(y,x)$ and $c(x,y)=c(y,x)$ whenever~$(x,y)=(y,x)\in\cmss E$.

Next we define some key notions of the theory. For any two distinct $u,v\in \cmss V$, let
\begin{equation}
\FF(u,v):=\bigl\{f \text{ function } \cmss V\to\R\colon\, f(u)=1,\,f(v)=0\bigr\}.
\end{equation}
We interpret~$f(x)$ as an assignment of a \emph{potential} to vertex~$x\in\cmss V$; each $f\in\FF(u,v)$ then has unit potential difference (a.k.a.~\textit{voltage}) between~$u$ and~$v$.
For any potential~$f\colon \cmss V\to\R$, define its Dirichlet energy by
\begin{equation}
\EE(f):=\sum_{(x,y)\in \cmss E}c(x,y)\bigl[f(y)-f(x)\bigr]^2,
\end{equation}
where, by our convention, each edge contributes only once. 

\begin{mydefinition}[Effective conductance] The infimum
\begin{equation}
\label{E:12.20}
C_\eff(u,v):=\inf\bigl\{\EE(f)\colon f\in\FF(u,v)\bigr\}
\end{equation}
is the \emph{effective conductance from~$u$ to~$v$}.
\end{mydefinition}

\noindent
Note that $C_\eff(u,v)>0$ since~$\mathfrak G$ is assumed connected and finite and the conductances are assumed to be strictly positive. 

Next we define the notion of (electric) current as follows:

\begin{mydefinition}[Current]
Let~$\vec{\cmss E}$ denote the set of oriented edges in~$\mathfrak G$, with both orientations included. A \emph{current from~$u$ to~$v$} is an assignment $i(e)$ of a real number to each~$e\in\vec{\cmss E}$ such that, writing $i(x,y)$ for~$i(e)$ with~$e=(x,y)$,
\begin{equation}
\label{E:12.21}
i(x,y)=-i(y,x),\quad (x,y)\in\vec{\cmss E}
\end{equation}
and
\begin{equation}
\label{E:12.22}
\sum_{y\colon (x,y)\in\vec{\cmss E}}i(x,y)=0,\quad x\in\cmss V\smallsetminus\{u,v\}.
\end{equation}
\end{mydefinition}

The first condition reflects on the fact that the current flowing along $(x,y)$ is the opposite of the current flowing along~$(y,x)$. The second condition then forces that the current be conserved at all vertices except~$u$ and~$v$. Next we observe:

\begin{mylemma}[Value of current]
For each current~$i$ from~$u$ to~$v$,
\begin{equation}
\label{E:12.23}
\sum_{x\colon (u,x)\in\vec{\cmss E}}i(u,x) = \sum_{x\colon (x,v)\in\vec{\cmss E}}i(x,v)\,.
\end{equation}
\end{mylemma}

\begin{proofsect}{Proof}
Conditions \twoeqref{E:12.21}{E:12.22} imply
\begin{equation}
\begin{aligned}
0=\sum_{(x,y)\in\vec{\cmss E}}i(x,y)&=\sum_{x\in\cmss V}\,\sum_{y\colon(x,y)\in\vec{\cmss E}}i(x,y)
\\
&=\sum_{y\colon(u,y)\in\vec{\cmss E}}i(u,y)+\sum_{y\colon(v,y)\in\vec{\cmss E}}i(v,y).
\end{aligned}
\end{equation}
Employing \eqref{E:12.21} one more time, we then get \eqref{E:12.23}.
\end{proofsect}

A natural interpretation of \eqref{E:12.23} is that the current incoming to the network at~$u$ equals the outgoing current at~$v$. (Note that this may be false in infinite networks.) We call the common value in \eqref{E:12.23} the \emph{value of current~$i$}, with the notation $\text{val}(i)$. It is natural to single out the currents with unit value into
\begin{equation}
\II(u,v):=\bigl\{i\colon \text{ current from~$u$ to~$v$ with }\text{val}(i)=1\bigr\}\,.
\end{equation}
For each current~$i$, its Dirichlet energy is then given by
\begin{equation}
\wt\EE(i):=\sum_{e\in \cmss E}r_e i(e)^2\,,
\end{equation}
where we again note that each edge contributes only one term to the sum.

\begin{mydefinition}[Effective resistance]
The infimum
\begin{equation}
\label{E:12.27}
R_\eff(u,v):=\inf\bigl\{\wt\EE(i)\colon i\in\II(u,v)\bigr\}
\end{equation}
is the \emph{effective resistance from~$u$ to~$v$}.
\end{mydefinition}

\noindent
Note that $\II(u,v)\ne\emptyset$ and thus also $R_\eff(u,v)<\infty$ thanks to the assumed finiteness and connectivity of~$\mathfrak G$. 

It is quite clear that the effective resistance and effective conductance must somehow be closely related. For instance, by \eqref{E:12.17} they are the reciprocals of each other in the network with two vertices and one edge. To address this connection in general networks, we first observe:

\begin{mylemma}
\label{lemma-12.11}
For any two distinct $u,v\in\cmss V$,
\begin{equation}
\label{E:12.28}
\EE(f)\wt\EE(i)\ge1,\quad f\in\FF(u,v),\,i\in\II(u,v).
\end{equation}
In particular, $R_\eff(u,v)C_\eff(u,v)\ge1$.
\end{mylemma}

\begin{proofsect}{Proof}
Let $f\in\FF(u,v)$ and $i\in\II(u,v)$. By a symmetrization argument and the definition of unit current,
\begin{multline}
\label{E:12.29}
\quad
\sum_{(x,y)\in\cmss E} i(x,y)\bigl[f(x)-f(y)\bigr]
=\frac12\sum_{x\in\cmss V}\,\sum_{y\colon (x,y)\in\cmss E} i(x,y)\bigl[f(x)-f(y)\bigr]
\\
=\sum_{x\in\cmss V}f(x)\sum_{y\colon (x,y)\in\cmss E}i(x,y)=f(u)-f(v)=1.
\quad
\end{multline}
On the other hand, \eqref{E:12.17} and the Cauchy-Schwarz inequality yield
\begin{multline}
\qquad
\sum_{(x,y)\in\cmss E} i(x,y)\bigl[f(x)-f(y)\bigr]
\\
=\sum_{(x,y)\in\cmss E} \sqrt{r(x,y)}\,i(x,y)\,\sqrt{c(x,y)}\,\bigl[f(x)-f(y)\bigr]
\\
\le\wt\EE(i)^{1/2}\EE(f)^{1/2}.
\qquad
\end{multline}
This gives \eqref{E:12.28}. The second part follows by optimizing over~$f$ and~$i$.
\end{proofsect}

We now claim:

\begin{mytheorem}[Electrostatic duality]
For any distinct $u,v\in\cmss V$,
\begin{equation}
C_\eff(u,v)=\frac1{R_\eff(u,v)}\,.
\end{equation}
\end{mytheorem}

\begin{proofsect}{Proof}
Since~$\II(u,v)$ can be identified with a closed convex subset of~$\R^{\cmss E}$ and $i\mapsto\wt\EE(i)$ with a strictly convex function on $\R^{\cmss E}$ that has compact level sets, there is a unique minimizer~$i_\star$ of~\eqref{E:12.27}. We claim that~$i_\star$ obeys the \emph{Kirchhoff cycle law}: For each $n\ge1$ and each~$x_0,x_1,\dots,x_n=x_0\in\cmss V$ with $(x_i,x_{i+1})\in\cmss E$ for each~$i=0,\dots,n-1$,
\begin{equation}
\label{E:12.32}
\sum_{k=1}^n r(x_k,x_{k+1})i_\star(x_k,x_{k+1})=0.
\end{equation}
To show this, let~$j$ be defined by $j(x_k,x_{k+1})=-j(x_{k+1},x_k):=1$ for~$k=1,\dots,n$ and $j(x,y):=0$ on all edges not belonging to the cycle $(x_0,\dots,x_n)$. Then $i_\star+a j\in\II(u,v)$ for any~$a\in\R$ and so, since~$i_\star$ is the minimizer,
\begin{equation}
\wt\EE(i_\star+a j)=\wt\EE(i_\star)+a\sum_{k=1}^nr(x_k,x_{k+1})i_\star(x_k,x_{k+1})+a^2\wt\EE(j)\ge\wt\EE(i_\star).
\end{equation}
Taking~$a\downarrow0$ then shows~``$\ge$'' in \eqref{E:12.32} and taking~$a\uparrow0$ then proves equality.

The fact that $e\mapsto r_ei_\star(e)$ obeys \eqref{E:12.32} implies that it is a gradient of a function. Specifically, we claim that there is~$f\colon\cmss V\to\R$ such that~$f(v)=0$ and
\begin{equation}
\label{E:12.34}
f(y)-f(x)=r(x,y)i_\star(x,y),\quad (x,y)\in\vec{\cmss E}.
\end{equation}
To see this, consider any path $x_0=v,x_1,\dots,x_n=x$ with $(x_k,x_{k+1})\in\cmss E$ for each~$k=0,\dots,n-1$ and let $f(x)$ be the sum of $-r_ei_\star(e)$ for edges along this path. The condition \eqref{E:12.32} then ensures that the value of~$f(x)$ is independent of the path chosen. Hence we get also \eqref{E:12.34}. 

Our next task is to compute~$f(u)$. Here we note that \eqref{E:12.34} identifies $\wt\EE(i_\star)$ with the quantity on the left of \eqref{E:12.29} and so
\begin{equation}
R_\eff(u,v)=\wt\EE(i_\star)=f(u)-f(v)=f(u).
\end{equation}
The function $\tilde f(x):=f(x)/R_\eff(u,v)$ thus belongs to~$\FF(u,v)$ and since, as is directly checked, $\EE(f)=\wt\EE(i_\star)=R_\eff(u,v)$, we get
\begin{equation}
C_\eff(u,v)\le\EE(\tilde f)=\frac1{R_\eff(u,v)^2}\,\EE(f)=\frac1{R_\eff(u,v)}.
\end{equation}
This gives~$C_\eff(u,v)R_\eff(u,v)\le1$, complementing the inequality from Lemma~\ref{lemma-12.11}. The claim follows.
\end{proofsect}

The above proof is based on optimizing over currents although one can also start from the minimizing potential. That this is more or less equivalent is attested by:

\begin{myexercise}[Ohm's law]
Prove that if~$f^\star$ is a minimizer of~$f\mapsto\EE(f)$ over~$f\in\FF(u,v)$, then 
\begin{equation}
\label{E:13.38nwt}
i(x,y):=c(x,y)\bigl[f^\star(y)-f^\star(x)\bigr],\quad (x,y)\in\vec{\cmss E},
\end{equation}
defines a current from~$u$ to~$v$ with~$\text{\rm val}(i)=C_\eff(u,v)$. (Compare \eqref{E:13.38nwt} with \eqref{E:12.34}.)
\end{myexercise}

\begin{figure}[t]
\vglue-1mm
\centerline{\includegraphics[width=0.7\textwidth]{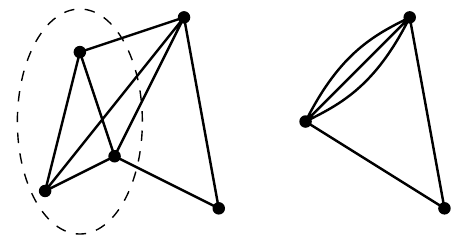}
}

\vglue0mm
\begin{quote}
\small 
\caption{
\label{fig-shorting}
\small
An example of shorting three vertices in a network (the one on the left) and collapsing them to one thus producing a new, effective network (the one on the right). All the outgoing edges are kept along with their original resistances.}
\normalsize
\end{quote}
\end{figure}

There is a natural extension of the effective resistance/conductance from pairs of vertices to pairs of sets. Indeed, for any pair of disjoint sets~$A,B\subset\cmss V$, we define $R_\eff(A,B)$ to be the effective resistance $R_\eff(\langle A\rangle,\langle B\rangle)$ in the network where all edges between the vertices in~$A$ as well as those between the vertices in~$B$ have been dropped and the vertices in~$A$ then merged into a single vertex~$\langle A\rangle$ and those in~$B$ merged into a vertex~$\langle B\rangle$. (The outgoing edges from~$A$ then emanate from~$\langle A\rangle$.) Similarly, we may define $C_\eff(A,B)$ as $C_\eff(\langle A\rangle,\langle B\rangle)$ or directly by
\begin{equation}
C_\eff(A,B):=\inf\Bigl\{\EE(f)\colon f|_A=1,\,f|_B=0\Bigr\}.
\end{equation}
Note that we simply set the potential to constants on~$A$ and~$B$. In the engineering vernacular, this amounts to \emph{shorting} the vertices in~$A$ and in~$B$; see Fig.~\ref{fig-shorting} for an illustration. The electrostatic duality still applies and so we have
\begin{equation}
C_\eff(A,B)=\frac1{R_\eff(A,B)}
\end{equation}
whenever~$A,B\subset\cmss V$ with~$A\cap B=\emptyset$.

\nopagebreak

\section{Markov chain connections and network reduction}
\noindent
With each electric network one can naturally associate a Markov chain on~$\cmss V$ with transition probabilities \begin{equation}
\cmss P(x,y):=\frac{c(x,y)}{\pi(x)}1_{(x,y)\in\cmss E}\quad\text{where}\quad\pi(x):=\sum_{y\colon (x,y)\in\cmss E}c(x,y).
\end{equation}
The symmetry condition $c(x,y)=c(y,x)$ then translates into
\begin{equation}
\pi(x)\cmss P(x,y)=\pi(y)\cmss P(y,x)
\end{equation}
 thus making~$\pi$ a reversible measure. Since~$\cmss G$ is connected, and the conductances are strictly positive, the Markov chain is also irreducible. Writing~$P^x$ for the law of the Markov chain started at~$x$, we then have:

\begin{myproposition}[Connection to Markov chain]
\label{prop-MC-conn}
The variational problem \eqref{E:12.20} has a unique minimizer~$f$ which is given by
\begin{equation}
\label{E:12.38}
f(x)=P^x(\tau_u<\tau_v)
\end{equation}
where~$\tau_z:=\inf\{n\ge0\colon X_n=z\}$.
\end{myproposition}

\begin{proofsect}{Proof}
Define an operator~$\LL$ on~$\ell^2(\cmss V)$ by
\begin{equation}
\label{E:12.40}
\LL f(x):=\sum_{y\colon (x,y)\in\cmss E}c(x,y)\bigl[f(y)-f(x)\bigr].
\end{equation}
(This is an analogue of the discrete Laplacian we encountered earlier in these lectures.) As is easy to check by differentiation, the minimizer of \eqref{E:12.20} obeys $\LL f(x)=0$ for all~$x\ne u,v$, with boundary values~$f(y)=1$ and~$f(v)=0$. The function $x\mapsto P^x(\tau_u<\tau_v)$ obeys exactly the same set of conditions; the claim thus follows from the next exercise.
\end{proofsect}

\begin{myexercise}[Dirichlet problem \&\ Maximum Principle]
\label{ex:max-principle}
Let~$\cmss U\subsetneq\cmss V$ and suppose~$f\colon\cmss V\to\R$ obeys $\LL f(x)=0$ for all~$x\in\cmss U$, where~$\LL$ is as in \eqref{E:12.40}. Prove that
\begin{equation}
f(x)=E^x\bigl(\,f(X_{\tau_{\cmss V\smallsetminus \cmss U}})\bigr),\quad x\in\cmss U,
\end{equation}
where $\tau_{\cmss V\smallsetminus \cmss U}:=\inf\{n\ge0\colon X_n\not\in\cmss U\}$. In particular, we have
\begin{equation}
\max_{x\in\cmss U}\bigl|f(x)\bigr|\le\max_{x\in\cmss V\smallsetminus\cmss U}\bigl|f(y)\bigr|
\end{equation}
and~$f$ is thus uniquely determined by its values on~$\cmss V\smallsetminus\cmss U$.
\end{myexercise}

We note that the above conclusion generally fails when~$\cmss U$ is allowed to be infinite. Returning to our main line of thought, from Proposition~\ref{prop-MC-conn} we immediately get:

\begin{mycorollary}
\label{cor-12.14}
Denoting~$\hat\tau_x:=\inf\{n\ge1\colon X_n=x\}$, for each~$u\ne v$ we have
\begin{equation}
\label{E:12.42}
\frac1{R_\eff(u,v)}=\pi(u)P^u(\hat\tau_u>\tau_v).
\end{equation}
\end{mycorollary}

\begin{proofsect}{Proof}
Let~$f$ be the minimizer of \eqref{E:12.20}. In light of $\LL f(x)=0$ for all~$x\ne u,v$, symmetrization arguments and~$f\in\FF(u,v)$ show
\begin{equation}
\EE(f)=-\sum_{x\in\cmss V}f(x)\LL f(x) = -f(u)\LL f(u)-f(v)\LL f(v)=-\LL f(u).
\end{equation}
The  representation \eqref{E:12.38} yields
\begin{equation}
\begin{aligned}
-\LL f(u)&=\pi(u)-\sum_{x\colon (u,x)\in\cmss E}c(u,x)P^x(\tau_x<\tau_v)
\\
&=\pi(u)\Bigl[1-\sum_{x\colon (u,x)\in\cmss E}\cmss P(u,x)P^x(\tau_x<\tau_v)\Bigr] =\pi(u)P^u(\hat\tau_u>\tau_v)\,,
\end{aligned}
\end{equation}
where we used the Markov property and the fact that $u\ne v$ implies $\hat\tau_u\ne\tau_v$. The claim follows from the Electrostatic Duality.
\end{proofsect}

The representation \eqref{E:12.42} leads to a criterion for recurrence/tran\-sience of a Markov chain~$X$ on an infinite, locally finite, connected electric network with positive resistances. Let $B(x,r)$ denote the ball in the graph-theoretical metric of radius~$r$ centered at~$x$. (The local finiteness ensures that~$B(x,r)$, as well as the set of all edges emanating from it, are finite.) First we note:

\begin{myexercise}
Denote by $C_\eff(x,B(x,r)^\cc)$, resp., $R_\eff(x,B(x,r)^\cc)$ the effective conductance, resp., resistance in the network where $B(x,r)^\cc$ has been collapsed to a single vertex. Prove, by employing a shorting argument, that $r\mapsto C_\eff(x,B(x,r)^\cc)$ is non-increasing.
\end{myexercise}

\noindent
This and the Electrostatic Duality ensure that
\begin{equation}
R_\eff(x,\infty):=\lim_{r\to\infty}R_\eff\bigl(x,B(x,r)^\cc\bigr)
\end{equation}
is well defined, albeit possibly infinite. The quantity $R_\eff(x,\infty)$, which we call effective resistance from~$x$ to infinity, depends on~$x$, although (by irreducibility) if~$R_\eff(x,\infty)$ diverges for one~$x$, then it diverges for all~$x$. We now have:

\begin{mycorollary}[Characterization of recurrence/transience]
\begin{equation}
\text{$X$\rm\ is recurrent}\quad\Leftrightarrow\quad R_\eff(\cdot,\infty)=\infty.
\end{equation}
\end{mycorollary}

\begin{proofsect}{Proof}
By Corollary~\ref{cor-12.14}, $P^x(\hat\tau_x>\tau_{B(x,r)^\cc})$ is proportional to~$C_\eff(u,B(x,r)^\cc)$. Since $\tau_{B(x,r)^\cc}\ge r$, $P^x(\hat\tau_x=\infty)$ is proportional to $R_\eff(x,\infty)^{-1}$.
\end{proofsect}

\nopagebreak
\begin{figure}[t]
\vglue-1mm
\centerline{\includegraphics[width=0.85\textwidth]{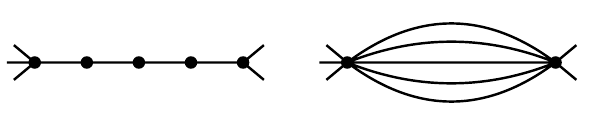}
}

\vglue0mm
\begin{quote}
\small 
\caption{
\label{fig-parallel-series}
\small
Examples of the situations handled by (left) the Series Law from Exercise~\ref{ex:series} and (right) the Parallel Law  from Exercise~\ref{ex:parallel}.}
\normalsize
\end{quote}
\end{figure}

The advantage of casting properties of Markov chains in electrostatic language is that we can manipulate networks using operations that do not always have a natural counterpart, or type of underlying monotonicity, in the context of Markov chains. We will refer to these operations using the (somewhat vague) term \emph{network reduction}.  Shorting part of the network serves as an example. Another example is:

\begin{myexercise}[Restriction to a subnetwork]
\label{ex:12.17}
Let~$\cmss V'\subset \cmss V$ and, for any function $f\colon \cmss V'\to\R$, set
\begin{equation}
\EE'(f):=\inf\bigl\{\EE(g)\colon g(x)=f(x)\,\,\forall x\in \cmss V'\bigr\}.
\end{equation}
Prove that~$\EE'(f)$ is still a Dirichlet energy of the form
\begin{equation}
\EE'(f)=\frac12\sum_{x,y\in \cmss V'}c'(x,y)\bigl[f(y)-f(x)\bigr]^2\,,
\end{equation}
where
\begin{equation}
\label{E:12.49}
c'(x,y):=\pi(x) P^x\bigl(X_{\hat\tau_{\cmss V'}}=y\bigr)
\end{equation}
with $\hat\tau_A:=\inf\{n\ge1\colon X_n\in A\}$.
\end{myexercise}

We remark that the restriction to a subnetwork does have a counterpart for the underlying Markov chain; indeed, it corresponds to observing the chain only at times when it is at~$\cmss V'$. (This gives rise to the formula \eqref{E:12.49}.)
The simplest instance is when $\cmss V'$ has only two vertices. Then, for any $u\ne v$, we have
\begin{equation}
\cmss V'=\{u,v\}\quad\Rightarrow\quad c'(u,v) = C_\eff(u,v).
\end{equation}
A another relevant example is the content of:

\begin{myexercise}[Series law]
\label{ex:series}
Suppose~$\mathfrak G$ contains a string of vertices $x_0,\dots,x_{n}$ such that $(x_{i-1},x_{i})\in\cmss E$ for each~$i=1,\dots,n$ and such that, for $i=1,\dots,n-1$, the vertex~$x_i$ has no other neighbors than~$x_{i-1}$ and~$x_{i+1}$. Prove that in the reduced network with $\cmss V':=\cmss V\smallsetminus\{x_1,\dots,x_{n-1}\}$ the string is replaced by an edge $(x_0,x_n)$ with resistance
\begin{equation}
r'(x_0,x_n):=\sum_{i=1}^n r(x_{i-1},x_i)\,.
\end{equation}
\end{myexercise}

\nopagebreak
\begin{figure}[t]
\vglue-1mm
\centerline{\includegraphics[width=0.6\textwidth]{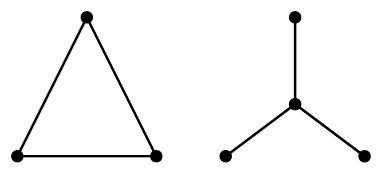}
}
\vglue0mm
\begin{quote}
\small 
\caption{
\label{fig-star-triangle}
\small
The ``triangle'' and ``star'' networks from Exercise~\ref{ex:12.20a}.}
\normalsize
\end{quote}
\end{figure}

There are other operations that produce equivalent networks which are not of the type discussed in Exercise~\ref{ex:12.17}:

\begin{myexercise}[Parallel law]
\label{ex:parallel}
Suppose $\mathfrak G$ contains~$n$ edges $e_1,\dots,e_n$ between vertices~$x$ and~$y$ with~$e_i$ having conductance~$c(e_i)$. Prove that we can replace these by a single edge~$e$ with conductance
\begin{equation}
c'(e):=\sum_{i=1}^n c(e_i).
\end{equation}
\end{myexercise}

\noindent
Here is another operation that does fall under the scheme of Exercise~\ref{ex:12.17}:

\begin{myexercise}[Star-triangle transformation]
\label{ex:12.20a}
Consider an electric network that contains (among others) three nodes $\{1,2,3\}$ and an edge between every pair of these nodes. Write~$c_{ij}$ for the conductance of edge $(i,j)$. Prove that an equivalent network is produced by replacing the ``triangle'' on these vertices by a ``star'' which consists of four nodes $\{0,1,2,3\}$ and three edges $\{(0,i)\colon i=1,2,3\}$, and no edges between the remaining vertices. See Fig.~\ref{fig-star-triangle}. [We call this replacement the ``star-triangle transformation.'']
\end{myexercise}

\noindent
The start-triangle transformation is a very powerful tool in network theory. Same tools as in the previous exercise can be used to conclude:

\begin{myexercise}
\label{ex:12.20}
For the network with $\cmss V:=\{1,2,3\}$ and $\cmss E:=\{(i,j)\colon 1\le i<j\le 3\}$, let~$c_{ij}$ denote the conductance of edge $(i,j)$. Denoting $R_{ij}:=R_\eff(i,j)$, prove that
\begin{equation}
\label{E:5.35b}
\frac{c_{12}}{c_{12}+c_{13}}=\frac{R_{13}+R_{23}-R_{12}}{2R_{23}}.
\end{equation}
\end{myexercise}

\noindent
The network reduction ideas will continue to be used heavily through the remaining lectures.


\chapter{Effective resistance control}
\label{lec-14}\noindent
The aim of this lecture is to develop further methods to control effective resistance/conductance in networks related to the DGFF as discussed above. We start by proving geometric representations which allow us to frame estimates of these quantities in the language of percolation theory. We then return to the model at hand and employ duality considerations to control the effective resistance across squares in~$\Z^2$. A Russo-Seymour-Welsh (type of) theory permits similar control for the effective resistance across rectangles. Finally, we use these to control the upper tail of the effective resistance between far-away points.

\section{Path-cut representations}
\noindent
The network reduction ideas mentioned at the end of the previous lecture naturally lead to representations on the effective resistance/conductance by quantities involving geometric objects such as paths and cuts. As these representations are the cornerstone of our approach, we start by explaining them in detail.

A path~$P$ from~$u$ to~$v$ is a sequence of edges~$e_1,\dots,e_n$, which we think of as oriented for this purpose, with~$e_1$ having the initial point at~$u$ and~$e_n$ the terminal endpoint at~$v$ and with the initial point of $e_{i+1}$ equal to terminal point of~$e_i$ for each~$i=1,\dots,n-1$. We will often identify~$P$ with the set of these edges; the notation $e\in P$ thus means that the path~$P$ crosses the unoriented edge~$e$. 
The following lemma is classical:
 
\begin{mylemma}
\label{lemma-12.19}
Suppose~$\PP$ is a finite set of edge disjoint paths --- i.e., those for which $P\ne P'$ implies $P\cap P'=\emptyset$ --- from~$u$ to~$v$. Then
\begin{equation}
\label{E:12.50}
R_\eff(u,v)\le\Bigl[\sum_{P\in\PP}\frac1{\sum_{e\in P}r_e}\Bigr]^{-1}.
\end{equation}
\end{mylemma}

\begin{proofsect}{Proof}
Dropping loops from paths decreases the quantity on the right of \eqref{E:12.50}, so we may and will assume that each path in~$\PP$ visits each edge at most once.
The idea is to route current along the paths in~$\PP$ to arrange for a unit net current to flow from~$u$ to~$v$.
Let $\overline R$ denote the quantity on the right of \eqref{E:12.50} and, for each~$P\in\PP$, set $i_P:=\overline R/\sum_{e\in P}r_e$. Define $i(e):=i_P$ if~$e$ lies in, and is oriented in the same direction as~$P$ (recall that the paths are edge disjoint and each edge is visited by each path at most once) and~$i(e):=0$ if~$e$ belongs to none of the paths in~$\PP$. From $\sum_{P\in\PP} i_P=1$ we then infer $i\in\II(u,v)$. A calculation shows $\wt\EE(i)=\overline R$ and so~$R_\eff(u,v)\le\overline R$.
\end{proofsect}

A natural question to ask is whether the upper bound \eqref{E:12.50} can possibly be sharp. As it turns out, all that stands in the way of this is the edge-disjointness requirement. This is overcome in:

\begin{myproposition}[Path representation of effective resistance]
\label{prop-12.20}
Let  $\mathfrak P_{u, v}$ denote the set of finite collections of paths from~$u$ to~$v$.
Then
\begin{equation}
\label{E:12.52}
R_\eff(u, v) = \inf_{\mathcal P \in \mathfrak P_{u, v}} \,\, \inf_{\{r_{e, P}\colon P\in \mathcal P,\,e\in P\} \in \mathfrak R_\mathcal P}\Bigl[\,\sum_{P\in \mathcal P} \frac{1}
{\sum_{e\in P} r_{e, P}}\Bigr]^{-1}\,,
\end{equation}
where $\mathfrak R_\mathcal P$ is the set of all collections of positive numbers  $\{r_{e, P}\colon P\in \mathcal P,\,e\in P\}$ such that 
\begin{equation}
\label{E:12.53}
\sum_{P\in \mathcal P} \frac{1}{r_{e, P}} \leq \frac{1}{r_e},\quad e \in \cmss E\,.
\end{equation}
The infima in~\eqref{E:12.52} are (jointly) achieved.
\end{myproposition}

\begin{proofsect}{Proof}
Pick a collection of paths $\PP\in \mathfrak P_{u, v}$ and positive numbers $\{r_{e,P}\colon P\in\PP\}$ satisfying \eqref{E:12.53}. Then split each edge~$e$ into a collection of edges $\{e_P\colon P\in\PP\}$ and assign resistance~$r_{e,P}$ to~$e_P$. If the inequality in \eqref{E:12.53} is strict, introduce also a dummy copy~$\tilde e$ of~$e$ and assign conductance $c_{\tilde e}:=1/r_e-\sum_{P\in\PP}1/r_{e,P}$ to~$\tilde e$. The Parallel Law shows that this operation produces an equivalent network in which, by way of interpretation, the paths in~$\PP$ are mutually edge disjoint. Lemma~\ref{lemma-12.19} then gives ``$\le$'' in~\eqref{E:12.52}. 

To get equality in \eqref{E:12.52}, let~$i_\star\in\II(u,v)$ be such that $\wt\EE(i_\star)=R_\eff(u,v)$. We will now recursively define a sequence of currents~$i_k$ (not necessarily of unit value) and paths~$P_k$ from~$u$ and~$v$.  First solve:

\begin{myexercise}
\label{ex:12.21}
Suppose~$e\mapsto i(e)$ is a current from~$u$ to~$v$ with~$\text{val}(i)>0$.  Show that there is a simple path~$P$ from~$u$ to~$v$ such that $i(e)>0$ for each~$e\in P$ (which we think of as oriented in the direction of~$P$).
\end{myexercise}

\noindent
We cast the recursive definition as an algorithm: INITIATE by setting $i_0:=i_\star$. Assuming that~$i_{k-1}$ has been defined for some~$k\ge1$, if~$\text{val}(i_{k-1})=0$ then STOP, else use Exercise~\ref{ex:12.21} to find a path~$P_k$ from~$u$ to~$v$ where $i_{k-1}(e)>0$ for each~$e\in P_k$ oriented along the path. Then set $\alpha_k:=\min_{e\in P_k}|i_{k-1}(e)|$, let
\begin{equation}
\label{E:14.4nwt}
i_k(e):=i_{k-1}(e)-\alpha_k\, \sgn(i_{k-1}(e))1_{\{e\in P_k\}}
\end{equation}
and, noting that~$i_k$ is a current from~$u$ to~$v$ with~$\text{val}(i)\ge0$, REPEAT. 

The construction ensures that $\{e\in E\colon i_k(e)\ne0\}$ is strictly decreasing in~$k$ and so the algorithm will terminate after a finite number of steps. Also $k\mapsto |i_k(e)|$ and~$k\mapsto\text{val}(i_k)$ are non-increasing and, as is checked from \eqref{E:14.4nwt}, 
\begin{equation}
\label{E:12.55}
\forall e\in\cmss E\colon\quad
\sum_{k\colon e\in P_k}\alpha_k \le |i_\star(e)|\quad\text{and}\quad\sum_k\alpha_k = \text{\rm val}(i_\star)=1.
\end{equation}
We now set $r_{e,P_k}:=|i_\star(e)| r_e/\alpha_k$ and note that \eqref{E:12.55} shows $\sum_k 1/r_{e,P_k}\le 1/r_e$ for each~$e\in\cmss E$. Moreover, from \eqref{E:12.55} we also get
\begin{equation}
\begin{aligned}
R_\eff(u,v)&=\sum_{e\in\cmss E}r_ei_\star(e)^2
\ge\sum_{e\in\cmss E}\,r_e|i_\star(e)|\sum_{k\colon e\in P_k}\alpha_k
\\&=\sum_{e\in\cmss E}\,\sum_{k\colon e\in P_k}r_{e,P_k}\alpha_k^2
=\sum_k\,\alpha_k^2\Bigl(\sum_{e\in P_k}r_{e,P_k}\Bigr).
\end{aligned}
\end{equation}
Denoting the quantity in the large parentheses by~$R_k$, among all positive~$\alpha_k$'s satisfying \eqref{E:12.55}, the right-hand side is minimized by~$\alpha_k:=\frac{1/R_k}{\sum_j 1/R_j}$. Plugging this in, we get ``$\ge$'' in~\eqref{E:12.52} and thus the whole claim.
\end{proofsect}

As it turns out, the effective conductance admits an analogous geometric variational characterization as well. Here one needs the notion of a \emph{cut}, or a \emph{cutset} form~$u$ to~$v$ which is a set of edges in~$\cmss E$ such that every path from~$u$ to~$v$ must contain at least one edge in this set. We again start with a classical bound:

\begin{mylemma}[Nash-Williams estimate]
\label{lemma-NW}
For any collection $\Pi$ of edge-disjoint cutsets  from~$u$ to~$v$,
\begin{equation}
\label{E:12.57}
C_\eff(u,v)\le\Bigl[\sum_{\pi\in\Pi}\frac1{\sum_{e\in\pi}c_e}\Bigr]^{-1}\,.
\end{equation}
\end{mylemma}

\begin{proofsect}{Proof}
Let~$i\in\II(u,v)$. The proof is based on:

\begin{myexercise}
For any cutset~$\pi$ from~$u$ to~$v$, $\sum_{e\in\pi}|i(e)|\ge1$.
\end{myexercise}

\noindent
Indeed, the Cauchy-Schwarz inequality and \eqref{E:12.17} tell us
\begin{equation}
1 \le \Bigl[\sum_{e\in\pi}|i(e)|\Bigr]^2\le \Bigl[\sum_{e\in\pi}r_ei(e)^2\Bigr]\Bigl[\sum_{e\in\pi}c_e\Bigr].
\end{equation}
The assumed edge-disjointness of the cutsets in~$\Pi$ then yields
\begin{equation}
\wt\EE(i)\ge\sum_{\pi\in\Pi}\sum_{e\in\pi}r_ei(e)^2\ge\sum_{\pi\in\Pi}\frac1{\sum_{e\in\pi}c_e}.
\end{equation}
The claim follows by the Electrostatic Duality.
\end{proofsect}

\begin{myremark}
Lemma~\ref{lemma-NW} is easy to prove by network reduction arguments when the cutsets are \emph{nested} meaning that they can be ordered in a sequence $\pi_1,\dots,\pi_n$ such that~$\pi_i$ separates~$\pi_{i-1}$ (as well as~$u$) from~$\pi_{i+1}$ (as well as~$v$). However, as the above proof shows, this geometric restriction is not needed (and, in fact, would be inconvenient to carry around).
\end{myremark}

The bound \eqref{E:12.57} was first proved by Nash-Williams~\cite{Nash-Williams} as a tool for proving recurrence of an infinite network. Indeed, to get $C_\eff(0,\infty):=1/R_\eff(0,\infty)=0$, it suffices to present a disjoint family of cutsets whose reciprocal total conductances add up to infinity. However, as far as the actual computation of $C_\eff(u,v)$ is concerned, \eqref{E:12.57} is generally not sharp; again, mostly due to the requirement of edge-disjointness. The following proposition provides the needed fix:

\begin{myproposition}[Cutset representation of effective conductance]
\label{prop-12.22}
Let $\mathfrak S_{u, v}$ be the set of all finite collections of cutsets between~$u$ and~$v$. Then 
\begin{equation}
\label{E:12.63}
C_\eff(u, v) = \inf_{\Pi \in \mathfrak S_{u, v}} \, \inf_{\{c_{e, \pi}\colon \pi\in \Pi,\,e\in\pi\} \in \mathfrak C_{\Pi}}\,\Bigl[\,\sum_{\pi\in \Pi} \frac{1}{\sum_{e\in \pi} c_{e,
\pi}}\Bigr]^{-1}\,,
\end{equation}
where $\mathfrak C_{\Pi}$ is the set of all families of positive numbers  $\{c_{e, \pi}\colon \pi\in \Pi,\,e\in\pi\}$ such that 
\begin{equation}
\label{E:12.64}
\sum_{\pi \in \Pi} \frac{1}{c_{e, \pi}} \leq \frac{1}{c_e}\,, \quad e \in \cmss E\,.
\end{equation}
 The infima in~\eqref{E:12.63} are (jointly) achieved.
\end{myproposition}

\begin{proofsect}{Proof}
Let~$\Pi$ be a family of cutsets from~$u$ to~$v$ and let $\{c_{e,\pi}\}\in\mathfrak C_\Pi$. We again produce an equivalent network as follows: Replace each edge~$e$ involved in these cutsets by a series of edges~$e_\pi$, $\pi\in\Pi$. Assign resistance $1/c_{e,\pi}$ to~$e_\pi$ and, should the inequality in \eqref{E:12.64} be strict for~$e$, add a dummy edge~$\tilde e$ with resistance~$r_{\tilde e}:=1/c_e-\sum_{\pi\in\Pi}1/c_{e,\pi}$. The cutsets can then be deemed edge-disjoint; and so we get~``$\le$'' in \eqref{E:12.63} by Lemma~\ref{lemma-NW}.

To prove equality in \eqref{E:12.63}, we consider the minimizer~$f^\star$ of the variational problem defining~$C_\eff(u,v)$. Using the operator in \eqref{E:12.40}, note that~$\LL f^\star(x)=0$ for~$x\ne u,v$. This is an important property in light of:

\begin{myexercise}
\label{ex:12.27}
Let~$f\colon\cmss V\to\R$ be such that~$f(u)>f(v)$ and abbreviate
\begin{equation}
\label{E:D-def}
D:=\bigl\{x\in \cmss V\colon f(x)=f(u)\bigr\}.
\end{equation}
 Assume $\LL f(x)=0$ for~$x\not\in D\cup\{v\}$. Prove that $\pi:=\partial D$ defines a cutset~$\pi$ from~$u$ to~$v$ such that $f(x)>f(y)$ holds for each edge~$(x,y)\in\pi$ oriented so that~$x\in D$ and~$y\not\in D$.
\end{myexercise}

\noindent
We will now define a sequence of functions $f_k\colon\cmss V\to\R$ and cuts~$\pi_k$ by the following algorithm: INITIATE by $f_0:=f^\star$. If $f_{k-1}$ is constant then STOP, else use Exercise~\ref{ex:12.27} with $D_{k-1}$ related to~$f_{k-1}$ as~$D$ is to~$f$ in \eqref{E:D-def} to define~$\pi_k:=\partial D_{k-1}$. Noting that $f_{k-1}(x)>f_{k-1}(y)$ for each edge $(x,y)\in\pi$ (oriented to point from~$u$ to~$v$), set~$\alpha_k:=\min_{(x,y)\in\pi}|f_{k-1}(x)-f_{k-1}(y)|$ and let
\begin{equation}
f_k(z):=f_{k-1}(z)-\alpha_k\,1_{D_{k-1}}(z).
\end{equation}
Then REPEAT.

As is checked by induction,~$k\mapsto D_k$ is strictly increasing and $k\mapsto f_k$ is non-increasing with $f_k=f^\star$ on~$\cmss V\smallsetminus D_k$. In particular, we have $\LL f_k(x)=0$ for all $x\not\in D_k\cup\{v\}$ and so Exercise~\ref{ex:12.27} can repeatedly be used. The premise of the strict inequality $f_k(u)>f_k(v)$ for all but the final step is the consequence of the Maximum Principle. 

Now we perform some elementary calculations. Let~$\Pi$ denote the set of the cutsets~$\pi_k$ identified above. For each edge~$(x,y)$ and each~$\pi_k\in\Pi$, define
\begin{equation}
\label{E:13.15}
c_{e,\pi_k} := \frac{|f^\star(y)-f^\star(x)|}{\alpha_k}\,c_e.
\end{equation}
The construction and the fact that $f^\star(u)=1$ and~$f^\star(v)=0$ imply 
\begin{equation}
\label{E:13.16}
\sum_k\alpha_k=1\,.
\end{equation}
For any $e=(x,y)$, we also get
\begin{equation}
\sum_{k\colon e\in\pi_k}\alpha_k = \bigl|f^\star(y)-f^\star(x)\bigr|.
\end{equation}
In particular, the collection $\{c_{e,\pi}\colon\pi\in\Pi,\,e\in\pi\}$ obeys \eqref{E:12.64}. Moreover, \eqref{E:13.15} also shows, for any $e=(x,y)$,
\begin{equation}
\sum_{k\colon e\in\pi_k} c_{e,\pi_k}\alpha_k^2 = c_e\bigl|f^\star(y)-f^\star(x)\bigr|\sum_{k\colon e\in\pi_k}\alpha_k = c_e\bigl|f^\star(y)-f^\star(x)\bigr|^2.
\end{equation}
Summing over~$e\in\cmss E$ and rearranging the sums yields
\begin{equation}
C_\eff(u,v) = \EE(f^\star) = \sum_{e\in\cmss E}\,\,\sum_{k\colon e\in\pi_k} c_{e,\pi_k}\alpha_k^2 =  
\sum_k\alpha_k^2\Bigl(\sum_{e\colon e\in\pi_k}c_{e,\pi_k}\Bigr)\,.
\end{equation}
Denoting the quantity in the large parentheses by~$C_k$, among all non-negative~$\alpha_k$'s satisfying \eqref{E:13.16}, the right-hand side is minimal for~$\alpha_k:=\frac{1/C_k}{\sum_j1/C_j}$. This shows ``$\ge$'' in \eqref{E:12.63} and thus finishes the proof.
\end{proofsect}

We note that the above derivations are closely related to characterizations of the effective resistance/conductance based on optimizing over \emph{random} paths and cuts. These are rooted in T.~Lyons' random-path method~\cite{TLyons} for proving finiteness of the effective resistance. Refined versions of these characterizations can be found in Berman and Konsowa~\cite{Berman-Konsowa}.

\section{Duality and effective resistance across squares}
\noindent
Although the above derivations are of independent interest, our main reason for giving their full account is to demonstrate the \emph{duality} ideas that underlie the derivations in these final lectures. This is based on the close similarity of the variational problems \twoeqref{E:12.52}{E:12.53} and \twoeqref{E:12.63}{E:12.64}. Indeed, these are identical provided we somehow manage to
\begin{enumerate}
\item[(1)] swap paths for cuts, and
\item[(2)] exchange resistances and conductances.
\end{enumerate}
Here (2) is achievable readily assuming symmetry of the underlying random field: For conductances (and  resistances) related to a random field~$h$ via \eqref{E:12.13},
\begin{equation}
\label{E:13.21}
h\,\,\,\laweq\,-h\quad\Rightarrow\quad \{c_e\colon e\in\cmss E\}\,\,\laweq\,\, \{r_e\colon e\in\cmss E\}.
\end{equation}
For~(1) we can perhaps hope to rely on the fact that paths on \emph{planar} graphs between the left and right sides of a rectangle are in one-to-one correspondence with cuts in the \emph{dual} graph separating the top and bottom sides of that rectangle. Unfortunately, swapping primal and dual edges seems to void the dependence \eqref{E:12.13} on the underlying field which we need for \eqref{E:13.21}. We will thus stay on the primal lattice and implement the path-cut duality underlying (1) only approximately.

\newcommand{\LR}{{\text{\tiny\rm LR}}}
\newcommand{\UD}{{\text{\tiny\rm UD}}}

Consider a rectangular box~$S$ of the form
\begin{equation}
B(M,N):=([-M,M]\times[N,N])\cap\Z^2
\end{equation}
and let $\partial^{\text{left}}S$, $\partial^{\text{down}}S$, $\partial^{\text{right}}S$, $\partial^{\text{up}}S$ denote the sets of vertices in~$S$ that have a neighbor in~$S^\cc$ to the left, downward, the right and upward thereof, respectively. Regarding~$S$, along with the edges with both endpoints in~$S$, as an electric network with conductances depending on a field~$h$ as in \eqref{E:12.13}, let
\begin{equation}
R^\LR_{S}(h):=R_\eff\bigl(\partial^{\text{left}}S,\partial^{\text{right}}S\bigr)
\end{equation}
and
\begin{equation}
R^\UD_{S}(h):=R_\eff\bigl(\partial^{\text{up}}S,\partial^{\text{down}}S\bigr)
\end{equation}
denote the left-to-right and up-down effective resistances across the box~$S$, respectively. To keep our formulas short, we will also abbreviate
\begin{equation}
B(N):=B(N,N).
\end{equation}
A key starting point of all derivations is the following estimate:

\begin{myproposition}[Effective resistance across squares]
\label{prop-13.9}
There is $\hat c>0$ and, for each $\epsilon>0$, there is~$N_0\ge1$ such that for all~$N\ge N_0$ and all~$M\ge 2N$, 
\begin{equation}
\label{E:13.24ai}
\BbbP\Bigl(R^\LR_{B(N)}(h^{B(M)})\le\texte^{\hat c\log\log M}\Bigr)\ge\frac12-\epsilon,
\end{equation}
where $h^{B(M)}:=$ DGFF in~$B(M)$.
\end{myproposition}

The proof will involve a sequence of lemmas. We begin by a general statement that combines usefully the \emph{electrostatic} duality along with approximate \emph{path-cut} duality. Consider a finite network~$\cmss G=(\cmss V,\cmss E)$ with resistances $\{r_e\colon e\in\cmss E\}$. A dual, or reciprocal, network is that on the same graph but with resistances $\{r_e^\star\colon e\in\cmss E\}$ where
\begin{equation}
\label{E:13.25uia}
r_e^\star:=\frac1{r_e},\quad e\in\cmss E.
\end{equation}
We will denote to the dual network as~$\cmss G^\star$. We will write~$R^\star_\eff$, resp., $C^\star_\eff$ to denote the effective resistances in~$\cmss G^\star$. We also say that edges~$e$ and~$e'$ are adjacent, with the notation~$e\sim e'$, if they share exactly one endpoint. We then have:

\begin{mylemma}
\label{lemma-dual}
Let~$(A,B)$ and~$(C,D)$ be pairs of non-empty disjoint subsets of~$\cmss V$ such that every path from~$A$ to~$B$ has a vertex in common with every path from~$C$ to~$D$. Then
\begin{equation}
R_\eff(A,B)R_\eff^\star(C,D)\ge\frac1{4\dmax^2\rhomax}\,,
\end{equation}
where $\dmax$ is the maximum vertex degree in~$\cmss G$ and
\begin{equation}
\rhomax:=\max\bigl\{r_e/r_{e'}\colon e\sim e'\bigr\}.
\end{equation}
\end{mylemma}

\begin{proofsect}{Proof}
The proof relies on the fact that every path~$P$ between $C$ and $D$ defines a cutset~$\pi_P$ between~$A$ and~$B$ by taking~$\pi_P$ to be the set of all edges adjacent to the edges in~$P$, but excluding the edges in $(A\times A)\cup(B\times B)$. In light of the Electrostatic Duality, we need to show 
\begin{equation}
\label{E:2.36ai}
C_\eff(A, B)\le 4 \dmax^2\rho_{\max} R_\eff^\star(C, D)\,.
\end{equation}
Aiming to use the variational characterizations \twoeqref{E:12.52}{E:12.53} and \twoeqref{E:12.63}{E:12.64},
fix a collection of paths $\mathcal P \in \mathfrak P_{C, D}$ and positive numbers $\{r'_{e,P}\colon P\in \mathcal P, e\in P\}$ such that
\begin{equation}
\label{eq-r-e3}
\sum_{P\in \mathcal P} \frac{1}{r'_{e, P}} \leq \frac{1}{r_e^\star}=\frac1{c_e},\quad e\in\cmss E\,,
\end{equation}
where the equality is a rewrite of \eqref{E:13.25uia}.
For each~$e\in\cmss E$ and each~$P\in\PP$, consider the cut~$\pi_P$ defined above and let
\begin{equation}
c_{e,\pi_P}:=2\dmax\rhomax\biggl(\,\sum_{e'\in P\colon e'\sim e}\,\frac1{r'_{e',P}}\biggr)^{-1}\,.
\end{equation}
Then, for each~$e\in\cmss E$, \eqref{eq-r-e3} along with the definitions of~$\rhomax$ and~$\dmax$ yield
\begin{equation}
\label{E:13.29ue}
\begin{aligned}
\sum_{P\in\PP}\frac1{c_{e,\pi_P}}&=\frac1{2\dmax\rhomax}\sum_{P\in\PP}\,
\sum_{e'\in P\colon e'\sim e}\,\frac1{r'_{e',P}}
\\
&\le\frac1{2\dmax\rhomax}\sum_{e'\sim e}\frac1{c_{e'}}
\le\frac1{2\dmax}\sum_{e'\sim e}\frac1{c_e}\le\frac1{c_e}\,.
\end{aligned}
\end{equation}
Using that $(\sum_ix_i)^{-1}\le\sum_i 1/x_i$ for any positive~$x_i$'s in turn shows
\begin{equation}
\sum_{e\in\pi_P}c_{e,\pi_P}
\le 2\dmax\rhomax \sum_{e\in\pi_P} \,\sum_{e'\in P\colon e'\sim e}\,r'_{e',P}
\le 4\dmax^2\rhomax\sum_{e'\in P}\,r'_{e',P}\,.
\end{equation}
In light of \eqref{E:13.29ue}, from \eqref{E:12.63} we thus have
\begin{multline}
\quad
C_\eff(A,B)\le \Bigl(\,\sum_{P\in\PP}\frac1{\sum_{e\in\pi_P}c_{e,\pi_P}}\Bigr)^{-1}
\\
\le 4\dmax^2\rhomax\Bigl(\,\sum_{P\in\PP}\frac1{\sum_{e'\in P}r'_{e',P}}\Bigr)^{-1}\,.
\quad
\end{multline}
This holds for all $\mathcal P$ and all positive $\{r'_{e, P}\colon P\in \mathcal P, e\in P\}$ subject to \eqref{eq-r-e3} and so \eqref{E:2.36ai}  follows from \eqref{E:12.52}.
\end{proofsect}

We also have the opposite inequality, albeit under somewhat different conditions on $(A,B)$ and~$(C,D)$:

\begin{mylemma}
\label{lemma-dual2}
For~$\dmax$ and~$\rhomax$ as in Lemma~\ref{lemma-dual}, 
let~$(A,B)$ and~$(C,D)$ be pairs of non-empty disjoint subsets of~$\cmss V$ such that for every cutset~$\pi$ between~$C$ and~$D$, the set of edges with one or both vertices in common with some edge in~$\pi$ contains a path from~$A$ to~$B$.~Then 
\begin{equation}
\label{E:13.33}
R_\eff(A , B) R_\eff^\star(C , D)\le 4{\dmax}^2\rho_{\max}\,.
\end{equation}
\end{mylemma}

\begin{proofsect}{Proof (sketch)}
For any cutset~$\pi$ between $C$ and $D$, the assumptions ensure the existence of a path $P_\pi$ from~$u$ to~$v$ that consists of edges that have one or both endpoints in common with some edge in~$\pi$. Pick a family of cutsets~$\Pi$ between~$C$ and~$D$ and positive numbers $\{c_{e, \pi}' \colon \pi \in \Pi,\, e\in P_\pi \}$ such that
\begin{equation}
\label{E:13.34}
\sum_{\pi \in \Pi}\frac{1}{c_{e, \pi}'} \leq \frac1{c_e^\star}=\frac{1}{r_e}\,.
\end{equation}
Following the exact same sequence of steps as in the proof of Lemma~\ref{lemma-dual}, the reader will readily construct $\{r_{e, P_\pi}\colon \pi \in \Pi,e\in P_\pi\}$ satisfying \eqref{E:12.53} such that
\begin{equation}
R_\eff(A,B)\le\Bigl(\,\sum_{\pi \in \Pi} \frac{1}{\sum_{e \in P_\pi} r_{e, 
P_\pi}}\Bigr)^{-1} \leq 4{\dmax}^2\rho_{\max}\Bigl(\,\sum_{\pi \in \Pi} \frac{1}{\sum_{e \in \pi} c_{e, \pi}'}\Bigr)^{-1}\,,
\end{equation}
where the first inequality is by Proposition~\ref{prop-12.20}. As this holds for all $\Pi$ and all positive $\{c_{e, \pi}' \colon \pi \in \Pi, e\in \pi \}$ satisfying \eqref{E:13.34}, we get \eqref{E:13.33}.
\end{proofsect}

The previous lemma shows that, in bounded degree planar graphs, primal and dual effective resistances can be compared as long as we can bound the ratio of the resistances on adjacent edges. Unfortunately, this would not work for resistances derived from the DGFF in~$B(M)$ via \eqref{E:12.13}; indeed, there~$\rhomax$ can be as large as a power of~$M$ due to the large \emph{maximal} local roughness of the field. We will resolve this by decomposing the DGFF into a smooth part, where the associated~$\rhomax$ is sub-polynomial in~$M$, and a rough part whose influence can be controlled directly. This is the content of:

\begin{mylemma}[Splitting DGFF into smooth and rough fields]
\label{lemma-split}
Recall that $h^{B(N)}$ denotes the DGFF in~$B(N)$. There is~$c>0$ and, for each~$N\ge1$, there are Gaussian fields $\varphi$ and~$\chi$ on~$B(N)$ such that
\begin{enumerate}
\item[(1)] $\varphi$ and~$\chi$ are independent with $\varphi+\chi\,\,\laweq\,\,h^{B(N)}$,
\item[(2)] $\Var(\chi_x)\le c\log\log N$ for each~$x\in B(N)$, and
\item[(3)] $\Var(\varphi_x-\varphi_y)\le c/\log N$ for every adjacent pair~$x,y\in B(N/2)$.
\end{enumerate}
Moreover, the law of~$\varphi$ is invariant under the rotations of~$B(N)$ by multiples of~$\ffrac\pi2$.
\end{mylemma}

\begin{proofsect}{Proof}
Let~$\{Y_n\colon n\ge0\}$ denote the discrete time simple symmetric random walk on~$\Z^2$ with holding probability~$1/2$ at each vertex and, given any finite~$\Lambda\subset\Z^2$, let~$\tau_{\Lambda^\cc}$ be the first exit time from~$\Lambda$. Writing~$P^x$ for the law of the walk started at~$x$ and denoting
\begin{equation}
\cmss Q(x,y):=P^x(Y_1=y,\,\tau_{\Lambda^\cc}\ge1),
\end{equation}
 we first ask reader to solve:

\begin{myexercise}
\label{ex:14.13}
Prove the following:
\begin{enumerate}
\item[(1)] the $n$-th matrix power of~$\cmss Q$ obeys $\cmss Q^n(x,y)=P^x(Y_n=y,\,\tau_{\Lambda^\cc}\ge n)$,
\item[(2)] the matrices~$\{\cmss Q^n(x,y)\colon x,y\in \Lambda\}$ are symmetric and positive semi-definite,
\item[(3)] the Green function in~$\Lambda$ obeys
\begin{equation}
\label{E:13.34a}
G^\Lambda(x,y) = \sum_{n\ge0}\frac12\cmss Q^n(x,y).
\end{equation}
\end{enumerate}
\end{myexercise}
We now apply the above for $\Lambda:=B(N)$. Writing $C_1(x,y)$ for the part of the sum in \eqref{E:13.34a} corresponding to~$n\le\lfloor\log N\rfloor^2$ and $C_2(x,y)$ for the remaining part of the sum, Exercise~\ref{ex:14.13}(2) ensures that the kernels~$C_1$ and~$C_2$ are symmetric and positive semidefinite with $G^D=C_1+C_2$. The fields
\begin{equation}
\chi:=\NN(0,C_1)\quad\text{and}\quad\varphi:=\NN(0,C_2)\quad\text{with}\quad\chi\independent\varphi
\end{equation}
then realize~(1) in the statement. To get~(2), we just sum the standard heat-kernel bound $\cmss Q^n(x,x)\le c/n$ (valid uniformly in~$x\in B(N)$). For~(3) we pick neighbors $x,y\in B(N/2)$ and use the Strong Markov Property to estimate
\begin{equation}
\begin{aligned}
\Bigl|E\bigl[&\varphi_x(\varphi_x-\varphi_y)\bigr]\Bigr|
\\
&\le\sum_{n>\lfloor\log N\rfloor^2}\Bigl|P^x(Y_n=y,\,\tau_{B(N)^\cc}>n)-P^x(Y_n=x,\,\tau_{B(N)^\cc}>n)\Bigr|
\\
&\le \sum_{n>\lfloor\log N\rfloor^2}\Bigl|P^x(Y_n=y)-P^x(Y_n=x)\Bigr|
\\
&\qquad\qquad+\sum_{n\ge1}E^x\Bigl|P^{X_{\tau_{B(N)^\cc}}}(Y_n=y)-P^{X_{\tau_{B(N)^\cc}}}(Y_n=x)\Bigr|.
\end{aligned}
\end{equation}
By \cite[Theorem~2.3.6]{Lawler-Limic} there is~$c>0$ such that for all~$(x',y')\in E(\Z^2)$ and all~$z\in\Z^2$,
\begin{equation}
\bigl|P^z(Y_n=y')-P^z(Y_n=x')\bigr| \le cn^{-3/2}.
\end{equation}
This shows that the first sum is~$O(1/\log N)$ and that, in light of~$x,y\in B(N/2)$, the second sum is $O(1/\sqrt N)$. This readily yields the claim.
\end{proofsect}

Let $R_{\eff,h}(u,v)$ mark the explicit $h$-dependence of the effective resistance from~$u$ to~$v$ in a network with conductances related to a field~$h$ as in \eqref{E:12.13}. In order to control the effect of the rough part of the decomposition of the DGFF from Lemma~\ref{lemma-split}, we will also need:

\begin{mylemma}
\label{lemma-remove}
For any fields~$\varphi$ and~$\chi$,
\begin{equation}
\label{E:13.38}
R_{\eff,\varphi+\chi}(u,v)\le R_{\eff,\varphi}(u,v)\max_{(x,y)\in\cmss E}\texte^{-\beta(\chi_x+\chi_y)}.
\end{equation}
Moreover, if~$\varphi\independent\chi$ then also
\begin{equation}
\label{E:13.39}
\E\bigl(R_{\eff,\varphi+\chi}(u,v)\,\big|\,\varphi\bigr)\le 
R_{\eff,\varphi}(u,v)\max_{(x,y)\in\cmss E}\E\bigl(\texte^{-\beta(\chi_x+\chi_y)}\bigr).
\end{equation}
\end{mylemma}

\begin{proofsect}{Proof}
Let~$i\in\II(u,v)$. Then \eqref{E:12.27} and \eqref{E:12.13} yield
\begin{equation}
R_{\eff,\varphi+\chi}(u,v)\le\sum_{e=(x,y)\in\cmss E}\texte^{-\beta(\varphi_x+\varphi_y)}\texte^{-\beta(\chi_x+\chi_y)}i(e)^2.
\end{equation}
Bounding the second exponential by its maximum and optimizing over~$i$ yields \eqref{E:13.38}. For \eqref{E:13.39} we first take the conditional expectation given~$\varphi$ and then proceed as before.
\end{proofsect}

Lemma~\ref{lemma-remove} will invariably be used through the following bound:

\begin{myexercise}[Removal of independent field]
\label{ex:13.15}
Suppose~$\varphi$ and~$\chi$ are independent Gaussian fields on~$\cmss V$. Denote $\sigma^2:=\max_{x\in\cmss V}\Var(\chi_x)$. Show that then for each~$a>0$ and each~$r>0$,
\begin{equation}
\Bigl|\BbbP\bigl(\,R_{\eff,\varphi+\chi}(u,v)\le ar\bigr)- \BbbP\bigl(\,R_{\eff,\varphi}(u,v)\le r\bigr)\Bigr|\le a^{-1}\texte^{2\beta^2\sigma^2}.
\end{equation}
\end{myexercise}

\noindent
We are ready to give:

\begin{proofsect}{Proof of Proposition~\ref{prop-13.9}}
Lemma~\ref{lemma-split} permits us to realize $h^{B(M)}$ as the sum of independent fields~$\varphi$ and~$\chi$ where, by the union bound and the standard Gaussian tail estimate (see Exercise~\ref{ex:2.2}), for each~$\epsilon>0$ there is~$c_1>0$ such that
\begin{equation}
\sup_{M\ge1}\,\BbbP\Bigl(\,\max_{\begin{subarray}{c}
x,y\in B(M/2)\\|x-y|_1\le 2
\end{subarray}}
|\varphi_x-\varphi_y|>c_1\Bigr)<\epsilon.
\end{equation}
Hence,~$\rhomax$ associated with field~$\varphi$ in the box $B(N)\subseteq B(M/2)$ via \eqref{E:12.13} obeys $\BbbP(\rhomax\le\texte^{c_1\beta})\ge1-\epsilon$. Next we observe that, abbreviating~$S:=B(N)$, the pairs
\begin{equation}
(A,B):=(\partial^{\text{left}}S,\partial^{\text{right}}S)\quad\text{and}\quad
(C,D):=(\partial^{\text{up}}S,\partial^{\text{down}}S)
\end{equation}
 obey the conditions in Lemma~\ref{lemma-dual2}. Since~$\mathfrak d=4$, \eqref{E:13.25uia} and \eqref{E:13.33} give
\begin{equation}
\BbbP\Bigl(R^\LR_{B(N)}(\varphi)R^\UD_{B(N)}(-\varphi)\le 64\texte^{c_1\beta}\Bigr)\ge1-\epsilon.
\end{equation}
The rotational symmetry of $B(N)$, $B(M)$ and~$\varphi$ along with the distributional symmetry~$\varphi\,\laweq\,-\varphi$ imply
\begin{equation}
R^\UD_{B(N)}(-\varphi)\,\,\laweq\,\,R^\LR_{B(N)}(\varphi).
\end{equation}
 A union bound then gives
\begin{equation}
\BbbP\bigl(R^\LR_{B(N)}(\varphi)\le 8\texte^{c_1\beta/2}\bigr)\ge\frac{1-\epsilon}2.
\end{equation}
Exercise~\ref{ex:13.15} with $r:=8\texte^{c_1\beta/2}$ and $a:=\texte^{\hat c\log\log M}(8\texte^{c_1\beta/2})^{-1}$ for any choice of $\hat c>2\beta^2 c$ in conjunction with Lemma~\ref{lemma-split} then imply the claim for~$N$ sufficiently large.\end{proofsect}

\section{RSW theory for effective resistance}
\noindent
The next (and the most challenging) task is to elevate the statement of Proposition~\ref{prop-13.9} from squares to rectangles. The desired claim is the content of:

\begin{myproposition}[Effective resistance across rectangles]
\label{prop-13.15}
There are $c,C\in(0,\infty)$ and~$N_1\ge1$ such that for all $N\ge N_1$, all $M\ge 16 N$ and any translate~$S$ of~$B(4N,N)$ satisfying~$S\subset B(M/2)$,
\begin{equation}
\label{E:13.46}
\BbbP\Bigl(R^\LR_{S}(h^{B(M)})\le C\texte^{\hat c\log\log M}\Bigr)\ge c.
\end{equation}
The same holds for $R^\UD_{S}(h^{B(M)})$ and translates~$S$ of $B(N,4N)$ with~$S\subset B(M/2)$. (The constant~$\hat c$ is as in Proposition~\ref{prop-13.9}.)
\end{myproposition}

The  proof (which we will only sketch) is rooted in the Russo-Seymour-Welsh (RSW) argument for critical percolation, which is a way to to bootstrap uniform lower bounds on the probability of an occupied crossing for squares to those for an occupied crossing for rectangles (in the ``longer'' direction) of a given aspect ratio. The technique was initiated in Russo~\cite{Russo78}, Seymour and Welsh~\cite{SW78} and Russo~\cite{Russo81} for percolation and later adapted to dependent models as well (e.g., Duminil-Copin, Hongler and Nolin~\cite{DHN11}, Beffara and Duminil-Copin~\cite{BDC12}). We will follow the version that Tassion~\cite{Tassion} developed for Voronoi percolation.

Proposition~\ref{prop-12.20} links the effective resistance to crossings by collections of paths. In order to mimic arguments from percolation theory, we will need a substitute for the trivial geometric fact that, if a path from~$u$ to~$v$ crosses a path from~$u'$ to~$v'$, then the union of these paths contains a path from~$u$ to~$v'$. For this, given a set~$\AA$ of paths, let~$R_\eff(\AA)$ denote quantity in \eqref{E:12.52} with the first infimum restricted to $\PP\subseteq\AA$ (not necessarily a subset of~$\mathfrak P_{u,v}$). We then have:

\begin{mylemma}[Subadditivity of effective resistance]
\label{lemma-13.16}
Let $\AA_1,\dots,\AA_n$ be sets of paths such that for each selection~$P_i\in\AA_i$, $i=1,\dots,n$, the graph union $P_1\cup\dots\cup P_n$ contains a path from~$u$ to~$v$. Then
\begin{equation}
R_\eff(u,v)\le\sum_{i=1}^n R_\eff(\AA_i).
\end{equation}
\end{mylemma}

Another property we will need can be thought of as a generalization of the Parallel Law: If we route current from~$u$ to~$v$ along separate families of paths, then $C_\eff(u,v)$ is at most the sum of the conductances of the individual families. This yields:

\begin{mylemma}[Subadditivity of effective conductance]
\label{lemma-13.17}
Let $\AA_1,\dots,\AA_n$ be sets of paths such that every path from~$u$ to~$v$ lies in $\AA_1\cup\dots\cup\AA_n$. Then
\begin{equation}
C_\eff(u,v)\le\sum_{i=1}^n R_\eff(\AA_i)^{-1}.
\end{equation}
\end{mylemma}

\noindent
We will not supply proofs of these lemmas as that amounts to further variations on the calculations underlying Propositions~\ref{prop-12.20} and~\ref{prop-12.22}; instead, we refer the reader to~\cite{BDG}. The use of these lemmas will be facilitated by the following observation:

\begin{myexercise}
\label{ex:13.18}
Show that, for any collection of paths~$\AA$ and for resistances related to field~$h$ as in \eqref{E:12.13}, $h\mapsto R_{\eff,h}(\AA)$ is decreasing (in each coordinate). Prove that, for~$h$ given as DGFF  and each~$a>0$, under the setting of Lemma~\ref{lemma-13.16} we have
\begin{equation}
\label{E:13.49a}
\BbbP\bigl(R_\eff(u,v)\le a\bigr)\ge \prod_{i=1}^n \BbbP\bigl(R_\eff(\AA_i)\le a/n\bigr)
\end{equation}
while under the setting of Lemma~\ref{lemma-13.17} we have
\begin{equation}
\label{E:13.50a}
\max_{i=1,\dots,n} \BbbP\bigl(R_\eff(\AA_i)\le a\bigr)\ge1-\Bigl[1-\BbbP\bigl(R_\eff(u,v)\le a/n\bigr)\Bigr]^{1/n}\,.
\end{equation}
Hint: Use the FKG inequality.
\end{myexercise}

\nopagebreak
\begin{figure}[t]
\vglue-1mm
\centerline{\includegraphics[width=0.45\textwidth]{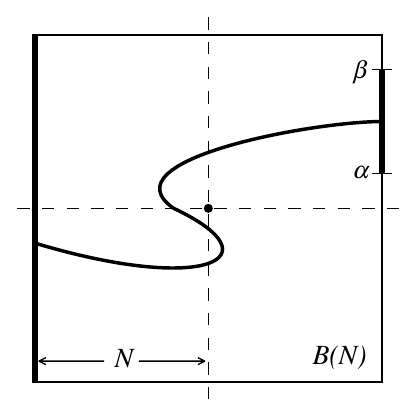}
}

\vglue0mm
\begin{quote}
\small 
\vglue-0.2cm
\caption{
\label{fig-alpha-beta}
\small
A path in the set~$\AA_{N,[\alpha,\beta]}$ underlying the definition of the effective resistivity $R^\LR_{N,[\alpha,\beta]}(h)$ between the left side of the square and the portion of the right-side marked by the interval~$[\alpha,\beta]$.}
\normalsize
\end{quote}
\end{figure} 

Tassion's version of the RSW argument proceeds by considering crossings of~$B(N)$ from the left-side of~$B(N)$ to (only) a \emph{portion} of the right side --- namely that corresponding to the interval~$[\alpha,\beta]$; see Fig.~\ref{fig-alpha-beta} for the geometric setting.
Writing~$\AA_{N,[\alpha,\beta]}$ for the set of all such crossing paths, we abbreviate
\begin{equation}
R^\LR_{N,[\alpha,\beta]}(h):=R_{\eff,h}(\AA_{N,[\alpha,\beta]}).
\end{equation}
We then have:

\begin{mylemma}
\label{lemma-14.20nwt}
For each~$\epsilon>0$ and all $M\ge 2N$ with~$N$ sufficiently large,
\begin{equation}
\label{E:13.50}
\BbbP\bigl(R^\LR_{N,[0,N]}(h^{B(M)})\le 2\texte^{\hat c\log\log M}\bigr)\ge1-\frac1{\sqrt 2}-\epsilon.
\end{equation}
\end{mylemma}

\begin{proofsect}{Proof}
Lemma~\ref{lemma-13.17} and the Electrostatic Duality show
\begin{equation}
\frac1{R^{\LR}_{B(N)}(h)}\le\frac1{R^\LR_{N,[0,N]}(h)}+\frac1{R^\LR_{N,[-N,0]}(h)}
\end{equation}
and the rotation (or reflection) symmetry yields
\begin{equation}
R^\LR_{N,[0,N]}(h^{B(M)})\,\laweq\,R^\LR_{N,[-N,0]}(h^{B(M)}).
\end{equation}
The bound \eqref{E:13.50a} along with Proposition~\ref{prop-13.9} then show \eqref{E:13.50}.
\end{proofsect}

\begin{figure}[t]
\vglue-2mm
\centerline{\includegraphics[width=0.75\textwidth]{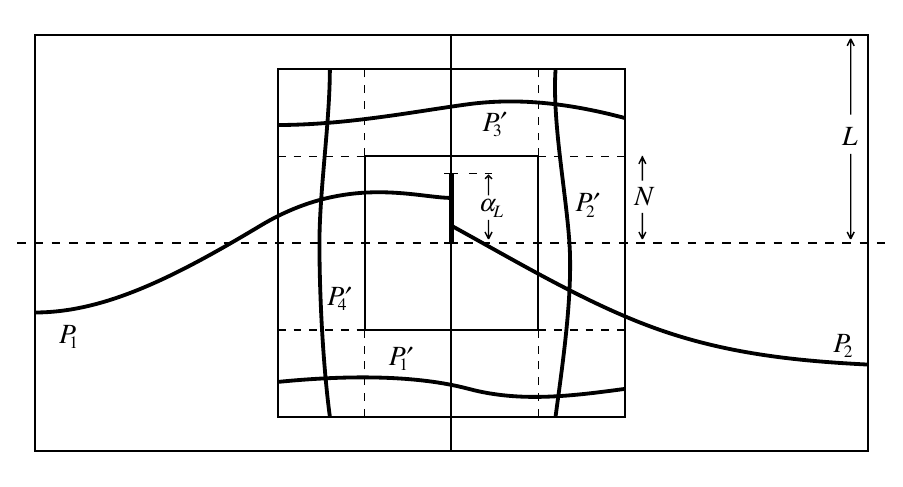}
}

\vglue-1mm
\begin{quote}
\small 
\vglue-0.2cm
\caption{
\label{fig-LN}
\small
A collection of paths that enforce a left-to-right crossing of $B(4L,L)$. The paths $P_1',\dots,P_4'$ arise from our assumption that \eqref{E:13.46} holds for~$N$ and the corresponding quantity $R^{\LR}_{B(4N,N)}$ is moderate. The paths $P_1$ and~$P_2$ arise from the assumption $\alpha_L\le N$. Note that since $\lfloor L/2\rfloor >N\ge\alpha_L$, we know that also $R^{\LR}_{L,[0,\alpha_L]}$ is moderate.}
\normalsize
\end{quote}
\end{figure}

Equipped with these techniques, we are ready to give:

\begin{proofsect}{Proof of Proposition~\ref{prop-13.15} (rough sketch)}
Since $\alpha\mapsto R^\LR_{N,[\alpha,N]}(h)$ is increasing on~$[0,N]$ with the value at~$\alpha=N$ presumably quite large, one can identify (with some degree of arbitrariness) a value~$\alpha_N$ where this function first exceeds a large multiple of~$\texte^{\hat c\log\log(2N)}$ with high-enough probability. More precisely, set
\begin{equation}
\phi_N(\alpha):=\BbbP\Bigl(R^\LR_{N,[\alpha,N]}(h^{B(2N)})>C\texte^{\hat c\log\log(2N)}\Bigr)
\end{equation}
for a suitably chosen~$C>2$ and note that,  by \eqref{E:13.50},~$\phi_N(0)<0.99$. Then set
\begin{equation}
\alpha_N:=\lfloor N/2\rfloor\wedge\min\Bigl\{\alpha\in\{0,\dots,\lfloor N/2\rfloor\}\colon \phi_N(\alpha)>0.99\Bigr\}\,.
\end{equation}
We now treat separately the cases~$\alpha_N=\lfloor N/2\rfloor$ and~$\alpha_N< \lfloor N/2\rfloor$ using the following sketch of the actual proof, for which we refer to the original paper:

\nopagebreak
\begin{figure}[t]
\vglue-1mm
\centerline{\includegraphics[width=0.95\textwidth]{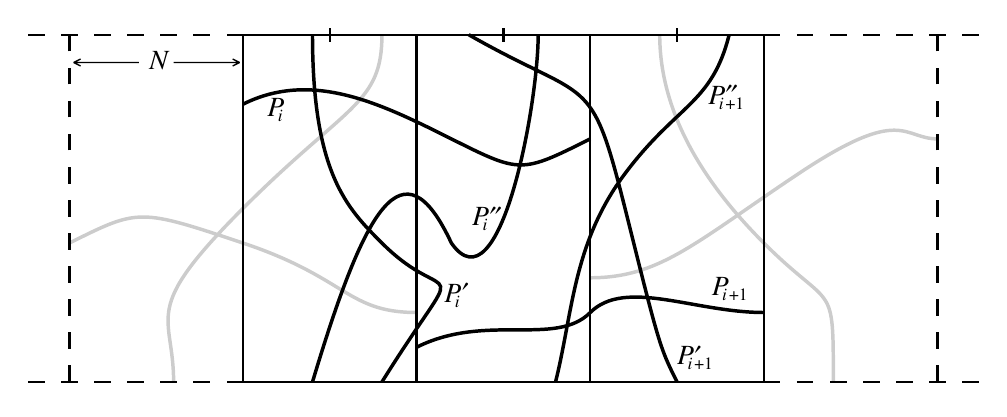}
}

\vglue0mm
\begin{quote}
\small 
\vglue-0.2cm
\caption{
\label{fig-5x2}
\small
An illustrations of paths~$P_i$, $P_i'$ and~$P_i''$, all of the form in Fig.~\ref{fig-alpha-beta} with~$\alpha:=\lfloor N/2\rfloor$ and~$\beta:=N$ (or just left-to-right crossing for~$P_i$) that ensure the existence a left-right crossing of the $4N\times N$ rectangle when~$\alpha_N=\lfloor N/2\rfloor$. The effective resistance in the rectangle is then bounded by the sum of (suitably shifted) copies of $R^\LR_{N,[0,N]}(h)$.}
\normalsize
\end{quote}
\end{figure}

\smallskip\noindent
\textit{CASE $\alpha_N=\lfloor N/2\rfloor$}: Here the fact that $\phi_N(\alpha_N-1)\le 0.99$ implies (via perturbation arguments that we suppress) that $\phi_N(\alpha_N)<0.999$ and so
\begin{equation}
\label{E:13.50uo}
\BbbP\bigl(R^\LR_{N,[\alpha_N,N]}\le C\texte^{\hat c\log\log(2N)}\bigr)\ge0.001.
\end{equation}
Invoking suitable shifts of the box~$B(N)$, Lemma~\ref{lemma-13.16} permits us to bound the left-to-right effective resistance $R^{\LR}_{B(4N,N)}$ by the sum of seven copies of effective resistances of the form $R^\LR_{N,[\alpha_N,N]}$ (and rotations thereof) and four copies of~$R^\LR_{B(N)}$ (and rotations thereof); see Fig.~\ref{fig-5x2}. The inequalities \eqref{E:13.49a}, \eqref{E:13.50uo} and \eqref{E:13.24ai} then yield \eqref{E:13.46}. A caveat is that these squares/rectangles are (generally) not centered at the same point as~$B(M)$, which we need in order to apply Proposition~\ref{prop-13.9} and the definition of~$\alpha_N$. This is remedied by invoking the Gibbs-Markov decomposition and removing the binding field via Exercise~\ref{ex:13.15}.

\smallskip\noindent
\textit{CASE $\alpha_N<\lfloor N/2\rfloor$}: In this case $\phi_N(\alpha_N)>0.99$. Since, by Lemma~\ref{lemma-13.17}, 
\begin{equation}
\frac1{R^\LR_{N,[0,N]}}\le \frac1{R^\LR_{N,[\alpha_N,N]}}+\frac1{R^\LR_{N,[0,\alpha_N]}}
\end{equation}
Lemma~\ref{lemma-14.20nwt} (and~$C>2$) show that $R^\LR_{N,[0,\alpha_N]}\le C'\texte^{\hat c\log\log(2N)}$ with a uniformly positive probability. Assuming in addition that 
\begin{equation}
\label{E:13.51}
\alpha_N\le 2\alpha_L\quad\text{for}\quad L:=\lfloor 4N/7\rfloor
\end{equation}
 another path crossing argument (see Fig.~\ref{fig-3x2}) is used to ensure a left-to-right crossing of the rectangle $B(2N-L,N)$. This is then readily bootstrapped to the crossing of~$B(4N,N)$, and thus a bound on the effective resistance, by way of Lemma~\ref{lemma-13.16} and the inequality \eqref{E:13.49a}.

\nopagebreak
\begin{figure}[t]
\vglue-1mm
\centerline{\includegraphics[width=0.6\textwidth]{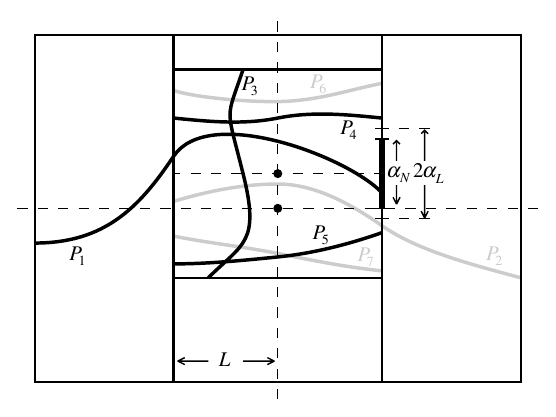}
}

\vglue0mm
\begin{quote}
\small 
\vglue-0.2cm
\caption{
\label{fig-3x2}
\small
Arranging a left-to-right crossing of $B(2N-L,N)$, with $L:=\lfloor 4N/7\rfloor$, under the assumption that~$\alpha_N\le2\alpha_L$. Some of the paths have been drawn in gray to make parsing of the connections easier.}
\normalsize
\end{quote}
\end{figure}

It remains to validate the assumption \eqref{E:13.51}. Here one proceeds by induction assuming that the statement \eqref{E:13.46} already holds for~$N$ and proving that, if~$\alpha_L\le N$, then the statement holds (with slightly worse constants) for~$L$ as well. This is based on a path-crossing argument whose geometric setting is depicted in Fig.~\ref{fig-LN}. This permits the construction of a sequence~$\{N_k\colon k\ge1\}$ such that $\alpha_{N_{k+1}}\le N_k$ and such that $N_{k+1}/N_k$ is bounded as~$k\to\infty$. (Another path crossing argument is required here for which we refer the reader to the original paper.) The claim is thus proved for~$N\in\{N_k\colon k\ge1\}$; general~$N$ are handled by monotonicity considerations.
\end{proofsect}

\nopagebreak

\section{Upper tail of effective resistance}
\noindent
With Proposition~\ref{prop-13.15} in hand, we are now ready to give the proof of:

\begin{mytheorem}[Effective resistance, upper tail]
\label{thm-13.19}
Given integers~$M\ge N\ge1$, let $R_\eff(u,v)$ denote the effective resistance from~$u$ to~$v$ in the network on~$B(N)$ with conductances related to~$h=h^{B(M)}$ via \eqref{E:12.13}. There are~$c_1,c_2\in(0,\infty)$ such that
\begin{equation}
\label{E:13.57}
\max_{u,v\in B(N)}\,\BbbP\Bigl(R_{\eff}(u,v)\ge c_1(\log M)\texte^{t\sqrt{\log M}}\Bigr)\le c_1(\log M)\texte^{-c_2 t^2}
\end{equation}
holds for all~$N\ge1$, all $M\ge 32 N$ and all~$t\ge1$.
\end{mytheorem}

\noindent
As we shall see, the proof is based on the following concentration estimate for the effective resistance:

\begin{mylemma}
\label{lemma-13.20}
For any collection of paths~$\AA$, let~$f(h):=\log R_{\eff,h}(\AA)$ (which entails that the resistances depend on~$h$ via \eqref{E:12.13}). Then
\begin{equation}
\sup_h\sum_x\Bigl|\frac{\partial f}{\partial h_x}(h)\Bigr|\le 2\beta.
\end{equation}
\end{mylemma}

\begin{proofsect}{Proof}
We will prove this under the simplifying assumption that~$\AA$ is the set of all paths from~$u$ to~$v$ and, therefore, $R_{\eff,h}(\AA)=R_{\eff,h}(u,v)$. Let~$i_\star$ be the current realizing the minimum in \eqref{E:12.27}. Then \eqref{E:12.13} yields
\begin{equation}
\begin{aligned}
\frac\partial{\partial h_x}\log R_{\eff,h}(u,v)
&=\frac1{R_{\eff,h}(u,v)}\sum_{e\in\cmss E}i_\star(e)^2\frac{\partial}{\partial h_x}r_e
\\
&=-\frac\beta{R_{\eff,h}(u,v)}\sum_{y\colon(x,y)\in\cmss E}i_\star(x,y)^2 \,r_e\,.
\end{aligned}
\end{equation}
Summing the last sum over all~$x$ yields $2R_{\eff,h}(u,v)$. The claim thus follows.
\end{proofsect}

\begin{myexercise}
Use the variational representation of~$R_\eff(\AA)$ as in \eqref{E:12.52} to prove Lem\-ma~\ref{lemma-13.20} in full generality.
\end{myexercise}

As a consequence we get:

\begin{mycorollary}
\label{cor-13.21}
There is a constant $c>0$ such that, for any $M\ge1$, any collection~$\AA$ of paths in~$B(M)$ and any~$t\ge0$,
\begin{equation}
\label{E:13.60}
\BbbP\Bigl(\bigl|\log R_{\eff,h^{B(M)}}(\AA)-\E\log R_{\eff,h^{B(M)}}(\AA)\bigr|>t\sqrt{\log M}\Bigr)\le2\texte^{-c t^2}.
\end{equation}
\end{mycorollary}

\begin{proofsect}{Proof}
Lemma~\ref{lemma-13.20} implies the premise \eqref{E:6.15a} (with~$M:=2\beta$) of the general Gaussian concentration estimate in Corollary~\ref{cor-6.6}. The claim follows by the union bound and the uniform control of the variance of the DGFF (cf, e.g., \eqref{E:2.3uai}).
\end{proofsect}

In order to localize the value of the expectation in \eqref{E:13.60}, we also need:

\begin{mylemma}[Effective resistance across rectangles, lower bound]
\label{lemma-13.22}
There are constants $c',C'>0$ such that for all~$N\ge1$, all $M\ge 16 N$ and any translate~$S$ of~$B(4N,N)$ with~$S\subset B(M/2)$,
\begin{equation}
\label{E:13.61}
\BbbP\Bigl(R^\LR_{S}(h^{B(M)})\ge C'\texte^{-\hat c\log\log M}\Bigr)\ge c'.
\end{equation}
The same holds for $R^\UD_{S}(h^{B(M)})$ and translates~$S$ of $B(N,4N)$ with~$S\subset B(M/2)$. (The constant~$\hat c$ is as in Proposition~\ref{prop-13.9}.)
\end{mylemma}

\begin{proofsect}{Proof (sketch)}
By swapping resistances for conductances (and relying on Lem\-ma~\ref{lemma-dual} instead of Lemma~\ref{lemma-dual2}), we find out that \eqref{E:13.61} holds for~$S:=B(N)$. Invoking the Gibbs-Markov property, Exercise~\ref{ex:13.15} yields a similar bound for all translates of~$B(N)$ contained in $B(M/4)$. Since every translate $S$ of~$B(4N,N)$ with~$S\subset B(M/2)$ contains a translate~$S'$ of~$B(N)$ contained in~$B(M/4)$, the claim  follows from the fact that $R^\LR_{S'}\le R^\LR_{S}$.
\end{proofsect}

\nopagebreak
\begin{figure}[t]
\vglue0mm
\centerline{\includegraphics[width=0.65\textwidth]{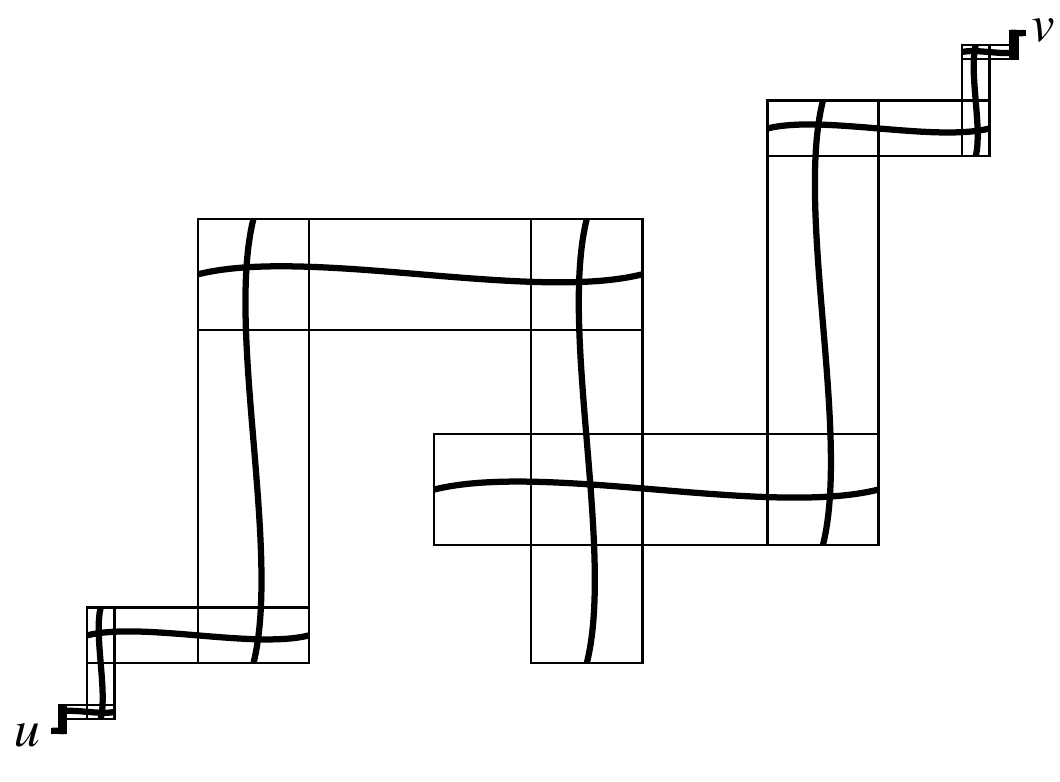}
}
\vglue-2mm
\begin{quote}
\small 
\vglue-0.2cm
\caption{
\label{fig-log-boxes}
\small
An illustration of a sequence of $4K\times K$ rectangles (for~$K$ taking values in powers of~$2$) such that the occurrence of a crossing of each of these rectangles implies the existence of a path from~$u$ to~$v$. For vertices within distance~$N$, order-$\log N$ rectangles suffice.}
\normalsize
\end{quote}
\end{figure}

With Lemma~\ref{lemma-13.22} in hand, we can finally prove the upper bound on the point-to-point effective resistances:

\begin{proofsect}{Proof of Theorem~\ref{thm-13.19}}
Using Corollary~\ref{cor-13.21} along with Proposition~\ref{prop-13.15} and Lemma~\ref{lemma-13.22} (and the fact that $\log\log M\le\sqrt{\log M}$ for~$M$ large) we get that
\begin{equation}
E\bigl|\log R^\LR_{S}(h^{B(M)})\bigr|\le \tilde c\sqrt{\log M}
\end{equation}
and
\begin{equation}
\label{E:13.63}
\BbbP\Bigl(\log R^\LR_{S}(h^{B(M)})\ge \tilde c(1+t)\sqrt{\log M}\Bigr)\le\texte^{-\tilde c' t^2},\quad t\ge0,
\end{equation}
hold with the same constants $\tilde c,\tilde c'\in(0,\infty)$ for all translates~$S$ of~$B(4N,N)$ or $B(N,4N)$ contained in~$B(M/2)$, uniformly in~$N\ge1$ and~$M\ge 32 N$.

Next we invoke the following geometric fact: Given $u,v\in B(N)$, there are at most order-$\log N$ translates of $B(4K,K)$ or~$B(K,4K)$ with $K\in\{1,\dots,N\}$ such that the existence of a crossing (in the longer direction) in each of these rectangles implies the existence of a path from~$u$ to~$v$; see Fig.~\ref{fig-log-boxes}. Assuming that the corresponding resistance (left-to-right in horizontal rectangles and top-to-bottom in vertical ones) is at most $\texte^{\tilde c(1+t)\sqrt{\log M}}$, Lemma~\ref{lemma-13.16} bounds $R_{\eff}(u,v)$ by a constant times $(\log M)\texte^{\tilde c(1+t)\sqrt{\log M}}$. The claim then follows from \eqref{E:13.63} and the union bound.
\end{proofsect}

Theorem~\ref{thm-13.19} now in turn implies the upper bound on the growth rate of resistances from~$0$ to~$B(N)^\cc$ in an infinite network on~$\Z^2$:

\begin{proofsect}{Proof of \eqref{E:1.11ua} in Theorem~\ref{thm-12.5} (sketch)}
First note that $R_\eff(0,B(N)^\cc)\le R_\eff(0,v)$ for any vertex~$v\in\partial B(N)$. The Gibbs-Markov property along with Exercise~\ref{ex:13.15} allow us to estimate the tail of $R_\eff(0,B(N)^\cc)$ for $h:=$ DGFF in~$\Z^2\smallsetminus\{0\}$ by essentially that for~$h:=$ DGFF in~$B(2N)$. Theorem~\ref{thm-13.19} gives\begin{equation}
\BbbP\Bigl(\log R_\eff(0,B(N)^\cc)\ge a(\log N)^{1/2}(\log\log N)^{1/2}\Bigr)\le c(\log N)^{1-c_2a^2}
\end{equation}
for all~$a>0$ as soon as~$N$ is sufficiently large.
Once~$a$ is sufficiently large, this is summable for~$N\in\{2^n\colon n\in\N\}$. The claim thus follows from the Borel-Cantelli lemma and the monotonicity of $N\mapsto R_\eff(0,B(N)^\cc)$.
\end{proofsect}

The lower corresponding lower bound \eqref{E:1.12ua} requires the use of the concentric decomposition and is therefore deferred to the next lecture.


\chapter{From resistance to random walk}
\label{lec-15}\noindent
We are left with the task to apply the technology for controlling the resistance to establish the desired conclusions about the random walk in DGFF landscape. It is here where we will tie the random walk problem to the overall theme of these lectures: the large values of the DGFF. We commence by some elementary identities from Markov chain theory which readily yield upper bounds on the expected exit time and the heat kernel. The control of the corresponding lower bounds is considerably more involved. We again focus on conveying the main ideas while referring the reader to the original paper for specific technical details.

\section{Hitting and commute-time identities}
\label{sec14.1}\noindent
We begin by considerations that apply to a general discrete-time Markov chain~$X$ on a countable state space~$\cmss V$ with transition probability~$\cmss P$ and a reversible measure~$\pi$. This problem always admits the electric network formulation by setting the conductances to $c(x,y):=\pi(x)\cmss P(x,y)$ --- whose symmetry is equivalent to reversibility of~$\pi$ --- and defining the resistances accordingly. The set of edges~$\cmss E$ is then the set of pairs~$(x,y)$ where~$c(x,y)>0$. Our main tool will be:

\begin{mylemma}[Hitting time identity]
\label{lemma-14.1}
Let~$A\subset\cmss V$ be finite and let us denote, as before, $\tau_{A^\cc}:=\inf\{n\ge0\colon X_n\not\in A\}$. Then for each~$x\in A$,
\begin{equation}
\label{E:14.1}
E^x(\tau_{A^\cc})=R_\eff(x,A^\cc)\sum_{y\in A}\pi(y)\phi(y)\,,
\end{equation}
where~$\phi(y):=P^x(\tau_x<\tau_{A^\cc})$ has the interpretation of the potential (a.k.a.\ voltage) that minimizes the variational problem defining $C_\eff(x,A^\cc)$. 
\end{mylemma}

\begin{proofsect}{Proof}
The Markov property ensures that the Green function associated with the Markov chain, defined by the formula \eqref{E:1.3au} or via~$G^A:=(1-\cmss P)^{-1}$, is $\cmss P$-harmonic in the first coordinate. The uniqueness of the solution to the Dirichlet problem in finite sets (see Exercise~\ref{ex:max-principle}) implies 
\begin{equation}
\phi(y)=\frac{G^A(y,x)}{G^A(x,x)}.
\end{equation}
Reversibility then yields
\begin{equation}
\begin{aligned}
\sum_{y\in A}\pi(y)\phi(y) &= G^A(x,x)^{-1}\sum_{y\in A}\pi(y)G^A(y,x)
\\&
=G^A(x,x)^{-1}\pi(x)\sum_{y\in A}G^A(x,y) 
\\
&= G^A(x,x)^{-1}\pi(x)E^x(\tau_{A^\cc}).
\end{aligned}
\end{equation}
The claim follows from
\begin{equation}
G^A(x,x)=\pi(x) R_\eff(x,A^\cc)
\end{equation}
as implied by Corollary~\ref{cor-12.14} and \eqref{E:1.15wk}.
\end{proofsect}

The above lemma is well known particularly in the form of the following corollaries. The first one of these is folklore:

\begin{mycorollary}[Hitting time estimate]
\label{cor-14.2}
For all finite~$A\subset\cmss V$ and all~$x\in A$,
\begin{equation}
E^x(\tau_{A^\cc})\le R_\eff(x,A^\cc)\pi(A).
\end{equation}
\end{mycorollary}

\begin{proofsect}{Proof}
Apply $\phi(y)\le 1$ for all~$y\in A$ in \eqref{E:14.1}.
\end{proofsect}

\noindent
The second corollary, which we include mostly for completeness of exposition, appeared first in~Chandra~\textit{et al}~\cite{CRRST}:

\begin{mycorollary}[Commute-time identity]
\label{cor-commute}
For~all $\cmss V$ finite and all distinct vertices $u,v\in\cmss V$,
\begin{equation}
E^u(\tau_v)+E^v(\tau_u) = R_\eff(u,v)\pi(\cmss V).
\end{equation}
\end{mycorollary}

\begin{proofsect}{Proof}
Let~$\phi$ be the minimizer of~$C_\eff(u,v)$ and~$\phi'$ the minimizer of~$C_\eff(v,u)$. Noting that $\phi(x)+\phi'(x)=1$ for all~$x=u,v$ the uniqueness of the solution to the Dirichlet problem implies that $\phi(x)+\phi'(x)=1$ all~$x\in\cmss V$. The claim follows from \eqref{E:14.1}. 
\end{proofsect}

\section{Upper bounds on expected exit time and heat kernel}
\noindent
The upshot of the above corollaries is quite clear: In order to bound the expected hitting time from above, we only need tight upper bounds on the effective resistance and the total volume of the reversible measure. Drawing heavily on our earlier work in these notes, this now permits to give:

\begin{proofsect}{Proof of Theorem~\ref{thm-exittime}, upper bound in \eqref{E:12.8}}
Let~$\delta>0$ and let~$h$ be the DGFF in $\Z^2\smallsetminus\{0\}$. The upper bound in Theorem~\ref{thm-12.5} ensures
\begin{equation}
R_{\eff,h}\bigl(0,B(N)^\cc\bigr)\le\texte^{(\log N)^{1/2+\delta}}
\end{equation}
with probability tending to one as~$N\to\infty$. The analogue of \eqref{E:2.17u} for the DGFF in~$\Z^2\smallsetminus\{0\}$ along with Markov inequality and Borel-Cantelli lemma yields
 \begin{equation}
\label{E:14.7}
\pi_h\bigl(B(N)\bigr)\le N^{\theta(\beta)}\texte^{(\log N)^\delta}
\end{equation}
with probability tending to one as~$N\to\infty$ when~$\beta\le\tilde\beta_\cc$; for~$\beta>\tilde\beta_\cc$ we instead use the Gibbs-Markov property and the exponential tails of the maximum proved in Lemma~\ref{lemma-8.2}. Corollary~\ref{cor-14.2} then gives
\begin{equation}
E^0_h(\tau_{B(N)^\cc})\le N^{\theta(\beta)}\texte^{2(\log N)^{1/2+2\delta}}
\end{equation}
with probability tending to one as $N\to\infty$.
\end{proofsect}

Next we will address the decay of the heat kernel:

\begin{proofsect}{Proof of Theorem~\ref{thm-heatkernel}, upper bound in~\eqref{E:12.5a}}
For~$T>0$ define the set
\begin{equation}
\Xi_T:=\{0\}\cup\Bigl\{x\in B(2T)\smallsetminus B(T)\colon \wt R_\eff(0,x)\le\texte^{(\log T)^{1/2+\delta}}\Bigr\},
\end{equation}
where $\wt R_\eff$ stands for the effective resistance in the network on~$B(2T)$.
From Theorem~\ref{thm-13.19} and the Markov inequality we conclude that~$|\Xi_T|$ will contain an overwhelming fraction of all vertices in the annulus $B(2T)\smallsetminus B(T)$. Let~$\wt X$ be the Markov chain on~$B(4T)$ defined by the same conductances as~$X$ except those corresponding to the jumps out of~$B(4T)$ which are set to zero (and these jumps are thus suppressed) and write~$\wt\pi_h$ for the correspondingly modified measure~$\pi_h$. Let~$Y_k$ be the position of the~$k$-th visit of~$\wt X$ to~$\Xi_T$, set~$\tau_0:=0$ and let $\tau_1$, $\tau_2$, etc be times of the successive visits of~$\wt X$ to~$0$. Let
\begin{equation}
\hat\sigma:=\inf\bigl\{k\ge1\colon \tau_k\ge T\bigr\}.
\end{equation} 
We ask the reader to verify:

\begin{myexercise}
Prove that $Y:=\{Y_k\colon k\ge0\}$ is a Markov chain on~$\Xi_T$ with stationary distribution
\begin{equation}
\nu(x):=\frac{\wt \pi_h(x)}{\wt\pi_h(\Xi_T)}.
\end{equation}
Prove also that~$\hat\sigma$ is a stopping time for the natural filtration associated with~$Y$ and that $E^x_h(\hat\sigma)<\infty$ (and thus $\hat\sigma<\infty$ $P^x_h$-a.s) hold for each~$x\in\Xi_T$. 
\end{myexercise}

\noindent
We also recall an exercise from a general theory of reversible Markov chains:

\begin{myexercise}
Prove that $T\mapsto P^0(\wt X_{2T}=0)$ is non-increasing.
\end{myexercise} 

\noindent
This permits us to write, for~$T$ even (which is all what matters),
\begin{equation}
\begin{aligned}
\frac12 TP^0_h(\wt X_T=0)&\le E^0_h\Bigl(\,\sum_{n=0}^{T-1}\1_{\{\wt X_n=0\}}\Bigr)
\\
&\le E^0_h\Bigl(\,\sum_{k=0}^{T-1}\1_{\{Y_k=0\}}\Bigr)\le E^0_h\Bigl(\,\sum_{k=0}^{\hat\sigma-1}\1_{\{Y_k=0\}}\Bigr),
\end{aligned}
\end{equation}
where the second inequality relies on the fact that~$0\in\Xi_T$. Next we observe:

\begin{myexercise}
Given a Markov chain~$Y$ with stationary distribution~$\nu$ and a state~$x$, suppose~$\sigma$ is a stopping time for~$Y$ such that~$Y_\sigma=x$ a.s. Then for each~$y$
\begin{equation}
E^x\Bigl(\,\sum_{k=0}^{\sigma-1}\1_{\{Y_k=y\}}\Bigr)=E^x(\sigma)\nu(y).
\end{equation}
Hint: The left-hand side is, as a function of~$y$, a stationary measure for~$Y$.
\end{myexercise}

By conditioning on~$Y_T$ we now have
\begin{equation}
E^0_h(\hat\sigma)\le T+E^0_h\bigl(E^{\wt X_T}(\tau_1)\bigr)\le T+\max_{u\in\Xi_T} E^u_h(\tau_1)
\end{equation}
and the commute-time identity (Corollary~\ref{cor-commute}) and the definition of~$\Xi_T$ give
\begin{equation}
E^u_h(\tau_1)\le\wt \pi_h(\Xi_T)\wt R_\eff(u,0)\le \pi_h(\Xi_T)\texte^{(\log T)^{1/2+\delta}}.
\end{equation}
The nearest-neighbor nature of the walk permits us to couple~$X$ and~$\wt X$ so that~$X$ coincides with~$\wt X$ at least up to time~$4T$. The above then bounds $P^0_h(X_T=0)$ by~$2\pi_h(0)T^{-1}\texte^{(\log T)^{1/2+\delta}}$.
\end{proofsect}

We can now also settle:

\begin{proofsect}{Proof of Corollary~\ref{cor-12.4}}
Using the Markov property, reversibility and the Cauchy-Schwarz inequality we get
\begin{equation}
\label{E:14.16au}
\begin{aligned}
P_h^0(X_{2T}=0)&\ge\sum_{x\in B(N)}P_h^0(X_T=x)P_h^x(X_T=0)
\\
&=\pi_h(0)\sum_{x\in B(N)}\frac{P_h^0(X_T=x)^2}{\pi_h(x)}
\ge\pi_h(0)\frac{P^0_h\bigl(X_T\in B(N)\bigr)^2}{\pi_h\bigl(B(N)\bigr)}.
\end{aligned}
\end{equation}
Invoking the (already proved) upper bound on the heat-kernel and the bound \eqref{E:14.7} we get that, with high probability,
\begin{equation}
P^0_h\bigl(X_T\in B(N)\bigr)
\le \Bigl[\frac 1 T \texte^{(\log T)^{1/2+\delta}}\,N^{\theta(\beta)}\texte^{(\log N)^\delta}\Bigr]^{1/2}\,.
\end{equation}
Setting $T:=N^{\theta(\beta)}\texte^{(\log N)^{1/2+2\delta}}$ gives the desired claim.
\end{proofsect}

\section{Bounding voltage from below}
\noindent
We  now turn to the lower bounds in Theorems~\ref{thm-heatkernel} and~\ref{thm-exittime}. As an inspection of \eqref{E:14.1} shows, this could be achieved by proving a lower bound on the potential difference (a.k.a.\ voltage) minimizing~$C_\eff(0,B(N)^\cc)$. This will be based on:

\begin{mylemma}
In the notation of Lemma~\ref{lemma-14.1}, for all~$x\in A$ and~$y\in A\smallsetminus\{x\}$,
\begin{equation}
\label{E:14.18}
2R_\eff(x,A^\cc)\phi(y)= R_\eff(x,A^\cc)+R_\eff(y,A^\cc)-R_{A,\eff}(x,y)\,,
\end{equation}
where~$R_{A,\eff}$ denotes the effective resistance in the network with~$A^\cc$ collapsed to a point.
\end{mylemma}

\begin{proofsect}{Proof}
We will apply the network reduction principle from Exercise~\ref{ex:12.17} and reduce the problem to the network with three nodes $1$, $2$ and~$3$ corresponding to $x$, $y$ and~$A^\cc$, respectively. Since~$\phi$ is harmonic, so is its restriction to the reduced network. But~$\phi(y)$ also has the interpretation of the probability that the reduced Markov chain at~$2$ jumps to~$1$ before~$3$. Writing the conductance between the~$i$th and the~$j$th node as $c_{ij}$, this probability is  given as $\frac{c_{12}}{c_{12}+c_{23}}$. Exercise~\ref{ex:12.20} equates this to the ratio of the right-hand side of~\eqref{E:14.18} and~$2R_\eff(0,A^\cc)$.
\end{proofsect}

For the setting at hand, the quantity on the right-hand side of \eqref{E:14.18} becomes:
\begin{equation}
D_N(x):=R_{\eff}\bigl(0,B(N)^\cc\bigr)+R_{\eff}\bigl(x,B(N)^\cc\bigr)-R_{B(N),\eff}(0,x)\,.
\end{equation}
We claim:

\begin{myproposition}
\label{prop-14.7}
For any~$\delta\in(0,1)$
\begin{equation}
\lim_{N\to\infty}\,\BbbP\Bigl(\,\min_{x\in B(N\texte^{-(\log N)^\delta})} D_N(x)\ge\log N\Bigr)=1.
\end{equation}
\end{myproposition}

We will only give a sketch of the main idea which is also key for the proof of the lower bound in Theorem~\ref{thm-12.5} that we will prove in some more detail. The proof relies on the concentric decomposition of the pinned DGFF from Exercise~\ref{ex:8.13} which couples the field to the Gaussian random walk $\{S_n\colon n\ge0\}$ defined in \eqref{E:8.42a} and some other Gaussian fields. We claim:

\begin{myproposition}
\label{prop-14.8}
Let~$h:=$ DGFF in~$\Z^2\smallsetminus\{0\}$ and let~$\{S_n\colon n\ge1\}$ be the Gaussian random walk related to~$h$ via Exercise~\ref{ex:8.13}. There is a constant~$c>0$ such that
\begin{equation}
R_{\eff,h}\bigl(0,B(N)^\cc\bigr)\ge \max_{k=1,\dots,\lfloor\log_4N\rfloor}\,\texte^{2\beta S_k-c(\log\log N)^2}
\end{equation}
fails for at most finitely many~$N$'s, $\BbbP$-a.s.
\end{myproposition}

For the proof, recall the definition of~$\Delta^k$ from \eqref{E:8.24ai} and write
\begin{equation}
A_k:=\Delta^k\smallsetminus\overline{\Delta^{k-1}},\qquad k\ge1,
\end{equation}
for the annuli used to define the concentric decomposition.
Invoking, in turn, the notation~$B(N)$ from \eqref{E:BN}, define also the thinned annuli
\begin{equation}
A_k':=B(3\cdot 2^{k-2})\smallsetminus B(2^{k-2})
\end{equation}
and note that $A_k'\subseteq A_k$ with $\dist_\infty(A_k',A_k^\cc)\ge 2^{k-2}$.
Network reduction arguments (underlying the Nash-Williams estimate) imply
\begin{equation}
\label{E:14.24}
R_{\eff,h}\bigl(0,B(N)^\cc\bigr)\ge R_{\eff,h}\bigl(\partial^{\text{\rm in}}A_k',\partial^{\text{\rm out}}A_k'\bigr),
\end{equation}
where $\partial^{\text{\rm in}}A$, resp., $\partial^{\text{\rm out}}A$ denote the inner, resp., outer external boundary of the annulus~$A$. We first observe:

\begin{mylemma}
\label{lemma-14.9}
There are~$c,C>0$ such that for each~$k\ge1$,
\begin{equation}
\label{E:14.25}
\BbbP\Bigl(R_{\eff,\,h^{A_k}}\bigl(\partial^{\text{\rm in}}A_k',\partial^{\text{\rm out}}A_k'\bigr)\ge C\texte^{-\hat c\log\log(2^k)}\Bigr)\ge c\,,
\end{equation}
where $h^{A_k}:=$ DGFF in~$A_k$.
\end{mylemma}

\begin{proofsect}{Proof}
Let $N:=2^{k-1}$ and let~$U_1,\dots,U_4$ denote the four (either $4N\times N$ or $N\times 4N$) rectangles that fit into~$A_k'$. We label these in the clockwise direction starting from the one on the right. Since every path from the inner boundary of~$A_k'$ to the outer boundary of~$A_k'$ crosses (and connects the longer sides of) one of these rectangles, Lemma~\ref{lemma-13.17} and the variational representation \twoeqref{E:12.52}{E:12.53} imply
\begin{equation}
R_{\eff,\,h}\bigl(\partial^{\text{\rm in}}A_k',\partial^{\text{\rm out}}A_k'\bigr)
\ge\frac14\min\Bigl\{R^\LR_{U_1}(h),\,R^\UD_{U_2}(h),\,R^\LR_{U_3}(h),\,R^\UD_{U_4}(h)\Bigr\}\,.
\end{equation}
(Note that the resistances are between the \emph{longer} sides of the rectangles.)
Invoking also the Gibbs-Markov property along with Lemma~\ref{lemma-remove}, it suffices to prove that, for some~$C,c>0$ and all~$N\ge1$,
\begin{equation}
\label{E:15.27uia}
\BbbP\Bigl(R^\LR_{U}(h^{B(16 N)})\ge C\texte^{-\hat c\log\log N}\Bigr)\ge c
\end{equation}
for any translate~$U$ of~$B(N,4N)$ contained in~$B(8N)$. 

We will proceed by the duality argument underlying the proof of Proposition~\ref{prop-13.9}. Abbreviate~$M:=16 N$ and consider the decomposition of~$h^{B(M)}=\varphi+\chi$ from Lemma~\ref{lemma-split}. Then for each~$r,a>0$, Exercise~\ref{ex:13.15} yields
\begin{equation}
\label{E:15.28uia}
\BbbP\bigl(R^\LR_U(h^{B(M)})\ge r\bigr)\ge \BbbP\bigl(R^\LR_U(\varphi)\ge r/a\bigr)-a^{-1}\texte^{\hat c\log\log M}\,.
\end{equation}
The path-cut approximate duality shows
\begin{equation}
\label{E:15.29uia}
\BbbP\Bigl(R^\LR_U(\varphi)R^{\UD}_U(-\varphi)\ge\texte^{-2\beta c_1}/64\Bigr)\ge1-\epsilon.
\end{equation}
For any~$r'>0$, Exercise~\ref{ex:13.15} also gives
\begin{equation}
\label{E:15.30uia}
\BbbP\bigl(R^{\UD}_U(\varphi)\le r'\bigr)\ge\BbbP\bigl(R^{\UD}_U(h^{B(M)})\le r'/a\bigr)-a^{-1}\texte^{\hat c\log\log M}.
\end{equation}
Setting~$r'/a:=\texte^{\hat c\log\log M}$ with $a:=C'\texte^{\hat c\log\log M}$ for some~$C'>0$ large enough, Lemma~\ref{lemma-13.22} bounds the probability on the right of \eqref{E:15.30uia} by a positive constant. Via \eqref{E:15.29uia} for~$\epsilon$ small, this yields a uniform lower bound on $\BbbP(R^\LR_U(\varphi)\ge r/a)$ for~$r/a:=(\texte^{-2\beta c_1}/64)/r'$. The bound \eqref{E:15.28uia} then gives \eqref{E:15.27uia}.
\end{proofsect}

\begin{proofsect}{Proof of Proposition~\ref{prop-14.8}}
Using the representation of $h:=$ DGFF in~$\Z^2\smallsetminus\{0\}$ from Exercise~\ref{ex:8.13}, the estimates from Lemmas~\ref{lemma-8.8}, \ref{lemma-8.10a}, \ref{lemma-8.11a} and~\ref{lemma-10.4a} on the various ``bits and pieces'' constituting the concentric decomposition yield
\begin{equation}
h\le -S_k+(\log n)^2+h_k'\quad\text{on }A_k'
\end{equation}
for all~$k=1,\dots,n$ as soon as~$n:=\lfloor\log_4 N\rfloor$ is sufficiently large. The inequality \eqref{E:14.24} then gives
\begin{equation}
R_{\eff,h}\bigl(0,B(N)^\cc\bigr)\ge \texte^{2\beta [S_k-(\log n)^2]}R_{\eff,\,h_k'}(\partial^{\text{\rm in}}A_k',\partial^{\text{\rm out}}A_k').
\end{equation}
The fields $\{h_k'\colon k\ge1\}$ are independent and hence so are the events
\begin{equation}
E_k:=\Bigl\{R_{\eff,\,h_k'}(\partial^{\text{\rm in}}A_k',\partial^{\text{\rm out}}A_k')\ge C\texte^{-\hat c\log\log(2^k)}\Bigr\},\quad k\ge1.
\end{equation}
Moreover, $h_k'\,\,\laweq\,\,h^{A_k}$ and so Lemma~\ref{lemma-14.9} shows $\BbbP(E_k)\ge c$ for all~$k\ge1$. 

A standard use of the Borel-Cantelli lemma implies that, almost surely for~$n$ sufficiently large, the longest interval of~$k\in\{1,\dots,n\}$ where~$E_k$ fails is of length at most order~$\log n$. The Gaussian nature of the increments of~$S_k$ permits us to assume that $\max_{1\le k\le n}|S_{k+1}-S_k|\le\log n$ for~$n$ large and so the value of~$S_k$ changes by at most another factor of order~$(\log n)^2$ over any interval of~$k$'s where~$E_k$ fails. Since~$\log\log(2^n)\ll(\log n)^2$ for~$n$ large, this and the fact that $\log n=\log\log N+O(1)$ readily implies the claim.
\end{proofsect}

We are now ready to complete our proof of Theorem~\ref{thm-12.5}:

\begin{proofsect}{Proof of \eqref{E:1.12ua} in Theorem~\ref{thm-12.5}}
Since~$\{S_k\colon k\ge1\}$ is a random walk with Gaussian increments of bounded and positive variance, Chung's Law of the Iterated Logarithm (see~\cite[Theorem~3]{Chung}) implies that, for some constant~$c>0$,
\begin{equation}
\liminf_{n\to\infty}\frac{\max_{k\le n}S_k}{\sqrt{n/\log\log n}}\ge c,\quad \BbbP\text{-a.s.}
\end{equation}
As $\log n=\log\log N+O(1)$, the claim  follows from Proposition~\ref{prop-14.8}.
\end{proofsect}

\nopagebreak
\begin{figure}[t]
\vglue-1mm
\centerline{\includegraphics[width=0.4\textwidth]{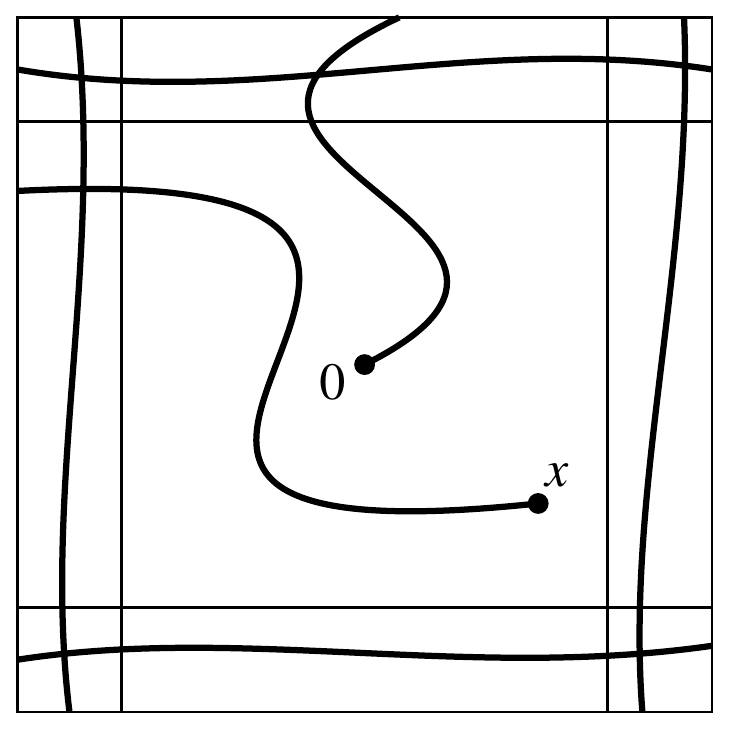}}
\vglue0mm
\begin{quote}
\small 
\caption{
\label{fig-path-to-annuli}
\small
A path-crossing event underlying the estimate \eqref{E:14.36}. Notice that the existence of a path from~$0$ to the outer square, another path from~$x$ to the outer square and the crossings of the four rectangles as shown imply the existence of a path from~$0$ to~$x$ within the outer square.}
\normalsize
\end{quote}
\end{figure}

It remains to give:

\begin{proofsect}{Proof of Proposition~\ref{prop-14.7}, main idea}
Consider the annuli~$A_k$ as above and notice the following consequences of the network reduction arguments discussed earlier. Fix~$n$ large and denote~$N:=2^n$. A shorting argument then implies that for any $k^\star\in\{1,\dots,n\}$ such that $x\in\Delta^{k^\star-1}$,
\begin{equation}
R_\eff\bigl(x,B(N)^\cc\bigr)\ge R_\eff\bigl(x,\partial^{\text{\rm in}}A_{k^\star}'\bigr)+R_{\eff}(\partial^{\text{\rm in}}A_{k^\star}',\partial^{\text{\rm out}}A_{k^\star}').
\end{equation}
Next consider the four rectangles constituting the annulus~$A_k'$ and let~$R^{\ssup i}_{A_k'}$ denote the effective resistance between the \emph{shorter} sides of the~$i$-th rectangle. Lemma~\ref{lemma-13.16} and the path-crossing argument from Fig.~\ref{fig-path-to-annuli} show that, for each~$k_\star\in\{1,\dots,n\}$ such that~$x\in\Delta^{k_\star-1}$,
\begin{equation}
\label{E:14.36}
R_{\eff,B(N)}(0,x)\le R_\eff\bigl(0,\partial^{\text{\rm out}}A_{k_\star}'\bigr)
+R_\eff\bigl(x,\partial^{\text{\rm out}}A_{k_\star}'\bigr)+\sum_{i=1}^4 R^{\ssup{i}}_{A_{k_\star}'}\,.\end{equation}
Assuming that~$k_\star<k_\star$, also have
\begin{equation}
\label{E:14.38}
R_\eff\bigl(x,\partial^{\text{\rm in}}A_{k^\star}'\bigr)\ge R_\eff\bigl(x,\partial^{\text{\rm out}}A_{k_\star}'\bigr)
\end{equation}
and so
\begin{equation}
\label{E:15.38nwt}
R_\eff\bigl(x,B(N)^\cc\bigr) - R_{\eff,B(N)}(0,x)
\ge R_{\eff}(\partial^{\text{\rm in}}A_{k^\star}',\partial^{\text{\rm out}}A_{k^\star}') 
-\sum_{i=1}^4 R^{\ssup{i}}_{A_{k_\star}'}\,.
\end{equation}
The point is now to show that, with probability tending to one as~$n\to\infty$, one can find integers~$k_\star,k^\star$ with
\begin{equation}
\label{E:14.39}
(\log n)^\delta<k_\star<k_\star\le n
\end{equation}
such that the first resistance on the right of \eqref{E:15.38nwt} significantly exceeds the sum of the four resistances therein. For this we note that, thanks to Lemma~\ref{lemma-14.9} and Proposition~\ref{prop-13.15}, the overall scale of these resistances is determined by the value of the random walk~$S_k$ at~$k:=k^\star$ (for the first resistance) and~$k:=k_\star$ (for the remaining resistances). In light of the observations made in the proof of Proposition~\ref{prop-14.8}, it will suffice to find~$k_\star,k^\star$ obeying \eqref{E:14.39} such that
\begin{equation}
S_{k^\star}\ge2(\log n)^2\quad\text{and}\quad S_{k_\star}\le -2(\log n)^2
\end{equation}
occur with overwhelming probability.
This requires a somewhat more quantitative version the Law of the Iterated Logarithm for which we refer the reader to the original paper.
\end{proofsect}

\section{Wrapping up}
\noindent
We will now finish by presenting the proofs of the desired lower bounds on the expected exit time and the heat kernel. 

\begin{proofsect}{Proof of Theorem~\ref{thm-exittime}, lower bound in \eqref{E:12.8}}
Let~$\delta\in(0,1)$. The hitting-time identity in Lemma~\ref{lemma-14.1} along with Proposition~\ref{prop-14.7} imply
\begin{equation}
\label{E:14.41}
E_h^0(\tau_{B(N)^\cc})\ge \pi_h\bigl(B(N\texte^{-(\log N)^\delta})\bigr)\log N
\end{equation}
with probability tending to one as~$N\to\infty$. Theorem~\ref{thm-intermediate} on the size of the inter\-mediate level sets and Theorem~\ref{thm-7.3} on the tightness of the absolute maximum yield
\begin{equation}
\pi_h\bigl(B(N\texte^{-(\log N)^\delta})\bigr)\ge N^{\theta(\beta)}\texte^{-(\log N)^{2\delta}}
\end{equation}
with probability tending to one as~$N\to\infty$. The claim follows.
\end{proofsect}

It remains to show:

\begin{proofsect}{Proof of Theorem~\ref{thm-heatkernel}, lower bound in~\eqref{E:12.5a}}
Let~$\Xi_T^\star$ be the union of $\{0\}\cup B(N)^\cc$ with the set of all~$x\in B(N\texte^{-(\log N)^\delta})$ such that
\begin{equation}
R_{\eff,B(N)}(0,x)\vee R_\eff\bigl(x,B(N)^\cc\bigr)\le\texte^{(\log N)^{1/2+\delta}}
\end{equation}
Abusing our earlier notation slightly, let~$Y_k$ be the $k$-th visit of~$X$ to~$\Xi_T^\star$ (counting the starting point as $k=0$ case). Denote $\hat\tau:=\inf\{k\ge0\colon Y_k\in B(N)^\cc\}$. Then
\begin{equation}
\begin{aligned}
E^0_h(\hat\tau)
&\le T P^0_h(\hat\tau\le T)+P^0_h(\hat\tau> T)\Bigl(T+\max_{x\in\Xi_N} E^x_h(\hat\tau)\Bigr)
\\
&=T+P^0_h(\hat\tau> T)\max_{x\in\Xi_N} E^x_h(\hat\tau).
\end{aligned}
\end{equation}
The hitting time estimate (Corollary~\ref{cor-14.2}) ensures
\begin{equation}
E^x_h(\hat\tau)\le\pi_h\bigl(\Xi_T^\star\cap B(N)\bigr)\texte^{(\log N)^{1/2+\delta}},\quad x\in\Xi^\star_T,
\end{equation}
implying
\begin{equation}
\label{E:14.46}
P^0_h(\hat\tau> T)\ge \pi_h\bigl(\Xi_T^\star\cap B(N)\bigr)^{-1}\texte^{-(\log N)^{1/2+\delta}}\bigl[E^0_h(\hat\tau)-T\bigr].
\end{equation}
By our choice of~$\Xi_T^\star$, the lower bound on the voltage from Proposition~\ref{prop-14.7} applies to all~$x\in\Xi^\star_T\cap B(N)$ and so 
\begin{equation}
\label{E:14.47}
E^0_h(\hat\tau)\ge \pi_h\bigl(\Xi_T^\star\cap B(N)\bigr)\log N.
\end{equation}
We now need to solve:

\begin{myexercise}
Prove that
\begin{equation}
\BbbP\Bigl(\pi_h\bigl(\Xi_N^\star\cap B(N)\bigr)\le N^{\theta(\beta)}\texte^{-(\log N)^{2\delta}}\Bigr)\,\underset{N\to\infty}\longrightarrow\,0.
\end{equation}
Hint: Estimate $\E\pi_h(B(N\texte^{(\log N)^\delta})\smallsetminus\Xi_N^\star)$ using the FKG inequality and the fact that $h\mapsto R_{\eff,h}(u,v)$ is decreasing.
\end{myexercise}

\noindent
For $N:=T^{1/\theta(\beta)}\texte^{(\log T)^\delta}$ this implies $E^0_h(\hat\tau)\ge 2T$ and via \twoeqref{E:14.46}{E:14.47} also
\begin{equation}
P^0_h(\hat\tau> T)\ge\texte^{-(\log N)^{1/2+\delta}}.
\end{equation}
But $\hat\tau\le\tau_{B(N)^\cc}:=\inf\{k\ge0\colon X_k\in B(N)^\cc\}$ and so
\begin{equation}
P^0\bigl(X_T\in B(N)\bigr)\ge P^0(\tau_{B(N)^\cc}>T\bigr)\ge \texte^{-(\log N)^{1/2+\delta}}.
\end{equation}
Plugging this in \eqref{E:14.16au} and invoking \eqref{E:14.7} then gives the claim.
\end{proofsect}


\chapter{Questions, conjectures and open problems}
\label{lec-16}\noindent
In this final lecture we discuss conjectures and open problems. These are all considerably harder than the exercises scattered throughout these notes. In fact, most of the problems below constitute non-trivial projects for future research. The presentation is at times (and sometimes deliberately) vague on detail and may even require concepts that lie outside the scope of these lecture notes.

\section{DGFF level sets}
\noindent
We open our discussion by questions that are related to the description of the DGFF level sets. Let~$h^{D_N}$ be the DGFF in~$D_N$, for a sequence~$\{D_N\colon N\ge1\}$ of admissible lattice approximation of domain~$D\in\mathfrak D$. We commence by:

\begin{myproblem}[Joint limit of intermediate level sets]
Find a way to extract a joint distributional limit of the level sets
\begin{equation}
\bigl\{x\in D_N\colon h_x\ge 2\sqrt g\lambda\log N\bigr\}
\end{equation}
or their associated point measures \eqref{E:2.30ua}, simultaneously for all $\lambda\in(0,1)$. Use this to design a coupling of the corresponding~$Z_\lambda^D$'s and show that they all arise (via the construction of the LQG measure) from the same underlying~CGFF.
\end{myproblem} 

\noindent
We remark that we see a way to control any finite number of the level sets by building on the moment calculations underpinning Theorem~\ref{thm-intermediate}. The problem is that, the more level sets we wish to include, the higher moments we seem to need. The key is to find a formulation that works simultaneously for all~$\lambda$'s and to show that that the limit measures can be linked to the same copy of the~CGFF. 

The insistence on connecting the limit process to the CGFF may in fact suggest a solution: Use the same sample of CGFF to define the DGFF on~$D_N$ for all~$N\ge1$ and thus all level sets simultaneously. Then show that this CGFF is what one gets (via the construction of Gaussian Multiplicative Chaos) as a limit measure. 
The same strategy will perhaps also make  the next question accessible:

\begin{myproblem}[Fluctuations and deviations]
Describe the rate of convergence in, e.g., Corollary~\ref{cor-2.6}. Identify the limit law of the fluctuations of the level set away from its limit value. Study (large) deviations estimates for the intermediate level sets.
\end{myproblem}

\noindent
By ``large deviations'' we mean, intentionally vaguely, any kind of deviation whose probability tends to zero as~$N\to\infty$. We note that this may be relevant for instance already in attempts to connect the limit law of the intermediate level sets to the construction of the LQG measure. Indeed, one is effectively asking to integrate the process~$\eta^D_N$ against the function
\begin{equation}
f(x,h):=\texte^{\beta h}\quad\text{where}\quad\beta:=\lambda\alpha.
\end{equation}
Plugging this formally into the limit (which is of course not mathematically justified) yields a diverging integral. This indicates that the problem requires some control of the deviations and/or the support of the LQG measure.

\section{At and near the absolute maximum}
\noindent
Moving to the extremal values, we start by a conjecture which, if true, might greatly simplify the existing proof of the uniqueness of the subsequential limits $\eta^D$ of $\{\eta^D_{N,r_N}\colon N\ge1\}$. Indeed, in the version of the proof presented in Lecture~\ref{lec-12}, a great deal of effort is spent on the verification of the Laplace transform tail \eqref{E:1.25ue} for the measures~$Z^D$ extracted from~$\eta^D$ via Theorem~\ref{prop-subseq}. This would not be necessary if we could resolve affimatively:

\begin{myconjecture}[By-passing the Laplace transform tail]
\label{conj-Laplace-tail}
Prove that any family of non-trivial random measures $\{Z^D\colon D\in\mathfrak D\}$ obeying (1-4) in Theorem~\ref{thm-10.14} satisfies also property~(5) for some~$\hat c>0$.
\end{myconjecture}

The reason why we believe this to be possible stems (somewhat vaguely) from the analogy of the Gibbs-Markov property along, say, partitions of squares into smaller squares, with fixed-point equations for the so-called \emph{smoothing transformations}. These have been classified in Durrett and Liggett~\cite{Durrett-Liggett}.

Our next point of interest is Theorem~\ref{thm-LLT} which, we recall, states the local limit theorem for the absolute maximum. The limiting position (scaled by~$N$) for the maximum restricted above~$m_N+t$ is then distributed according to the density $x\mapsto\rho^D(x,t)$ in \eqref{E:density}. From the estimate \eqref{E:1.25ue} we are able to conclude
\begin{equation}
\label{E:rhoD-asymptotic}
\rho^D(x,t)\sim c_\star t\texte^{-\alpha t}r_D(x)^2,\quad t\to\infty,
\end{equation}
but we have not been able to characterize~$\rho^D$ explicitly for any finite~$t$. Still,~$\rho^D$ is determined by~$Z^D$, which should be universal (up to a constant multiple) for a whole class of logarithmically correlated fields in~$d=2$. By analogy with the role of the KPP equation plays in the analysis of the Branching Brownian Motion (cf Bramson~\cite{Bramson}), we hope the following is perhaps achievable: 

\begin{myproblem}[Nailing density explicitly]
Determine~$\rho^D$ explicitly or at least characterize it via, e.g., a PDE that has a unique solution.
\end{myproblem}

Another question related to the extremal values is the crossover between the regime of intermediate level sets and those within order-unity of the absolute maximum. We state this as:

\begin{myconjecture}[Size of deep extremal level set]
\label{cor-12.13ua}
There is~$c\in(0,\infty)$ such that, for each open~$A\subseteq D$ (with $\Leb(\partial A)=0$)
\begin{equation}
\label{E:12.1ua}
\frac1t\texte^{-\alpha t}\,\#\bigl\{x\in\Gamma_N^D(t)\colon x/N\in A\bigr\}\,\,\,\underset{\begin{subarray}{c} N\to\infty\\t\to\infty\end{subarray}}\lawarrow\,\,\, cZ^D(A),
\end{equation}
with convergence in law valid even jointly for any finite number of disjoint open sets $A_1,\dots,A_k\subseteq D$ (with $\Leb(\partial A_j)=0$ for all~$j$).
\end{myconjecture}

Note that this constitutes a stronger version of
Theorem~\ref{thm-DZ1} because \eqref{E:12.1ua} implies
\begin{equation}
\label{E:12.2ua}
\lim_{r\downarrow0}\,\liminf_{t\to\infty}\,\liminf_{N\to\infty}\, P\Bigl(rt\texte^{\alpha t}\le|\Gamma_N^D(t)|\le r^{-1}t\texte^{\alpha t}\Bigr)=1.
\end{equation}
Writing $\Gamma_N(A,t)$ for the set on the left-hand side of \eqref{E:12.1ua}, Theorem~\ref{thm-DZ2} and truncation of the absolute maximum to lie below~$m_N+\sqrt t$ implies, for each~$\lambda>0$,
\begin{multline}
E\bigl(\texte^{-\lambda|\Gamma_N^D(A,t)|}\bigr)
\\=E\biggl(\exp\Bigl\{-Z^D(A)\texte^{\alpha t}\int_0^{t+\sqrt t}\textd h\,\texte^{-\alpha h}E_\nu(1-\texte^{-\lambda f_h(\phi)})\Bigr\}\biggr)+o(1)\,,
\end{multline}
where $o(1)\to0$ as~$N\to\infty$ followed by $t\to\infty$ and
\begin{equation}
f_h(\phi):=\sum_{z\in\Z^2}1_{[0,h]}(\phi_z).
\end{equation}
Using the above for~$\lambda:=t^{-1}\texte^{-\alpha t}$, Conjecture~\ref{cor-12.13ua} seems closely linked to:

\begin{myconjecture}
For~$f_h$ as above,
\begin{equation}
\lim_{h\to\infty}\texte^{-\alpha h} E_\nu(f_h) \text{ exists in } (0,\infty)\,.
\end{equation}
\end{myconjecture}

\noindent
Note that that \eqref{E:12.1ua} is similar to the limit law for the size of Daviaud's level sets proved in Corollary~\ref{cor-2.6}. 

Another conjecture concerning deep extremal level sets is motivated by the question on at what levels is the critical LQG measure mainly supported and how is the DGFF distributed thereabout. The Seneta-Heyde normalization indicates that this happens at levels order $\sqrt{\log N}$ below the absolute maximum. Discussions of this problem with O.~Louidor during the summer school led to:

\begin{myconjecture}[Profile of near-extremal level sets]
\label{conj-15.5}
There is a constant~$c>0$ such that (writing $h$ for the DGFF in~$D_N$),
\begin{equation}
\eta^D_N:=\frac1{\log N}\sum_{x\in D_N}\texte^{\alpha (h_x-m_N)}\,\delta_{x/N}\otimes\delta_{\frac{m_N-h_x}{\sqrt{\log N}}},
\end{equation}
where~$m_N$ is as in \eqref{E:7.8}, obeys
\begin{equation}
\eta^D_N\,\,\underset{N\to\infty}\lawarrow\,\,c\,Z^D(\textd x)\otimes\,1_{[0,\infty)}(h)\,h\,\texte^{-\frac{h^2}{2g}}\textd h\,.
\end{equation}
Here~$Z^D$ is the measure in Theorem~\ref{thm-extremal-vals}.
\end{myconjecture}

\noindent
In particular, the density
\begin{equation}
h\mapsto 1_{[0,\infty)}(h)\,h\,\texte^{-\frac{h^2}{2g}}
\end{equation}
gives the asymptotic ``profile'' of the values of the DGFF that contribute to the critical LQG measure at scales order-$\sqrt{\log N}$ below~$m_N$. 

\section{Universality}
\noindent
A very natural question to ask is to what extent are the various results reported here universal with respect to various changes of the underlying setting. As a starter, we pose:

\begin{myproblem}[Log-correlated Gaussian fields in~$d=2$]
\label{pb-15.6}
For each~$N\ge1$, consider a Gaussian field in~$D_N$ whose covariance matrix~$C_N$ scales, in the sense described in Theorem~\ref{thm-1.17}, to the continuum Green function in the underlying (continuum) domain and such that for all~$z\in\Z^2$, the limit
\begin{equation}
\lim_{N\to\infty} \bigl[C_N(x,x+z)-C_N(x,x)\bigr]
\end{equation}
exists, is independent of~$x$ and is uniform in~$x$ with~$\dist_\infty(x,D_N^\cc)>\epsilon N$, for every~$\epsilon>0$. Prove the analogues of Theorems~\ref{thm-intermediate} and~\ref{thm-extremal-vals}.
\end{myproblem}

Some instances of this can be resolved based on the results already presented here. For instance, if~$C_N$ is given as
\begin{equation}
C_N(x,y)=G^N(x,y)+F(x-y),
\end{equation}
where~$F\colon\Z^2\to\R$ a bounded, positive semi-definite function that decays fast to zero (as the argument tends to infinity), one can realize $\NN(0,C_N)$ as the sum of the DGFF plus a Gaussian field with fast-decaying correlations. The effect of the additional Gaussian field can be handled within the structure of Theorems~\ref{thm-intermediate} and~\ref{thm-extremal-vals}. The point is thus to go beyond such elementary variations.

Another question is that of generalizations to log-correlated Gaussian processes in dimensions~$d\ne 2$. Naturally, here the choice of the covariance is considerably less canonical. Notwithstanding, convergence of the centered maximum to a non-degenerate limit law has already been established in a class of such fields (in any~$d\ge1$) by Ding, Roy and Zeitouni~\cite{DRZ}. (This offers some headway on Problem~\ref{pb-15.6} as well.) The next goal is an extension to the full extremal process as stated for the DGFF in Theorem~\ref{thm-extremal-vals}.

The question of universality becomes perhaps even more profound once we realize that there are other non-Gaussian models where one expects the results reported in these notes to be relevant. The first one of these is the class of \emph{Gradient Models}. These are random fields with law in finite sets~$V\subset\Z^2$ given~by
\begin{equation}
\label{E:1.1b}
P(h^V\in A):=\frac1{\text{\rm norm.}}\int_A\texte^{-\sum_{(x,y)\in\cmss E(\Z^d)}\VV(h_x-h_y)}\prod_{x\in V}\textd h_x\prod_{x\not\in V}\delta_0(\textd h_x)\,,
\end{equation}
where~$\VV\colon\R\to\R$ is the potential which is assumed even, continuous, bounded from below and with superlinear growth at infinity. As an inspection of \eqref{E:1.1} shows, the DGFF is included in this class via
\begin{equation}
\VV(h):=\frac1{4d}h^2.
\end{equation}
The formula \eqref{E:1.1b} defines the field with zero boundary conditions, although other boundary conditions can be considered as well (although they cannot be easily reduced to zero boundary conditions outside the Gaussian case).

Much is known about these models when~$\VV$ is uniformly strictly convex (i.e., for~$\VV''$ positive and bounded away from zero and infinity). Here through the work of Funaki and Spohn~\cite{FS}, Naddaf and Spencer~\cite{NS}, Giacomin, Olla and Spohn~\cite{GOS} we know that the field tends to a linear transform of CGFF in the thermodynamic limit (interpreted in the sense of gradients in $d=1,2$), and by Miller~\cite{Miller} also in the scaling limit in~$d=2$. Recently, Belius and Wu~\cite{Belius-Wu} have proved the counterpart of Theorem~\ref{thm-2.1} by identifying the leading-order growth of the absolute maximum. Wu and Zeitouni~\cite{WZ} have then extended the Dekking-Host subsequential tightness argument from Lemma~\ref{lemma-DH} to this case as well. 

The situation seems ripe to tackle  more complicated questions such as those discussed in these notes. We believe that the easiest point of entry is via:

\begin{myproblem}[Intermediate level sets for gradient fields]
\label{pb-15.7}
Prove the scaling limit of intermediate level sets for gradient models with uniformly strictly convex potentials.
\end{myproblem}

\noindent
This should be manageable as the main technical input in these are moment calculations that lie at the heart of~\cite{Belius-Wu} which deals with the case where they should, in fact, be hardest to carry out. 

Perhaps even closer to the subject of these notes is the problem of DGFF generated by a \emph{Random Conductance Model} (see~\cite{Biskup-RCM}): Assign to each edge $e$ in~$\Z^2$ a non-negative conductance~$c(e)$. We may in fact assume that the conductances are \emph{uniformly elliptic}, meaning that
\begin{equation}
\exists\lambda\in(0,1)\,\,\forall e\in E(\Z^2)\colon\qquad c(e)\in[\lambda,\lambda^{-1}]\quad\text{a.s.}
\end{equation}
Given a realization of the conductances, we can define an associated DGFF in domain~$D\subsetneq\Z^2$ as~$\NN(0,G^D)$ where~$G^D:=(1-\cmss P)^{-1}$ is, for~$\cmss P$ related to the conductances as discussed at the beginning of Section~\ref{sec14.1}, the Green function of the random walk among random conductances. We then pose:

\begin{myproblem}[DGFF over Random Conductance Model]
\label{pb-15.8}
Assume that the conductances are uniformly elliptic and their law is invariant and ergodic with respect to the shifts of~$\Z^2$. Prove that for almost every realization of the conductances, the analogues of Theorems~\ref{thm-intermediate} and~\ref{thm-extremal-vals} hold.
\end{myproblem}

This could in principle be easier than Problem~\ref{pb-15.7} due to the underlying Gaussian nature of the field. Notwithstanding, the use of homogenization theory, a key tool driving many studies of gradient models, can hardly be avoided. Incidentally, Problem~\ref{pb-15.8} still falls under the umbrella of gradient models, albeit with non-convex interactions. The connection to the DGFF over random conductances yields much information about the corresponding gradient models as well; see the papers of the author with R.~Koteck\'y~\cite{BK} and H.~Spohn~\cite{BS}.

Another problem of interest that should be close to the DGFF is that of \emph{local time} of the simple random walk. The connection to the DGFF arises via the \emph{Dynkin isomorphism} (usually attributed to Dynkin~\cite{Dynkin} though going back to Symanzik~\cite{Symanzik}) or, more accurately, the \textit{second Ray-Knight theorem} (see, e.g., Eisenbaum~\textit{et al}~\cite{EKMRS}).
This tool has been useful in determining the asymptotic behavior of the \emph{cover time} (namely, in the studies of the leading-order behavior by Dembo, Peres, Rosen and Zeitouni~\cite{DPRZ}, Ding~\cite{Ding-cover-time}, Ding, Lee and Peres~\cite{Ding-Lee-Peres}, in attempts to nail the subleading order in Belius and Kistler~\cite{Belius-Kistler} and Abe~\cite{Abe2} and, very recently, in a proof of tightness under proper centering in Belius, Rosen and Zeitouni~\cite{BRZ}).

We wish to apply the connection to the study of the local time $\ell_t$ for a continuous-time (constant-speed) random walk on $N\times N$ torus $\T_N$ in~$\Z^2$ started at~$0$. We choose the  parametrization by the actual time the walk spends at~$0$; i.e.,~$\ell_t(0)=t$ for all~$t\ge0$. For~$\theta>0$, define the time scale
\begin{equation}
t_\theta:=\theta [g\log N]^2.
\end{equation}
 Then the total time of the walk, $\sum_x\ell_{t_\theta}(x)$, is very close to a $\theta$-multiple of the cover time. It is known that $(\ell_t-t)/\sqrt{t}$ scales to a DGFF on $\T_N\smallsetminus\{0\}$ in the limit~$t\to\infty$. This motivated Abe~\cite{Abe1} to show that, as~$N\to\infty$, the maximum of~$\ell_{t_\theta}$ scales as
\begin{equation}
\frac{\max_{x\in \T_N}\ell_{t_\theta}(x)-t_\theta}{\sqrt{t_\theta}}=\Bigl(1+\frac1{2\sqrt\theta}+o(1)\Bigr)2\sqrt g\log N\,.
\end{equation}
He also studied the cardinality of the level set
\begin{equation}
\label{level}
\Bigl\{x\in \T_N\colon \frac{\ell_{t_\theta}(x)-t_\theta}{\sqrt{t_\theta}}\ge2\eta\sqrt g\log N\Bigr\}
\end{equation}
and proved a result analogous, albeit with \emph{different} exponents, to Daviaud's~\cite{Daviaud} (see Theorem~\ref{thm-2.2}). Some correspondence with the DGFF at the level of intermediate level sets and maxima thus exists but is far from straightforward.

The proofs of \cite{Abe1} are based on moment calculations under a suitable truncation, which is also the method used for proving \eqref{E:2.7}. It thus appears that a good starting point would be to solve:

\begin{myproblem}
Show that a suitably normalized level set in \eqref{level}, or a point measure of the kind \eqref{E:2.30ua}, admits a scaling limit for all~$0<\eta<1+\frac1{2\sqrt\theta}$.
\end{myproblem}

We remark that Abe~\cite{Abe2} studied the corresponding problem also for the random walk on a homogeneous tree of depth~$n$. By embedding the tree ultrametrically into a unit interval, he managed to describe the full scaling limit of the process of extremal local maxima, much in the spirit of Theorem~\ref{prop-subseq}. Remarkably, the limit process coincides (up to a constant multiple of the intensity measure) with that for the DGFF on the tree. A description of the cluster process is in the works~\cite{Abe-Biskup}.

\textit{Update in revision:} As shown in a recent posting by Abe and the author~\cite{Abe-Biskup2}, the intermediate level sets for the local time have now been shown to have the same scaling limit as the DGFF, albeit in a different (and considerably more convenient) parametrization than suggested in \eqref{level}. Concerning the distributional limit of the cover time; this has recently been achieved for the random walk on a homogeneous tree by Cortines, Louidor and Saglietti~\cite{CLS}.

\section{Random walk in DGFF landscape}
\noindent
The next set of problems we wish to discuss deals with the random walk driven by the DGFF. Here the first (and obvious) follow-up question is the complementary inequality to that stated in Corollary~\ref{cor-12.4}:

\begin{myproblem}[Subdiffusive upper bound]
Using the notation from \eqref{E:12.9}, show that, for each~$\beta>0$ and each~$\delta>0$,
\begin{equation}
P_h^0\Bigl(|X_T|\le T^{1/\psi(\beta)}\texte^{(\log T)^{1/2+\delta}}\Bigr)
\,\,\underset{T\to\infty}\longrightarrow\,\,1\,,\quad\text{\rm in $\BbbP$-probability}.
\end{equation}
\end{myproblem}

\noindent
A key issue in here is to show that the exit time~$\tau_{B(N)^\cc}$ is (reasonably) concentrated around its expectation. We believe that this can perhaps be resolved by having better control of the concentration of the effective resistance. 

\begin{figure}[t]
\vglue-1mm
\centerline{\includegraphics[width=0.27\textwidth]{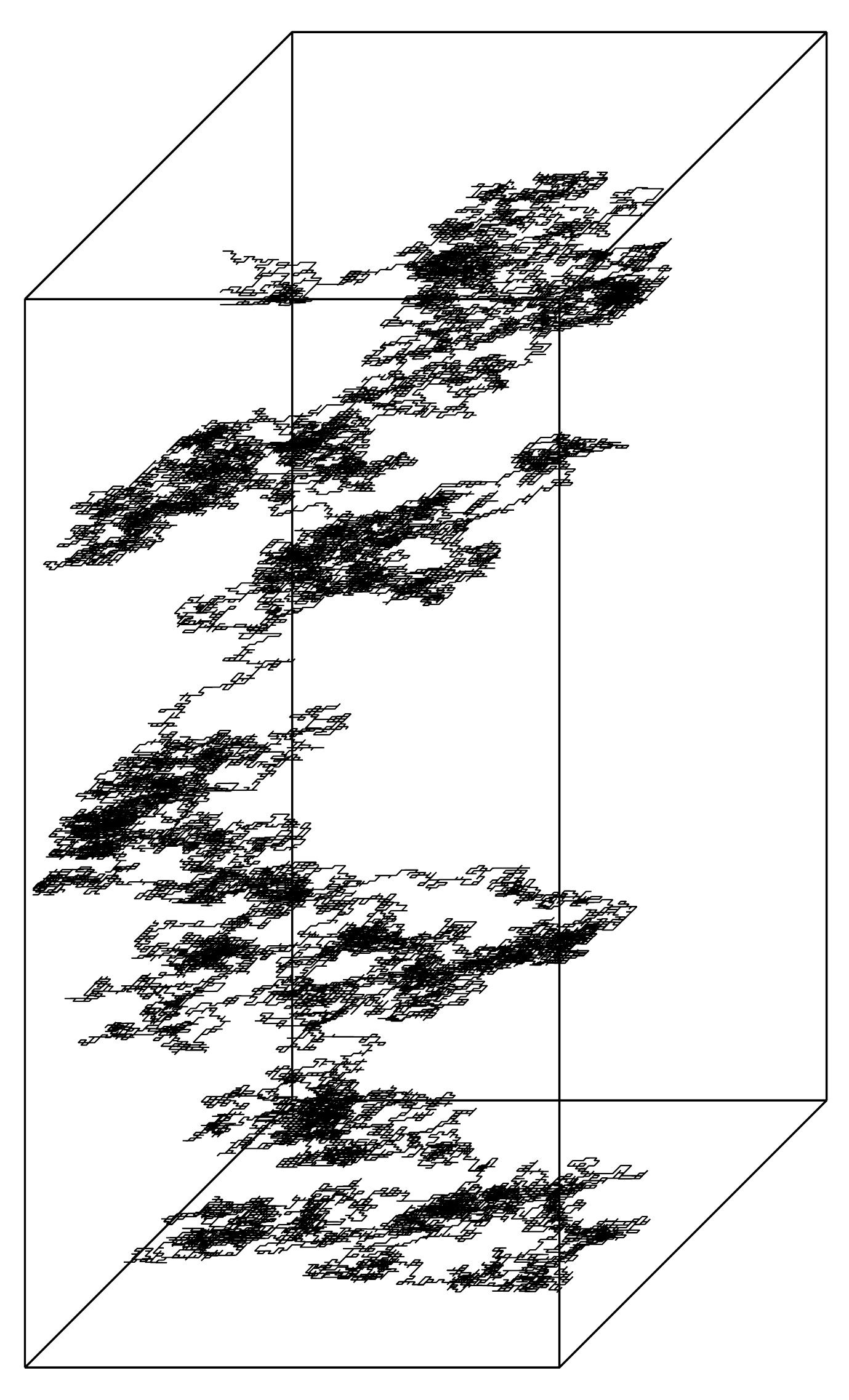}
\includegraphics[width=0.27\textwidth]{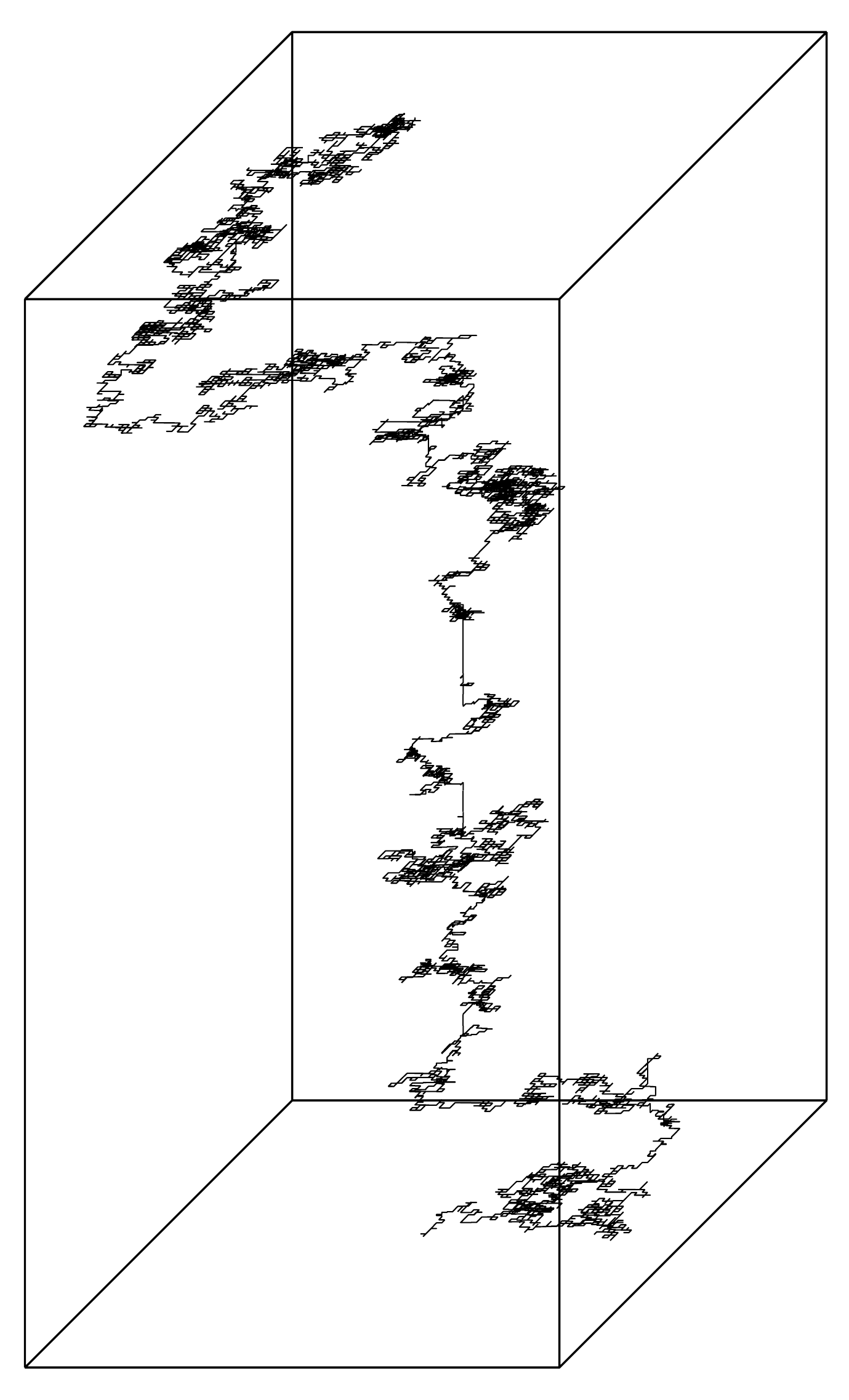}
\includegraphics[width=0.27\textwidth]{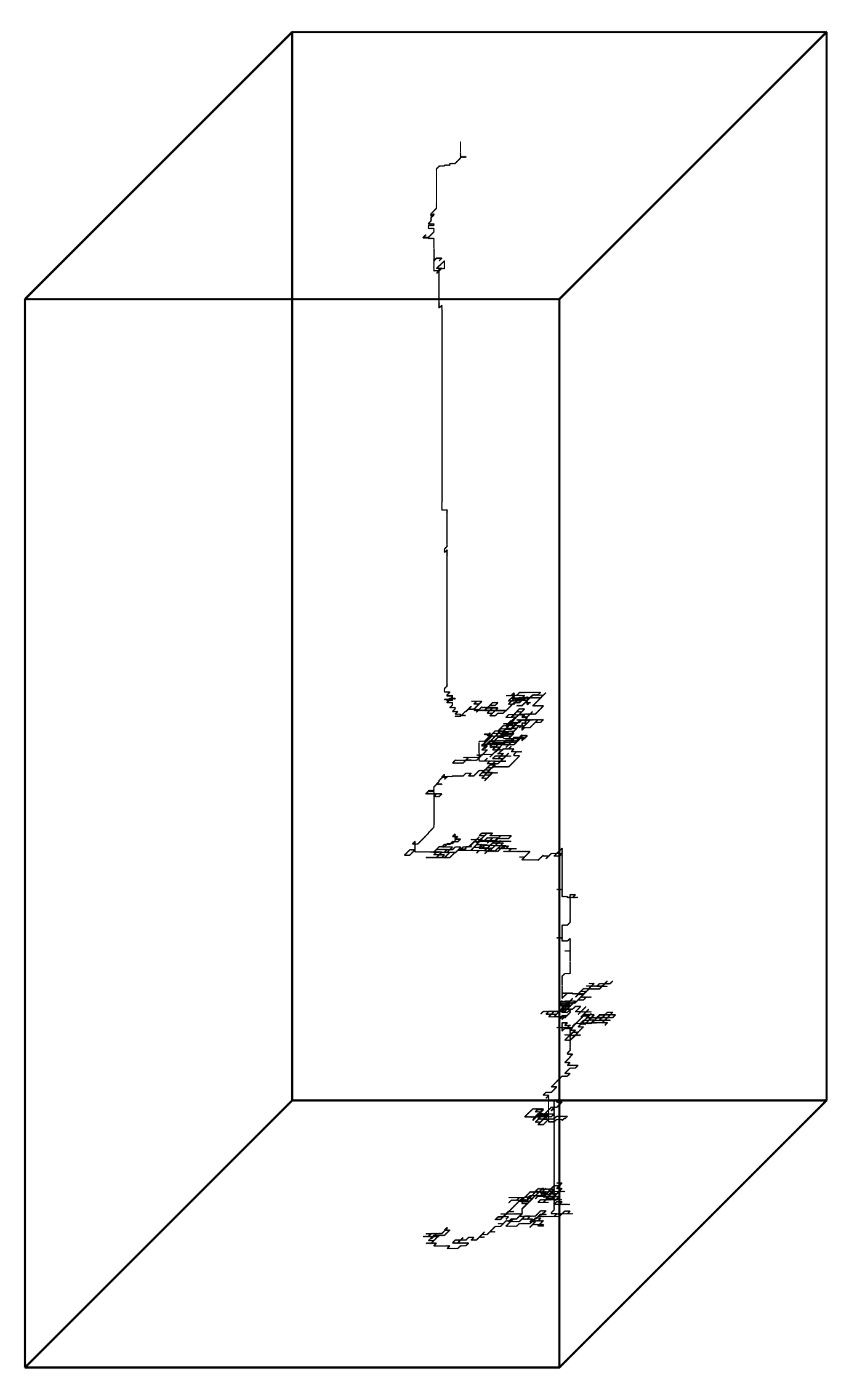}
}

\vglue0mm
\begin{quote}
\small 
\caption{
\label{fig-lw-walks}
\small
Runs of $100000$ steps of the Liouville random walk with~$\beta$ equal to~$0.2$, $0.6$ and~$1.2$ multiples of~$\beta_\cc$. Time runs upwards along the vertical axis. For~$\beta>\beta_\cc$ the walk becomes trapped for a majority of its time.}
\normalsize
\end{quote}
\end{figure}

Another question concerns a possible scaling limit of the random walk. Heuristic considerations make it fairly clear that the limit process cannot be the \emph{Liouville Brownian Motion} (LBM), introduced in Berestycki~\cite{Berestycki-LBM} and Garban, Rhodes and Vargas~\cite{GRV}. Indeed, this process exists only for~$\beta\le\beta_\cc$ (with the critical version constructed by Rhodes and Vargas~\cite{RV-LBM}) and simulations of the associated Liouville Random Walk, which is a continuous-time simple symmetric random walk with exponential holding time with parameter~$\texte^{\beta h_x}$ at~$x$, show a dramatic change in the behavior as~$\beta$ increases through~$\beta_\cc$; see Fig.~\ref{fig-lw-walks}. (The supercritical walk is trapped for a majority of its time with transitions described by a $K$-process; cf Cortines, Gold and Louidor~\cite{CGL}.) No such dramatic change seems to occur for the random walk with the transition probabilities in \eqref{eq-def-transition}, at any~$\beta>0$.

\nopagebreak

The reason why (we believe) the LBM is not the limit process is because its paths, once stripped of their specific parametrization, are those of the two-dimensional Brownian motion. In particular, their law is completely decoupled from the underlying CGFF. On the other hand, the limit of our random walk must feel the drift towards larger ``values'' of the CGFF which will in turn couple the law of the paths (still regarded as geometric objects) to the (limiting) CGFF. 

We thus have to define the limit process through the important features of the discrete problem. Here we observe that the random walk is reversible with respect to the measure~$\pi_h$ defined in \eqref{E:12.12}. As~$\pi_h(x)\approx \texte^{2\beta h_x}$, whose limit is (for~$\beta<\tilde\beta_\cc=\beta_\cc/2$) a LQG measure, one might expect that the resulting process should be reversible with respect to the LQG measure. Although this remains  problematic because of the restriction on~$\beta$, we can perhaps make sense of this by imposing reversibility only for a regularized version of the CGFF and taking suitable limits.

Consider a sample of a smooth, centered Gaussian process $\{h_\epsilon(x)\colon x\in D\}$ such that, for $N:=\lfloor\epsilon^{-1}\rfloor$,
\begin{equation}
\Cov\bigl(h_\epsilon(x),h_\epsilon(y)\bigr)= G^{D_N}\bigl(\lfloor Nx\rfloor, \lfloor Ny\rfloor\bigr)+o(1),\quad \epsilon\downarrow0.
\end{equation}
Then $h_\epsilon$ tends in law to the CGFF as~$\epsilon\downarrow0$. Define a diffusion~$X^\epsilon$ via the Langevin equation
\begin{equation}
\label{E:15.18}
\textd X^\epsilon_t=\beta\nabla h_\epsilon(X_t^\epsilon)\textd t+\sqrt 2\,\textd B_t\,,
\end{equation}
where~$\{B_t\colon t\ge0\}$ is a standard Brownian motion. The (unnormalized) Gibbs measure $\texte^{\beta h_\epsilon(x)}\textd x$ is then stationary and reversible for~$X^\epsilon$. Moreover, under a suitable time change, the process $X^\epsilon$ mimics closely the dynamics of the above random walk on boxes of side-length~$\epsilon^{-1}$ (and~$\beta$ replaced by~$2\beta$). We then pose:

\begin{myproblem}[Liouville Langevin Motion]
\label{pb-15.11}
Under a suitable time change for~$X^\epsilon$, prove that the pair $(h_\epsilon, X^\epsilon)$ converges as~$\epsilon\downarrow0$ jointly in law to $(h,X)$, where~$h$ is a CGFF and~$X$ is a process with continuous (non-constant) paths (correlated with~$h$).
\end{myproblem}

We propose to call the limit process~$X$ the \emph{Liouville Langevin Motion}, due to the connection with the Langevin equation \eqref{E:15.18} and LQG measure. One strategy for tackling the above problem is to prove  characterizing all possible limit processes through their behavior under conformal maps and/or restrictions to subdomain via the Gibbs-Markov property.

\section{DGFF electric network}
\noindent
The final set of problems to be discussed here concern the objects that we used to control the random walk in DGFF landscape. Consider first the effective resistivity $R_\eff(x,B(N)^\cc)$. Proposition~\ref{prop-14.8} suggests that, for~$x:=0$, the logarithm thereof is to the leading order the maximum of the random walk~$S_n$ associated with the concentric decomposition. Upon scaling by~$\sqrt n$, this maximum tends in law to absolute value of a normal random variable. The question is:

\nopagebreak
\begin{figure}[t]
\vglue-1mm
\centerline{\includegraphics[width=0.45\textwidth]{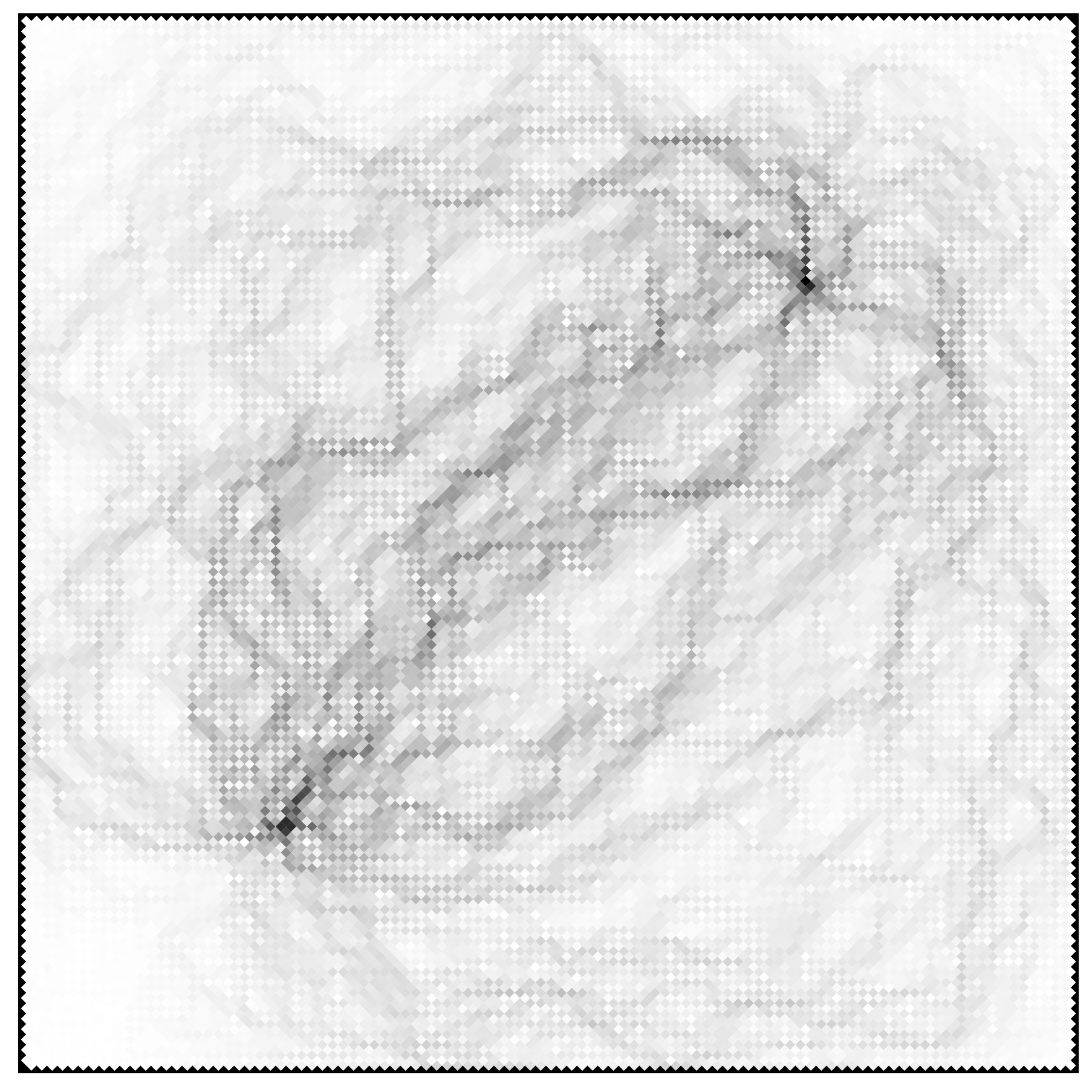}
\includegraphics[width=0.45\textwidth]{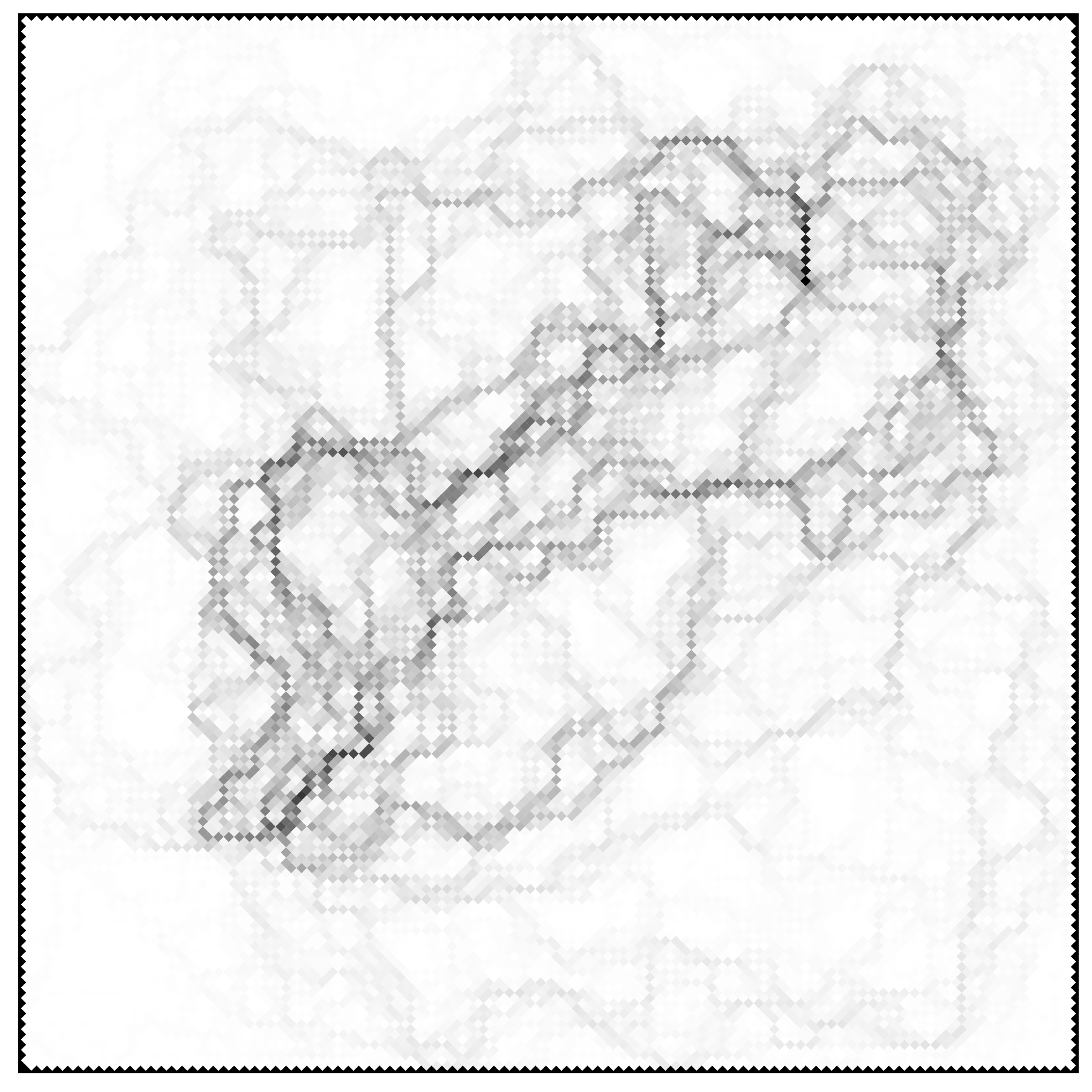}
}

\vglue0mm
\begin{quote}
\small 
\caption{
\label{fig-currents}
\small
Samples of the current minimizing the effective resistance between two points on the diagonal of the square (the dark spots in the left figure), for the network on the square with conductances related to the DGFF on the square as in \eqref{E:12.13} at inverse temperatures~$\beta=0.4\,\beta_\cc$ (left) and~$\beta=\beta_\cc$ (right). The shade of each pixel is proportional to intensity of the current. The current seems to be carried along well defined paths, particularly, in the figure on the right.}
\normalsize
\end{quote}
\end{figure}

\begin{myproblem}
For any sequence~$\{D_N\colon N\ge1\}$ of admissible approximations of domain~$D\in\mathfrak D$, characterize the limit distribution of
\begin{equation}
x\mapsto (\log N)^{-1/2}\log R_\eff(\lfloor x N\rfloor,D_N^\cc)
\end{equation}
 as $N\to\infty$.
\end{myproblem}

\noindent
In light of the concentric decomposition, a natural first guess for the limit process is the CGFF but proving this explicitly seems far from clear.
 
Next consider the problem of computing~$R_\eff(u,v)$ for the network in a square $[-N,N]^2\cap\Z^2$, the resistances associated to the DGFF via \eqref{E:12.13} and~$u$ and~$v$ two generic points in the square. The question we ask: Is there a way to describe the scaling limit of the minimizing current~$i_\star$? Tackling this directly may end up to be somewhat hard, because (as seen in Fig.~\ref{fig-currents}) the current is supported on a sparse (fractal) set.

A possible strategy to address this problem is based on the path-representation of the effective resistance from Proposition~\ref{prop-12.20}. Indeed, the proof contains an algorithm that identifies a set of simple paths~$\PP^\star$ from~$u$ to~$v$ and resistances $\{r^\star_{e,P}\colon e\in\cmss E,P\in\PP^\star\}$ that achieve the infima in \eqref{E:12.52}. Using these objects, $i_\star$ can be decomposed into a family of currents~$\{i_P\colon P\in\PP^\star\}$ from~$u$ to~$v$, such that~$i_P(e)=0$ for~$e\not\in P$ and with $i_P(e)$ equal to
\begin{equation}
i_P:=\frac{R_\eff(u,v)}{\sum_{e\in P}r_{e,P}^\star}
\end{equation}
on each edge~$e\in P$ oriented in the direction of~$P$ (recall that~$P$ is simple).
Noting that
\begin{equation}
\sum_{P\in\PP^\star}i_P=\text{val}(i_\star)=1\,,
\end{equation}
this recasts the computation of, say, the net current flowing through a linear segment~$[a,b]$ between two points $a,b\in D$ as the difference of the probability that a \emph{random path} crosses~$[a,b]$ in the given direction and the probability that it crosses $[a,b]$ in the opposite direction, in the distribution where path~$P$ (from~$u$ to~$v$) is chosen with probability~$i_P$. We pose:

\nopagebreak
\begin{figure}[t]
\vglue-1mm
\centerline{\includegraphics[width=0.45\textwidth]{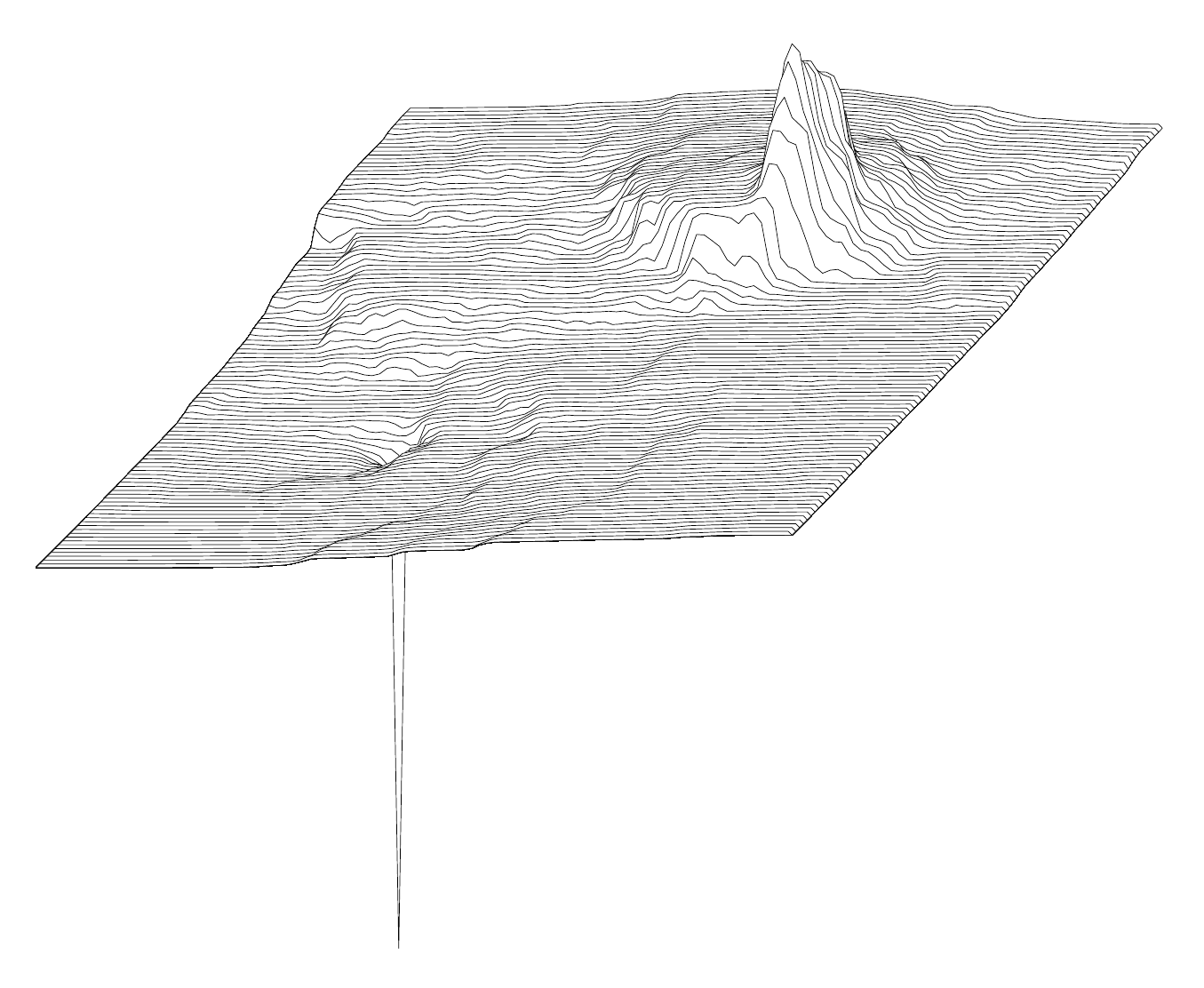}
\includegraphics[width=0.45\textwidth]{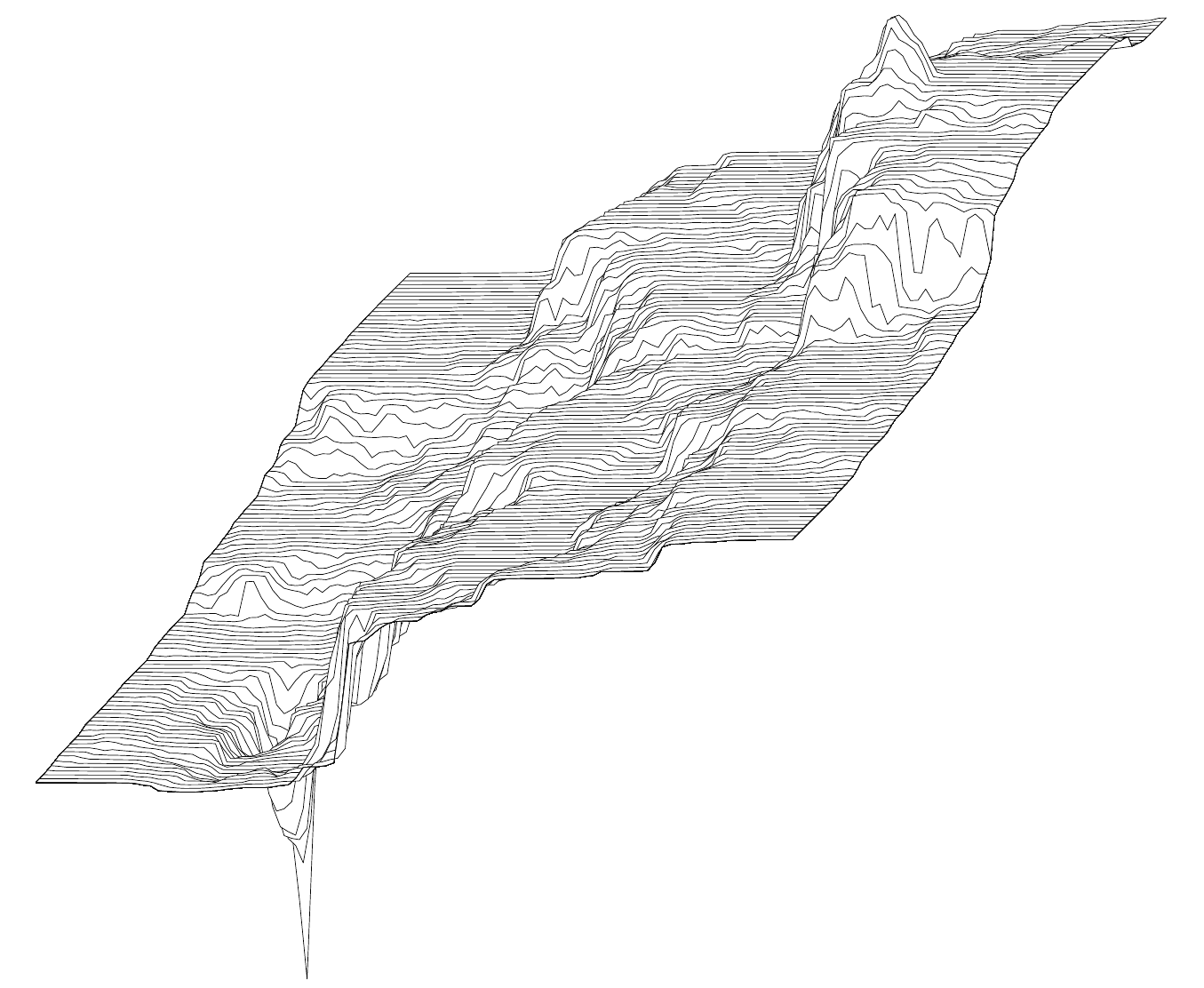}
}

\vglue0mm
\begin{quote}
\small 
\caption{
\label{fig-voltage}
\small
Two samples of the voltage realizing the effective conductances $C_\eff(u,v)$ in a box of side~$100$ and~$u$ and~$v$ two points on the diagonal. The conductances are related to the DGFF in the box via \eqref{E:12.13}.}
\normalsize±
\end{quote}
\end{figure}

\begin{myproblem}[Scaling limit of minimizing current]
\label{pb-15.12}
Prove that the joint law of the DGFF and the random path~$P$ admits a non-degenerate scaling limit as~$N\to\infty$.
\end{myproblem}

We envision tackling this using the strategy that we have used already a few times throughout these notes: First prove tightness, then extract subsequential limits and then characterize the limit law uniquely by, e.g., use of the Gibbs-Markov property and/or conformal maps. We suspect that the limit path law (in Problem~\ref{pb-15.12}) will be closely related to the scaling limit of the actual random walk process (in Problem~\ref{pb-15.11}); e.g., through a suitable loop erasure.

Once the continuum version of the above picture has been elucidated, the following additional questions (suggested by J.~Ding) come to mind:

\begin{myproblem}
What is the Hausdorff dimension of the (continuum) random path? And, say, for the electric current between opposite boundaries of a square box, what is the dimension of the support for the energy?
\end{myproblem}

\noindent
Similar geometric ideas can perhaps be used to address the computation of the voltage minimizing the effective conductance~$C_\eff(u,v)$; see Fig.~\ref{fig-voltage} for some illustrations. There collections of paths will be replaced by collections of nested cutsets (as output by the algorithm in the proof of Proposition~\ref{prop-12.22}). 

An alternative approach to the latter problem might thus be based on the observation that the cutsets are nested, and thus can be thought of as time snapshots of a \emph{growth process}. The use of Loewner evolution for the description of the kinematics of the process naturally springs to mind.

\newpage

\specialtocentry{References}

\mytocentry{\mytocfilestring}

\end{document}